\documentclass[oneside,12pt]{memoir}
\pdfoutput=1
 \let\subfloat\relax
\usepackage{graphicx,amsfonts, amssymb, amsthm, mathrsfs,  tikz, eucal, color,caption,float,enumitem, tabto,chngcntr,array,MnSymbol, rotating,tikz-cd,tabto,chngcntr,array,mathtools,caption, booktabs,footnote,subcaption,colortbl,hhline,breakcites}
\usepackage[style=alphabetic,backend=bibtex]{biblatex}
\bibliography{references}
\TabPositions{7.4cm}
\usetikzlibrary{decorations.shapes,positioning,arrows,calc,fit,shapes.geometric,decorations.pathmorphing}
\def\clap#1{\hbox to 0pt{\hss#1\hss}}

\renewbibmacro{in:}{}
 \restylefloat{figure}
 \counterwithin{table}{chapter}
 \counterwithin{equation}{chapter}
\theoremstyle{definition}
\newtheorem{theorem}{Theorem}  
\newtheorem{conjecture}[theorem]{Conjecture}
\newtheorem*{conjecture*}{Conjecture}
\newtheoremstyle{mytheorem}
{3pt}
{3pt}
{}
{}
{\bfseries}
{.}
{.5em}
{}

\theoremstyle{mytheorem}

\newtheorem{lemma}[theorem]{Lemma}
\newtheorem{proposition}[theorem]{Proposition}
\newtheorem{corollary}[theorem]{Corollary}

\newtheorem{fact}[theorem]{Fact}
\newtheorem{definition}[theorem]{Definition}
\newtheorem{remark}[theorem]{Remark}
\newtheorem{example}[theorem]{Example}

\newtheorem{convention}[theorem]{Convention}

\newcommand{\cat}[1]{\mathbf{#1}}
\DeclareMathAlphabet{\mathpzc}{OT1}{pzc}{m}{it}
\DeclareMathOperator{\Ty}{typ}
\DeclareMathOperator{\Fl}{Fl}
\DeclareMathOperator{\Orb}{Orb}
\DeclareMathOperator{\Sym}{Sym}

\DeclareMathOperator{\Fun}{Fun}
\DeclareMathOperator{\Hom}{Hom}
\newcommand{\op}{{\mathrm{op}}}
\DeclareMathOperator{\CalC}{\mathcal{C}}
\DeclareMathOperator{\CalD}{\mathcal{D}}

\numberwithin{theorem}{section}
\newcommand{\wh}{\widehat}

\newcommand{\image}{\mathrm{Im}\,}

\newcommand{\lra}{\longrightarrow}
\newcommand{\ra}{\rightarrow}

\newcommand{\Bic}{\mathrm{Bic}\,}
\newcommand{\Id}{\mathrm{Id}}
\newcommand{\ID}{\mathbb{I}\mathrm{d}}
\newcommand{\id}{\mathrm{id}}

\newcommand{\trhd}{\widetilde{\rhd}}
\newcommand{\tlhd}{\widetilde{\lhd}}
\newcommand{\lact}{\filledmedtriangleleft}
\newcommand{\colim}{\mbox{colim}}
\newcommand{\ract}{\filledmedtriangleright}
\newcommand{\uact}{\filledmedtriangleup}
\newcommand{\dact}{\filledmedtriangledown}
\newcommand{\cosk}{\mathbf{cosk}}
\newcommand{\tr}{\mathbf{tr}}
\newcommand{\Cosk}{\mathbf{Cosk}}
\newcommand{\sk}{\mathbf{sk}}
\newcommand{\Sk}{\mathbf{Sk}}
\newcommand{\D}{\mathcal{D}}
\newcommand{\C}{\mathcal{C}}
\newcommand{\K}{\mathcal{K}}
\newcommand{\Set}{\cat{Set}}

\newcommand{\Vdc}{\mathrm{Vdc}}
\newcommand{\ES}{\mathrm{ES}}
\newcommand{\B}{\mathcal{B}}
\newcommand{\BBT}{\mathbb{T}}

\newcommand{\n}{\mathbf{n}}
\newcommand{\m}{\mathbf{m}}
\newcommand{\e}{\mathbf{e}}
\newcommand{\fk}{\mathfrak{h}}
\newcommand{\frakf}{\mathfrak{f}}
\newcommand{\squish}{\mathrm{Sq}}
\newcommand{\fold}{\mathrm{Fold}}
\newcommand{\trunc}{\mathrm{Trunc}}
\newcommand{\Hcal}{\mathcal{H}}
\newcommand{\Vcal}{\mathcal{V}}
\newcommand{\T}{\mathcal{T}}
\newcommand{\blank}{\cellcolor{gray!60}}
\newcommand{\onecell}{\cellcolor{blue!15}}
\newcommand{\zzcell}{\cellcolor{yellow!27}}
\newcommand{\diff}{\mathrm{diff}}
\newcommand{\TBic}{\mathrm{FBic}}
\newcommand{\alddiv}{\!\! \upspoon \! \!}
\newcommand{\spoon}{\mathrm{sp}_f}
\newcommand{\vbar}{\ ; \ }
\newcommand{\fraks}{\mathfrak{s}}
\newcommand{\frakH}{\mathfrak{H}}
\newcommand{\frakh}{\mathfrak{h}}
\newcommand{\frakp}{\mathfrak{p}}
\newcommand{\frakq}{\mathfrak{q}}
\newcommand{\frakr}{\mathfrak{r}}
\newcommand{\frakt}{\mathfrak{t}}
\newcommand{\frakm}{\mathfrak{m}}

\newcommand{\TCat}{\cat{FCat}}
\newcommand{\wtb}{\widetilde{\B}}
\newcommand{\olb}{\overline{\B}}
\newcommand{\wtc}{\widetilde{\C}}
\newcommand{\olc}{\overline{\C}}
\newcommand{\A}{\mathcal{A}}
\newcommand{\wtf}{\widetilde{F}}
\newcommand{\olf}{\overline{F}}
\newcommand{\Vdcth}{\mathrm{TVdc}}
\newcommand{\SFun}{\mathrm{StFun}}
\newcommand{\sfancya}{f^{\uparrow}}
\newcommand{\cfancya}{f^{\downarrow}}

\newcommand{\fancyhom}{{\mathcal{H}\mbox{om}}}
\newcommand{\fancyadj}{\mathcal{A}\mbox{dj}}
\newcommand{\Adj}{\mathrm{Adj}}
\newcommand{\Adjst}{\mathrm{Adj}^{\mathrm{st}}}
\newcommand{\Sing}{\mathrm{Sing}}
\DefineBibliographyStrings{english}{%
bibliography = {References},
}

\newcolumntype{z}{@{\hskip\tabcolsep\vrule width 1.5pt\hskip\tabcolsep}}

\newcounter{nodemarkers}
\newcommand\circletext[1]{%
    \tikz[overlay,remember picture] 
        \node (marker-\arabic{nodemarkers}-a) at (0,1.5ex) {};%
    #1%
    \tikz[overlay,remember picture]
        \node (marker-\arabic{nodemarkers}-b) at (0,0){};%
    \tikz[overlay,remember picture,inner sep=2pt]
        \node[draw,ellipse,fit=(marker-\arabic{nodemarkers}-a.center) (marker-\arabic{nodemarkers}-b.center)] {};%
    \stepcounter{nodemarkers}%
}

\newcounter{sarrow}

\def\presuper#1#2%
  {\mathop{}%
   \mathopen{\vphantom{#2}}^{#1}%
   \kern-\scriptspace%
   #2}

\setlrmarginsandblock{1.0in}{1.0in}{*}
\setulmarginsandblock{1.6in}{1.0in}{*}

\checkandfixthelayout
\definecolor{comment}{rgb}{0,0,.73}

\begin{document}

\title{Non-Simplicial Nerves for Two-Dimen\-sional Categorical Structures}
\author{Nathaniel Watson\thanks{nathaniel.g.watson@gmail.com}}




\maketitle

\begin{abstract}


The most natural notion of a simplicial nerve for a (weak) bicategory was given by Duskin in \cite{Dus02}. Duskin showed that a simplicial set is isomorphic to the nerve of a $(2,1)$-category (i.e. a bicategory with invertible $2$-morphisms) if and only if it is a quasicategory which has unique fillers for inner horns of dimension $3$ and greater. Using Duskin's technique, we show how his nerve applies to $(2,1)$-category functors, making it a fully faithful inclusion of $(2,1)$-categories into simplicial sets. Then we consider analogues of this extension of Duskin's result for several different two-dimensional categorical structures, defining and analysing nerves valued in presheaf categories based on $\Delta^2$, on Segal's category $\Gamma$, and Joyal's category $\Theta_2$. In each case, our nerves yield exactly those presheaves meeting a certain ``horn-filling'' condition, with unique fillers for high-dimensional horns. Generalizing our definitions to higher dimensions and relaxing this uniqueness condition, we get proposed models for several different kinds higher-categorical structures, with each of these models closely analogous to quasicategories. Of particular interest, we conjecture that our ``inner-Kan $\Gamma$-sets'' are a combinatorial model for symmetric monoidal $(\infty,0)$-categories, i.e.  $E_\infty$-spaces.

\end{abstract}


\tableofcontents
\clearpage



\pagestyle{headings}

\chapter{Introduction}
The \emph{nerve} of a small category $\CalC$, first defined by Grothendieck in the 1960's, is a simplicial set whose $0$-cells are the objects of $\CalC$ and whose $n$-cells are composable sequences of $n$ morphisms in $\CalC$. The face and degeneracy maps of $N(\CalC)$ are defined by: 
\begin{align*}
d_0\left( c_0 \stackrel{f_1}{\lra} c_1 \stackrel{f_2}{\lra} c_2 \cdots c_{n-1}\stackrel{f_n}{\lra} c_n \right) &=  c_1 \stackrel{f_2}{\lra} c_2 \cdots c_{n-1}\stackrel{f_n}{\lra} c_n\\
d_i\left(  c_0 \stackrel{f_1}{\lra} c_1 \stackrel{f_2}{\lra} c_2 \cdots c_{n-1}\stackrel{f_n}{\lra} c_n \right) &=   c_0 \stackrel{f_1}{\lra} c_1\cdots  c_{i-1}\stackrel{f_{i+1}\circ f_i}{\lra} c_{i+1} \cdots c_{n-1}\stackrel{f_n}{\lra} c_n  \ \ \ \  \ \mbox{for $0<i<n$  } \\
d_n\left( c_0 \stackrel{f_1}{\lra} c_1 \stackrel{f_2}{\lra} c_2 \cdots c_{n-1}\stackrel{f_n}{\lra} c_n \right) &=  c_0 \stackrel{f_1}{\lra} c_1 \cdots c_{n-2}\stackrel{f_{n-1}}{\lra} c_{n-1}\\
s_i\left( c_0 \stackrel{f_1}{\lra} c_1 \stackrel{f_2}{\lra} c_2 \cdots c_{n-1}\stackrel{f_n}{\lra} c_n \right) &=   c_0 \stackrel{f_1}{\lra} c_1\cdots  c_{i-1}\stackrel{f_i}{\lra}c_{i}\stackrel{\id_{c_i}}{\lra} c_{i}\stackrel{f_{i+1}}{\lra}c_{i+1} \cdots c_{n-1}\stackrel{f_n}{\lra} c_n  
\end{align*}
The simplicial set $N(\CalC)$ always has the \emph{inner-Kan condition}, meaning whenever $0<i<n$, a map from the horn $\Lambda^n_i$ (which is obtained from removing the $i$th face from the simplex $\Delta[n]$) to $N(\CalC)$ always has a \emph{filler}, i.e. an extension along the inclusion $\Lambda^n_i \ra \Delta[n].$ According to \cite{Dus02}, it was Ross Street who first observed that these fillers are always unique for $N(\CalC)$, and if a simplicial set has unique fillers for all inner horns, it is the nerve of some small category.\footnote{A similar characterization of the nerve of a category goes back to Grothendieck \cite{Gro60}, however.}  The following slightly stronger statement summarizes Street's observation:
\begin{theorem}[Street] \label{catcase} $N$ is an equivalence of categories from the category of small categories and functors to the category of simplicial sets which have unique fillers for every inner horn.
\end{theorem} 
See for instance Chapter 4 of \cite{Dus02} for details. Boardman and Vogt in \cite{BV73} were the first to suggest studying simplicial sets satisfying the inner-Kan condition, calling them \emph{quasicategories}, and quasicategory theory has been advanced by Joyal (for instance in \cite{Joy02a} and \cite{Joy08}) and later by Lurie in \cite{Lur09}. Quasicategories are used as a model for $(\infty,1)$-categories, i.e. weak $\infty$-categories whose morphisms are invertible for $n>1.$

\subsection{Chapter~\ref{bicchapter}: The Duskin nerve} 
Many generalizations of Grothendieck's nerve have been given, but the nerve that is generalized in this thesis is the Duskin nerve, defined in \cite{Dus02}. This is a simplicial set $N(\B)$ defined from a small bicategory $\B$, with $N(\B)_0$ and $N(\B)_1$ respectively defined to be the objects and morphisms of $\B$, and a $2$-cell of $N(\B)$ defined to be a quadruple $(h,g,f \vbar \eta)$ where $h,g,f$ are $1$-morphisms of $\B$ and $\eta: g \Rightarrow h \circ f $ is a $2$-morphism. Details of the construction are provided in Section~\ref{duskinnerve}. Duskin establishes a characterization of those simplicial sets which are isomorphic to the Duskin nerve of a small bicategory in Theorem $8.6$ of \cite{Dus02}.

In case $\B$ is a small $(2,1)$-category, i.e. the $2$-morphisms of $\B$ are invertible, Duskin's characterization is very simple: a simplicial set $X$ is isomorphic to the nerve of a small $(2,1)$-category if it meets the inner-Kan condition, and every inner horn in $X$ of dimension $3$ or greater has a unique filler. We call a simplicial set meeting this condition a \emph{$2$-reduced inner-Kan simplicial set}. In Chapter~\ref{bictheorem}, we recapitulate Duskin's proof in this special case, leading to Theorem~\ref{daggertheorem}, which is part of Duskin's Theorem $8.6$. We then slightly extend Duskin's construction, defining the Duskin nerve for functors between bicategories, proving Theorem~\ref{bicsummary}, which shows that the Duskin nerve is an equivalence of categories, generalizing Theorem~\ref{catcase}.

\subsection{Overview of generalizations of the Duskin nerve}
Given the success of quasicategory theory, several generalizations and analogues of quasicategories have been defined. Notably, Dominic Verity's \emph{weak complicial sets}, defined in \cite{Ver}, are simplicial sets equipped with some extra structure and horn-filling conditions based on this structure, modelling general $\infty$-categories.

In this thesis, we consider three analogues of the inner-Kan condition, for categories of presheaves of sets on $\Delta\times \Delta$, Segal's category $\Gamma$, and Joyal's globular category $\Theta_2$. In each case, we show an equivalence of ``inner-Kan'' presheaves meeting a uniqueness condition for horn-fillers, which are called \emph{$2$-reduced inner-Kan presheaves}, with a certain algebraically defined categorical structure:
\begin{table}[H]
\begin{center}
\begin{tabular}{rlcll}
$N:$& $\! 2$-red. inner-Kan $\Delta\times\Delta$-sets &$\leftrightarrow$  & Verity double categories & $:\Vdc$\\
$\Orb \circ N:$& $\! 2$-red. inner-Kan $\Gamma$-sets             &$\leftrightarrow$ & Symmetric monoidal groupoids & $:\Sym\circ \Fl_2$ \\
$ N_\theta:$&$\! 2$-red. inner-Kan $\Theta_2$-sets	&	$\leftrightarrow$ 	& Fancy bicategories &  $:\TBic$
\end{tabular}\vspace{-.1in}\end{center}\end{table}
Symmetric monoidal groupoids are symmetric monoidal categories whose morphisms are invertible. See below in this introduction for discussion of Verity double categories and fancy bicategories.
\subsection{Chapter~\ref{vdcchapter}: Verity double categories and their bisimplicial nerves}
A \emph{double category} captures the notion of a category with two disjoint classes of $1$-morphisms between a shared class of objects. A \emph{strict double category} is  a category internal to the category $\cat{Cat}.$  The most common type of double category is called a \emph{pseudo-double category}, defined to be a category internal the strict bicategory $\cat{Cat}$, in the appropriate sense of being internal to a bicategory. A pseudo-double category consists of two categories $\mathcal{C}_1$, $\mathcal{C}_0$ with functors $s,t: \mathcal{C}_1\ra \mathcal{C}_0$ and $\id:\mathcal{C}_0 \ra \mathcal{C}_1$, together with a \emph{composition} functor: $$c:\mathcal{C}_1 \underset{\mathcal{C}_0}{\times}\mathcal{C}_1\ra \mathcal{C}_1$$ and a certain \emph{associator} natural transformation ensuring that this composition functor is associative up to isomorphism. We call the morphisms of $\mathcal{C}_0$ the \emph{vertical $1$-morphisms}, and the objects of $\mathcal{C}_1$ the \emph{horizontal $1$-morphisms}. Note that in this definition, the vertical $1$-morphisms have a strictly associative composition, whereas the horizontal $1$-morphisms have a composition which is only associative up to isomorphism. While pseudo-double categories are the most common type of double category in applications, this asymmetry complicates their relationship to bisimplicial sets.

Instead, we consider as our object of study in Chapter~\ref{vdcchapter} a notion of a double category which has weakly associative composition in both directions, defined by Verity in his 1992 Ph.D. thesis \cite{Ver11}. A Verity double category consists of two $(2,1)$-categories $H$ and $V$ sharing a set of objects $O$, and a set of squares $\mbox{Sq}(f,f',p,p')$ where $f,f'$ are $1$-morphisms in $H$, and  $p,p'$ are $1$-morphisms in $V$, with a square $\Theta$ pictured as shown:
\begin{center}
\begin{tikzpicture}[scale=1.8,auto]

\begin{scope}

\node (10) at (1,1) {};
\node (00) at (0,1) {};
\node (11) at (1,0) {};
\node (01) at (0,0) {};
\node[rotate=-45] at (.5,.5){$\Rightarrow$};
\node[scale=.8] at (.6,.65){$\Theta$};

\path[->] (00) edge node[midway]{$f$}(10);
\path[->] (00) edge node[midway,swap]{$p$}(01);
\path[->] (01) edge node[midway,swap]{$f'$}(11);
\path[->] (10) edge node[midway]{$p'$}(11);
\end{scope}

\end{tikzpicture}
\end{center}
Definition~\ref{Veritydef} provides a complete definition of Verity double categories. 

On the other hand, the generalization of the inner-Kan condition to bisimplicial sets is not difficult. In fact, Jardine defines in \cite{Jar10} a \emph{Kan condition} for a bisimplicial set by requiring a filler for every \emph{horn}, defined by removing any face from the boundary of a bisimplex, $d\Delta^2[m,n]$.  Each such face corresponds to a face either of the simplex $\Delta[m]$ or the simplex $\Delta[n]$, and we say a face is \emph{inner} if it is associated to an inner simplex face. We define an \emph{inner horn} to be a horn obtained by removing an inner face from $d\Delta^2[m,n]$, and we say a bisimplicial set $X$ is \emph{inner-Kan} if every inner horn in $X$ has a filler. If additionally such fillers are always unique when $m+n>2$, we say $X$ is \emph{$2$-reduced}.

In Chapter~\ref{vdcchapter}, we construct a bisimplicial nerve $N(\D)$ of a small Verity double category $\D$, which is a $2$-reduced inner-Kan bisimplicial set. Conversely, we show how a $2$-reduced inner-Kan bisimplicial set $X$ can be used to construct a small Verity double category $\Vdc(X)$, and we show these are constructions are in fact inverse equivalences of categories as the main theorem of Chapter~\ref{vdcchapter}, Theorem~\ref{vdcsummary}.

The crucial technique used in proving this theorem, as well as the main theorems of Chapter~\ref{gammachapter} and Chapter~\ref{thetachapter} is a generalized version of so-called ``Glenn tables'', which are extensively used by Duskin in \cite{Dus02} to construct horns and spheres in a simplicial set. Section~\ref{genglennsec} provides a setting in which this technique can be defined and justified.
\subsection{Chapter~\ref{bicatchapter}: Two bisimplicial nerves for fancy bicategories} 
Ehresmann first observed in \cite{Ehr631} that a strict double category can be constructed from a bicategory in two different ways. Either we can make a 2-category into a \emph{vertically trivial} strict double category having only identity vertical $1$-morphisms, or we can make a strict $2$-category $\mathcal{B}$ into double category $\ES(\B)$ whose vertical $1$-morphisms are identical to its horizontal $1$-morphisms, with both identical to the $1$-morphisms of $\B$. Then a square:
\begin{center}
\begin{tikzpicture}[scale=1.4,auto]
\begin{scope}
\node (10) at (1,1) {};
\node (00) at (0,1) {};
\node (11) at (1,0) {};
\node (01) at (0,0) {};
\node[rotate=45] at (.5,.5){$\Rightarrow$};
\path[->] (00) edge node[midway]{$f'$}(10);
\path[->] (00) edge node[midway,swap]{$f$}(01);
\path[->] (01) edge node[midway,swap]{$g$}(11);
\path[->] (10) edge node[midway]{$g'$}(11);
\end{scope}
\end{tikzpicture}
\end{center}
of $\ES(\mathcal{B})$ is a $2$-morphism $g \circ f \Rightarrow g'\circ f'$. 

Since Verity double categories are rare in nature, and bicategories are much more common, we consider in this chapter the extent to which the constructions above can be generalized, yielding a model for bicategories within the setting of Verity double categories. The natural generalization takes as input bicategories with some extra data:
\begin{definition} A \emph{fancy bicategory} $\B$ consists of a bicategory $\wtb$ together with a $(2,1)$-category $\olb$ and a strict functor $t_\B:  \olb \ra \wtb$ such that $t_\B$ is an isomorphism on objects and $1$-morphisms.
\end{definition}
We show there is a generalization of $\ES$ and of the vertically trivial construction, both making a Verity double category, and consider the structure and conditions on a Verity double category needed to invert this construction. We also consider the equivalent theory for bisimplicial sets, leading to two bisimplicial definitions of fancy bicategories.

\subsection{Chapter~\ref{gammachapter}: Symmetric monoidal groupoids and the $\Gamma$-nerve. }
In this chapter, we define and analyse a certain inner-Kan condition for $\Gamma$-sets, i.e. presheaves on Segal's category $\Gamma$. As before, the inner-Kan condition is a filler condition for \emph{inner horns}, which are subsheaves of representable presheaves obtained by removing certain faces. The appropriate definition of these inner horns is less obvious in this case, and is given in Definition~\ref{innerkangammadef}. We say an inner-Kan $\Gamma$-set $X$ is \emph{$2$-reduced} if it has unique fillers for inner horns of dimension $3$ or greater (as a note of caution, the appropriate general definition of an $n$-reduced inner-Kan $\Gamma$-set is slightly more complicated, see Definition~\ref{nreddef}).

We relate $2$-reduced inner-Kan $\Gamma$-sets to symmetric monoidal groupoid through an equivalent intermediate category. A monoidal groupoid (i.e. a monoidal category with invertible morphisms) is equivalent to a $(2,1)$-category with one object, or equivalently (as shown in Chapter~\ref{bicchapter}) a $2$-reduced inner-Kan simplicial set with one object. A \emph{$2$-reduced symmetric quasimonoid}, defined in Definition~\ref{symquasidef}, is a $2$-reduced inner-Kan simplicial set with one object, equipped with an involution $\sigma$  its $2$-cells which satisfies certain properties. 

In Section~\ref{orbsec}, we give an equivalence between symmetric quasimonoids and inner-Kan $\Gamma$-sets:
\begin{table}[H]
\begin{center}
\begin{tabular}{rlcll}
$\Orb:$&$2$-reduced symmetric quasimonoids &$\leftrightarrow$  & $2$-reduced inner-Kan $\Gamma$-sets & $:\Fl_2$
\end{tabular}\vspace{-.1in}\end{center}\end{table}
Then in Section~\ref{Symsection}  we show that if we take a small symmetric monoidal groupoid $\CalC$ viewed as a $(2,1)$-category with one object, then the Duskin nerve $N(\CalC)$ naturally has the structure of a $2$-reduced symmetric quasimonoid. Furthermore, if $X$ is a $2$-reduced symmetric quasimonoid then the associated  $(2,1)$-category $\Bic(X)$ can be given the structure of a symmetric monoidal category, which we call $\Sym(X)$. We show there is an equivalence:
\begin{table}[H]
\begin{center}
\begin{tabular}{rlcll}
$N:$&Small symmetric monoidal groupoids &$\leftrightarrow$  &$2$-red. symmetric quasimonoids  & $:\Sym$
\end{tabular}\vspace{-.1in}\end{center}\end{table}
Theorem~\ref{gammasummary}, noting that $\Orb \circ N$ and $\Sym\circ \Fl_2$ are inverse equivalences of categories, is the main result of this chapter.

\subsection{Chapter~\ref{thetachapter}: The globular nerve}
The \emph{globular category} $\Theta_n$ was defined by Joyal in an unpublished note \cite{Joy97}, which suggested using $\Theta_n$-sets meeting an inner-Kan condition as a definition of $n$-category. Leinster discusses this definition in \cite{Lei02}. The category $\Theta_n$ was fruitfully used by Rezk, who used certain $\Theta_n$-spaces (i.e. functors $\Theta_n^{\op} \ra \cat{Set}_\Delta$) to model $(\infty,n)$-categories in \cite{Rez10}, generalizing the model category of \emph{complete Segal spaces} as a model for $(\infty,1)$-categories, which is also due to Rezk (\cite{Rez97}).

We use Joyal's concept of an \emph{inner horn} and an \emph{inner-Kan} $\Theta_2$-set, defining a $\Theta_2$-set having unique fillers for inner horns of dimension $3$ and greater to be \emph{$2$-reduced inner Kan}. This $2$-reduced condition is similar in spirit but not exactly equivalent to the condition that Joyal proposes to put on inner-Kan $\Theta_2$-sets to model bicategories, but it the correct condition to use if we wish to get an equivalence to a non-strict $2$-dimensional categorical structure. 

We relate \emph{$2$-reduced inner-Kan} simplicial sets not to bicategories but to fancy bicategories, yielding an equivalence:
\begin{table}[H]
\begin{center}
\begin{tabular}{rlcll}
$ N_\theta:$&$2$-red. inner-Kan $\Theta_2$-sets	&	$\leftrightarrow$ 	& Small fancy bicategories &  $:\TBic$
\end{tabular}\vspace{-.1in}\end{center}\end{table}
Unlike in the other chapters, considerable work is necessary to show that the Glenn table technique works for $\Theta_2$-sets, and to this end a significant fraction of this chapter is devoted to combinatorial study $\Theta_2$.
\subsection{Chapter~\ref{epilogue}: Extensions and generalizations}
In the final chapter we consider inner-Kan $\Delta^n$-sets, inner-Kan $\Gamma$-sets, inner-Kan $\Theta_n$ sets, and inner-Kan $\Gamma \times \Delta^n$-sets. We suggest three ways of defining fancy $(\infty, n)$-categories and a way of defining fancy symmetric monoidal $(\infty,n)$-categories. 

If $\C$ is a small category, let $\cat{Set}_\C$ denote the category $\Fun(\C^{\op},\cat{Set})$ of presheaves of sets on $\C$. Similarly $\cat{Space}_\C$ denotes $\Fun(\C^{\op},\cat{Set}_\Delta)$, the category of $\C$-spaces, i.e. presheaves of simplicial sets on $\C$. We conjecture the existence of certain model structures on the presheaf categories $\cat{Set}_\Gamma$ and $\cat{Set}_{\Theta_n}$ which are conjecturally Quillen equivalent to known model structures on $\cat{Space}_\Gamma$ and $\cat{Space}_{\Theta_n}$. The $\Gamma$ case is particularly interesting:
\begin{conjecture*} There is a model structure $M_{\Gamma}$ on $\cat{Set}_\Gamma$ whose fibrant objects are the inner-Kan $\Gamma$-sets and which is Quillen equivalent to the model structure on $\cat{Space}_\Gamma$ whose fibrant objects are the special $\Gamma$-spaces, which was first defined in \cite{BF77}.
\end{conjecture*}
If this conjecture is true, $\Gamma$-sets provide a fairly simple combinatorial model of $E_\infty$-spaces, in the same way that quasicategories provide a simple combinatorial model for $(\infty,1)$-categories.

Since there are three combinatorial models for fancy bicategories in this thesis, we consider fancy bicategory theory in Section~\ref{fancysec}. After all, it is of little use to generalize fancy bicategories if they are themselves of no interest! We find that not only can bicategory concepts be extended to fancy bicategory theory, something is actually gained by making this generalization, as fancy bicategory theory unifies ``weak'' and ``strict'' notions in bicategory theory. As an illustration of this concept, we show that we can define a notion of $2$-limit for fancy bicategories that generalizes both the usual (weak) notion of a $2$-limit and the (fully) strict $2$-limit.

\subsection{Acknowledgements}
I wish to thank the Max Planck Institute for Mathematics for support and hospitality, and the Simons Foundation and the National Science Foundation for support during my graduate studies. Dan Berwick-Evans, Julie Bergner, and Chris Schommer-Pries provided comments and suggestions for this project. I also thank Constantin Teleman and Ori Ganor for serving on my committee and reviewing this thesis. I owe a special intellectual debt to John Duskin, since his bicategory nerve is the chief inspiration for this thesis. Special thanks is due to my advisor Peter Teichner, whose advice was invaluable. Without his encouragement, this project would never have been started, let alone finished. Finally, I hardly know how to thank my parents, Charles and Nancy Watson, enough. I hope they know how much I love and admire them.

\chapter{The bicategory $\mbox{Bic}(X)$ and the Duskin nerve\label{bicchapter}}
\begin{definition} An inner-Kan simplicial set will be called \emph{$2$-reduced} if it meets the inner-Kan condition uniquely for $n > 2$. \end{definition}

This chapter concerns two constructions due to Duskin in \cite{Dus02}. A $2$-reduced inner-Kan simplicial set  together with data $\chi$ specifying preferred fillers for $\Lambda^2_1$-horns in $X$ can be used to construct a $(2,1)$-category $\Bic(X,\chi)$ (i.e. a bicategory with invertible $2$-morphisms). In the other direction, a $(2,1)$-category $\B$ has a ``Duskin nerve'' $N(\B)$ which is a $2$-reduced inner-Kan simplicial set. Duskin showed these constructions are inverse to each other, thereby concluding that every $2$-reduced inner-Kan simplicial set is the Duskin nerve of a $(2,1)$-category. The main construction of this chapter takes a $2$-reduced inner-Kan simplicial set $X$ together with data $\chi$ specifying preferred fillers for $\Lambda^2_1$-horns in $X$ and gives a $(2,1)$-category $\Bic(X,\chi)$.

\begin{remark}In \cite{Dus02}, Duskin defines his $\Bic(X)$ for an inner-Kan simplicial set meeting a more general condition, but proves that, if $X$ is $2$-reduced, the $2$-morphisms of $\mbox{Bic}(X)$ are invertible. The construction of $\Bic(X,\chi)$ and the nerve $N$ in this section is merely a recapitulation of Duskin's work, with some shortcuts made possible by the fact we are working in this special case. They are included for the sake of making our exposition comprehensible and self-contained. The extension of $\Bic$ in $N$ to functors given in Sections~\ref{functorsection} and Section~\ref{summarysection} is anticipated by Duskin, but not explicated. However, in \cite{Gur09}, Gurski does show how $N$ may be extended to a full and faithful functor from the category of bicategories to the category of weak complicial sets, which are a kind of simplicial sets equipped with extra structure and properties. The constructions given in Sections~\ref{functorsection} and Section~\ref{summarysection} are essentially a special case of Gurski's construction.

\end{remark}

\section{Bicategories and $(2,1)$-categories}
Because we will build $\Bic(X)$ ``from scratch'' we will not be able to get around explicitly checking every bicategory axiom. For this reason, it will be helpful to write the axioms in a simple, explicit way. Using this definition sacrifices brevity and the heuristic force of more traditional definitions, but it will hopefully aid the reader in folowing the steps. These axioms are similiar to the axioms Duskin employs, in that they use ``whiskering'' instead of full horizontal composition of 2-morphisms, and this will simplify our definitions.  Since we are only interested in the case of a bicategory with invertible $2$-morphisms, we will include this assumption in the axioms, which slightly reduces their complexity.
\begin{definition} A $(2,1)$-category $\B$ consists of the following:

Data:
\begin{itemize}
\item A class of objects $O$. If $O$ is a set, then $\B$ is called \emph{small}
\item A set of \emph{$1$-morphisms} $\Hom_1(a,b)$ for any two objects 
\item For any $a,b$ and any morphisms $f,g \in \Hom_1(a,b)$, a set of \emph{$2$-morphisms} $\Hom_2(f,g)$
\end{itemize}
Structure:
\begin{enumerate}
\item For each object $a$, a \emph{pseudo-identity} $\id_{a}:a\ra a $ 
\item For objects $a,b,c$ and $f:a\ra b$ and $g:b \ra c$ a \emph{composite} $g \circ f: a \lra c$ 
\item For $f: a\ra b$ a $2$\emph{-identity}	 $\Id_f: f \Rightarrow f$ 
\item For $f,g,h: a \ra b$ and $\eta: f\Rightarrow g$ and  $\theta: g\Rightarrow h$ a \emph{vertical composite} $\theta \bullet \eta$ 
\item For $f,g: a\ra b$ and $h:b \ra c$ and $\eta:f \Rightarrow g$ a $\emph{right whiskering}$ $h \rhd \eta: f \circ h \Rightarrow g \circ h$
\item For $f: a \ra b$ and $g,h: b \ra c$ and  $\eta: g\Rightarrow h$ a $\emph{left whiskering}$ $\eta \lhd f:g \circ f \Rightarrow h \circ f$
\item For $f,g: a \ra b$ and $\eta: f\Rightarrow g$ an \emph{inverse}  $\eta^{-1} : g \Rightarrow f$
\item For $f:a\ra b,$ a \emph{right unitor} $\rho_f:f \Rightarrow f\circ \id_a $
\item For $f:a\ra b,$ a \emph{left unitor} $\lambda_f:f \Rightarrow\id_b \circ f$ 
\item For $f: a\ra b,$ and $g: b \ra c,$ and $h: c \ra d$ an \emph{associator} $\alpha_{h,g,f}: h \circ (g \circ f) \ra (h \circ g) \circ f$ 
\end{enumerate}
Axioms:
\begin{itemize}

\item \emph{category axioms for vertical composition}

\begin{enumerate}
\item For all $f,g$ and $\eta:f\Rightarrow g$,  \tab $\eta \bullet \Id_f =\Id_g \bullet \eta = \eta$.
\item For all $f \stackrel{\eta}{\Rightarrow} g \stackrel{\theta}{\Rightarrow} h \stackrel{\iota}{\Rightarrow}i,$ \tab $\iota \bullet (\theta \bullet \eta) = (\iota \bullet \theta) \bullet \eta.$
\item For all $f,g$ and $\eta:f \Rightarrow g$, \tab $\eta \bullet \eta^{-1}= \Id_g$ and $\eta^{-1} \bullet \eta= \Id_f$ 
\end{enumerate}

\item \emph{interchange of whiskering and vertical composition and identity}
\begin{enumerate}[resume]
\item  For all $a \stackrel{f}\ra b \stackrel{g}\to c$, \tab $g \rhd \Id_f=\Id_g \lhd f=\Id_{g \circ f }$
\item  For all $f \stackrel{\eta}\Rightarrow g \stackrel{\theta}\Rightarrow h: a \to b$ and $i:b \to c$,  \tab $(i \rhd \theta) \bullet (i \rhd \eta)=i \rhd (\theta \bullet \eta)$
\item  For all $f: a \to b$ and $g \stackrel{\eta}\Rightarrow h \stackrel{\theta}\Rightarrow i: b \to c$, \tab $(\theta \lhd f) \bullet (\eta \lhd f)=(\theta \bullet \eta) \lhd f$
\end{enumerate}
\item \emph{naturality of the unitors}
\begin{enumerate}[resume]
\item For all $f\stackrel{\eta}{\Rightarrow}g:a\ra b$ \tab  $ (\eta \lhd \id_a) \bullet \rho_f =\rho_g \bullet\eta $
\item For all $f\stackrel{\eta}{\Rightarrow}g:a \ra b$ \tab $  (\id_b \rhd \eta) \bullet\lambda_f=   \lambda_g\bullet\eta$

\end{enumerate}
\item\emph{naturality of the associator}
\begin{enumerate}[resume]
\item For all $ f\stackrel{\eta}{\Rightarrow} g:a \to b$ and $b \stackrel{h}\to c \stackrel{i}\to d$, \tab $\alpha_{i,h,g} \bullet (i \rhd (h \rhd \eta))=((i \circ h) \rhd \eta) \bullet \alpha_{i,h,f}$ 
\item For all $a\stackrel{f}{\ra}b$ and $g\stackrel{\eta}{\Rightarrow}h:b\ra c$ and $c\stackrel{i}{\ra}d,$ \tab $\alpha_{i,h,f} \bullet (i \rhd (\eta \lhd f))=((i \rhd \eta) \lhd f) \bullet \alpha_{i,g,f}$
\item For all $a\stackrel{f}{\ra}b\stackrel{g}{\ra}c$ and $h\stackrel{\eta}{\Rightarrow}i:c\ra d,$ \tab $\alpha_{i,g,f} \bullet (\eta \lhd (g \circ f))=((\eta \lhd g) \lhd f) \bullet \alpha_{h,g,f}$
\end{enumerate}
\item\emph{compatibility of the unitors and the pseudo-identity}
\begin{enumerate}[resume]
\item For all $a$ \tab $\lambda_{\id_a}=\rho_{\id_a}$
\end{enumerate}
\item\emph{compatibility of the unitors and the associator}
\begin{enumerate}[resume]
\item For all $a\stackrel{f}{\ra}b\stackrel{g}{\ra}c$  \tab $\alpha_{g,f,\id_a}\bullet(g\rhd 
\rho_f)=\rho_{g\circ f}$
\item For all $a\stackrel{f}{\ra}b\stackrel{g}{\ra}c$  \tab $\rho_g \lhd f  =\alpha_{g,\id_b,f}\bullet(g\rhd 
\lambda_f)$
\item For all $a\stackrel{f}{\ra}b\stackrel{g}{\ra}c$ \tab $\lambda_g \lhd f  =\alpha_{\id_c,g,f}\bullet\lambda_{g\circ f}$
\end{enumerate}
\item\emph{full interchange}
\begin{enumerate}[resume]
\item For all $f\stackrel{\eta}{\Rightarrow} g:a\ra b$ and $h\stackrel{\theta}{\Rightarrow} i:b\ra c$ \tab $(i\rhd \eta) \bullet (\theta \lhd f)=(\theta \lhd g)\bullet (h \rhd \eta)$
\end{enumerate}
\item\emph{pentagon identity}
\begin{enumerate}[resume]
\item For all $a \stackrel{f}\to b \stackrel{g}\to c \stackrel{h}\to d \stackrel{i}\to e$, \tab $(\alpha_{i,h,g} \lhd f) \bullet (\alpha_{i,h \circ g,f} \bullet (i \rhd \alpha_{h,g,f}))$ \\  ~ \tab \tab $\quad =\alpha_{i \circ h,g,f}\bullet \alpha_{i,h,g \circ f} $ 
\end{enumerate} 
\end{itemize}


\end{definition}

\begin{definition} A \emph{bicategory} consists of the same data and structure satisfying the same axioms, except with Structure 7 (the inverse of a 2-morphism) and Axiom 3 omitted, and axioms added to assert the invertibility of the unitors and the associator.
\end{definition}

\section{Glenn tables}

Following Duskin, we will use the simplicial table technique of Glenn described in \cite{Gle82}. An $d\Delta[n]$-\emph{sphere} in a simplicial set $Y$ will mean a map $d( \Delta[n]) \ra Y$, or equivalently a list $[y_0, \ldots, y_{n}]$  of $(n-1)$-simplices in $y$ that meet the  face relations they would need to have if they were the (ordered) boundary faces of a $n$-cell of $Y.$ \footnote{A is $d\Delta[n]$-sphere is usually called an $(n-1)$-sphere, since the topological realization $|d\Delta[n]|$ an $(n-1)$-sphere. We use our more cumbersome notation to avoid confusion arising from the fact that a $d\Delta[n]$-sphere is also a $n$-horn in the generalized sense of Definition~\ref{mincompdimdef}.   } This relation is simply \begin{equation}\label{simp} d_j(y_i)= d_{i}(y_{j+1}), \ \ 0 \leq i \leq j \leq {n-1}.\end{equation} To write down a specific $d\Delta[n]$-sphere $(y_0, \ldots, y_{n})$ in $Y$, we write a $(n+1) \times n$ table, where we list the $(n-2)$-cell $d_j(y_i)$ in the $i$th row, $j$th column position. Then condition ($\ref{simp}$) is easily verified by checking that the $(i,j)$ entry of this table on or above the main diagonal matches the $(j+1,i)$ entry below the main diagonal. Visually, this corresponds to checking the table is ``almost-symmetric'' in that the triangle of above-diagonal entries of the table can be flipped to correspond to the triangle of below-diagonal entries. This is the chief purpose of Glenn tables, to facilitate a quick visual check that a given list of cells meets the condition ($\ref{simp}$) which must hold in order for the cells to fit together as faces of a sphere.

In Section~\ref{genglennsec}, a generalization of the Glenn table technique is described. Readers seeking a more rigorous discussion of Glenn tables are referred to this section.

By itself this table only conveys the names of the $(n-2)$-cells of our $d\Delta[n]$-sphere, so we will give the names of each $(n-1)$-simplex in our sphere next to its corresponding row. An example of a sphere in $Y$ is pictured in Figure \ref{Ysphere}:
\begin{figure}[H]
\begin{center}
\begin{tikzpicture}[scale=1.8,auto]

\begin{scope}[shift=(54:3)]

\node (00) at (0,0) {$y_0$};
\node (10) at (90:1.5) {$y_1$};
\node (20) at (210:1.5) {$y_2$};
\node (30) at (330:1.5) {$y_3$};

\node[shape=coordinate, shift={(90:.3)}] (10out) at (10){};
\node[shape=coordinate, shift={(210:.3)}] (20out) at (20){};
\node[shape= coordinate, shift={(330:.3)}] (30out) at (30){};

\coordinate[shift={(255:.4)}](10e2) at (10){};
\coordinate[shift={(285:.4)}](10e3) at (10){};
\coordinate[shift={(45:.4)}](20e1) at (20){};
\coordinate[shift={(15:.4)}](20e3) at (20){};
\coordinate[shift={(165:.4)}](30e2) at (30){};
\coordinate[shift={(135:.4)}](30e1) at (30){};

\coordinate[shift={(150:.15)}](00up) at (00){};
\coordinate[shift={(270:.15)}](00left) at (00){};
\coordinate[shift={(30:.15)}](00right) at (00){};

\filldraw[gray!40] (10e2)--(20e1)--(00up)--(10e2);
\filldraw[blue!20] (10e3)--(30e1)--(00right)--(10e3);
\filldraw[yellow!35] (20e3)--(30e2)--(00left)--(20e3);
%
\draw[brown] (10out) to node[black,swap]{$y_{12}$} (20out);
\draw[brown](20out) to node[black,swap]{$y_{23}$} (30out);
\draw[brown] (30out) to node[black,swap]{$y_{13}$}(10out);

\path[->] (00) edge node[sloped]{$y_{02}$}(20);
\path[->] (00) edge node[midway,shift={(270:.3)}]{$y_{01}$} (10);  
\path[->]  (00) edge node[sloped,shift={(180:.6)}]{$y_{03}$} (30);

\path[->]
(10) edge (20)
		 edge (30)
(20) edge (30);

\node[shift=(150:.7)] at (0,0){$y_{012}$};
\node[shift=(270:.7)] at (0,0){$y_{023}$};
\node[shift=(45:.75)] at (0,0){$y_{013}$};

back lable 123
\node (back0) at (5:1.06){};
\node[shift=(30:.7)] (backlabel0) at (back0){$y_{123}$};
\draw [<-](back0) .. controls +(0:.3) and +(260:.3) .. (backlabel0);

\end{scope}

\end{tikzpicture}
\end{center}
\caption{A sphere in $Y$ with standard labels for the faces, edges, and vertices }\label{Ysphere}
\end{figure}
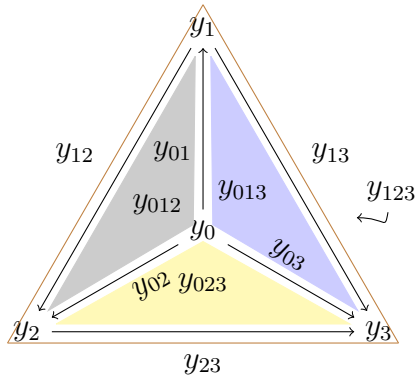

The sphere pictured in Figure \ref{Ysphere} has the following Glenn table:  

\begin{table}[H]\begin{center}
 \begin{tabular}{ r | l || l | l | l | }
     \cline{2-5}
        & $y_{123}$    &$y_{23}$      & $y_{13}$                   & $y_{12}$                    \\   \cline{2-5} 
         & $y_{023}$    &$y_{23}$      & $y_{03}$     & $y_{02}$         \\ \cline{2-5}
          & $y_{013}$  &$y_{13}$      & $y_{03}$     & $y_{01}$           \\   \cline{2-5}
            & $y_{012}$   &$y_{12}$      & $y_{02}$                   & $y_{01}$                    \\   \cline{2-5}
    \cline{2-5}
    \end{tabular}\end{center}
    \end{table}
The above table also illustrates a convention we will use for giving a naming $n$-cell in a simplicial set $Y$ and all its (iterated) faces: $y_{01\ldots n}$ denotes a generic $n$-cell, and its $i$th face is indicated by removing the $i$th index, e.g. $d_1 y_{01234}= y_{0234}$   and $d_2d_1 y_{01234}= y_{024}$. 

We will often wish to list the faces of a given $n$-cell, i.e. give the sphere that is its boundary. We denote this using the symbol $d$, for instance:
$$d y_{0123} = [y_{123},\ y_{023},\ y_{013},\ y_{012}].$$

Our main tool for this chapter will be horns, especially horns of the shape $\Lambda^3_1,\Lambda^3_2, \Lambda^4_1, \Lambda^4_2,$ and $\Lambda^4_3$-horns in a simplicial set also have simplicial tables.  We will list the faces of the ``missing face'' in the Glenn table of a horn as if it were present, but we will indicate which face is missing by a $\Lambda$ sign next to the appropriate row.
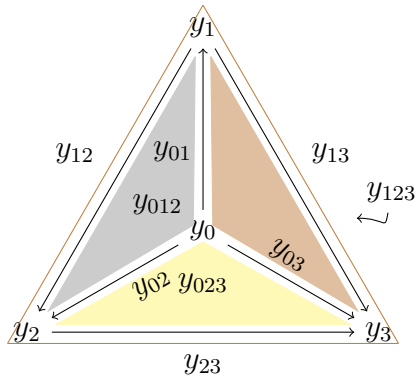
\begin{figure}[H]
\begin{center}
\begin{tikzpicture}[scale=1.8,auto]

\begin{scope}[shift=(54:3)]

\node (00) at (0,0) {$y_0$};
\node (10) at (90:1.5) {$y_1$};
\node (20) at (210:1.5) {$y_2$};
\node (30) at (330:1.5) {$y_3$};

\node[shape=coordinate, shift={(90:.3)}] (10out) at (10){};
\node[shape=coordinate, shift={(210:.3)}] (20out) at (20){};
\node[shape= coordinate, shift={(330:.3)}] (30out) at (30){};

\coordinate[shift={(255:.4)}](10e2) at (10){};
\coordinate[shift={(285:.4)}](10e3) at (10){};
\coordinate[shift={(45:.4)}](20e1) at (20){};
\coordinate[shift={(15:.4)}](20e3) at (20){};
\coordinate[shift={(165:.4)}](30e2) at (30){};
\coordinate[shift={(135:.4)}](30e1) at (30){};

\coordinate[shift={(150:.15)}](00up) at (00){};
\coordinate[shift={(270:.15)}](00left) at (00){};
\coordinate[shift={(30:.15)}](00right) at (00){};

\filldraw[gray!40] (10e2)--(20e1)--(00up)--(10e2);
\filldraw[brown!50] (10e3)--(30e1)--(00right)--(10e3);
\filldraw[yellow!35] (20e3)--(30e2)--(00left)--(20e3);
%
\draw[brown] (10out) to node[black,swap]{$y_{12}$} (20out);
\draw[brown](20out) to node[black,swap]{$y_{23}$} (30out);
\draw[brown] (30out) to node[black,swap]{$y_{13}$}(10out);

\path[->] (00) edge node[sloped]{$y_{02}$}(20);
\path[->] (00) edge node[midway,shift={(270:.3)}]{$y_{01}$} (10);  
\path[->]  (00) edge node[sloped,shift={(180:.6)}]{$y_{03}$} (30);

\path[->]
(10) edge (20)
		 edge (30)
(20) edge (30);

\node[shift=(150:.7)] at (0,0){$y_{012}$};
\node[shift=(270:.7)] at (0,0){$y_{023}$};

back lable 123
\node (back0) at (5:1.06){};
\node[shift=(30:.7)] (backlabel0) at (back0){$y_{123}$};
\draw [<-](back0) .. controls +(0:.3) and +(260:.3) .. (backlabel0);

\end{scope}

\end{tikzpicture}
\end{center} \caption{A labelled $\Lambda^3_2$-horn in $Y$.\label{lambda23horn}}
\end{figure}
The horn in Figure \ref{lambda23horn} has the following Glenn table: 

 \begin{table}[H] \begin{center}
 \begin{tabular}{ r | l || l | l || l | }
     \cline{2-5}
         & $y_{123}$    &$y_{23}$      & $y_{13}$                       & $y_{12}$                   \\   \cline{2-5} 
          & $y_{023}$  &$y_{23}$      & $y_{03}$                        & $y_{02}$         \\ \cline{2-5}
$\Lambda$    &  &$y_{13}$      & $y_{03}$                        & $y_{01}$         \\   \cline{2-5}
             & $y_{012}$      &$y_{12}$      & $y_{02}$                   & $y_{01}$                \\   \cline{2-5}
    \cline{2-5}
    \end{tabular}\end{center}
    
    \end{table}
When we define 3-horns and 4-horns in a $2$-reduced inner-Kan set, it will be in order to consider its unique filler. We will usually wish to give a name to the cells that fill the interior and empty face of our horn. We will put the name of the filling face in the empty space in the table, which is next to the row that names this face's edges. The name of the horn itself will be indicated if needed at the top of the table. We will also have a convention for naming the cell that fills such a horn:
\begin{convention} The name of the filling cell of a horn that has been given a name with a $\Lambda$ in it, (e.g. $\Lambda_{\wedge}\eta$) can be derived by changing the $\Lambda$ to a $\Delta$ (e.g. $\Delta_{\wedge}\eta$).
\end{convention}


\section{The data of $\mbox{Bic}(X)$}


Given only the simplicial set which is the Duskin nerve $N(\CalC)$ of a small bicategory $\CalC$ (see Section~\ref{duskinnerve}), we can recover the set of objects and morphisms of $\CalC$, but we can only tell the composition of two morphisms up to equivalence. This is because a $2$-cell in $N(\C)$ is a $2$-morphism $g \Rightarrow h \circ f$, but we can't tell using only the structure of $N(\CalC)$ whether or not a such a $2$-cell comes from an identity.
 
Instead, if we know the composition of $1$-morphisms in $\CalC$, we can use it to define a preferred filler of every $\Lambda_1^2$ in $N(\CalC)$, namely we fill \vspace{10pt}
 \begin{center}
 \begin{tabular}{ r| l ||  l | l | }
     \cline{2-4}
         & $x_{12}$    &$x_2$      & $x_1$                                \\   \cline{2-4} 
$\Lambda$  &       &$x_2$      & $x_0$                                    \\ \cline{2-4}
         & $x_{01}$         &$x_1$      & $x_0$                            \\   \cline{2-4}
   
    \cline{2-4}
    \end{tabular}
\end{center}\vspace{10pt}
by the identity morphism on the $1$-cell $x_{12}\circ x_{01}.$ 

This explains why the bicategory $\mbox{Bic}(X)$ we construct will depend not only on $X$ but on extra data $\chi$ consisting of a filler for every $\Lambda_1^2$ in $X$. A $\Lambda_1^2$-horn in $x$ is determined by two faces $f$ and $g,$ and we write $\chi(g, f)$ for our preferred filling 2-cell of this horn. 
\begin{definition}Following Nikolaus~\cite{Nik11} we call a inner-Kan simplicial set together with a set of preferred filling cells for every inner  horn a \emph{algebraic inner-Kan simplicial set}. A morphism $F:(X,\chi)\ra (Y,\chi')$ will be called \emph{strict} if $F(\chi(g,f))=\chi'(F(g),F(f))$ for all $1$-cells $a \stackrel{f}{\ra} b \stackrel{g}{\ra} c$ in $X.$ \end{definition}
For a $2$-reduced inner-Kan simplicial set, an algebraic structure $\chi$ is uniquely determined for all horns except $\Lambda_1^2$-horns, so we treat an algebraic  $2$-reduced inner-Kan simplicial set as having only this data. We define $\Bic(X)$ for an algebraic  $2$-reduced inner-Kan simplicial set, and if we wish to emphasize the dependence of $\mbox{Bic}$ on the algebraic data $\chi$, we will write $\mbox{Bic}(X,\chi).$

Let $(X,\chi)$ be an algebraic $2$-reduced inner-Kan simplicial set.

\begin{definition} The objects of $\mbox{Bic}(X)$ are the $0$-cells of $X$. A morphism from $a$ to $b$ in $\Bic(X)$ is a $1$-cells $f$ of $X$ with $df=[b,a]$.
\end{definition} 
\begin{definition} Let $f$ and $g$ be two $1$-cells in $X$ with $df=dg=[b,a]$. The 2-morphisms of $\Bic(X)$ with orientation $f\Rightarrow g$  are the $2$-cells $\eta$ of $X$ with faces $$d\eta=[g,f,s_0(a)].$$ 
\end{definition}

There is an arbitrary choice we have made in giving this definition, because we could have instead insisted that the zeroth face of $\eta$ is degenerate. We will have occasion to consider $2$-cells of this opposite type:
\begin{definition} \label{altmorphisms} Given $a,b$ and $f,g$ as above, an \emph{alt-2-morphism} in $X$ from  $f$ to $g$ is a $2$-cells $\eta$ of $X$ with faces  $$d\eta =[s_0(b),f,g].$$ 
\end{definition}

\begin{definition}\label{hat} Let $\eta$ be a 2-morphism from $f$ to $g$ as defined above. Then we can make the \emph{alt version} of $\eta$ by:

\begin{table}[H]\caption*{$\Lambda_{\wedge}(\eta)$, defining $\widehat{\eta}$}\begin{center}
\begin{tabular}{ r | l || l | l | l | }
     \cline{2-5}
          &  $s_1(g)$   &         $s_0(b)$      &  $g$             &    $g$                \\  \cline{2-5} 
   $\Lambda$   &   $=:\widehat{\eta}$     &  $s_0(b)$     &        $f$        &  $g$                       \\ \cline{2-5}
   &         $\eta$    &   $g$  & $f$               &     $s_0(a)$     \\   \cline{2-5}
       &     $s_0(g)$     &    $g$         &                  $g$       &  $s_0(a)$                       \\   \cline{2-5}
    \cline{2-5}
    \end{tabular} \end{center}
    \end{table}

If instead $\eta'$ is an alt-2-morphism with faces $d\eta'=[s_0(b),f,g]$ we define $\widehat{\eta'}$ by:

\begin{table}[H]\caption*{$\Lambda_{\wedge}(\eta')$, defining $\widehat{\eta'}$}\begin{center}
    \begin{tabular}{ r | l || l | l | l | }
     \cline{2-5}
       &  $s_1(g)$   &         $_0s(b)$      &  $g$             &       $g$                   \\  \cline{2-5} 
         &    $\eta'$        &         $s_0(b)$     &        $f$        &     $g$           \\ \cline{2-5}
$\Lambda$ & $=:\widehat{\eta'}$  &        $g$           & $f$               &     $s_0(a)$                 \\   \cline{2-5}
         &     $s_0(g)$      &         $g$         &      $g$       &  $s_0(a)$                    \\   \cline{2-5}
    \cline{2-5}
    \end{tabular}\end{center}
    \end{table}

\end{definition} 
\begin{lemma}\label{hatsonhats}If $\eta$ is a 2-morphism or alt-2-morphism, $\widehat{\widehat{\eta}}=\eta.$ This gives a bijection between 2-morphisms and alt-2-morphisms.
\end{lemma}
\begin{proof} $\Lambda_{\wedge}{\eta}$ is evidently a filler for the horn defining $\widehat{\widehat{\eta}}.$ By uniqueness of fillers, $\eta = \widehat{\widehat{\eta}}.$
\end{proof}


\section{The structure of $\Bic(X).$}


\begin{definition}[Structure 1 and 2] For $a\in \Bic(X),$ we define $\id_{a}=s_0(a).$ The composition of 1-morphisms in $\mbox{Bic}(X,\chi)$ is defined by $\chi$ in the obvious way: $$g \circ f := d_1( \chi(g, f))$$
\end{definition}

Often, we will have a $2$-cell $x$ with $dx =[h,g,f]$, which we wish to turn into a $2$-morphism $\underline{x}:g\Rightarrow h\circ f .$
\begin{definition}\label{underline}~

\begin{table}[H] \caption*{$\Lambda_{-} (x)$ }\begin{center}
    \begin{tabular}{ r | l || l | l | l | }
     \cline{2-5}

                          &$\chi(h, f)$       &   $h$   &     $h\circ f$    &  $f$          \\   \cline{2-5}

                              & $x$         &  $h$    & $g $                &  $f$             \\   \cline{2-5}
        
          $\Lambda$     &    $=:\underline{x}$     &   $h\circ f$   &   $g$    &    $\id_a$              \\ \cline{2-5}
           
                          & $s_0(f)$        &   $f$    & $f$        &  $\id_a$               \\   \cline{2-5} 
    \cline{2-5}
    \end{tabular}\end{center}
    \end{table}

\end{definition}

\begin{definition}[Structure 3 and 4] \label{bullet} We take $\Id_{f}= s_0(f).$  Let $f,g,h: a \ra b$ We define vertical composition $\theta \bullet \eta$ where $f\stackrel{\eta}{\Rightarrow}g$ and $g\stackrel{\theta}{\Rightarrow}f$ by:

\begin{table}[H]\caption*{$\Lambda_{ \bullet}(\theta,\eta)$, defining $\theta\bullet \eta$}\begin{center}
    \begin{tabular}{ r | l || l | l | l | }
     \cline{2-5}
      
         & $\theta$     &   $h$    & $g$                   &  $\id_a$                   \\   \cline{2-5}

$\Lambda$    & $=:\theta \bullet \eta$  &  $h$  & $f$     &    $\id_a$        \\   \cline{2-5}
     
         & $\eta$  &   $g$  & $f$     &    $\id_a$           \\ \cline{2-5}
     
        & $\Id_{\id_{a}}$      &$\id_a$      & $\id_a$                   & $\id_a$                   \\   \cline{2-5}

    \cline{2-5}
    \end{tabular}\end{center}
    \end{table}
\end{definition}

\begin{definition}[Structure 5 and 6]\label{triangle}  
For left whiskering,  suppose we have  $f:a\ra b,$ $g, h:b \ra c$ and $\eta: g\Rightarrow h.$ Then we can define $\wh{\eta\lhd f}$ with the following horn:

\begin{table}[H] \caption*{$\Lambda_{\lhd}(\eta,f)$}\begin{center}
    \begin{tabular}{ r | l || l | l || l | }
     \cline{2-5}

             & $\wh{\eta}$     &$\id_c$      & $g$                   & $h$               \\   \cline{2-5}

 $\Lambda$  &  $=:       \wh{\eta\lhd f}$  & $\id_c$    & $g \circ f $     & $h\circ f$                      \\ \cline{2-5}

            & $\chi(g, f)$    &$g$      & $g \circ f$     & $f$                           \\   \cline{2-5}  

            & $\chi(h, f)$   &$h$      &  $h\circ f$                  & $f$                \\   \cline{2-5}
 
    \cline{2-5}
    \end{tabular}\end{center}
    \end{table}
We of course define $\eta\lhd f:=\wh{\wh{\eta\lhd f}}$ using Lemma~\ref{hatsonhats} to justify our notation.

For right whiskering, suppose we have $f, g:a \ra b$ and $\eta: f\Rightarrow g$   and $h:b\ra c.$ We make the definition:

\begin{table}[H] \caption*{$\Lambda_{\rhd}( h,\eta)$}\begin{center}
    \begin{tabular}{ r | l || l | l | l | }
     \cline{2-5}
         
        & $\chi(h,g)$      &$h$            & $h\circ g$                   & $g$                  \\   \cline{2-5} 

            & $\chi(h, f)$   &$h$           & $h\circ f$       & $f$        \\   \cline{2-5}

$\Lambda$    & $=: h\rhd\eta$  &$h\circ g$            & $h\circ f$             & $\id_a$            \\ \cline{2-5}

             & $\eta$    &$g$      & $f$                   & $\id_a$                  \\   \cline{2-5}    \cline{2-5}
    \end{tabular}\end{center}
    \end{table}
\end{definition} 
\begin{definition}[Structure 7]\label{invdef} Let $f,g:a \ra b$ and $\eta: f\Rightarrow g.$ We define:

\begin{table}[H]\caption*{$\Lambda_{\mbox{inv}}(\eta)$, defining $\eta^{-1}$ in $\Bic(X)$}\begin{center}
    \begin{tabular}{ r | l || l | l | l | }
     \cline{2-5}

                  & $\eta$          &      $g$       & $f$                   &    $\id_a$              \\   \cline{2-5}        
  
                   &   $\Id_{g}$   &    $g$         & $g$     &  $\id_a$ \\   \cline{2-5}
  
  $\Lambda$       & $=:\eta^{-1}$    &   $f$    & $g$     &   $\id_a$      \\ \cline{2-5}
    
                   &  $\Id_{\id_{a}}$  &$\id_a$      & $\id_a$                   & $\id_a$                    \\   \cline{2-5} 
 
    \cline{2-5}
    \end{tabular}\end{center}
    \end{table}
\end{definition}

\begin{definition}[Structure 8 and 9] Let $f: a\ra b$. Then $$d(\underline{\Id_f})=[f \circ \id_a,f, \id_a]$$ which makes  $\underline{\Id_f}$ a 2-morphism $f \Rightarrow f \circ \id_a.$ So we define $\rho_f :=\underline{\Id_f}.$ Likewise $$d(\underline{(\wh{\Id_f})})=[\id_b \circ f, f, \id_a]$$ so $\underline{(\wh{\Id_f})}$ is a 2-morphism  $f \Rightarrow \id_b \circ f,$ and we define $\lambda_f=\underline{(\wh{\Id_f})}. $
\end{definition}
\begin{definition}[Structure 10] Let $f:a\ra b$ and $f: b \ra c$ and $f: c \ra d$. Then define:

\begin{table}[H]\caption*{$\Lambda_{\widetilde{\alpha}}(h,g,f)$, defining $\widetilde{\alpha}_{h,g,f}$}\begin{center}
    \begin{tabular}{ r | l || l | l | l | }
     \cline{2-5}
              &  $\chi(h,g)$            &$h$              & $h\circ g$                    & $g$                \\   \cline{2-5} 
          
                & $\chi(h, g \circ f)$      &$h$              & $h \circ (g \circ f) $        & {$g \circ f$}   \\ \cline{2-5}
     
         $\Lambda$  & $=:\widetilde{\alpha}_{h,g,f}$  &$h \circ g$      & $h \circ (g \circ f)$          & $f$            \\   \cline{2-5}
            
                     & $\chi(g,f)$        &$g$              & $g \circ f$                   & $f$                \\   \cline{2-5}
    \cline{2-5}
    \end{tabular}\end{center}
    \end{table}

\end{definition} 
Then we define $\alpha_{h,g,f}:= \underline{ \widetilde{\alpha}_{h,g,f}}.$

\section{Verification of the bicategory axioms for $\Bic(X)$.}

\subsection{Preliminaries}

As before, let $(X,\chi)$ be a algebraic $2$-reduced inner-Kan simplicial set.  

We begin by making setting some convenient terminology and giving a name to some obvious facts we will be repeatedly appealing to.
\begin{definition}An $d\Delta[n]$-sphere in a $X$ is said to be \emph{commutative} if it the boundary of some $n$-cell. 
\end{definition}
\begin{fact}[horn verification of commutativity] A $d\Delta[n]$-sphere in a inner-Kan complex is commutative if it forms the boundary of the missing face in an inner $n+1$-horn. 
\end{fact}
\begin{fact}[substitution] If $[a_0,a_1, \ldots, a_n]$ is an $d\Delta[n]$-sphere in a simplicial set, and $da_i= db,$ then $[a_0,a_1, \ldots,b,\ldots a_n]$ is also an $d\Delta[n]$-sphere.
\end{fact} 
\begin{lemma}[Matching Lemma] \label{lemma1} Suppose $[a_0,a_1,\ldots,a_n]$ is a commutative $d\Delta[n]$-sphere in $X$, with $n \geq 2$, and $d a_i = db$ with $0<i<n$. Then if $[a_0,a_1, \ldots, a_{i-1},b,a_{i+1},\ldots, a_n]$ is commutative, for instance, if it is the missing face of an inner $n+1$-horn, then $a_i=b.$
\end{lemma}
\begin{proof} By the uniqueness of fillers for the horn $[a_0,a_1,\ldots, a_{i-1},-,a_{i+1},\ldots,a_n]$.
\end{proof}

\begin{definition} If $\widehat{\eta}:f\Rightarrow g$ and $\widehat{\theta}:g\Rightarrow h$ are alt-2-morphisms, then we define their \emph{composition} $\widehat{\theta}\widehat{\bullet} \widehat{\eta}$ by the following horn:

\begin{table}[H]\caption*{$\Lambda_{ \widehat{\bullet}}(\widehat{\theta},\widehat{\eta})$}\begin{center}
    \begin{tabular}{ r | l || l | l | l | }
     \cline{2-5}
            & $\Id_{\id_{b}}$      &$\id_b$      & $\id_b$                   & $\id_b$              \\   \cline{2-5} 
         & $\widehat{\eta}$    &$\id_b$      & $f$     & {$g$}        \\ \cline{2-5}
$\Lambda$  & $=:\widehat{\theta} \widehat{\bullet} \widehat{\eta}$ &$\id_b$      & $f$     & $h$            \\   \cline{2-5}
          & $\widehat{\theta}$    &$\id_b$      & $g$                   & $h$                     \\   \cline{2-5}
    \cline{2-5}
    \end{tabular}\end{center}
    \end{table}

\end{definition}
\begin{lemma}\label{altbullet}Let $f\stackrel{\eta}{\Rightarrow}g\stackrel{\theta}{\Rightarrow}h$ be 2-morphisms, then $\widehat{\theta \bullet \eta}=\widehat{\theta}\widehat{\bullet}\widehat{\eta} $
\end{lemma}
\begin{proof} By the Matching Lemma, this equivalent to the statement that $$[\Id_{\id_b},\widehat{\eta},\widehat{\theta\bullet\eta},\widehat{\theta}]$$ is commutative. We show this with a two-step table proof. Consider the following Glenn table:

\begin{table}[H]\caption{~} \begin{center}
    \begin{tabular}{ r | l || l | l | l | l | }
     \cline{2-6}

     & $s_1(\Id_h)$    & $\Id_{\id_b}$   &    $\widehat{\Id_h}$     &  $\widehat{\Id_h}$            & $\widehat{\Id_h}$         \\ \cline{2-6}
     
   $\Lambda$     &      & $\Id_{\id_b}$  &  $\widehat{\eta}$   & $\widehat{\theta\bullet\eta}$   &       $\widehat{\theta}$      \\   \cline{2-6}  

     $\odot$     &  (Tabel~\ref{altbulletpt2})     &  $\widehat{\Id_h}$    &   $\widehat{\eta}$                     &    $\theta\bullet \eta$   &  $\theta$            \\   \cline{2-6}  
     
        &   $\Delta_{\wedge}(\theta\bullet\eta)$ & $\widehat{\Id_h}$&  $\widehat{\theta\bullet\eta}$      & $\theta\bullet\eta$        &  $\Id_h$ \\ \cline{2-6}

        &    $\Delta_{\wedge}(\theta)$ &  $\widehat{\Id_h}$ &   $\widehat{\theta}$ &     $\theta$    &   $\Id_h$      \\   \cline{2-6} 
    \end{tabular}\end{center}
    \end{table}
    
If we can show that the sphere which is marked $\odot$ is commutative, then if we fill it we get a horn verifying the desired commutativity. So it now suffices to show this face is commutative. This commutativity is verified by the following horn:
 
\begin{table}[H]\caption{~}\label{altbulletpt2} \begin{center}
    \begin{tabular}{ r | l || l | l | l | l | }
     \cline{2-6}

      & $s_2\theta$  & $\widehat{\Id_h}$   &    $\widehat{\Id_g}$     &  $\theta$            & $\theta$          \\ \cline{2-6}
     
   $\Lambda$  &       & $\widehat{\Id_h}$  &  $\widehat{\eta}$   & $\theta\bullet\eta$   &       $\theta$        \\   \cline{2-6}  

       & $ \Delta_{\wedge}(\eta)$    &  $\widehat{\Id_g}$  &   $\widehat{\eta}$                     &    $ \eta$   &  $\Id_g$           \\   \cline{2-6}  
     
      &   $\Delta_{\bullet}(\theta,\eta)$    & $\theta$&  $\theta\bullet\eta$      & $\eta$        &  $\Id_{\id_a}$  \\ \cline{2-6}

      &   $s_0\theta$   &  $\theta$ &   $\theta$ &     $\Id_g$    &  $\Id_{\id_a}$     \\   \cline{2-6} 
    \end{tabular}\end{center}  
    \end{table}
  
\end{proof}
\begin{remark}This method of proof will be needed to verify some of the axioms. We will mark faces whose commutativity must be verified with $\odot$, and refer to the table which verifies the commutativity in the left column. In this way, we will use a sequence of tables to prove the commutativity of a sphere.\end{remark}

We now define an alt version of $\Delta_{-}:$
\begin{definition} For any $2$-cell $x$ with $dx=[h,g ,f]$ define the alt-2-morphism $\overline{x}$ by:

\begin{table}[H] \caption*{$\Lambda_{\stackrel{\wedge}{-}} (x)$ }\begin{center}
    \begin{tabular}{ r | l || l | l | l | }
     \cline{2-5}

                             &$\widehat{\Id_h}$     &   $\id_c$   &     $h$    &  $h$           \\   \cline{2-5}
      
                $\Lambda$     &   $ =:\overline{x}$    &    $\id_c$   &   $g$    &  $h\circ f$                \\ \cline{2-5}

                               & $x$            &  $h$    & $g $                &  $f$                \\   \cline{2-5}
  
                                   & $\chi(h,f)$       &   $h$    & $h\circ f$        &  $f$              \\   \cline{2-5} 
    \cline{2-5}
    \end{tabular}\end{center}
    \end{table}

\end{definition}
\begin{lemma}\label{overlinelemma}For any $2$-cell $x$ in $X$ we have $\overline{x}=\widehat{(\underline{x})}$
\end{lemma}
\begin{proof} Let $dx=[h,g,f]$ as above. We have $$d\left(\Delta_{\wedge} (\underline{x})\right)=[s_1(h\circ f),\widehat{(\underline{x})}, \underline{x}, s_0(h\circ f)] $$ so by the Matching Lemma it suffices to show the sphere  $$[s_1(h\circ f),\overline{x}, \underline{x}, s_0(h\circ f)]$$ is commutative. The following horn verifies this commutativity:

\begin{table}[H] \begin{center}
    \begin{tabular}{ r | l || l | l | l | l | }
     \cline{2-6}

    & $s_2(\chi(h,f))$        & $s_1(h)$   &    $s_1(h\circ f)$     &  $\chi(h,f)$            & $\chi(h,f)$      \\ \cline{2-6}
     
     & $\Delta_{\stackrel{\wedge}{-}}(x) $    & $s_1(h)$  &  $\overline{x}$   & $x$   &       $\chi(h,f)$          \\   \cline{2-6}  

  $\Lambda$        &        &  $s_1(h\circ g)$  &   $\overline{x}$      &    $ \underline{x}$   &  $s_0(h\circ f)$       \\   \cline{2-6}  
     
    & $\Delta_{-}(x)$     & $\chi(h,f)$&  $x$      & $\underline{x}$        &  $s_0(f)$ \\ \cline{2-6}

    &   $s_0(\chi(h,f))$   &  $\chi(h,f)$ &   $\chi(h,f)$ &     $s_0(h\circ f)$    &  $s_0(f)$       \\   \cline{2-6} 
    \end{tabular}\end{center}
    \end{table}

\end{proof}

Up to now we have defined composition of 2-morphism and alt-2-morphisms. We now define an operation which generalizes composition of 2-morphisms, precomposing any 2-cell with a 2-morphism: 
\begin{definition}\label{bullety}Let $x$ be a 2-cell in $X$ with $dx= [i,  g,  h]$ and $\eta: f \rightarrow g,$ then $x \bullet \eta$ is the 2-cell defined by the following horn:

\begin{table}[H] \caption*{$\Lambda_{\bullet} (x,\eta)$ }\begin{center}
    \begin{tabular}{ r | l || l | l | l | }
     \cline{2-5}

                              &$x$   &   $i$   &     $g$    &  $h$            \\   \cline{2-5}
      
                $\Lambda$    &   $ =:x\bullet \eta$    &    $i$   &   $f$    &  $h$                 \\ \cline{2-5}

                                  & $\eta$            &  $g$    & $f $                &  $\id_a$             \\   \cline{2-5}
  
                                & $\Id_h$         &   $h$    & $h$        &  $\id_a$               \\   \cline{2-5} 
    \cline{2-5}
    \end{tabular}\end{center}
    \end{table}

\end{definition}
Note that this definition agrees with the previous one in the case that $x$ is a 2-morphism. 
\begin{lemma}\label{bulletylemma} With $x, \eta$ as above, $\underline{x \bullet \eta}= \underline{x} \bullet\eta$
\end{lemma}
\begin{proof} We have \begin{align*}d\left(\Delta_{-}(x\bullet\eta)\right) &=[\chi(i,h),x\bullet\eta,\underline{x\bullet\eta},\Id_h] \end{align*} so by the Matching Lemma it suffices to show $$[\chi(i,h),x\bullet\eta,\underline{x}\bullet\eta,\Id_h]$$ is commutative. The following horn verifies this commutativity:

\begin{table}[H] \begin{center}
    \begin{tabular}{ r | l || l | l | l | l | }
     \cline{2-6}

     & $\Delta_{-}(x)$    & $\chi(i,h)$   &    $x$     &  $\underline{x}$            & $\Id_h$         \\ \cline{2-6}
     
    $\Lambda$   &           & $\chi(i,h)$ &  $x\bullet \eta$   & $\underline{x}\bullet\eta$   &       $\Id_h$      \\   \cline{2-6}  

    &   $\Delta_{\bullet}(x,\eta)$    &  $x$  &   $x\bullet \eta$      &    $\eta$   &  $\Id_h$           \\   \cline{2-6}  
     
     &   $\Delta_{\bullet}(\underline{x},\eta)$     & $\underline{x}$&  $\underline{x}\bullet \eta$      & $\eta$        &  $\Id_{\id_a}$  \\ \cline{2-6}

     &   $s_1(h)$    &  $\Id_h$ &   $\Id_h$ &     $\Id_h$    &  $\Id_{\id_a}$     \\   \cline{2-6} 
    \end{tabular}\end{center}

    \end{table}
\end{proof}

\subsection{Category axioms for vertical composition}

\begin{proposition}[Axiom 1] For all $f,g: a\ra b$ in $\Bic(X)$ and $\eta:f\Rightarrow g$ we have  $\eta \bullet \Id_f=\eta$ and $\Id_g \bullet \eta = \eta$.
\end{proposition}
\begin{proof} First showing $\eta \bullet \Id_f=\eta$, recall from Definition~\ref{bullet}:
$$d\left(\Delta_{ \bullet}(\eta,\Id_f)\right) = [\eta, \eta \bullet \Id_f, \Id_f,\Id_{\id_a}].$$  On the other hand applying the simplicial identities $$d(s_0\eta)=[ \eta,\eta, s_0 d_1 \eta, s_0 d_2\eta ]=[\eta,\eta,\Id_f,\Id_{\id_a}].$$ Applying the Matching Lemma to these commutative spheres at face 1, we conclude $$\eta \bullet \Id_f=\eta.$$

Likewise for $\Id_g\bullet\eta$, we have $$d\left(\Lambda_{ \bullet}(\Id_g,\eta)\right)=[\Id_g,\Id_g \bullet \eta,\eta,\Id_{\id_{a}}]$$
Whereas $$d(s_1\eta) = [s_0 d_0\eta, \eta, \eta, s_1 d_2\eta]=[s_1\id_b, \eta,\eta,s_1 g]= [\Id_g,\eta,\eta,\Id_{\id_a}].$$ Applying the Matching Lemma to these commutative spheres at face 1, we see $\Id_g \bullet \eta=\eta$.\end{proof}

\begin{proposition}[Axiom 2]
 Let $f \stackrel{\eta}{\Rightarrow} g \stackrel{\theta}{\Rightarrow} h \stackrel{\iota}{\Rightarrow}i:a\ra b.$ Then $\iota \bullet (\theta \bullet \eta) = (\iota \bullet \theta) \bullet \eta.$\end{proposition}
\begin{proof}In order to check this, we consider  $$d\left(\Lambda_{\bullet}(\iota\bullet \theta,\eta)\right)= [\iota\bullet \theta,(\iota\bullet \theta)\bullet \eta,\eta,\Id_{\id_{a}}].$$ Then $$[\iota\bullet \theta,\iota\bullet(\theta\bullet\eta),\eta,\Id_{\id_{a}}]$$ is a sphere by substitution, with $d(\iota\bullet(\theta\bullet\eta))=d((\iota\bullet \theta)\bullet \eta)$ because both are 2-morphisms $f \Rightarrow i$.  Commutativity of this sphere is verified by filling the following horn, showing $$\iota\bullet(\theta\bullet\eta)=(\iota\bullet \theta)\bullet \eta$$ by the Matching Lemma:

\begin{table}[H]\begin{center}
    \begin{tabular}{ r | l || l | l | l | l | }
     \cline{2-6}

        & $\Delta_{\bullet}(\iota,\theta)$     &$\iota$                & $\iota\bullet \theta$               & $\theta$       &   $ \Id_{\id_{a}}$         \\ \cline{2-6}
    
       & $\Delta_{\bullet}(\iota,\theta\bullet \eta)$    &$\iota$               & $\iota\bullet (\theta\bullet \eta)$        & $\theta\bullet \eta$ &   $\Id_{\id_{a}} $  \\ \cline{2-6}
 
 $\Lambda$           &            &$\iota\bullet \theta$     & $\iota\bullet (\theta\bullet \eta)$      &$\eta$   &   $\Id_{\id_{a}} $              \\   \cline{2-6}
     
    &    $\Delta_{\bullet}(\theta,\eta)$      & $\theta$    &   $\theta\bullet \eta$                       &   $\eta$       &   $\Id_{\id_{a}}$                  \\   \cline{2-6}  
      
    & $s_0(\Id_{\id_{a}})$         &  $  \Id_{\id_{a}} $   &$\Id_{\id_{a}} $     & $\Id_{\id_{a}} $                & $\Id_{\id_{a}}  $            \\   \cline{2-6} 
    \end{tabular}\end{center}
    \end{table}

\end{proof}
\begin{proposition}[Axiom 3]Let $f,g:a\ra b$ and $\eta:f \Rightarrow g.$ Then $\eta \bullet \eta^{-1}= \Id_g$ and $\eta^{-1} \bullet \eta= \Id_f$ 
\end{proposition}
\begin{proof}Note $d\left(\Delta_{\bullet}(\eta,\eta^{-1})\right)=[\eta, \eta\bullet\eta^{-1},\eta^{-1},\Id_{\id_a}]$ and $$d\left(\Delta_{\mbox{inv}}(\eta)\right)=[\eta, \Id_g, \eta^{-1}, \Id_{\id_a}]$$
showing that $\eta \bullet \eta^{-1}=\Id_{g}$ by the Matching Lemma. To show the other identity, 
$$d\left( \Delta_{\bullet}(\eta^{-1},\eta)\right)= [\eta^{-1},\eta^{-1}\bullet \eta,\eta,\Id_{\id_{a}}],$$ so by substitution $ [\eta^{-1},\Id_{f},\eta,\Id_{\id_{a}}]$  is a sphere. Commutativity of this sphere is verified by the following horn, showing  $\eta^{-1} \bullet \eta =\Id_f$ by the Matching Lemma:

\begin{table}[H] \begin{center}
    \begin{tabular}{ r | l || l | l | l | l | }
     \cline{2-6}

              & $\Delta_{\mbox{inv}}(\eta)$    &$\eta$                & $\Id_g$               & $\eta^{-1}$       &   $ \Id_{\id_{a}}$        \\ \cline{2-6}
     
   (Axiom 1)   & $\Delta_{\bullet}(\eta,\Id_f)$   &$\eta$               & $\eta$        & $\Id_f$ &   $\Id_{\id_{a}} $  \\ \cline{2-6}
 
    (Axiom 1)   &   $\Delta_{\bullet}(\Id_g,\eta)$        &$\Id_g$     & $ \eta$      &$\eta$   &   $\Id_{\id_{a}} $                  \\   \cline{2-6}
     
   $\Lambda$    &          & $\eta^{-1}$    &   $\Id_f$                       &   $\eta$       &   $\Id_{\id_{a}}$           \\   \cline{2-6}  
      
          & $s_0(\Id_{\id_{a}})$  &  $  \Id_{\id_{a}} $   &$\Id_{\id_{a}} $     & $\Id_{\id_{a}} $                & $\Id_{\id_{a}}  $             \\   \cline{2-6} 
    \end{tabular}\end{center}
    \end{table}

Note that we have applied Axiom 1 to the marked rows, to show they are indeed the boundary of the named cells.
\end{proof}
\begin{lemma}
\label{lemma2} Let $a\stackrel{f}{\ra} b \stackrel{g}{\ra} c$, then $\underline{\chi(g,f)}=\Id_{g\circ f}.$
\end{lemma}
\begin{proof}From Definition~\ref{underline}, we see $$d\left( \Delta_{-}(\chi(g,f))\right)=[\chi(g,f),\chi(g,f),\underline{\chi(g,f)}, \Id_f].$$ On the other hand we can compute $$d\left( s_0\chi(g,f)\right)=[\chi(g,f),\chi(g,f),s_0d_1\chi(g,f),s_0d_2\chi(g,f)]=[\chi(g,f),\chi(g,f),\Id_{g\circ f},\Id_f]$$ so by the Matching Lemma, we conclude $\underline{\chi(g,f)}=\Id_{g\circ f}.$
\end{proof}
\begin{lemma}\label{lemma3} Let $f:a \ra b$ then $\widehat{\Id_f}=s_1f$.
\end{lemma}
\begin{proof} By Definition~\ref{hat}, we get $$d\left( \Delta_{\wedge}(\Id_f)\right)=[s_1f,\widehat{\Id_f},\Id_f,\Id_f],$$ whereas we can compute \begin{align*}d(s_0 s_1 f)&=[s_1 f,s_1 f, s_0 d_1 s_1 f,s_0 d_2 s_1 f] \\ &=[s_1 f,s_1 f,s_0 f,s_0 f]\\&=[s_1 f,s_1 f,\Id_f,\Id_f].\end{align*} We conclude $\widehat{\Id_f}=s_1f$ by the Matching Lemma.
\end{proof}

\subsection{Interchange laws}
\begin{proposition}[Axiom 4]Let $a\stackrel{f}{\ra} b \stackrel{g}{\ra} c$, then $g \rhd \Id_f=\Id_{g \circ f }$ and $\Id_g \lhd f=\Id_{g\circ f}.$
\end{proposition}
\begin{proof} First, $$d\left(\Delta_{\rhd}(g,\Id_f) \right)=[\chi(g,f), \chi(g,f),g \rhd \Id_f,\Id_f],$$ whereas we can compute \begin{align*} d\left(s_0\chi(g,f)\right)&=[\chi(g,f),\chi(g,f),s_0d_1\chi(g,f),s_0d_2\chi(g,f)] \\&=[\chi(g,f), \chi(g,f),\Id_{g\circ f},\Id_f].\end{align*} Then we apply the Matching Lemma to conclude $g \rhd \Id_f=\Id_{g\circ f}.$

For $\Id_g \lhd f $ by Lemma~\ref{hatsonhats} it suffices to show  $\wh{\Id_g \lhd f}=\wh{\Id_{g\circ f}}.$   $$d\left(\Delta_{\lhd}(\Id_g,f) \right)=[\wh{\Id_g}, \wh{\Id_{g}\lhd f},\chi(g,f),\chi(g,f)],$$ whereas \begin{align*} d\left(s_2\chi(g,f)\right)&=[s_1d_0\chi(g,f),s_1d_1\chi(g,f), \chi(g,f),\chi(g,f)] \\&=[\wh{\Id_g},\wh{\Id_{g\circ f}},\chi(g,f), \chi(g,f)].\end{align*} Note we have applied Lemma~\ref{lemma3} in the last step. By the  Matching Lemma, $$\wh{\Id_g \lhd f}=\wh{\Id_{g\circ f}}.$$ 
\end{proof}
\begin{proposition}[Axiom 5]Let $f \stackrel{\eta}\Rightarrow g \stackrel{\theta}\Rightarrow h: a \to b$ and $i:b \to c.$ Then $$(i \rhd \theta) \bullet (i \rhd \eta)=i \rhd (\theta \bullet \eta).$$
\end{proposition}
\begin{proof} $$d\left(\Delta_\bullet (i\rhd \theta, i \rhd \eta) \right)=[i\rhd \theta,\ (i\rhd \theta)\bullet (i\rhd \eta),\ i\rhd \eta,\ \Id_{\id_a}]$$ so by the Matching Lemma it suffices to show $$[i\rhd \theta,\ i\rhd (\theta \bullet \eta) ,\ i\rhd \eta,\ \Id_{\id_a}]$$ is commutative. The following horn proves this commutativity:

\begin{table}[H] \begin{center}
    \begin{tabular}{ r | l || l | l | l | l | }
     \cline{2-6}

        & $\Delta_{\rhd}(i,\theta)$    &$\chi(i,h)$   &  $\chi(i,g)$       & $i\rhd\theta$               & $\theta $         \\ \cline{2-6}
     
      &  $\Delta_{\rhd}(i,\theta\bullet\eta)$    &$\chi(i,h)$&   $\chi(i,f)$      & $i\rhd(\theta\bullet\eta)$        &  $\theta\bullet\eta $ \\ \cline{2-6}

       &  $\Delta_{\rhd}(i,\eta)$     &  $\chi(i,g)$    &   $\chi(i,f)$                      &   $i\rhd \eta$       &      $\eta$            \\   \cline{2-6}

  $\Lambda$        &        &  $i\rhd \theta$   &  $i\rhd(\theta\bullet\eta)$    & $i\rhd\eta $  & $\Id_{\id_a} $         \\   \cline{2-6}      
  
      &   $\Delta_{\bullet}(\theta,\eta)$    &$\theta $   &$\theta\bullet\eta $   &  $\eta$             &   $\Id_{\id_a} $     \\   \cline{2-6} 
    \end{tabular}\end{center}
    \end{table}

\end{proof}
\begin{proposition}[Axiom 6]Let $f: a \to b$ and $g \stackrel{\eta}\Rightarrow h \stackrel{\theta}\Rightarrow i: b \to c$. Then $$(\theta \lhd f) \bullet (\eta \lhd f)=(\theta \bullet \eta) \lhd f.$$
\end{proposition}
\begin{proof} 
It suffices to show $$\widehat{(\theta \lhd f) \bullet (\eta \lhd f)}=\wh{(\theta \lhd f)}\widehat{\bullet}\wh{(\eta \lhd f)}=\wh{(\theta \bullet \eta) \lhd f}.$$
We have
$$d\left(\Delta_{\widehat{\bullet}} (\wh{(\theta \lhd f)},\wh{(\eta \lhd f)})\right)=[\Id_{\id_c}, \ \wh{(\theta \lhd f)},\  \wh{(\eta \lhd f)}\widehat{\bullet} \wh{(\theta \lhd f)}   , \ \wh{(\eta \lhd f)}]$$
so by the Matching Lemma, it suffices to show $$[\Id_{\id_c}, \ \wh{(\eta \lhd f)},\  \wh{(\theta \bullet \eta) \lhd f}   , \ \wh{(\theta \lhd f)}]$$ is commutative. The following horn verifies this commutativity:

\begin{table}[H]\begin{center}
    \begin{tabular}{ r | l || l | l | l | l | }
     \cline{2-6}

   (Lemma~\ref{altbullet})      & $\Delta_{\widehat{\bullet}}(\wh{\theta},\wh{\eta})$     &$\Id_{\id_c}$                & $\wh{\eta}$               & $\wh{\theta\bullet\eta}$       &   $ \wh{\theta}$         \\ \cline{2-6}
     
$\Lambda$ &  &$\Id_{\id_c}$               & $\wh{(\eta \lhd f)}$        & $\wh{(\theta \bullet \eta) \lhd f}$ &   $\wh{(\theta \lhd f)} $   \\ \cline{2-6}
 
       &  $ \Delta_{\lhd}(\eta,f)  $     &$\wh{\eta}$     & $\wh{(\eta \lhd f)}$      &$\chi(g,f)$   &   $\chi(h,f)$                \\   \cline{2-6}
     
   &   $  \Delta_{\lhd}(\theta\bullet\eta,f)$        & $\wh{\theta\bullet\eta}$    &   $\wh{(\theta \bullet \eta) \lhd f}$                       &   $\chi(g,f)$       &   $\chi(i,f)$               \\   \cline{2-6}  
      
       &  $\Delta_{\lhd}(\theta,f) $   &   $ \wh{\theta}$    & $ \wh{(\theta \lhd f)}$       & $\chi(h,f) $                & $\chi(i,f)$         \\   \cline{2-6} 
    \end{tabular}\end{center}
    \end{table}

\end{proof}

\subsection{Naturality of unitors and the associator}

\begin{definition}
Let $f, g:a \ra b$ and $\eta: f\Rightarrow g$   and $h:b\ra c.$ We define $h \trhd \eta$ by:

\begin{table}[H] \caption*{$\Lambda_{\widetilde{\rhd}}( h,\eta)$}\begin{center}
    \begin{tabular}{ r | l || l | l | l | }
     \cline{2-5}
         
           & $\Id_{h}$         &$h$            & $h$                    & $\id_b$            \\   \cline{2-5} 

         & $\chi(h, f)$   &$h$      & $h\circ f$       & $f$           \\   \cline{2-5}

$\Lambda$ & $=: h\widetilde{\rhd}\eta$  &$h$            & $h\circ f$             & $g$               \\ \cline{2-5}

             & $\widehat{\eta}$     &$\id_b$      & $f$                   & $g$                 \\   \cline{2-5}    \cline{2-5}
    \end{tabular}\end{center}
    \end{table}

\end{definition}

Now let $f:a\ra b$ and  $g, h:b \ra c$ and $\eta: g\Rightarrow h.$  We define $\eta \widetilde{\lhd}f$ by:

\begin{table}[H] \caption*{$\Lambda_{\widetilde{\lhd}}(\eta,f)$}\begin{center}
    \begin{tabular}{ r | l || l | l | l | }
     \cline{2-5}

        & $\eta$      &$h$      & $g$                   & $\id_b$                   \\   \cline{2-5}

 $\Lambda$  & $=:\eta \tlhd f$  & {$h$}    & $g \circ f $     & $f$         \\ \cline{2-5}

           & $\chi(g, f)$  &$g$      & $g \circ f$     & $f$            \\   \cline{2-5}  

           & $\widehat{\Id_f}$      &$\id_b$      & $f$                   & $f$                 \\   \cline{2-5}
 
    \cline{2-5}
    \end{tabular}\end{center}
    \end{table}

\begin{lemma}\label{trianglelemma} $h \rhd \eta =\underline{h \widetilde{\rhd}\eta} $ and $\eta \lhd f=\underline{\eta \tlhd f}. $
\end{lemma}

\begin{proof} First for $h \rhd \eta =\underline{h \widetilde{\rhd}\eta} $ we have $$\Delta_{-}(h\trhd \eta)=[\chi(h,g),\ h\trhd \eta, \ \underline{h\trhd \eta},\ \Id_g]$$ so by the Matching Lemma it suffices to show $$[\chi(h,g),\ h\trhd \eta, \ h\rhd\eta,\ \Id_g]$$ is commutative. The following horn verifies this commutativity:

\begin{table}[H] \begin{center}

    \begin{tabular}{ r | l || l | l | l | l | }
     \cline{2-6}

     & $s_1(\chi(h,g))$    & $\Id_{h}$   &    $\chi(h,g)$     &$\chi(h,g)$           &$\wh{\Id_g}$      \\ \cline{2-6}
     
       &  $\Delta_{\trhd}(h,\eta)$    & $\Id_{h}$   &  $\chi(h,f)$   & $h\trhd \eta$   &     $\wh{\eta}$   \\   \cline{2-6}  

      & $\Delta_{\rhd}(h,\eta)$  &  $\chi(h,g)$ &   $\chi(h,f)$      &    $h\rhd \eta $   &  $\eta$       \\   \cline{2-6}  
     
   $\Lambda$    &    &$\chi(h,g)$  &   $h\trhd \eta $  &   $h\rhd \eta $    &  $\Id_g$  \\ \cline{2-6}

       & $\Delta_{\wedge}(\eta)$& $\wh{\Id_g}$ &  $\wh{\eta}$  &  $\eta$      &  $\Id_g$\\   \cline{2-6} 
    \end{tabular}\end{center}
    \end{table}

On the other hand, to show  $\eta \lhd f=\underline{\eta \tlhd f},$ by Lemma~\ref{overlinelemma} it suffices to show $\widehat{\eta \lhd f} =\overline{\eta \tlhd f}.$ We have  $$\Delta_{\stackrel{\wedge}{-}}( \eta\tlhd f)=[\wh{\Id_h},\ \overline{\eta\tlhd f}, \ \eta \tlhd f,\ \chi(h,f)]$$ so by the Matching Lemma it suffices to show  $$ [\wh{\Id_h},\ \wh{\eta\lhd f}, \ \eta \tlhd f,\ \chi(h,f)]$$ is commutative. The following horn verifies this commutativity:

\begin{table}[H] \begin{center}

    \begin{tabular}{ r | l || l | l | l | l | }
     \cline{2-6}

      &$\Delta_{\wedge}(\eta)$   & $\wh{\Id_{h}}$   &    $\wh{\eta}$     &$\eta$           &$\Id_h$     \\ \cline{2-6}
     
  $\Lambda$  &      & $\wh{\Id_{h}}$   &  $\wh{\eta\lhd f}$   & $\eta\tlhd f$   &     $\chi(h,f)$       \\   \cline{2-6}  

     & $\Delta_{\lhd}(\eta,f)$  &  $\wh{\eta}$ &   $\wh{\eta\lhd f}$       &    $\chi(g,f) $  & $\chi(h,f)$         \\   \cline{2-6}  
     
     &   $\Delta_{\tlhd}(\eta,f)$    &$\eta$  &   $\eta \tlhd f $  &   $\chi(g,f) $    &    $\wh{\Id_f}$\\ \cline{2-6}

    & $s_1(\chi(h,f))$   & $\Id_h$ &  $\chi(h,f)$   & $\chi(h,f)$    &  $\wh{\Id_f}$ \\   \cline{2-6} 
    \end{tabular}\end{center}
    \end{table}

\end{proof}

\begin{proposition}[Axiom 7]Let $\eta:f\Rightarrow g:a \ra b$ in $\Bic(X).$ Then $ (\eta \lhd \id_a) \bullet \rho_f =\rho_g \bullet\eta  .$\end{proposition}
\begin{proof}
We will show $$(\eta \lhd \id_a) \bullet \rho_f=\underline{\eta} =\rho_g \bullet\eta.$$ To show $\underline{\eta} =\rho_g \bullet\eta=\underline{\Id_g}\bullet \eta$ we apply Lemma~\ref{bulletylemma} which shows $\underline{\Id_g}\bullet \eta=\underline{\Id_g\bullet \eta}=\underline{\eta}$ by Axiom~2.

To show $(\eta \lhd \id_a) \bullet \rho_f=\underline{\eta}$ note that by Lemma~\ref{bulletylemma} we have $$\underline{(\eta \tlhd \id_a) \bullet \rho_f}=\underline{(\eta \tlhd \id_a)} \bullet \rho_f=(\eta \lhd \id_a) \bullet \rho_f$$ so it suffices to show  $(\eta \tlhd \id_a) \bullet \rho_f=\eta.$ We have $$d\left(\Delta_{\bullet}( \eta \tlhd \id_a,  \rho_f) \right)=[  \eta \tlhd \id_a, \  (\eta \tlhd \id_a)    \bullet  \rho_f   , \   \rho_f, \ \Id_{\id_a}       ]$$ so by the Matching Lemma it is enough to show $$[  \eta \tlhd \id_a, \  \eta   , \   \rho_f, \ \Id_{\id_a}       ]$$ is commutative. The following horn verifies this commutativity:

\begin{table}[H] \begin{center}

    \begin{tabular}{ r | l || l | l | l | l | }
     \cline{2-6}

     &$\Delta_{\tlhd}(\eta,\id_a)$   & $\eta$   &    $\eta\tlhd \id_a$     &$\chi(f,\id_a)$           &$\Id_{\id_a}$      \\ \cline{2-6}
     
    &  $s_0(\eta)$ &$\eta$   &  $\eta$   & $\Id_f$   &     $\Id_{\id_a}$        \\ \cline{2-6}  

   $\Lambda$     &     &  $\eta\tlhd \id_a$ &   $\eta$       &    $\rho_f= \underline{\Id_f} $  & $\Id_{\id_a}$        \\   \cline{2-6}  
     
     &   $\Delta_{-}(\Id_f)$   &$\chi(f,\id_a)$   &  $\Id_f$   &    $\underline{\Id_f} $    &    $\Id_{\id_a}$ \\ \cline{2-6}

    & $s_0(\Id_{\id_a})$    & $\Id_{\id_a}$ &  $\Id_{\id_a}$    & $\Id_{\id_a}$     & $\Id_{\id_a}$  \\   \cline{2-6} 
    \end{tabular}\end{center}
    \end{table}

\end{proof}

\begin{proposition}[Axiom 8]Let $\eta:f{\Rightarrow}g:a\ra b$ in $\Bic(X).$ Then $(\id_b \rhd \eta) \bullet \lambda_f=\lambda_g \bullet \eta.$
\end{proposition}
\begin{proof} We will show $(\id_b \rhd \eta) \bullet \lambda_f=\underline{(\wh{\eta})} = \lambda_g\bullet\eta.$ We have $$d\left(\Delta_{\bullet}(\lambda_g,\eta)\right)=[\lambda_g,\lambda_g\bullet\eta,\eta, \Id_{\id_a}]$$ so by the Matching Lemma we can show $\lambda_g\bullet\eta=\underline{(\wh{\eta})}$ by showing $$[\lambda_g,\underline{(\wh{\eta})},\eta, \Id_{\id_a}]$$ is commutative. The following horn verifies this:

\begin{table}[H] \begin{center}

    \begin{tabular}{ r | l || l | l | l | l | }
     \cline{2-6}

      &$\Delta_{-}(\wh{\Id_g})$    & $\chi(\id_b,g)$   &    $\wh{\Id_g}$     &$\underline{(\wh{\Id_g})}$           &$\Id_g$    \\ \cline{2-6}
     
   &  $\Delta_{-}(\wh{\eta})$   &$\chi(\id_b,g)$   &  $\wh{\eta}$   & $\underline{(\wh{\eta})}$   &     $\Id_{g}$      \\ \cline{2-6}  

  &  $\Delta_{\wedge}(\eta)$   &  $\wh{\Id_g}$ &   $\wh{\eta}$       &    $\eta$  & $\Id_{g}$         \\   \cline{2-6}  
     
     $\Lambda$    &       &$\lambda_g=\underline{(\wh{\Id_g})}$   &  $\underline{(\wh{\eta})}$ &    $\eta $    &    $\Id_{\id_a}$\\ \cline{2-6}

     & $s_1{\Id_g}$  & $\Id_g$ &  $\Id_g$    & $\Id_g$     & $\Id_{\id_a}$   \\   \cline{2-6} 
    \end{tabular}\end{center}
    \end{table}

To show the other equality,
$(\id_b \rhd \eta) \bullet \lambda_f=\underline{(\wh{\eta})}$ by Lemma~\ref{bulletylemma} and Lemma~\ref{trianglelemma} it is enough to show $(\id_b \trhd \eta) \bullet \lambda_f=\wh{\eta}.$ We have $$d\left(\Delta_{\bullet}(\id_b \trhd \eta,\lambda_f )\right)=[\id_b \trhd \eta,(\id_b \trhd \eta)\bullet\lambda_f,\lambda_f, \Id_g]$$ so by the Matching Lemma it suffices to show $$d\left(\Delta_{\bullet}(\id_b \trhd \eta,\lambda_f)\right)=[\id_b \trhd \eta,\wh{\eta},\lambda_f, \Id_g].$$ The following horn verifies this commutativity:

\begin{table}[H] \begin{center}

    \begin{tabular}{ r | l || l | l | l | l | }
     \cline{2-6}

      &$\Delta_{\trhd}(\id_b,\eta)$   & $\Id_{\id_b}$   &    $\chi(\id_b,f)$     &$\id_b \trhd \eta$           &$\wh{\eta}$     \\ \cline{2-6}
     
  &  $s_2(\wh{\eta})$   &$\Id_{\id_b}$   &  $\wh{\Id_f}$   & $\wh{\eta}$   &    $\wh{\eta}$       \\ \cline{2-6}  

   &  $\Delta_{-}(\wh{\Id_f})$   &  $\chi(\id_b,f)$ &   $\wh{\Id_f}$       &    $\underline{(\wh{\Id_f})}$  & $\Id_{f}$        \\   \cline{2-6}  
     
     $\Lambda$     &    &$\id_b \trhd \eta$   &  $\wh{\eta}$ &    $\lambda_f=\underline{(\wh{\Id_f})}$    &    $\Id_g$  \\ \cline{2-6}

     & $s_0(\wh{\eta})$     & $\wh{\eta}$& $\wh{\eta}$    &$\Id_{f}$   &$\Id_{g}$\\   \cline{2-6} 
    \end{tabular}\end{center}
    \end{table}

\end{proof}

\begin{proposition}[Axiom 9]Let $\eta: f \Rightarrow g:a \to b$ and $b \stackrel{h}\to c \stackrel{i}\to d$ in $\Bic(X).$ Then $$\alpha_{i,h,g} \bullet (i \rhd (h \rhd \eta))=((i \circ h) \rhd \eta) \bullet \alpha_{i,h,f}.$$
\end{proposition}

\begin{proof}
$$d\left(\Delta_{\bullet}((i\circ h)\rhd \eta,\alpha_{i,h,f})\right)=[(i\circ h)\rhd \eta, \ ((i\circ h)\rhd \eta) \bullet \alpha_{i,h,f} ,\ \alpha_{i,h,f},\ \Id_{\id_a}]$$
By the Matching Lemma, it suffices to show $$[(i\circ h)\rhd \eta, \ \alpha_{i,h,g} \bullet (i\rhd (h\rhd \eta)),\ \alpha_{i,h,f},\ \Id_{\id_a}]$$ is commutative. The following table proof verifies this commutativity:

\begin{table}[H] \caption{~ }\label{ax9pt1}\begin{center}

    \begin{tabular}{ r | l || l | l | l | l | }
     \cline{2-6}

      & $\Delta_{\rhd}(i\circ h, \eta) $ & $\chi(i\circ h, g)$   &    $\chi(i\circ h, f)$     &  $(i\circ h)\rhd \eta$            &$\eta$         \\ \cline{2-6}
     
   $\odot$   &  (Table~\ref{ax9pt2})  & $\chi(i\circ h, g)$    &  $\widetilde{\alpha}_{i,h,f}$   & $\alpha_{i,h,g} \bullet (i\rhd (h\rhd \eta))$   &     $\eta$         \\   \cline{2-6}  

   &   $\Delta_{-}(\widetilde{\alpha}_{i,h,f} )$   &   $\chi(i\circ h, f)$&  $\widetilde{\alpha}_{i,h,f}$     &    $\alpha_{i,h,f}$   &  $\Id_f$        \\   \cline{2-6}  
     
    $\Lambda$  &       & $(i\circ h)\rhd \eta$  & $\alpha_{i,h,g} \bullet (i\rhd (h\rhd \eta))$  &  $\alpha_{i,h,f}$     &  $\Id_{\id_a}$    \\ \cline{2-6}

      & $s_0\eta$  &$\eta$  &       $\eta$    &  $\Id_f$      &   $\Id_{\id_a}$ \\   \cline{2-6} 
    \end{tabular}\end{center}
    \end{table}

\begin{table}[H] \caption{~ }\label{ax9pt2}\begin{center}

    \begin{tabular}{ r | l || l | l | l | l | }
     \cline{2-6}

     & $\Delta_{-}(\widetilde{\alpha}_{i,h,g})$   & $\chi(i\circ h, g)$   &    $\widetilde{\alpha}_{i,h,g}$     &  $\alpha_{i,h,g}$            &$\Id_g$        \\ \cline{2-6}
     
   $\Lambda$     &            & $\chi(i\circ h, g)$    &  $\widetilde{\alpha}_{i,h,f}$   & $\alpha_{i,h,g} \bullet (i\rhd (h\rhd \eta))$   &     $\eta$     \\   \cline{2-6}  

    $\odot$     &   (Table~\ref{ax9pt3})     &  $\widetilde{\alpha}_{i,h,g}$ &  $\widetilde{\alpha}_{i,h,f}$     &    $i\rhd (h\rhd \eta)$    &  $\eta$         \\   \cline{2-6}  
     
  &   $\Delta_{\bullet}(\alpha_{i,h,g}, i\rhd (h\rhd \eta))$    &  $\alpha_{i,h,g}$   & $\alpha_{i,h,g} \bullet (i\rhd (h\rhd \eta))$  &  $i\rhd (h\rhd \eta)$     &  $\Id_{\id_a}$ \\ \cline{2-6}

      & $s_1\eta$       &$\Id_g$ &       $\eta$    &  $\eta$      &   $\Id_{\id_a}$ \\   \cline{2-6} 
    \end{tabular}\end{center}

    \end{table}

\begin{table}[H] \caption{~ }\label{ax9pt3}\begin{center}

    \begin{tabular}{ r | l || l | l | l | l | }
     \cline{2-6}

      & $\Delta_{\widetilde{\alpha}}(i,h,g)$   & $\chi(i, h)$   &    $\chi(i,h\circ g)$    &  $\widetilde{\alpha}_{i,h,g}$           &$\chi(h,g)$       \\ \cline{2-6}
     
        & $\Delta_{\widetilde{\alpha}}(i,h,f)$     & $\chi(i,h)$    &   $\chi(i,h\circ f)$   & $\widetilde{\alpha}_{i,h,f}$   &    $\chi(h,f)$      \\   \cline{2-6}  

&   $\Delta_{\rhd}(i, h\rhd \eta)$    &  $\chi(i,h\circ g)$  & $\chi(i,h\circ f)$  &  $i\rhd (h\rhd \eta)$     &  $h\rhd \eta$\\ \cline{2-6}

    $\Lambda$     &      &  $\widetilde{\alpha}_{i,h,g}$ &  $\widetilde{\alpha}_{i,h,f}$     &    $i\rhd (h\rhd \eta)$    &  $\eta$        \\   \cline{2-6}

        & $\Delta_{\rhd}(h, \eta)$    &$\chi(h,g)$ &       $\chi(h,f)$    &  $h\rhd \eta$      &   $\eta$  \\   \cline{2-6} 
    \end{tabular}\end{center}
    \end{table}
\end{proof}

\begin{proposition}[Axiom 10]Let $a\stackrel{f}{\ra}b$ and $g\stackrel{\eta}{\Rightarrow}h:b\ra c$ and $c\stackrel{i}{\ra}d$ be in $\Bic(X).$ Then  $$\alpha_{i,h,f} \bullet (i \rhd (\eta \lhd f))=((i \rhd \eta) \lhd f) \bullet \alpha_{i,g,f}.$$ \end{proposition}
\begin{proof} By Lemmas~\ref{bulletylemma} and \ref{trianglelemma},  \begin{align*}\alpha_{i,h,f} \bullet (i \rhd (\eta \lhd f))&=\underline{\widetilde{\alpha}_{i,h,f} \bullet (i \rhd (\eta \lhd f))}\\((i \rhd \eta) \lhd f) \bullet \alpha_{i,g,f} &= \underline{((i \rhd \eta) \tlhd f) \bullet \alpha_{i,g,f}}.\end{align*} So it suffices to show $$\widetilde{\alpha}_{i,h,f} \bullet (i \rhd (\eta \lhd f))=((i \rhd \eta) \tlhd f) \bullet \alpha_{i,g,f}.$$ We have
$$d\left(\Delta_{\bullet}((i \rhd \eta) \tlhd f,\alpha_{i,g,f})\right)=[(i \rhd \eta) \tlhd f,\ ((i \rhd \eta) \tlhd f) \bullet \alpha_{i,g,f}, \ \alpha_{i,g,f} , \ \Id_f],$$ so by the Matching Lemma, it suffices to show $$[(i \rhd \eta) \tlhd f,\ \widetilde{\alpha}_{i,h,f} \bullet (i \rhd (\eta \lhd f)), \ \alpha_{i,g,f} , \ \Id_f]$$ is commutative. The following table proof verifies this commutativity:

\begin{table}[H] \caption{~ }\label{ax10pt1}\begin{center}

    \begin{tabular}{ r | l || l | l | l | l | }
     \cline{2-6}

    & $\Delta_{\tlhd}(i \rhd \eta,f)$        & $i\rhd \eta$   &    $(i\rhd \eta)\tlhd f$     &  $\chi(i\circ g,f)$            & $\widehat{\Id_f}$      \\ \cline{2-6}
     
    $\odot$      &  (Table~\ref{ax10pt2})   & $i\rhd \eta$ &  $\widetilde{\alpha}_{i,h,f}\bullet(i\rhd(\eta \lhd f))$   & $\widetilde{\alpha}_{i,g,f}$   &       $\widehat{\Id_f}$      \\   \cline{2-6}  

   $\Lambda$   &       &  $(i\rhd \eta)\tlhd f$  &   $\widetilde{\alpha}_{i,h,f}\bullet (i\rhd(\eta \lhd f))$      &    $\alpha_{i,g,f}$   &  $\Id_f$           \\   \cline{2-6}  
     
      &   $\Delta_{-}(\widetilde{\alpha}_{i,g,f})$    & $\chi(i\circ g,f)$&  $\widetilde{\alpha}_{i,g,f}$      &$\alpha_{i,g,f}$        &  $\Id_f$  \\ \cline{2-6}

      &   $s_2 \Id_f$   & $\widehat{\Id_f}$ &   $\widehat{\Id_f}$ &     $\Id_f$    &  $\Id_f$     \\   \cline{2-6} 
    \end{tabular}\end{center}
    \end{table}

\begin{table}[H] \caption{~ }\label{ax10pt2}\begin{center}

    \begin{tabular}{ r | l || l | l | l | l | }
     \cline{2-6}

                        & $\Delta_{\rhd}(i,\eta)$     & $\chi(i,h)$   &    $\chi(i,g)$    &   $i\rhd \eta$           & $\eta$          \\ \cline{2-6}
     
       $\odot$   &  (Table~\ref{ax10pt3})   & $\chi(i,h)$ & $\chi(i,g\circ f)$     &  $\widetilde{\alpha}_{i,h,f}\bullet(i\rhd(\eta \lhd f))$    &      $\eta \tlhd f$        \\   \cline{2-6}  

    &   $\Delta_{\widetilde{\alpha}}(i,g,f)$    & $\chi(i,g)$ &    $\chi(i,g\circ f)$     &    $\widetilde{\alpha}_{i,g,f}$   &  $\chi(g,f)$          \\   \cline{2-6}  
     
    $\Lambda$  &    &  $i\rhd \eta$&  $\widetilde{\alpha}_{i,h,f}\bullet(i\rhd(\eta \lhd f))$      & $\widetilde{\alpha}_{i,g,f}$         &  $\Id_f$  \\ \cline{2-6}

    &   $\Delta_{\tlhd}(\eta,f)$   & $\eta$ &   $\eta\tlhd f$ & $\chi(g,f)$    &  $\Id_f$       \\   \cline{2-6} 
    \end{tabular}\end{center}

    \end{table}

\begin{table}[H] \caption{~ }\label{ax10pt3}\begin{center}

  \scalebox{0.97}{  \begin{tabular}{ r | l || l | l | l | l | }
     \cline{2-6}

                         & $\Delta_{\widetilde{\alpha}}(i,h,f)$        & $\chi(i,h)$   &    $\chi(i,h\circ f)$    &   $\widetilde{\alpha}_{i,h,f}$           & $\chi(h,f)$      \\ \cline{2-6}

    $\Lambda$      &        & $\chi(i,h)$ & $\chi(i,g\circ f)$     &  $\widetilde{\alpha}_{i,h,f}\bullet(i\rhd(\eta \lhd f))$    &      $\eta \tlhd f$     \\   \cline{2-6}  

     &   $\Delta_{\rhd}(i,\eta\lhd f)$  &$\chi(i,h\circ f)$ &    $\chi(i,g\circ f)$     &    $i\rhd(\eta\lhd f)$   &  $\eta \lhd f$           \\   \cline{2-6}  
    &   $\Delta_{\bullet}(\widetilde{\alpha}_{i,h,f},i\rhd(\eta\lhd f))$  &   $\widetilde{\alpha}_{i,h,f}$&  $\widetilde{\alpha}_{i,h,f}\bullet(i\rhd(\eta \lhd f))$      &  $i\rhd(\eta\lhd f)$          &  $\Id_f$ \\ \cline{2-6}

       &   $\Delta_{-}(\eta\tlhd f)$   & $\chi(h,f)$ &   $\eta\tlhd f$ & $\eta \lhd f$   &  $\Id_f$    \\   \cline{2-6} 
    \end{tabular}}\end{center}
    \end{table}
\end{proof}

\begin{proposition}[Axiom 11]Let $a\stackrel{f}{\ra}b\stackrel{g}{\ra}c$ and $h\stackrel{\eta}{\Rightarrow}i:c\ra d$ in $\Bic(X).$ Then $$\alpha_{i,g,f} \bullet (\eta \lhd (g \circ f))=((\eta \lhd g) \lhd f) \bullet \alpha_{h,g,f}$$
\end{proposition}
\begin{proof} This identity is equivalent to $$\wh{\alpha_{i,g,f}} \widehat{\bullet} \wh{(\eta \lhd (g \circ f))}=\wh{((\eta \lhd g) \lhd f)} \widehat{\bullet} \wh{\alpha_{h,g,f}}.$$ We have

$$d\left(\Delta_{\widehat{\bullet}}(\wh{(\eta\lhd g)\lhd f},\wh{\alpha_{h,g,f}})\right)=[\Id_{\id_d}, \ \wh{\alpha_{h,g,f}}  ,\ \wh{((\eta\lhd g)\lhd f)}\widehat{\bullet} \wh{\alpha_{h,g,f}}, \  \wh{(\eta\lhd g)\lhd f}],$$ so by the Matching Lemma it suffices to show $$[\Id_{\id_d}, \ \wh{\alpha_{h,g,f}}  ,\   \wh{\alpha_{i,g,f}} \widehat{\bullet} \wh{(\eta \lhd (g \circ f))} , \  \wh{(\eta\lhd g)\lhd f}]$$ is commutative. The following table proof verifies this commutativity:

\begin{table}[H] \caption{~ }\label{ax11pt1}\begin{center}

\scalebox{.95}{  \begin{tabular}{ r | l || l | l | l | l | }
     \cline{2-6}

       & $s_2(\wh{\eta\lhd g})$     & $\Id_{\id_d}$   &    $\wh{\Id_{h\circ g}}$     &  $\wh{\eta \lhd g}$            &$\wh{\eta \lhd g}$     \\ \cline{2-6}
     
       $\Lambda$   &   & $\Id_{\id_d}$   &  $\wh{\alpha_{h,g,f}}$   & $\wh{\alpha_{i,g,f}} \widehat{\bullet} \wh{(\eta \lhd (g \circ f))} $   &     $\wh{(\eta\lhd g)\lhd f}$     \\   \cline{2-6}  

   (Lem.~\ref{overlinelemma})   &   $\Delta_{\stackrel{\wedge}{-}}(\widetilde{\alpha}_{h,g,f} )$  &  $\wh{\Id_{h\circ g}}$  &  $\wh{\alpha_{h,g,f}}$     &    $\widetilde{\alpha}_{h,g,f}$   &  $\chi(h\circ g, f)$        \\   \cline{2-6}  
     
      $\odot$  &   (Table~\ref{ax11pt2})       & $\wh{\eta \lhd g}$  & $\wh{\alpha_{i,g,f}}\widehat{\bullet} \wh{(\eta \lhd (g \circ f))} $  &  $\widetilde{\alpha}_{h,g,f}$     &  $\chi(i\circ g, f)$\\ \cline{2-6}

    & $\Delta_{\lhd}(\eta\lhd g,f)$   &$\wh{\eta \lhd g}$  &       $\wh{(\eta\lhd g)\lhd f}$    &  $\chi(h\circ g, f)$      &  $\chi(i \circ g, f)$ \\   \cline{2-6} 
    \end{tabular}}\end{center}
    \end{table}

\begin{table}[H] \caption{~ }\label{ax11pt2}\begin{center}

  \scalebox{.95}{ \begin{tabular}{ r | l || l | l | l | l | }
     \cline{2-6}

       &  $s_1(\wh{\eta \lhd g})$ & $\Id_{\id_d}$   &   $\wh{\eta \lhd g}$   &  $\wh{\eta \lhd g}$            & $\wh{\Id_{i\circ g}}$       \\ \cline{2-6}
     
     &  $\Delta_{\widehat{\bullet}}(  \wh{\alpha_{i,g,f}},\wh{\eta \lhd (g \circ f)})$    &  $\Id_{\id_d}$ &  $\wh{\eta \lhd (g \circ f)}$   &  $\wh{\alpha_{i,g,f}} \widehat{\bullet}\wh{(\eta \lhd (g \circ f))} $  &     $ \wh{\alpha_{i,g,f}}$       \\   \cline{2-6}  

 $\odot$    & (Table~\ref{ax11pt3})    & $\wh{\eta \lhd g}$  &   $\wh{\eta \lhd (g \circ f)}$     &    $\widetilde{\alpha}_{h,g,f}$   & $\widetilde{\alpha}_{i,g,f}$             \\   \cline{2-6}  
     
     $\Lambda$   &       & $\wh{\eta \lhd g}$  & $\wh{\alpha_{i,g,f}} \widehat{\bullet} \wh{(\eta \lhd (g \circ f))} $  &  $\widetilde{\alpha}_{h,g,f}$     &  $\chi(i\circ g, f)$ \\ \cline{2-6}

 (\ref{overlinelemma}) &   $\Delta_{\stackrel{\wedge}{-}}(\widetilde{\alpha}_{i,g,f})$    & $\wh{\Id_{i\circ g}}$ &    $\wh{\alpha_{i,g,f}}$&      $\widetilde{\alpha}_{i,g,f}$    &  $\chi(i\circ g, f)$   \\   \cline{2-6} 
    \end{tabular}}\end{center}
    \end{table}

\begin{table}[H] \caption{~ }\label{ax11pt3}
\begin{center}

    \begin{tabular}{ r | l || l | l | l | l | }
     \cline{2-6}

      & $\Delta_{\lhd}(\eta,g)$     & $\wh{\eta}$   &    $\wh{\eta \lhd g}$     &   $\chi(h,g) $           &$\chi(i,g) $        \\ \cline{2-6}
     
      &    $\Delta_{\lhd}(\eta,g\circ f)$      & $\wh{\eta}$ & $\wh{\eta \lhd (g \circ f)}$   & $\chi(h,g\circ f) $  &        $\chi(i,g\circ f) $     \\   \cline{2-6}  

   $\Lambda$    &    & $\wh{\eta \lhd g}$  &   $\wh{\eta \lhd (g \circ f)}$     &    $\widetilde{\alpha}_{h,g,f}$   & $\widetilde{\alpha}_{i,g,f}$         \\   \cline{2-6}  
     
    & $\Delta_{\widetilde{\alpha}}(h,g,f)$ &  $\chi(h,g) $ & $\chi(h,g\circ f)$ &   $\widetilde{\alpha}_{h,g,f}$      &  $\chi(g,f)$   \\ \cline{2-6}

    &  $\Delta_{\widetilde{\alpha}}(i,g,f)$   & $\chi(i,g) $  &    $\chi(i,g\circ f) $   & $\widetilde{\alpha}_{i,g,f}$      &   $\chi(g,f)$    \\   \cline{2-6} 
    \end{tabular}\end{center}
    \end{table}

\end{proof}
\subsection{Compatibility of the unitors and the pseudo-identity and associator}

\begin{proposition}[Axiom 12] For any object $a$ in $\Bic(X)$, $\lambda_{\id_a}=\rho_{\id_a}.$\end{proposition}
\begin{proof} $$\lambda_{\id_a}=\underline{\widehat{\Id_{\id_a}}}=\underline{s_1s_0(a)}=\underline{s_0s_0(a)}=\underline{\Id_{\id_a}}=\rho_{\id_a}.$$ Note we have applied  Lemma~\ref{lemma3} for the second equality.
\end{proof}
\begin{proposition}[Axiom 13]Let $a\stackrel{f}{\ra}b\stackrel{g}{\ra}c$ in $\Bic(X)$ then $$\alpha_{g,f,\id_a}\bullet(g\rhd 
\rho_f)=\rho_{g\circ f}.$$
 \end{proposition}
 \begin{proof} By definition $\rho_{g\circ f}=\underline{\Id_{g\circ f}}.$ By Lemma~\ref{bulletylemma} we have $$\underline{\widetilde{\alpha}_{g,f,\id_a}\bullet(g\rhd 
\rho_f)}=\alpha_{g,f,\id_a}\bullet(g\rhd \rho_f)$$ so it is enough for us to show $$\widetilde{\alpha}_{g,f,\id_a}\bullet(g\rhd \rho_f)=\Id_{g\circ f}$$
we have $$d\left(  \Delta_{\bullet}(\widetilde{\alpha}_{g,f,\id_a} ,g\rhd \rho_f ) \right)=[\widetilde{\alpha}_{g,f,\id_a},\ \widetilde{\alpha}_{g,f,\id_a}\bullet(g\rhd \rho_f),\ g\rhd \rho_f ,\ \Id_{\id_a} ] $$ so by the Matching Lemma it is enough to show $$[\widetilde{\alpha}_{g,f,\id_a},\ \Id_{g\circ f},\ g\rhd \rho_f ,\ \Id_{\id_a}]$$ is commutative. The following horn verifies this commutativity:
\begin{table}[H] \begin{center}

    \begin{tabular}{ r | l || l | l | l | l | }
     \cline{2-6}

    & $\Delta_{\widetilde{\alpha}}(g,f,\id_a)$      & $\chi(g,f)$   &    $\chi(g,f\circ\id_a)$     &   $\widetilde{\alpha}_{g,f,\id_a} $           &$\chi(f,\id_a) $         \\ \cline{2-6}
     
      &    $s_0 (\chi(g,f))$      & $\chi(g,f)$ &  $\chi(g,f)$    & $\Id_{g\circ f} $  &        $\Id_f $     \\   \cline{2-6}  

    &  $\Delta_{\rhd}(g,\underline{\Id_f})$  &  $\chi(g,f\circ \id_a)$ &    $\chi(g,f)$     &    $g\rhd\underline{\Id_f}$   & $\underline{\Id_f}$           \\   \cline{2-6}  
     
  $\Lambda$  &      &$\widetilde{\alpha}_{g,f,\id_a} $    &  $\Id_{g\circ f} $    &$g\rhd\underline{\Id_f}$       &   $\Id_{\id_a}$\\ \cline{2-6}

      &  $\Delta_{-}(\Id_f)$    & $\chi(f,\id_a) $  &     $\Id_f$   &  $\underline{\Id_f}$       &   $\Id_{\id_a}$ \\   \cline{2-6} 
    \end{tabular}\end{center}
    \end{table}

 \end{proof}

\begin{proposition}[Axiom 14]Let $f\stackrel{\eta}{\Rightarrow} g: a \ra b$ and $h\stackrel{\theta}{\Rightarrow}i : b \ra c$ then $$\rho_g \lhd f  =\alpha_{g,\id_b,f}\bullet (g\rhd \lambda_f).$$
\end{proposition}
\begin{proof}By Lemmas~\ref{bulletylemma} and \ref{trianglelemma},  \begin{align*}\alpha_{g,\id_b,f}\bullet (g \rhd \lambda_f)  &=\underline{\widetilde{\alpha}_{g,\id_b,f}\bullet (g \rhd \lambda_f)}\\\rho_g \lhd f. &= \underline{\rho_g \tlhd f.}.\end{align*} So it suffices to show $$\widetilde{\alpha}_{g,\id_b,f}\bullet (g \rhd \lambda_f)=\rho_g \tlhd f.$$ We have
$$d\left(\Delta_{\bullet}(\widetilde{\alpha}_{g,\id_b,f},g \rhd \lambda_f)\right)=[\widetilde{\alpha}_{g,\id_b,f},\ \widetilde{\alpha}_{g,\id_b,f} \bullet (g \rhd \lambda_f), \ g \rhd \lambda_f , \ \Id_f],$$ and by the Matching Lemma it is enough to show $$[\widetilde{\alpha}_{g,\id_b,f},\ \rho_g\tlhd f, \  g \rhd \lambda_f , \ \Id_f]$$ is commutative. The following table proof verifies this commutativity:

\begin{table}[H] \caption{~ }\label{ax12pt1}\begin{center}

    \begin{tabular}{ r | l || l | l | l | l | }
     \cline{2-6}

        & $\Delta_{\widetilde{\alpha}}(g,\id_b,f)$  & $\chi(g,\id_b)$   &    $\chi(g,\id_b\circ f)$     &   $\widetilde{\alpha}_{g,\id_b,f}$           &$\chi(\id_b,f) $         \\ \cline{2-6}
     
 $\odot$     &  (Table~\ref{ax12pt2})  &  $\chi(g,\id_b)$ & $\chi(g,f)$   &$\rho_g\tlhd f$  &        $\wh{\Id_f} $      \\   \cline{2-6}  

       & $\Delta_{\rhd}(g,\lambda_f)$     &  $\chi(g,\id_b\circ f)$ &  $\chi(g,f)$    &    $g \rhd \lambda_f$   &$\lambda_f=\underline{(\wh{\Id_f})}$          \\   \cline{2-6}  
     
    $\Lambda$ &    & $\widetilde{\alpha}_{g,\id_b,f}$   &  $\rho_g\tlhd f$   &   $g\rhd \lambda_f$     & $\Id_f$   \\ \cline{2-6}

   &  $\Delta_{-}(\wh{\Id_f})$      & $\chi(\id_b,f) $    & $\wh{\Id_f} $     &  $\underline{(\wh{\Id_f})}$ & $\Id_f$   \\   \cline{2-6} 
    \end{tabular}\end{center}
    \end{table}

\begin{table}[H] \caption{~ }\label{ax12pt2}\begin{center}

    \begin{tabular}{ r | l || l | l | l | l | }
     \cline{2-6}

         & $\Delta_{-}(\Id_g)$   & $\chi(g,\id_b)$   &    $\Id_g$     &   $\rho_g=\underline{\Id_g}$           &$\Id_{\id_b} $       \\ \cline{2-6}
     
    $\Lambda$  &     &  $\chi(g,\id_b)$ & $\chi(g,f)$   &$\rho_g\tlhd f$  &        $\wh{\Id_f} $       \\   \cline{2-6}  

        & $s_1(\chi(g,f) )$    & $\Id_g$       &  $\chi(g,f)$          &    $\chi(g,f)$   &   $\wh{\Id_f}$           \\   \cline{2-6}  
      
    & $\Delta_{\tlhd}(\rho_g,f)$     & $\rho_g$   &  $\rho_g\tlhd f$   &   $\chi(g,f)$     &    $\wh{\Id_f}$  \\ \cline{2-6}

      & $s_2(\wh{\Id_f})$   & $\Id_{\id_b} $    & $\wh{\Id_f} $     &  $\wh{\Id_f}$  & $\wh{\Id_f}$   \\   \cline{2-6} 
    \end{tabular}\end{center}
    \end{table}

\end{proof}
\begin{proposition}[Axiom 15]Let $a\stackrel{f}{\ra}b\stackrel{g}{\ra}c$ in $\Bic(X)$. Then $\lambda_g \lhd f  =\alpha_{\id_c,g,f}\bullet\lambda_{g\circ f}.$
\end{proposition}
By Lemma~\ref{trianglelemma} we have $\underline{\lambda_g \tlhd f   } = \lambda_g \lhd f.$ Also, by Lemma~\ref{bulletylemma}, $$\underline{\widetilde{\alpha}_{\id_c,g,f}\bullet\lambda_{g\circ f}}=\alpha_{\id_c,g,f}\bullet\lambda_{g\circ f}$$ so it suffices to show $$\lambda_g \tlhd f=\widetilde{\alpha}_{\id_c,g,f}\bullet\lambda_{g\circ f}.$$ We have $$d\left( \Delta_{\bullet}(\widetilde{\alpha}_{\id_c,g,f} , \lambda_{g\circ f} ) \right)= [  \widetilde{\alpha}_{\id_c,g,f} ,\  \widetilde{\alpha}_{\id_c,g,f}\bullet \lambda_{g\circ f} ,\   \lambda_{g\circ f} ,\ \Id_f ]$$ so by the Matching Lemma we need only show $$[  \widetilde{\alpha}_{\id_c,g,f} ,\ \lambda_g \tlhd f ,\   \lambda_{g\circ f} ,\ \Id_f ]$$ is commutative. The following table proof verifies this commutativity:

\begin{table}[H] \caption{~ }\label{ax15pt1}\begin{center}

    \begin{tabular}{ r | l || l | l | l | l | }
     \cline{2-6}

          & $\Delta_{\widetilde{\alpha}}(\id_c,g,f)$   & $\chi(\id_c,g)$   &    $\chi(\id_c,g\circ f)$     &   $\widetilde{\alpha}_{\id_c,g,f}$           &$\chi(g,f) $      \\ \cline{2-6}
     
  $\odot$   & (Table \ref{ax15pt2})    &  $\chi(\id_c,g)$ & $\widehat{\Id_{g\circ f}}$   & $\lambda_g\tlhd f$  &        $\chi(g,f) $   \\   \cline{2-6}  

        & $\Delta_{-}(\widehat{\Id_{g\circ f}})$     &$\chi(\id_c,g\circ f)$        &  $\widehat{\Id_{g\circ f}}$          &    $\underline{\widehat{\Id_{g\circ f}}}$  &   $\Id_{g\circ f}$          \\   \cline{2-6}  
      
     $\Lambda$    &    & $\widetilde{\alpha}_{\id_c,g,f}$   &  $\lambda_g\tlhd f$   &  $\lambda_{g\circ f}=\underline{\widehat{\Id_{g\circ f}}}$     &    $\Id_f $    \\ \cline{2-6}

      & $s_0(\chi(g,f))$ &$\chi(g,f) $   & $\chi(g,f) $     &    $\Id_{g\circ f}$  & $\Id_f$     \\   \cline{2-6} 
    \end{tabular}\end{center}
    \end{table}

\begin{table}[H] \caption{~ }\label{ax15pt2}\begin{center}

    \begin{tabular}{ r | l || l | l | l | l | }
     \cline{2-6}

       & $\Delta_{-}(\widehat{\Id_g})$   & $\chi(\id_c,g)$   &    $\widehat{\Id_g}$     &   $\underline{(\widehat{\Id_g})}$           &$\Id_g $         \\ \cline{2-6}
     
  $\Lambda$      &  &  $\chi(\id_c,g)$ & $\widehat{\Id_{g\circ f}}$   & $\lambda_g\tlhd f$  &        $\chi(g,f) $    \\   \cline{2-6}  

          & $s_2(\chi(g,f) )$     & $\widehat{\Id_g}$       &  $\widehat{\Id_{g\circ f}}$          &    $\chi(g,f)$  &  $ \chi(g,f)$         \\   \cline{2-6}  
      
      & $\Delta_{\tlhd}(\lambda_g,f)$    &  $\lambda_g=\underline{(\widehat{\Id_g})}$    &  $\lambda_g\tlhd f$    &  $\chi(g,f)$     &    $\widehat{\Id_f}$     \\ \cline{2-6}

       & $s_1(\chi(g,f))$    &$\Id_g $  & $\chi(g,f) $  &  $\chi(g,f) $& $\widehat{\Id_f}$ \\   \cline{2-6} 
    \end{tabular}\end{center}
    \end{table}

\subsection{Full interchange}
\begin{proposition}[Axiom 16] Let $f\stackrel{\eta}{\Rightarrow} g:a\ra b$ and $h\stackrel{\theta}{\Rightarrow} i:c\ra d$ in $\Bic(X).$ Then $$(i\rhd \eta) \bullet (\theta \lhd f)=(\theta \lhd g)\bullet (h \rhd \eta).$$
\end{proposition}
\begin{proof}$$d\left(\Delta_{\bullet}(i\rhd \eta,\theta \lhd f)\right)=[i\rhd \eta, (i\rhd \eta) \bullet (\theta \lhd f), \theta \lhd f,\Id_{\id_a}]$$ so by the Matching Lemma it suffices to show $$[i\rhd \eta, (\theta \lhd g)\bullet (h \rhd \eta), \theta \lhd f,\Id_{\id_a}]$$ is commutative. The following table proof verifies this commutativity:

\begin{table}[H] \caption{~ }\label{ax13pt1}\begin{center}

    \begin{tabular}{ r | l || l | l | l | l | }
     \cline{2-6}

      & $\Delta_{\rhd}(i,\eta)$  & $\chi(i,g)$   &    $\chi(i,f)$     &   $i\rhd \eta$           &$\eta$           \\ \cline{2-6}
     
    $\odot$  & (Table~\ref{ax13pt2}) &$\chi(i,g)$  & $\theta \tlhd f$  & $(\theta \lhd g)\bullet (h \rhd \eta)$  &        $\eta$     \\   \cline{2-6}  

  & $\Delta_{-}(\theta \tlhd f)$     &$\chi(i,f)$     &$\theta \tlhd f$      &   $\theta \lhd f$    &$\Id_f$          \\   \cline{2-6}  
     
    $\Lambda$  &    & $i\rhd \eta$   &  $(\theta \lhd g)\bullet (h \rhd \eta)$   &   $\theta\lhd f$     &    $\Id_{\id_a}$ \\ \cline{2-6}

       &  $s_0\eta$  & $\eta$  & $\eta $     &  $\Id_f$& $\Id_{\id_a}$  \\   \cline{2-6} 
    \end{tabular}\end{center}
    \end{table}

\begin{table}[H] \caption{~ }\label{ax13pt2}\begin{center}

    \begin{tabular}{ r | l  || l | l | l| l | }
     \cline{2-6}

                  & $\Delta_{-}(\theta\tlhd g)$     & $\chi(i,g)$   &   $\theta\tlhd g$      &   $\theta\lhd g$           &$\Id_g$          \\ \cline{2-6}
     
  $\Lambda$         &         &$\chi(i,g)$  & $\theta \tlhd f$  & $(\theta \lhd g)\bullet (h \rhd \eta)$  &        $\eta$   \\   \cline{2-6}  

      $\odot$     &  (Table~\ref{ax13pt3})        &$\theta\tlhd g$    &$\theta \tlhd f$      &   $h\rhd\eta$    &$\eta$           \\   \cline{2-6}  
     
    & $\Delta_{\bullet}(\theta\lhd g,h\rhd\eta)$    & $\theta\lhd g$   &  $(\theta \lhd g)\bullet (h \rhd \eta)$   &   $h\rhd \eta$     &    $\Id_{\id_a}$\\ \cline{2-6}

      &  $s_1\eta$  & $\Id_g$  & $\eta $     &  $\eta$& $\Id_{\id_a}$   \\   \cline{2-6} 
    \end{tabular}\end{center}
    \end{table}

\begin{table}[H] \caption{~ }\label{ax13pt3}\begin{center}

    \begin{tabular}{ r | l || l | l | l | l | }
     \cline{2-6}

                   & $\Delta_{\tlhd}(\theta, g)$      & $\theta$   &   $\theta\tlhd g$      &   $\chi(h,g)$           &$\wh{\Id_g}$        \\ \cline{2-6}
     
                     &  $\Delta_{\tlhd}(\theta, f)$       &$\theta$  & $\theta \tlhd f$  & $\chi(h,f)$     &        $\wh{\Id_f}$    \\   \cline{2-6}  

      $\Lambda$           &       &$\theta\tlhd g$    &$\theta \tlhd f$      &   $h\rhd\eta$    &$\eta$        \\   \cline{2-6}  
     
    & $\Delta_{\rhd}(h,\eta)$    & $\chi(h,g)$   &  $\chi(h,f)$   &  $h\rhd\eta$    &    $\eta$\\ \cline{2-6}

     &  $s_2\eta$   & $\wh{\Id_g}$ &   $\wh{\Id_f}$  &  $\eta$& $\eta$   \\   \cline{2-6} 
    \end{tabular}\end{center}
    \end{table}

\end{proof}
\subsection{Pentagon identity}
\begin{proposition}[Axiom 17] Let $a \stackrel{f}\to b \stackrel{g}\to c \stackrel{h}\to d \stackrel{i}\to e$ in $\Bic(X).$  Then $$(\alpha_{i,h,g} \lhd f) \bullet (\alpha_{i,h \circ g,f} \bullet (i \rhd \alpha_{h,g,f}))=\alpha_{i \circ h,g,f}\bullet \alpha_{i,h,g \circ f}.$$
\end{proposition}
\begin{proof} By Lemmas~\ref{bulletylemma} and \ref{trianglelemma}

\begin{align*}(\alpha_{i,h,g} \lhd f) \bullet (\alpha_{i,h \circ g,f} \bullet (i \rhd \alpha_{h,g,f}))&=\underline{(\alpha_{i,h,g} \tlhd f) \bullet (\alpha_{i,h \circ g,f} \bullet (i \rhd \alpha_{h,g,f}))}\\ \alpha_{i \circ h,g,f}\bullet \alpha_{i,h,g \circ f}&= \underline{\widetilde{\alpha}_{i \circ h,g,f}\bullet \alpha_{i,h,g \circ f}}.\end{align*}
Thus it suffices to show

$$ (\alpha_{i,h,g} \tlhd f) \bullet (\alpha_{i,h \circ g,f} \bullet (i \rhd \alpha_{h,g,f}))=\widetilde{\alpha}_{i \circ h,g,f}\bullet \alpha_{i,h,g \circ f}.$$
We have \begin{multline*} d \left(\Delta_{\bullet} ( \alpha_{i,h,g} \tlhd f , \alpha_{i,h \circ g,f} \bullet (i \rhd \alpha_{h,g,f}))\right)= \\ [\alpha_{i,h,g} \tlhd f,\  (\alpha_{i,h,g} \tlhd f) \bullet (\alpha_{i,h \circ g,f} \bullet (i \rhd \alpha_{h,g,f})),\alpha_{i,h \circ g,f} \bullet (i \rhd \alpha_{h,g,f}), \ \Id_f] \end{multline*} so by the Matching Lemma it is enough to show  $$ [\alpha_{i,h,g} \tlhd f,\ \widetilde{\alpha}_{i \circ h,g,f}\bullet \alpha_{i,h,g \circ f},\ \alpha_{i,h \circ g,f} \bullet (i \rhd \alpha_{h,g,f}), \ \Id_f]$$
is commutative. The following table proof verifies this commutativity:

\begin{table}[H] \caption{~ }\label{ax14pt1}\begin{center}

 \scalebox{.93}{   \begin{tabular}{ r | l ||  l | l | l |l | }
     \cline{2-6}

                    & $\Delta_{\tlhd}(\alpha_{i,h,g}, f)$        & $\alpha_{i,h,g}$   &   $\alpha_{i,h,g}\tlhd f$      &   $\chi(i\circ(h\circ g),f)$           &$\wh{\Id_f}$     \\ \cline{2-6}
     
            $\odot$  &  (Table \ref{ax14pt2}) & $\alpha_{i,h,g}$  & $\widetilde{\alpha}_{i\circ h,g,f} \bullet \alpha_{i,h,g\circ f}$  &    $\widetilde{\alpha}_{i,h \circ g,f} \bullet (i \rhd \alpha_{h,g,f})$       &        $\wh{\Id_f}$   \\   \cline{2-6}  

      $\Lambda$        &        &$\alpha_{i,h,g} \tlhd f$    &$\widetilde{\alpha}_{i \circ h,g,f}\bullet \alpha_{i,h,g \circ f}$      &   $\alpha_{i,h \circ g,f} \bullet (i \rhd \alpha_{h,g,f})$    &$ \Id_f$          \\   \cline{2-6}  
     
    & $\Delta_{-}(\widetilde{\alpha}_{i,h \circ g,f} \bullet (i \rhd \alpha_{h,g,f}))$    & $\chi(i\circ (h\circ g),f)$   &   $\widetilde{\alpha}_{i,h \circ g,f} \bullet (i \rhd \alpha_{h,g,f})$   &    $\alpha_{i,h \circ g,f} \bullet (i \rhd \alpha_{h,g,f})$   &    $ \Id_f$\\ \cline{2-6}

     &  $s_2(\Id_f)$   & $\wh{\Id_f}$ &   $\wh{\Id_f}$  &  $\Id_f$& $\Id_f$   \\   \cline{2-6} 
    \end{tabular}}\end{center}
    \end{table}

\begin{table}[H] \caption{~ }\label{ax14pt2}
\begin{center}

    \begin{tabular}{ r | l || l | l | l | l | }
     \cline{2-6}

               & $\Delta_{-}(\widetilde{\alpha}_{i,h,g})$             & $\chi(i\circ h,g)$   &   $\widetilde{\alpha}_{i,h,g}$      &  $\alpha_{i,h,g}$          &$\Id_g$     \\ \cline{2-6}

      $\odot_1$       &  (Table  \ref{ax14pt3})      &$\chi(i\circ h,g)$    &$\widetilde{\alpha}_{i,h,g\circ f}$      &   $\widetilde{\alpha}_{i\circ h,g,f} \bullet \alpha_{i,h,g\circ f}$    & $\chi(g,f)$         \\   \cline{2-6}  
     
    $\odot_2$  &  (Table \ref{ax14pt4})& $\widetilde{\alpha}_{i,h,g}$   &   $\widetilde{\alpha}_{i,h,g\circ f}$   &    $\widetilde{\alpha}_{i,h \circ g,f} \bullet (i \rhd \alpha_{h,g,f})$   &     $\chi(g,f)$ \\ \cline{2-6}

     $\Lambda$   &  & $\alpha_{i,h,g}$  & $\widetilde{\alpha}_{i\circ h,g,f} \bullet \alpha_{i,h,g\circ f}$  &    $\widetilde{\alpha}_{i,h \circ g,f} \bullet (i \rhd \alpha_{h,g,f})$       &        $\wh{\Id_f}$     \\   \cline{2-6}    
  
      &  $s_1(\chi(g,f))$   & $\Id_g$ &  $\chi(g,f)$  &  $\chi(g,f)$& $\wh{\Id_f}$ \\   \cline{2-6} 
    \end{tabular}
\end{center}
    \end{table}

\begin{table}[H] \caption{Commutativity of  $\odot_1$  in Table~\ref{ax14pt2}}\label{ax14pt3}\begin{center}

    \begin{tabular}{ r | l || l | l | l | l | }
     \cline{2-6}

          & $\Delta_{\widetilde{\alpha}}(i\circ h,g,f)$              & $\chi(i\circ h,g)$   &   $\chi(i\circ h,g\circ f)$      &  $\alpha_{i\circ h,g,f}$          &$\chi(g,f)$         \\ \cline{2-6}

      $\Lambda$   &     &$\chi(i\circ h,g)$    &$\widetilde{\alpha}_{i,h,g\circ f}$      &   $\widetilde{\alpha}_{i\circ h,g,f} \bullet \alpha_{i,h,g\circ f}$    & $\chi(g,f)$          \\   \cline{2-6}  
     
  & $\Delta_{-}(\widetilde{\alpha}_{i, h, g\circ f})$ &  $\chi(i\circ h,g\circ f)$   &   $\widetilde{\alpha}_{i,h,g\circ f}$   &    $\alpha_{i,h,g\circ f}$   &     $\Id_{g\circ f}$\\ \cline{2-6}

   &  $\Delta_{\bullet}( \widetilde{\alpha}_{i\circ h,g,f},\alpha_{i, h,g\circ f} )  $    & $\widetilde{\alpha}_{i\circ h,g,f}$  & $\widetilde{\alpha}_{i\circ h,g,f} \bullet \alpha_{i,h,g\circ f}$  &    $\alpha_{i, h,g\circ f}$   &     $\Id_f$     \\   \cline{2-6}    
  
    &  $s_1(\chi(g,f))$   & $\Id_g$ &  $\chi(g,f)$  &  $\chi(g,f)$& $\wh{\Id_f}$   \\   \cline{2-6}
    \end{tabular} \end{center}
    \end{table}

\begin{table}[H] \caption{Commutativity of  $\odot_2$ in Table~\ref{ax14pt2}, part 1}\label{ax14pt4}\begin{center}

    \begin{tabular}{ r | l || l | l | l | l | }
     \cline{2-6}

        & $\Delta_{\widetilde{\alpha}}(i,h,g)$    &  $\chi(i,h)$  &   $\chi(i,h\circ g)$      & $\widetilde{\alpha}_{i,h,g}$           &$\chi(h,g)$        \\ \cline{2-6}

          &   $\Delta_{\widetilde{\alpha}}(i,h,g\circ f)$      &  $\chi(i,h)$   &$\chi(i,h\circ (g\circ f))$     &  $\widetilde{\alpha}_{i,h,g\circ f}$    & $\chi(h,g\circ f)$     \\   \cline{2-6}  
 
     $\odot$   &  (Table  \ref{ax14pt5})   & $\chi(i,h\circ g)$ & $\chi(i,h\circ (g\circ f))$ &    $\widetilde{\alpha}_{i,h \circ g,f} \bullet (i \rhd \alpha_{h,g,f})$       &      $\widetilde{\alpha}_{h,g,f}$     \\   \cline{2-6}     
        
    $\Lambda$&    & $\widetilde{\alpha}_{i,h,g}$   &   $\widetilde{\alpha}_{i,h,g\circ f}$   &    $\widetilde{\alpha}_{i,h \circ g,f} \bullet (i \rhd \alpha_{h,g,f})$   &     $\chi(g,f)$\\ \cline{2-6}

     &  $\Delta_{\widetilde{\alpha}}( h,g,f)$  & $\chi(h,g)$ &  $\chi(h,g\circ f)$  &  $\widetilde{\alpha}_{h,g,f}$   & $\chi(g,f)$  \\   \cline{2-6} 
    \end{tabular}\end{center}
    \end{table}

\begin{table}[H] \caption{Commutativity of $\odot_2$ in Table~\ref{ax14pt2}, part 2}\label{ax14pt5}\begin{center}

 \scalebox{.90}{   \begin{tabular}{ r | l || l | l | l | l | }
     \cline{2-6}

      & $\Delta_{\widetilde{\alpha}}(i,h\circ g,f)$      &  $\chi(i,h\circ g)$  &   $\chi(i,(h\circ g)\circ f)$      & $\widetilde{\alpha}_{i,h\circ g,f}$           &$\chi(h\circ g,f)$        \\ \cline{2-6}

   $\Lambda$      &     & $\chi(i,h\circ g)$ & $\chi(i,h\circ (g\circ f))$ &    $\widetilde{\alpha}_{i,h \circ g,f} \bullet (i \rhd \alpha_{h,g,f})$       &      $\widetilde{\alpha}_{h,g,f}$   \\   \cline{2-6}  

          &   $\Delta_{\rhd}(i,\alpha_{h,g,f})$  &   $\chi(i,(h\circ g)\circ f)$  & $\chi(i,h\circ (g\circ f))$     &   $i \rhd \alpha_{h,g,f}$     & $ \alpha_{h,g,f}$         \\   \cline{2-6}

          &  $\Delta_{\bullet}(\widetilde{\alpha}_{i,h\circ g,f},i \rhd \alpha_{h,g,f} )$    & $\widetilde{\alpha}_{i,h\circ g,f}$   &  $\widetilde{\alpha}_{i,h \circ g,f} \bullet (i \rhd \alpha_{h,g,f})$     &    $i \rhd \alpha_{h,g,f}$   & $\Id_f$   \\ \cline{2-6}

     &  $\Delta_{-}(\widetilde{\alpha}_{h,g,f})$   & $\chi(h\circ g,f)$ &   $\widetilde{\alpha}_{h,g,f}$ &  $\alpha_{h,g,f}$   &$\Id_f$ \\   \cline{2-6} 
    \end{tabular}}\end{center}
    \end{table}

\end{proof}


\section{Coskeleta of simplicial sets}


 Having finished the proof that $\Bic(X)$ is a bicategory, we now work toward defining the Duskin nerve of a bicategory, again following \cite{Dus02}. This construction will be a crucial part of each of the nerve constructions given in subsequent chapters. We first make some remarks about the coskeleton functor for simpicial sets. For the omitted proofs, refer to Subsection~\ref{bicoskeletasec} where the proofs are given in a more general setting, or refer to \cite{May67} or \cite{GJ99}.

\begin{definition}Let $\Delta |^k_0$ denote the full subcategory of $\Delta$ on the objects $\{[0],[1],\ldots,[k]\}.$ Then let $\cat{Set}_{\Delta |^k_0}$ denote the category of presheaves of sets on $\Delta |^k_0$. Recalling that a simplicial set is a presheaf of sets on $\Delta,$ we get the \emph{truncation} functors $$\mathbf{tr}^k: \cat{Set}_\Delta \ra \cat{Set}_{\Delta |^k_0}$$ by restriction.
\end{definition}
\begin{definition}
The functor $\mathbf{tr}^k: \cat{Set}_\Delta  \ra \cat{Set}_{\Delta |^k_0}$ has adjoints on both sides. The left adjoint is called the \emph{$k$-skeleton}, $\mathbf{sk}^k$ and the right adjoint the \emph{$k$-coskeleton}, \emph{$\mathbf{cosk}^k$}. Composing with $\mathbf{tr}^k$ we get endofunctors of $\cat{Set}_\Delta$, $\mathbf{Sk}^k=\mathbf{sk}^k\mathbf{tr}^k$ and $\mathbf{Cosk}^k=\mathbf{cosk}^k\mathbf{tr}^k,$ which we also refer to as the $k$-skeleton and $k$-coskeleton functors. We will focus our attention on the coskeleton, which can be computed by the formula: $$(\mathbf{cosk}^k(X))_n= \Hom(\mathbf{tr}^k(\Delta[n]), X).$$
\end{definition}
\begin{proposition} Let $Y\in \cat{Set}_\Delta$. Then the canonical map $Y \ra \mathbf{Cosk}^k(Y)$ is an isomorphism if and only if every $d\Delta[n]$-sphere has a unique filler for $n>k$, that is, if every map $d\Delta[n]\ra Y$ extends uniquely to a map $\Delta[n]\ra Y$ for $n>k$. If this is the case, we say $Y$ is \emph{$k$-coskeletal}.
\end{proposition}
\begin{remark}
If $Y$ is Kan, then $Y\ra \Cosk^k(Y)$ is a model for $k$th Postnikov section, inducing isomorphisms $\pi_r(Y) \cong \pi_r(\Cosk^k(Y))$ for $r<k$ and with $\pi_r(\Cosk^k(Y))$ being trivial for $r\geq k.$ There is a diagram of natural maps:
\begin{equation*}\label{coskdiagram1} \cdots \ra  \Cosk^{n+2}Y \ra  \Cosk^{n+1} Y\ra \Cosk^n Y \end{equation*} with indirect limit $Y$, giving a model for the Postnikov tower of $Y$. See \cite{DK84} for details.
\end{remark}
\begin{proposition}\label{spherefiller} Suppose $Y\in  \cat{Set}_\Delta$ has the property that every inner $n$-horn has a unique filler for $n \geq k,$ with $k \geq 2$. Then every $d\Delta[k]$-sphere in $Y$ has at most one filler, and for $m>k$ every $m$-sphere has a unique filler, thus $Y$ is $k+1$-coskeletal.
\end{proposition}
\begin{proof}
First, let $n\geq k$, consider an inner horn inside a sphere $\Lambda^n_p \ra d\Delta[n]\ra X.$ Then any filler of the sphere fills this inner horn, so the filler must be unique if it exists. This shows that any $d\Delta[n]$-sphere has at most one filler for $k\geq n$

Now let $n>k$ and consider an inner horn $H$ inside a sphere $S$, $$\Lambda^n_i \ra d\Delta[n]\ra X.$$ Consider the unique filler $F$ of the horn $H$. Clearly every face of $F$, other than perhaps the $i$th face, agrees with the corresponding face of $S.$ From this, the boundary of the $i$th face of $F$ matches the boundary of the $i$th face of $S$, so they are both fillers of the same $d\Delta[n-1]$-sphere, and are identical by what we proved above. So $F$ matches our sphere at every face, $dF=S$ , thus $F$ fills $S.$ 
\end{proof} 
\begin{proposition} \label{spherefiller1}
\label{coskeletaliskan0} Suppose $Y \in \cat{Set}_{\Delta |^k_0}$ with $k \geq 0$. Then:
\begin{itemize}
 \item $\cosk^k Y$ meets the inner Kan condition (uniquely) for $n$ if and only if every inner horn $\tr^k(\Lambda^n_i ) \ra Y$ has a (unique) extension to a map $\tr^k(\Delta[n])\ra Y$ along the map $\tr^k(\Lambda^n_i)\ra \tr^k(\Delta[n]).$ 
 \item $\cosk^k Y$ meets the inner Kan condition uniquely for  $n \geq k+2$.
\end{itemize}
\end{proposition}
\begin{proof} The first statement is immediate from the fact that $\cosk^k$ is right adjoint to $\tr^k$. The second statement can be seen to be a special case of the first since $$\tr^k (\Lambda^n_i)\ra \tr^k(\Delta[n])$$ is an isomorphism if $n\geq k+2$.
\end{proof}


\section{The Duskin nerve of a small bicategory \label{duskinnerve}}


Let $\B$ be a small bicategory. We first define the truncated Duskin nerve $N(\B)|^3_0\in \cat{Set}_{\Delta |^k_0}.$ The 0-cells of $N(\B)|^3_0$ are objects of $\B$, and the $1$-cells are the 1-morphisms. The face maps $d_1,d_0:N(\B)|^3_0(1)\ra N(\B)|^3_0(0)$ are defined by $d_0(f)=\mbox{target}(f)$ and $d_1(f)=\mbox{source}(f).$ The degeneracy map $s_0: N(\B)|^3_0(0)\ra N(\B)|^3_0(1)$ is given by $s_0(a)=\id_a.$

We define a 2-cell in $N(\B)|^3_0$  to be a triple $(h,g,f)$ together with an \emph{interior} 2-morphism $\eta:g \ra h \circ f,$ which we use the notation $(h,g,f   \vbar \eta)$ to specify. Such a 2-cell has boundary $d(h,g,f \vbar \eta)=[h,g,f].$ When no confusion is possible, we will use the same name for a $2$-cell and its interior $2$-morphism.  The two degeneracy maps $s_0,s_1:N(\B)|^3_0(1)\ra N(\B)|^3_0(2)$ are given for $f:a\ra b$ by letting $s_0(f)$ be the triple $(f,f,\id_a)$  with the 2-morphism  $\rho_f:f \Rightarrow f\circ \id_a$ and $s_1(f)$ be the triple $(\id_b,f,f)$ with the 2-morphism  $\lambda_f:f \Rightarrow \id_b\circ f.$
To define a 3-cell $x_{0123}$ in $N(\B)|^3_0$ is a quadruple of 2-cells $(x_{123},x_{023},x_{013}, x_{012})$ satisfying appropriate face relations so that $[x_{123},x_{023},x_{013},x_{012}]$ forms a sphere, and also satisfying a condition saying the appropriate respective compositions of the even and odd faces agree up to the associator. That is:
 \begin{equation}\label{3cellcondition}\alpha_{x_{23},x_{12},x_{01}}\bullet(x_{23}\rhd x_{012})\bullet x_{023}=(x_{123}\lhd x_{01})\bullet x_{013}\end{equation}
where $x_{23}=d_0(x_{123})=d_0(x_{023})$ and $x_{01}=d_2(x_{012})=d_2(x_{013}).$  The face maps for 3-cells are the obvious maps built into our definition.

Our definitions of the  three degeneracy maps $s_0,s_1,s_2:N(\B)|^3_0(2)\ra N(\B)|^3_0(3)$ are determined by the simplicial identities for $d_is_j$: \begin{align}  s_0(h,g,f \vbar \eta)&=[\eta,\eta,\rho_g,\rho_f]\\s_1(h,g,f \vbar \eta)&= [\rho_h,\eta,\eta,\lambda_f]  \\ s_2(h,g,f \vbar \eta)&= [\lambda_h,\lambda_g,\eta,\eta]  \end{align}
In order to show these definitions give 3-cells in $N(\B)|^3_0,$ we must check that Equation~\ref{3cellcondition} holds for these spheres. First for $s_0$, we must show
\begin{equation}\label{s0condition}\alpha_{h,f,\id_a}\bullet (h\rhd \rho_f) \bullet \eta=(\eta \lhd \id_a)\bullet\rho_g.\end{equation} 
By Axiom 7, i.e. the naturality of $\rho$, we have $$(\eta \lhd \id_a)\bullet\rho_g=\rho_{h\circ f} \bullet \eta.$$ By Axiom 13, the first compatibility axiom of the associator and the unitors, $$\alpha_{h,f,\id_a}\bullet (h\rhd \rho_f)=\rho_{h\circ f},$$ and combining these two we get Equation~\ref{s0condition}. Similarly, for $s_1$ the commutativity condition
$$\alpha_{h,\id_b,f}\bullet(h \rhd \lambda_f)\bullet \eta=(\rho_h\lhd f)\bullet \eta$$ follows directly from Axiom 14, the second  compatibility axiom of the associator and the unitors. Finally for $s_2$ the commutativity condition is $$\alpha_{\id_c,h,f}\bullet(\id_c \rhd \eta)\bullet \lambda_g=(\lambda_h\lhd f)\bullet \eta.$$ This follows, similarly to the $s_0$ case, from Axiom 8 and Axiom 15, i.e. the naturality of $\lambda$ and the third compatibility axiom for the unitors and the associator.
 
Our definition of $N(\B)|^3_0$ makes the simplicial identities for $d_id_j$ and $d_is_j$ hold ipso facto. However, the identities for $s_is_j$ need to be checked. First for a 0-cell, an object $a\in \Bic(X)$, we need to check $s_1s_0a=s_0s_0a,$  that is, $\lambda_{\id_a}=\rho_{\id_a}$. This is Axiom 12, the compatability of the unitors and the pseudo-identity. Now let $f:a\ra b$ be a $1$-cell, a morphism in $\Bic(X).$ We have \begin{align*}s_1s_0 f&= [\rho_f,\ \rho_f,\ \rho_f ,\ \lambda_{\id_a}]\\s_0s_0f&= [\rho_f,\ \rho_f,\ \rho_f ,\ \rho_{\id_a}]. \end{align*} So this again hold by Axiom 12. The case $s_2s_1f=s_1s_1f$ also follows by Axiom 12. The case $s_2s_0f=s_0s_1f$ can be seen to hold immediately.

This completes the definition of $N(\B)|^3_0.$
\begin{definition} \label{nervedef} Let $\B$ be a small bicategory. Then the \emph{Duskin nerve}, or simply the \emph{nerve} of $\B$, is defined by $N(\B)=\mathbf{cosk}^3(N(\B)|^3_0).$    
\end{definition}


\subsection{The inner-Kan conditions for $N(\B).$}


Let $B$ be a small $(2,1)$-category. 
\begin{theorem}[Duskin] \label{daggertheorem} $N(\B)$ is $2$-reduced inner-Kan, that is, $N(\B)$ satisfies the inner-Kan condition, and satisfies it uniquely for $n > 2.$
\end{theorem}\begin{proof}
By Proposition~\ref{spherefiller1} it is enough to check the filler conditions for truncated horns of dimension $2$, $3$, and $4$ in $N(\B)|^3_0$

A horn of type $\Lambda^2_1 \ra N(\B)|^3_0$ is equivalent to a pair of composable morphisms in $\B,$ $a\stackrel{f}{\ra} b\stackrel{g}{\ra}c.$ We gave in Definition~\ref{nervedef} a preferred filler for these horns, but this filler in general is not unique.

The horns of type $\Lambda^3_1 \ra N(\B)|^3_0$ and $\Lambda^3_2$ do not necessarily have fillers for an arbitrary bicategory.  The statement that a horn of type $\Lambda^3_1 \ra N(\B)|^3_0$ has a unique filler in $N(\B)|^3_0$ is equivalent to the statement, that given $2$-cells $x_{012},x_{123},x_{013},$ there is a unique $2$-morphism $\eta$ such that $$\alpha_{x_{23},x_{12},x_{01}}\bullet(x_{23}\rhd x_{012})\bullet \eta=(x_{123}\lhd x_{01})\bullet x_{013}.$$ If $\B$ is a small $(2,1)$-category, this follows from Axiom 3, the invertibility of 2-morphisms. Note that Axiom~2, associativity for 2-morphisms, is also needed here since it is needed to show the uniqueness of inverses for 2-morphisms. The case of $\Lambda^3_2$-horns is similar.

Finally we consider horns of type $\Lambda^4_1,\Lambda^4_2,$, and $\Lambda^4_3.$ Note that by Proposition~\ref{spherefiller1} we are looking for fillers of the form $\tr^3(\Lambda^4_i) \ra \mathbf{tr}^3 \Delta[4]$, which corresponds to filling in the ``missing face'' of the horn. For instance, consider a horn $$[x_{1234},\ - ,\ x_{0134},\ x_{0124},\ x_{0123}]$$ in $N(B)|^3_0.$ The filling face has faces $[x_{234},x_{034},x_{024},x_{023}].$
This uniquely specifies it as a 3-cell according to our definition, we need only show that these faces satisfy the commutativity condition. The commutativity conditions for the five faces are:
 \begin{align}
 \alpha_{x_{34},x_{23},x_{12}}\bullet(x_{34}\rhd x_{123})\bullet x_{134}&=(x_{234}\lhd x_{12})\bullet x_{124} \\ 
 \alpha_{x_{34},x_{23},x_{02}}\bullet(x_{34}\rhd x_{023})\bullet x_{034}&=(x_{234}\lhd x_{12})\bullet x_{024} \label{cond1} \\
 \alpha_{x_{34},x_{13},x_{01}}\bullet(x_{34}\rhd x_{013})\bullet x_{034}&=(x_{134}\lhd x_{01})\bullet x_{013} \label{cond2}\\
 \alpha_{x_{24},x_{12},x_{01}}\bullet(x_{24}\rhd x_{012})\bullet x_{024}&=(x_{124}\lhd x_{01})\bullet x_{014} \label{cond3} \\
 \alpha_{x_{23},x_{12},x_{01}}\bullet(x_{23}\rhd x_{012})\bullet x_{023}&=(x_{123}\lhd x_{01})\bullet x_{013} 
 \end{align}
For the horns $\Lambda^4_1,\Lambda^4_2,$, and $\Lambda^4_3.$ we need to show that each of the middle three conditions, \ref{cond1}, \ref{cond2}, and \ref{cond3} hold, given in each case the other four conditions. This follows formally from the bicategory axioms, by a proof that uses all cases of interchange, the naturality of the associator in each of its variables, and the pentagon identity. We refer the reader to Section~6.7 of \cite{Dus02} for details.
\end{proof}
 Given a $\Lambda_1^2$-horn given by two composable morphisms $f,g$ in $\B$, we have a preferred filler $(g,g\circ f, f \vbar \Id_{g \circ f})$ in $N(\B)$, and we let these fillers define $N(\B)$ as an algebraic $2$-reduced inner-Kan simplicial set.

\section{The isomorphism $(X,\chi) \cong N(\Bic(X,\chi))$ \label{bicisosec}}
In this section, we give a strict (i.e. algebraic-structure-preserving) isomorphism $$u:(X,\chi) \ra N(\Bic(X,\chi)),$$ and finish the characterization of the simplicial sets which are nerves of bicategories. Because $N(\Bic(X))$ and $X$ are $3$-coskeletal, that is \begin{align*}X&\cong \mathbf{Cosk}^3(X)=\mathbf{cosk}^3\mathbf{tr}^3(X) \\ N(\Bic(X))&\cong \mathbf{Cosk}^3(N(\Bic(X)))=\mathbf{cosk}^3\mathbf{tr}^3(N(\Bic(X)))\end{align*} it suffices to give the isomorphism $u:\mathbf{tr}^3(X) \cong \mathbf{tr}^3(N(\Bic(X)))$. 

For $0$-cells and $1$-cells, by definition $X$ and $N(\Bic(X))$ are identical, so we take $u(x)=x.$ For a $2$ cell $x$ with $dx = [h,g,f]$ we let $u(x)= (h,g,f \vbar \underline{x}).$ We construct the inverse map $u^{-1}$ for a $2$ cell $(h,g,f \vbar \eta)$ in $N(\Bic(X))$ by the following horn: 
\begin{table}[H]\caption*{$\Lambda_{u^{-1}}(h,g,f \vbar \eta)$}\begin{center}
    \begin{tabular}{ r | l || l | l | l | }
     \cline{2-5}
              &  $\chi(h,f)$            &$h$              & $h\circ f$                    & $f$                \\   \cline{2-5} 
          
          $\Lambda$        & $=:u^{-1}(h,g,f \vbar \eta)$      &$h$              & $g$        &  $f$  \\ \cline{2-5}
    
                        & $\eta$       &$h\circ f$      & $g$          & $\id_a$            \\   \cline{2-5}
            
                     & $\Id_f$        &$f$              & $f$                   & $\id_a$                \\   \cline{2-5}
    \cline{2-5}
    \end{tabular}\end{center}
    \end{table}
So we have \begin{align*}
d\left(\Delta_{u^{-1}}(h,g,f \vbar \underline{x})\right)&=[\chi(h,f),\ u^{-1}(h,g,f \vbar \underline{x})=u^{-1}(u(x)),\  \underline{x},\ \Id_f] \\ 
d\left(\Delta_{-}(x)\right)&=[\chi(h,f),\ x,\  \underline{x},\ \Id_f].
\end{align*}
 So by the Matching Lemma $uu^{-1}$ is the identity. By a similar argument, we see $$\underline{u^{-1}(h,g,f \vbar \eta)}=\eta$$ and so $u^{-1}u$ is the identity, showing $u$ is an isomorphism for $2$-cells.
 
 For $3$ cells, the map $u$ is given by $$u(x_{0123})=[ u(x_{123}),u(x_{023}),u(x_{013}),u(x_{012})].$$ This map is injective by Proposition~\ref{spherefiller}, but we must check that it is well-defined and surjective. This amounts to showing the sphere \begin{equation}\label{xsphere}[ x_{123},\ x_{023},\ x_{013},\ x_{012}]\end{equation} is commutative if and only if that the $3$-cell condition holds, that is 
 \begin{equation*}\alpha_{x_{23},x_{12},x_{01}}\bullet(x_{23}\rhd \underline{x_{012}})\bullet \underline{x_{023}}=(\underline{x_{123}}\lhd x_{01})\bullet \underline{x_{013}}.\end{equation*}
 By Lemma~\ref{bulletylemma}, this is equivalent to:
 \begin{equation*}\widetilde{\alpha}_{x_{23},x_{12},x_{01}}\bullet(x_{23}\rhd \underline{x_{012}})\bullet \underline{x_{023}}=(\underline{x_{123}}\tlhd x_{01})\bullet \underline{x_{013}}.\end{equation*}
which is equivalent to the commutativity of the sphere \begin{equation}\label{condsphere}[  \underline{x_{123}}\tlhd x_{01}    ,\    \widetilde{\alpha}_{x_{23},x_{12},x_{01}}\bullet(x_{23}\rhd \underline{x_{012}})\bullet \underline{x_{023}}    ,\     \underline{x_{013}}   ,\     \Id_{\id_{x_0}} ].\end{equation} 
 The following table proof shows the commutativity of \ref{condsphere}, given the commutativity of \ref{xsphere}:
 \begin{table}[H] \caption{~}\label{isopt1a}\begin{center}

    \begin{tabular}{ r | l || l | l | l | l | }
     \cline{2-6}

      & $\Delta_{\tlhd}(\underline{x_{123}},x_{01})$    &   $\underline{x_{123}}$   &  $\underline{x_{123}}\tlhd x_{01}$       & $\chi(x_{13},x_{01})$           &$\wh{\Id_{x_{01}}}$        \\ \cline{2-6}

   $\odot$                        &      (Table~\ref{isopt2a})                 &  $\underline{x_{123}}$                        &          $\widetilde{\alpha}_{x_{23},x_{12},x_{01}}\bullet(x_{23}\rhd \underline{x_{012}})\bullet \underline{x_{023}}$    &                 $x_{013}$                  &      $\wh{\Id_{x_{01}}}$   \\   \cline{2-6}  

      $\Lambda$                         &          &    $\underline{x_{123}}\tlhd x_{01}$     & $\widetilde{\alpha}_{x_{23},x_{12},x_{01}}\bullet(x_{23}\rhd \underline{x_{012}})\bullet \underline{x_{023}}$                  &   $\underline{x_{013}}$                                   & $\Id_{x_{01}}$         \\   \cline{2-6}

                                 &  $\Delta_{-}(\chi(x_{13},x_{01}))$                                       & $\chi(x_{13},x_{01})$                       &    $x_{013}$      &    $\underline{x_{013}}$    &$\Id_{x_{01}}$      \\ \cline{2-6}

     &  $s_2(\Id_{x_{01}})$                                 &   $\wh{\Id_{x_{01}}}$                        &    $\wh{\Id_{x_{01}}}$          &$\Id_{x_{01}}$     &$\Id_{x_{01}}$    \\   \cline{2-6} 
    \end{tabular}\end{center}
    \end{table}
 \begin{table}[H] \caption{~}\label{isopt2a}\begin{center}
\scalebox{.97}{
    \begin{tabular}{ r | l || l | l | l | l | }
     \cline{2-6}

      & $\Delta_{-}(x_{123})$    &   $\chi(x_{23},x_{12})$   &  $x_{123}$       &$\underline{x_{123}}$          &$\Id_{x_{12}}$        \\ \cline{2-6}

   $\odot$      &      (Table~\ref{isopt3})         &  $\chi(x_{23},x_{12})$       &        $Q$  &        $\widetilde{\alpha}\bullet(x_{23}\rhd \underline{x_{012}})\bullet \underline{x_{023}}$  &      $\chi(x_{12},x_{01})$   \\   \cline{2-6}  

           &    (Def.)  \ \   $=:\Delta_{Q}   $           &    $x_{123}$                   & $=:Q$                  & $x_{013}$        &                  $\chi(x_{12},x_{01})$        \\  \cline{2-6}

       $\Lambda$          &      &  $\underline{x_{123}}$       &          $\widetilde{\alpha}\bullet(x_{23}\rhd \underline{x_{012}})\bullet \underline{x_{023}}$    &     $x_{013}$      &        $\wh{\Id_{x_{01}}}$ \\   \cline{2-6}  
  
     &  $s_1(\chi(x_{12},x_{01}))$                                 &   $\Id_{x_{12}}$                        &      $\chi(x_{12},x_{01})$        &  $\chi(x_{12},x_{01})$  &  $\wh{\Id_{x_{01}}}$  \\   \cline{2-6} 
    \end{tabular}}\end{center}
    \end{table}
$Q$ in  Table~\ref{isopt2} above is defined by face 2 of the table, that is, as the filler of the horn $$[x_{123},\  - ,\ x_{013},\ \chi(x_{12},x_{01})].$$ The verification that this does indeed give a horn is left to the reader. The next three tables are found sideways on the next page, and we continue after these with the final table below.
\begin{sidewaystable}[h,p] \caption{~}\label{isopt3}\begin{center}

    \begin{tabular}{ r | l || l | l | l | l | }
     \cline{2-6}

      & $\Delta_{\widetilde{\alpha}}(x_{23},x_{12},x_{01}  )$    &   $\chi(x_{23},x_{12})$   &  $\chi(x_{23},x_{12}\circ x_{01})$       &$\widetilde{\alpha}_{x_{23},x_{12},x_{01}}$         & $\chi(x_{12},x_{01})$       \\ \cline{2-6}

         $\Lambda$                    &                     &  $\chi(x_{23},x_{12})$                        &        $Q$  &               $\widetilde{\alpha}\bullet(x_{23}\rhd \underline{x_{012}})\bullet \underline{x_{023}}$                   &      $\chi(x_{12},x_{01})$   \\   \cline{2-6}  

    $\odot$                     &        (Table~\ref{isopt4})          &   $\chi(x_{23},x_{12}\circ x_{01})$                                                      & $Q$                  & $(x_{23}\rhd \underline{x_{012}})\bullet \underline{x_{023}}$                    &                  $\Id_{x_{12}\circ x_{01}}$        \\   \cline{2-6}

                            &        $\Delta_{\bullet}(   \widetilde{\alpha}  , (x_{23}\rhd \underline{x_{012}})\bullet \underline{x_{023}}  )$         &   $\widetilde{\alpha}_{x_{23},x_{12},x_{01}}$                        &          $\widetilde{\alpha}\bullet(x_{23}\rhd \underline{x_{012}})\bullet \underline{x_{023}}$    &                 $(x_{23}\rhd \underline{x_{012}})\bullet \underline{x_{023}}$                     &        $\Id_{x_{01}}$ \\   \cline{2-6}  
  
     &  $s_1(\chi(x_{12},x_{01}))$                                 &   $\chi(x_{12},x_{01})$                        &      $\chi(x_{12},x_{01})$        &    $\Id_{x_{12}\circ x_{01}}$  &   $\Id_{x_{01}}$\\   \cline{2-6} 
    \end{tabular}\end{center}
\begin{center} \caption{~}\label{isopt4}
    \begin{tabular}{ r | l || l | l | l | l | }
     \cline{2-6}

      & $\Delta_{-}(x_{23}\trhd \underline{x_{012}})$    &   $\chi(x_{23},x_{12}\circ x_{01})$     &  $x_{23}\trhd \underline{x_{012}}$       &$x_{23}\rhd \underline{x_{012}}$        & $\Id_{x_{12}\circ x_{01}}$       \\ \cline{2-6}

   $\Lambda$                  &                &   $\chi(x_{23},x_{12}\circ x_{01})$                  & $Q$         & $(x_{23}\rhd \underline{x_{012}})\bullet \underline{x_{023}}$              &          $\Id_{x_{12}\circ x_{01}}$  \\\cline{2-6}  

         $\odot$       &     (Table~\ref{isopt5})           &   $x_{23}\trhd \underline{x_{012}}$              & $Q$       & $\underline{x_{023}} $                                             &      $\Id_{x_{12}\circ x_{01}}$         \\   \cline{2-6}

      &        $\Delta_{\bullet}(   x_{23}\rhd \underline{x_{012}} , \underline{x_{023}}  )$         &   $x_{23}\rhd \underline{x_{012}}$    &          $(x_{23}\rhd \underline{x_{012}})\bullet \underline{x_{023}}$       &                 $\underline{x_{023}}$     &        $\Id_{\id_{x_0}}$ \\   \cline{2-6}  
  
     &  $s_1(\Id_{x_{12}\circ x_{01}})$       & $\Id_{x_{12}\circ x_{01}}$                    &    $\Id_{x_{12}\circ x_{01}}$        &    $\Id_{x_{12}\circ x_{01}}$    &   $\Id_{\id_{x_0}}$\\   \cline{2-6} 
    \end{tabular}\end{center}
\begin{center} \caption{~}\label{isopt5}
    \begin{tabular}{ r | l || l | l | l | l | }
     \cline{2-6}

      & $\Delta_{\trhd}(x_{23},\underline{x_{012}})$    &   $\Id_{x_{23}}$     &  $\chi(x_{23},x_{02})$       &$x_{23}\trhd \underline{x_{012}}$        & $\wh{(\underline{x_{012}})}$       \\ \cline{2-6}

        $\odot$         &    (Table~\ref{isopt6}   )                      &   $\Id_{x_{23}}$                  &  $x_{023}$       &    $Q$            &          $\wh{(\underline{x_{012}}) }$  \\\cline{2-6}  

                    &   $\Delta_{-}(x_{023})$      &   $\chi(x_{23},x_{02})$              &  $x_{023}$         & $\underline{x_{023}} $          &      $\Id_{x_{02}}$         \\   \cline{2-6}

    $\Lambda$  &               &   $x_{23}\trhd \underline{x_{012}}$    &          $Q$       &                 $\underline{x_{023}}$     &        $\Id_{x_{12}\circ x_{01}}$ \\   \cline{2-6}  
  
     &  $s_0(\wh{(\underline{x_{012}})})$       & $\wh{(\underline{x_{012}}) }$                   &   $\wh{(\underline{x_{012}}) }$      &    $\Id_{x_{02}}$    &   $\Id_{x_{12}\circ x_{01}}$ \\   \cline{2-6} 
    \end{tabular}
\end{center}
    \end{sidewaystable}
\begin{table}[H]
\begin{center} \caption{~}\label{isopt6}
     \begin{tabular}{ r | l || l | l | l | l | }
     \cline{2-6}

      & $s_1(x_{123})$    &   $\Id_{x_{23}}$     &  $x_{123}$       &$x_{123}$        & $\wh{\Id_{x_{12}}}$       \\ \cline{2-6}

         $\Lambda$       &        &   $\Id_{x_{23}}$                  &  $x_{023}$       &    $Q$            &          $\wh{(\underline{x_{012}}) }$  \\ \cline{2-6}  

           $\odot$           & (hypothesis) & $x_{123}$              &$x_{023}$           &$x_{013}$            &      $x_{012}$         \\   \cline{2-6}

                       &     $\Delta_{Q}$   &   $x_{123}$     &          $Q$       &                $x_{013}$    &       $\chi(x_{12},x_{01})$\\   \cline{2-6}  
  
                      &  $\Delta_{\stackrel{\wedge}{-}}(x_{012})$       &$\wh{\Id_{x_{12}}}$                 & $\wh{(\underline{x_{012}}) }$       &   $x_{012}$      &   $\chi(x_{12},x_{01})$ \\   \cline{2-6} 
    \end{tabular}\end{center}

\end{table}
This proves the commutativity of \ref{condsphere}, given the commutativity of $[x_{123},\ x_{023},\ x_{013},\ x_{012}].$ For the other direction, if we instead interpret the last table as a horn in which we wish to fill face 2, then the reading the proof backwards, switching the symbols $\odot$ and $\Lambda$, we see the face $[x_{123},\ x_{023},\ x_{013},\ x_{012}]$ is commutative, given the hypothesis that \ref{condsphere} is commutative. This completes the proof that the canonical map $u:X \ra N(\Bic(X))$ is an isomorphism. 
\begin{theorem}[Duskin] \label{bictheorem} A simplicial set $X$ is $2$-reduced inner-Kan, i.e. meets the inner Kan condition for all $n$ and meets it uniquely for $n>2$, if and only if it is isomorphic to the nerve of a small $(2,1)$-category.
\end{theorem}
\begin{proof} We showed that $N(B)$ is $2$-reduced inner-Kan as Theorem~\ref{daggertheorem} and we have shown in this section that if $X$ is $2$-reduced inner-Kan, then $X \cong  N(\Bic(X))$.

Finally, to show that the isomorphism $u$ is strict, we must show $u$ sends the algebraic structure $\chi$ of $X$ to the natural algebraic structure of $N(\Bic(X))$. This is equivalent to showing $\underline{\chi(g,f)}= \Id_{g\circ f}$, which is done above as Lemma~\ref{lemma2}.\end{proof}
%
\section{Promoting $N$ and $\Bic$ to functors \label{functorsection}}
We now consider how $N$ and $\Bic$ can be defined, respectively, for either  strictly identity preserving functors between $(2,1)$-categories or for morphisms of algebraic $2$-reduced inner-Kan simplicial sets. Our constructions are a special case of those given in \cite{Gur09}.
First recall the definition of a (weak) functor between bicategories. These are also often called \emph{pseudo-functors} or \emph{homomorphisms}.
\begin{definition} If $\B$ and $\B'$ are bicategories, then a functor $(F,\phi, \upsilon):\B \ra \B'$ consists of:
\begin{itemize}
\item A mapping $F$ from objects, $1$-morphisms, and $2$-morphisms of $\B$ and to those of $\B'$:
\item For all $a \stackrel{f}{\ra} b \stackrel{g}{\ra} c$ in $\B$, an invertible $2$-morphism $\phi_{g,f}: F(g\circ f)\Rightarrow F(g) \circ F(f) $, collectively called the \emph{distributor} of $F$.
\item For every object $a$ of $\B$, an invertible $2$-morphism $\upsilon_a : F(\id_a)\Rightarrow \id_{F(a)} $, collectively called the \emph{unitor} of $F$.
\end{itemize}
Such that the following axioms are satisfied ($\alpha, \rho, \gamma$ denote associator and unitors of $\B,$ and $\alpha', \rho', \gamma'$ are the associator and unitors of $\B'$):
\begin{itemize}
\item \emph{F is strictly functorial with respect to $\bullet$}
\begin{description}
\item[BFun1.]  For all $f: a \ra b,$ in $\B,$ \tab  $F(\Id_f)=\Id_{F(f)}$ 
\item[BFun2.]  For all $f \stackrel{\eta}{\Rightarrow} g \stackrel{\theta}{\Rightarrow} h$, \tab $F(\theta)\bullet F(\eta) = F(\theta \bullet \eta)$
\end{description}
\item \emph{naturality of the distributor}
\begin{description}
\item[BFun3.]  For all $f\stackrel{\eta}{\Rightarrow}g:a\ra b$ and $h:b \ra c$ in $\B,$

 \tab  $ \phi_{h,g} \bullet  F(h \rhd \eta) = (F(h) \rhd F(\eta)) \bullet \phi_{h,f} $
\item[BFun4.]  For all $f:a \ra b$ and $g \stackrel{\eta}{\Rightarrow} h: b \ra c$  in $\B,$

 \tab      $ \phi_{h,f} \bullet  F(\eta\lhd f) = (F(\eta) \lhd F(f)) \bullet \phi_{g,f}$
\end{description}
\item\emph{compatibility of the distributor with the associators}
\begin{description}
\item[BFun5.] For all $a \stackrel{f}{\ra} b \stackrel{g}{\ra} c \stackrel{h}{\ra} d$ in $\B$ \tab $  \alpha'_{F(h),F(g),F(f)}  \bullet ( F(h)\rhd\phi_{g,f} ) \bullet  \phi_{h, g\circ f}$ 

\tab $\quad =    (\phi_{h,g} \lhd F(f)) \bullet\phi_{h\circ g,f}    \bullet  F(\alpha_{h,g,f}) $
\end{description}
\item\emph{compatibility of the unitors of $F$ and the unitors of $\B$ and $\B'$}
\begin{description}
\item[BFun6.]  For all $f:a \ra b$ in $\B$, \tab    $ \rho'_{F(f)} = (F(f)\rhd \upsilon_a) \bullet \phi_{  f ,\id_a} \bullet F(\rho_f) $
\item[BFun7.]  For all $f:a \ra b$ in $\B$, \tab $ \lambda'_{F(f)} = (\upsilon_b \lhd F(f)) \bullet \phi_{\id_b  ,f } \bullet F(\lambda_f)  $
\end{description}
\end{itemize}
A functor $(F, \phi,\upsilon)$ between bicategories is called \emph{strict} if $\phi$ and $\upsilon$ are identities. In this case, the axioms \textbf{BFun1}-\textbf{BFun7} are equivalent $F$ strictly preserving all compositions, both identities, the associator, and the unitors. If instead only $\upsilon$ is required to consist of identities, we say $F$ is \emph{strictly identity-preserving}.
\end{definition}
The composition $GF$ of functors $(F,\phi,\upsilon)$ and $(G,\phi', \upsilon')$ is given by composing the action of $G$ and $F$ on objects, $1$-morphisms, and $2$-morphisms. The distributor of $GF$ is given by $G(\phi_{h,f})\bullet \phi'_{F(h), F(f)}:$   $$GF(h\circ f)\stackrel{G(\phi_{h,f})}{\ra} G(F(h)\bullet F(f))\stackrel{\phi'_{F(h), F(f)}}{\ra} GF(h)\circ GF(f).$$
The unitor of $GF$ is similarly given by $\upsilon'_{F(a)}\circ G(\upsilon_a).$

\subsection{The bicategory functor $\Bic(F)$}
First let  $F:(X,\chi)\ra (Y,\chi')$ be a (not necessarily strict) morphism between two algebraic $2$-reduced inner-Kan simplicial sets.  We define a functor of $(2,1)$-categories $\Bic(F):\Bic(X)\ra \Bic(Y)$ as follows:
\begin{itemize}
\item $\Bic(F)$ is defined by $F$ for objects and $1$-morphisms, since the objects and $1$-morphisms of $\Bic(X)$ and $\Bic(Y)$ are the same as the $0$-cells and $1$-cells of $X$ and $Y$.
\item Clearly $F$ sends those $2$-cells of $X$ which are $2$-morphisms of $\Bic(X)$ to $2$-cells of $Y$ which are $2$-morphisms of $\Bic(Y)$. This defines $\Bic(F)$ on $2$-morphisms.
\item Note that $\Bic(F)$ is strictly identity-preserving since $F$ respects degeneracy maps. Thus we can define the unitor of $\Bic(F)$ to be the identity.
\item Let $a \stackrel{f}{\ra} b \stackrel{g}{\ra} c$ be $1$-morphisms in $\Bic(X)$. Then the distributor $$\phi_{g,f}: F(g \circ f)\Rightarrow F(g) \circ F(f) $$ of $\Bic(Y)$ is defined by $$\phi_{g,f} =\underline{F(\chi(g,f))}$$
By comparing  $\Delta_{-} (F(\chi(g,f)))$ with $\Delta_\bullet (\chi'(F(g),F(f)),\phi_{g,f})$ and applying the Matching Lemma it is easy to see $$F(\chi(g,f)) = \chi'(F(g),F(f))\bullet \phi_{g,f}.$$
\end{itemize}
\begin{theorem} $\Bic(F)$ as defined above is a functor of $(2,1)$-categories. 
\end{theorem} 
\begin{proof}~
\begin{itemize}
\item \textbf{BFun1} follows from the fact that $F$ preserves degeneracy maps.
\item For \textbf{BFun2}, we can easily see $F(\Lambda_\bullet(\theta, \eta))= \Lambda_\bullet ( F(\theta), F(\eta) )$. Then $F(\Delta_\bullet (\theta, \eta)$ fills $\Lambda_\bullet ( F(\theta), F(\eta) )$ from which we see $F(\Delta_\bullet (\theta,\eta))=\Delta_\bullet (F(\theta),F(\eta))$ and in particular $F(\theta)\bullet F(\eta) = F(\theta \bullet \eta)$.
\item For \textbf{BFun3} let  $f\stackrel{\eta}{\Rightarrow}g:a\ra b$ and $h:b \ra c$ in $\Bic(X).$ We must show $$\phi_{h,g} \bullet  F(h \rhd \eta) = (F(h) \rhd F(\eta)) \bullet \phi_{h,f}.$$ We have $$\Delta_\bullet (\phi_{h,f}, F(h \rhd \eta))=[\phi_{h,g}, \ \phi_{h,g} \bullet F(h \rhd \eta), \  F(h \rhd \eta)  ,\ \Id_{\id_a}]$$ so by the Matching Lemma it suffices to show that the sphere $$\Delta_\bullet (\phi_{h,g}, F(h \rhd \eta))=[\phi_{h,g}, \ (F(h) \rhd F(\eta)) \bullet \phi_{h,f}, \  F(h \rhd \eta)  ,\ \Id_{\id_a}]$$ is commutative. The following Glenn table proof verifies this commutativity:

\begin{table}[H] \caption{~\label{Bfun3proof1}}\begin{center}
\scalebox{.96}{    \begin{tabular}{ r | l || l | l | l | l | }
     \cline{2-6}

    & $\Delta_{-} (F(\chi(g,f)))$       & $\chi'(F(h),F(g))$   &    $F(\chi(h,g))$     &  $\phi_{h,g}$            & $\Id_{F(g)}$      \\ \cline{2-6}
     
  $\odot$   & (Table~\ref{Bfun3proof2})   &  $\chi'(F(h),F(g))$   &  $F(\chi(h,f))$   & $(F(h)\rhd F(\eta))\bullet \phi_{h,f}$   &       $F(\eta)$          \\   \cline{2-6}  

    & $F(\Delta_{\rhd}(h,\eta))$       &   $F(\chi(h,g))$    &    $F(\chi(h,f))$     &    $ F(h \rhd \eta) $   &  $F(\eta)$       \\   \cline{2-6}  
     
  $\Lambda$  &    & $\phi_{h,g}$&  $(F(h) \rhd F(\eta))\bullet \phi_{h,f}$      & $ F(h \rhd \eta)$        &  $\Id_{\id_a}$ \\ \cline{2-6}

    &   $s_1(F(\eta))$   &  $\Id_{F(g)}$&     $F(\eta)$&      $F(\eta)$   & $\Id_{\id_a}$    \\   \cline{2-6} 
    \end{tabular}}\end{center}
    \end{table}

\begin{table}[H] \caption{~\label{Bfun3proof2}}\begin{center}\scalebox{0.88}{
    \begin{tabular}{ r | l || l | l | l | l | }
     \cline{2-6}

         & $\Delta_\rhd(F(h),F(g))$  & $\chi'(F(h),F(g))$   &    $\chi'(F(h),F(f))$     &  $F(h)\rhd F(\eta)$            & $F(\eta)$      \\ \cline{2-6}
     
  $\Lambda$                 &    &  $\chi'(F(h),F(g))$   &  $F(\chi(h,f))$   & $(F(h)\rhd F(\eta))\bullet \phi_{h,f}$   &       $F(\eta)$          \\   \cline{2-6}  

    & $\Delta_{-} (F(\chi(h,f)))$     &   $\chi'(F(h),F(f))$    &    $F(\chi(h,f))$     &    $ \phi_{h,f} $   &  $\Id_{F(g)}$       \\   \cline{2-6}  
     
  & $\Delta_\bullet(F(h)\rhd F(\eta), \phi_{h,f})$   & $F(h)\rhd F(\eta)$   &  $(F(h) \rhd F(\eta))\bullet \phi_{h,f}$      &  $ \phi_{h,f} $         &  $\Id_{\id_a}$ \\ \cline{2-6}

    &   $s_0(F(\eta))$   &  $F(\eta)$  &     $F(\eta)$&       $\Id_{F(g)}$    & $\Id_{\id_a}$    \\   \cline{2-6} 
    \end{tabular}   }\end{center}
    \end{table}
 
\item For \textbf{BFun4}, note that for a $2$-morphism $\eta$ in $X$ we have $F(\Lambda_\wedge (\eta))= \Lambda_\wedge (F(\eta))$ so $F(\Delta_\wedge (\eta))$ fills $\Lambda_\wedge (F(\eta))$ and by uniqueness of fillers in $Y$ we have $F(\Delta_\wedge (\eta))=\Delta_\wedge (F(\eta))$ and in particular $F(\wh{\eta})=\wh{F(\eta)}.$
Now to verify \textbf{BFun4}, suppose we have  $f:a \ra b$ and $g \stackrel{\eta}{\Rightarrow} h: b \ra c$  in $\Bic(X).$

We must show $ \phi_{h,f} \bullet  F(\eta\lhd f) = (F(\eta) \lhd F(f)) \bullet \phi_{g,f}$. We have $$\Delta_\bullet (\phi_{h,f},  F(\eta\lhd f))=[\phi_{h,f}, \ (F(\eta) \lhd F(f)) \bullet \phi_{g,f}, \  F(\eta\lhd f)  ,\ \Id_{\id_a}]$$ so by the Matching Lemma it suffices to show the commutativity of the sphere $$\Delta_\bullet (\phi_{h,g}, F(h \rhd \eta))=[\phi_{h,f}, \ (F(\eta) \lhd F(f)) \bullet \phi_{g,f}, \  F(\eta\lhd f)  ,\ \Id_{\id_a}].$$  The following Glenn table proof verifies this commutativity:

\begin{table}[H] \caption{~\label{Bfun4proof1}}\begin{center}
  \scalebox{.94}{  \begin{tabular}{ r | l || l | l | l | l | }
     \cline{2-6}

    & $\Delta_{-} (F(\chi(h,f)))$        & $\chi'(F(h),F(f))$   &    $F(\chi(h,f))$     &  $\phi_{h,f}$            & $\Id_{F(f)}$      \\ \cline{2-6}
     
  $\odot$   & (Table~\ref{Bfun4proof2})   &  $\chi'(F(h),F(f))$   &  $F(\eta\tlhd f)$   &$(F(\eta) \lhd F(f)) \bullet \phi_{g,f}$   &       $\Id_{F(f)}$          \\   \cline{2-6}  

    & $F(\Delta_{-}(\eta\tlhd f))$       &   $F(\chi(h,f))$   &   $F(\eta\tlhd f)$    &   $F(\eta\lhd f)$  &  $F(\Id_f)$       \\   \cline{2-6}  
     
  $\Lambda$  &    & $\phi_{h,f}$&  $(F(\eta) \lhd F(f)) \bullet \phi_{g,f}$      & $F(\eta\lhd f)$        &  $\Id_{\id_a}$ \\ \cline{2-6}

    &     $s_0(\Id_{F(f)})$   &    $\Id_{F(f)}$ &      $\Id_{F(f)}$ &        $F(\Id_f)=\Id_{F(f)}$   & $\Id_{\id_a}$    \\   \cline{2-6} 
    \end{tabular}}\end{center}
    \end{table}

\begin{table}[H] \caption{~\label{Bfun4proof2}}\begin{center}
   \scalebox{.90}{  \begin{tabular}{ r | l || l | l | l | l | }
     \cline{2-6}

    & $\Delta_{-}( F(\eta)\tlhd F(f))$        & $\chi'(F(h),F(f))$   &    $F(\eta)\tlhd F(f)$     &  $F(\eta)\lhd F(f)$            & $\Id_{F(f)}$      \\ \cline{2-6}
     
  $\Lambda$   &    &  $\chi'(F(h),F(f))$   &  $F(\eta\tlhd f)$   &$(F(\eta) \lhd F(f)) \bullet \phi_{g,f}$   &       $\Id_{F(f)}$          \\   \cline{2-6}  

     $\odot$   &(Table~\ref{Bfun4proof3})     &   $F(\eta)\tlhd F(f)$    &   $F(\eta\tlhd f)$    &    $\phi_{g,f}$   &  $\Id_{F(f)}$       \\   \cline{2-6}  
     
                                               & $\Delta_{\bullet}(F(\eta) \lhd F(f),  \phi_{g,f})$  & $F(\eta) \lhd F(f)$&  $(F(\eta) \lhd F(f)) \bullet \phi_{g,f}$      & $\phi_{g,f}$        &  $\Id_{\id_a}$ \\ \cline{2-6}
 
                                              &     $s_0(\Id_{F(f)})$   &    $\Id_{F(f)}$ &      $\Id_{F(f)}$ &        $\Id_{F(f)}$   & $\Id_{\id_a}$    \\   \cline{2-6} 
    \end{tabular}}\end{center}
    \end{table}
 
\begin{table}[H] \caption{~\label{Bfun4proof3}}\begin{center}
    \begin{tabular}{ r | l || l | l | l | l | }
     \cline{2-6}

    & $\Delta_{\tlhd}(F(\eta), F(f))$                      & $F(\eta)$   &    $F(\eta)\tlhd F(f)$     &  $\chi'(F(g),F(f))$            & $\wh{\Id_{F(f)}}$      \\ \cline{2-6}
     
                             &   $F(  \Delta_{\tlhd}(\eta, f) )$            &  $F(\eta)$   &  $F(\eta\tlhd f)$       &$F(\chi(g,f))$   &       $F(\wh{\Id_f})$          \\   \cline{2-6}  

       $\Lambda$      &                             &   $F(\eta)\tlhd F(f)$    &   $F(\eta\tlhd f)$    &    $\phi_{g,f}$   &  $\Id_{F(f)}$       \\   \cline{2-6}  
     
                                              &$\Delta_{-} (F(\chi(g,f)))$   &  $\chi'(F(g),F(f))$  &  $F(\chi(g,f))$    & $\phi_{g,f}$        &   $\Id_{F(f)}$     \\ \cline{2-6}
 
                                              &  $s_2(\Id_{F(f)})$   &     $\wh{\Id_{F(f)}}$ &       $F(\wh{\Id_f})=\wh{\Id_{F(f)}}$      &        $\Id_{F(f)}$   & $\Id_{F(f)}$    \\   \cline{2-6} 
    \end{tabular}\end{center}
    \end{table}   
    
\item For \textbf{BFun5}, take $a \stackrel{f}{\ra} b \stackrel{g}{\ra} c \stackrel{h}{\ra}$ in $\Bic(X).$

We must show    $$\alpha'_{F(h),F(g),F(f)}\bullet(F(h)\rhd \phi_{g,f})\bullet \phi_{h, g\circ f}= (\phi_{h,g}\lhd F(f))\bullet \phi_{h\circ g,f}\bullet F(\alpha_{h,g,f}). $$ By Lemma~\ref{bulletylemma} it suffices to show  $$\widetilde{\alpha}'_{F(h),F(g),F(f)}\bullet(F(h)\rhd \phi_{g,f})\bullet \phi_{h,g\circ f }= (\phi_{h,g}\tlhd F(f))\bullet \phi_{h\circ g,f}\bullet F(\alpha_{h,g,f}). $$ We have 
\begin{multline*} \Delta_\bullet(\widetilde{\alpha}'_{F(h),F(g),F(f)}, (F(h)\rhd \phi_{g,f})\bullet \phi_{h, g\circ f}) = \\ [ \widetilde{\alpha}'_{F(h),F(g),F(f)},\  \widetilde{\alpha}'_{F(h),F(g),F(f)}\bullet (F(h)\rhd \phi_{g,f})\bullet \phi_{h, g\circ f}  ,\ (F(h)\rhd \phi_{g,f})\bullet \phi_{h,g\circ f } ,\ \Id_{F(f)}       ],
\end{multline*}
so by the Matching Lemma it suffices to show commutativity for the sphere 
\begin{multline*} \Delta_\bullet(\widetilde{\alpha}'_{F(h),F(g),F(f)}, (F(h)\rhd \phi_{g,f})\bullet \phi_{h,g\circ f }) = \\ [ \widetilde{\alpha}'_{F(h),F(g),F(f)},\  (\phi_{h,g}\tlhd F(f))\bullet \phi_{h\circ g,f}\bullet F(\alpha_{h,g,f})  ,\ (F(h)\rhd \phi_{g,f})\bullet \phi_{h,g\circ f } ,\ \Id_{F(f)}].
\end{multline*}
We show this using a Glenn table proof, given in Tables~\ref{Bfun5proof1}--\ref{Bfun5proof7}. For the purpose of abbreviation, we define $$ A =  (\phi_{h,g}\tlhd F(f))\bullet \phi_{h\circ g,f}\bullet F(\alpha_{h,g,f}).$$

\begin{sidewaystable}[h,p] \caption{~\label{Bfun5proof1}}\begin{center}
\scalebox{.95}{    \begin{tabular}{ r | l || l | l | l | l | }
     \cline{2-6}
     
    & $\Delta_{\widetilde{\alpha}}(F(h), F(g), F(f))$                      & $\chi'(F(h),F(g))$   &    $\chi'(F(h), F(g)\circ F(f))$     &  $\widetilde{\alpha}'_{F(h),F(g),F(f)}$            & $\chi'(F(g),F(f))$      \\ \cline{2-6}
      
    $\odot$       &    (Table~\ref{Bfun5proof2})         &  $\chi'(F(h),F(g))$   &   $(F(h)\trhd \phi_{g,f} )\bullet \phi_{h,g\circ f}$       &   $A$  &      $\chi'(F(g),F(f))$          \\   \cline{2-6}  

           &   $\Delta_{-}((F(h)\rhd \phi_{g,f})\bullet \phi_{h,g\circ f } )$                &   $\chi'(F(h), F(g)\circ F(f))$  &    $(F(h)\trhd \phi_{g,f})\bullet \phi_{h, g\circ f}$      &   $(F(h)\rhd \phi_{g,f})\bullet \phi_{h, g\circ f}$   &  $\Id_{F(g)\circ F(f)}$       \\   \cline{2-6}  
     
   $\Lambda$                   &      &   $\widetilde{\alpha}'_{F(h),F(g),F(f)}$  & $A$    & $(F(h)\rhd \phi_{g,f})\bullet \phi_{ h,g\circ f}$        &   $\Id_{F(f)}$     \\ \cline{2-6}
 
                     &  $s_0(\chi'(F(g),F(f)))$   &   $\chi'(F(g),F(f))$       &       $\chi'(F(g),F(f))$      &        $\Id_{F(g)\circ F(f)}$   & $\Id_{F(f)}$    \\   \cline{2-6} 
    \end{tabular}}\end{center}
    \begin{center} \caption{~\label{Bfun5proof2}}
    \begin{tabular}{ r | l || l | l | l | l | }
     \cline{2-6}

      & $\Delta_{-}(F(\chi(h,g)))$    &  $\chi'(F(h),F(g))$     &  $F(\chi(h,g))$       &$\phi_{h,g}$        & $\Id_{F(g)}$       \\ \cline{2-6}

   $\Lambda$                  &                &   $\chi'(F(h),F(g))$        & $(F(h)\trhd \phi_{g,f})\bullet \phi_{h,g\circ f}$         & $A$             &          $\chi'(F(g), F(f))$  \\\cline{2-6}  

         $\odot_1$       &     (Table~\ref{Bfun5proof3})           &    $F(\chi(h,g))$               & $(F(h)\trhd \phi_{g,f})\bullet \phi_{h,g\circ f}$       & $F(\widetilde{\alpha}_{h,g,f}) $                                             &      $\chi'(F(g), F(f))$        \\   \cline{2-6}

    $\odot_2$      &       (Table~\ref{Bfun5proof6})         &  $\phi_{h,g}$   &     $A$    &                  $F(\widetilde{\alpha}_{h,g,f}) $     &        $\widehat{\Id_{F(f)}}$ \\   \cline{2-6}  
  
     &  $s_1(\chi'(F(g), F(f)))$       & $\Id_{F(g)}$                   &    $\chi'(F(g), F(f))$       &     $\chi'(F(g), F(f))$    &  $\widehat{\Id_{F(f)}}$ \\   \cline{2-6} 
    \end{tabular}\end{center}
    \begin{center} \caption{Commutativity of $\odot_1$ in Table~\ref{Bfun5proof2} \label{Bfun5proof3}}
  \scalebox{.97}{  \begin{tabular}{ r | l || l | l | l | l | }
     \cline{2-6}
     
      & $s_1(F(\chi(h,g)))$  &  $\Id_{F(h)}$     &  $F(\chi(h,g))$       &$F(\chi(h,g))$        & $\wh{\Id_{F(g)}}$       \\ \cline{2-6}

   $\odot_3$      &  (Table~\ref{Bfun5proof4})            &  $\Id_{F(h)}$       & $F(\chi(h,g\circ f))$         & $(F(h)\trhd \phi_{g,f})\bullet \phi_{h,g\circ f}$             &          $\wh{\phi_{g,f}}$  \\\cline{2-6}  

               &  $F(\Delta_{\wh{\alpha}}(h,g,f)) $         &   $F(\chi(h,g))$               & $F(\chi(h,g\circ f))$       & $F(\widetilde{\alpha}_{h,g,f}) $            &      $F(\chi(g,f))$        \\   \cline{2-6}  
        
     $\Lambda$       &        &  $F(\chi(h,g))$    &     $(F(h)\trhd \phi_{g,f})\bullet \phi_{h,g\circ f}$    &                  $F(\widetilde{\alpha}_{h,g,f}) $     &        $\chi'(F(g),F(f))$ \\   \cline{2-6}  
  
  $\odot_4$    &   (Table~\ref{Bfun5proof5})       & $\wh{\Id_{F(g)}}$                   &    $\wh{\phi_{g,f}}$       &     $F(\chi(g,f))$   &   $\chi'(F(g),F(f))$  \\   \cline{2-6} 
    \end{tabular}}\end{center}
     \end{sidewaystable}      
\begin{sidewaystable}[h,p] 
    \begin{center} \caption{Commutativity of $\odot_3$ in Table~\ref{Bfun5proof3} \label{Bfun5proof4}}
    \begin{tabular}{ r | l || l | l | l | l | }
     \cline{2-6}
     
      & $\Delta_{\trhd} (F(h),\phi_{g,f} )$  &  $\Id_{F(h)}$     &  $\chi'(F(h), F(g\circ f))$       &$\chi'(F(h)\trhd \phi_{g,f})$        & $\wh{\phi_{g,f}}$       \\ \cline{2-6}

    $\Lambda$       &                &  $\Id_{F(h)}$       & $F(\chi(h,g\circ f))$         & $(F(h)\trhd \phi_{g,f})\bullet \phi_{h,g\circ f}$             &          $\wh{\phi_{g,f}}$  \\\cline{2-6}  

               &   $\Delta_{-} (F(\chi(h,g\circ f)))$         &  $\chi'(F(h), F(g\circ f))$               & $F(\chi(h,g\circ f))$       & $\phi_{h, g\circ f} $      &      $\Id_{F(g\circ f)}$      \\   \cline{2-6}  
        
       &  $\Delta_{\bullet}( F(h)\trhd \phi_{g,f}, \phi_{h,g\circ f} ) $     &  $\chi'(F(h)\trhd \phi_{g,f})$   &     $(F(h)\trhd \phi_{g,f})\bullet \phi_{h,g\circ f}$    &       $\phi_{ h,g\circ f} $     &   $\Id_{F(g)\circ F(f)}$  \\   \cline{2-6}  
  
     &  $s_0  ( \wh{\phi_{g,f}}) $ & $\wh{\phi_{g,f}}$                      &   $\wh{\phi_{g,f}}$         &     $\Id_{F(g\circ f)}$   &   $\Id_{F(g)\circ F(f)}$  \\   \cline{2-6} 
    \end{tabular}\end{center}    
    \begin{center} \caption{Commutativity of $\odot_4$ in Table~\ref{Bfun5proof3}\label{Bfun5proof5}}
    \begin{tabular}{ r | l || l | l | l | l | }
     \cline{2-6}
     
      & $s_2(\chi'(F(g),F(f))))$  &  $\wh{\Id_{F(g)}}$     &   $\wh{\Id_{F(g)\circ F(f)}}$     &$\chi'(F(g),F(f))$        & $\chi'(F(g),F(f))$        \\ \cline{2-6}

    $\Lambda$       &                &  $\wh{\Id_{F(g)}}$       & $\wh{\phi_{g,f}}$         & $F(\chi(g,f))$             &          $\chi'(F(g),F(f))$  \\\cline{2-6}  

               &  $\Delta_{\wedge}(\phi_{g,f}) $         &   $\wh{\Id_{F(g)\circ F(f)}}$               & $\wh{\phi_{g,f}}$        & $\phi_{g,f}$      &      $\Id_{F(g)\circ F(f)}$      \\   \cline{2-6}  
        
       &     $\Delta_{-} (F(\chi(g, f)))$     &  $\chi'(F(g),F(f))$   &     $F(\chi(g,f))$    &       $\phi_{g,f}$      &   $\Id_{F(f)}$  \\   \cline{2-6}  
  
     &  $s_0  ( \chi'(F(g),F(f))) $ & $\chi'(F(g),F(f))$                 &    $\chi'(F(g),F(f))$        &     $\Id_{F(g)\circ F(f)}$   &   $\Id_{F(f)}$  \\   \cline{2-6} 
    \end{tabular}\end{center}    
        \begin{center} \caption{Commutativity of $\odot_2$ in Table~\ref{Bfun5proof2}, part 1 \label{Bfun5proof6}}
\scalebox{.96}{ \begin{tabular}{ r | l || l | l | l | l | }
     \cline{2-6}

     $\odot$     &   (Table~\ref{Bfun5proof7})  &  $\phi_{h,g}$     &  $(\phi_{h,g}\tlhd F(f))\bullet \phi_{h\circ g,f}$       &$F(\chi(h\circ g,f))$        & $\wh{\Id_{F(f)}}$       \\ \cline{2-6}

   $\Lambda$        &                &   $\phi_{h,g}$         & $(\phi_{h,g}\tlhd F(f))\bullet \phi_{h\circ g,f}\bullet F(\alpha_{h,g,f})$         & $F(\widetilde{\alpha}_{h,g,f})$             &         $\wh{\Id_{F(f)}}$     \\   \cline{2-6}  

            &     $\Delta_{\bullet}((\phi_{h,g}\tlhd F(f))\bullet \phi_{h\circ g,f} ,\ F(\alpha_{h,g,f}) )$         &   $(\phi_{h,g}\tlhd F(f))\bullet \phi_{h\circ g,f}$                 &  $(\phi_{h,g}\tlhd F(f))\bullet \phi_{h\circ g,f}\bullet F(\alpha_{h,g,f})$         & $F(\alpha_{h,g,f}) $                                             &      $\Id_{F(f)}$        \\   \cline{2-6}

       &    $F(\Delta_{-}(\widetilde{\alpha}_{h,g,f}))$      &  $F(\chi(h\circ g,f))$   &    $F(\widetilde{\alpha}_{h,g,f})$     &      $F(\alpha_{h,g,f}) $    &       $\Id_{F(f)}$   \\   \cline{2-6}  
  
     &  $s_2(\Id_{F(f)})$       &$\wh{\Id_{F(f)}}$                     &  $\wh{\Id_{F(f)}}$         &    $\Id_{F(f)}$  &   $\Id_{F(f)}$ \\   \cline{2-6} 
    \end{tabular}}\end{center}
    \end{sidewaystable}   
  \begin{table}[h,p,t,b]
     \begin{center} \caption{Commutativity of $\odot_2$ in Table~\ref{Bfun5proof2}, part 2 \label{Bfun5proof7}}
    \scalebox{.90}{\begin{tabular}{ r | l || l | l | l | l | }
     \cline{2-6}

         &   $\Delta_{\tlhd}(\phi_{h,g}, F(f))$ &  $\phi_{h,g}$     &  $\phi_{h,g} \tlhd F(f)$      &$\chi'(F(h\circ g), F(f))$        & $\wh{\Id_{F(f)}}$       \\ \cline{2-6}

   $\Lambda$        &                &   $\phi_{h,g}$         &  $(\phi_{h,g}\tlhd F(f))\bullet \phi_{h\circ g,f}$          & $F(\chi(h\circ g, f))$             &         $\wh{\Id_{F(f)}}$     \\   \cline{2-6}  

            &     $\Delta_{\bullet}(\phi_{h,g} \tlhd F(f), \phi_{h\circ g,f}) $         &    $\phi_{h,g} \tlhd F(f)$                     &$(\phi_{h,g}\tlhd F(f))\bullet \phi_{h\circ g,f}$      & $\phi_{h\circ g,f}$      &      $\Id_{F(f)}$        \\   \cline{2-6}

       &    $\Delta_{-}(F(\chi(h\circ g,f)))$      & $\chi'(F(h\circ g), F(f))$  &     $F(\chi(h\circ g, f))$     &      $\phi_{h\circ g,f}$    &       $\Id_{F(f)}$   \\   \cline{2-6}  
  
     &  $s_2(\Id_{F(f)})$       &$\wh{\Id_{F(f)}}$                     &  $\wh{\Id_{F(f)}}$         &    $\Id_{F(f)}$  &   $\Id_{F(f)}$ \\   \cline{2-6} 
    \end{tabular}}\end{center}  \end{table} 
\item For \textbf{BFun6}  let $f:a \ra b$ in $\Bic(X)$. We must show   $$ \rho'_{F(f)} = (F(f)\rhd \upsilon_a) \bullet \phi_{ f, \id_a } \bullet F(\rho_f).$$ Recall that we have defined $\upsilon_a = \Id_{\id_{F(a)}}$ and we have $$\rho'_{F(f)} = \underline{\Id_{F(f)}}  \ \ \ \ \ \ \ \ \ \ \ \  \ F(\rho_f)=F(\underline{\Id_f}). $$ By Lemma~\ref{bulletylemma} it suffices to show  $$ \Id_{F(f)} = (F(f)\trhd \upsilon_a) \bullet \phi_{ f, \id_a} \bullet F(\underline{\Id_f}).$$ We have
\begin{multline*} \Delta_{\bullet}(  F(f)\trhd \Id_{\id_{F(a)}} ,\  \phi_{ f ,\id_a } \bullet F(\underline{\Id_f}) ) =  \\ [F(f)\trhd  \Id_{\id_{F(a)}},\ (F(f)\trhd  \Id_{\id_{F(a)}}) \bullet \phi_{ f ,\id_a } \bullet F(\underline{\Id_f})  ,\  \phi_{ f,\id_a  } \bullet F(\underline{\Id_f})    ,\  \Id_{\id_{F(a)}} ]\end{multline*} so by the Matching Lemma, it suffices to show the sphere \begin{equation}\label{BFun6eq} [F(f)\trhd  \Id_{\id_{F(a)}},\  \Id_{F(f)}  ,\  \phi_{f ,\id_a  } \bullet F(\underline{\Id_f})    ,\  \Id_{\id_{F(a)}} ]\end{equation} is commutativity. First we observe that
 \begin{align*} 
 d(\Delta_{\trhd} (F(f),\Id_{\id_{F(b)}})) &=  [ \Id_{F(f)}   ,\   \chi'(F(f), \id_{F(a)})   ,\   F(f) \trhd \Id_{\id_{F(a)}} ,\ \wh{ \Id_{\id_{F(a)}}}] \\
 s_1(  \chi'(F(f), \id_{F(a)}) ) &= [ \Id_{F(f)},\   \chi'(F(f), \id_{F(a)})   ,\    \chi'(F(f), \id_{F(a)})       ,\ \wh{ \Id_{\id_{F(a)}}}].
 \end{align*}
Thus by the Matching Lemma $ F(f) \trhd \Id_{\id_{F(a)}}= \chi'(F(f), \id_{F(a)})$. So making this substitution into sphere~\ref{BFun6eq}, we must show the commutativity of \begin{equation*} [\chi'(F(f), \id_{F(a)}),\  \Id_{F(f)}  ,\  \phi_{f ,\id_a  } \bullet F(\underline{\Id_f})    ,\  \Id_{\id_{F(a)}} ].\end{equation*} The horn below verifies this commutativity.

  \begin{table}[H]
     \begin{center} 
   \scalebox{.96}{\begin{tabular}{ r | l || l | l | l | l | }
     \cline{2-6}

                    &   $\Delta_{-}(F(\chi( f,\id_a)))$      &  $\chi'(F(f),F(\id_a))$     &  $F(\chi(f,\id_a))$      &$\phi_{f,\id_a}$        & $\Id_{F(\id_{a})}$       \\ \cline{2-6}

           $\Lambda$         &                          &    $\chi'(F(f),F(\id_a))$         &  $\Id_{F(f)}$          &$\phi_{f,\id_a}\bullet F(\underline{\Id_f})$            &         $ \Id_{\id_{F(a)}}$     \\   \cline{2-6}  

                       &    $ F(\Delta_{-}(\Id_f))$      &  $F(\chi(f,\id_a))$                  &$F(\Id_f)$     & $F(\underline{\Id_f})$      &        $F(\Id_{\id_a})$         \\   \cline{2-6}  
        
                        &   $\Delta_{\bullet}( \phi_{f,\id_a},\ F(\underline{\Id_f}) ) $             &$\phi_{f,\id_a}$   &    $\phi_{f,\id_a}\bullet F(\underline{\Id_f})$    &     $F(\underline{\Id_f})$      &      $ \Id_{\id_{F(a)}}$     \\   \cline{2-6}  
  
                        & $s_0( \Id_{\id_{F(a)}})$       & $ \Id_{\id_{F(a)}}$                       & $ \Id_{\id_{F(a)}}$             &  $ \Id_{\id_{F(a)}}$     &   $ \Id_{\id_{F(a)}}$     \\   \cline{2-6} 
    \end{tabular}}\end{center}  \end{table} 
 
  \item For \textbf{BFun7}  let $f:a \ra b$ in $\Bic(X)$. We must show   $$  \lambda'_{F(f)} = (\upsilon_b \lhd F(f)) \bullet \phi_{ \id_b,f  } \bullet F(\lambda_f).$$ We have  $$\upsilon_b = \Id_{\id_{F(b)}}\ \ \ \ \ \ \ \ \ \ \ \ \ \ \lambda'_{F(f)} = \underline{(\wh{\Id_{F(f)}})}  \ \ \ \ \ \ \ \ \ \ \ \  \ F(\lambda_f)=F(\underline{(\wh{\Id_f})}). $$ By Lemma~\ref{bulletylemma} it suffices to show  $$   \wh{\Id_{F(f)}}= (\Id_{\id_{F(b)}} \tlhd F(f)) \bullet \phi_{\id_b ,f  } \bullet F(\underline{(\wh{\Id_f})}).$$ We have
\begin{multline*} \Delta_{\bullet}( \Id_{\id_{F(b)}} \tlhd F(f) ,\  \phi_{\id_b  ,f } \bullet F(\underline{(\wh{\Id_f})}) ) =  \\ [ \Id_{\id_{F(b)}} \tlhd F(f)  ,\ (\Id_{\id_{F(b)}} \tlhd F(f)) \bullet  \phi_{ \id_b ,f } \bullet F(\underline{(\wh{\Id_f})})        ,\  \phi_{ \id_b ,f } \bullet F(\underline{(\wh{\Id_f})})   ,\  \Id_{F(f)} ]\end{multline*} so by the Matching Lemma, it suffices to show \begin{equation} \label{BFun7eq} [ \Id_{\id_{F(b)}} \tlhd F(f)  ,\   \wh{\Id_{F(f)}}   ,\  \phi_{\id_b  ,f } \bullet F(\underline{(\wh{\Id_f})})   ,\  \Id_{F(f)} ].\end{equation} First note that
 \begin{align*} 
 d(\Delta_{\tlhd} (\Id_{\id_{F(b)}},F(f))) &=  [ \Id_{\id_{F(b)}}   ,\   \Id_{\id_{F(b)}}\tlhd F(f) ,\  \chi(\id_{F(b)},\ F(f))  ,\ \wh{ \Id_{F(f)}}] \\
 s_1( \chi(\id_{F(b)}, F(f)) )  &= [ Id_{\id_{F(b)}} ,\   \chi(\id_{F(b)},\ F(f))  ,\     \chi(\id_{F(b)},\ F(f))     ,\ \wh{ \Id_{F(f)}}]
 \end{align*} so by the Matching Lemma $$\Id_{\id_{F(b)}}\tlhd F(f)=  \chi(\id_{F(b)},\ F(f)).$$ Substituting into sphere~\ref{BFun7eq} we must show the commutativity of $$[ \chi(\id_{F(b)},\ F(f))  ,\   \wh{\Id_{F(f)}}   ,\  \phi_{  \id_b,f } \bullet F(\underline{(\wh{\Id_f})})   ,\  \Id_{F(f)} ].$$  The horn below verifies this commutativity.
  \begin{table}[H]
     \begin{center} 
   \scalebox{.90}{\begin{tabular}{ r | l || l | l | l | l | }
     \cline{2-6}

         &   $\Delta_{-}(\chi( \id_b,f))$ &  $\chi'(F(\id_b),F(f))$     &  $F(\chi(\id_b,f))$      &$\phi_{\id_b,f}$        & $\Id_{F(f)}$       \\ \cline{2-6}

   $\Lambda$         &       &  $\chi(\id_{F(b)},F(f))$         &  $\wh{\Id_{F(f)}}$          &$\phi_{\id_b,f}\bullet F(\underline{(\wh{\Id_f})})$            &         $ \Id_{F(f)}$     \\   \cline{2-6}  

          &    $ F(\Delta_{-}(\wh{\Id_{f}}))$      &  $F(\chi(\id_b,f))$                &$F( \wh{\Id_f} )=\wh{\Id_{F(f)}}$     & $ F(\underline{(\wh{\Id_f})})$      &         $F(\Id_f)=\Id_{F(f)}$        \\   \cline{2-6}  
        
   &   $\Delta_{\bullet}(\phi_{\id_b,f},\ F(\underline{(\wh{\Id_f})}) ) $  &$\phi_{\id_b,f}$  &  $\phi_{\id_b,f}\bullet F(\underline{(\wh{\Id_f})})$  &  $F(\underline{(\wh{\Id_f})})$   &   $ \Id_{\id_{F(a)}}$  \\   \cline{2-6}  
  
                        & $s_0( \Id_{F(f)})$       &$\Id_{F(f)}$        &$\Id_{F(f)}$          &  $\Id_{F(f)}$       &   $ \Id_{\id_{F(a)}}$     \\   \cline{2-6} 
    \end{tabular}}\end{center}  \end{table}
\end{itemize}\end{proof}
\begin{proposition} $\Bic$ preserves strictness in the sense that if $F:(X,\chi)\ra (Y,\chi')$ is a strict morphism of algebraic $2$-reduced inner-Kan simplicial sets, then $\Bic(F)$ is a strict functor of $(2,1)$-categories.
\end{proposition}
\begin{proof}If $F$ is strict, then $F(\chi(g,f))= \chi'(F(g),F(f))$ therefore for $\Bic(F)$ we have $$\phi_{g,f}=\underline{F(\chi(g,f))}=\underline{\chi'(F(g),F(f))}=\Id_{F(g)\circ F(f)} $$ by Lemma~\ref{lemma2}.
\end{proof}
\begin{proposition}
$\Bic$ is functorial, i.e. $\Bic(G F)=\Bic(G)\circ \Bic(F)$.
\end{proposition}
\begin{proof} It is clear that $\Bic(G F)$ and $\Bic(G)\circ \Bic(F)$ are the same on objects, $1$-morphisms and $2$-morphisms, and both are strictly identity preserving, having a trivial unitor. We must check that the distributors are the same. The distributor of $\Bic(G F)$ is $\underline{GF(\chi(g,f))}$, whereas the distributor of $\Bic(G)\circ \Bic(F)$ is $\underline{G(\chi(F(g),F(f)))} \bullet G(\underline{F(\chi(g,f))}).$
To show $$\underline{GF(\chi(g,f))}=\underline{G(\chi(F(g),F(f)))} \bullet G(\underline{F(\chi(g,f))})$$ by Lemma~\ref{bulletylemma} it suffices to show 
$$GF(\chi(g,f))=G(\chi(F(g),F(f))) \bullet G(\underline{F(\chi(g,f))}).$$ We have \begin{align*} &\Delta_\bullet( G(\chi(F(g),F(f))),\ G(\underline{F(\chi(g,f))})     ) \\ &\quad =[ G(\chi(F(g),F(f))),\  G(\chi(F(g),F(f)))\bullet G(\underline{F(\chi(g,f))})     ,\ G(\underline{F(\chi(g,f))}),\ \Id_{\id_{GF(a)}}   ].\end{align*}
By the Matching Lemma, it suffices to show the commutativity of: $$[ G(\chi(F(g),F(f))),\ GF(\chi(g,f))    ,\ G(\underline{F(\chi(g,f))}),\ \Id_{\id_{GF(a)}}   ].$$
This sphere is easily seen to be $d\left(G\left(\Delta_{-}(  F(\chi(g,f))    )\right)\right).$
\end{proof}
\subsection{The nerve of a functor of $(2,1)$-categories}
Besides being anticipated by Duskin in \cite{Dus02}, the construction defined in this section is described in \cite{BFB04} in the case of (weak) functors between strict bicategories, and a proof that the Duskin nerve is full and faithful in this case is also sketched there.
\begin{definition}
Let $\B$ and $\B'$ be a $(2,1)$-categories, and $F:\B \ra \B'$ be a strictly identity-preserving functor. 
The \emph{nerve} $N(F)$ of $F$ will be defined as a map of simplicial sets $N(F): N(\B)\ra N(\B')$. First we define $N(F)|^3_0: N(\B)|^3_0 \ra N(\B')|^3_0$ as a map of truncated simplicial sets as follows:
\begin{itemize} 
\item $N(F)$ is identical to $F$ on objects and $1$-morphisms. 
\item For a $2$-cell $(h, g, f \vbar \eta)$ in $N(\B)|^3_0$, where $\eta: g \Rightarrow h\circ f$ we define $$N(F)|^3_0(h, g, f \vbar ,\eta):=(F(h),F(g),F(f) \vbar \phi_{h,f} \bullet F(\eta)).$$
\item For a $3$-cell $[x_{123},\, x_{023},\, x_{013},\, x_{012}]$ in $N(\B)|^3_0$, which by definition is a $4$-tuple of $2$-cells meeting the commutativity condition expressed in Equation~\ref{3cellcondition}, we define
\begin{multline*}N(F)|^3_0 ( [x_{123},\, x_{023},\, x_{013},\, x_{012}] ) := \\ [N(F)|^3_0(x_{123}), \,  N(F)|^3_0(x_{023} ), \ N(F)|^3_0(x_{013}),\,  N(F)|^3_0(x_{012})].\end{multline*} In order for this map to be well-defined, we must show the $4$-tuple  $$[N(F)|^3_0(x_{123}), \,  N(F)|^3_0(x_{023} ), \ N(F)|^3_0(x_{013}),\,  N(F)|^3_0(x_{012})]$$ is indeed a $3$-cell of $N(\B')|^3_0$, i.e. that it meets the condition given in Equation~\ref{3cellcondition}. This is shown below in Proposition~\ref{F3cellcond}.
\end{itemize}
In order to check that $N(F)|^3_0$ is a morphism of truncated simplicial sets in $\cat{Set}_{\Delta|^3_0}$, we must check that it respects the face and degeneracy maps. $N(F)|^3_0$ can be seen to respect the face maps immediately from its definition, and the fact that it respects the degeneracy maps is shown in Proposition~\ref{NFrespectsdegeneracies} below. 

Then we can define $N(F):N(\B)\ra N(\B')$ by $$N(F)=\cosk^3(N(F)|^3_0):\cosk^3 (N(\B)|^3_0)=N(\B) \ra \cosk^3(N(\B')|^3_0)= N(\B').$$
\end{definition}
\begin{proposition} \label{NFrespectsdegeneracies} $N(F)|^3_0: N(\B)|^3_0 \ra N(\B')|^3_0$ respects each of the degeneracy maps. 
\end{proposition}
\begin{proof}
First let $a\in N(\B)_0$ be an object of $\B$. Then $$s_0 (N(F)|^3_0(a))=s_0 F(a) = \id_{F(a)} = F(\id_a)= N(F)|^3_0(\id_a)=N(F)|^3_0(s_0 a)$$ since $F$ is assumed to be strictly identity-preserving.

Now let $f \in N(\B)_1$ be a $1$-morphism of $\B.$ Then we have $$N(F)|^3_0 (s_0 f) = N(F)|^3_0 (\rho_f)=\phi_{f,\id_a} \bullet F(\rho_f).$$ Since we have assumed $F$ is strictly identity-preserving, we have $\upsilon_a=\Id_{F(a)}$. Therefore applying \textbf{BFun6} we have $$N(F)|^3_0 (s_0 f) =\phi_{f,\id_a} \bullet F(\rho_f) = (F(f)\rhd \upsilon_a)\bullet \phi_{f,\id_a} \bullet F(\rho_f)=\rho'_{F(f)}=s_0(F(f))\!=s_0(N(F)|^3_0 (f)).$$ By a similar argument using \textbf{BFun7} we have $N(F)|^3_0(s_1 f)= s_1(N(F)|^3_0 (f))$, so $N(F)|^3_0 $ respects the degeneracy maps in this dimension.

For a $2$-cell $\eta \in N(\B)_0$, the fact that $N(F)|^3_0(s_i \eta) = s_i N(F)|^3_0 (\eta)$ can be seen as an immediate consequence of the definitions of these degeneracy maps  and the fact that $N(F)|^3_0$ preserves face maps and the degeneracy maps for $1$-cells.
\end{proof}
\begin{proposition}\label{F3cellcond} Let $x_{0123}=[x_{123},\, x_{023},\, x_{013},\, x_{012}]$ be a $3$-cell of $N(\B)|^3_0$, i.e. $4$-tuple of $2$-morphisms of $N(\B)$ with appropriate sources and targets (for instance $x_{012}: x_{02}\ra x_{12}\circ x_{01} $) meeting the commutativity condition of Equation~\ref{3cellcondition}.
Then also $N(F)|^3_0 (x_{0123})$ meets the commutativity condition of Equation~\ref{3cellcondition} and is thus a $3$-cell of $N(\B')|^3_0.$\end{proposition}
\begin{proof}
The commutativity condition for $3$-cells (Equation~\ref{3cellcondition}) for $x_{0123}$ asserts:
\begin{equation}\label{x0123condition}\alpha_{x_{23},x_{12},x_{01}}\bullet(x_{23}\rhd x_{012})\bullet x_{023}=(x_{123}\lhd x_{01})\bullet x_{013}
\end{equation}
Now we unwrap our definition of $N(F)|^3_0 (x_{0123}).$
\begin{align*}N(F)|^3_0 (x_{0123}) &= [N(F)|^3_0(x_{123}), \,  N(F)|^3_0(x_{023} ), \ N(F)|^3_0(x_{013}),\,  N(F)|^3_0(x_{012})] \\ 
							&= [\phi_{x_{23},x_{12}} \bullet F(x_{123}),\ \phi_{x_{23},x_{02}} \bullet F(x_{023}),\ \phi_{x_{13},x_{01}} \bullet F(x_{013}),\ \phi_{x_{12},x_{01}} \bullet F(x_{012})]
\end{align*}
so the commutativity condition for $3$-cells for $N(F)|^3_0 (x_{0123})$ is 
\begin{multline} \label{NF3cellcond} \alpha'_{F(x_{23}),F(x_{12}),F(x_{01})}\bullet [F(x_{23})\rhd (\phi_{x_{12},x_{01}} \bullet F(x_{012}))] \bullet \phi_{x_{23},x_{02}} \bullet F(x_{023})= \\ [(\phi_{x_{23},x_{12}} \bullet F(x_{123}))\lhd F(x_{01})]\bullet \phi_{x_{13},x_{01}} \bullet F(x_{013}).
\end{multline}
We begin by manipulating the left hand side of Equation~\ref{NF3cellcond}
\begin{align}&\alpha'_{F(x_{23}),F(x_{12}),F(x_{01})}\bullet [F(x_{23})\rhd (\phi_{x_{12},x_{01}} \bullet F(x_{012}))] \bullet \phi_{x_{23},x_{02}} \bullet F(x_{023}) \nonumber \\
&\quad = \alpha'_{F(x_{23}),F(x_{12}),F(x_{01})}\bullet (F(x_{23})\rhd \phi_{x_{12},x_{01}}) \bullet (F(x_{23})\rhd  F(x_{012})) \bullet \phi_{x_{23},x_{02}} \bullet F(x_{023}) \nonumber \\
&\quad = \alpha'_{F(x_{23}),F(x_{12}),F(x_{01})}\bullet (F(x_{23})\rhd \phi_{x_{12},x_{01}}) \bullet \phi_{x_{23},x_{12}\circ x_{01} } \bullet F(x_{23}\rhd  x_{012}) \bullet F(x_{023}) \nonumber \\
&\quad = \label{NF3cellcond2} (\phi_{x_{23},x_{12}} \lhd F(x_{01})) \bullet\phi_{x_{23}\circ x_{12},x_{01}}    \bullet  F(\alpha_{x_{23},x_{12},x_{01}}) \bullet F(x_{23}\rhd  x_{012}) \bullet F(x_{023}) 
\end{align} 
using \textbf{B5}, \textbf{BFun3}, and \textbf{BFun5}. Now apply $F$ to both sides of Equation~\ref{x0123condition} and distribute $F$ over $\bullet$ using \textbf{BFun2} to yield
$$F(\alpha_{x_{23},x_{12},x_{01}})\bullet F(x_{23}\rhd x_{012})\bullet F(x_{023})=F(x_{123}\lhd x_{01})\bullet F(x_{013}).$$
Substitute this equation into Equation~\ref{NF3cellcond2} to get:
\begin{align*}&(\phi_{x_{23},x_{12}} \lhd F(x_{01})) \bullet \phi_{x_{23}\circ x_{12},x_{01}}    \bullet  F(x_{123}\lhd x_{01})\bullet F(x_{013}) \\ 
& \quad = (\phi_{x_{23},x_{12}} \lhd F(x_{01})) \bullet   (F(x_{123})\lhd F(x_{01}))\bullet \phi_{x_{13},x_{01}}   \bullet F(x_{013})\\
& \quad = [(\phi_{x_{23},x_{12}} \bullet F(x_{123}))\lhd F(x_{01})]\bullet \phi_{x_{13},x_{01}} \bullet F(x_{013})
\end{align*} using \textbf{B6} and \textbf{BFun4}. This shows Equation~\ref{NF3cellcond} holds.
\end{proof}
\begin{proposition} If $F:\B\ra\C $ is a strict functor, then $N(F):N(\B)\ra N(\C)$ is strict.
\end{proposition}
\begin{proof}If $\phi_{g,f}=\Id_{F(g)\circ F(f)}$ then $N(F)$ sends $(g,g\circ f, f \vbar \Id_{g\circ f})$ to $$(F(g), F(g)\circ F(f),F(f) \vbar F(\Id_{g\circ f})  )=\Id_{F(g)\circ F(f)}$$ thus $N(F)$ preserves the natural algebraic structures of $N(\B)$ and $N(\C)$. 
\end{proof}
\begin{proposition}
$N$ is functorial, i.e. $N(G)\circ N(F)=N(GF)$.
\end{proposition}
\begin{proof}It is clear that $N(G)\circ N(F)=N(GF)$ for $0$ and $1$-cells. For a $2$-cell $((h,g,f),\eta)$ we have 
\begin{align*}
N(GF)(h,g,f \vbar \eta)&=(GF(h),GF(g),GF(f) \vbar (\phi'_{F(h), F(f)}\bullet G(\phi_{h,f})) \bullet GF(\eta))\\
N(G)\circ N(F)(h,g,f\vbar \eta)&=(GF(h),GF(g),GF(f) \vbar \phi'_{F(h), F(f)}\bullet G(\phi_{h,f}\bullet F(\eta)))
\end{align*}
The equality follows from the functoriality of $G$ with respect to $\bullet$.
\end{proof}
\section{The natural isomorphism $\Bic(N(\B))\cong \B$ and summary \label{summarysection}} 

In Section~\ref{bicisosec} we defined an isomorphism $u:N(\Bic(X))\ra X.$ Now that we have defined $N$ and $\Bic$ as functors, we show that $u$ is a natural transformation. Unwrapping the definitions, this works out to showing that if  $x$ is a $2$-cell in $X$ with $dx = [ h, g, f],$  and $F:(X,\chi)\ra (Y,\chi')$ is a map of algebraic $2$-reduced inner-Kan simplicial sets  then $$\underline{F(x)} = \underline{F(\chi(h,f))} \bullet F( \underline{x})$$ we have $$\Delta_\bullet(\underline{F(\chi(h,f))}  , F( \underline{x})) = [ \underline{F(\chi(h,f))} ,\    \underline{F(\chi(h,f))}\bullet F( \underline{x})  ,\  F( \underline{x})    ,\  \Id_{\id_{F(a)}  }] $$ so by the Matching Lemma it suffices to show the commutativity of the sphere: $$\Delta_\bullet(\underline{F(\chi(h,f))}  , F( \underline{x})) = [ \underline{F(\chi(h,f))} ,\   \underline{F(x)}  ,\  F( \underline{x})    ,\  \Id_{\id_{F(a)}  }]. $$ The following horn shows this commutativity:
\begin{table}[H]
     \begin{center} 
   \scalebox{.91}{\begin{tabular}{ r | l || l | l | l | l | }
     \cline{2-6}

         &   $\Delta_{-}(F(\chi(h,f)))$ &  $\chi'(F(h),F(f))$     &  $F(\chi(h,f))$      &$\underline{ F(\chi(h,f))}$        & $\Id_{F(f)}$       \\ \cline{2-6}

        &  $\Delta_{-} (F(x))$     &  $\chi'(F(h),F(f))$       &  $F(x)$          &$\underline{F(x)} $            &         $ \Id_{F(f)}$     \\   \cline{2-6}  

          &    $ F(\Delta_{-}(x))$      &   $F(\chi(h,f))$    & $F(x)$      & $ F(\underline{x})$      &         $F(\Id_f)=\Id_{F(f)}$        \\   \cline{2-6}  
        
   $\Lambda$   	 &     &$\underline{F(\chi(h,f))}$  &  $\underline{F(x)}$  &  $ F(\underline{x})$   &   $ \Id_{\id_{F(a)}}$  \\   \cline{2-6}  
  
                        & $s_0( \Id_{F(f)})$       &$\Id_{F(f)}$        &$\Id_{F(f)}$          &  $\Id_{F(f)}$       &   $ \Id_{\id_{F(a)}}$     \\   \cline{2-6} 
    \end{tabular}}\end{center}  \end{table}

In the opposite direction there is a strict isomorphism of $(2,1)$-categories $\Bic ( N ( \B),\chi) \cong \B$, where $\chi$ is the natural algebraic structure on $N(\B)$, as given in Definition~\ref{nervedef}. First we use the definitions to describe data of $\Bic ( N ( \B))$ explicitly: 
\begin{itemize}
\item The objects and morphisms of $\Bic(N(\B))$ are identical to those in $\B.$
\item A $2$-morphism $f \Rightarrow g:a \ra b$ in $\Bic(N(\B))$ is given by a $2$ morphism $ f \Rightarrow g \circ \id_a$ in $\Bic(N(\B)).$
\end{itemize}
Following \cite{Dus02}, we construct a functor of $(2,1)$-category $U: \B \Rightarrow \Bic(N(\B))$ as follows:
\begin{itemize}
\item $U$ is the identity on objects and morphisms.
\item For a $2$-morphism $\eta : f \Rightarrow g$ we define $U(\eta)=\rho'_g \bullet \eta.$
\item $U$ is strict, i.e., the unitor of $U$ consists of identity $2$-morphisms.
\end{itemize}
Since $U$ is strict, the functor axioms of $U$ assert that $U$ strictly preserves each of the bicategory structures of $\B$. The proof of this fact involves checking  for each piece of structure of $\B$ that a certain identity in $\B$ follows from the bicategory axioms. See Lemma 8.4 and Theorem 8.5 in \cite{Dus02} for details.

The naturality of $U$ is also easily checked. Suppose $F:\B \ra \C$ is a functor of $(2,1)$-categories. The functor $U_{\C} F:\B \ra \Bic(N(X))$ is the identity on objects and morphisms, and has a trivial unitor. For a $2$-morphism $\eta$ in $\B$ we have $$U_{\C} F (\eta) = \rho'_{F(g)} \bullet F(\eta)$$ (where $\rho'$ and $\bullet$ indicate the structures of $\C$, not of  $\Bic ( N(X))$). The distributor of $U_{\C} F$ is given by $U_{\C}(\phi_{g,f})= \rho'_{F(g)\circ F(f)} \bullet \phi_{g,f}. $

The functor $\Bic ( N(F)) U_{\B}$ again is the identity on objects and morphisms and has a trivial unitor. For a $2$-morphism $\eta:f\Rightarrow g : a\ra b$ in $\B$ we compute $$\Bic ( N(F)) U_{\B} (\eta) = \Bic ( N(F))(\rho_g \bullet \eta) = \phi_{g,\id_a} \bullet F(\rho_g) \bullet F(\eta).$$ From \textbf{BFun6} for $F$ we have $$  \rho'_{F(g)} = (F(g)\rhd \upsilon_a) \bullet \phi_{ g,\id_a  } \bullet F(\rho_g)  $$ and since $\upsilon_a=\Id_{F(a)}$ we have $\rho'_{F(g)} =  \phi_{ g ,\id_a} \bullet F(\rho_g).$ This ensures $$U_{\C} F (\eta)= \rho'_{F(g)} \bullet F(\eta)    =  \phi_{g,\id_a} \bullet F(\rho_g) \bullet F(\eta)= \Bic ( N(F)) U_{\B} (\eta) .$$

To compute the distributor of $\Bic ( N(F)) U_{\B},$ we must first compute the distributor of $\Bic ( N(F)),$ which for $f,g$ is $\underline{N(F)(\chi(g,f))}.$ By the definition of the natural algebraic structure $\chi$ for $N(\B)$, we have $\chi(g,f)=(g,g\circ f,f \vbar \Id_{g\circ f})$ so $$\underline{N(F)(\chi(g,f))}=\underline{N(F)(g,g\circ f, f \vbar \Id_{g\circ f})}=\underline{\phi_{g,f}}.$$ Applying the definition of $\Delta_{-}(\phi_{g,f})$ and the $3$-cell condition in $N(\B)$, Equation~\ref{3cellcondition}, we get 
$$\alpha'_{F(g),F(f),\id}\bullet F(g) \rhd \rho'_{F(f)} \bullet \phi_{g,f}=(\Id_{F(g)\circ F(f)}\lhd \id)\bullet \underline{\phi_{g,f}}.$$ Applying the compatibility of the associator with $\rho'$, axiom $\textbf{B13}$, we get that the distributor of  $\Bic ( N(F)) U_{\B}$ is $\underline{\phi_{g,f}}=\rho'_{F(g)\circ F(f)} \bullet \phi_{g,f}.$ This matches the distributor of $U_{\C} F$, so we have shown that $U$ is natural.

We have now shown:
\begin{theorem} 
\label{bicsummary} The functors $N$, $\Bic$ are inverse equivalences of categories between the category of $2$-reduced inner-Kan algebraic simplicial sets and the category of small $(2,1)$-categories and strictly identity preserving functors. Furthermore, $N$ and $\Bic$ preserve strictness, and the natural isomorphisms $u:N \Bic\cong \Id $ and $U: \Bic N \cong \Id$ exhibiting the equivalence are strict, thus $N$ and $\Bic$ are also inverse equivalences of categories between the category of $2$-reduced inner-Kan algebraic simplicial sets and strict morphisms and the category of small $(2,1)$-categories and strict functors.
\end{theorem}

\begin{remark}
The functor that forgets algebraic structure from an inner-Kan simplicial set is clearly fully faithful and essentially surjective, thus $N$ gives an equivalence of categories between the category of small $(2,1)$-categories and the category of (non-algebraic) $2$-reduced inner-Kan simplicial sets. The algebraic structure is necessary only to construct an inverse to this functor. In particular, if $\chi$ and $\chi'$ are algebraic structures on $X$, then $\Bic(X,\chi)$ and $\Bic(X,\chi')$ are (weakly) isomorphic as $(2,1)$-categories. 
\end{remark}

\begin{remark}
Theorem~\ref{bicsummary} could also be obtained with minimal effort as a corollary of Theorem~3.17 in \cite{Gur09}
\end{remark}



\chapter{The bisimplicial nerve of a Verity double category \label{vdcchapter}}




\section{Dimensional categories \label{genglennsec}}
In this section, we give a general setup in which we can define the concepts of spheres, horns, the coskeleton and skeleton functors, and Glenn tables. For $n$-simplicial sets, the generalizations are quite straightforward, but we will also use these concepts in subsequent chapters, so we work in a more general setting.

\begin{definition} A \emph{dimensional category} is a small category $C$ equipped with a non-negative grading $\dim:\mbox{Ob}(C) \ra \mathbb{N},$ called the \emph{dimension}, such that all isomorphisms preserve $\dim.$\end{definition}
We can define the coskelton and skeleton functors and the concept of a ``sphere'' for the presheaf category of any dimensional category.
\subsection{Coskeleta of dimensional categories\label{bicoskeletasec}}
Let $C$ be a dimensional category. We denote the category of \emph{$C$-sets}  (presheaves of sets on $C$) by $\cat{Set}_{C}$. The Yoneda embedding $\Hom(-,c)$ of an object $c \in \C$ in $\cat{Set}_{C}$ will be denoted by $C[c]$, which we sometimes refer to as the \emph{$c$-simplex}. If $X$ is a $C$-set, we will call the elements of $X(c)$ (also denoted $X_c$) the \emph{$c$-cells} of $X$.
\begin{definition} Let $C|^n_0$ denote the full subcategory of $C$ of objects of dimension $\leq n.$  We have a \emph{truncation} functor, given by restriction, that is, the presheaf pullback along the inclusion $C|^n_0\ra C$, which is a map  $$\mathbf{tr}^n:\cat{Set}_{C}\ra \cat{Set}_{C |^n_0}.$$
\end{definition}
\begin{definition}\label{coskeletaldef} The functor $\mathbf{tr}^n$ has adjoints on both the left and the right. The right adjoint is the presheaf extension associated with $\mathbf{tr}^n$ which we denote $\mathbf{cosk}^n,$ computed by the formula 
$$(\mathbf{cosk}^n  Y)(c) = \Hom( \mathbf{tr}^n c, Y).$$
For a $C $-set $X$, we define $\mathbf{Cosk}^nX:=\mathbf{cosk}^n\mathbf{tr}^n X.$ A $C $-set $X$ is called $n$-coskeletal if the natural unit map $X\ra \mathbf{Cosk}^n X$ is an isomorphism. If this map is a monomorphism, we will call $X$ $n$-subcoskeletal. \end{definition}

Similarly, the left adjoint to $\tr^n$ is denoted $\sk^n$. We can construct this adjoint as follows:
\begin{definition}For a $C $-set $X$, let $\Sk^n X$ denote the subcomplex of $X$ made up of cells $X(c)$ whose classifying map $C [c] \ra X$ factor through a cell of dimension $\leq n$. For $Y \in \cat{Set}_{C |^n_0} $ let $\sk^nY $ be a colimit-preserving functor such that $\sk^n(C [c])\cong \Sk^n(c)$, i.e., $$\sk^n Y :=\int^{c \in C |^n_0}\Sk^n(c)  \times Y(c).$$
\end{definition}
\begin{proposition}\label{skelcomp} $\sk^n$ is left adjoint to $\tr^n$. Furthermore, we have a natural isomorphism $\Sk^n \cong \sk^n\tr^n$ with the counit $\sk^n \tr^n X \ra X$ corresponding to the inclusion $\Sk^n X \ra X.$
\end{proposition}
\begin{proof} We will give the unit and counit of the adjunction $\sk^n : \cat{Set}_{C |^n_0} \leftrightarrows \cat{Set}_{C }: \tr^n$. For the unit $\sk^n\tr^n \ra \Id$ first note that $\sk^n$ preserves colimits by its definition, and $\tr^n$ preserves colimits by the fact that colimits are preserved pointwise in a presheaf category. So $\sk^n \tr^n$ and $\Id$ both preserve colimits, so it suffices to construct the counit on representable presheaves. For $c \in C $ we have $$\sk^n \tr^n (C [c]) \cong \sk^n(C [c])\cong \Sk^n(C [c])$$ so we define the counit using the natural inclusion $\Sk^n(C [c])\ra C [c]$. 

For the unit, let $c \in C |^n_0$ we have $$\tr^n \sk^n(C [c])\cong \tr^n \Sk^n(C [c])\cong C [c]$$ with the last equality following since $\dim(c) \leq n.$ Again using the fact $\tr^n\sk^n$ preserves colimits, this isomorphism can be used to define a natural unit isomorphism $\Id \ra \tr^n\sk^n$.

The zig-zag identities for the unit and counit can likewise be checked on representable presheaves, where they follow immediately from the definitions.

The isomorphism $\sk^n\tr^n \cong \Sk^n$ follows since both are colimit-preserving which agree on representable presheaves. The colimit agrees with the natural inclusion by definition.
\end{proof}
It follows directly that $\Sk$ is left adjoint to $\Cosk$, with the natural map $X\ra \Cosk^n X$ adjoint to the map $\Sk^n X \ra X$, since both are adjoint to $\id: \tr^n X \ra \tr^n X. $
\begin{lemma} \label{trunklemma} Let $m \geq n$ The natural map $\tr^n\Sk^mX \ra \tr^n X $ and $\tr^n X \ra \tr^n\Cosk^m X $ are isomorphisms. In other words $\Cosk^mX$ matches $X$ up to cells of dimension $m.$ \end{lemma}
\begin{proof} This follows for $\Sk^m$ directly from its definition. For $\Cosk, $
restricting to $\dim(c) \leq n$, the unit map $X_{c} \ra (\mathbf{Cosk}^mX)_{c}$ is given by 
\begin{align*} \tr^m: X_{c}\cong \Hom(C [c],X) \ra \Hom(\tr^m C [c], \tr^m(X)) & \cong \Hom(C [c], \tr^m X) \cong (\Cosk^mX)_{c} \\ &\cong (\tr^m X)_{c} \end{align*} (using the Yoneda lemma), which by the naturality of the Yoneda embedding is the canonical map $X_{c}\ra (\tr^m X)_{c}$, which is an isomorphism since $\dim(c) \leq n \leq m$.  
\end{proof}
\begin{lemma} \label{iterlemma}Let $m\geq n$. Then the natural (given by the counit  $\Sk^m\ra \Id$) map $$\Sk^m \Sk^nX \ra \Sk^mX$$ is an isomorphism. Similarly, the natural (given by the unit $\Id \ra \Cosk^m$) map $$\Cosk^nX \ra\Cosk^m\Cosk^nX$$ is an isomorphism.
\end{lemma}
\begin{proof} The statement for $\Sk^m$ is immediate from the definition. For the statement for $\Cosk$ it is then clear that we have a natural isomorphism $$\Cosk^n \cong \Cosk^m\Cosk^n$$ since they are respectively right adjoint to $\Sk^m$ and $\Sk^n\Sk^m$. It is not immediately obvious that this isomorphism is the same map as given by the unit $\Id \ra \Cosk^m$, but this follows formally from our definitions, and we leave the verification to the reader (for instance by a ``string diagram chase''). 
\end{proof}

Using Lemma~\ref{iterlemma} and the unit $\Id \ra \Cosk^n$ we have a map $\Cosk^m \ra \Cosk^m \Cosk^n\cong \Cosk^n$, and it's easy to see this map is compatible with the unit maps, that is, the following diagram commutes:
\begin{center}
\begin{tikzpicture}[scale=1.5,auto]
\node (10) at (1.5,1) {$\Cosk^mX$};
\node (00) at (0,1) {$X$};
\node (11) at (1.5,0) {$\Cosk^nX$};

\draw[->] (00)--(10);
\draw[->] (00)--(11);
\draw[->] (10)--(11);
\end{tikzpicture}
\end{center}
By Lemma~\ref{trunklemma} $X \ra \Cosk^n X$ is an isomorphism up to cells of dimension $n$, so it follows that $X$ is the (indirect) limit of the diagram:
\begin{equation}\label{coskdiagram} \cdots \ra  \Cosk^{n+2} X\ra  \Cosk^{n+1} X \ra \Cosk^n X. \end{equation}

\begin{definition}\label{spheredef} For an object $c$ of dimension $n$, we define the \emph{universal sphere} $d C [c]:=\Sk^{n-1}(C[c]).$  A $dC[c]$-\emph{sphere} in $X$ is a map $d C[c]\ra X.$ A filler of a sphere is an extension of the sphere to a map $\C[c]\ra X$ along the natural (counit) map $d C[c] \ra C[c]$. \end{definition}
\begin{proposition} \label{coskeletallemma} A $C $-set $X$ is $n$-coskeletal if and only if every $dC[c]$-sphere in $X$ has a unique filler when $\dim(c)>n.$  Similarly $X$ is $n$-subcoskeletal if and only if every $dC[c]$-sphere in $X$ has at most one filler when $\dim(c)>n.$
\end{proposition}
\begin{proof} Note that the filler conditions mentioned above equivalently state that the natural map $$\Hom(  C [c], X) \ra \Hom(d C [c],X)$$ induced by the inclusion $d C [c]\ra  C [c]$ is an isomorphism, or respectively, an injective map. First assume $X$ is sub-coskeletal, i.e. $X \ra  \Cosk^n(X)$ is a monomorphism. Then we have a commutative diagram of natural maps:

\begin{center}
\begin{tikzcd}
\Hom(  C [c], X) \arrow{r}{\subseteq} \arrow{d}{} & \Hom( C [c],   \Cosk^n X)  \arrow{r}{\cong} \arrow{d}{} &    \Hom(\Sk^n  ( C [c]), X)       \arrow{d}{\cong}            \\
\Hom(d C [c],X)        \arrow{r}{\subseteq}     & \Hom( d C [c], \Cosk^n X)  \arrow{r}{\cong}             &    \Hom(\Sk^n(d C [c]), \Cosk^n X) 
\end{tikzcd}
\end{center}

The map on the far right $$\Hom(\Sk^n  ( C [c]), X) \ra \Hom(\Sk^n(d  C [c]), \Cosk^n X) $$ is induced by the unit map $$ \Sk^n \Sk^{\dim(c)-1}( C [c])  \ra  \Sk^n( C [c]) $$ which is an isomorphism by Lemma~\ref{iterlemma}. The commutativity of this diagram then guarantees the natural map $\Hom(  C [c], X) \ra \Hom(d C [c],X)$ is injective. 

Alternatively if $X$ is coskeletal and thus $X \ra \Cosk^n(X)$ is an isomorphism, the horizontal maps in the above diagram are isomorphisms, so by the commutativity of the diagram we conclude the map $$\Hom(  C [c], X) \ra \Hom(d C [c],X)$$ is bijective.

For the other direction suppose every sphere in $X$ of dimension $\geq n$ has a unique filler (has at most one filler).  Let $m \geq n$ and consider the unit map $\tr^{m+1}X \ra \tr^{m+1}\Cosk^m(X)$. This map is an isomorphism on cells of dimension $\leq m$ by Lemma~\ref{trunklemma}, and for cells of dimension $m+1$, the  map is given by $\Hom( C [c],X)\ra \Hom(\Sk^m{ C [c]}= d C [c],  X)$ which is bijective (injective) by the hypothesis. Thus the map $\tr^{m+1}X \ra \tr^{m+1}\Cosk^n(X)$ is an isomorphism, and applying $\cosk^{m+1}$ we deduce that the natural maps  $$\Cosk^{m+1}(X) \ra \Cosk^{m+1}\Cosk^m(X)\cong \Cosk^m(X)$$ in Diagram~\ref{coskdiagram} are all isomorphisms (injective). Thus the map $X \ra \Cosk^n X$ is an isomorphism (is injective), since it is a component map of the limit.
\end{proof}
\subsection{Good and excellent dimensional categories}
\begin{definition} A dimensional category will be called \emph{simple} if it has no nonidentity isomorphisms.\end{definition}
\begin{definition}\label{cofacesystemdef} For a dimensional category $C $, a \emph{coface} map in $C $ is a monic map which increases dimension by $1$  We say $C $ is a \emph{good} if the following holds:
\begin{itemize}
\item For every morphism $f:a\ra c$ in $C $ with $\dim a < \dim c$, there is a coface map $d'$ and a morphism $f'$ such that $f=d'\circ f'$ 
\end{itemize}
If furthermore the following property holds, then we say $C $ is \emph{excellent}:
\begin{itemize}[resume]
\item Let $f:a\ra c$ be a morphism in $C $ with $\dim a< \dim c$. Suppose $f$ factors through coface maps in two ways $f=d'\circ f'=d''\circ f''$, where $d'$ and $d''$ are not related by precomposition with an isomorphism. Then there are coface maps $e'$ and $e''$ and a morphism $g$ such that $d'e'=d''e''$ and such that $f$ factors as $f=d'e'g=d''e''g.$
\end{itemize}
\end{definition}
\begin{example}The category $\Delta^k$, can be given the structure of an excellent simple dimensional category, with $$\dim([\n^1,\n^2, \ldots, \n^k])\ =\sum_{1\leq i\leq k} \n^i.$$ Segal's category $\Gamma$ (defined in Section~\ref{gammasection}) with $\dim(\underline{n})=n$ is also an excellent dimensional category.
\end{example}
Recall for any presheaf category $C $, we have a coend formula expressing any presheaf $X$ on $C $ as a colimit of representable presheaves:
\begin{equation*}X\cong\int^{c\in C } X(c)\times  C [c]\cong\mbox{colim}_{x\in X(c)}  C [c] \end{equation*}
Where the colimit is taken over the \emph{category of objects of $X$}, whose objects are elements of any $X(c)$, with a morphism from $x \in X(c)$ to $x' \in X(c')$ for each $f: c \ra c'$ with $X(f)(x')= x.$

\subsection{Glenn categories and tables for simple dimensional categories.}
Let $C $ be a dimensional category.
\begin{definition}\label{mincompdimdef} Let $c$ be an object of $C $, and let $S$ be a subset of the set of coface maps with target $c$.
We will say a set of coface maps $S$ \emph{respects isomorphisms}  if whenever $ s=s'\varphi $ where $\varphi$ is an isomorphism, then $s\in S$ if and only if $s' \in S$.

If $S$ is a set of coface maps with target $c$ which respects isomorphisms, then the \emph{ universal $S$-horn} $\Lambda_S^c$ of $c$ is the subsheaf of $ C [c]$ consisting of maps $b\ra c$ which factor through some face map which is not in $S$. The \emph{dimension} of the horn $\Lambda_S^c$ is $\dim{c}$, while the \emph{minimal complementary dimension} of $\Lambda_S^c$ is the smallest $n$ for which the natural map $\tr^n \Lambda_S^c\ra  C [c]$ is not an isomorphism. A horn is called \emph{nice} if its minimal complementary dimension is equal to one less than its dimension.

If $X$ is a $C$-set, a $\Lambda^c_S$-horn in $X$ (also referred to as a $c$-horn in $X$) is a map $\Lambda^c_S \ra X$. A \emph{filler} of such a horn in $X$ is an extension of this map along the natural inclusion $\Lambda^c_S \ra C[c]$.
\end{definition}
\begin{definition}
A \emph{special coface system} for $C $ to consists of a choice of representative for each equivalence class of coface maps with target $c$ under the equivalence relation of differing by precomposition with an isomorphism. Our chosen representatives are called \emph{special coface maps}. Note that if $C $ is \emph{simple}, i.e. has no non-identity isomorphisms, then this choice is trivial.
\end{definition}

For a universal horn $\Lambda_S^c$ consider the subcategory $G(c,S)$ of the category of objects $O(c,S)$ of $\Lambda_S^c$ having as objects those morphisms $a \ra c$ which are either special coface maps, or can be written as compositions of two special coface maps. The morphisms are given by precompositions with coface maps and isomorphisms. We call this subcategory the \emph{ Glenn category of $\Lambda_S^c$}.
\begin{proposition}
If $C $ is a good dimensional category, then $\Lambda_{\emptyset} (c)=d C [c].$
\end{proposition} 
\begin{proof} The definition of a good dimensional categories ensures these are the same subsheaf of $ C [c]$.
\end{proof}

%
%
%

We seek to explicate a precise sense in which we can say that $G(c,S)$ contains enough information to compute colimits over $O(c,S)$. We recall the following standard notion from category theory:
\begin{definition}
Suppose we have a functor $F:C \ra \mathcal{D}$ and an object $d\in \mathcal{D}$. The \emph{comma category} $(d/ F)$ is the category whose objects are arrows $f:F(c) \ra d$ and whose morphisms from $f:F(c)\ra d$ to $g:F(c')\ra d$ are morphisms $c \ra c'$ in $C$ which make the obvious triangle commutative. $F$ is called \emph{cofinal} if for every $d \in \mathcal{D}$ the category $(d/F)$ is non-empty and connected.
\end{definition}
\begin{lemma}\label{cofinallemma}
If $F$ is cofinal if and only if for any category $\mathcal{E}$ and any functor $G: \mathcal{D} \ra \mathcal{E}$ the natural map $$\colim(G\circ F)\ra  \colim(G)$$ is an isomorphism.
\end{lemma}

\begin{definition}
Let $C $ be a good dimensional category with a special coface system. We say the universal horn $\Lambda_S^c$ is \emph{tabular} if the natural inclusion $I:G(c,S)\ra O(c,S)$ is cofinal, ensuring by Lemma~\ref{cofinallemma} that the colimit for $\Lambda_S^c$ can be computed by \begin{equation}\Lambda_S^c\cong \mbox{colim}_{g\in G(c,S)} \ \underline{\mbox{source}(g)}. \label{coendeq}\end{equation}
\end{definition}
\begin{proposition}
If $C $ is an excellent dimensional category, then for any special coface system on $C $, every universal horn $\Lambda_S^c$ is tabular.
\end{proposition}
\begin{proof}
Let $f:d \ra c$ be an object of $O(c,S)$. The category $(f/I)$ has the following description:
\begin{itemize}
\item An object of $(f/I)$ is a factorization $f=h\circ g$ where $h$ is either a special coface map not in $S$ or a composition of two special coface maps not in $S$. 
\item A morphism $h\circ g \ra h'\circ g'$ is a factorization $h'= h \circ h''$ where $h''$  is a special coface map or an isomorphism. Note that since $h$ is monic we have that $g=h''\circ g'$ holds automatically.
\end{itemize}
Note that since the coface maps of $C $ are monic, the factorization $f = h \circ g$ is uniquely determined by $h$. Thus we will think of the objects of $(f/I)$ as being coface maps not in $S$ or composites of two coface maps not in $S$ through which $f$ factors, without mention of $g$.  There is a morphism from $h$ to every composite $h \circ h''$ through which $f$ also factors. For brevity, we call a map through which $f$ factors \emph{$f$-valid}.

The fact that $(f/I)$ is non-empty is almost tautological, by definition $f$ factors through some coface map $h'$ not in $S$, so $f$ factors through the special coface map which has the form $h=h'\varphi$ where $\varphi$ is an isomorphism.

To show $(f/I)$ is connected, let $\smileeq$ be the equivalence relation on objects of $(f/I)$ generated by the morphisms. We must show $\smileeq$ is the complete relation. Clearly if $h$ is a composition of two special coface maps $h=h_0\circ h_1$, then $h\smileeq h_0$. So it suffices to show $h \smileeq h'$  in the case where $h$ and $h'$ are special $f$-valid coface maps. Since they are special, if $h\neq h'$ then they are not related by precomposition with an isomorphism so we apply the definition of  s excellent dimensional category to get coface maps $e,e'$ such that $h \circ e = h' \circ e'$ and this is a $f$-valid map. Then take the special coface maps representing $e$ and $e'$, which we denote $\hat{e}=e\circ \varphi$ and $\hat{e'}= e' \circ \varphi'$ such that $\varphi,\varphi'$ are isomorphisms. Then clearly $h \circ \hat{e}$ and $h' \circ \hat{e'}$ are $f$-valid and related by precomposition with $\varphi^{-1} \circ \varphi'.$ Thus $h \smileeq h \circ \hat{e} \smileeq h \circ \hat{e'} \smileeq h'.$ 
\end{proof}

Suppose $C $ has a locally finite set of coface maps, meaning there are finitely many coface maps with any given target. Suppose further we have chosen for every object $b$ of $C $ an ordering for every set of coface maps with target $b$. In this case, we say we have an \emph{ordered special coface system}.

If we have an ordered special coface system for a good dimensional category $C $ can give a graphical notation explicitly describing the Glenn category of a tabular horn $\Lambda_S^c$. First we consider the case where $S=\emptyset$ so that $\Lambda_{\emptyset} (c)=d C [c].$ We make a table, each row representing a special coface map $f$ with target $c$. In the row for $f$ we list $f$ itself, followed by the maps which are obtained by precomposing $f$ with each possible special coface map, in order. This table is called the \emph{universal Glenn table} for $c$. Note that it need not be a rectangular array. Table~\ref{univtableex} in the next section gives an example of a universal Glenn table for the category $\Delta^2.$
 
A map $F: d C [c]\ra X$ is equivalent to making a cone to $X$ from the canonical functor $g\ra\underline{\mbox{source}(g)}$ along which the colimit $$\Lambda_S^c\cong \mbox{colim}_{g\in G(c,S)} \ \underline{\mbox{source}(g)}$$ from Equation~\ref{coendeq} giving $\Lambda_{\emptyset} (c)$ is computed. For each object $g$ in the  Glenn category of $ C [c]$, we must give a map from the representable presheaf $\underline{\mbox{source}(g)}$ to $X$, or equivalently, an object of $X(\mbox{source}(g))$. To describe this cone, we make a new table from the universal  Glenn table in which we replace $g$ with this cell $X(\mbox{source}(g)).$ We call this table a Glenn table for $c$ in $X$.

To ensure our Glenn table gives a unique component of the cone for each object of $G(c,S)$, if two entries of the universal table ``match'', i.e., give the same object in the  Glenn category, then we must have the same entries in the corresponding places of the Glenn table.
The property that these components have the ``cone property'' with respect to maps in $G(c,S)$ given by precomposition with coface maps, is equivalent to the statement that, if the entry for a special coface $f$ of $c$ at the beginning of a row is replaced with the cell $y$, and $z$ is an entry in this row spot corresponding to the (iterated) special coface $g$, then $X(g)(y)=z$. That is, the entries of a row beginning with the cell $y$ constitute an ordered list of the \emph{special faces} of $y$. If $C $ is simple, precomposition with a coface map is the only kind of non-identity morphism in the  Glenn category, so this is the only property we must check. Otherwise, we must consider the cone property with respect to isomorphisms in $G(c,S)$. For the Glenn table, this corresponds to the condition that if we have $f = f'\varphi $ for two entries $f, f'$ in the universal table and an isomorphism $\varphi$, then if $A$ and $A'$ are the corresponding entry in the Glenn table, then $ A = X(\varphi) A'.$ 

We give an example which depends on definitions from Chapter~\ref{gammachapter}, the universal Glenn  table for $\underline{4}\in \Gamma$:
\begin{table}[H]  \begin{center}\caption{ The universal Glenn  table for $\underline{4}$}

    \begin{tabular}{ r | c || c | c | c | c z c | c | }
     \cline{2-8}

     &$123$   & $23$  &    $(12)3$     &$1(23)$           &$12$   & $13$ & $(13)2$  \\ \cline{2-8}
     
    &  $(01)23$ &$23$   &  $(012)3$   & $(01)(23) $   &     $(01)2$    &$(01)3$  & $(013)2$   \\ \cline{2-8}  

    & $0(12)3$     &   $(12)3$  &     $(012)3$      &    $0(123)$  & $0(12)$   &\circletext{$03$} & $(03)(12)$    \\   \cline{2-8}  
    
     & $01(23)$     &   $1(23)$  &     $(01)(23)$      &    $0(123)$  & $01$   &$0(23)$ & $(023)1$    \\   \cline{2-8}   
     
     & $012$     &   $12$  &     $(01)2$      &    $0(12)$  & $01$   &$02$ & $(02)1$    \\   \cmidrule[1pt]{2-8}  
   
     &$302$   & $02$  &    $(03)2$     &$3(02)$           &\circletext{$30$}  & $32$ & $(23)0$  \\ \cline{2-8}
     
    &  $(13)02$ &$02$   &  $(013)2$   & $(13)(02) $   &     $(13)0$    &$(13)2$  & $(123)0$   \\ \cline{2-8}  

    & $1(03)2$     &   $(03)2$  &     $(013)2$      &    $1(023)$  & $1(03)$   &$12$ & $(12)(03)$    \\   \cline{2-8}  
    
     & $13(02)$     &   $3(02)$  &     $(13)(02)$      &    $1(023)$  & $13$   &$1(02)$ & $(012)3$    \\   \cline{2-8}   
     
     & $130$     &  {$30$}  &     $(13)0$      &    $1(03)$  & $13$   &$10$ & $(01)3$    \\   \cline{2-8}   
 
    \end{tabular}\end{center}
    \end{table}
 The circled entries $03$ and $30$ are related by precomposition with the map $10:\underline{2} \ra \underline{2}.$ Thus, in the Glenn  table giving a sphere in $X$, the entries in the corresponding places must be related by the operator $X(10)=\gamma_{10}$.

It is not hard to see by a similar argument that for any tabular universal horn $\Lambda_S^c$ we can define a $\Lambda_S^c$-horn in $X$ by giving part of a Glenn table, omitting the rows which correspond to special coface maps which are in $S$. Even though these omitted rows of the Glenn table are not needed to define the horn, we will still list them in our table, but we indicate they are ``officially'' supposed to be omitted by marking them with a $\Lambda$ symbol. In these rows, we list any entries in these rows that can inferred from the rest of the table, which give us information about the faces of any potential filler of the horn.

\subsection{Fillers for spheres vs. fillers for horns}
Let $X$ be a $C $-set.
\begin{lemma} \label{fillerstosubcoskeletal} Suppose for every $c\in C $ with $\dim{c}> n$, there is a universal horn $\Lambda_S^c$ such that every horn of this type in $X$ has a unique filler. Then $X$ is $n$-subcoskeletal.
\end{lemma}
\begin{proof}Any filler of $d C [c]\ra X$ is also a filler of some horn of every type $\Lambda_S^c\ra d C [c] \ra X.$ The condition in the lemma thus ensures such sphere fillers are unique if $\dim c > n,$ thus $X$ is $n$-subcoskeletal.
\end{proof}
\begin{lemma} \label{fillerstocoskeletal1} Suppose $X$ is $(n-1)$-subcoskeletal and for every $c\in C $ with $\dim{c}\geq n$, there is a nice universal horn $\Lambda_S^c$ such that every horn of this type in $X$ has a filler. Then $X$ is $n$-coskeletal.
\end{lemma}
\begin{proof}
Suppose we have a sphere $s:d C [c] \ra X$ with $\dim c >n$.  Consider a nice horn $\Lambda_S^c$ such that every horn of this type has a filler. Consider the restriction $ h:\Lambda_S^c\ra X$ of $s$ to $\Lambda_S^c$, and let $\hat{h}$ be the filler of this horn. By the fact that $\Lambda_S^c$ is nice, we have $$\Sk^{n-1} \Lambda_S^c \cong \Sk^{n-1}  C [c] \cong \Sk^{n-1} d C [c],$$ with the equivalences given by the natural maps. So for the restricted maps, we have $$s = \hat{h}=h:\Sk^{n-1} d C [c]\ra X.$$ Applying the adjointness of $\Sk^{n-1}$ and $\Cosk^{n-1}$, we see that the compositions \begin{align*}s:d C [c]\ra X\ra \Cosk^{n-1} X \\ \hat{h}:d C [c] \ra X \ra \Cosk^{n-1} X
\end{align*} are equal. Since $X$ is $(n-1)$-subcoskeletal, the map $X\ra \Cosk^n X$ is monic, thus $s$ and $\hat{h}$ are equal as maps $d C [c]\ra X.$ Thus $\hat{h}$ is a filler of $s$.
\end{proof} 
\begin{corollary} \label{fillerstocoskeletal} Suppose for every $c\in C $ with $\dim{c}\geq n$, there is a nice universal horn $\Lambda_S^c$ such that every horn of this type in $X$ has a unique filler. Then $X$ is $n$-coskeletal.
\end{corollary}
\begin{proof}  $X$ is $(n-1)$-subcoskeletal by Lemma~\ref{fillerstosubcoskeletal} and so Lemma~\ref{fillerstocoskeletal1} applies. 
\end{proof}
\begin{lemma}\label{coskeletaliskan} Let $Y$ be a $C $-Set. Then every inner horn of type $c$ in $\cosk^k Y$ has a filler (has at most one filler) if and only if every inner horn $\tr^k(\Lambda) \ra Y$ of this type has a (has at most one) extension to a map $\tr^k( C [c])\ra Y$ along the map $\tr^k(\Lambda)\ra \tr^k( C [c]).$ 
\end{lemma}
\begin{proof} This statement is immediate from the fact that $\cosk^k$ is right adjoint to $\tr^k$. 
\end{proof}

\begin{lemma}\label{coskeletaltofiller} If $X$ is $n$-subcoskeletal, every horn in $X$ of minimal complementary dimension $n+1$ or greater has at most one filler. If $X$ is $n$-coskeletal, every horn in $X$ of minimal complementary dimension $n+1$ or greater has a unique filler. 
\end{lemma}
\begin{proof} Let $H:=\Lambda_S (c)$ be a universal horn of minimal complementary dimension $n+1$ or greater. Then $\Sk^n H \ra \Sk^n   C [c]$ is an isomorphism. So a filler of a $H$ in $X$ is equivalent to the extension of a map $f:\Sk^n   C [c]\ra X$ to $ C [c]$. Applying adjunction of $\Cosk^n$ with $\Sk^n$, this is equivalent to finding a map $f'$ in the following diagram:
\begin{center}
\begin{tikzcd}
 C [c] \arrow[dashed]{d}{f'}\arrow{dr}{f'} &       \\
X   \arrow{r}{can.} & \Cosk^n X  
\end{tikzcd}
\end{center}
If $X$ is $n$-subcoskeletal, the bottom canonical map is monic, so there is at most one such $f'$. If $X$ is $n$-coskeletal, the bottom map is an isomorphism, so $f'$ exists and is unique.
\end{proof}


\section{Multi-simplicial sets}

\begin{definition} A \emph{$k$-simplicial set} is a presheaf of sets on $\Delta^k$, i.e. a $\Delta^k$-set.
\end{definition} 
 be a multi-index, and let $\e^1, \e^2, \ldots, \e^k$ denote the standard basis of $\mathbb{N}^k$, viewed as elements of $\Delta^k$, i.e.  $\e^j=[0,0,\ldots, 0, 1, 0, \ldots, 0]$ where the $1$ is in the $j$th position.

A $k$-simplicial set $X$ can equivalently be viewed as a collection of sets $X_{\n}$ for all multi-indices $\n$,  with face maps $d^\alpha_i:X_{\n}\ra X_{\n - \e^\alpha}$  and degeneracy maps $s^\alpha_i:X_{\n} \ra X_{\n +\e^\alpha}$ for $0\leq \alpha \leq k$ and $0 \leq i \leq n_\alpha,$ all satisfying the following $k$-simplicial identities:
\begin{itemize}
\item $d^\alpha_id^\alpha_j=d^\alpha_{j-1}d_i.$ 
\item For $d^\alpha_is^\alpha_j$ : 
\begin{itemize}
\item $d^\alpha_is^\alpha_j=s^\alpha_{j-1}d^\alpha_i$ if $i<j.$
\item $d^\alpha_is^\alpha_j=\id$ \ if $i=j$ or $i=j+1.$
\item $d^\alpha_is^\alpha_j=s^\alpha_jd^\alpha_{i-1}$  if $i>j+1.$
\end{itemize}
\item $s^\alpha_is^\alpha_j=s^\alpha_js^\alpha_{i-1}.$
\item $d^\alpha_i$ and $s^\alpha_j$ each commute with $d^{\alpha'}_i$ and $d^{\alpha'}_j$ if $\alpha \neq \alpha'$.
\end{itemize}

\begin{definition} \label{dimensiondef} Let $\n=[\n^1,\ldots,\n^k]\in \Delta^k.$ The \emph{dimension} $\dim(\n)$ of $\n$ is $\sum\limits_{1\leq \alpha\leq k } \n^\alpha.$ This definition makes $\Delta^k$ a dimensional category.
\end{definition}

\begin{definition} Recall that a coface map in $\Delta^k$ is a monic map which increases dimension by $1$. The coface maps with target $\n$ have the form $$\hat{d}^\alpha_i :=[\id, \ldots, \id , \hat{d}^\alpha_i,\id, \ldots, \id]$$
where $\hat{d}^\alpha_i$ is a coface map in $\Delta$ with target $[\n^\alpha].$ We call $\alpha$ the \emph{direction} of $\alpha$.
\end{definition}
\begin{proposition} The dimensional category $\Delta^k$ is excellent.
\end{proposition}
\begin{proof} We leave the straightforward verification of this fact to the reader.
\end{proof}
We make this an ordered coface system by ordering the $\hat{d}^\alpha_i$ lexically by $\alpha$ and then by $i$.
\begin{definition} \label{multihorns}
The universal horn in $\Delta^k$ formed by removing $\hat{d}^\alpha_i$ from $\Delta^k[\n]$ is denoted $\Lambda^{\n}_{\{d^\alpha_i\}},$ in accordance with Definition~\ref{mincompdimdef}. If $0<i <\n^\alpha$, we say this horn is \emph{inner}. Similarly, an inner horn in a $\Delta^k$-set $X$ is a map from an inner universal horn to $X$. We call $X$ \emph{inner-Kan} if every inner horn in $X$ has a filler.
\end{definition}
It's not hard to see that every horn $\Lambda^{\n}_{\{d^\alpha_i\}}$ is nice. We call an inner-Kan $\Delta^k$-set $X$ \emph{$m$-reduced} if every inner horn in $X$ of minimal complementary dimension $m$ or greater has a unique filler, or equivalently if a $\Lambda^{\n}$-horn has a unique filler if $\dim(\n)>m.$

\subsection{The box product for multi-simplicial set}
For simplicial sets, we have a familiar adjunction $$\mathbf{const}: \cat{Set} \leftrightarrow \cat{Set}_{\Delta}:\mathbf{res}$$ where $\mathbf{Const}(X)_i=X$ and $\mathbf{res}(Y_\bullet)=Y_0$. Constant simplicial sets in the image of $\mathbf{const}$ are thought of as being discrete, since they are $0$-dimensional in the sense of have only degenerate simplices in dimensions greater than $0$. The functor $\mathbf{const}$ tells us how a set may be thought of as being a simplicial set, and so is often suppressed in notation, for instance if $X\in \cat{Set}$ and $Y \in \cat{Set}_\Delta$, then $X \times Y$ means $\mathbf{const}(X)\times Y$. For multi-simplicial sets, we can apply this adjunction for any index, or combination of indices.
\begin{definition} \label{constdef} Let $m\geq n$ and let $(i_1,\ldots, i_n)$ be an ordered subset of $(1,\ldots, m)$.  Then we define an adjunction $$\mathbf{const}^m_{i_1,\ldots,i_n}:\Set_{\Delta^{n}} \leftrightarrow \Set_{\Delta^{m}}:\mathbf{res}^m_{i_1,\ldots,i_n}$$ by the formula

\begin{itemize} \item $\mathbf{const}^m_{i_1,\ldots,i_n}(Y)_{q_1,\ldots, q_{m}}= Y_{q_{i_1},q_{i_2},\ldots, q_{i_n}}.$
 \item $\mathbf{res}^m_{i_1,\ldots,i_n}(X)_{p_1,\ldots, p_{n}}= X_{j_1, \ldots, j_m}$ where $j_l=\left\{
	\begin{array}{ll}
		0 & \mbox{if } l \notin (i_1,\ldots, i_n) \\
	p_k& \mbox{if } l = i_k \end{array}  \right. $    
\end{itemize} 
\end{definition}
\begin{definition} \label{boxdef} Let $X\in \Set_{\Delta^{k}}$ and $Y \in \Set_{\Delta^{l}}$ then define the \emph{box product} $X\square Y$ to be the $(k+l)$-simplicial set:
$$ X \square Y := \mathbf{const}^{k+l}_{1,\ldots, k}(  X  )\times \mathbf{const}^{k+l}_{k+1, \ldots, k+l} ( Y )$$
\end{definition} 

We now state a few basic observations about $\square$, the verifications of which we leave to the reader. 
\begin{proposition} The box product is associative (up to a canonical natural isomorphism).
\end{proposition}
\begin{proposition} \label{boxysimplicies} There is a canonical natural isomorphism $$\Delta^k[n_1,\ldots, n_k]\cong \Delta[n_1]\square \cdots \square \Delta[n_k]$$.
\end{proposition}
\begin{proposition} \label{boxytimes} If $X$ is a $k$-simplicial set, the functor $X \square -$ preserves colimits. Also, if $X, X'$ are $k$-simiplicial sets and $Y,Y'$ are $l$-simplicial sets, we have a natural isomorphism $$(X\square Y)\times (X'\square Y') \cong (X \times X') \square (Y\times Y') $$
\end{proposition}


\subsection{Notation and Glenn tables for bisimplicial sets}


In this chapter, we will be concerned specifically with bisimplicial sets, so we will adopt some helpful notation to use in this case. First, we write $d_i$ for $d^0_i$ and $s_i$ for $s^0_i.$ Likewise we write $\delta_i$ and $\varsigma_i$ for $d^1_i$ and $s^1_i$ respectively. Similarly we write $\hat{\delta}_i$ for $\hat{d}^1_i$, and etc. 

To write a general map in $\Delta^2$ with target $[m,n]$, we write the standard notation for the first component of the map above the standard notation for the second component. For instance, $$ \hat{d}_1\hat{\varsigma}_0 :[1,3]\ra [2,2]$$ is denoted ${\binom{02}{0012}}.$

Now we adopt a notation allowing us to name a generic cell of a bisimplicial set $X$ along with all of its faces and iterated faces.  $x_{01\ldots n}^{01\ldots m}$ indicates a $[m,n]$-cell of $X$, and we write for instance $$x^{02}_{0012}:=  X\binom{02}{0012}  (x_{01\ldots n}^{01\ldots m}).$$ As before, we will often wish to make a list of all the faces of a cell, its \emph{boundary}, though now a cell in general has two kinds of faces, those given by applying a $d_i$, the \emph{horizontal faces} and those given by applying a $\delta_i$, the \emph{vertical faces}. We will always list horizontal faces first, and divide horizontal faces from vertical faces by a $\mathbf{|}$ symbol, for instance: $$d(x_{012}^{012})=[x^{12}_{012} ,\ x^{02}_{012},\  x^{01}_{012}  \ \  | \ \ x^{012}_{12},\ x^{012}_{02},\ x^{012}_{01}].$$ Such a list of faces satisfying appropriate face relations  necessary to be a boundary of an $[m,n]$-cell in $X$ is called an \emph{$d\Delta^2[m,n]$-sphere} in $X$. Note that such a $d\Delta^2[m,n]$-sphere is equivalent to a sphere as defined in Definition~\ref{spheredef}, a map $d\Delta^2[m,n]\ra X$. An \emph{filler} of an $d\Delta^2[m,n]$-sphere $d\Delta^2[m,n]\ra X$ is an extension of the sphere to a map $\Delta^2[m,n] \ra X$, or equivalently, a $[m,n]$-cell that has the specified sphere as its boundary. A $d\Delta^2[m,n]$-sphere in a bisimplicial set is called $\emph{commutative}$ if it has a filler.

As before, we will often wish to construct spheres and especially horns from individual cells of a bisimplicial set $X$. To do this, we will use the technique of Glenn tables described above in Section~\ref{genglennsec}, using the ordered system of cofaces given above, where the faces of a $[m,n]$-cells are ordered $$d_0, d_1,\ldots, d_m, \delta_0, \delta_1, \ldots, \delta_n.$$ To illustrate, the universal Glenn table for $[3,3]\in \Delta^2$ is given below:

\begin{table}[H]\begin{center} \caption{The universal Glenn table for $[3,3]$ \label{univtableex}}

    \begin{tabular}{ r | l || l | l !{\vrule width 2pt} l | l | l |}
    \cline{2-7}

        &  $  \binom{12}{012} $  &  $  \binom{2}{012} $  &   $  \binom{1}{012} $   &     $  \binom{12}{12}  $    &    $  \binom{12}{02} $   &  $  \binom{12}{01}  $        \\ \cline{2-7}

        &  $  \binom{02}{012} $  &  $  \binom{2}{012} $  &   $  \binom{0}{012} $   &     $  \binom{02}{12}  $    &    $  \binom{02}{02} $   &  $  \binom{02}{01}  $       \\ \cline{2-7}

        &  $  \binom{01}{012} $  &  $  \binom{1}{012} $  &   $  \binom{0}{012} $   &     $  \binom{01}{12}  $    &    $  \binom{01}{02} $   &  $  \binom{01}{01}  $         \\ \cline{2-7} 
          
    \end{tabular} \vspace{.5cm}
   
 \begin{tabular}{    r | l || l| l |  l !{\vrule width 2pt} l | l  | }\cline{2-7} 
  &  $  \binom{012}{12} $&   $  \binom{12}{12} $  &   $  \binom{02}{12} $   &     $  \binom{01}{12}  $    &    $  \binom{012}{2} $   &  $  \binom{012}{1}  $    \\ \cline{2-7} 
 &  $  \binom{012}{02} $&   $  \binom{12}{02} $  &   $  \binom{02}{02} $   &     $  \binom{01}{02}  $    &    $  \binom{012}{2} $   &  $  \binom{012}{0}  $     \\ \cline{2-7}
 &  $  \binom{012}{01} $  &   $  \binom{12}{01} $  &   $  \binom{02}{01} $   &     $  \binom{01}{01}  $    &    $  \binom{012}{1} $   &  $  \binom{012}{0}  $    \\ \cline{2-7}
\end{tabular}   
\end{center}
  \end{table}

We have seperated this table into two parts, the top for the $d$ faces and the bottom for the $\delta$ face. Both are parts of one Glenn table, with the separation done merely to make the table easier to read.

The relations \begin{align*} d_jd_i&= d_{i}d_{j+1},  & 0 \leq i \leq j \leq {m-1}. \\ \delta_j\delta_i &= \delta_{i}\delta_{j+1},  & 0 \leq i \leq j \leq {n-1}.   \end{align*} correspond to the agreement of the $(i,j)$th entry on or above the main diagonal to the $(j+1,i)$th entry below the diagonal for the left part of the top and the right part of the bottom. This is the same symmetry relation that Glenn tables for simplicial sets satisfy. The relation $d_j\delta_i=\delta_id_j$ corresponds to agreement between the $(i,j)$th entry of the right part of the top table with the $(j,i)$th entry of the left part of the bottom table. That is, the ``middle'' two parts of our array are transpose matrices, because of the fact that $d$ and $\delta$ commute with each other.  The same pattern describes the ``matching pattern'' for every universal Glenn  table for $\Delta^2$, and is thus the pattern we must check in every Glenn table for $\Delta^2$.

Often, we will consider spheres or horns where $m=1$ or $n=1$, in which case the far left or far right parts of the array will be missing. For instance, below is the Glenn  table for the sphere $d(x_{012}^{01}):$ 

\begin{table}[H]\begin{center}

    {\begin{tabular}{ r | l || l | l !{\vrule width 2pt} l | l  |}
    \cline{2-6}

        &  $x^{12}_{01}$  &  $x^{2}_{01}$  &   $x^{1}_{01}$   &     $x^{12}_{1} $    &    $x^{12}_{0}$       \\ \cline{2-6} 
     
        &  $x^{02}_{01}$  &  $x^{2}_{01}$  &   $x^{0}_{01}$   &     $x^{02}_{1} $    &    $x^{02}_{0}$     \\ \cline{2-6} 
 
        &  $x^{01}_{01}$  &  $x^{1}_{01}$  &   $x^{0}_{01}$   &     $x^{01}_{1} $    &    $x^{01}_{0}$      \\ \cline{2-6}  \end{tabular}}
       \vspace{.5cm}
       
\begin{tabular}{r | l || l | l | l !{\vrule width 2pt}}
    \cline{2-5}

       &  $x^{012}_{1}$ &     $x^{12}_{1} $    &    $x^{02}_{1}$   &  $x^{01}_{1} $        \\ \cline{2-5} 
     
        &  $x^{012}_{1}$     &     $x^{12}_{0} $    &    $x^{02}_{0}$   &  $x^{01}_{0} $     \\ \cline{2-5} 
 
  \end{tabular}

 \end{center}   \end{table}

\begin{definition}
As described in Definition~\ref{multihorns}, we can make a horn from the sphere $d \Delta^2[m,n]$ we can make \emph{horns} by removing a face.  $\Lambda[m^{\, i}, n]$ will serve as a more convenient notation for the horn given by removing the face $d_i$ from  $d \Delta^2[m,n],$ which is the horn $\Lambda^{[m,n]}_{\{d^0_i\}}$ as defined in Definition~\ref{mincompdimdef}. Similarly  $\Lambda[m, n^{\, j}]$ will be our new notation for $\Lambda^{[m,n]}_{\{d^1_j\}}.$ Again, a horn is said to be \emph{inner} if $0<i<m$ in the first case, or $0<j<n$ in the second case.
\end{definition}

As before, we will use Glenn tables to build horns inside a bisimplicial set $X$. A Glenn  table provides an easy way to visually check that a given list of horizontal and vertical faces satisfies the appropriate simplicial identities to make a horn, and gives the boundary of the ``missing face'' that is to be filled in. We will signify the missing face of a horn in its Glenn  table by the symbol $\Lambda.$

\begin{lemma}\label{kaniscoskeletal}If a inner-Kan bisimplicial set $X$ is $2$-reduced, then $X$ is $2$-subcoskeletal and  $3$-coskeletal.
\end{lemma}
\begin{proof} By Lemma~\ref{fillerstosubcoskeletal} and Corollary~\ref{fillerstocoskeletal}.
\end{proof}

\section{Verity double categories}


\begin{definition} A \emph{strict double category} is  a category internal to the category $\cat{Cat}.$  A \emph{pseudo-double category} is defined to be a category internal the bicategory $\cat{Cat}$, in the appropriate sense of being internal to a bicategory. A double category therefore consists of two categories $\mathcal{C}_1$, $\mathcal{C}_0$ with functors $s,t: \mathcal{C}_1\ra \mathcal{C}_0$ and $\id:\mathcal{C}_0 \ra \mathcal{C}_1$, together with a \emph{composition}: $$c:\mathcal{C}_1 \underset{\mathcal{C}_0}{\times}\mathcal{C}_1\ra \mathcal{C}_1.$$ We call the morphisms of $\mathcal{C}_0$ the \emph{vertical $1$-morphisms}, and the objects of $\mathcal{C}_1$ the \emph{horizontal $1$-morphisms}. In a strict double category, we assume this functor $c$ satisfies an associativity condition, while in a \emph{pseudo-double category} we assume the existence of an \emph{associator} natural transformation. 
\end{definition}
Note that there is an asymmetry in the definition of a pseudo double category between the two directions of $1$-morphisms. Namely, composition for vertical $1$-morphisms, which is given by composition in $\mathcal{C}_0,$ is strictly associative, while the composition of horizontal $1$-morphisms, given by $c$, is only associative up to the associator, whose components are morphism of $\mathcal{C}_1$. This asymmetry is present in many naturally occurring double categories, for instance the double category where $\mathcal{C}_0$ is the category of commutative rings and ring morphisms, and $\mathcal{C}_1$ is the category of bimodules and bimodule morphisms. However, this asymmetry is problematic from the point of view of relating these double categories to bisimplicial sets, which are symmetrical in their two directions. This is our immediate motive for considering Verity double categories, a kind of double category that is weak in both directions, introduced by Dominic Verity as ``double bicategories'' in his 1992 Ph.D. thesis, republished as \cite{Ver11}. We loosely follow Jeffery Morton's exposition of the definition \cite{Mor09}, which is responsible for bringing Verity's definition to the attention of the author. 
\begin{definition}[Verity] \label{Veritydef} A Verity double category (also called a VDC) consists of the following data:
\begin{itemize}
\item $(2,1)$-categories $H$ and $V$ sharing a class of objects $O$. (For Verity, $H$ and $V$ are only assumed to be bicategories, but we restrict to the case where the $2$-morphisms of $H$ and $V$ are invertible.) If $O$ is a set, we say our VDC is \emph{small}.
\item A set of squares $\mbox{Sq}(f,f',p,p')$ for all $f,f',p,p'$ where $f: a \ra b$ and $f':a'\ra b'$ are $1$-morphisms in $H$, and  $p:a \ra a'$ and $p':b\ra b'$ are $1$-morphisms in $V$. We picture such a square $\Theta$ as being ``inside'' the square made by the four morphisms, as shown:

\begin{center}
\begin{tikzpicture}[scale=1.8,auto]

\begin{scope}

\node (10) at (1,1) {$b$};
\node (00) at (0,1) {$a$};
\node (11) at (1,0) {$b'$};
\node (01) at (0,0) {$a'$};
\node[rotate=-45] at (.5,.5){$\Rightarrow$};
\node[scale=.8] at (.6,.65){$\Theta$};

\path[->] (00) edge node[midway]{$f$}(10);
\path[->] (00) edge node[midway,swap]{$p$}(01);
\path[->] (01) edge node[midway,swap]{$f'$}(11);
\path[->] (10) edge node[midway]{$p'$}(11);
\end{scope}

\end{tikzpicture}
\end{center}
\end{itemize}
The squares can be composed and acted on by $2$-morphisms of $H$ and $V$ in the following ways:
\begin{enumerate}\item Horizontal composition of squares $\Pi \boxvert \Theta$, (a mnemonic for $\boxvert$ is that the symbol is a little picture of two side-by-side horizontally composible squares):

\begin{center}
\begin{tikzpicture}[scale=1.8,auto]

\begin{scope}

\node (10) at (1,1) {$b$};
\node (00) at (0,1) {$a$};
\node (11) at (1,0) {$b'$};
\node (01) at (0,0) {$a'$};
\node[rotate=-45] at (.5,.5){$\Rightarrow$};
\node[scale=.8] at (.6,.65){$\Theta$};

\node (20) at (2,1) {$c$};
\node (21) at (2,0) {$c'$};
\node[rotate=-45] at (1.5,.5){$\Rightarrow$};
\node[scale=.8] at (1.6,.65){$\Pi$};

\path[->] (00) edge node[midway]{$f$}(10);
\path[->] (00) edge node[midway,swap]{$p$}(01);
\path[->] (01) edge node[midway,swap]{$f'$}(11);
\path[->] (10) edge node[midway]{$p'$}(11);
\path[->] (10) edge node[midway]{$g$}(20);
\path[->] (11) edge node[midway,swap]{$g'$}(21);
\path[->] (20) edge node[midway]{$p''$}(21);
\node at (2.75,.5){$\leadsto$};
\end{scope}
\begin{scope}[shift={(3.5,0)}]
\node (10') at (1,1) {$c$};
\node (00') at (0,1) {$a$};
\node (11') at (1,0) {$c'$};
\node (01') at (0,0) {$a'$};
\path[->] (00') edge node[midway]{$g\circ f$}(10');
\path[->] (00') edge node[midway,swap]{$p$}(01');
\path[->] (01') edge node[midway,swap]{$g'\circ f'$}(11');
\path[->] (10') edge node[midway]{$p''$}(11');
\node[rotate=-45] at (.5,.43){$\Rightarrow$};
\node[scale=.8] at (.55,.63){$\Pi \boxvert \Theta$};

\end{scope}
\end{tikzpicture}

\end{center}
\item Vertical composition of squares, $\Pi \boxminus \Theta$:

\begin{center}
\begin{tikzpicture}[scale=1.8,auto]

\begin{scope}

\node (10) at (1,1) {$b$};
\node (00) at (0,1) {$a$};
\node (11) at (1,0) {$b'$};
\node (01) at (0,0) {$a'$};
\node[rotate=-45] at (.5,.5){$\Rightarrow$};
\node[scale=.8] at (.6,.65){$\Theta$};

\node (11') at (1,-1) {$b''$};
\node (01') at (0,-1) {$a''$};
\node[rotate=-45] at (0.5,-.58){$\Rightarrow$};
\node[scale=.8] at (0.6,-.43){$\Pi$};

\path[->] (00) edge node[midway]{$f$}(10);
\path[->] (00) edge node[midway,swap]{$p$}(01);
\path[->] (01) edge node[midway,swap]{$f'$}(11);
\path[->] (10) edge node[midway]{$p'$}(11);
\path[->] (11) edge node[midway]{$q'$}(11');
\path[->] (01) edge node[midway,swap]{$q$}(01');
\path[->] (01') edge node[midway,swap]{$f''$}(11');
\node at (1.75,0){$\leadsto$};
\end{scope}
\begin{scope}[shift={(2.7,-.5)}]
\node (10') at (1,1) {$b$};
\node (00') at (0,1) {$a$};
\node (11') at (1,0) {$b''$};
\node (01') at (0,0) {$a''$};
\path[->] (00') edge node[midway]{$f$}(10');
\path[->] (00') edge node[midway,swap]{$q\circ p$}(01');
\path[->] (01') edge node[midway,swap]{$f''$}(11');
\path[->] (10') edge node[midway]{$q' \circ p'$}(11');
\node[rotate=-45] at (.5,.43){$\Rightarrow$};
\node[scale=.8] at (.55,.63){$\Pi \boxminus \Theta$};

\end{scope}
\end{tikzpicture}
\end{center}
\item For every horizontal $1$-morphism $f$, a \emph{pseudo-identity square for vertical composition} $\ID_f$ and for every vertical $1$-morphism $p$, a \emph{pseudo-identity square for horizontal composition} $\ID_p$:

\begin{center}
\begin{tikzpicture}[scale=1.8,auto]

\begin{scope}

\node (10) at (1,1) {$b$};
\node (00) at (0,1) {$a$};
\node (11) at (1,0) {$b$};
\node (01) at (0,0) {$a'$};
\node[rotate=-45] at (.5,.5){$\Rightarrow$};
\node[scale=.8] at (.6,.65){$\ID_f$};

\path[->] (00) edge node[midway]{$f$}(10);
\path[->] (00) edge node[midway,swap]{$\id_a$}(01);
\path[->] (01) edge node[midway,swap]{$f$}(11);
\path[->] (10) edge node[midway]{$\id_b$}(11);
\end{scope}
\begin{scope}[shift={(2.5,0)}]
\node (10') at (1,1) {$a$};
\node (00') at (0,1) {$a$};
\node (11') at (1,0) {$a'$};
\node (01') at (0,0) {$a'$};
\path[->] (00') edge node[midway]{$\id_a$}(10');
\path[->] (00') edge node[midway,swap]{$p$}(01');
\path[->] (01') edge node[midway,swap]{$\id_{a'}$}(11');
\path[->] (10') edge node[midway]{$p$}(11');
\node[rotate=-45] at (.5,.43){$\Rightarrow$};
\node[scale=.8] at (.55,.63){$\ID_p$};

\end{scope}
\end{tikzpicture}

\end{center}

\item Action of horizontal 2-morphisms on the top, $\Theta \filledmedtriangleup \beta$:

\begin{center}
\begin{tikzpicture}[scale=1.8,auto]

\begin{scope}

\node (10) at (1,1) {$b$};
\node (00) at (0,1) {$a$};
\node (11) at (1,0) {$b'$};
\node (01) at (0,0) {$a'$};
\node[rotate=-45] at (.5,.5){$\Rightarrow$};
\node[scale=.8] at (.6,.65){$\Theta$};

\path[->] (00) edge node (g)[midway,below]{$g$}(10);
\draw[->] (00) to [out=60,in=120] node(f)[midway]{$f$} (10) ;
\path[->] (00) edge node[midway,swap]{$q$}(01);
\path[->] (01) edge node[midway,swap]{$g'$}(11);
\path[->] (10) edge node[midway]{$q'$}(11);
\node[rotate=-90] (areta) at ($(f)+(0,-.38)$){$\Rightarrow$};
\node (eta) at ($(areta)+(.15,0)$){$\beta$};
\node at (1.75,.5){$\leadsto$};
\end{scope}
\begin{scope}[shift={(2.5,0)}]
\node (10') at (1,1) {$b$};
\node (00') at (0,1) {$a$};
\node (11') at (1,0) {$b'$};
\node (01') at (0,0) {$a'$};

\path[->] (00') edge node[midway]{$f$}(10');
\path[->] (00') edge node[midway,swap]{$q$}(01');
\path[->] (01') edge node[midway,swap]{$g'$}(11');
\path[->] (10') edge node[midway]{$q'$}(11');
\node[rotate=-45] at (.5,.43){$\Rightarrow$};
\node[scale=.8] at (.55,.63){$\Theta\filledmedtriangleup\beta  $};

\end{scope}
\end{tikzpicture}

\end{center}
\item Action of horizontal 2-morphisms on the bottom, $ \Theta \filledmedtriangledown \beta' $:

\begin{center}
\begin{tikzpicture}[scale=1.8,auto]

\begin{scope}

\node (10) at (1,1) {$b$};
\node (00) at (0,1) {$a$};
\node (11) at (1,0) {$b'$};
\node (01) at (0,0) {$a'$};
\node[rotate=-45] at (.5,.5){$\Rightarrow$};
\node[scale=.8] at (.6,.65){$\Theta$};

\path[->] (00) edge node (g)[midway]{$g$}(10);
\draw[->] (01) to [out=-60,in=-120] node(f)[midway,below]{$f'$} (11) ;
\path[->] (00) edge node[midway,swap]{$q$}(01);
\path[->] (01) edge node[midway]{$g'$}(11);
\path[->] (10) edge node[midway]{$q'$}(11);
\node[rotate=90] (areta) at ($(f)+(0,+.4)$){$\Rightarrow$};
\node (eta) at ($(areta)+(.15,0)$){$\beta'$};
\node at (1.75,.5){$\leadsto$};
\end{scope}
\begin{scope}[shift={(2.5,0)}]
\node (10') at (1,1) {$b$};
\node (00') at (0,1) {$a$};
\node (11') at (1,0) {$b'$};
\node (01') at (0,0) {$a'$};

\path[->] (00') edge node[midway]{$g$}(10');
\path[->] (00') edge node[midway,swap]{$q$}(01');
\path[->] (01') edge node[midway,swap]{$f'$}(11');
\path[->] (10') edge node[midway]{$q'$}(11');
\node[rotate=-45] at (.5,.43){$\Rightarrow$};
\node[scale=.8] at (.55,.63){$\Theta\filledmedtriangledown \beta' $};

\end{scope}
\end{tikzpicture}
\end{center}

\item Action of vertical 2-morphisms on the left, $ \Theta\filledmedtriangleleft\eta $:

\begin{center}
\begin{tikzpicture}[scale=1.8,auto]

\begin{scope}

\node (10) at (1,1) {$b$};
\node (00) at (0,1) {$a$};
\node (11) at (1,0) {$b'$};
\node (01) at (0,0) {$a'$};
\node[rotate=-45] at (.5,.5){$\Rightarrow$};
\node[scale=.8] at (.6,.65){$\Theta$};

\path[->] (00) edge node (g)[midway]{$g$}(10);
\draw[->] (00) to [out=-150,in=150] node(f)[midway,left]{$p$} (01) ;
\path[->] (00) edge node[midway]{$q$}(01);
\path[->] (01) edge node[midway,below]{$g'$}(11);
\path[->] (10) edge node[midway]{$q'$}(11);

\node[rotate=0] (areta) at ($(f)+(.35,0)$){$\Rightarrow$};
\node (eta) at ($(areta)+(0,.15)$){$\eta$};
\node at (1.75,.5){$\leadsto$};
\end{scope}
\begin{scope}[shift={(2.5,0)}]
\node (10') at (1,1) {$b$};
\node (00') at (0,1) {$a$};
\node (11') at (1,0) {$b'$};
\node (01') at (0,0) {$a'$};

\path[->] (00') edge node[midway]{$g$}(10');
\path[->] (00') edge node[midway,swap]{$p$}(01');
\path[->] (01') edge node[midway,swap]{$g'$}(11');
\path[->] (10') edge node[midway]{$q'$}(11');
\node[rotate=-45] at (.5,.43){$\Rightarrow$};
\node[scale=.8] at (.55,.63){$\Theta\filledmedtriangleleft \eta $};

\end{scope}
\end{tikzpicture}
\end{center}
\item Action of vertical 2-morphisms on the right, $ \Theta \filledmedtriangleright\eta'$:

\begin{center}
\begin{tikzpicture}[scale=1.8,auto]
\begin{scope}
\node (10) at (1,1) {$b$};
\node (00) at (0,1) {$a$};
\node (11) at (1,0) {$b'$};
\node (01) at (0,0) {$a'$};
\node[rotate=-45] at (.5,.5){$\Rightarrow$};
\node[scale=.8] at (.6,.65){$\Theta$};

\path[->] (00) edge node (g)[midway]{$g$}(10);
\draw[->] (10) to [out=-30,in=30] node(f)[midway,right]{$p'$} (11) ;
\path[->] (00) edge node[midway,left]{$q$}(01);
\path[->] (01) edge node[midway,below]{$g'$}(11);
\path[->] (10) edge node[midway,left]{$q'$}(11);
\node[rotate=180] (areta) at ($(f)+(-.35,0)$){$\Rightarrow$};
\node (eta) at ($(areta)+(0,.15)$){$\eta'$};
\node at (2,.5){$\leadsto$};
\end{scope}
\begin{scope}[shift={(2.75,0)}]
\node (10') at (1,1) {$b$};
\node (00') at (0,1) {$a$};
\node (11') at (1,0) {$b'$};
\node (01') at (0,0) {$a'$};

\path[->] (00') edge node[midway]{$g$}(10');
\path[->] (00') edge node[midway,swap]{$q$}(01');
\path[->] (01') edge node[midway,swap]{$g'$}(11');
\path[->] (10') edge node[midway]{$p'$}(11');
\node[rotate=-45] at (.5,.43){$\Rightarrow$};
\node[scale=.8] at (.55,.63){$\Theta \filledmedtriangleright \eta'$};

\end{scope}
\end{tikzpicture}
\end{center}
\end{enumerate}
\end{definition}
\begin{remark} We have defined $\Theta\filledmedtriangledown \beta' $ and $ \Theta\filledmedtriangleright\eta' $ using $2$-morphisms going the opposite direction from those used to define the corresponding actions in \cite{Ver11}. Each of our actions by 2-morphisms will therefore be right (contravariant) actions on the set of squares, whereas for Verity the bottom and right actions are covariant. However, since we assume $H$ and $V$ are $(2,1)$-categories and thus have invertible 2-morphisms, the two kinds of actions are equivalent. This convention will make our constructions easier later, but Verity's convention is perhaps more natural.
\end{remark}
\begin{convention} When necessary, we will distinguish the structures of $H$ and $V$ using subscripts $h$ and $v$, for instance $\bullet_h$ is $2$-morphism composition in $H$. However, we will almost always supress this notation. To help avoid confusion, we will always name horizontal $1$-morphisms $e,f,g,h,$ vertical $1$-morphisms $p,q,r,s,$ horizontal $2$-morphisms $\alpha,\beta,\gamma,$ and vertical $2$-morphsims $\eta,\theta,\pi.$
\end{convention}

These operations are subject to the following conditions:
\begin{itemize}
\item  $\filledmedtriangleup$ and $\filledmedtriangledown$ are actions of the $2$-morphisms of $H$ and $\filledmedtriangleleft$ and $\filledmedtriangleright$ are actions of the $2$-morphisms of $V$. This means:
\begin{itemize} 
\item \textbf{VDC1}: \ The actions are unital, e.g. $ \Theta\filledmedtriangleup \Id_g=\Theta.$
\item \textbf{VDC2}: \  The actions are associative, e.g. $  \Theta \filledmedtriangleup(\beta\bullet \alpha)= (\Theta \filledmedtriangleup \beta)\filledmedtriangleup \alpha.$
\end{itemize}
\item \textbf{VDC3}: \ The four actions $\filledmedtriangleup,$ $\filledmedtriangledown,$ $\filledmedtriangleleft,$ and $\filledmedtriangleright$ commute with each other, e.g. $ (\Theta \filledmedtriangleup \beta)\filledmedtriangleleft\eta =( \Theta\filledmedtriangleleft\eta ) \filledmedtriangleup \beta .$
\item \textbf{VDC4}: \  The four actions are compatible with composition of squares in the appropriate direction, e.g. $(\Pi \boxvert\Theta ) \filledmedtriangleleft \eta=\Pi \boxvert(\Theta  \filledmedtriangleleft \eta).$
\item \textbf{VDC5}: \  The pseudo-identity is compatible with the identities for objects. For any object $a$, $\ID_{(\Id_a)_h}=\ID_{(\Id_a)_v}.$
\item \textbf{VDC6}: \ The pseudo-identity is compatible with composition of $1$ morphisms. If $a\stackrel{f}{\ra}b \stackrel{g}{\ra}c$ are horizontal $1$-morphisms, we have $\ID_g \boxvert \ID_f =\ID_{g\circ f},$ and similarly for vertical 2-morphisms.
\item \textbf{VDC7}: \  The pseudo-identity is compatible with the actions. For a horizontal $2$-morphism $\beta: f \Rightarrow g$ this means $\ID_g \uact \beta = \ID_f \dact \beta^{-1}$ or equivalently $(\ID_g \uact \beta)\dact \beta = \ID_f.$ Similarly for a vertical $2$-morphism $\eta: p \Rightarrow q$ we have $(\ID_q \lact\eta)\ract \eta=\ID_p.$ 
\item \textbf{VDC8}: \ Interchange laws for $\filledmedtriangleup,$ $\filledmedtriangledown,$ $\filledmedtriangleleft,$ and $\filledmedtriangleright$ and whiskering are satisfied. To give the interchange laws for $\uact$, suppose we have squares $\Theta$ and $\Pi$ along with horizontal $2$-morphism $\beta:e\Rightarrow g$ as pictured below:

\begin{center}
\begin{tikzpicture}[scale=1.8,auto]

\begin{scope}

\node (10) at (1,1) {$b$};
\node (00) at (0,1) {$a$};
\node (11) at (1,0) {$b'$};
\node (01) at (0,0) {$a'$};
\node[rotate=-45] at (.5,.45){$\Rightarrow$};
\node[scale=.8] at (.6,.6){$\Theta$};

\node (20) at (2,1) {$c$};
\node (21) at (2,0) {$c'$};
\node[rotate=-45] at (1.5,.5){$\Rightarrow$};
\node[scale=.8] at (1.6,.65){$\Pi$};

\path[->] (00) edge node[midway,below]{$g$}(10);
\path[->] (00) edge node[midway,swap]{$q$}(01);
\path[->] (01) edge node[midway,swap]{$g'$}(11);
\path[->] (10) edge node[midway]{$q'$}(11);
\path[->] (10) edge node[midway]{$h$}(20);
\path[->] (11) edge node[midway,swap]{$h'$}(21);
\path[->] (20) edge node[midway]{$q''$}(21);

\draw[->] (00) to [out=60,in=120] node(f)[midway]{$e$} (10) ;

\node[rotate=-90] (areta) at ($(f)+(0,-.32)$){$\Rightarrow$};
\node (eta) at ($(areta)+(.15,0)$){$\beta$};
\end{scope}
\end{tikzpicture}
\end{center}
then $\Pi \boxvert (\Theta \uact \beta)= (\Pi\boxvert \Theta)\uact (h \rhd \beta).$ Similarly, if $\Theta$ and $\Pi$ are as above and we instead have a horizontal $2$-morphism $\gamma:f \Rightarrow h$: 
\begin{center}
\begin{tikzpicture}[scale=1.8,auto]

\begin{scope}

\node (10) at (1,1) {$b$};
\node (00) at (0,1) {$a$};
\node (11) at (1,0) {$b'$};
\node (01) at (0,0) {$a'$};
\node[rotate=-45] at (.5,.5){$\Rightarrow$};
\node[scale=.8] at (.6,.65){$\Theta$};

\node (20) at (2,1) {$c$};
\node (21) at (2,0) {$c'$};
\node[rotate=-45] at (1.5,.45){$\Rightarrow$};
\node[scale=.8] at (1.6,.6){$\Pi$};

\path[->] (00) edge node[midway,above]{$g$}(10);
\path[->] (00) edge node[midway,swap]{$q$}(01);
\path[->] (01) edge node[midway,swap]{$g'$}(11);
\path[->] (10) edge node[midway]{$q'$}(11);
\path[->] (10) edge node[midway,below]{$h$}(20);
\path[->] (11) edge node[midway,swap]{$h'$}(21);
\path[->] (20) edge node[midway]{$q''$}(21);

\draw[->] (10) to [out=60,in=120] node(f)[midway]{$f$} (20) ;

\node[rotate=-90] (areta) at ($(f)+(0,-.38)$){$\Rightarrow$};
\node (eta) at ($(areta)+(0.15,0)$){$\gamma$};

\end{scope}
\end{tikzpicture}
\end{center}
then $(\Pi \uact \gamma) \boxvert \Theta = (\Pi\boxvert \Theta)\uact (\gamma\lhd g).$ Similar interchange laws hold for $\dact, \lact,$ and $\ract.$

\item \textbf{VDC9}:  \  Horizontal and vertical composition of squares are compatible with the associators of $H$ and $V$ respectively. In other words, horizontal (vertical) composition of squares is associative up to the associator of $H$ ($V$). For horizontal composition, suppose we have three horizontally composible squares $\Theta,$ $\Pi$ and $\Sigma$:

\begin{center}
\begin{tikzpicture}[scale=1.8,auto]

\begin{scope}

\node (10) at (1,1) {$b$};
\node (00) at (0,1) {$a$};
\node (11) at (1,0) {$b'$};
\node (01) at (0,0) {$a'$};
\node[rotate=-45] at (.5,.5){$\Rightarrow$};
\node[scale=.8] at (.6,.65){$\Theta$};

\node (20) at (2,1) {$c$};
\node (21) at (2,0) {$c'$};
\node (30) at (3,1) {$d$};
\node (31) at (3,0) {$d'$};
\node[rotate=-45] at (1.5,.5){$\Rightarrow$};
\node[scale=.8] at (1.6,.65){$\Pi$};
\node[rotate=-45] at (2.5,.5){$\Rightarrow$};
\node[scale=.8] at (2.6,.65){$\Sigma$};

\path[->] (00) edge node[midway,above]{$f$}(10);
\path[->] (00) edge node[midway,swap]{$q$}(01);
\path[->] (01) edge node[midway,swap]{$f'$}(11);
\path[->] (10) edge node[midway]{$q'$}(11);
\path[->] (10) edge node[midway]{$g$}(20);
\path[->] (11) edge node[midway,swap]{$g'$}(21);
\path[->] (20) edge node[midway]{$q''$}(21);
\path[->] (20) edge node[midway]{$h$}(30);
\path[->] (21) edge node[midway,swap]{$h'$}(31);
\path[->] (30) edge node[midway]{$q'''$}(31);

\end{scope}
\end{tikzpicture}
\end{center}
then  $((\Sigma \boxvert \Pi) \boxvert \Theta)\uact \alpha_{h,g,f}\dact \alpha_{h',g',f'}=\Sigma \boxvert (\Pi \boxvert \Theta)$, (where $\alpha_{-,-,-}$ is the associator of $H$):

\begin{center}
\begin{tikzpicture}[scale=1.8,auto]
\begin{scope}[shift={(3.5,0)}]

\node (00) at (0,1) {$a$};

\node (01) at (0,0) {$a'$};
\node[rotate=-45] at (.9,.38){$\Rightarrow$};
\node[scale=.8] at (1,.57){$\Sigma \boxvert (\Pi\boxvert \Theta)$};

\node (20) at (2,1) {$d$};
\node (21) at (2,0) {$d'$};

\path[->] (00) edge node[midway,above]{\scriptsize $h\circ (g\circ f)$}(20);
\path[->] (00) edge node[midway,swap]{$q$}(01);
\path[->] (01) edge node[midway,below]{\scriptsize $h'\circ (g'\circ f')$}(21);
\path[->] (20) edge node[midway]{$q'''$}(21);

\end{scope}\node (eq) at (2.75,.5){$=$};
\begin{scope}

\node (00) at (0,1) {$a$};

\node (01) at (0,0) {$a'$};
\node[rotate=-45] at (.9,.38){$\Rightarrow$};
\node[scale=.8] at (1,.57){$(\Sigma \boxvert \Pi)\boxvert \Theta$};

\node (20) at (2,1) {$d$};
\node (21) at (2,0) {$d'$};

\path[->] (00) edge node[midway,below]{\scriptsize $(h\circ g)\circ f$}(20);
\path[->] (00) edge node[midway,swap]{$q$}(01);
\path[->] (01) edge node[midway,above]{\scriptsize $(h'\circ g')\circ f'$}(21);
\path[->] (20) edge node[midway]{$q'''$}(21);

\draw[->] (01) to [out=-45,in=-135] node(f)[midway,below]{\scriptsize $h'\circ (g' \circ f')$} (21) ;

\node[rotate=+90] (areta) at ($(f)+(-.15,+.45)$){$\Rightarrow$};
\node (eta) at ($(areta)+(.42,0)$){\scriptsize $\alpha_{h',g',f'}$};

\draw[->] (00) to [out=45,in=135] node(f)[midway]{\scriptsize $h\circ (g \circ f)$} (20) ;

\node[rotate=-90] (areta) at ($(f)+(0,-.42)$){$\Rightarrow$};
\node (eta) at ($(areta)+(.35,0)$){\scriptsize $\alpha_{h,g,f}$};
\end{scope}
\end{tikzpicture}
\end{center}
and similarly, the same compatability holds for vertical composition of squares and the associator of $V.$
\item \textbf{VDC10}: \  The pseudo-identity squares for horizontal and vertical composition are compatible with the unitors of $H$ and $V$ respectively. That is, the pseudo-identity squares they act as identities up to the appropriate unitors. For $\Theta \in \mbox{Sq}(f,f',p,p'),$ the condition is $$(\Theta\boxvert \ID_{p})\uact \rho_{f}\dact \rho_{f'}=\Theta= (\ID_{p'}\boxvert\Theta)\uact \lambda_f\dact \lambda_{f'}:$$

\begin{center}
\begin{tikzpicture}[scale=1.8,auto]

\begin{scope}

\node (01) at (0,0) {$a'$};
\node (00) at (0,1) {$a$};
\node (20) at (2,1) {$b$};
\node (21) at (2,0) {$b'$};
\node[rotate=-45] at (1,.4){$\Rightarrow$};
\node[scale=.8] at (1.1,.55){$\Theta\boxvert\ID_{p}$};

\path[->] (00) edge node[midway,below]{\scriptsize $f\circ \id_a$}(20);
\path[->] (00) edge node[midway,left]{$p$}(01);
\path[->] (01) edge node[midway,above]{\scriptsize $f'\circ \id_{a'}$}(21);
\path[->] (20) edge node[midway]{$p'$}(21);

\draw[->] (01) to [out=-45,in=-135] node(f)[midway,below]{ $f'$} (21) ;

\node[rotate=+90] (areta) at ($(f)+(0,+.45)$){$\Rightarrow$};
\node (eta) at ($(areta)+(.22,0)$){ $\rho_{f'}$};

\draw[->] (00) to [out=45,in=135] node(f)[midway]{$ f$} (20) ;

\node[rotate=-90] (areta) at ($(f)+(0,-.42)$){$\Rightarrow$};
\node (eta) at ($(areta)+(.22,0)$){ $\rho_f$};

\node(eq1) at (2.67,.5){$=$};

\end{scope}

\begin{scope}[shift={(3.3,0)}]
\node (10') at (1,1) {$b$};
\node (00') at (0,1) {$a$};
\node (11') at (1,0) {$b'$};
\node (01') at (0,0) {$a'$};
\path[->] (00') edge node[midway]{$f$}(10');
\path[->] (00') edge node[midway,swap]{$p$}(01');
\path[->] (01') edge node[midway,swap]{$ f'$}(11');
\path[->] (10') edge node[midway]{$p'$}(11');
\node[rotate=-45] at (.5,.43){$\Rightarrow$};
\node[scale=.8] at (.55,.63){$\Theta$};

\node(eq1) at (1.70,.5){$=$};

\end{scope}

\begin{scope}[shift={(5.7,0)}]
\node (01) at (0,0) {$a'$};
\node (00) at (0,1) {$a$};
\node (20) at (2,1) {$b$};
\node (21) at (2,0) {$b'$};
\node[rotate=-45] at (1,.38){$\Rightarrow$};
\node[scale=.8] at (1.1,.59){$\ID_{p'}\boxvert\Theta$};

\path[->] (00) edge node[midway,below]{\scriptsize $ \id_b\circ f$}(20);
\path[->] (00) edge node[midway,left]{$p$}(01);
\path[->] (01) edge node[midway,above]{\scriptsize $ \id_{b'}\circ f'$}(21);
\path[->] (20) edge node[midway]{$p'$}(21);

\draw[->] (01) to [out=-45,in=-135] node(f)[midway,below]{ $f'$} (21) ;

\node[rotate=+90] (areta) at ($(f)+(0,+.45)$){$\Rightarrow$};
\node (eta) at ($(areta)+(.22,0)$){ $\lambda_{f'}$};

\draw[->] (00) to [out=45,in=135] node(f)[midway]{$ f$} (20) ;

\node[rotate=-90] (areta) at ($(f)+(0,-.42)$){$\Rightarrow$};
\node (eta) at ($(areta)+(.22,0)$){ $\lambda_f$};
\end{scope}
\end{tikzpicture}
\end{center}
\item \textbf{VDC11}: \ Horizontal and vertical composition of squares satisfy the interchange law for the two ways of composing a 2-by-2 grid, $$(\Pi'\boxvert\Theta')\boxminus(\Pi \boxvert \Theta)=(\Pi'\boxminus\Pi)\boxvert(\Theta' \boxminus \Theta),$$ for $\Theta,\Pi,\Theta',\Pi'$ as shown:

\begin{center}
\begin{tikzpicture}[scale=1.8,auto]

\begin{scope}

\node (10) at (1,1) {$b$};
\node (00) at (0,1) {$a$};
\node (11) at (1,0) {$b'$};
\node (01) at (0,0) {$a'$};
\node[rotate=-45] at (.5,.5){$\Rightarrow$};
\node[scale=.8] at (.6,.65){$\Theta$};

\node (20) at (2,1) {$c$};
\node (21) at (2,0) {$c'$};
\node(22) at (2,-1){$c''$};
\node[rotate=-45] at (1.5,.5){$\Rightarrow$};
\node[scale=.8] at (1.6,.65){$\Pi$};

\node (12) at (1,-1) {$b''$};
\path[->] (00) edge node[midway]{$f$}(10);
\path[->] (00) edge node[midway,swap]{$p$}(01);
\path[->] (01) edge node[midway,swap]{$f'$}(11);
\path[->] (10) edge node[midway]{$p'$}(11);
\path[->] (10) edge node[midway]{$g$}(20);
\path[->] (11) edge node[midway,swap]{$g'$}(21);
\path[->] (20) edge node[midway]{$p''$}(21);
\path[->](21) edge node[midway]{$q''$}(22);
\path[->](12) edge node[midway,below]{$g''$}(22);

\node (01') at (0,-1) {$a''$};
\node[rotate=-45] at (0.5,-.58){$\Rightarrow$};
\node[scale=.8] at (0.6,-.43){$\Theta'$};

\path[->] (11) edge node[midway]{$q'$}(12);
\path[->] (01) edge node[midway,swap]{$q$}(01');
\path[->] (01') edge node[midway,swap]{$f''$}(12);
\node[rotate=-45] at (1.5,-.58){$\Rightarrow$};
\node[scale=.8] at (1.6,-.43){$\Pi'$};

\end{scope}
\end{tikzpicture}
\end{center}
\end{itemize}


\section{The Verity double category $\mbox{Vdc}(X)$.}


\begin{definition} An \emph{algebraic} inner-Kan bisimplicial set is an inner-Kan bisimplicial set $X$ equipped with a set $\chi$ of preferred fillers for its inner horns. If $X$ is $2$-reduced inner-Kan, the only case in which this choice is non-trivial is for horns having the form  $\Lambda[2^{\, 1}, 0]$ and $\Lambda[0,2^{\, 1}]$, so the data of $\chi$ is equivalent to giving an algebraic structures $\chi_h$ and $\chi_v$ to both $X_{0\bullet}$ and $X_{\bullet 0}.$  A map of algebraic inner-Kan bisimplicial sets is called \emph{strict} if it preserves the algebraic structure $\chi.$     \end{definition}

Let $(X,\chi)$ be a $2$-reduced inner-Kan bisimplicial set. Let $X_h:=X_{0 \bullet}$ and $X_v:=X_{\bullet0}$ denote the simplicial sets that make up the $0$th row and $0$th column of $X$. Note that the $(m,0)$ Kan conditions for $X$ is equivalent to the $m$th Kan condition for $X_h$, and the $(0,n)$ Kan conditions for $X$ is equivalent to the $n$th Kan condition for $X_v$, so $X_h$ and $X_v$, so $X_h$ and $X_v$ satisfy the inner Kan condition for all $n$, and satisfy the inner Kan condition uniquely for $n>2,$ i.e. they are $2$-reduced inner-Kan.

In this section, we will define from $X$ the structures of a Verity double category $\mbox{Vdc}(X)$, then, in the next section, we will check that these structures meet the axioms for a VDC.
\begin{definition} The horizontal $(2,1)$-category $H$ of $\mbox{Vdc}(X)$ is defined to be $\Bic(X_h,\chi_h)$, as defined in Chapter~\ref{bicchapter}. The vertical $(2,1)$-category $V$ is $\Bic(X_v,\chi_v)$.
\end{definition}
\begin{definition} The squares of $\mbox{Vdc}(X)$ are defined to be the elements of $X_{11}$. Let $\Theta \in X_{11}$ and let $d\Theta=[g',\ g\ \ | \ \ q',\ q].$ Note that $g,g'\in X_{01}$ so they are $1$-morphisms in $H$, and  $q,q'\in X_{10}$ are $1$-morphisms in $V$.  So we define this element $\Theta$ to be a square in $\mbox{Sq}(g',g,q',q).$ 
\end{definition}

\begin{definition}[Horizontal composition of squares]\label{horzcompdef} Let $\Theta$, $\Pi$ be two squares as shown:

\begin{center}
\begin{tikzpicture}[scale=1.8,auto]

\begin{scope}

\node (10) at (1,1) {$b$};
\node (00) at (0,1) {$a$};
\node (11) at (1,0) {$b'$};
\node (01) at (0,0) {$a'$};
\node[rotate=-45] at (.5,.5){$\Rightarrow$};
\node[scale=.8] at (.6,.65){$\Theta$};

\node (20) at (2,1) {$c$};
\node (21) at (2,0) {$c'$};
\node[rotate=-45] at (1.5,.5){$\Rightarrow$};
\node[scale=.8] at (1.6,.65){$\Pi$};

\path[->] (00) edge node[midway]{$g$}(10);
\path[->] (00) edge node[midway,swap]{$q$}(01);
\path[->] (01) edge node[midway,swap]{$g'$}(11);
\path[->] (10) edge node[midway]{$q'$}(11);
\path[->] (10) edge node[midway]{$h$}(20);
\path[->] (11) edge node[midway,swap]{$h'$}(21);
\path[->] (20) edge node[midway]{$q''$}(21);
\end{scope}
\end{tikzpicture}
\end{center}
then the horizontal composition $\Pi\boxvert\Theta$ in  $\mbox{Vdc}(X)$ is defined by the following horn:

\begin{table}[H]\caption*{$\Lambda_{\boxvert}(\Pi,\Theta),$ defining $\Pi\boxvert \Theta$}\begin{center}

\begin{tabular}{ r | l || !{\vrule width 2pt}  l | l | l |  }
    \cline{2-5}

   &  $\chi(h',g')$    &     $h'$    &    $h'\circ g'$   &  $g'$          \\ \cline{2-5} 
     
 &  $\chi(h,g)$ &     $h $    &    $h \circ g$   &  $g$         \\ \cline{2-5} 
 
  \end{tabular}\vspace{12pt}

    \begin{tabular}{ r |   l ||  l | l !{\vrule width 2pt} l | l |}
    \cline{2-6}

       &  $\Pi$     &  $h'$  &                        $h$                  &     $q'' $    &    $q'$     \\ \cline{2-6} 
     
 $ \Lambda$       &      $=:\Pi\boxvert\Theta$       &    $h'\circ g'$     &     $h \circ g$     &     $q'' $    &    $q$         \\ \cline{2-6} 
 
    &  $\Theta$ &              $g'$  &                       $g$                &     $q'$    &    $q$      \\ \cline{2-6}  \end{tabular}

\end{center}
\end{table}

\end{definition}
\begin{definition}[Vertical composition of squares] Let $\Theta$, $\Pi$ be two squares as shown:

\begin{center}
\begin{tikzpicture}[scale=1.8,auto]

\begin{scope}

\node (10) at (1,1) {$b$};
\node (00) at (0,1) {$a$};
\node (11) at (1,0) {$b'$};
\node (01) at (0,0) {$a'$};
\node[rotate=-45] at (.5,.5){$\Rightarrow$};
\node[scale=.8] at (.6,.65){$\Theta$};

\node (11') at (1,-1) {$b''$};
\node (01') at (0,-1) {$a''$};
\node[rotate=-45] at (0.5,-.58){$\Rightarrow$};
\node[scale=.8] at (0.6,-.43){$\Pi$};

\path[->] (00) edge node[midway]{$g$}(10);
\path[->] (00) edge node[midway,swap]{$q$}(01);
\path[->] (01) edge node[midway,swap]{$g'$}(11);
\path[->] (10) edge node[midway]{$q'$}(11);
\path[->] (11) edge node[midway]{$r'$}(11');
\path[->] (01) edge node[midway,swap]{$r$}(01');
\path[->] (01') edge node[midway,swap]{$g''$}(11');
\end{scope}
\end{tikzpicture}
\end{center}
then the vertical composition $\Pi\boxminus \Theta$ in $\mbox{Vdc}(X)$ is defined by the following horn:

\begin{table}[H]\caption*{$\Lambda_{\boxminus}(\Pi,\Theta),$ defining $\Pi\boxminus \Theta$}\begin{center}
   \begin{tabular}{ r | l || l | l !{\vrule width 2pt} l | l  |}
    \cline{2-6}

        &  $\Pi$  &  $g''$  &   $g'$   &     $r' $    &    $r$       \\ \cline{2-6} 
     
 $\Lambda$       &  $=:\Pi\boxminus \Theta$  &  $g''$  &   $g$   &   $r'\circ q' $    &    $r\circ q$     \\ \cline{2-6} 
 
             &  $\Theta$  &  $g'$  &   $g$   &     $q'$    &    $q$      \\ \cline{2-6}  \end{tabular}
\vspace{12pt}

\begin{tabular}{ r | l || l | l | l !{\vrule width 2pt} }
    \cline{2-5}

       &  $\chi(r',q')$   &     $r' $    &    $r'\circ q'$   &  $q' $       \\ \cline{2-5} 
     
        &  $\chi(r,q)$     &     $r $    &    $r\circ q$   &  $q$      \\ \cline{2-5} 
 
  \end{tabular}
  \end{center}  \end{table}

\end{definition}

\begin{definition}[Pseudo-identity squares] Let $f$ be a horizontal $1$-morphism of $\mbox{Vdc}(X),$ which is an element of $X_{01}.$ Then we define the pseudo-identity square for vertical composition $\ID_f:= s_0 f\in X_{11}.$ Similarly, for a vertical $1$-morphism $p$ in $X_{10}$ we define $\ID_p:= \varsigma_0  p\in X_{11}$.
\end{definition}
\begin{definition}[Top action]\label{uactdef} Let $\Theta$ be a square of $\mbox{Vdc}$ and $\beta$ be a horizontal $2$-morphism as shown:

\begin{center}
\begin{tikzpicture}[scale=1.8,auto]

\begin{scope}

\node (10) at (1,1) {$b$};
\node (00) at (0,1) {$a$};
\node (11) at (1,0) {$b'$};
\node (01) at (0,0) {$a'$};
\node[rotate=-45] at (.5,.5){$\Rightarrow$};
\node[scale=.8] at (.6,.65){$\Theta$};

\path[->] (00) edge node (g)[midway,below]{$g$}(10);
\draw[->] (00) to [out=60,in=120] node(f)[midway]{$f$} (10) ;
\path[->] (00) edge node[midway,swap]{$q$}(01);
\path[->] (01) edge node[midway,swap]{$g'$}(11);
\path[->] (10) edge node[midway]{$q'$}(11);
\node[rotate=-90] (areta) at ($(f)+(0,-.38)$){$\Rightarrow$};
\node (eta) at ($(areta)+(.15,0)$){$\beta$};

\end{scope}
\end{tikzpicture}
\end{center}
then the top action $\Theta \uact \beta$ is defined by the following horn:

\begin{table}[H]\caption*{$\Lambda_{\uact}(\Theta,\beta),$ defining $\Theta \uact \beta$}\begin{center}

\begin{tabular}{ r | l || !{\vrule width 2pt}  l | l | l |  }
    \cline{2-5}

  &  $\Id_{g'}$    &  $g'$    &     $g'$    &    $\id_{a'}$          \\ \cline{2-5} 
     
    &  $\beta$ &     $g $    &    $f$   &  $\id_{a}$         \\ \cline{2-5} 
 
  \end{tabular}\vspace{12pt}

\begin{tabular}{ r | l ||  l | l !{\vrule width 2pt} l | l |}
    \cline{2-6}

          &  $\Theta$ &  $g'$  &   $g$   &     $q' $    &    $q$       \\ \cline{2-6} 
     
 $ \Lambda$      &      $=:\Theta \uact \beta$       &    $g'$     &     $f$     &     $q' $    &    $q$          \\ \cline{2-6} 
 
            &  $\ID_q$  &              $\id_{a'}$  &   $\id_a$   &                                              $q$    &    $q$      \\ \cline{2-6}  \end{tabular}

\end{center}
\end{table}

\end{definition}
\begin{definition}[Bottom action]Let $\Theta$ be a square of $\mbox{Vdc}$ and $\beta'$ be a horizontal $2$-morphism as shown:
\item Action of horizontal 2-morphisms on the bottom, $ \Theta \filledmedtriangledown \beta' $:

\begin{center}
\begin{tikzpicture}[scale=1.8,auto]

\begin{scope}

\node (10) at (1,1) {$b$};
\node (00) at (0,1) {$a$};
\node (11) at (1,0) {$b'$};
\node (01) at (0,0) {$a'$};
\node[rotate=-45] at (.5,.5){$\Rightarrow$};
\node[scale=.8] at (.6,.65){$\Theta$};

\path[->] (00) edge node (g)[midway]{$g$}(10);
\draw[->] (01) to [out=-60,in=-120] node(f)[midway,below]{$f'$} (11) ;
\path[->] (00) edge node[midway,swap]{$q$}(01);
\path[->] (01) edge node[midway]{$g'$}(11);
\path[->] (10) edge node[midway]{$q'$}(11);
\node[rotate=90] (areta) at ($(f)+(0,+.4)$){$\Rightarrow$};
\node (eta) at ($(areta)+(.15,0)$){$\beta'$};
\end{scope}
\end{tikzpicture}
\end{center}
then the bottom action $\Theta \dact \beta'$ is defined by the following horn:

\begin{table}[H]\caption*{$\Lambda_{\dact}(\Theta,\beta'),$ defining $\Theta \dact \beta'$}\begin{center}
 \subfloat

\begin{tabular}{ r | l || !{\vrule width 2pt}  l | l | l |  }
    \cline{2-5}

  & $\beta'$      &  $g'$    &     $f'$    &    $\id_{a'}$          \\ \cline{2-5} 
     
    & $\Id_{g}$ &     $g $    &    $g$   &  $\id_{a}$         \\ \cline{2-5} 
 
  \end{tabular}\vspace{12pt}

    \begin{tabular}{ r |   l ||  l | l !{\vrule width 2pt} l | l |}
    \cline{2-6}

          &  $\Theta$ &  $g'$  &   $g$   &     $q' $    &    $q$      \\ \cline{2-6} 
     
 $ \Lambda$     &      $=:\Theta \dact \beta'$      &    $f'$     &     $g$     &     $q' $    &    $q$       \\ \cline{2-6} 
 
               &  $\ID_q$     &               $\id_{a'}$  &   $\id_a$   &                                              $q$    &    $q$     \\ \cline{2-6}  \end{tabular}

\end{center}
\end{table}

\end{definition}
\begin{definition}[Left action] Let $\Theta$ be a square of $\mbox{Vdc}$ and $\eta$ be a vertical $2$-morphism as shown: \label{lactdef}
\begin{center}
\begin{tikzpicture}[scale=1.8,auto]

\begin{scope}
\node (10) at (1,1) {$b$};
\node (00) at (0,1) {$a$};
\node (11) at (1,0) {$b'$};
\node (01) at (0,0) {$a'$};
\node[rotate=-45] at (.5,.5){$\Rightarrow$};
\node[scale=.8] at (.6,.65){$\Theta$};

\path[->] (00) edge node (g)[midway]{$g$}(10);

\path[->] (00) edge node[midway]{$q$}(01);
\path[->] (01) edge node[midway,below]{$g'$}(11);
\path[->] (10) edge node[midway]{$q'$}(11);
\draw[->] (00) to [out=-150,in=150] node(f)[midway,left]{$p$} (01) ;
\node[rotate=0] (areta) at ($(f)+(.35,0)$){$\Rightarrow$};
\node (eta) at ($(areta)+(0,.15)$){$\eta$};
\end{scope}
\end{tikzpicture}
\end{center}
then the left action $\Theta \lact \eta$ is defined by the following horn:

\begin{table}[H]\caption*{$\Lambda_{\lact}(\Theta,\eta),$ defining $\Theta \lact \eta$}\begin{center}
\begin{tabular}{ r | l || l | l !{\vrule width 2pt} l | l  |}
    \cline{2-6}

                &  $\Theta$                  &  $g'$  &   $g$   &        $q' $    &    $q$       \\ \cline{2-6} 
     
 $\Lambda$       &  $=:\Theta \lact \eta$    &  $g'$  &   $g$   &        $q'$    &    $p$     \\ \cline{2-6} 
 
                 &  $\ID_g$                 &   $g$  &   $g$   &       $\id_b$    &    $\id_a$      \\ \cline{2-6}  \end{tabular}\vspace{12pt}
                 
\begin{tabular}{ r | l || l | l | l !{\vrule width 2pt} }
    \cline{2-5}

        &  $\Id_{q'}$  &     $q' $    &    $q'$   &  $\id_b$       \\ \cline{2-5} 
     
          &  $\eta$    &     $q $    &    $p$   &  $\id_a$        \\ \cline{2-5} 
 
  \end{tabular}
 \end{center}   \end{table}

\end{definition}
\begin{definition}[Right action]\label{ractdef}Let $\Theta$ be a square of $\mbox{Vdc}$ and $\eta'$ be a vertical $2$-morphism as shown:

\begin{center}
\begin{tikzpicture}[scale=1.8,auto]
\begin{scope}
\node (10) at (1,1) {$b$};
\node (00) at (0,1) {$a$};
\node (11) at (1,0) {$b'$};
\node (01) at (0,0) {$a'$};
\node[rotate=-45] at (.5,.5){$\Rightarrow$};
\node[scale=.8] at (.6,.65){$\Theta$};

\path[->] (00) edge node (g)[midway]{$g$}(10);

\path[->] (00) edge node[midway,left]{$q$}(01);
\path[->] (01) edge node[midway,below]{$g'$}(11);
\path[->] (10) edge node[midway,left]{$q'$}(11);

\draw[->] (10) to [out=-30,in=30] node(p')[midway,right]{$p'$} (11) ;
\node[rotate=180] (aretar) at ($(p')+(-.35,0)$){$\Rightarrow$};
\node (eta) at ($(aretar)+(0,.15)$){$\eta'$};

\end{scope}
\end{tikzpicture}
\end{center}then the right action $\Theta \ract \eta'$ is defined by the following horn:

\begin{table}[H]\begin{center}\caption*{$\Lambda_{\ract}(\Theta,\eta'),$ defining $\Theta \ract \eta'$}
   \begin{tabular}{ r | l || l | l !{\vrule width 2pt} l | l  |}
    \cline{2-6}

                &  $\Theta$                  &  $g'$  &   $g$   &        $q' $    &    $q$       \\ \cline{2-6} 
     
 $\Lambda$       &  $=:\Theta \ract \eta'$    &  $g'$  &   $g$   &        $p'$    &    $q$     \\ \cline{2-6} 
 
                 &  $\ID_g$                 &   $g$  &   $g$   &       $\id_b$    &    $\id_a$      \\ \cline{2-6}  \end{tabular}
\vspace{12pt}

\begin{tabular}{ r | l || l | l | l !{\vrule width 2pt} }
    \cline{2-5}

         &  $\eta'$ &     $q' $    &    $p'$   &  $\id_b$       \\ \cline{2-5} 
     
         &  $\Id_q$   &     $q $    &    $q$   &  $\id_a$       \\ \cline{2-5} 
 
  \end{tabular}
 \end{center}   \end{table}

\end{definition}

\section{Verification of the VDC axioms for $\mbox{Vdc}(X).$}
For the this section we will always have $\Theta,\beta,\beta',\eta,\eta'$ as in Definitions~\ref{uactdef}--\ref{ractdef}, as shown in Figure~\ref{Pipic}.

\begin{figure}[h,b,t]\begin{center}
\begin{tikzpicture}[scale=1.8,auto]
\begin{scope}
\node (10) at (1,1) {$b$};
\node (00) at (0,1) {$a$};
\node (11) at (1,0) {$b'$};
\node (01) at (0,0) {$a'$};
\node[rotate=-45] at (.5,.5){$\Rightarrow$};
\node[scale=.8] at (.6,.65){$\Theta$};

\path[->] (00) edge node (g)[midway,below]{$g$}(10);

\path[->] (00) edge node[midway,right]{$q$}(01);
\path[->] (01) edge node[midway,above]{$g'$}(11);
\path[->] (10) edge node[midway,left]{$q'$}(11);

\draw[->] (10) to [out=-30,in=30] node(p')[midway,right]{$p'$} (11) ;
\node[rotate=180] (areta) at ($(p')+(-.35,0)$){$\Rightarrow$};
\node (eta) at ($(areta)+(0,.15)$){$\eta'$};

\draw[->] (00) to [out=-150,in=150] node(p)[midway,left]{$p$} (01) ;
\node[rotate=0] (aretal) at ($(p)+(.35,0)$){$\Rightarrow$};
\node (eta) at ($(aretal)+(0,.15)$){$\eta$};

\draw[->] (01) to [out=-60,in=-120] node(f')[midway,below]{$f'$} (11) ;
\node[rotate=90] (aretab) at ($(f')+(0,+.4)$){$\Rightarrow$};
\node (eta) at ($(aretab)+(.15,0)$){$\beta'$};

\draw[->] (00) to [out=60,in=120] node(f)[midway]{$f$} (10) ;
\node[rotate=-90] (aretau) at ($(f)+(0,-.38)$){$\Rightarrow$};
\node (eta) at ($(aretau)+(.15,0)$){$\beta$};

\end{scope}
\end{tikzpicture}\caption{~\label{Pipic}}\end{center}
\end{figure}
\subsection{Preliminaries}
Let $X$ be a $2$-reduces inner-Kan bisimplicial set. The following analogue of Lemma~\ref{lemma1} is the main tool we use for checking the Verity double category axioms for $\Vdc(X):$
\begin{lemma}[Matching Lemma] Suppose we have two commutative $d\Delta[n,m]$-spheres $S$, $S'$ in $X$, such that every corresponding face, except for possibly a single inner face, of these spheres match. Then $S=S'$ and in particular they indeed match on this putatively non-matching inner face.
\end{lemma}
\begin{proof} By the uniqueness of fillers for inner horns in $X$.
\end{proof}
 In the previous section, we defined the data of a (putative) Verity double category $\mbox{Vdc}(X)$. In this section, we check that the Verity double axioms hold for $\mbox{Vdc}(X).$ Recall that in the verification the bicategory axioms for $\Bic(X)$ in Chapter~\ref{bicchapter}, even axioms that were in a sense ``symmetric'' to each other, for instance the naturality of the associator in its first vs. its third argument, the proofs of the two cases were not neatly symmetrical. Instead, because of the asymmetry inherent in the definition of the $2$-morphisms in $\Bic(X)$, one case was usually more complex, with more alt 2-morphisms and ``hats'' involved. However, in the verification of the VDC axioms for $\mbox{Vdc}(X)$, there is an obvious symmetry between the horizontal and vertical directions for both the axioms and our definitions. Moreover there is a close symmetry between the definition of $\uact$ and $\dact$ (and also between $\lact$ and $\ract$) that will make the proofs for the axioms that only involve one kind of $2$-morphism action (action uniticity, action associativity, compatibility with square composition, and interchange for actions) exactly mirror each other in these cases.

Still, we cannot completely avoid dealing with alt $2$-morphisms. We defined actions on each side of a square by a $2$-morphisms, but it equally possible to define the actions for alt $2$-morphisms.
\begin{definition} Let $\Theta \in \mbox{Sq}(g',g, q',q)$, and $\widehat{\beta}$ be an alt $2$-morphism $f \Rightarrow g'$ (see Definition~\ref{altmorphisms}). Then we define $\Theta \uact \wh{\beta}$ by the following horn:
\begin{table}[H]\caption*{$\Lambda_{\uact}(\Theta,\wh{\beta}),$ defining $\Theta \uact \wh{\beta}$}\begin{center}

\begin{tabular}{ r | l || !{\vrule width 2pt}  l | l | l |  }
    \cline{2-5}

  &  $\wh{\Id_{g'}}$    &  $\id_{b'}$    &     $g'$    &   $g'$           \\ \cline{2-5} 
     
    &  $\wh{\beta}$ &   $\id_{b}$  &  $f$    &    $g$            \\ \cline{2-5} 
 
  \end{tabular}\vspace{12pt}

\begin{tabular}{ r |   l ||  l | l !{\vrule width 2pt} l | l |}
    \cline{2-6}

                    &  $\ID_{q'}$    &  $\id_{b'}$   &  $\id_b$  &     $q' $    &    $q'$      \\ \cline{2-6} 
     
 $ \Lambda$        &      $=:\Theta \uact \wh{\beta}$     &    $g'$       &     $f$     &     $q' $    &    $q$          \\ \cline{2-6} 
 
                      &    $\Theta$   &   $g'$     &   $g$   &     $q'$    &    $q$     \\ \cline{2-6}  \end{tabular}

\end{center} 
\end{table}
We define the other actions for alt 2-morphisms similarly by horns as follows, with $\Theta,$ $\beta'$, $\eta$, and $\eta'$ defined as in Figure~\ref{Pipic}:
\begin{align*}
\Delta_{\dact}(\Theta,\wh{\beta'})&= [\wh{\beta'},\ \wh{\Id_g} \ \ | \ \      \ID_{q'}, \   =:\Theta \dact \wh{\beta},\ \Theta                                    ]     \\ 
\Delta_{\lact}(\Theta,\wh{\eta})&=  [ \ID_{g'}   , \   =:\Theta \lact \wh{\eta},\ \Theta  \ \ | \ \ \wh{\Id_{q'}} ,\      \wh{\eta}                                   ]   \\
\Delta_{\ract}(\Theta,\wh{\eta'})&= [ \ID_{g'},\        =:\Theta \ract \wh{\eta},\  \Theta  \ \ | \ \  \wh{\eta'}   , \  \wh{\Id_{q}}                                  ]   
\end{align*} 
\end{definition}
\begin{lemma} \label{altaction} Let $\Theta \in \mbox{Sq}(g',g, q',q)$, as above and $\beta$ be an $2$-morphism $f \Rightarrow g'.$ Then $\Theta \uact \beta=\Theta\uact\wh{\beta}.$
\end{lemma}
\begin{proof} By the Matching Lemma, it suffices to show the following sphere is commutative:
$$[  \wh{\Id_{g'}},\ \wh{\beta}           \ \  |  \ \   \ID_{q'},\ \Theta\uact\beta,\ \Theta         ].$$ The following horn verifies this commutativity:
\begin{table}[H]\begin{center}
    \begin{tabular}{ r | l ||!{\vrule width 2pt} l | l | l | l |}
    \cline{2-6} 
        &        $\varsigma_2(\Id_{g'})$       &$\wh{\Id_{g'}}$ & $\wh{\Id_{g'}}$    &  $\Id_{g'}$   &  $\Id_{g'}$          \\ \cline{2-6}
     
       &  $\Delta_{\wedge}(\beta)$    &    $\wh{\Id_g}$  &  $\wh{\beta}$    & $\beta$       &     $\Id_{g}$           \\ \cline{2-6}
                   
    \end{tabular} 

\vspace{12pt}

 \begin{tabular}{    r | l || l | l  !{\vrule width 2pt}  l | l | l | }\cline{2-7} 
              &  $\varsigma_1(\Theta)$ &   $\wh{\Id_{g'}}$      &      $\wh{\Id_g}$             &     $\ID_{q'}$            &    $\Theta$                      & $\Theta$           \\ \cline{2-7} 
      $\Lambda$    &   & $\wh{\Id_{g'}}$   &     $\wh{\beta}$              &     $\ID_{q'}$            &   $\Theta \uact \beta$           &  $\Theta$        \\ \cline{2-7}
                        &  $\Delta_{\uact}(\Theta,\beta)$  &   $\Id_{g'}$  &     $\beta$                   &     $\Theta$              &   $\Theta \uact \beta$           &  $\ID_q $          \\ \cline{2-7}
                    &  $\varsigma_0(\Theta)  $    &  $\Id_{g'}$     &     $\Id_{g}$             &     $\Theta$              &       $\Theta$                   &  $\ID_q $         \\ \cline{2-7}
\end{tabular}   
 \end{center}
    \end{table}
\end{proof}
\subsection{Contravariant action of $2$-morphisms}
\begin{proposition}[Uniticity of actions (\textbf{VDC1})] Let $\Theta\in \mbox{Sq}(g',g,q',q)$ be in $\mbox{Vdc}(X),$ as shown in Figure~\ref{Pipic}. Then $\Theta\uact\Id_g=\Theta\dact\Id_g'=\Theta$, and $\Theta\lact\Id_q=\Theta\ract\Id_q'=\Theta.$
\end{proposition}
\begin{proof} We have by definition $$d\left(\Delta_{\uact}(\Theta,\Id_g)\right)=d\left(\Delta_{\dact}(\Theta,\Id_{g'})\right)= [ \Id_{g'},\ \Id_g,\ \ | \ \ \Theta,\ \Theta\uact \Id_g=\Theta \dact \Id_{g'},\ \ID_q ], $$ whereas we compute 
\begin{align*}
d(\varsigma_0 \Theta)&= [ d_0 \varsigma_0 \Theta,\  d_1 \varsigma_0 \Theta,\ \ | \ \  \delta_0 \varsigma_0 \Theta,\ \delta_1 \varsigma_0 \Theta,\ \delta_2 \varsigma_0 \Theta ] \\
                         &= [ \varsigma_0  d_0\Theta,\   \varsigma_0 d_1\Theta,\ \ | \ \  \Theta,\  \Theta,\ \varsigma_0\delta_1  \Theta ] \\
                         &= [  \varsigma_0 g',\   \varsigma_0 g,\ \ | \ \  \Theta,\  \Theta,\ \varsigma_0 q ] \\
                         &= [  \Id_{g'},\   \Id_g,\ \ | \ \  \Theta,\  \Theta,\ \ID_q ]. 
\end{align*}
Thus, by the Matching Lemma,   $\Theta\uact\Id_g=\Theta\dact\Id_g'=\Theta$. By symmetry, $\Theta\lact\Id_q=\Theta\ract\Id_q'=\Theta.$
\end{proof}
\begin{proposition}[Associativity of actions (\textbf{VDC2})] Let  $\Theta\in \mbox{Sq}(g',g,q',q)$, and horizontal $2$-morphisms \mbox{$e\stackrel{\alpha}{\Rightarrow}f \stackrel{\beta}{\Rightarrow}g$} be in $\mbox{Vdc}(X).$ Then $(\Theta\uact \beta)\uact \alpha=\Theta \uact (\beta \bullet\alpha).$ A similar axiom holds for $\dact$, $\lact$, and~$\ract$.
\end{proposition}
\begin{proof} $$d \left( \Delta_{\uact}( \Theta \uact \beta, \alpha)\right)= [ \Id_{g'},\ \alpha,\ \ | \ \ \Theta\uact \beta,\ (\Theta\uact \beta) \uact \alpha,\ \ID_q ].$$ Thus by the Matching Lemma it suffices to show $$[ \Id_{g'},\ \alpha,\ \ | \ \ \Theta\uact \beta,\ \Theta\uact (\beta \bullet \alpha),\ \ID_q ]$$ is commutative. The following horn of type $\Lambda[1,3^{\ 2} ]$ verifies this commutativity:
\begin{table}[H]\begin{center}
    \begin{tabular}{ r | l ||!{\vrule width 2pt} l | l | l | l |}
    \cline{2-6} 
        &  $\varsigma_1(\Id_{g'})$  & $\Id_{g'}$     &    $\Id_{g'}$   &  $\Id_{g'}$     &   $\Id_{\id_{a'}}$         \\ \cline{2-6}
     
        &  $\Delta_{\bullet}(\beta, \alpha )$  &  $\beta$  & $\beta\bullet \alpha$   &   $\alpha$      &  $\Id_{\id_a}$        \\ \cline{2-6}
                   
    \end{tabular} 
\vspace{12pt}

 \begin{tabular}{    r | l || l | l  !{\vrule width 2pt}  l | l | l |}\cline{2-7} 
               &  $\Delta_{\uact}(\Theta,\beta)$ &   $\Id_{g'}$     &      $\beta$                 &     $\Theta$               &    $\Theta \uact \beta$                      &  $\ID_q $        \\ \cline{2-7} 
              &  $\Delta_{\uact}(\Theta,\beta\bullet\alpha)$ &   $\Id_{g'}$     &     $\beta\bullet \alpha$   &     $\Theta$                &   $\Theta \uact (\beta\bullet\alpha)$       &  $\ID_q $           \\ \cline{2-7}
 $\Lambda$      &   &   $\Id_{g'}$     &     $\alpha$                &     $\Theta\uact \beta $    &    $\Theta \uact (\beta\bullet\alpha)$      &  $\ID_q $       \\ \cline{2-7}
                 &  $\varsigma_0(\ID_q)$  &   $\Id_{\id_{a'}}$  &     $\Id_{\id_a}$             &     $\ID_q $    &        $\ID_q $                                   &  $\ID_q $               \\ \cline{2-7}
\end{tabular}   
 \end{center}
    \end{table}

The associativity action for $\dact$ is similar, and the axioms for $\lact$ and $\ract$ follow by symmetry. 
\end{proof}
\begin{proposition}[Commutativity of actions (\textbf{VDC3})] \label{commute} The actions $\uact$ $\dact$, $\lact$, and $\ract$ commute with each other.
\end{proposition}
\begin{proof} Let $\Theta\in \mbox{Sq}(g',g,q',q)$ and $\beta,\beta',\eta,\eta'$ as in Figure~\ref{Pipic}. We first show $$(\Theta \dact \beta')\uact \beta= (\Theta \uact \beta)\dact \beta'.$$ We have $$d \left( \Delta_{\uact}(\Theta\dact \beta',\beta )\right)= [ \Id_{f'},\ \beta,\ \ | \ \ \Theta \dact \beta',\   (\Theta \dact \beta') \uact \beta,\  \ID_q]$$ so by the Matching Lemma it is enough to show  $$[ \Id_{f'},\ \beta,\ \ | \ \ \Theta \dact \beta',\  (\Theta \uact \beta)\dact \beta ,\  \ID_q]$$ is commutative. The following  table proof verifies this commutativity:
\begin{table}[H]\begin{center}\caption{~}
    \begin{tabular}{ r | l ||!{\vrule width 2pt} l | l | l | l |}
    \cline{2-6} 
        &  $\varsigma_0(\beta')$ & $\beta'$     &    $\beta'$   &  $\Id_{f'}$     &   $\Id_{\id_{a'}}$         \\ \cline{2-6}
     
        &  $\varsigma_1(\beta)$ &  $\Id_g$  & $\beta$   &   $\beta$      &  $\Id_{\id_a}$        \\ \cline{2-6}
                   
    \end{tabular} 

\vspace{12pt}
 \begin{tabular}{    r | l || l | l  !{\vrule width 2pt}  l | l | l |}\cline{2-7} 
               &  $\Delta_{\dact}(\Theta,\beta')$ &   $\beta'$     &      $\Id_g$                 &     $\Theta$               &    $\Theta \dact \beta'$                      &  $\ID_q $        \\ \cline{2-7} 
       $\odot$   & (Table~\ref{com1pt2})  &   $\beta'$     &     $\beta$                  &     $\Theta$                &   $(\Theta \uact \beta)\dact\beta'$       &  $\ID_q $         \\ \cline{2-7}
   $\Lambda$     &    &   $\Id_{f'}$     &     $\beta$                &     $\Theta\dact \beta' $    &    $(\Theta \uact \beta)\dact \beta'$      &  $\ID_q $       \\ \cline{2-7}
                 &  $\varsigma_0(\ID_q)  $   &   $\Id_{\id_{a'}}$  &     $\Id_{\id_a}$             &     $\ID_q $    &        $\ID_q $                                   &  $\ID_q $        \\ \cline{2-7}
\end{tabular}   
 \end{center}
    \end{table}
\begin{table}[H]\begin{center}\caption{~}\label{com1pt2}
    \begin{tabular}{ r | l ||!{\vrule width 2pt} l | l | l | l |}
    \cline{2-6} 
        &  $\varsigma_1(\beta')$ &  $\Id_{g'}$     &    $\beta'$   & $\beta'$     &   $\Id_{\id_{a'}}$         \\ \cline{2-6}
     
        &  $\varsigma_0(\beta)$ &  $\beta$    & $\beta$   &  $\Id_g$     &  $\Id_{\id_a}$        \\ \cline{2-6}
                   
    \end{tabular} 

\vspace{12pt}

 \begin{tabular}{    r | l || l | l  !{\vrule width 2pt}  l | l | l | }\cline{2-7} 
               &  $\Delta_{\uact}(\Theta,\beta)$  &   $\Id_{g'}$     &      $\beta$                 &     $\Theta$               &    $\Theta \uact \beta$                      &  $\ID_q $       \\ \cline{2-7} 
   $\Lambda$     &   &   $\beta'$     &     $\beta$                  &     $\Theta$                &   $(\Theta \uact \beta)\dact\beta'$       &  $\ID_q $        \\ \cline{2-7}
                     & $\Delta_{\dact}(\Theta\uact \beta,\beta')$     &   $\beta'$     &     $\Id_g$                &     $\Theta\uact \beta $    &    $(\Theta \uact \beta)\dact \beta'$      &  $\ID_q $        \\ \cline{2-7}
                &  $\varsigma_0(\ID_q)  $   &   $\Id_{\id_{a'}}$  &     $\Id_{\id_a}$             &     $\ID_q $    &        $\ID_q $                                   &  $\ID_q $         \\ \cline{2-7}
\end{tabular}   
\end{center}
\end{table}
Next we show  $(\Theta \lact \eta)\uact\beta=(\Theta \uact\beta)\lact \eta.$ $$d\left( \Delta_\lact( \Theta \uact \beta, \eta )  \right)=[\Theta \uact \beta,\ (\Theta \uact\beta)\lact \eta,\ \ID_f \ \ \ | \ \ \Id_{q'},\ \eta     ] $$ so by the Matching Lemma it is enough to show that  $$[\Theta \uact \beta,\ (\Theta \lact \eta)\uact\beta,\ \ID_f \ \ \ | \ \ \Id_{q'},\ \eta     ] $$ is commutative. The following horn verifies this commutativity:
\begin{table}[H]\begin{center}
    \begin{tabular}{ r | l || l | l !{\vrule width 2pt} l | l | l |}
    \cline{2-7}

        &  $\Delta_{\uact}(\Theta,\beta)$  &     $\Id_{g'}$        &   $\beta$   &  $\Theta$              &   $\Theta \uact \beta$               &  $\ID_q$                    \\ \cline{2-7} 
        &  $\Delta_{\uact}(\Theta\lact \eta,\beta)$ &      $\Id_{g'}$       &   $\beta$   &  $\Theta \lact \eta$   &   $(\Theta \lact \eta)\uact\beta$    &  $\ID_p$               \\ \cline{2-7}       
        &  $ s_0(\beta) $   &    $\beta$              &   $\beta$   &     $\ID_g$            &   $\ID_f$                             &   $\ID_{\id_a}$       \\ \cline{2-7} 
          
    \end{tabular} 
    \vspace{12pt}

 \begin{tabular}{    r | l || l | l |  l !{\vrule width 2pt} l | l | }\cline{2-7} 
           &  $\Delta_{\lact}(\Theta,\eta)$ &   $\Theta$              &   $\Theta \lact \eta$              &     $\ID_g$           &    $\Id_{q'}$        &  $\eta$     \\ \cline{2-7} 
 $\Lambda$ &   &   $\Theta \uact \beta$ &   $(\Theta \lact \eta)\uact\beta$  &     $\ID_f$           &    $\Id_{q'}$        &  $\eta$      \\ \cline{2-7}
          &  $\varsigma_0(\eta)$    &   $\ID_q$              &   $\ID_p$                          &     $\ID_{\id_a}$     &    $\eta$            &  $\eta$    \\ \cline{2-7}
\end{tabular}   
\end{center}  \end{table}
Next we must check  $(\Theta \ract \eta')\uact\beta=(\Theta \uact\beta)\ract \eta'.$ $$d\left( \Delta_\ract( \Theta \uact \beta, \eta' )  \right)=[\Theta \uact \beta,\ (\Theta \uact\beta)\ract \eta',\ \ID_f \ \ \ | \ \ \eta' ,\   \Id_q ] $$ so by the Matching Lemma it is enough to show that  $$[\Theta \uact \beta,\ (\Theta \ract \eta')\uact\beta,\ \ID_f \ \ \ | \ \ \eta',\  \Id_q   ]$$ is commutative. The following horn verifies this commutativity:
\begin{table}[H]\begin{center}
    \begin{tabular}{ r | l || l | l !{\vrule width 2pt} l | l | l |}
    \cline{2-7}

        &  $\Delta_{\uact}(\Theta,\beta)$  &     $\Id_{g'}$        &   $\beta$   &  $\Theta$                       &   $\Theta \uact \beta$                &  $\ID_q$                    \\ \cline{2-7} 
        &  $\Delta_{\uact}(\Theta\lact \eta',\beta)$ &      $\Id_{g'}$       &   $\beta$   &  $\Theta \lact \eta'$   &   $(\Theta \lact \eta')\uact\beta$    &  $\ID_q$               \\ \cline{2-7}       
        &  $ s_0(\beta) $   &    $\beta$              &   $\beta$   &     $\ID_g$            &   $\ID_f$                             &   $\ID_{\id_a}$       \\ \cline{2-7} 
          
    \end{tabular} 
    \vspace{12pt}

 \begin{tabular}{    r | l || l | l |  l !{\vrule width 2pt} l | l |}\cline{2-7} 
          &  $\Delta_{\lact}(\Theta,\eta)$ &   $\Theta$              &        $\Theta \lact \eta'$        &     $\ID_g$           &    $\eta'$          &  $\Id_q$         \\ \cline{2-7} 
 $\Lambda$   &  &   $\Theta \uact \beta$  &   $(\Theta \lact \eta')\uact\beta$ &     $\ID_f$           &    $ \eta'$         &  $\Id_q $       \\ \cline{2-7}
            &  $s_0(\ID_q)$  &    $\ID_q$                &     $\ID_q$                      &     $\ID_{\id_a}$     &    $\Id_q$          &  $\Id_q $      \\ \cline{2-7}
\end{tabular}   
\end{center}  \end{table}
The final case we must check is $(\Theta \ract \eta')\dact\beta'=(\Theta \dact\beta')\ract \eta'.$  $$d\left( \Delta_\ract( \Theta \dact \beta', \eta' )  \right)=[\Theta \dact \beta',\ (\Theta \dact\beta')\ract \eta',\ \ID_g \ \ \ | \ \ \eta' ,\  \Id_{q}   ] $$ so by the Matching Lemma it is enough to show that  $$[\Theta \dact \beta',\ (\Theta \ract \eta')\dact\beta',\ \ID_g \ \ \ | \ \ \eta' ,\  \Id_{q}   ]    $$ is commutative. The following horn verifies this commutativity::
\begin{table}[H]\begin{center}
    \begin{tabular}{ r | l || l | l !{\vrule width 2pt} l | l | l |}
    \cline{2-7}

        &  $\Delta_{\dact}(\Theta,\beta')$  &              $\beta'$         &   $\Id_{g}$                &  $\Theta$              &   $\Theta \dact \beta'$               &  $\ID_q$          \\ \cline{2-7} 
        &  $\Delta_{\dact}(\Theta\ract \eta',\beta')$ &    $\beta'$         &   $\Id_{g}$    &  $\Theta \ract \eta'$   &   $(\Theta \ract \eta')\dact\beta'$    &    $\ID_q$                    \\ \cline{2-7}       
        &  $ s_0(\Id_g) $   &                              $\Id_{g}$             &   $\Id_{g}$     &     $\ID_{g}$            &   $\ID_g$                             &   $\ID_{\id_a}$     \\ \cline{2-7} 
          
    \end{tabular} 
    \vspace{12pt}

 \begin{tabular}{    r | l || l | l |  l !{\vrule width 2pt} l | l | }\cline{2-7} 
            &  $\Delta_{\ract}(\Theta,\eta')$ &   $\Theta$              &   $\Theta \ract \eta'$              &     $\ID_g$           &   $ \eta' $       &  $\Id_q$    \\ \cline{2-7} 
 $\Lambda$&     &   $\Theta \dact \beta'$ &   $(\Theta \ract \eta')\dact\beta'$  &     $ \ID_g $           &    $ \eta' $        &  $\Id_{q} $     \\ \cline{2-7}
         &  $\varsigma_0(\Id_q)$  &   $\ID_q$              &   $\ID_q$                          &     $\ID_{\id_a}$     &    $\Id_q$            &  $\Id_q$       \\ \cline{2-7}
\end{tabular}   
\end{center}  \end{table}
\end{proof}
As a corollary to the above proof, we note that we have shown there is a single horn that defines  ``action on both the top and bottom of a square at the same time.''
\begin{definition}\label{updownact} Let $\Theta\in \mbox{Sq}(g',g,q',q)$ and $\beta,\beta',\eta,\eta'$ as in Figure~\ref{Pipic}  $\Lambda_{\uact\dact}(\Theta,\beta,\beta') $ is defined to be face $1$ in Table~\ref{com1pt2}:
$$\Lambda_{\uact\dact}(\Theta,\beta,\beta') :=     [ \beta'  ,\      \beta     \ \  | \ \            \Theta ,\    -   ,\ \ID_q ].$$ We showed in the proof of Proposition~\ref{commute} that the filler of this horn is:
\begin{align*}\Delta_{\uact\dact}(\Theta,\beta,\beta') &=     [ \beta'  ,\      \beta     \ \  | \ \            \Theta ,\    (\Theta \uact \beta)\dact\beta'  ,\ \ID_q ] \\         
                                                       &=     [ \beta'  ,\      \beta     \ \  | \ \            \Theta ,\    (\Theta\dact \beta')  \uact \beta  ,\ \ID_q ] .\end{align*} 
                                                       
Similarly we can define   $\Lambda_{\lact\ract}(\Theta,\eta,\eta') :=     [   \Theta ,\    -   ,\ \ID_g   \ \  | \ \   \eta'  ,\      \eta          ]$   so that         
\begin{align*}\Delta_{\lact\ract}(\Theta,\eta,\eta') &=       [   \Theta ,\  (\Theta \lact \eta)\ract\eta'     ,\ \ID_g   \ \  | \ \   \eta'  ,\      \eta          ] \\         
                                                     &=       [   \Theta ,\   (\Theta\ract \eta')\lact \eta  ,\ \ID_g   \ \  | \ \   \eta'  ,\      \eta          ].      \end{align*} 
                                       
\end{definition}

\begin{proposition}[Compatibility of actions and square composition (\textbf{VDC4})] Let $\Theta,$ $\Pi$, $\eta$, and $\pi'$ as shown below. Then $(\Pi \boxvert \Theta)\lact \eta=\Pi \boxvert (\Theta\lact \eta)$ and $(\Pi\ract \pi') \boxvert \Theta=(\Pi \boxvert \Theta)\ract \pi',$ and similar laws hold for $\uact$ and $\dact$.
\begin{center}
\begin{tikzpicture}[scale=1.8,auto]

\begin{scope}

\node (10) at (1,1) {$b$};
\node (00) at (0,1) {$a$};
\node (11) at (1,0) {$b'$};
\node (01) at (0,0) {$a'$};
\node[rotate=-45] at (.5,.5){$\Rightarrow$};
\node[scale=.8] at (.6,.65){$\Theta$};

\node (20) at (2,1) {$c$};
\node (21) at (2,0) {$c'$};
\node[rotate=-45] at (1.5,.5){$\Rightarrow$};
\node[scale=.8] at (1.6,.65){$\Pi$};

\path[->] (00) edge node[midway]{$g$}(10);
\path[->] (00) edge node[midway,right]{$q$}(01);
\path[->] (01) edge node[midway,swap]{$g'$}(11);
\path[->] (10) edge node[midway]{$q'$}(11);
\path[->] (10) edge node[midway]{$h$}(20);
\path[->] (11) edge node[midway,swap]{$h'$}(21);
\path[->] (20) edge node[midway,left]{$q''$}(21);
\draw[->] (00) to [out=-150,in=150] node(f)[midway,left]{$p$} (01) ;
\node[rotate=0] (areta) at ($(f)+(.35,0)$){$\Rightarrow$};
\node (eta) at ($(areta)+(0,.15)$){$\eta$};

\draw[->] (20) to [out=-30,in=30] node(f)[midway,right]{$p'$} (21) ;
\node[rotate=180] (areta) at ($(f)+(-.35,0)$){$\Rightarrow$};
\node (eta) at ($(areta)+(0,.15)$){$\pi'$};
\end{scope}
\end{tikzpicture}
\end{center}
\end{proposition}
\begin{proof} We first show $(\Pi \boxvert \Theta)\lact \eta=\Pi \boxvert (\Theta\lact \eta).$  $$d\left(\Delta_{\boxvert}(\Pi,\Theta\lact \eta ) \right)=[ \chi(h',g') ,\ \chi(h,g)  \ \  | \ \   \Pi,\   \Pi \boxvert (\Theta\lact \eta)  ,\   \Theta\lact \eta] $$ so by the Matching Lemma it suffices to show $$[ \chi(h',g') ,\ \chi(h,g)  \ \  | \ \   \Pi,\   (\Pi \boxvert \Theta)\lact \eta  ,\   \Theta\lact \eta].$$
The following horn verifies this commutativity:
\begin{table}[H]\begin{center}
    \begin{tabular}{ r | l || l | l !{\vrule width 2pt} l | l | l |}
    \cline{2-7}

                 &  $\Delta_{\boxvert}(\Pi,\Theta)$  &            $\chi(h',g')$        &  $\chi(h,g)$            &   $\Pi$              &   $\Pi\boxvert\Theta$               &  $\Theta$                    \\ \cline{2-7} 
   $\Lambda$     &   &      $\chi(h',g')$       &   $\chi(h,g)$        &  $\Pi$          &   $(\Pi \boxvert \Theta)\lact\eta$      &  $\Theta\lact \eta$               \\ \cline{2-7}       
                 &  $ s_0(\chi(h,g)) $            &              $\chi(h,g)$                &  $\chi(h,g)$     &      $\ID_h$           &   $\ID_{h\circ g}$        &   $\ID_g$       \\ \cline{2-7} 
          
    \end{tabular} 
    \vspace{12pt}

 \begin{tabular}{    r | l || l | l |  l !{\vrule width 2pt} l | l | }\cline{2-7} 
           &  $s_0(\Pi)$ &    $\Pi$                 &  $\Pi$                                      &             $\ID_h$             &    $\Id_{q'}$                         &    $\Id_{q'}$        \\ \cline{2-7} 
          &    $\Delta_{\lact}(\Pi \boxvert \Theta, \eta)$   &   $\Pi \boxvert \Theta$  &   $(\Pi \boxvert \Theta)\lact\eta$          &     $\ID_{h\circ g}$                    &    $\Id_{q'}$        &  $\eta$      \\ \cline{2-7}
            &  $\Delta_{\lact}(\Theta, \eta)$  &   $\Theta$               &   $\Theta \lact \eta$                       &     $\ID_{g}$     &    $\Id_{q'}$         &  $\eta$    \\ \cline{2-7}
\end{tabular}   
\end{center}  \end{table}
This the compatibility $(\Pi\ract \pi') \boxvert \Theta=(\Pi \boxvert \Theta)\ract \pi'$ follows similarly. The compatibility of $\uact$ and $\dact$ with vertical composition follows by symmetry.
\end{proof} 

\begin{proposition}[Compatibility of $\ID$ with identities for objects (\textbf{VDC5})] If $a$ is an object of $\mbox{Vdc}(X)$, then $ \ID_{(\Id_a)_h}=\ID_{(\Id_a)_v}.$
\end{proposition}
\begin{proof} Applying definitions, this statement works out to $s_0 \varsigma_0 a = \varsigma_0 s_0 a$ which is a bisimplicial identity.\end{proof}
\begin{proposition}[Compatibility of $\ID$ with 1-morphism composition (\textbf{VDC6})] Let $a\stackrel{f}{\ra}b \stackrel{g}{\ra}c$ be horizontal $1$-morphisms in $\mbox{Vdc}(X)$, we have $\ID_g \boxvert \ID_f =\ID_{g\circ f}.$ The same axiom holds for vertical $1$-morphisms.
\end{proposition} 
\begin{proof}\begin{align*} d\left(s_0(\chi(g,f))\right)&=[\chi(g,f),\ \chi(g,f) \ \ | \ \ \ID_g, \     \ID_{g\circ f},\    \ID_f]\\
d\left(\Delta_{\boxvert}(\ID_g, \ID_f ) \right)&=[\chi(g,f),\ \chi(g,f) \ \ | \ \ \ID_g, \     \ID_g\boxvert \ID_f,\    \ID_f].\end{align*}
By the Matching Lemma, we conclude $ \ID_{g\circ f}=\ID_g\boxvert \ID_f.$ The axiom for vertical composition follows by symmetry.

\end{proof}
\begin{proposition}[Compatibility of the $\ID$ with the actions (\textbf{VDC7})]
Let $\beta:f\Rightarrow g$ be a horizontal $2$-morphism. Then $\ID_g \uact\beta\dact\beta =\ID_f.$ Similarly if $\eta:p\Rightarrow q$ is a vertical $2$-morphism, $\ID_q \lact\eta\ract \eta =\ID_p$.
\end{proposition}
\begin{proof}
\begin{align*}
d(s_0(\beta))&=[ \beta,\ \beta     \ \ | \ \   \ID_g,\ \ID_f,  \Id_{\id_a}]\\ 
d\left(\Delta_{\uact\dact}(\ID_g,\beta,\beta)\right)&=[\beta,\beta \ \ | \ \ \ID_g,\ \ID_g\uact\beta \dact \beta,\ \ID_{\id_a} ]
\end{align*}
So by the Matching Lemma, $\ID_g \uact\beta\dact\beta =\ID_f.$ The vertical case $\ID_q \lact\eta\ract \eta =\ID_p$ follows by symmetry.
\end{proof}
\begin{proposition}[Interchange of whiskering and square composition (\textbf{VDC8})] Let $\Pi, \Theta$ and $\beta, \gamma$ be as shown below.
\begin{center}
\begin{tikzpicture}[scale=1.8,auto]

\begin{scope}

\node (10) at (1,1) {$b$};
\node (00) at (0,1) {$a$};
\node (11) at (1,0) {$b'$};
\node (01) at (0,0) {$a'$};
\node[rotate=-45] at (.5,.45){$\Rightarrow$};
\node[scale=.8] at (.6,.6){$\Theta$};

\node (20) at (2,1) {$c$};
\node (21) at (2,0) {$c'$};
\node[rotate=-45] at (1.5,.5){$\Rightarrow$};
\node[scale=.8] at (1.6,.65){$\Pi$};

\path[->] (00) edge node[midway,below]{$g$}(10);
\path[->] (00) edge node[midway,swap]{$q$}(01);
\path[->] (01) edge node[midway,swap]{$g'$}(11);
\path[->] (10) edge node[midway]{$q'$}(11);
\path[->] (10) edge node[midway,below]{$h$}(20);
\path[->] (11) edge node[midway,swap]{$h'$}(21);
\path[->] (20) edge node[midway]{$q''$}(21);

\draw[->] (00) to [out=60,in=120] node(e)[midway]{$e$} (10) ;

\node[rotate=-90] (areta1) at ($(e)+(0,-.32)$){$\Rightarrow$};
\node (eta) at ($(areta1)+(.15,0)$){$\beta$};

\draw[->] (10) to [out=60,in=120] node(f)[midway]{$f$} (20) ;

\node[rotate=-90] (areta) at ($(f)+(0,-.38)$){$\Rightarrow$};
\node (eta) at ($(areta)+(0.15,0)$){$\gamma$};

\end{scope}
\end{tikzpicture}
\end{center}
then $\Pi \boxvert (\Theta \uact \beta)= (\Pi\boxvert \Theta)\uact (h \rhd \beta)$ and $(\Pi \uact \gamma) \boxvert \Theta = (\Pi\boxvert \Theta)\uact (\gamma\lhd g).$ Similar laws hold for $\dact, \lact,$ and $\ract.$ 
\end{proposition}\begin{proof} First we show $\Pi \boxvert (\Theta \uact \beta)= (\Pi\boxvert \Theta)\uact (h \rhd \beta).$
$$d\left(\Delta_{\boxvert}( \Pi\boxvert\Theta  ,  h \rhd \beta ) \right)=[\Id_{h'\circ g'},\ h\rhd\beta \ \  | \ \ \Pi\boxvert\Theta,\ (\Pi\boxvert\Theta)\uact (h\rhd \beta),\  \Id_q  ]$$
so by the Matching Lemma, it suffices to show the sphere $$[\Id_{h'\circ g'},\ h\rhd\beta \ \  | \ \ \Pi\boxvert\Theta,\ \Pi \boxvert (\Theta \uact \beta),\  \ID_q  ]$$ is commutative. The following horn verifies this commutativity:
\begin{table}[H]\begin{center}
    \begin{tabular}{ r | l ||!{\vrule width 2pt} l | l | l | l |}
    \cline{2-6} 
     &   $\varsigma_0( \chi(h',g')) $      &   $\chi(h',g')$  & $\chi(h',g')$      &  $\Id_{h'\circ g'}$   &  $\Id_{g'}$           \\ \cline{2-6}
     
     &  $\Delta_{\rhd}(h,\beta)$     &   $\chi(h,g)$   &   $\chi(h,f)$     & $h\rhd \beta$     &    $\beta$                \\ \cline{2-6}
                   
    \end{tabular} 

\vspace{12pt}

 \begin{tabular}{    r | l || l | l  !{\vrule width 2pt}  l | l | l | }\cline{2-7} 
            &  $\Delta_{\boxvert}(\Pi,\Theta)$  &  $\chi(h',g')$     &      $\chi(h,g)$              &     $\Pi$                &    $\Pi\boxvert\Theta$                  &  $\Theta $       \\ \cline{2-7} 
            &  $\Delta_{\boxvert}(\Pi,\Theta\uact\beta)$   & $\chi(h',g')$       &     $\chi(h,f)$     &       $\Pi$              &   $\Pi \boxvert (\Theta \uact \beta)$   &  $\Theta\uact\beta $  \\ \cline{2-7}
 $\Lambda$  &   &   $\Id_{h'\circ g'}$     &     $h\rhd \beta$       &   $\Pi\boxvert\Theta$    &    $\Pi \boxvert (\Theta \uact \beta)$  &  $\ID_q $        \\ \cline{2-7}
             &  $\Delta_{\uact}(\Theta,\beta) $   &   $\Id_{ g'}$           &     $\beta$             &     $\Theta $            &          $\Theta\uact\beta $                       &  $\ID_q $      \\ \cline{2-7}
\end{tabular}   
\end{center}
\end{table}

Next we show $(\Pi \uact \gamma) \boxvert \Theta = (\Pi\boxvert \Theta)\uact (\gamma\lhd g).$ By Lemma~\ref{altaction} $\Pi \uact \gamma=\Pi \uact \wh{\gamma},$ so we show $(\Pi \uact \wh{\gamma}) \boxvert \Theta = (\Pi\boxvert \Theta)\uact \wh{\gamma\lhd g}.$  By the Matching Lemma applied to $\Delta_{\uact}(\Pi\boxvert\Theta,\wh{\gamma\lhd g})$, it suffices to show that $$[\widehat{\Id_{h'\circ g'}},\ \widehat{\gamma\lhd g} \ \  | \ \ \ID_{q''},\ (\Pi \uact \wh{\gamma}) \boxvert \Theta,\  \Pi\boxvert\Theta  ]$$ is commutative. The following horn verifies this commutativity:
\begin{table}[H]\begin{center}
    \begin{tabular}{ r | l ||!{\vrule width 2pt} l | l | l | l |}
    \cline{2-6} 
     &   $s_2( \chi(h',g')) $      &   $\wh{\Id_{h'}}$  &$\wh{\Id_{h'\circ g'}}$    &   $\chi(h',g')$    &   $\chi(h',g')$           \\ \cline{2-6}
     
     &  $\Delta_{\lhd}(\gamma,g)$     &   $\wh{\gamma}$    &   $\wh{\gamma \lhd g}$      & $\chi(f,g)$      &   $\chi(h,g)$                \\ \cline{2-6}
                   
    \end{tabular} 

\vspace{12pt}

 \begin{tabular}{    r | l || l | l  !{\vrule width 2pt}  l | l | l | }\cline{2-7} 
           &  $\Delta_{\uact}(\Pi,\wh{\gamma})$ &  $\wh{\Id_{h'}}$     &      $\wh{\gamma}$                    &     $\ID_{q''}$                     &   $\Pi \uact \wh{\gamma}$                   &    $ \Pi $         \\ \cline{2-7} 
   
      $\Lambda$ &  & $\wh{\Id_{h'\circ g'}}$      &     $\wh{\gamma \lhd g}$     &       $\ID_{q''}$      & $(\Pi \uact \wh{\gamma}) \boxvert \Theta $ &      $\Pi\boxvert\Theta $    \\ \cline{2-7}

               &    $\Delta_{\boxvert}(\Pi \uact \wh{\gamma},\Theta)$   &    $\chi(h',g')$    &     $\chi(f,g)$       &   $\Pi \uact \wh{\gamma}$          & $(\Pi \uact \wh{\gamma}) \boxvert \Theta $    &       $\Theta $    \\ \cline{2-7}

            &  $\chi(h',g')$          &    $\chi(h,g)$            &     $\Pi $                   &          $\Pi\boxvert\Theta $         &  $\Theta  $       & $\Delta_{\boxvert}(\Pi,\Theta)$   \\ \cline{2-7}
\end{tabular}   
\end{center}
The interchange law for $\lact$ follows similarly, and then the laws for  $\dact$ and $\ract$ follow by symmetry.
\end{table} \end{proof}
\begin{proposition}[Compatibility of composition and the associators (\textbf{VDC9})] Horizontal (vertical) composition of squares is associative up to the associator of $H$ ($V$). For horizontal composition, let $\Theta,$ $\Pi$ and $\Sigma$ be as shown below:
\begin{center}
\begin{tikzpicture}[scale=1.8,auto]

\begin{scope}

\node (10) at (1,1) {$b$};
\node (00) at (0,1) {$a$};
\node (11) at (1,0) {$b'$};
\node (01) at (0,0) {$a'$};
\node[rotate=-45] at (.5,.5){$\Rightarrow$};
\node[scale=.8] at (.6,.65){$\Theta$};

\node (20) at (2,1) {$c$};
\node (21) at (2,0) {$c'$};
\node (30) at (3,1) {$d$};
\node (31) at (3,0) {$d'$};
\node[rotate=-45] at (1.5,.5){$\Rightarrow$};
\node[scale=.8] at (1.6,.65){$\Pi$};
\node[rotate=-45] at (2.5,.5){$\Rightarrow$};
\node[scale=.8] at (2.6,.65){$\Sigma$};

\path[->] (00) edge node[midway,above]{$f$}(10);
\path[->] (00) edge node[midway,swap]{$q$}(01);
\path[->] (01) edge node[midway,swap]{$f'$}(11);
\path[->] (10) edge node[midway]{$q'$}(11);
\path[->] (10) edge node[midway]{$g$}(20);
\path[->] (11) edge node[midway,swap]{$g'$}(21);
\path[->] (20) edge node[midway]{$q''$}(21);
\path[->] (20) edge node[midway]{$h$}(30);
\path[->] (21) edge node[midway,swap]{$h'$}(31);
\path[->] (30) edge node[midway]{$q'''$}(31);

\end{scope}
\end{tikzpicture}
\end{center}
then  $((\Sigma \boxvert \Pi) \boxvert \Theta)\uact \alpha_{h,g,f}\dact \alpha_{h',g',f'}=\Sigma \boxvert (\Pi \boxvert \Theta)$, (where $\alpha_{-,-,-}$ is the associator of $H$). A similar compatibility law holds for vertical composition.
\end{proposition}
\begin{proof}From Definition~\ref{updownact}: 
\begin{align*}d(\Delta_{\uact\dact}((\Sigma \boxvert \Pi)& \boxvert \Theta, \alpha_{h,g,f},  \alpha_{h',g',f'})) \\ \hspace{40pt} &= [\alpha_{h',g',f'},\  \alpha_{h,g,f} \ \ | \ \ (\Sigma \boxvert \Pi) \boxvert \Theta,\ ((\Sigma \boxvert \Pi) \boxvert \Theta)\uact \alpha_{h,g,f}\dact \alpha_{h',g',f'},\ \ID_q].\end{align*} By the Matching Lemma applied to this sphere it suffices to show $$[\alpha_{h',g',f'},\  \alpha_{h,g,f} \ \ | \ \ (\Sigma \boxvert \Pi) \boxvert \Theta,\ \Sigma \boxvert (\Pi \boxvert \Theta),\ \ID_q]$$ is commutative. The following  table proof verifies this commutativity:
\begin{table}[H]\begin{center}\caption{~}
    \begin{tabular}{ r | l ||!{\vrule width 2pt} l | l | l | l |}
    \cline{2-6} 
     &  $\Delta_{-}(\widetilde{\alpha}_{h',g',f'})$      & $\chi(h'\circ g', f') $     &  $\widetilde{\alpha}_{h',g',f'}$      &  $\alpha_{h',g',f'}$   & $\Id_{f'}$          \\ \cline{2-6}
     
     &  $\Delta_{-}(\widetilde{\alpha}_{h,g,f})$     &  $\chi(h\circ g, f)$      & $\widetilde{\alpha}_{h,g,f}$               & $\alpha_{h,g,f}$              &  $\Id_{f}$             \\ \cline{2-6}      
    \end{tabular} 

\vspace{12pt}

 \begin{tabular}{    r | l || l | l  !{\vrule width 2pt}  l | l | l | }\cline{2-7} 
           &  $\Delta_{\boxvert}(\Sigma \boxvert \Pi,\Theta)$    &  $\chi(h'\circ g', f')$     &     $\chi(h\circ g, f)$                                          &     $\Sigma \boxvert \Pi$                     &   $(\Sigma \boxvert \Pi) \boxvert \Theta$                   &    $  \Theta$      \\ \cline{2-7} 
   
$\odot $        &    (Table~\ref{asscomp2})     &    $\widetilde{\alpha}_{h',g',f'}$    &     $\widetilde{\alpha}_{h,g,f}$       &  $\Sigma\boxvert\Pi$          & $\Sigma \boxvert (\Pi \boxvert \Theta)$    &       $\Theta $        \\ \cline{2-7}
               
      $\Lambda$  & & $\alpha_{h',g',f'}$      &     $\alpha_{h,g,f} $                                 &       $(\Sigma \boxvert \Pi) \boxvert \Theta$      & $\Sigma \boxvert (\Pi \boxvert \Theta)$ &            $ \ID_q$    \\ \cline{2-7}

            &  $\varsigma_0(\Theta)$ &  $\Id_{f'}$          &    $\Id_f$            &       $\Theta $                 &      $\Theta  $           &   $ \ID_q$     \\ \cline{2-7}
\end{tabular}   
\end{center}\end{table}

\begin{table}[H]\begin{center}\caption{~}\label{asscomp2}
    \begin{tabular}{ r | l ||!{\vrule width 2pt} l | l | l | l |}
    \cline{2-6} 
     &  $\Delta_{-}(\widetilde{\alpha}_{h',g',f'})$  &$\chi(h',g') $   & $\chi(h'\circ g', f') $     &  $\widetilde{\alpha}_{h',g',f'}$                  & $\chi(g',f')$          \\ \cline{2-6}
     
     &  $\Delta_{-}(\widetilde{\alpha}_{h,g,f})$         &$\chi(h,g)$      &  $\chi(h\circ g, f)$      & $\widetilde{\alpha}_{h,g,f}$                    &  $\chi(g,f)$             \\ \cline{2-6}      
    \end{tabular} 

\vspace{12pt}

 \begin{tabular}{    r | l || l | l  !{\vrule width 2pt}  l | l | l |}\cline{2-7} 
               &  $\Delta_{\boxvert}(\Sigma , \Pi)$   &  $\chi(h',g')$     &     $\chi(h, g)$                                          &     $\Sigma$                     &   $\Sigma\boxvert\Pi$             &    $  \Pi$   \\ \cline{2-7}

                   & $\Delta_{\boxvert}(\Sigma , \Pi \boxvert \Theta)$ & $\chi(h'\circ g',f')   $      &     $\chi(h\circ g,f) $                                     &        $\Sigma$           & $\Sigma \boxvert (\Pi \boxvert \Theta)$ &            $\Pi \boxvert \Theta$ \\ \cline{2-7}
       
$\Lambda $        &             &    $\widetilde{\alpha}_{h',g',f'}$    &     $\widetilde{\alpha}_{h,g,f}$       &  $\Sigma\boxvert\Pi$          & $\Sigma \boxvert (\Pi \boxvert \Theta)$    &       $\Theta $           \\ \cline{2-7}    
          
           &  $\Delta_{\boxvert}(\Pi,\Theta)$  &  $\chi(g',f') $        &    $\chi(g,f) $          &     $  \Pi$                  &     $\Pi \boxvert \Theta$        &   $\Theta $       \\ \cline{2-7}
\end{tabular}   
\end{center}\end{table} The axiom for vertical composition follows by symmetry.
\end{proof}
\begin{proposition}[Compatibility of the pseudo-identity squares and unitors (\textbf{VDC10})] Let $\Theta\in\mbox{Sq}(q',q,g',g)$ as in Figure~\ref{Pipic}. Then  \begin{align*}(\Theta\boxvert \ID_{q})\uact (\rho_{g})_h\dact (\rho_{g'})_h=\Theta= (\ID_{q'}\boxvert\Theta)\uact (\lambda_g)_h \dact (\lambda_{g'})_h \\  (\Theta\boxminus \ID_{g})\lact (\rho_{q})_v\ract (\rho_{q'})_v=\Theta= (\ID_{g'}\boxminus\Theta)\lact (\lambda_q)_v\ract (\lambda_{q'})_v.
\end{align*}\end{proposition}
\begin{proof} First we show $(\Theta\boxvert \ID_{q})\uact (\rho_{g})_h\dact (\rho_{g'})_h=\Theta.$ Applying the Matching Lemma to $\Delta_{\uact\dact}(\Theta\boxvert\ID_q,\rho_{g},\rho_{g'})$, it is enough to show  the sphere $$[\rho_{g'}=\underline{\Id_{g'}},\ \rho_g=\underline{\Id_g} \  \  |  \ \ \Theta \boxvert \ID_q,\ \Theta,\ \ID_q]$$ is commutative. The following horn verifies this commutativity: 
\begin{table}[H]\begin{center}
    \begin{tabular}{ r | l ||!{\vrule width 2pt} l | l | l | l |}
    \cline{2-6} 
     &  $\Delta_{-}(\Id_{g'})$  &$\chi(g',\id_{a'})$   & $\Id_{g'}$     &  $\underline{\Id_{g'}}$                  & $\Id_{\id_{a'}}$          \\ \cline{2-6}
     
     &  $\Delta_{-}(\Id_g)$         &$\chi(g,\id_a)$      &  $\Id_g$      & $\underline{\Id_g}$                    &  $\Id_{\id_a}$             \\ \cline{2-6}      
    \end{tabular} 

\vspace{12pt}

 \begin{tabular}{    r | l || l | l  !{\vrule width 2pt}  l | l | l | }\cline{2-7} 
           &  $\Delta_{\boxvert}(\Theta , \ID_q)$   &  $\chi(g',\id_{a'})$     &             $\chi(g,\id_a)$                                    &     $\Theta$                     &  $\Theta\boxvert \ID_q$            &    $  \ID_q$       \\ \cline{2-7} 
               
                & $\varsigma_0(\Theta)$   & $\Id_{g'}$      &     $\Id_g$                                     &        $\Theta$           & $\Theta$ &            $\ID_q$   \\ \cline{2-7}
       
$\Lambda $      &     & $\underline{\Id_{g'}}$  &   $\underline{\Id_g}$                      &  $\Theta\boxvert \ID_q$          & $\Theta$    &       $ \ID_q$           \\ \cline{2-7}    
          
             &  $\varsigma_0(\ID_q)$ &  $\Id_{\id_{a'}} $        &    $\Id_{\id_a}$                       &     $  \ID_q$                  &     $  \ID_q$          &      $  \ID_q$        \\ \cline{2-7}
\end{tabular}   
\end{center}\end{table}
For $\Theta= (\ID_{q'}\boxvert\Theta)\uact \lambda_g \dact \lambda_{g'} $, by the Matching Lemma applied to $\Delta_{\uact\dact}(\Theta\boxvert\ID_q, \lambda_g,\lambda_{g'})$ it is enough to show $$[\lambda_{g'}=\underline{(\wh{\Id_{g'}})},\ \lambda_g=\underline{(\wh{\Id_g})} \  \  |  \ \  \ID_{q'}\boxvert\Theta ,\ \Theta,\ \ID_q].$$ The following horn verifies this commutativity:
\begin{table}[H]\begin{center}
    \begin{tabular}{ r | l ||!{\vrule width 2pt} l | l | l | l |}
    \cline{2-6} 
     &  $\Delta_{-}(\wh{\Id_{g'}})$  &$\chi(\id_{b'},g')$   & $\wh{\Id_{g'}}$     &  $\underline{(\wh{\Id_{g'}})}$                  & $\Id_{g'} $          \\ \cline{2-6}
     
     &  $\Delta_{-}(\wh{\Id_g})$         &$\chi(\id_b,g)$      &  $\wh{\Id_g}$      & $\underline{(\wh{\Id_g})}$                    &  $\Id_g$             \\ \cline{2-6}      
    \end{tabular} 

\vspace{12pt}

 \begin{tabular}{    r | l || l | l  !{\vrule width 2pt}  l | l | l | }\cline{2-7} 
             &  $\Delta_{\boxvert}(\ID_{q'} ,\Theta )$  &  $\chi(\id_{b'},g')$      &             $\chi(\id_b,g)$                                    &     $\ID_{q'}$                     &  $\ID_{q'}\boxvert\Theta $            &    $  \Theta$      \\ \cline{2-7} 
               
                 & $\varsigma_1(\Theta)$   & $\wh{\Id_{g'}}$      &     $\wh{\Id_g}$                                     &         $\ID_{q'}$        & $\Theta$ &     $\Theta$           \\ \cline{2-7}
       
$\Lambda $   &   & $\underline{(\wh{\Id_{g'}})}$  &   $\underline{(\wh{\Id_g})}$                      &  $\ID_q'\boxvert \Theta$          & $\Theta$    &       $ \ID_q$                \\ \cline{2-7}    
          
            &  $\varsigma_0(\Theta)$ &  $\Id_{g'} $         &    $\Id_g$                       &     $  \Theta$                  &     $  \Theta$          &      $  \ID_q$         \\ \cline{2-7}
\end{tabular}   
\end{center}\end{table}
The compatibility of the unitors of $V$ with vertical composition follows by symmetry.
\end{proof}
\begin{proposition}[Interchange for squares (\textbf{VDC11})] Let $\Theta, \Theta',\Pi,\Pi'$ be squares as shown below. Then $$(\Pi'\boxvert\Theta')\boxminus(\Pi \boxvert \Theta)=(\Pi'\boxminus\Pi)\boxvert(\Theta' \boxminus \Theta).$$ 
\begin{center}
\begin{tikzpicture}[scale=1.8,auto]

\begin{scope}

\node (10) at (1,1) {$b$};
\node (00) at (0,1) {$a$};
\node (11) at (1,0) {$b'$};
\node (01) at (0,0) {$a'$};
\node[rotate=-45] at (.5,.5){$\Rightarrow$};
\node[scale=.8] at (.6,.65){$\Theta$};

\node (20) at (2,1) {$c$};
\node (21) at (2,0) {$c'$};
\node(22) at (2,-1){$c''$};
\node[rotate=-45] at (1.5,.5){$\Rightarrow$};
\node[scale=.8] at (1.6,.65){$\Pi$};

\node (12) at (1,-1) {$b''$};
\path[->] (00) edge node[midway]{$f$}(10);
\path[->] (00) edge node[midway,swap]{$p$}(01);
\path[->] (01) edge node[midway,swap]{$f'$}(11);
\path[->] (10) edge node[midway]{$p'$}(11);
\path[->] (10) edge node[midway]{$g$}(20);
\path[->] (11) edge node[midway,swap]{$g'$}(21);
\path[->] (20) edge node[midway]{$p''$}(21);
\path[->](21) edge node[midway]{$q''$}(22);
\path[->](12) edge node[midway,below]{$g''$}(22);

\node (01') at (0,-1) {$a''$};
\node[rotate=-45] at (0.5,-.58){$\Rightarrow$};
\node[scale=.8] at (0.6,-.43){$\Theta'$};

\path[->] (11) edge node[midway]{$q'$}(12);
\path[->] (01) edge node[midway,swap]{$q$}(01');
\path[->] (01') edge node[midway,swap]{$f''$}(12);
\node[rotate=-45] at (1.5,-.58){$\Rightarrow$};
\node[scale=.8] at (1.6,-.43){$\Pi'$};

\end{scope}
\end{tikzpicture}
\end{center}
\end{proposition}
\begin{proof} By the Matching Lemma applied to $\Delta_{\boxminus}(\Pi'\boxvert\Theta',\Pi \boxvert \Theta),$ it suffices to show the sphere
$$[ \Pi'\boxvert\Theta',\ (\Pi' \boxminus \Pi)\boxvert (\Theta'\boxminus\Theta),\ \Pi\boxvert\Theta \ \ | \ \ \chi(q'',p''),\ \chi(q,p)]$$
is commutative. The following horn verifies this commutativity:

\begin{table}[H]\begin{center}
    \begin{tabular}{ r | l || l | l !{\vrule width 2pt} l | l | l |}
    \cline{2-7}        &  $\Delta_{\boxvert}(\Pi',\Theta')$   &    $\chi(g'',f'')$        &  $\chi(g',f')$     &   $\Pi'$              &   $\Pi'\boxvert\Theta'$              &     $\Theta'$          \\ \cline{2-7} 
                      &    $\Delta_{\boxvert}(\Pi' \boxminus \Pi,\Theta'\boxminus\Theta)$        &     $\chi(g'',f'')$            &      $\chi(g,f)$         & $\Pi'\boxminus\Pi$       &   $(\Pi' \boxminus \Pi)\boxvert (\Theta'\boxminus\Theta)$  &     $\Theta'\boxminus\Theta$         \\ \cline{2-7}
        &  $\Delta_{\boxvert}(\Pi,\Theta)$                   &  $\chi(g',f')$                &   $\chi(g,f)$                        &     $\Pi$     &    $\Pi\boxvert\Theta$            &  $\Theta$      \\ \cline{2-7}

    \end{tabular} 
    \vspace{12pt}

 \begin{tabular}{    r | l || l | l |  l !{\vrule width 2pt} l | l | }\cline{2-7} 
         
                  &   $\Delta_{\boxminus}(\Pi',\Pi)$            &  $\Pi'$                  &   $\Pi' \boxminus \Pi$               &  $\Pi$                     &  $\chi(q'',p'')$  &     $\chi(q',p')$                       \\ \cline{2-7} 
 $\Lambda$   &          &  $\Pi'\boxvert\Theta'$   &   $(\Pi' \boxminus \Pi)\boxvert (\Theta'\boxminus\Theta)$    &  $\Pi\boxvert\Theta$  &  $\chi(q'',p'')$  &      $\chi(q,p)$                    \\ \cline{2-7}       
                 &   $\Delta_{\boxminus}(\Theta',\Theta)$     &   $\Theta'$              &   $\Theta'\boxminus\Theta$              &   $\Theta$               &  $ \chi(q',p') $   &    $\chi(q,p)$                 \\ \cline{2-7} 
          
\end{tabular}   
\end{center}  \end{table}
\end{proof}
\section{The bisimplicial nerve of a VDC}\subsection{The definition of $N(\D)$}
Let $\D$ be a small VDC. In this section, we will construct a bisimplicial nerve $N(\D)$, generalizing the Duskin nerve of a small $(2,1)$-category. As in the definition of the Duskin nerve, we first define the truncated nerve $N(\D)|^3_0,$ which is a presheaf of sets on $\Delta \times \Delta|^3_0$, i.e. an object of $(\Delta\times \Delta|^3_0)^{\op}\cat{-Set}.$ 

The two ``edges'' of $N(\D)$ are given by the Duskin nerve of the horizontal and vertical small $(2,1)$-categories of $\D.$ Specifically, $N(\D)_{0i}=N(\D_h)_i$ and $N(\D)_{i0}=N(\D_v)_i$ for $0 \leq i \leq 3$, where $N(\B)$ denotes the Duskin nerve of a small $(2,1)$-category $\B,$ as defined in Section~\ref{duskinnerve}. Recalling the definition of the Duskin nerve, this means the objects of  $N(\D)_{00}$ are the shared objects of $\D_h$ and $\D_v$, the objects of $N(\D)_{01}$ (respectively $N(\D)_{10}$) are the horizontal (vertical) $1$-morphisms, and the objects of $N(\D)_{02}$ consist of a triple $(h,g,f)$ of horizontal $1$-morphism with a horizontal $2$-morphism $\beta:g \Rightarrow h \circ f.$

We define the elements of $N(\D)_{11}$ to be the squares of $\D$, with face maps given by $d(S)=[g',\ g,\ \ | \ \ q',\ q]$ for $S \in \mbox{Sq}(g',g,q',q).$ An element $x \in N(\D)_{12}$ is a sphere $$[(h',g',f'\vbar\beta' ),\ (h,g,f\vbar \beta) \ \ | \ \ \Sigma,\ \Pi,\ \Theta]$$ in $N(\D)|^2_0$ (which we have now defined) such that  \begin{equation}\label{12condition}(\Sigma\boxvert\Theta)\uact \beta \dact \beta'=\Pi.\end{equation} The faces of $x$ are (of course) given by  $$d(x)=[(h',g',f'\vbar \beta'),\  ( h,g,f \vbar \beta)  \ \    |   \  \   \Sigma    ,\  \Pi    ,\ \Theta  ].$$ Similarly an element $y \in N(\D)_{21}$ is a sphere $$[ \Sigma,\ \Pi,\ \Theta  \ \ | \ \  (r',q',p'\vbar \eta'),\ (r,q,p \vbar \eta,) ]$$ in $N(\D)|^2_0$  such that  \begin{equation}\label{21condition}(\Sigma\boxminus \Theta)\lact \eta\ract\eta'=\Pi.\end{equation} The faces of $y$ are given by  $$d(y)=[ \Sigma    ,\  \Pi    ,\ \Theta   \ \    |   \  \  (r',q',p'\vbar \eta') ,\ (r,q,p \vbar \eta ) ].$$

We have defined each of the sets in $N(\D)|^3_0,$ and their face maps, so we now define the degeneracy maps. The degeneracy maps in the first row and column of $N(\D)|^3_0$ are given by the corresponding degeneracy maps in the Duskin nerves $N(\D_h)_i$ and $N(\D_v)_i$, leaving us to define:
\begin{itemize}
\item $s_0: N(\D)_{01}\ra N(\D)_{11}$ and $\varsigma_0 : N(\D)_{10}\ra N(\D)_{11}.$ These are defined by $s_0 f:=\ID_f$ and $\varsigma_0 p:=\ID_p$.
\item $s_0: N(\D)_{02}\ra N(\D)_{12}$ and $\varsigma_0 : N(\D)_{20}\ra N(\D)_{21}.$ These are defined by $$s_0 ((h,g,f),\beta)= [ (h,g,f\vbar \beta)  ,\  (h,g,f\vbar\beta) \ \ | \ \ \ID_h ,\ \ID_g ,\ \ID_f ]$$ and $$\varsigma_0(r,q,p\vbar\eta)=  [ \ID_r ,\ \ID_q ,\ \ID_p  \ \  | \ \ (r,q,p\vbar \eta)  ,\  (r,q,p\vbar\eta)].  $$ Then the $(1,2)$-cell condition \ref{12condition} for $s_0 (h,g,f\vbar \beta)$ follows from the compatibility of the pseudo-identity with composition and with the action of $\beta$: $$(\ID_h \boxvert  \ID_f)\uact \beta \dact \beta=\ID_{h\circ f}\uact \beta \dact \beta=\ID_g.$$  The $(2,1)$-cell condition for $\varsigma_0(r,q,p\vbar \eta)$ follows similarly.
\item $s_0,s_1: N(\D)_{11}\ra N(\D)_{21}$ and $\varsigma_0,\varsigma_1 : N(\D)_{11}\ra N(\D)_{12}.$ These are defined for $\Theta \in \mbox{Sq}(g',g,q',q)$ by 
\begin{align*}s_0\Theta&=[\Theta,\ \Theta,\ \ID_g \ \  |  \  \    (\rho_{q'})_v,\ (\rho_q)_v  ]   \\   
 s_1\Theta&=[\ID_{g'},\ \Theta,\ \Theta \ \  |  \  \    (\lambda_{q'})_v,\ (\lambda_q)_v  ]  \\
  \varsigma_0\Theta&=[ (\rho_{g'})_h,\ (\rho_g)_h\ \  |  \  \  \Theta,\ \Theta,\ \ID_q ]  \\
   \varsigma_1\Theta&=[ (\lambda_{g'})_h,\ (\lambda_g)_h \ \  |  \  \ \ID_{q'},\ \Theta,\ \Theta  ].
      \end{align*}
      The $(2,1)$-cell conditions for $s_0 \Theta$ and $s_1\Theta$ follow directly from the compatibility of the unitors of $V$ with vertical composition, and the $(1,2)$-cell conditions for $\varsigma_0\Theta$ and $\varsigma_1\Theta$ follow from the compatibility of the unitors of $H$ with horizontal composition.
\end{itemize}

We have left to check that the bisimplicial identities hold for $N(\D)|^3_0.$ The identities involving face maps $d_i$ or $\delta_i$ can all be checked directly from our definitions. These verifications are left to the reader. We are left to check the following identities for $\varsigma_i\varsigma_j$ and $\varsigma_i s_j$ and $s_is_j:$
\begin{itemize}
\item  For $a\in N(\D)_{00}$,  $\varsigma_0 s_0 a = s_0\varsigma_0 a.$ Applying definitions, this is equivalent to $$\ID_{(\id_a)_v}=\ID_{(\id_a)_h},$$ which states the compatibility of $\ID$ with the identities for objects in $\D$.
\item For $f:a \ra b$  in $ N(\D)_{01}$, $\varsigma_0s_0 f=s_0\varsigma_0 f$ and $\varsigma_1s_0 f=s_0\varsigma_1 f.$ This follows directly from the definitions, \begin{align*}\varsigma_0s_0 f=s_0\varsigma_0 f&=[ (\rho_{f})_h,\ (\rho_f)_h\ \  |  \  \  \ID_f,\ \ID_f,\ \ID_{\id_a} ]   \\    \varsigma_1s_0 f=s_0\varsigma_1 f&=[ (\lambda_{f})_h,\ (\lambda_f)_h \ \  |  \  \ \ID_{\id_{b}},\ \ID_f,\ \ID_f  ]  \end{align*}
\item For $p:a \ra b$  in $ N(\D)_{10}$, $s_0\varsigma_0 p=\varsigma_0s_0 p$ and $s_1\varsigma_0 p=\varsigma_0s_1 p.$ This follows from the definitions as above.
\item For $f:a\ra b$  in  $ N(\D)_{01}$, $s_0s_0f=s_1s_0f.$ Applying the definitions, we get
\begin{align*} 
s_0s_0f=s_0\ID_f &=[\ID_f,\ \ID_f,\ \ID_f \ \  |  \  \    (\rho_{\id_b})_v,\ (\rho_{\id_a})_v \\ 
s_1s_0=s_1\ID_f &= [\ID_f,\ \ID_f,\ \ID_f \ \  |  \  \    (\lambda_{\id_b})_v,\ (\rho_{\id_a})_v  ]
\end{align*} 
The equality follows from bicategory axiom $12$ for $V$ which we called ``compatibility of the unitors with the identity for objects'' in Chapter~\ref{bicchapter}, which states $(\rho_{\id_a})_v=(\lambda_{\id_a})_v$ for all objects $a$ of $V$.
\item For $p:a\ra b$ in   $ N(\D)_{10}$, $\varsigma_0\varsigma_0p=\varsigma_1\varsigma_0p.$ This follows from the compatibility of the unitors of $H$ with the identity for objects in $H.$
\end{itemize} 

This completes the definition of $N(\D)|^3_0.$ We define $N(\D):= \mathbf{cosk}^3(N(\D)|^3_0 ) .$

\subsection{Horn-filling conditions for $N(\D).$}
The following lemma expresses a kind of ``associativity'' for actions and square composition that is not directly assumed by the VDC axioms, but follows from the other axioms.
\begin{lemma}[middle associativity] \label{midass} Let $\D$ be a small VDC, $\Theta, \Pi$ be squares of $\D,$ and $\theta,\pi$ be vertical $2$-morphisms as pictured below:
\begin{center}
\begin{tikzpicture}[scale=2.0,auto]

\begin{scope}

\node (10) at (1,1) {$b$};
\node (00) at (0,1) {$a$};
\node (11) at (1,0) {$b'$};
\node (01) at (0,0) {$a'$};
\node[rotate=-45] at (.35,.5){$\Rightarrow$};
\node[scale=.8] at (.45,.65){$\Theta$};

\node (20) at (2,1) {$c$};
\node (21) at (2,0) {$c'$};
\node[rotate=-45] at (1.5,.5){$\Rightarrow$};
\node[scale=.8] at (1.6,.65){$\Pi$};

\path[->] (00) edge node[midway]{$g$}(10);
\path[->] (00) edge node[midway,right]{}(01);
\path[->] (01) edge node[midway,swap]{$g'$}(11);
\path[->] (10) edge node[midway,shift={(0,-.3)}]{$p'$}(11);
\path[->] (10) edge node[midway]{$h$}(20);
\path[->] (11) edge node[midway,swap]{$h'$}(21);
\path[->] (20) edge node[midway,left]{}(21);
\draw[->] (10) to [out=-150,in=150] node(f)[midway,left,shift={(.1,-.3)}]{$q'$} (11) ;
\node[rotate=0] (areta) at ($(f)+(.25,.15)$){$\Leftarrow$};
\node (eta) at ($(areta)+(0,.15)$){$\theta$};

\draw[->] (10) to [out=-30,in=30] node(f)[midway,right]{} (11) ;
\node[rotate=180] (areta) at ($(f)+(-.25,0)$){$\Leftarrow$};
\node (eta) at ($(areta)+(0,.15)$){$\pi$};

\end{scope}
\end{tikzpicture}
\end{center}
 Then $$(\Pi\lact \pi) \boxvert (\Theta \ract \theta)= (\Pi\lact (\pi \bullet \theta^{-1}))\boxvert \Theta=\Pi \boxvert  ( \Theta \ract (\pi^{-1} \bullet \theta)).$$ A similar statement holds for vertical composition.
\end{lemma}
\begin{proof}~~

\begin{tabular}{p{10cm} p{5cm}}
$(\Pi\lact \pi) \boxvert (\Theta \ract \theta)$                                       &  \\ 
$=(\Pi \lact \pi) \boxvert [ ((\ID_{q'}\boxvert\Theta)\uact\lambda_g\dact\lambda_{g'})\ract \theta ]$ & (comp. of $\boxvert$ with unitors) \\ 
$=(\Pi \lact \pi) \boxvert[ ((\ID_{q'}\boxvert\Theta)\ract \theta )\uact\lambda_g\dact\lambda_{g'}]$ & (commutativity of actions)\\
$=(\Pi \lact \pi) \boxvert[ ((\ID_{q'}\ract \theta)\boxvert\Theta) \uact\lambda_g\dact\lambda_{g'}]$ & (compatibility of $\ract$ and $\boxvert$)) \\
$=[(\Pi \lact \pi) \boxvert ((\ID_{q'}\ract \theta)\boxvert\Theta)] \uact(h\rhd \lambda_g)\dact (h' \rhd \lambda_{g'})$ & (interchange)\\
$=[(\Pi \lact \pi) \boxvert ((\ID_{p'}  \lact \theta^{-1})\boxvert\Theta)] \uact(h\rhd \lambda_g)\dact (h' \rhd \lambda_{g'})$ & (comp. of actions with $\ID$)\\
$=[[((\Pi \lact \pi) \boxvert ( \ID_{p'}\lact \theta^{-1} ))\boxvert\Theta]  \uact \alpha_{h,\id_b,g} \dact \alpha_{h',\id_{b'},g'} ]$ &  (comp. of $\boxvert$ with $\alpha$)\\  $\uact(h\rhd \lambda_g)\dact (h' \rhd \lambda_{g'})$ & \\
$=[((\Pi \lact \pi) \boxvert ( \ID_{p'}\lact \theta^{-1} ))\boxvert\Theta]  \uact [\alpha_{g,\id_b,h}\bullet (h\rhd \lambda_g)]$ &  (associativity of $\uact$ and $\dact$)\\  $\dact[\alpha_{h',\id_{b'},g'}\bullet (h' \rhd \lambda_{g'})]$ &\\
$=[((\Pi \lact \pi) \boxvert ( \ID_{p'}\lact \theta^{-1} ))\boxvert\Theta]  \uact (\rho_h\lhd g)\dact (\rho_{h'}\lhd g')$ &  (comp. of $\alpha$ and the unitors)\\  
$=[(((\Pi \lact \pi) \boxvert  \ID_{p'})\lact \theta^{-1} )\boxvert\Theta]  \uact (\rho_h\lhd g)\dact (\rho_{h'}\lhd g')$ &  (comp. of $\lact$ and $\boxvert$)\\ 
$=[(((\Pi \lact \pi) \boxvert  \ID_{p'} )\lact \theta^{-1} ) \uact \rho_h\dact \rho_{h'}]\boxvert\Theta$ &  (interchange)\\   
$=[(((\Pi \lact \pi) \boxvert  \ID_{p'} ) \uact \rho_h\dact \rho_{h'}) \lact \theta^{-1}]\boxvert\Theta$ &  (commutativity of actions))\\   
$=[(\Pi \lact \pi)  \lact \theta^{-1}]\boxvert\Theta$ & (comp. of $\ID$ and $\rho$) \\
$=(\Pi \lact (\pi \bullet \theta^{-1}))\boxvert \Theta$ & (associativity of $\lact$) \\
\end{tabular}

The equality  $(\Pi\lact \pi) \boxvert (\Theta \ract \theta)=\Pi \boxvert  ( \Theta \ract (\pi^{-1} \bullet \theta))$ follows by a similar argument, and the statement for vertical composition follows by symmetry.
\end{proof}
\begin{corollary}[middle cancellation] \label{midcancel} With $\Theta,\Pi, \theta$ as in Lemma~\ref{midass} and $\pi=\theta$, $$(\Pi\lact \theta) \boxvert (\Theta \ract \theta)= \Pi\boxvert \Theta.$$ Similarly for vertical composition.
\end{corollary}

 \begin{theorem}\label{horntheorem} If $\D$ is a small VDC, $N(\D)$ is $2$-reduced inner-Kan.
 \end{theorem}
 \begin{proof}
 Using Lemma~\ref{coskeletaliskan}, it suffices to check that (truncated) horns of dimension $\leq 4$ in $N(\D)|^3_0$ have the appropriate filling conditions. The conditions for horns of type $\Lambda[0,k^{\ i}]$ and $\Lambda[k^{ \ i},0]$ follow from the fact that $N(H)$ and $N(V)$ are $2$-reduced inner-Kan, which is part of the main theorem of Chapter~\ref{bicchapter}, Theorem~\ref{bictheorem}. 
 
For the remaining horns, first consider horns of type $\Lambda[ 1,\ 2^{\ 1}]$ and $\Lambda[2^{\ 1},\ 1].$  A horn of type $\Lambda[ 1,\ 2^{\ 1}]$ has the form:
 $$d(x)=[(h',g',f'\vbar \beta'),\  (h,g,f\vbar \beta)  \ \    |   \  \   \Sigma    ,\  -    ,\ \Theta  ].$$ A filler of this horn is a sphere $$d(x)=[(h',g',f'\vbar \beta'),\  (h,g,f\vbar \beta )  \ \    |   \  \   \Sigma    ,\  \Pi    ,\ \Theta  ]$$ satisfying Condition~\ref{12condition}, $$(\Sigma\boxvert\Theta)\uact \beta \dact \beta'=\Pi,$$ so obviously this filler exists and is unique. Similarly for horns of the type  $\Lambda[2^{\ 1},\ 1].$
 
Next we consider a horn of type $\Lambda[2,\ 2^{\ 1}]$ in $N(\D)|^3_0:$ 
\begin{table}[H]\begin{center}
    \begin{tabular}{ r | l || l | l !{\vrule width 2pt} l | l | l |}
    \cline{2-7}

      $ \quad$ &  $x^{12}_{012}$  &  $x^{2}_{012}$  &   $x^{1}_{012}$   &     $x^{12}_{12} $    &    $x^{12}_{02}$   &  $x^{12}_{01} $        \\ \cline{2-7}
     
        &  $x^{02}_{012}$  &  $x^{2}_{012}$  &   $x^{0}_{012}$   &     $x^{02}_{12} $    &    $x^{02}_{02}$   &  $x^{02}_{01} $       \\ \cline{2-7}

        &  $x^{01}_{012}$  &  $x^{1}_{012}$  &   $x^{0}_{012}$   &     $x^{01}_{12} $    &    $x^{01}_{02}$   &  $x^{01}_{01} $         \\ \cline{2-7} 
          
    \end{tabular} 
\vspace{12pt}

 \begin{tabular}{    r | l|| l | l |  l !{\vrule width 2pt} l | l  | }\cline{2-7} 
 &  $x^{012}_{12}$ &   $x^{12}_{12}$  &   $x^{02}_{12}$   &     $x^{01}_{12} $    &    $x^{012}_{2}$   &  $x^{012}_{1} $     \\ \cline{2-7} 
$\Lambda$  &  &   $x^{12}_{02}$  &   $x^{02}_{02}$   &     $x^{01}_{02} $    &    $x^{012}_{2}$   &  $x^{012}_{0} $   \\ \cline{2-7}
&  $x^{012}_{01}$  &   $x^{12}_{01}$  &   $x^{02}_{01}$   &     $x^{01}_{01} $    &    $x^{012}_{1}$   &  $x^{012}_{0} $     \\ \cline{2-7}
\end{tabular}   
  \end{center} \end{table}
We must check that the $(2,1)$-cell condition \ref{21condition} holds if the $(2,1)$-cell and $(1,2)$-cell conditions hold for all the other faces.  That is we must show 
\begin{equation} 
x^{02}_{02} = (x^{12}_{02}\boxvert x^{01}_{02})\uact x^{012}_0 \dact x^{012}_2 \label{sphereface13}
\end{equation}
given the following:
\begin{align}
x^{02}_{12} &= (x^{12}_{12}\boxvert x^{01}_{12})\uact x^{012}_1 \dact x^{012}_2 \label{sphereface14}\\
x^{02}_{01} &= (x^{12}_{01}\boxvert x^{01}_{01})\uact x^{012}_0 \dact x^{012}_1 \label{sphereface15}
\end{align}
\begin{align}
x^{12}_{02} &= (x^{12}_{12}\boxminus x^{12}_{01})\lact x^1_{012} \ract x^2_{012} \label{sphereface16}\\
x^{02}_{02} &= (x^{02}_{12}\boxminus x^{02}_{01})\lact x^0_{012} \ract x^2_{012} \label{sphereface17}\\
x^{01}_{02} &= (x^{01}_{12}\boxminus x^{01}_{01})\lact x^0_{012} \ract x^1_{012} \label{sphereface18} .
\end{align}

We begin by substituting from Equations~\ref{sphereface14} and \ref{sphereface15} into Equation~\ref{sphereface17}:

$$x^{02}_{02} = ([ (x^{12}_{12}\boxvert x^{01}_{12})\uact x^{012}_1 \dact x^{012}_2 ]\boxminus [(x^{12}_{01}\boxvert x^{01}_{01})\uact x^{012}_0 \dact x^{012}_1 ])\lact x^0_{012} \ract x^2_{012} $$
Next apply commutativity of $\uact$ and $\dact$ and then the vertical case of middle cancellation (Corollary~\ref{midcancel}) to cancel $x^{012}_1,$ yielding:
$$x^{02}_{02} = ([ (x^{12}_{12}\boxvert x^{01}_{12}) \dact x^{012}_2 ]\boxminus [(x^{12}_{01}\boxvert x^{01}_{01})\uact x^{012}_0  ])\lact x^0_{012} \ract x^2_{012} $$
Apply the compatability of $\uact$ and $\dact$ with $\boxminus$ to give:
$$x^{02}_{02} = [ (x^{12}_{12}\boxvert x^{01}_{12})\boxminus (x^{12}_{01}\boxvert x^{01}_{01})  ]\uact x^{012}_0\dact x^{012}_2 \lact x^0_{012} \ract x^2_{012}. $$
Square interchange yields:
$$x^{02}_{02} = [ (x^{12}_{12}\boxminus x^{12}_{01} )\boxvert (x^{01}_{12}\boxminus x^{01}_{01})  ]\uact x^{012}_0\dact x^{012}_2 \lact x^0_{012} \ract x^2_{012}. $$
Apply commutativity of actions and the compatability of  $\lact$ and $\ract$ with $\boxvert$ to get
$$x^{02}_{02} = [ ((x^{12}_{12}\boxminus x^{12}_{01}) \ract x^2_{012})\boxvert ((x^{01}_{12}\boxminus x^{01}_{01}) \lact x^0_{012})  ]\uact x^{012}_0\dact x^{012}_2.$$
Apply the horizontal case of middle cancellation (Corollary~\ref{midcancel}) to add in two instances of $x^1_{012}$:
$$x^{02}_{02} = [ ((x^{12}_{12}\boxminus x^{12}_{01}) \ract x^2_{012}\lact x^1_{012})\boxvert ((x^{01}_{12}\boxminus x^{01}_{01}) \lact x^0_{012}\ract x^1_{012} )  ]\uact x^{012}_0\dact x^{012}_2.$$
Finally commute the $\ract$ and $\lact$ on the left then substitute using Equations~\ref{sphereface16} and \ref{sphereface18}, giving:
$$x^{02}_{02}= (x^{12}_{02} \boxvert x^{01}_{02})\uact x^{012}_0\dact x^{012}_2,$$
which was to be shown.  The $\Lambda[2^{\ 1},\ 2]$ case is symmetric.

Next consider a horn of type $\Lambda[1, 3^{\ 1}]$: 
\begin{table}[H]\begin{center}
    \begin{tabular}{ r | l ||!{\vrule width 2pt} l | l | l | l |}
    \cline{2-6} 
   $\quad$  &  $x_{0123}^1$  &$x_{123}^1$   & $x_{023}^1$     &  $x_{013}^1$                  & $x_{012}^1 $          \\ \cline{2-6}
     
     &  $x_{0123}^0$         &$x_{123}^0$      &  $x_{023}^0$      & $x_{013}^0$                    &  $x_{012}^0$             \\ \cline{2-6}      
    \end{tabular} 

\vspace{12pt}

 \begin{tabular}{    r | l || l | l  !{\vrule width 2pt}  l | l | l | }\cline{2-7} 
                  &  $x_{123}^{01}$   &  $x_{123}^1$               &   $x_{123}^0$                &    $x_{23}^{01}$        &  $x_{13}^{01}$    &   $x_{12}^{01}$       \\ \cline{2-7} 
               
   $\Lambda $      & & $x_{023}^{1}$      &   $x_{023}^0$                             &    $x_{23}^{01}$     &    $x_{03}^{01}$     &     $x_{02}^{01}$       \\ \cline{2-7}
       
                    &     $x_{013}^{01}$  & $x_{013}^{1}$     &   $x_{013}^{0}$                                         &  $x_{13}^{01}$        &$x_{03}^{01}$   &      $x_{01}^{01}$           \\ \cline{2-7}    
          
                    &  $x_{012}^{01}$   &   $x_{012}^1 $         &   $x_{012}^0$                                     &  $x_{12}^{01}$       &    $x_{02}^{01}$    &   $x_{01}^{01}$  \\ \cline{2-7}
\end{tabular}   
\end{center}\end{table}
We must show that the $(1,2)$-cell condition holds for vertical face $1$, marked $\Lambda$, if the $(1,2)$-cell and $3$-cell conditions hold for the other faces. Explicitly, we must show
\begin{equation}x^{01}_{03} = (x^{01}_{23}\boxvert x^{01}_{02}) \uact x^0_{023}\dact x^1_{023}\label{sphereface19} \end{equation} given:
\begin{align}
\alpha_{x^1_{23},x^1_{12},x^1_{01}}\bullet(x^1_{23}\rhd x^1_{012})\bullet x^1_{023}&=(x^1_{123}\lhd x^1_{01})\bullet x^1_{013} \label{sphereface20}\\
 \alpha_{x^0_{23},x^0_{12},x^0_{01}}\bullet(x^0_{23}\rhd x^0_{012})\bullet x^0_{023}&=(x^0_{123}\lhd x^0_{01})\bullet x^0_{013} \label{sphereface21}\\
x_{13}^{01}  &= (x_{23}^{01} \boxvert x_{12}^{01}) \uact x_{123}^0 \dact x_{123}^1 \label{sphereface22} \\
x^{01}_{03}  &= (x^{01}_{13}\boxvert x^{01}_{01})   \uact x^0_{013}\dact x^1_{013} \label{sphereface23} \\
x_{02}^{01}  &= (x_{12}^{01} \boxvert x_{01}^{01}) \uact x_{012}^0 \dact x_{012}^1 \label{sphereface24}
\end{align}
Similarly for the $\Lambda[1, 3^{\ 2}]$ condition we must show Equation \ref{sphereface23} follows from the other equations above. So we can show both horn-filling conditions by showing Equations \ref{sphereface19} and \ref{sphereface23} are equivalent, given the other four equations. Beginning with Equation~\ref{sphereface23}:

$$x^{01}_{03}  = (x^{01}_{13}\boxvert x^{01}_{01})   \uact x^0_{013}\dact x^1_{013}. $$
Substitute for $x^{01}_{13}$ from Equation~\ref{sphereface22}:
$$x^{01}_{03}  = [((x_{23}^{01} \boxvert x_{12}^{01}) \uact x_{123}^0 \dact x_{123}^1) \boxvert x^{01}_{01} ]   \uact x^0_{013}\dact x^1_{013}.$$
Apply interchange of $\uact$ and $\dact$ with whiskering:
$$x^{01}_{03}  = [((x_{23}^{01} \boxvert x_{12}^{01})  \boxvert x^{01}_{01})\uact (x_{123}^0\lhd x^0_{01} ) \dact (x_{123}^1 \lhd x^1_{01})]   \uact x^0_{013}\dact x^1_{013}. $$
Apply commutativity and associativity of actions:
$$x^{01}_{03}  = ((x_{23}^{01} \boxvert x_{12}^{01})  \boxvert x^{01}_{01})\uact ((x_{123}^0\lhd x^0_{01} )\bullet   x^0_{013}   )\dact ((x_{123}^1 \lhd x^1_{01})\bullet x^1_{013} ) .$$
Substitute from \ref{sphereface20} and \ref{sphereface21}:
$$x^{01}_{03}  = ((x_{23}^{01} \boxvert x_{12}^{01})  \boxvert x^{01}_{01})\uact ( \alpha_{x^0_{23},x^0_{12},x^0_{01}}\bullet(x^0_{23}\rhd x^0_{012})\bullet x^0_{023} )\dact (\alpha_{x^1_{23},x^1_{12},x^1_{01}}\bullet(x^1_{23}\rhd x^1_{012})\bullet x^1_{023}) .$$
Apply associativity of $\uact$ and $\dact$ then the compatibility of $\boxvert$ with the associator to get:
$$x^{01}_{03}  = (x_{23}^{01} \boxvert (x_{12}^{01}  \boxvert x^{01}_{01}))\uact ((x^0_{23}\rhd x^0_{012})\bullet x^0_{023} )\dact ((x^1_{23}\rhd x^1_{012})\bullet x^1_{023}) .$$
Apply the associativity of $\uact$ and $\dact$ and interchange:
$$x^{01}_{03}  = [x_{23}^{01} \boxvert (( x_{12}^{01}  \boxvert x^{01}_{01}  )\uact x^0_{012} \dact x^1_{012} )]\uact  x^0_{023} \dact x^1_{023}.$$
Substitute from \ref{sphereface24}:
$$x^{01}_{03}  = (x_{23}^{01} \boxvert x_{02}^{01} )\uact  x^0_{023} \dact  x^1_{023}.$$
This is Equation~\ref{sphereface19}. Each step is an equivalence, so we have shown the horn filling conditions for $\Lambda[1, 3^{\ 1}]$ and $\Lambda[1, 3^{\ 2}]$. The conditions for $\Lambda[ 3^{\ 1},1]$ and $\Lambda[3^{\ 2},1]$ follow by symmetry. 
\end{proof}
\subsection{The isomorphism $u:X\ra N(\mbox{Vdc}(X))$}
Let $X$ be a $2$-reduced inner-Kan bisimplicial set.
Since both $X$ and $N(\mbox{Vdc}(X))$ are $3$-coskeletal, it is enough to give an isomorphism $X|^3_0 \ra N(\mbox{Vdc}(X))|^3_0.$
Recall the horizontal and vertical bicategories of $\mbox{Vdc}(X)$ are given by $\Bic(X_h)$ and $\Bic(X_v),$ and the first row and column of $N(\D)$ are given by the Duskin nerve of the horizontal and vertical bicategories of $D$. So the first row and column of the isomorphism $u:X|^3_0\ra N(\mbox{Vdc}(X))|^3_0$ are given by the isomorphism $u$ constructed in Section~\ref{bicisosec}. Recall that an element of $N(\D)_{02}$ is a triple of horizontal $1$-morphisms $(h,g,f)$ together with a horizontal $2$-morphism $\beta: g \Rightarrow h\circ f$.  So for $x_{012}\in X_{02}$ we have $u(x_{012})= (x_{12},x_{02},x_{01}\vbar   \underline{x_{012}}).$ 

By definition, $X_{11}= N(\mbox{Vdc}(X))_{11}$ so $u(x)=x$ for $x\in X_{11}$. For an element $x^{01}_{012}\in X_{12}$ we define $u(x^{01}_{012})$ to be:
 \begin{equation} \label{usphere} [  u(x^1_{012})= (x^1_{12},x^1_{02},x^1_{01} \vbar  \underline{x^1_{012}}),\ \  u(x^0_{012})=(x^0_{12},x^0_{02},x^0_{01}\vbar  \underline{x^0_{012}}) \ \ | \ \  x^{01}_{12},\ x^{01}_{02},\ x^{01}_{01} ]. \end{equation} 
This map is injective by Lemma~\ref{kaniscoskeletal} and the injectivity of $u$ on $X_{02}$. We must show that it is well-defined and that it is surjective. Since $u$ is also surjective, this amounts to showing \begin{equation}\label{12sphere}[ x^1_{012},\ x^0_{012}  \ \ | \ \  x^{01}_{12},\ x^{01}_{02},\ x^{01}_{01}  ]\end{equation} is commutative if and only if \ref{usphere} meets the $(1,2)$-cell condition:
\begin{equation}  x^{01}_{02} =( x^{01}_{12}\boxvert x^{01}_{01} ) \uact  \underline{x^0_{012}} \dact \underline{x^1_{012}}
\end{equation}
which is equivalent to the commutativity of the following sphere in X:
\begin{equation} \label{showisosphere}  [   \underline{x^1_{012}},\ \underline{x^0_{012}} \ \ | \ \   x^{01}_{12}\boxvert x^{01}_{01},\ x^{01}_{02},\   \ID_{x^{01}_0} ].\end{equation}
Consider the following horn in $X$:

\begin{table}[H]\begin{center}\caption{~}\label{uiso}
    \begin{tabular}{ r | l ||!{\vrule width 2pt} l | l | l | l |}
    \cline{2-6} 
     &  $x_{0123}^1$  &$x_{123}^1$   & $x_{023}^1$     &  $x_{013}^1$                  & $x_{012}^1 $          \\ \cline{2-6}
     
     &  $x_{0123}^0$         &$x_{123}^0$      &  $x_{023}^0$      & $x_{013}^0$                    &  $x_{012}^0$             \\ \cline{2-6}      
    \end{tabular} 

\vspace{12pt}

 \begin{tabular}{    r | l || l | l  !{\vrule width 2pt}  l | l | l | }\cline{2-7} 
                 &  $\Delta_{\boxvert}(x^{01}_{12}, x^{01}_{01})$    &  $\chi(x^1_{12}, x^1_{01})$               &   $\chi(x^0_{12}, x^0_{01})$                &    $x_{12}^{01}$        &  $x_{12}^{01}\boxvert x^{01}_{01}$    &   $ x^{01}_{01}$       \\ \cline{2-7} 
               
   $\Lambda $       & (sphere~\ref{12sphere}) & $x_{012}^{1}$      &   $x_{012}^0$                             &    $x_{12}^{01}$     &    $x_{02}^{01}$     &     $x_{01}^{01}$     \\ \cline{2-7}
       
      $\odot$           &     (sphere~\ref{showisosphere})       & $ \underline{x^1_{012}}$     &   $ \underline{x^0_{012}}$      &  $x^{01}_{12}\boxvert x^{01}_{01}$        &$x_{01}^{02}$   &      $\ID_{x^{01}_0}$        \\ \cline{2-7}    
          
                  &  $\varsigma_0 (x^{01}_{01}) $  &   $\ID_{x_{01}^1}$         &   $\ID_{x_{01}^0}$                            &  $x_{01}^{01}$       &    $x_{01}^{01}$    &  $\ID_{x^{01}_0}$       \\ \cline{2-7}
\end{tabular} \end{center}\end{table}  
As written, the horn in Table~\ref{uiso} shows the commutativity of sphere~\ref{12sphere} given the commutativity of sphere~\ref{showisosphere}. If we instead take the commutativity of  sphere~\ref{12sphere} as a hypothesis, then we can view the above table as a horn of type $\Lambda[1,3^{\ 2}]$, (switch the symbols $\odot$ and $\Lambda$), and this horn then verifies the commutativity of sphere~\ref{showisosphere}. So the commutativity of these spheres is equivalent, which was to be shown.

With $u:X_{21}\ra N(\mbox{Vdc}(X))_{21}$ constructed symmetrically, we have a canonical isomorphism $u:X|^3_0 \ra N(\mbox{Vdc}(X))|^3_0$, and applying $\mathbf{cosk}^3$ we get a canonical isomorphism  $u:X \ra N(\mbox{Vdc}(X)).$ It is trivial to infer that this isomorphism is natural from the fact that the isomorphism $u:X\cong N(\Bic(X))$ constructed in Chapter~\ref{bicchapter} is natural.
\begin{theorem}Let $X$ be a bisimplicial set. Then $X$ is isomorphic to the nerve of a small VDC if and only if $X$ is $2$-reduced inner-Kan.
\end{theorem}
\begin{proof} The nerve of a small VDC is $2$-reduced inner-Kan by Theorem~\ref{horntheorem}, and if $X$ is $2$-reduced inner-Kan, then $X$ is isomorphic to $N(\mbox{Vdc}(X))$ by the isomorphism constructed in this section.
\end{proof}
\section{Promoting $N$ and $\mbox{Vdc}$ to functors}
First we define the notion of a functor of Verity double categories:
\begin{definition} Let $\D$, $\D'$ be Verity double categories. Then a functor $(F, \phi^v,\upsilon^v,\phi^h,\upsilon^h):\D\ra \D'$ consists of
\begin{itemize}
\item A map $F$ from each set of data of $\D$, (objects, vertical/horizontal $1$-morphisms, vertical/horizontal $2$-morphisms, and squares) to the corresponding set of data of $\D'$.
\item  A \emph{horizontal distributor} $\phi^h$ and \emph{vertical unitor}  $ \upsilon^h$, such that the restriction $F|_{\mathcal{H}}:\mathcal{H} \ra \mathcal{H'}$ together with $\phi^h$ and $\upsilon^h$ is a functor of $(2,1)$-categories. 
\item A \emph{vertical distributor} $\phi^v$ and \emph{vertical unitor}  $ \upsilon^v$, such that the restriction $F|_{\mathcal{V}}:\mathcal{V} \ra \mathcal{V'}$ together with $\phi^v$ and $\upsilon^v$ is a functor of $(2,1)$-categories. 
\end{itemize} 
In addition to the axioms that $(F|_{\mathcal{H}},\phi^h,\upsilon^h)$ and  $(F|_{\mathcal{V}},\phi^v,\upsilon^v)$ satisfy by virtue of being $(2,1)$-category functors, $F$ satisfies the following axioms:
\begin{description}
\item[VFun1.] For squares $\Theta, \Pi$ in $\D$ as shown below, $[F(\Pi) \boxvert F(\Theta)] \uact \phi^h_{g,f}\dact \phi^h_{g',f'}= F(\Pi \boxvert \Theta).$ A similar axiom holds for $\boxminus$ and $\phi^v$
\item[VFun2.] For a horizontal $1$-morphism $f:a\ra b$,    $\ID_{F(f)}\lact \upsilon^v_a \ract \upsilon^v_b=F(\ID_f).$ A similar axiom holds for a vertical $1$-morphism. 
\item[VFun3.] For a square $\Theta$ and a morphism $\eta$ such that $\Theta \lact \eta$ is defined, then $F(\Theta \lact \eta)= F(\Theta)\lact F(\eta).$ Similarly for $\ract, \uact$ and $\dact$. 
\end{description}

Such a functor $F$ is called \emph{strict} if $F|_{\mathcal{H}}$ and $F|_{\mathcal{V}}$ are strict, and \emph{strictly identity-preserving} if $F|_{\mathcal{H}}$ and $F|_{\mathcal{V}}$ are strictly identity-preserving. Note that \textbf{VFun2} ensures that $F$ strictly preserves pseudo-idenity squares if $F$ is strictly identity-preserving.
\end{definition}

We now describe how $N$ and $\mbox{Vdc}$ can be applied respectively to functors of VDC's and to morphisms of bisimplicial sets. \begin{definition}Let $F:X\ra Y$ be a map between  algebraic $2$-reduced inner-Kan bisimplical sets. Define $\mbox{Vdc}(F): \mbox{Vdc}(X) \ra \mbox{Vdc}(Y)$ as follows:
\begin{itemize}
\item $\mbox{Vdc}(F)|_{\mathcal{H}}$ is defined by $\Bic(F|_{\mathcal{H}})$ and $\mbox{Vdc}(F)|_{\mathcal{V}}$ by $\Bic(F|_{\mathcal{V}}).$
\item For a square $S$ in $\mbox{Vdc}(X)$ we take $\mbox{Vdc}(F)(S)=F(S),$ which is an element in $Y_{11}$ and thus a square of $\mbox{Vdc}(Y).$
\end{itemize}\end{definition}
\begin{proposition} $\mbox{Vdc}(F)$ is a strictly identity-preserving functor of VDC's.
\end{proposition}\begin{proof}
We must check the axioms for a functor of VDC's. 
\begin{itemize}
\item For \textbf{VDC1} let $\Theta$ and $\Pi$ be as in Definition~\ref{horzcompdef}. Then we have  $$d F(\Delta_{\boxvert}(\Pi,\Theta))=[F(\pi),\ F(\Pi \boxvert \Theta), F(\Theta) \ \ | \ \ \chi(r',q'),\ \chi(r,q)].$$ The commutativity condition Equation~\ref{12condition} for this cell is  $$F(\Pi \boxvert \Theta)=F(\Pi) \boxvert F(\Theta) \uact \underline{F(\chi(r,q))} \dact \underline{F(\chi(r',q'))} =F(\Pi) \boxvert F(\Theta) \uact \phi_{r,q} \dact \phi_{r',q'},$$ which gives one case of \textbf{VDC1}, the other case following by symmetry.
\item \textbf{VDC2} follows from the fact that $F|_{\mathcal{H}}$ and $F|_{\mathcal{V}}$ are strictly identity-preserving and $F$ preserves pseudo-identities squares. Note that this also assures that $F$ is strictly identity-preserving.
\item For \textbf{VDC3} let $\Theta$ and $\eta$ be as in Definition~\ref{lactdef}. We have: \begin{align*}dF(\Delta_{\lact}(\Theta,\eta)) &= [F(\Theta),\ F(\Theta \lact \eta),\ F(\ID_g) \ \ | \ \ F(\Id_{q'}),\ F(\eta)]  \\  &=   [F(\Theta),\ F(\Theta \lact \eta),\ \ID_{F(g)} \ \ | \ \ \Id_{F(q')},\ F(\eta)] \\ \\ d\left( \Delta_{\lact}(F(\Theta),F(\eta))\right)&=   [F(\Theta),\ F(\Theta) \lact F(\eta),\ \ID_{F(g)} \ \ | \ \ \Id_{F(q')},\ F(\eta)].\end{align*} By the Matching Lemma, we conclude $F(\Theta \lact \eta)= F(\Theta)\lact F(\eta).$ The other cases of \textbf{VFun3} follow by a similar argument. \end{itemize}\end{proof}

We now turn to defining the nerve of a VDC functor.
\begin{definition} Let $F:\D \ra \D'$ be a strictly identity-preserving functor of small VDC's. We first define $N(F):N(\D)|^3_0 \ra N(\D')|^3_0.$ For the edges ($0$th row and column) of $N(\D)|^3_0$ we define $N(F)|^3_0$ using $N(F|_{\Hcal})$ and $N(F|_{\Vcal})$ as defined in Chapter~\ref{bicchapter}. For elements of $N(\D)_{11}$, $N(F)|^3_0$ is defined by $F$, the since the elements of $N(\D)_{11}$ correspond to the squares of $\D$. An element of $N(\D)_{12}$ is given by a 5-tuple $$x^{01}_{012}=[x^1_{012},\ x^0_{012} \ \ | \ \ x^{01}_{12}, \ x^{01}_{02}, \ x^{01}_{01}]$$ of cells of dimension $2$, meeting the 3-cell condition of Equation~\ref{12condition}. So we define $$N(F)|^3_0(x^{01}_{012})=[N(F)|^3_0 (x^1_{012}), \ N(F)|^3_0 (x^0_{012}) \ \ | \ \ N(F)|^3_0( x^{01}_{12}),\ N(F)|^3_0 (x^{01}_{02}),\ N(F)|^3_0 (x^{01}_{01})].$$ $N(F)|^3_0$ is defined in a similar manner for cells in $N(\D)_{21}.$ To show this is well defined, we must show that the 3-cell conditions of Equations~\ref{12condition} and \ref{21condition} for $N(F)|^3_0(x^{01}_{012})$ and $N(F)|^3_0 (x^{012}_{01})$ follow from the same conditions for $x^{01}_{012}$ and $x^{012}_{01}.$ This is shown in Proposition~\ref{3cellconditioncheck}
below. The fact that $N(F)|^3_0$ commutes with face maps is immediate from our definition. The nontrivial cases of commutativity with degeneracy maps follows from the fact that $F$ preserves pseudo-identity squares (because it is strictly identity-preserving). 
\end{definition}
\begin{proposition} \label{3cellconditioncheck} For $x^{01}_{012}\in N(\D)_{12}$, the 5-tuple $N(F)|^3_0(x^{01}_{012} )$ meets the 3-cell condition of Equation~\ref{12condition}, and thus is an element of $N(\D')_{12}$. Similarly for an element of $x^{012}_{01}\in N(\D)_{21}.$
\end{proposition}
\begin{proof} The $3$-cell condition of Equation~\ref{12condition} for $x^{01}_{012}$, we have \begin{equation}\label{12cellcondforx}(x^{01}_{12}\boxvert x^{01}_{01})\uact x^0_{012} \dact x^1_{012}= x^{01}_{02}.\end{equation} We wish to show the $3$-cell condition for $N(F)|^3_0(x^{01}_{012} ),$ which asserts $$(N(F)(x^{01}_{12})\boxvert N(F)(x^{01}_{01}))\uact N(F)(x^0_{012}) \dact N(F)(x^1_{012})= N(F)(x^{01}_{02}).$$ Applying the definition of $N(F)$, this is equivalent to: 
\begin{align*}(F(x^{01}_{12})\boxvert F(x^{01}_{01}))\uact (\phi_{x^0_{12},x^0_{01}} \bullet F(x^0_{012})) \dact (\phi_{x^1_{12},x^1_{01}}\bullet F(x^1_{012}))& = F(x^{01}_{02}) \\   
[(F(x^{01}_{12})\boxvert F(x^{01}_{01}))\uact (\phi_{x^0_{12},x^0_{01}}) \dact (\phi_{x^1_{12},x^1_{01}})]\uact F(x^0_{012}) \dact  F(x^1_{012})     &=F(x^{01}_{02})   \\
F(x^{01}_{12}\boxvert x^{01}_{01})\uact F(x^0_{012}) \dact  F(x^1_{012})     &=F(x^{01}_{02}) \\
F(x^{01}_{12}\boxvert x^{01}_{01}\uact x^0_{012} \dact  x^1_{012})           &=F(x^{01}_{02})
\end{align*} where we have applied \textbf{VFun1} and \textbf{VFun3}. The last equality follows by applying $F$ to both sides of Equation~\ref{12cellcondforx}. The corresponding statement for cells in $N(\D)_{21}$ follows by symmetry.
\end{proof}
\section{The natural isomorphism $\D\cong \mbox{Vdc}(N(\D)) $ and summary} 
In this section, we construct an isomorphism $$U:\D \stackrel{\cong}{\lra} \mbox{Vdc}(N(\D)) .$$ This isomorphism is given on the horizontal and vertical bicategories by the isomorphism $U: \B \Rightarrow \Bic(N(\B))$ as constructed in Section~\ref{summarysection}. On squares, $U$ acts as the identity. From this definition, and the fact that the $(2,1)$-categorical version of $U$ from Section~\ref{summarysection} is a strict natural isomorphism, it is easy to see that $U$ is strict natural isomorphism if it can be shown that $U$ satisfies the axioms for a functor of VDC's. 

To check the axioms \textbf{VFun1}-\textbf{VFun3}, we first describe the square composition and $2$-morphism actions for the $\mbox{Vdc}(N(\D))$, which we denote for instance by $\widetilde{\boxvert}$ and $\widetilde{\uact},$ in terms of the structure of $\D$. Note that the squares of $\mbox{Vdc}(N(\D))$ are identical to those of $\D$. For squares $\Pi$ and $\Theta$ as in the definition of $\boxvert$ we have $$\Pi \widetilde{\boxvert} \Theta = (\Pi \boxvert \Theta)\uact \Id_{g\circ f} \dact \Id_{g' \circ f'} = \Pi \boxvert \Theta.$$ Similarly $\widetilde{\boxminus}$ is identical to $\boxminus$. For $\Theta$ and $\beta$ as Figure~\ref{Pipic} in $\mbox{Vdc}(N(\D))$, recall that the horizontal morphism $\beta : f \Rightarrow g$ is given by a horizontal morphism $f \Rightarrow g \circ \id_a$ in $\D$. Unwrapping the definition of $\widetilde{\uact}$, we get the formula: $$\Theta \widetilde{\uact} \beta = (\Theta \boxvert \ID_q) \uact \beta \dact \rho_{g'}.$$
The formulas for the other actions in  $\mbox{Vdc}(N(\D))$ are similar (again using $\beta'$, $\eta$, and $\eta'$ as in Figure~\ref{Pipic})
\begin{align*}
\Theta \widetilde{\dact} \beta' &= (\Theta \boxvert \ID_q)  \uact \rho_{g}\dact \beta. \\
\Theta \widetilde{\lact} \eta &= (\Theta \boxminus \ID_g) \lact \eta \ract \rho_{q'}. \\
\Theta \widetilde{\ract} \eta' &= (\Theta \boxminus \ID_g) \lact  \rho_{q} \ract \eta'. \\							
\end{align*}

Since $U$ is taken to be strict and thus has trivial distributors and unitors, the fact that $U$ strictly preserves $\boxvert$ and $\boxminus$ as shown above verifies \textbf{VFun1}. Also $U$ can be easily be seen to preserve pseudo-identity squares directly from the definitions of $N$ and $\mbox{Vdc}$, showing \textbf{VFun2}. For \textbf{VFun3} we must show for instance: $$F(\Theta) \widetilde{\uact} F(\beta)=(\Theta \boxvert \Id_q) \uact (\rho_g \bullet \beta) \dact \rho_{g'}\stackrel{?}{=}F(\Theta \uact \beta)=\Theta \uact \beta.$$ Manipulating the left-hand side using \textbf{VDC2} and \textbf{VDC3} we get $$((\Theta \boxvert \Id_q) \uact \rho_g  \dact \rho_{g'} )\uact \beta$$ which is equal to the right-hand side $\Theta \uact \beta$ by \textbf{VDC10}. This verifies that $U$ is a functor of VDC's. 

\begin{theorem} \label{vdcsummary} The functors $N$, $\Vdc$ are inverse equivalences of categories between the category of algebraic $2$-reduced inner-Kan bisimplicial sets and the category of small VDC's and strictly identity-preserving functors. Furthermore, $N$ and $\Vdc$ preserve strictness, and the natural isomorphisms $u:N \Vdc\cong \Id $ and $U: \Vdc N \cong \Id$ exhibiting the equivalence are strict thus $N$ and $\Vdc$ are also inverse equivalences of categories between the category of $2$-reduced inner-Kan algebraic bisimplicial sets and strict morphisms and the category of small Verity double categories and strict functors.
\end{theorem}
\begin{proof} We have already shown each part of this theorem except for the fact that $N$ and $\Vdc$ preserve strictness, which follows immediately from the fact that $N$ and $\Bic$ from Chapter~\ref{bicchapter} preserve strictness, which is given as part of Theorem~\ref{bicsummary}.
\end{proof}
\begin{remark}
The functor that forgets algebraic structure from an inner-Kan bisimplicial set is an equivalence of categories, thus $N$ gives an equivalence of categories from the category of small Verity double categories to the category of (non-algebraic) $2$-reduced inner-Kan bisimplicial sets. 
\end{remark}

\chapter{Bisimplicial nerves for fancy bicategories \label{bicatchapter}}

\section{Edge-symmetric VDC's \label{edgesymsec}}
\subsection{The edge-symmetric strict double category}
Ehresmann first observed in \cite{Ehr631} that a strict 2-category gives rise to a strict double category in two different ways. The first way is to consider a 2-category as a strict double category with only identity vertical $1$-morphisms.

On the other hand, we can make a strict $2$-category $\mathcal{B}$ into double category $\ES(\B)$ whose vertical $1$-morphisms are identical to its horizontal $1$-morphisms, with both identical to the $1$-morphisms of $\B$. Then a square:
\begin{center}
\begin{tikzpicture}[scale=1.4,auto]
\begin{scope}
\node (10) at (1,1) {};
\node (00) at (0,1) {};
\node (11) at (1,0) {};
\node (01) at (0,0) {};
\node[rotate=45] at (.5,.5){$\Rightarrow$};
\path[->] (00) edge node[midway]{$f'$}(10);
\path[->] (00) edge node[midway,swap]{$f$}(01);
\path[->] (01) edge node[midway,swap]{$g$}(11);
\path[->] (10) edge node[midway]{$g'$}(11);
\end{scope}
\end{tikzpicture}
\end{center}
of $\ES(\mathcal{B})$ is a $2$-morphism $g \circ f \Rightarrow g'\circ f'$. Square composition is given by whiskering:
\begin{center}
\begin{tikzpicture}[scale=1.8,auto]

\begin{scope}

\node (10) at (1,1) {};
\node (00) at (0,1) {};
\node (11) at (1,0) {};
\node (01) at (0,0) {};
\node[rotate=45] at (.5,.55){$\Rightarrow$};
\node[scale=.8] at (.6,.45){$\theta$};

\node (20) at (2,1) {};
\node (21) at (2,0) {};
\node[rotate=45] at (1.5,.55){$\Rightarrow$};
\node[scale=.8] at (1.6,.45){$\pi$};

\path[->] (00) edge node[midway]{$f'$}(10);
\path[->] (00) edge node[midway,swap]{$f$}(01);
\path[->] (01) edge node[midway,swap]{$g$}(11);
\path[->] (10) edge node[midway]{$g''$}(11);
\path[->] (10) edge node[midway]{$g'$}(20);
\path[->] (11) edge node[midway,swap]{$h$}(21);
\path[->] (20) edge node[midway]{$h'$}(21);
\node at (2.75,.5){$\leadsto$};
\end{scope}
\begin{scope}[shift={(3.5,0)}]
\node (10') at (1.5,1) {};
\node (00') at (0,1) {};
\node (11') at (1.5,0) {};
\node (01') at (0,0) {};
\path[->] (00') edge node[midway]{$g'\circ f'$}(10');
\path[->] (00') edge node[midway,swap]{$f$}(01');
\path[->] (01') edge node[midway,swap]{$h\circ g$}(11');
\path[->] (10') edge node[midway]{$h'$}(11');
\node[rotate=45] at (.7,.43){$\Rightarrow$};
\node[scale=.8] at (.75,.63){$ (\pi \lhd f') \bullet (h \rhd \theta )$};
\end{scope}
\end{tikzpicture}
\end{center}
The other direction of composition is similar. Horizontal and vertical identity squares are given by identity $2$-morphisms in $\B$. 

What structure can be given to a strict double category will allow us to invert this construction? The answer was suggested by Spencer \cite{Spe77}, with details and proofs given by Brown and Mosa. We summarize the theorems of \cite{BM99} in a slightly tweaked form.

 Let $\mathbb{D}=\mathbb{D}_1 \rightrightarrows \mathbb{D}_0$ be a strict double category. The category $\ES(\mathbb{D}_0 )$, with $\mathbb{D}_0$ considered as a $2$-category with only trivial morphisms, has $\ES(\mathbb{D}_0 )_0 = \mathbb{D}_0,$ that is, its objects and vertical morphisms are by definition to those of  $\mathbb{D}_0$. A square of $\ES(\mathbb{D}_0 )$ is the same as a commutative square of morphisms of $\mathbb{D}_0$. 
\begin{definition} \label{strictthinstructure} A \emph{thin structure} on a strict double category $\mathbb{D}$ is a functor of double categories $c:\ES(\mathbb{D}_0 )\ra \mathbb{D}$ which is the identity on the vertical category, and is an isomorphism on the horizontal categories. 
\end{definition}
In particular, a thin structure identifies the vertical category of $\mathbb{D}$ with the horizontal category. If $\B$ is a strict 2-category, then let $\B_0$ be the 2-category of objects and $1$-morphisms of $\B$, with only identity $2$-morphisms. The embedding $\B_0 \ra \B$ as $2$ categories induces a map $\ES(\B_0)\ra \ES(\B)$, which is easily seen to be a thin structure. We call this the \emph{canonical thin structure} on $\ES(\B)$.
\begin{theorem}[Brown-Mosa, Spencer] \label{BMSthm} The map that sends $\B\ra \ES(\B)$ with the canonical connection is an equivalence of categories between the category of strict $2$-categories and the category of double categories with thin structure.
\end{theorem}
\subsection{The edge-symmetric VDC}
\begin{definition} \label{fancydef} A \emph{(small) fancy bicategory} $\B$ consists of an \emph{underlying} (small) bicategory $\widetilde{\B}$ together with a (small) $(2,1)$-category $\overline{\B}$, called the \emph{thin structure $(2,1)$-category}, and strict functor $t_{\B}: \overline{\B} \ra \widetilde{\mathcal{B}}$ called the \emph{thin structure map} such that $t_\B$ is an isomorphism on objects and $1$-morphisms. We treat the objects and $1$-morphisms of $\olb$ and $\wtb$ as being identified by $t_\B$, and refer to them as just the objects and $1$-morphisms of $\B$

A \emph{functor of fancy bicategories} $F:\B \ra \C$ consists of (weak) functors of bicategories $\olf:\olb \ra \olc$ and $\wtf:\wtb \ra \wtc$, such that the following diagram commutes strictly:
\begin{center}
\begin{tikzcd}
\olb \arrow{r}{\olf} \arrow{d}{t_\B} & \olc \arrow{d}{t_\C} \\
\wtb \arrow{r}{\wtf}                  & \wtc
\end{tikzcd}
\end{center} If $\olf$ and $\wtf$ are strict functors, we say $F$ is strict. If they are strictly identity preserving, we say $F$ is strictly identity preserving.\end{definition}
\begin{definition} \label{injectivedefinition} A fancy bicategory $\B$ is called \emph{injective} if $t_\B$ is injective on $2$-morphisms, making $\olb$ into a subcategory of $\wtb.$ In this case, we call a $2$-morphism of $\wtb$ \emph{thin} if it is in the image of $\olb.$ Note that the strictness of $t_\B$ ensures that the components of the associator and unitors of $\wtb$ are thin. A functor between injective fancy bicategories is equivalent a functor $F:\wtb \ra \wtc$ taking thin $2$-morphisms to thin $2$-morphisms, such that the distributor and unitors of $F$ have thin components. 
\end{definition}

\begin{definition} \label{fancificationdef}
Given a bicategory $\B$, we can construct the \emph{complete fancification $\llcorner\B\lrcorner$}  by taking $\widetilde{\llcorner\B\lrcorner}:=\B$ and $\overline{\llcorner\B\lrcorner}$ to be the subcategory of $\B$ consisting of all objects and $1$-morphisms and all invertible $2$-morphisms. This operation is a full and faithful left adjoint to the functor $\B \ra \olb. $  \footnote{Note that this statement can be interpreted in a way that removes the elicit reference to the``category'' of bicategories} A fancy bicategory will be called \emph{complete} if it is isomorphic to a fancy bicategory image of the complete fancification functor.

Given a strict bicategory $\B$, the \emph{sparse fancification $\ulcorner\B\urcorner$ } is formed by taking $\widetilde{\ulcorner\B\urcorner}:=\B$  and $\overline{\ulcorner\B\urcorner}$ to be the objects and $1$-morphisms of $\B$ together with identity $2$-morphisms. This is a full and faithful functor from the ``category'' of strict bicategories and strict functors to the ``category'' of fancy bicategories. A fancy bicategory $\B$  be called \emph{sparse} if it is isomorphic to a fancy bicategory in the image of this functor, or equivalently if $\overline{\B}$ has only identity $2$-morphisms.
\end{definition}

$\ES$ can be generalized to a construction that takes a fancy bicategories to a Verity double category $\ES(\B)$. This construction is very closely related to Verity's $\mathrm{Sq}$ construction, which makes a VDC from a profunctor equipment, and is the chief motivation for the introduction of VDC's in \cite{Ver11}. In the case $\wtb$ is a $(2,1)$-category, and $t_\B$ is the identity, this construction is in fact a special case of Verity's $\mathrm{Sq}$ construction.

\begin{itemize}
\item The horizontal and vertical $(2,1)$-category of $\ES(\B)$ are each identical to $\olb$
\item A square of $\ES(\B)$ is a $2$-morphism $\Theta:g \circ f \Rightarrow g'\circ f'$ in $\B$: 
\begin{figure}
\begin{center}
\begin{tikzpicture}[scale=1.4,auto]
\begin{scope}
\node (10) at (1,1) {};
\node (00) at (0,1) {};
\node (11) at (1,0) {};
\node (01) at (0,0) {};
\node[rotate=45] at (.5,.5){$\Rightarrow$};
\path[->] (00) edge node[midway]{$f'$}(10);
\path[->] (00) edge node[midway,swap]{$f$}(01);
\path[->] (01) edge node[midway,swap]{$g$}(11);
\path[->] (10) edge node[midway]{$g'$}(11);
\end{scope}
\end{tikzpicture}
\end{center}
\caption{A square in $\ES(\B)$. \label{essquare}}
\end{figure}

\item Actions of $2$-morphisms of $\olb$ on a such a square $\Theta$ are given by:

\begin{itemize}
\item  $\Theta \lact \beta := \Theta \bullet (g \rhd t(\beta)) $
\item $\Theta \dact \beta := \Theta \bullet (t(\beta) \lhd f)$
\item $\Theta \uact \beta := (g' \rhd t(\beta^{-1})) \bullet \Theta$
\item $\Theta \ract \beta := (t(\beta^{-1}) \lhd f') \bullet \Theta$
\end{itemize}


\item Horizontal composition uses whiskering:
\begin{figure}[H]
\begin{center}
\begin{tikzpicture}[scale=1.4,auto]

\begin{scope}

\node (10) at (1,1) {};
\node (00) at (0,1) {};
\node (11) at (1,0) {};
\node (01) at (0,0) {};
\node[rotate=45] at (.5,.55){$\Rightarrow$};
\node[scale=.8] at (.6,.42){$\Theta$};

\node (20) at (2,1) {};
\node (21) at (2,0) {};
\node[rotate=45] at (1.5,.55){$\Rightarrow$};
\node[scale=.8] at (1.6,.42){$\Pi$};

\path[->] (00) edge node[midway]{$f'$}(10);
\path[->] (00) edge node[midway,swap]{$f$}(01);
\path[->] (01) edge node[midway,swap]{$g$}(11);
\path[->] (10) edge node[midway]{$g''$}(11);
\path[->] (10) edge node[midway]{$g'$}(20);
\path[->] (11) edge node[midway,swap]{$h$}(21);
\path[->] (20) edge node[midway]{$h'$}(21);
\node at (2.75,.5){$\leadsto$};
\end{scope}
\begin{scope}[shift={(3.5,0)}]
\node (10') at (1,1) {};
\node (00') at (0,1) {};
\node (11') at (1,0) {};
\node (01') at (0,0) {};
\path[->] (00') edge node[midway]{$g'\circ f'$}(10');
\path[->] (00') edge node[midway,swap]{$f$}(01');
\path[->] (01') edge node[midway,swap]{$h\circ g$}(11');
\path[->] (10') edge node[midway]{$h'$}(11');
\node[rotate=45] at (.5,.43){$\Rightarrow$};

\node[scale=.8] at (.50,.63){$\Pi\boxvert\Theta$};
\node[scale=.8] at (4.50,.63){$\Pi\boxvert\Theta:= \alpha^{-1}_{h',g',f'}  \bullet (\Pi \lhd f') \bullet \alpha_{h,g'',f'} \bullet (h \rhd \Theta )\bullet \alpha^{-1}_{h,g,f} $};
\end{scope}
\end{tikzpicture}
\end{center}\caption{Horizontal square composition in $\ES(\B)$ \label{horzcompfig}}
\end{figure}
\item Vertical composition is similar:
\begin{figure}[H]
\begin{center}
\begin{tikzpicture}[scale=1.4,auto]

\begin{scope}

\node (10) at (1,1) {};
\node (00) at (0,1) {};
\node (11) at (1,0) {};
\node (01) at (0,0) {};
\node[rotate=45] at (.5,.55){$\Rightarrow$};
\node[scale=.8] at (.6,.42){$\Theta$};

\node (02) at (0,-1) {};
\node (12) at (1,-1) {};
\node[rotate=45] at (.5,-.5){$\Rightarrow$};
\node[scale=.8] at (.6,-.63){$\Pi$};

\node at (1.75,.0){$\leadsto$};
\path[->] (00) edge node[midway]{$f'$}(10);
\path[->] (00) edge node[midway,swap]{$f$}(01);
\path[->] (01) edge node[midway,swap]{$g''$}(11);
\path[->] (10) edge node[midway]{$g'$}(11);
\path[->] (01) edge node[midway,left]{$g$}(02);
\path[->] (11) edge node[midway,right]{$h'$}(12);
\path[->] (02) edge node[midway,below]{$h$}(12);
\end{scope}
\begin{scope}[shift={(2.8,-.5)}]
\node (10') at (1,1) {};
\node (00') at (0,1) {};
\node (11') at (1,0) {};
\node (01') at (0,0) {};
\path[->] (00') edge node[midway]{$f'$}(10');
\path[->] (00') edge node[midway,swap]{$g\circ f$}(01');
\path[->] (01') edge node[midway,swap]{$h$}(11');
\path[->] (10') edge node[midway]{$h'\circ g'$}(11');
\node[rotate=45] at (.5,.43){$\Rightarrow$};

\node[scale=.8] at (.50,.63){$\Pi\boxminus\Theta$};
\node[scale=.8] at (4.75,.63){$\Pi\boxminus\Theta:= \alpha_{f',g',h'}  \bullet (h' \rhd \Theta) \bullet \alpha^{-1}_{f,g'',h'} \bullet (\Pi \lhd f )\bullet \alpha_{f,g,h} $};
\end{scope}
\end{tikzpicture}
\end{center} \caption{Vertical square composition in $\ES(\B)$ \label{vertcompfig}}
\end{figure}
\item Pseudo-identity squares:

\begin{center}
\begin{tikzpicture}[scale=1.4,auto]

\begin{scope}[shift={(2.8,-.5)}]
\node (10') at (1,1) {};
\node (00') at (0,1) {};
\node (11') at (1,0) {};
\node (01') at (0,0) {};
\path[->] (00') edge node[midway]{$f$}(10');
\path[-] (00') edge[thick,double] (01');
\path[->] (01') edge node[midway,swap]{$f$}(11');
\path[-] (10') edge[thick,double] (11');
\node[rotate=45] at (.5,.39){$\Rightarrow$};

\node[scale=.8] at (.50,.65){$\ID^h_f$};
\node at (2,.5){$:=  \lambda_f \bullet \rho^{-1}_f$};
\end{scope}
\begin{scope}[shift={(6.8,-.5)}]
\node (10') at (1,1) {};
\node (00') at (0,1) {};
\node (11') at (1,0) {};
\node (01') at (0,0) {};
\path[-] (00') edge[thick,double] (10');
\path[->] (00') edge node[midway,left]{$f$} (01');
\path[-] (01') edge[thick,double] (11');
\path[->] (10') edgenode[midway,right]{$f$} (11');
\node[rotate=45] at (.5,.39){$\Rightarrow$};

\node[scale=.8] at (.50,.65){$\ID^v_f$};
\node at (2,.5){$:=  \rho_f \bullet \lambda^{-1}_f$};
\end{scope}

\end{tikzpicture}
\end{center}
\end{itemize}

Note that if $\wtb$ (and thereby $\olb$) have only the identity $2$-morphisms, this reduces to the definition given above for $\ES(\B)$. We now check that $\ES(\B)$ satisfies the axioms for a Verity double category. Consider, for instance, the square interchange axiom, \textbf{VDC11}.
Suppose we have squares $\Theta,\Pi,\Pi'$, and $\Sigma$ in $\ES(\B)$ as shown in Figure~\ref{interchangefigure}. 
\begin{figure}[h]
\begin{center}
\begin{tikzpicture}[scale=1.5,auto]

\begin{scope}

\node (10) at (1,1) {};
\node (00) at (0,1) {};
\node (11) at (1,0) {};
\node (01) at (0,0) {};
\node[rotate=45] at (.5,.6){$\Rightarrow$};
\node[scale=.8] at (.6,.45){$\Theta$};

\node (20) at (2,1) {};
\node (21) at (2,0) {};
\node(22) at (2,-1){};
\node[rotate=45] at (1.5,.6){$\Rightarrow$};
\node[scale=.8] at (1.6,.45){$\Pi'$};

\node (12) at (1,-1) {};
\path[->] (00) edge node[midway]{\scriptsize $f'$}(10);
\path[->] (00) edge node[midway,swap]{\scriptsize $ f$}(01);
\path[->] (01) edge node[midway,above]{\scriptsize $ g''$}(11);
\path[->] (10) edge node[midway,right]{\scriptsize $g'''$}(11);
\path[->] (10) edge node[midway]{\scriptsize $g'$}(20);
\path[->] (11) edge node[midway,swap]{\scriptsize $h'''$}(21);
\path[->] (20) edge node[midway]{\scriptsize $h'$}(21);
\path[->](21) edge node[midway]{\scriptsize $i'$}(22);
\path[->](12) edge node[midway,below]{$i$}(22);

\node (01') at (0,-1) {};
\node[rotate=45] at (0.5,-.43){$\Rightarrow$};
\node[scale=.8] at (0.6,-.58){$\Pi$};

\path[->] (11) edge node[midway]{\scriptsize $h''$}(12);
\path[->] (01) edge node[midway,swap]{\scriptsize $g$ }(01');
\path[->] (01') edge node[midway,swap]{\scriptsize $h$}(12);
\node[rotate=45] at (1.5,-.43){$\Rightarrow$};
\node[scale=.8] at (1.6,-.58){$\Sigma$};

\end{scope}
\end{tikzpicture}
\end{center} \caption{\label{interchangefigure} }\end{figure}
The square interchange axiom asserts that $$(\Sigma \boxvert \Pi) \boxminus  (\Pi' \boxvert \Theta)= (\Sigma \boxminus \Pi')\boxvert (\Pi \boxminus \Theta). $$ Applying our definitions, this works out to the following monstrous equality of $2$-morphisms in $\B$:

\begin{equation}
 \begin{split}
&\alpha_{i', h',g'\circ f'} \bullet (i'\rhd [  \alpha^{-1}_{h',g',f'}  \bullet (\Pi' \lhd f') \bullet \alpha_{h,g''',f'} \bullet (h''' \rhd \Theta )\bullet \alpha^{-1}_{h''',g'',f}  ])  \bullet \alpha^{-1}_{i',h'''\circ g'',f } \\
 &\quad \bullet  ([   \alpha^{-1}_{i',h''',g''}  \bullet (\sigma \lhd g'') \bullet \alpha_{i,h'',g''} \bullet (i \rhd \Pi )\bullet \alpha^{-1}_{i,h,g}   ] \lhd f)\bullet  \alpha_{i\circ h,g,f}  \\
&= \alpha^{-1}_{i'\circ h',g',f'}   \bullet ([\alpha_{i',h',g'}  \bullet  (i'\rhd  \Pi') \bullet    \alpha^{-1}_{i',h''',g'''} \bullet    (\sigma \lhd g''')      \bullet   \alpha_{i,h'',g'''} ] \lhd f')  \bullet \alpha_{i,h'' \circ g''',f'}  \\
 &\quad\bullet (i \rhd [\alpha_{h'',g''',f'}    \bullet    (h'' \rhd \Theta )  \bullet \alpha^{-1}_{h'',g'',f} \bullet (\Pi  \lhd f) \bullet\alpha_{h,g,f}])\bullet \alpha^{-1}_{i,h,g\circ f}. 
 \end{split} \label{biginterchangeequation}
\end{equation}

We could of course proceed to prove this directly using the bicategory axioms, repeatedly applying the naturality of the associator, the pentagon identity, and interchange. However, if we view Figure~\ref{interchangefigure} not as a diagram of squares in $\ES(\B)$ but as a ``pasting diagram'' of $2$-morphisms in $\wtb$, then Equation~\ref{biginterchangeequation} asserts the equality of two particular ways of interpreting the diagram as a $2$-morphism from (a particular bracketing of) the bottom ``edge'' of the diagram to (a particular bracketing of) the top edge.

In fact, it is possible to state a coherence result that covers this situation. 
\begin{lemma} \label{coherencelemma} In general, given a pasting diagram and a bracketing of the $1$-morphisms making up the top and bottom edges of the diagram, every way of sprinkling in associators and unitors to compose the diagram into a $2$-morphism between the composed bottom \mbox{$1$-morphism} and the composed top \mbox{$1$-morphism} gives the same result.
\end{lemma}
This sort of statement will be familiar folklore to some readers, and of course it has something to do with MacLane's coherence theorem. However, it takes some effort to state it precisely. A rigorous version of Lemma~\ref{coherencelemma} is developed in appendix A of \cite{Ver11} (see in particular the discussion on pages 261-262), for exactly the same purposes as we now apply it to.  We will use the coherence of pasting without giving a precise statement and refer interested readers to \cite{Ver11}.

It is important to note that this coherence lemma also allows for the cases involving identities under the rubric of ``sprinkling in unitors''. For instance, in \textbf{VDC10} we wish to show
 \begin{equation*}\label{unitproof1}\begin{split}
&(\Theta\boxminus \ID_{f'})\lact \rho_f \ract \rho_{g'} \\
&= (\rho^{-1}_{g'} \lhd f') \bullet [\alpha_{g',\id,f'}  \bullet (g' \rhd [\lambda_{f'} \bullet \rho^{-1}_{f'}]) \bullet \alpha^{-1}_{\id,f',g'} \bullet (\Theta \lhd \id )\bullet \alpha_{\id,f,g}]\bullet (g\rhd \rho_f).  \\
&= \Theta.
\end{split}
\end{equation*}
This again follows from Lemma~\ref{coherencelemma}, since both sides are possible ways of composing the diagram below.
\begin{center}
\begin{tikzpicture}[scale=1.4,auto]
\begin{scope}
\node (10) at (1,1) {};
\node (00) at (0,1) {};
\node (11) at (1,0) {};
\node (01) at (0,0) {};
\node[rotate=45] at (.5,.5){$\Rightarrow$};
\path[->] (00) edge node[midway]{$f'$}(10);
\path[->] (00) edge node[midway,swap]{$f$}(01);
\path[->] (01) edge node[midway,swap]{$g$}(11);
\path[->] (10) edge node[midway]{$g'$}(11);
\end{scope}
\end{tikzpicture}
\end{center}

The other axioms similarly follow from coherence, so $\ES(\B)$ is a VDC. Now we aim to promote $\ES$ to a functor from the ``category'' of fancy bicategories to the ``category'' of VDC's. 

\begin{definition} Let $F:\B\ra\C$ be a functor of fancy bicategories. Then $\ES(F):\ES(\B)\ra \ES(\C)$ is defined as follows:
\begin{itemize} \item $\ES(F)$ is defined on objects, horizontal and vertical $1$-morphisms, and horizontal and vertical $2$-morphisms by $\overline{F}$. 
\item The vertical and horizontal distributors and unitors of $\ES(F)$ are given by the distributors and unitors $\phi$ and $\upsilon$ of $\olf$.
\item For a square $\Theta:g\circ f \ra g'\circ f'$ in $\ES(B)$ define $$\ES(F)(\Theta):= t_\C(\phi_{g',f'}) \bullet F(\Theta)\bullet t_\C(\phi^{-1}_{g,f}).$$  
\end{itemize}
\end{definition}
We know by definition that $\ES(F)$ gives a functor on the vertical and horizontal $(2,1)$-categories of $\ES(\B)$ and $\ES(\B)$, both given by $\olf$. We need only check the axioms, \textbf{VFun1}, \textbf{VFun2}, and \textbf{VFun3}. These can be seen to follow from Corollary~A.0.12 in \cite{Ver11}, which states the sense in which pasting is preserved by a bicategory functor. Alternatively, a direct verification is provided below. We denote the bicategory axioms 1-17 from Chapter~\ref{bicchapter} by \textbf{B1}-\textbf{B17}.

\begin{itemize}
\item For \textbf{VFun1}, let $\Theta,\Pi$ be squares of $\ES(\B)$ as shown below: 

\begin{center}
\begin{tikzpicture}[scale=1.4,auto]

\begin{scope}

\node (10) at (1,1) {};
\node (00) at (0,1) {};
\node (11) at (1,0) {};
\node (01) at (0,0) {};
\node[rotate=45] at (.5,.55){$\Rightarrow$};
\node[scale=.8] at (.6,.41){$\Theta$};

\node (20) at (2,1) {};
\node (21) at (2,0) {};
\node[rotate=45] at (1.5,.55){$\Rightarrow$};
\node[scale=.8] at (1.6,.41){$\Pi$};

\path[->] (00) edge node[midway]{$f'$}(10);
\path[->] (00) edge node[midway,swap]{$f$}(01);
\path[->] (01) edge node[midway,swap]{$g$}(11);
\path[->] (10) edge node[midway]{$g''$}(11);
\path[->] (10) edge node[midway]{$g'$}(20);
\path[->] (11) edge node[midway,swap]{$h$}(21);
\path[->] (20) edge node[midway]{$h'$}(21);
\end{scope}\end{tikzpicture}
\end{center}

Writing down the definitions, we have:
 \begin{equation}\label{VFun1proof1}
 \begin{split}
&[\ES(F)(\Pi)\boxvert \ES(F)(\Theta)]\dact \phi_{h,g}\uact \phi_{g',f'}\\
&\quad = (F(h')\rhd \phi^{-1}_{g',f'} ) \bullet [{\alpha'}^{-1}_{F(h'),F(g'),F(f')} \bullet ( [ \phi_{h',g'} \bullet F(\Pi) \bullet \phi^{-1}_{h,g''}]\lhd  F(f')   ) 
\\ &\quad \bullet {\alpha'}_{F(h),F(g''),F(f')} \bullet (F(h)\rhd [ \phi_{g'',f'} \bullet F(\Theta) \bullet \phi^{-1}_{g,f}]  ) \bullet {\alpha'}^{-1}_{F(h),F(g),F(f)}] 
\\ &\quad \bullet (\phi_{h,g} \lhd F(f))
 \end{split}
\end{equation}
Distribute $\lhd F(f')$ and $F(h)\rhd$ in Equation~\ref{VFun1proof1} using interchange (\textbf{B5} and \textbf{B6}) to get: 
 \begin{equation} \label{VFun1proof2}
 \begin{split}
&[\ES(F)(\Pi)\boxvert \ES(F)(\Theta)]\dact \phi_{h,g}\uact \phi_{g',f'}\\
&=\! (F(h')\rhd \phi^{-1}_{g',f'} )\! \bullet \! {\alpha'}^{-1}_{F(h'),F(g'),F(f')} \bullet ( \phi_{h',g'}\lhd  F(f') ) \bullet (F(\Pi)\lhd  F(f'))  \bullet (\phi^{-1}_{h,g''}\lhd  F(f') ) 
\\ &\quad \bullet {\alpha'}_{F(h),F(g''),F(f')} \bullet (F(h)\rhd \phi_{g'',f'}) \bullet (F(h)\rhd F(\Theta)) \bullet (F(h)\rhd\phi^{-1}_{g,f})   \bullet {\alpha'}^{-1}_{F(h),F(g),F(f)}
\\ &\quad \bullet (\phi_{h,g} \lhd F(f))
 \end{split}
 \end{equation}
Next we give identities using the functor axiom \textbf{BFun5}:
\begin{align*}
(F(h')\rhd \phi^{-1}_{g',f'} ) \bullet {\alpha'}^{-1}_{F(h'),F(g'),F(f')} \bullet ( \phi_{h',g'}\lhd  F(f') )&= \phi_{h',g'\circ f'}  \bullet F(\alpha_{h',g',f'})^{-1}\bullet \phi^{-1}_{h'\circ g',f'}\\
(\phi^{-1}_{h,g''}\lhd  F(f') ) \bullet {\alpha'}_{F(h),F(g''),F(f')} \bullet (F(h)\rhd \phi_{g'',f'})&=\phi_{ h\circ g'',f' } \bullet F(\alpha_{h,g'',f'})\bullet \phi^{-1}_{h,g''\circ f'} \\
 (F(h)\rhd\phi^{-1}_{g,f})   \bullet {\alpha'}^{-1}_{F(h),F(g),F(f)} \bullet (\phi_{h,g} \lhd F(f)) &= \phi_{ h,g\circ f} \bullet  F(\alpha_{h,g,f})^{-1} \bullet \phi^{-1}_{h\circ g, f}  \\
\end{align*}
Substituting these into Equation~\ref{VFun1proof2}
 \begin{equation} \label{VFun1proof3}
 \begin{split}
&[\ES(F)(\Pi)\boxvert \ES(F)(\Theta)]\dact \phi_{h,g}\uact \phi_{g',f'}\\
&\quad = \phi_{h',g'\circ f'}  \bullet F(\alpha_{h',g',f'})^{-1}\bullet \phi^{-1}_{h'\circ g',f'} \bullet (F(\Pi)\lhd  F(f'))  \bullet \phi_{ h\circ g'' ,f'}
\\ &\quad \bullet F(\alpha_{h,g'',f'})\bullet \phi^{-1}_{h,g''\circ f'} \bullet (F(h)\rhd F(\Theta)) \bullet \phi_{ h,g\circ f} \bullet  F(\alpha_{h,g,f})^{-1}\bullet \phi^{-1}_{h\circ g, f}
 \end{split}
 \end{equation}
Next note the following identities from functor axioms \textbf{BFun4} and \textbf{BFun3}:
\begin{align*}
 \phi^{-1}_{h'\circ g',f'} \bullet (F(\Pi)\lhd  F(f'))  \bullet \phi_{h\circ g'', f' }&=  F(\Pi \lhd f')           \\
\phi^{-1}_{h,g''\circ f'} \bullet (F(h)\rhd F(\Theta)) \bullet \phi_{h, g\circ f}     &=  F(h \rhd \Theta)            \\
\end{align*}
Substitute these into Equation~\ref{VFun1proof3}:
 \begin{equation} \label{VFun1proof4}
 \begin{split}
&[\ES(F)(\Pi)\boxvert \ES(F)(\Theta)]\dact \phi_{h,g}\uact \phi_{g',f'}\\
&\quad = \phi_{h',g'\circ f'}  \bullet F(\alpha_{h',g',f'})^{-1}\bullet  F(\Pi \lhd f')    
\bullet F(\alpha_{h,g'',f'})\bullet F(h \rhd \Theta)  \bullet  F(\alpha_{h,g,f})^{-1}\bullet \phi^{-1}_{h\circ g, f}
 \end{split}
 \end{equation}
 Finally \textbf{BFun1} and \textbf{BFun2} state that $F$ is functorial with respect to $\bullet$. Apply this to Equation~\ref{VFun1proof4} to get:
 \begin{equation*}
 \begin{split}
&[\ES(F)(\Pi)\boxvert \ES(F)(\Theta)]\dact \phi_{h,g}\uact \phi_{g',f'}\\
&\quad = \phi_{h',g'\circ f'}  \bullet F(\alpha_{h',g',f'}^{-1}\bullet  (\Pi \lhd f')    
\bullet \alpha_{h,g'',f'}\bullet (h \rhd \Theta)  \bullet  \alpha_{h,g,f}^{-1})\bullet \phi^{-1}_{h\circ g, f} \\
& \quad = \ES(F)(\Pi \boxvert \Theta).
 \end{split}
 \end{equation*}
The $\boxminus$ case is similar.
\item For \textbf{VFun2}, let $f:a \ra b$ be a $1$-morphism in $\B$. Then considering $f$ as a horizontal $1$-morphism in $\ES(F)$, we have 
 \begin{align*}
\ID^h_{\ES(F)(f)}  \lact \upsilon^v_a \ract \upsilon^v_b &=  (\upsilon^{-1}_b\lhd F(f))  \bullet  {\lambda'}_{F(f)} \bullet {\rho'}^{-1}_{F(f)}  \bullet (F(f)\rhd \upsilon_a) \\
&=  \phi_{\id_b,f}\bullet F(\lambda_f)\bullet F(\rho^{-1}_f) \bullet \phi^{-1}_{ f,\id_a} \\ 
&= \phi_{\id_b,f}\bullet F(\lambda_f \bullet \rho^{-1}_f) \bullet \phi^{-1}_{f, \id_a} \\
&= \ES(F)(\ID_f)
 \end{align*}
 where we use \textbf{BFun6} and \textbf{BFun7} for the first simplification. The  $\ID^v$ case is similar.
 \item For \textbf{VFun3}, let $\Theta: g\circ f \Rightarrow g'\circ f'$ be a square in $\ES(\B)$ and let $\eta:i\Rightarrow f$ be a $2$-morphism in $\overline{\B}$, viewed as a vertical $2$-morphism of $\ES(\B)$. We have:
 \begin{align*} \ES(F)(\Theta) \lact \ES(F)(\eta)& =  \phi_{g',f'} \bullet F(\Theta) \bullet \phi^{-1}_{g,f} \bullet ( F(g) \rhd F(\eta) )\\
 &=  \phi_{g',f'} \bullet F(\Theta) \bullet F(g\rhd \eta) \bullet \phi^{-1}_{i,g} \\
 &=  \phi_{g',f'} \bullet F(\Theta \bullet (g\rhd \eta)) \bullet \phi^{-1}_{i,g} \\ 
 &= \ES(F)(\Theta \lact \eta).
 \end{align*} 
where we use \textbf{BFun3} for the first simplification. The other cases follow similarly.
\end{itemize}

\subsection{Thin structures for VDC's}

\begin{definition}
Let $\D$ be a Verity double category.  A \emph{thin structure} on $\D$ consists of a $(2,1)$-category $\T$ together with a functor $\frakt:\ES(\llcorner\T \lrcorner )\ra \D$ inducing isomorphisms on horizontal and vertical $(2,1)$-categories, $\T \cong \mathcal{H} \cong \mathcal{V}$. 
\end{definition}

If $\B$ is a fancy bicategory, we have an embedding $\llcorner\olb \lrcorner\ra \B.$ Applying $\ES$, we get a functor $\ES(\llcorner\olb \lrcorner)\ra \ES(\B)$, giving a canonical thin structure on $\ES(\B).$ We now give an inverse construction, taking a Verity double category to a bicategory:

\begin{definition}  Let $(\D,\T, \frakt)$ be a Verity double category with thin structure. We construct a fancy bicategory structure $\fold(\D, \T,\frakt)$ as follows: 
\begin{itemize}
\item The objects and morphisms of  $\fold(\D, \T,\frakt)$ are the objects and horizontal $1$-morphisms of $\D.$ $1$-morphism composition and pseudo-identities for objects are defined as in $\D$.
\item A $2$-morphism $f\Rightarrow g$ of $\widetilde{ \fold(\D, \T,\frakt)   }$ is a horizontal-globular square of $\D$ from $f$ to $g$.
\item If $f \stackrel{\theta}{\Rightarrow} g \stackrel{\eta}{\Rightarrow} h: a \ra b$ are $2$-morphisms given by horizontal-globular squares, we define $$\theta \bullet \eta = (\theta \boxminus \eta) \lact \lambda^v_{\id_a} \ract\lambda^v_{\id_b}.$$ Recall that $\lambda^v_{\id_a}=\rho^v_{\id_b}$, resolving the apparent asymmetry of this definition.
\item The identity for a $1$-morphism $f$ is given by the pseudo-identity square $\ID_f$.
\item If $f:a\ra b$ is a $1$-morphism of $  \fold(\D, \T,\frakt),$ given by a horizontal $1$-morphism of $\D$, and $g\stackrel{\eta}{\Rightarrow} h$ is a $2$-morphism given by a horizontal-globular square, $$\eta \lhd f := \eta \boxvert \ID^h_f.$$ Similarly, if $f \stackrel{\eta}{\Rightarrow} g:a \ra b$ and $h:b\ra c$ then  $$h\rhd \eta := \ID^h_h \boxvert \eta.$$
\item We take $\overline{ \fold(\D, \T,\frakt)   }:=\mathcal{H}$ and the functor $t:\overline{ \fold(\D, \T,\frakt)   } \ra  \widetilde{\fold(\D, \T,\frakt)}$ to be given by the identity on objects and $1$-morphisms and for $\theta:f \ra g$ by  $t(\theta):= \Id_g \uact \theta$

\item The associator and unitors of $\widetilde{\fold(\D, \T,\frakt)}$  are defined from those of $\overline{\fold(\D, \T,\frakt)},$ for instance  $t(\alpha_{h,g,f})$ is the associator of $\widetilde{\fold(\D, \T,\frakt)}.$
\end{itemize}
\end{definition}

\begin{proposition} $\fold( \D, \T, \frakt )$ satisfies the axioms for a fancy bicategory.
\end{proposition} 
\begin{proof}
If we had an appropriate generalization of Lemma~\ref{coherencelemma} that asserted coherence for pasting diagrams in VDC's, then we could use that to coherence to easily check the axioms. However, we will not develop this idea, so instead it is necessary to check the axioms manually. We check that  $t:\overline{ \fold(\D, \T,\frakt)   } \ra  \widetilde{\fold(\D, \T,\frakt)}$ is a strict functor.  It's immediate from definitions that $t$ preserves identity $2$-morphisms, the associator, and the unitors. Let us check that $t$ preserves $\bullet,$ $\rhd,$ and $\lhd$:
\begin{itemize}
\item For $\bullet$, take $f \stackrel{\theta}{\Rightarrow} g \stackrel{\pi}{\Rightarrow} h:a\ra b$ in $\mathcal{H}$  
\begin{align*} 
T(\pi)\bullet T(\theta)&= [(\ID_h \uact \pi) \boxminus (\ID_g\uact \theta )]\ract \lambda_{\id_b}  \lact  \lambda_{\id_a}  \\
                       &= [\ID_h \boxminus (\ID_g\uact \theta\dact \pi^{-1})]\ract \lambda_{\id_b}  \lact  \lambda_{\id_a}  \\
\end{align*}
by mid-associativity (Lemma~\ref{midass}). Next from \textbf{VDC7} we have $$\ID_g=  \ID_h\uact \pi \dact \pi.$$ Making this substitution, we have:

\begin{tabular}{p{1.8cm} p{8.5cm} p{3.8cm}}
$t(\pi)\bullet t(\theta)$ &  $ =[\ID_h \boxminus (\ID_h\uact \pi \dact \pi \uact \theta\dact \pi^{-1})]\ract \lambda_{\id_b} \lact  \lambda_{\id_a}$ &  \\
  &$=  [\ID_h \boxminus (\ID_h\uact \pi \uact \theta)]\ract \lambda_{\id_b}  \lact  \lambda_{\id_a}$ & (\textbf{VDC2} and \textbf{VDC3})  \\
  &$=  [\ID_h \boxminus (\ID_h\uact (\pi \bullet\theta) )] \ract \lambda_{\id_b} \lact  \lambda_{\id_a}$ & (\textbf{VDC2})  \\
  &$=\ID_h\uact (\pi \bullet\theta) $& (\textbf{VDC10})\\
  &$= t(\pi\bullet \theta)$&
\end{tabular}

\item For $\rhd,$ take

 $f \stackrel{\theta}{\Rightarrow}  g: a \ra b$ and $h: b\ra c$ in $\mathcal{H}$
 
\begin{tabular}{p{2cm} p{8.5cm} p{2cm}}

$t(h) \rhd t(\theta)$& $=  [\ID_h \boxvert (\ID_g\uact \theta)]$ &  \\
& $=  [\ID_h\boxvert \ID_g]\uact (h \rhd \theta)$ & (\textbf{VDC8})  \\
& $=  \ID_{h\circ g}\uact (h \rhd \theta)$ & (\textbf{VDC6})  \\
& $ =t(h \rhd \theta) $& \end{tabular}
\item $t$ preserves $\lhd$ by a similar argument.

\end{itemize}
We can use this functor to immediately deduce the pentagon identity (\textbf{B17}), the compatibility of the unitors and the associator (\textbf{B13}, \textbf{B14}, and \textbf{B15}), the compatibility of the unitors and the pseudo-identity (\textbf{B12}), and the interchange of whiskering and the identity (\textbf{B4}) for  $ \widetilde{\fold(\D, \T,\frakt)}$, which all follow from the corresponding axioms for $\mathcal{H}=  \overline{ \fold(\D, \T,\frakt)   } $. We now check the other axioms:

\begin{itemize}
\item \textbf{B1} follows directly from the compatibility of $\lambda$ and $\boxminus$.
\item  For \textbf{B2} we have:

\begin{tabular}{p{1.7cm} p{10.6cm} p{1.8cm}}
$\Sigma \bullet (\Pi \bullet \Theta)$  &   $=(\Sigma \boxminus [ (\Pi\boxminus \Theta) \lact \lambda^v_{\id} \ract \lambda^v_{\id}])\lact \lambda^v_{\id} \ract \lambda^v_{\id}$ &  \\
  &$=[(\Sigma \boxminus  (\Pi\boxminus \Theta)) \lact (\id \rhd \lambda^v_{\id}) \ract (\id \rhd \lambda^v_{\id})]\lact \lambda^v_{\id} \ract \lambda^v_{\id}$  & (\textbf{VDC8}) \\
  & $=[[((\Sigma \boxminus  \Pi)\boxminus \Theta) \lact \alpha_{\id,\id,\id} \ract \alpha_{\id,\id,\id} ]$& (\textbf{VDC9}) \\& \quad $\lact (\id \rhd \lambda^v_{\id}) \ract (\id \rhd \lambda^v_{\id})]\lact \lambda^v_{\id} \ract \lambda^v_{\id}$             & \\
  & $=[((\Sigma \boxminus  \Pi)\boxminus \Theta) \lact (\alpha_{\id,\id,\id} \bullet (\id \rhd \lambda^v_{\id}))\ract (\alpha_{\id,\id,\id} \bullet (\id \rhd \lambda^v_{\id}))]$& (\textbf{VDC2})\\& \quad $\lact \lambda^v_{\id} \ract \lambda^v_{\id}$   &   \\
 & $=[((\Sigma \boxminus  \Pi)\boxminus \Theta) \lact (\rho^v_{\id} \lhd \id )\ract (\rho^v_{\id} \lhd \id )]\lact \lambda^v_{\id} \ract \lambda^v_{\id}$   & (\textbf{B14})   \\
  & $=[((\Sigma \boxminus  \Pi)\boxminus \Theta) \lact (\lambda^v_{\id} \lhd \id )\ract (\lambda^v_{\id} \lhd \id )]\lact \lambda^v_{\id} \ract \lambda^v_{\id}$   & (\textbf{B12})  \\
 &   $=[((\Sigma \boxminus \Pi)\lact \lambda^v_{\id} \ract \lambda^v_{\id})\boxminus  \Theta] \lact \lambda^v_{\id} \ract \lambda^v_{\id} $ &  \\ 
    & $=(\Sigma \bullet \Pi) \bullet \Theta$   &   
\end{tabular}

\item For \textbf{B5}:

\begin{tabular}{p{1.8cm} p{10cm} p{2.1cm}}
$i \rhd (\Pi \bullet \Theta)$  &    $=\Id^h_i \boxvert [(\Pi \boxminus \Theta)\lact \lambda^v_{\id} \ract \lambda^v_{\id} ]$  &  \\
  &  $=[(\Id^h_i\boxminus \Id^h_i)\lact \lambda^v_{\id} \ract \lambda^v_{\id} ] \boxvert [(\Pi \boxminus \Theta)\lact \lambda^v_{\id} \ract \lambda^v_{\id}] $  & (\textbf{VDC10}) \\
   &  $=[(\Id^h_i\boxminus \Id^h_i)\ract \lambda^v_{\id} ] \boxvert [(\Pi \boxminus \Theta)\lact \lambda^v_{\id} ] $          &  (Cor.~\ref{midcancel})        \\
   &  $=[(\Id^h_i\boxminus \Id^h_i) \boxvert (\Pi \boxminus \Theta)]\lact \lambda^v_{\id}  \ract \lambda^v_{\id} $          &      (\textbf{VDC4})      \\
      &  $=[(\Id^h_i\boxvert \Pi) \boxminus ( \Id^h_i \boxvert \Theta)]\lact \lambda^v_{\id}  \ract \lambda^v_{\id} $          &   (\textbf{VDC11})        \\
            &  $=(i \rhd \Pi) \bullet ( i \rhd \Theta) $          & 
\end{tabular}

\textbf{B6} is similar.
\item For \textbf{B7}, let $\Theta:f\Rightarrow g$ in  $\widetilde{\fold(\D, \T,\frakt)}$.

\begin{tabular}{p{2.1cm} p{9.7cm} p{2.1cm}}
$(\Theta \lhd \id )\bullet \rho_f $  &    $=((\Theta \boxvert \ID_\id) \boxminus (\ID_{f \circ \id} \uact \rho_f))\lact \lambda^v_{\id} \ract \lambda^v_{\id}$  &  \\
  &  $=([(\Theta \boxvert \ID_\id)\uact \rho_f ] \boxminus (\ID_{f \circ \id}  \uact \rho_f\dact \rho_f ))\lact \lambda^v_{\id} \ract \lambda^v_{\id}$  & (Cor.~\ref{midcancel})  \\
  &  $=([(\Theta \boxvert \ID_\id)\uact \rho_f ] \boxminus \ID_{f})\lact \lambda^v_{\id} \ract \lambda^v_{\id}$  & (\textbf{VDC7})  \\
  &  $=([(\Theta \boxvert \ID_\id)\uact \rho_f ] \boxminus \ID_{f})\lact \rho^v_{\id} \ract \rho^v_{\id}$  &  (\textbf{B12})  \\
  &  $=(\ID_g \boxminus [(\Theta \boxvert \ID_\id)\uact \rho_f ])\lact \lambda^v_{\id} \ract \lambda^v_{\id}$  &  (\textbf{VDC10})    \\
  &  $=((\ID_g \uact \rho_g) \boxminus [(\Theta \boxvert \ID_\id)\uact \rho_f\dact \rho_g ])\lact \lambda^v_{\id} \ract \lambda^v_{\id}$  &   (Cor.~\ref{midcancel})      \\
  &  $=((\ID_g \uact \rho_g) \boxminus \Theta )\lact \lambda^v_{\id} \ract \lambda^v_{\id}$  &   (\textbf{VDC10})      \\  
  &  $=\rho_g \bullet \eta$  &           
\end{tabular}

\textbf{B8} is similar.

\item The verifications of \textbf{B9}, \textbf{B10}, and \textbf{B11} are similar to \textbf{B7}. We leave the details to the reader.

\item For \textbf{B16}, let $\Theta: f\Rightarrow g: a \ra b $ and $\Pi: h \Rightarrow i$ in $\widetilde{\fold(\D, \T,\frakt)}$

\begin{tabular}{p{2.9cm} p{9.2cm} p{2.1cm}}
$(i\rhd \Theta) \bullet (\Pi \lhd f) $  &    $=[(\ID_i \boxvert \Theta)\boxminus (\Pi \boxvert \ID_f)]\lact \lambda_{\id}\ract \lambda_{\id} $  &  \\
     &    $=[(\ID_i \boxminus \Pi )\boxvert (\Theta \boxminus \ID_f)]\lact \lambda_{\id}\ract \lambda_{\id} $  & (\textbf{VDC11})  \\
     &    $=((\ID_i \boxminus \Pi )\ract \lambda_{\id} \lact \lambda_{\id}) \boxvert ((\Theta \boxminus \ID_f)\lact \lambda_{\id}\ract \lambda_{\id}) $ & (Cor.~\ref{midcancel}) \\
     &    $=((\ID_i \boxminus \Pi )\ract \lambda_{\id} \lact \lambda_{\id}) \boxvert ((\Theta \boxminus \ID_f)\lact \rho_{\id}\ract \rho_{\id}) $ &  (\textbf{B12}) \\   
     &    $=(( \Pi \boxminus \ID_h)\ract \rho_{\id} \lact \rho_{\id}) \boxvert ((\ID_g \boxminus \Theta )\lact \lambda_{\id}\ract \lambda_{\id}) $ & (\textbf{VDC10}) \\     
     &    $=(( \Pi \boxminus \ID_h)\ract \lambda_{\id} \lact \lambda_{\id}) \boxvert ((\ID_g \boxminus \Theta )\lact \lambda_{\id}\ract \lambda_{\id}) $ &  (\textbf{B12}) \\    
     &    $=(( \Pi \boxminus \ID_h)\ract \lambda_{\id} ) \boxvert ((\ID_g \boxminus \Theta )\lact \lambda_{\id}) $ & (Cor.~\ref{midcancel})  \\ 
     &    $=[( \Pi \boxminus \ID_h) \boxvert (\ID_g \boxminus \Theta )]\lact \lambda_{\id}\ract \lambda_{\id}  $ &  (\textbf{VDC4}) \\ 
     &    $=[( \Pi \boxvert \ID_g) \boxminus (\ID_h \boxvert \Theta )]\lact \lambda_{\id}\ract \lambda_{\id}  $ & (\textbf{VDC11})  \\   
     &    $=( \Pi \lhd g) \bullet (h \rhd \Theta ) $ &  \      
           \end{tabular}

\end{itemize}
\end{proof}
\begin{conjecture}\label{foldesconj}
$\mbox{Fold}$ and $\ES$ are inverse equivalences of categories between the category of small fancy bicategories and strictly identity-preserving functors, and the category of small Verity double categories with thin structure and strictly identity-preserving functors. 
\end{conjecture}
The above conjecture can be proved in a straightforward but very tedious manner by dressing up the proof given in \cite{BM99} of Theorem~\ref{BMSthm} with distributors and associators. This generalization of \cite{BM99} would be made quite easy if we had a coherence theorem along the lines of Lemma~\ref{coherencelemma} for Verity double categories.

\section{Homotopy categories for $k$-fold quasicategories \label{homotopysection}}
In this section, we will show how to construct from a inner-Kan $k$-simplicial set $X$ an ``underlying'' $m$-reduced inner-Kan $k$-simplicial set $h_m X$, which is obtained by applying a homotopy relation to the cells of dimension $m$ in $X$. This construction is given in the case $k=1$ in Subsection~2.3.4 of \cite{Lur09}, and we follow Lurie's construction to some extent, though his exposition is expedited by the use of some homotopy theory developed in other parts of his book, which we do not generalize in this thesis.
\subsection{The relation of homotopy rel. boundary.}
\begin{definition}   Let $x$ be a $\n$-cell in $X$ with $\dim(\n)>0$. The \emph{standard $\alpha$-squished horn of $x$}, denoted $\squish^\alpha(x)$ is the subhorn of type $\fk:\Lambda^{\n+\e^\alpha}_{\{ d^\alpha_1\}}\ra X$  of the degenerate cell $s^\alpha_0 x.$
\end{definition}
\begin{proposition} \label{daggericonditions} Suppose $X$ is inner-Kan $k$-simplicial set, and $x$ and $x'$ are $\n$-cells of $X$, with $\dim (\n) > 0$. Let $\alpha$ be an index for which $ \n^\alpha >0$. Then the following are equivalent:
\begin{enumerate}
\item There exists a filler $H$  of the standard squished horn $\squish^\alpha (x)$ such that $d^\alpha_1 H = x'$
\item For every (not necessarily inner) horn $\fk$ of type $\Lambda^{\n+\e^\beta}_{\{d^\beta_i\}}$, the following conditions are equivalent:
\begin{enumerate}
\item $\fk$ has a filler $\hat{\fk}$  with $d^\beta_i\hat{\fk} = x$ 
\item $\fk$ has a filler $\hat{\fk}'$  with $d^\beta_i\hat{\fk}' = x'$
\end{enumerate} 
\item There is some $\beta$ and some \emph{inner} horn $\fk$ of type $\Lambda^{\n+\e^\beta}_{\{d^\beta_i\}}$ such that $\fk$ has fillers $\hat{\fk}$ and $\hat{\fk}'$ with $d^\beta_i\hat{\fk} = x$ and $d^\beta_i\hat{\fk}' = x'$
\end{enumerate}
If these equivalent conditions are met, we say \emph{$x$ is homotopic to $x'$ rel. boundary} which we write $x\sim x'$  (in particular $x \sim x'$ implies $dx =d x'$). 
\end{proposition}
\begin{proof} We first show that if $(1)$, then $(\mbox{a}) \Rightarrow (\mathrm{b}).$

We consider three cases, based on the type of $\fk$. First, suppose $\beta \neq \alpha$. 

Take $\hat{\fk}$ to be a filler of $\fk$ with $d^\beta_i \hat{\fk}= x,$ whose existence is guaranteed by the hypothesis $(\mathrm{a}).$  Consider the cell $s^{\alpha}_0\hat{\fk},$ whose faces in the $\alpha$ direction are $$d^\alpha(\hat{\fk})=[\hat{\fk},\ \hat{\fk},\ s^\alpha_0 d^\alpha_1 \hat{\fk},\ s^\alpha_0d^\alpha_2\hat{\fk},\ \ldots,\ s^\alpha_0d^\alpha_m \hat{\fk}].$$ Below we give the $d^\beta d^\alpha$ and $d^\alpha d^\beta$ portion of the Glenn table of this cell.

\begin{table}[H]\caption{$d^\beta d^\alpha$ part of the Glenn table for $s^{\alpha}_0\hat{\fk}$}\label{simequivproof1} \begin{center}
    \begin{tabular}{ r | l || !{\vrule width 2pt} l | l | l | l | l | l | l  !{\vrule width 2pt} }
    \cline{2-5} \cline{7-7} \cline{9-9}

                      & $\hat{\fk}$  & $d^\beta_0 \hat{\fk}$   &   $d^\beta_1 \hat{\fk}$    & $d^\beta_2 \hat{\fk}$             & $\cdots$   &   $d^\beta_i \hat{\fk}=x$      & $\cdots$ &  $d^\beta_{\n^\beta+1} \hat{\fk}$     \\ \cline{2-5} \cline{7-7} \cline{9-9}
     
                      & $\hat{\fk}$  &  $d^\beta_0 \hat{\fk}$   & $d^\beta_1 \hat{\fk}$   & $d^\beta_2 \hat{\fk}$        &   &   $d^\beta_i \hat{\fk}=x$       &        &  $d^\beta_{\n^\beta+1} \hat{\fk}$    \\ \cline{2-5} \cline{7-7} \cline{9-9}

                      & $s^\alpha_0 d^\alpha_1 \hat{\fk}$    &  $s^\alpha_0 d^\beta_0 d^\alpha_1 \hat{\fk}$  &  $s^\alpha_0 d^\beta_1 d^\alpha_1 \hat{\fk}$         &  $s^\alpha_0 d^\beta_2 d^\alpha_1 \hat{\fk}$    &       &$s^\alpha_0 d^\beta_i d^\alpha_1 \hat{\fk}$   &     & $s^\alpha_0 d^\beta_{\n^\beta+1} d^\alpha_1 \hat{\fk}$             \\  \cline{2-5} \cline{7-7} \cline{9-9}

    \multicolumn{1}{r}{}   &   \multicolumn{1}{c}{$\vdots$}     & \multicolumn{1}{r}{} &  \multicolumn{1}{r}{}     & \multicolumn{1}{c}{}       &  \multicolumn{1}{c}{$\ddots$}    & \multicolumn{1}{r}{} & \multicolumn{1}{r}{$\ddots$}  & \multicolumn{1}{c}{$\vdots$}         \\ \cline{2-5} \cline{7-7} \cline{9-9}

      &   $s^\alpha_0 d^\alpha_{\n^\alpha} \hat{\fk}$   &  $s^\alpha_0 d^\beta_0 d^\alpha_{\n^\alpha} \hat{\fk}$  &  $s^\alpha_0 d^\beta_1 d^\alpha_{\n^\alpha} \hat{\fk}$         &  $s^\alpha_0 d^\beta_2 d^\alpha_{\n^\alpha} \hat{\fk}$    &  $\cdots$      &$s^\alpha_0 d^\beta_i d^\alpha_{\n^\alpha} \hat{\fk}$   &     $\cdots$ & $s^\alpha_0 d^\beta_{\n^\beta+1} d^\alpha_{\n^\alpha} \hat{\fk}$             \\  \cline{2-5} \cline{7-7} \cline{9-9}  
    \end{tabular}\end{center}  
    \end{table}

\begin{table}[H]\caption{$d^\alpha d^\beta$ part of the Glenn table for $s^{\alpha}_0\hat{\fk}$}\label{simequivproof2} \begin{center}
    \begin{tabular}{ r | l || !{\vrule width 2pt} l | l | l  | l | l  !{\vrule width 2pt} }
    \cline{2-5} \cline{7-7}

                      & $d^\beta_0 s^{\alpha}_0 \hat{\fk}$  &  $d^\beta_0 \hat{\fk}$    &  $d^\beta_0 \hat{\fk}$  & $s^\alpha_0  d^\alpha_1 d^\beta_0 \hat{\fk}$          & $\cdots$ &  $s^\alpha_0 d^\alpha_{\n^\alpha} d^\beta_0  \hat{\fk}$     \\ \cline{2-5} \cline{7-7} 
     
                      & $d^\beta_1 s^{\alpha}_0 \hat{\fk}$  &  $d^\beta_1 \hat{\fk}$    &  $d^\beta_1 \hat{\fk}$  & $s^\alpha_0  d^\alpha_1 d^\beta_1 \hat{\fk}$           & &  $s^\alpha_0 d^\alpha_{\n^\alpha} d^\beta_1  \hat{\fk}$     \\ \cline{2-5} \cline{7-7}

    \multicolumn{1}{r}{}   &   \multicolumn{1}{c}{$\vdots$}     & \multicolumn{1}{r}{} &  \multicolumn{1}{r}{}       & \multicolumn{1}{r}{} & \multicolumn{1}{r}{}  & \multicolumn{1}{c}{$\vdots$}        \\ \cline{2-5} \cline{7-7}

      &   $d^\beta_i s^\alpha_0  \hat{\fk}=s^\alpha_0 x$   &  $x$ &   $x$ &    $s^\alpha_0  d^\alpha_1 x$&    &     $s^\alpha_0 d^\alpha_{\n^\alpha} x$     \\   \cline{2-5} \cline{7-7} 
      
    \multicolumn{1}{r}{}   &   \multicolumn{1}{c}{$\vdots$}     & \multicolumn{1}{r}{} &  \multicolumn{1}{l}{}    & \multicolumn{1}{r}{} & \multicolumn{1}{r}{$\ddots$}  & \multicolumn{1}{c}{$\vdots$}         \\ \cline{2-5} \cline{7-7}

                      & $d^\beta_{\n^\beta+1} s^{\alpha}_0 \hat{\fk}$  &  $d^\beta_{\n^\beta+1} \hat{\fk}$    &  $d^\beta_{\n^\beta
                      1} \hat{\fk}$  & $s^\alpha_0  d^\alpha_1 d^\beta_{\n^\beta+1} \hat{\fk}$           & $\cdots$ &  $s^\alpha_0 d^\alpha_{\n^\alpha} d^\beta_{\n^\beta+1}  \hat{\fk}$     \\ \cline{2-5} \cline{7-7}  
    \end{tabular}\end{center}  
    \end{table}

 Now consider the $ \Lambda^{\n+\e^\beta+\e^\alpha}_{\{d^\alpha_{1}\}}$-subhorn $L$ contained inside $s^\alpha_0 \hat{\fk}$, obtained from removing the face $d^\alpha_1 s^\alpha_0 \hat{\fk} = \hat{\fk} $ in Table~\ref{simequivproof1}. We can make a new horn $L'$ by replacing the $d^\beta_i$ face of $L$, (which is $s^\alpha_0 x$ in Table~\ref{simequivproof2}), with the filler $H$ for $\squish^\alpha(x)$ whose existence is guaranteed by the hypothesis. Since $H$ agrees with $s^\alpha_0 x$ except on $d^\alpha_1$, and since $d^{\alpha}_1 d^\beta_i=d^\beta_i d^\alpha_1 $, this non-matching face lies along the ``missing face'' of the horn. Below we give the $d^\beta d^\alpha$ portion of the Glenn table of $L'$.

\begin{table}[H]\caption{$d^\beta d^\alpha$ part of the Glenn table for $L'$}\label{simequivproof3} \begin{center}
    \begin{tabular}{ r | l || !{\vrule width 2pt} l | l | l | l | l | l | l  !{\vrule width 2pt} }
    \cline{2-5} \cline{7-7} \cline{9-9}

                          & $\hat{\fk}$  & $d^\beta_0 \hat{\fk}$   &   $d^\beta_1 \hat{\fk}$    & $d^\beta_2 \hat{\fk}$             & $\cdots$   &   $x$      & $\cdots$ &  $d^\beta_{\n^\beta+1} \hat{\fk}$     \\ \cline{2-5} \cline{7-7} \cline{9-9}
     
       $\Lambda$          &   $=:\hat{\fk}'$ &  $d^\beta_0 \hat{\fk}$   & $d^\beta_1 \hat{\fk}$   & $d^\beta_2 \hat{\fk}$        &   &   $x'$       &        &  $d^\beta_{\n^\beta+1} \hat{\fk}$    \\ \cline{2-5} \cline{7-7} \cline{9-9}

                      & $s^\alpha_0 d^\alpha_1 \hat{\fk}$    &  $s^\alpha_0 d^\beta_0 d^\alpha_1 \hat{\fk}$  &  $s^\alpha_0 d^\beta_1 d^\alpha_1 \hat{\fk}$         &  $s^\alpha_0 d^\beta_2 d^\alpha_1 \hat{\fk}$    &       &$s^\alpha_0 d^\beta_i d^\alpha_1 \hat{\fk}$   &     & $s^\alpha_0 d^\beta_{\n^\beta+1} d^\alpha_1 \hat{\fk}$             \\  \cline{2-5} \cline{7-7} \cline{9-9}

    \multicolumn{1}{r}{}   &   \multicolumn{1}{c}{$\vdots$}     & \multicolumn{1}{r}{} &  \multicolumn{1}{r}{}     & \multicolumn{1}{c}{}       &  \multicolumn{1}{c}{$\ddots$}    & \multicolumn{1}{r}{} & \multicolumn{1}{r}{$\ddots$}  & \multicolumn{1}{c}{$\vdots$}         \\ \cline{2-5} \cline{7-7} \cline{9-9}

      &   $s^\alpha_0 d^\alpha_{\n^\alpha} \hat{\fk}$   &  $s^\alpha_0 d^\beta_0 d^\alpha_{\n^\alpha} \hat{\fk}$  &  $s^\alpha_0 d^\beta_1 d^\alpha_{\n^\alpha} \hat{\fk}$         &  $s^\alpha_0 d^\beta_2 d^\alpha_{\n^\alpha} \hat{\fk}$    &  $\cdots$      &$s^\alpha_0 d^\beta_i d^\alpha_{\n^\alpha} \hat{\fk}$   &     $\cdots$ & $s^\alpha_0 d^\beta_{\n^\beta+1} d^\alpha_{\n^\alpha} \hat{\fk}$             \\  \cline{2-5} \cline{7-7} \cline{9-9}  
    \end{tabular}\end{center}  
    \end{table}
We define $ \hat{\fk}'$ as the $d^\alpha_1$ face of a filler of $L'$, as shown in Table~\ref{simequivproof3}. Clearly $\hat{\fk}'$ agrees with $\hat{\fk}$ except on $d^\beta_i$, so $\hat{\fk}'$ is a filler of $\fk$ with $d^\beta_i\hat{\fk}' = x',$ verifying $(\mathrm{b})$.

Next we consider the case where $\beta=\alpha$ and $i>0$, so that $\fk$ is of type $\Lambda^{\n+\e^\alpha}_{\{d^\alpha_i\}},$ and take $\hat{\fk}$ to be a filler of $\fk$ with $d^\alpha_i \hat{\fk}=x.$  We proceed similarly to above, taking $L$ to be the $ \Lambda^{\n+\e^\alpha+\e^\alpha}_{\{d^\alpha_1\}}$-subhorn contained inside $s^\alpha_0 \hat{\fk}$, obtained from removing the face $d^\alpha_1 s^\alpha_0 \hat{\fk} = \hat{\fk} $ from $L$. We define a new horn $L'$ by replacing the $d^\alpha_{i+1}$ face of $L$, which is the cell $$d^\alpha_{i+1} s^\alpha_0 \fk= s^\alpha_0 d^\alpha_i \fk = s^\alpha_0 x,$$ with the filler $H$ for $\squish^\alpha(x)$. $H$ matches $s^\alpha_0 x$ except on $d^\alpha_1$, and since $i>0$ we have $ d^{\alpha}_1 d^\alpha_{i+1}=d^\alpha_i d^\alpha_1$, and so this non-matching face lies along the ``missing face'' of the horn. We give the $d^\alpha d^\alpha$ portion of the Glenn table of $L'$.

\begin{table}[H]\caption{$d^\alpha d^\alpha$ part of the Glenn table for $L'$}\label{daggerconditionsproof4} \begin{center}
    \begin{tabular}{ r | l || !{\vrule width 2pt} l | l | l | l | l | l | l  !{\vrule width 2pt} }
    \cline{2-5} \cline{7-7} \cline{9-9}

                      & $\hat{\fk}$  & $d^\alpha_0 \hat{\fk}$   &   $d^\alpha_1 \hat{\fk}$    & $d^\alpha_2 \hat{\fk}$             & $\cdots$   &   $x$      & $\cdots$ &  $d^\alpha_{\n^\alpha+1} \hat{\fk}$     \\ \cline{2-5} \cline{7-7} \cline{9-9}
     
    $\Lambda$         & $=:\hat{\fk}'$  &  $d^\alpha_0 \hat{\fk}$   & $d^\alpha_1 \hat{\fk}$   & $d^\alpha_2 \hat{\fk}$        &   &   $x'$       &        &  $d^\alpha_{\n^\alpha+1} \hat{\fk}$    \\ \cline{2-5} \cline{7-7} \cline{9-9}

                      & $s^\alpha_0 d^\alpha_1 \hat{\fk}$    &  $d^\alpha_1 \hat{\fk}$  &  $d^\alpha_1 \hat{\fk}$         &  $s^\alpha_0 d^\alpha_1 d^\alpha_1 \hat{\fk}$     &     &$s^\alpha_0 d^\alpha_i d^\alpha_1 \hat{\fk}$  & & $s^\alpha_0 d^\alpha_{\n^\alpha} d^\alpha_1 \hat{\fk}$               \\  \cline{2-5} \cline{7-7} \cline{9-9}
     
    \multicolumn{1}{r}{}   &   \multicolumn{1}{c}{$\vdots$}     & \multicolumn{1}{r}{} &  \multicolumn{1}{r}{}     & \multicolumn{1}{r}{}       &  \multicolumn{1}{l}{$\ddots$}    & \multicolumn{1}{r}{} & \multicolumn{1}{r}{}  & \multicolumn{1}{c}{$\vdots$}        \\ \cline{2-5} \cline{7-7} \cline{9-9}

             & $H$     &  $x$ &   $x'$ &    $s^\alpha_0 d^\alpha_1 d^\alpha_i \hat{\fk}$ &     &  $s^\alpha_0 d^\alpha_i d^\alpha_i \hat{\fk}$ & &    $s^\alpha_0 d^\alpha_{\n^\alpha} d^\alpha_i \hat{\fk}$      \\   \cline{2-5} \cline{7-7} \cline{9-9}
      
    \multicolumn{1}{r}{}   &   \multicolumn{1}{c}{$\vdots$}     & \multicolumn{1}{r}{} &  \multicolumn{1}{r}{}     & \multicolumn{1}{r}{}       &  \multicolumn{1}{l}{}    & \multicolumn{1}{r}{} & \multicolumn{1}{r}{$\ddots$}  & \multicolumn{1}{c}{$\vdots$}         \\ \cline{2-5} \cline{7-7} \cline{9-9}

      &   $s^\alpha_0 d^\alpha_{\n^\alpha+1} \hat{\fk}$   &  $d^\alpha_{\n^\alpha+1} \hat{\fk}$ &  $d^\alpha_{\n^\alpha+1} \hat{\fk}$  &     $s^\alpha_0 d^\alpha_1 d^\alpha_{\n^\alpha+1} \hat{\fk}$     &  $\cdots$    & $s^\alpha_0 d^\alpha_i d^\alpha_{\n^\alpha+1} \hat{\fk}$ &  $\cdots$  &   $s^\alpha_0 d^\alpha_{\n^\alpha} d^\alpha_{\n^\alpha+1} \hat{\fk}$      \\   \cline{2-5} \cline{7-7} \cline{9-9}      
    \end{tabular}\end{center}  
    \end{table}
Define $ \hat{\fk}'$ as the $d^\alpha_1$ face of a filler of $L'$, as shown in Table~\ref{daggerconditionsproof4}. $\hat{\fk}'$ matches $\hat{\fk}$ except on $d^\alpha_i$, so $\hat{\fk}'$ is a filler of $\fk$ with $d^\alpha_i\hat{\fk}' = x',$ showing $(\mathrm{b})$. 

In the final case $\hat{\fk}$ is a face of type $d^\alpha_0.$ In this case we proceed as above, except taking $L$ to be the $ \Lambda^{\m+\e^\alpha+\e^\alpha}_{\{d^\alpha_{2}\}}$-subhorn contained inside $s^\alpha_1 \hat{\fk}.$ In this case we replace the face 
$d^\alpha_0 L$, which is the cell $$d^\alpha_0 s^\alpha_1 \fk=s^\alpha_0 d^\alpha_0 \fk =s^\alpha_0 x,$$ with the filler $H$ for $\squish^\alpha(x).$ $H$ matches $s^\alpha_0 x$ except on $d^\alpha_1$, and we have $ d^{\alpha}_1 d^\alpha_{0}=d^\alpha_0 d^\alpha_2$, and so this non-matching face lies along the ``missing face'' of the horn. Similarly to above, define $\hat{\fk}'$ as the $d^\alpha_2$ face of a filler of $L'$. $\hat{\fk}'$ matches $\hat{\fk}$ except on $d^\alpha_0$, so $\hat{\fk}'$ is a filler of $\fk$ with $d^\alpha_0\hat{\fk}'=d^\alpha_1 H = x',$ showing $(\mathrm{b})$.

We have shown $(1) \Rightarrow ( (\mathrm{a}) \Rightarrow (\mathrm{b})).$ Since $dx = dx'$, note that $s^\alpha_0 x$  matches $s^\alpha_0 x'$  except on $d^\alpha_0$ and $d^\alpha_1.$ Thus replacing $d^\alpha_0 \squish^\alpha(x)$ by $x'$ yields $\squish^\alpha(x')$. Therefore if we have $(\mathrm{a}) \Rightarrow (\mathrm{b}),$ we can apply it with $\fk$ taken to be the $\Lambda^{\n+\e^\alpha}_{\{d^\alpha_0\}}$ subhorn of $\squish^\alpha(x)$ and $\hat{\fk}=s^\alpha_0 x$, yielding a cell $\hat{\fk}'$ which is also a filler of $\squish^\alpha(x')$ with $d^\alpha_1 \hat{\fk}'=x.$ That is, assuming $(\mathrm{a}) \Rightarrow (\mathrm{b}),$ we have the reverse of $(1)$, with the roles of $x$ and $x'$ switched. Since $(1) \Rightarrow ( (\mathrm{a}) \Rightarrow (\mathrm{b})),$ this reverse of $(1)$ implies that $(b) \Rightarrow (a)$. The argument used above can be applied to show $(b) \Rightarrow (a)$ implies $(1)$. So we have  the following implications
\begin{itemize}
\item $(1)$ implies $(\mathrm{a}) \Rightarrow (\mathrm{b})$
\item $(\mathrm{a}) \Rightarrow (\mathrm{b})$ implies the reverse  of $(1)$, with the roles of $x$ and $x'$ switched.
\item The reverse of $(1)$ implies $ (\mathrm{b}) \Rightarrow (\mathrm{a})$
\item  $ (\mathrm{b}) \Rightarrow (\mathrm{a})$ implies $(1).$
\end{itemize}
Thus these four conditions are equivalent and therefore $(1) \Leftrightarrow (2).$  This is enough to show that condition $(1)$ is equivalent for different values of $\alpha$. 

$(1)$ is a special case of $(3)$ thus $(1) \Rightarrow (3).$ Now we show $(3) \Rightarrow (1)$. Suppose we have $(3)$ and $\fk$ is a horn of type $\Lambda^{\n+\e^\beta}_{\{d^\beta_i\}}$ such that $\fk$ has fillers $\hat{\fk}$ and $\hat{\fk}'$ with $d^\beta_i\hat{\fk} = x$ and $d^\beta_i\hat{\fk}' = x'.$  Since $(1)$ is equivalent for different values of $\alpha$, and since we must have $\n^\beta>0$ from the fact that there is an inner horn of this type, we can assume $\alpha = \beta$.

Consider the cell $s^{\alpha}_0\hat{\fk},$ whose faces in the $\alpha$ direction are $$d^\alpha(\hat{\fk})=[\hat{\fk},\ \hat{\fk},\ s^\alpha_0 d^\alpha_1 \hat{\fk},\ s^\alpha_0d^\alpha_2\hat{\fk},\ \ldots,\ s^\alpha_0d^\alpha_{\n^\alpha+1} \hat{\fk}].$$ 
Take $L$ to be the $ \Lambda^{\n+\e^\alpha+\e^\alpha}_{\{d^\alpha_{i+1}\}}$-subhorn contained inside $s^\alpha_0 \hat{\fk}$, obtained from removing the face $s^\alpha_0 d^\alpha_i\hat{\fk} = s^\alpha_0 x $ in the list of faces above. We can make a new horn $L'$ by replacing the $d^\alpha_1$ face of $L$, (which is $\hat{\fk}$), with $\hat{\fk}'$, since $d\hat{\fk}'$ agrees with $d\hat{\fk}$ except on $d^\alpha_i$, and since $d^\alpha_id^\alpha_1 =d^{\alpha}_1d^\alpha_{i+1}$ (using $i\geq 1$), this non-matching face lies along the ``missing face'' of the horn. Below we give the $d^\alpha_i d^\alpha_j$ portion of the Glenn table of $L'$. 

\begin{table}[H]\caption{The $d^\alpha d^\alpha$ part of the Glenn table for $L'$}\label{daggerconditionsproof6} \begin{center}
    \begin{tabular}{ r | l || !{\vrule width 2pt} l | l | l | l | l | l | l  !{\vrule width 2pt} }
    \cline{2-5} \cline{7-7} \cline{9-9}

                      & $\hat{\fk}$  & $d^\alpha_0 \hat{\fk}$   &   $d^\alpha_1 \hat{\fk}$    & $d^\alpha_2 \hat{\fk}$             & $\cdots$   &   $x$      & $\cdots$ &  $d^\alpha_{\n^\alpha+1} \hat{\fk}$     \\ \cline{2-5} \cline{7-7} \cline{9-9}
     
                      & $\hat{\fk}'$  &  $d^\alpha_0 \hat{\fk}$   & $d^\alpha_1 \hat{\fk}$   & $d^\alpha_2 \hat{\fk}$        &   &   $x'$       &        &  $d^\alpha_{\n^\alpha+1} \hat{\fk}$    \\ \cline{2-5} \cline{7-7} \cline{9-9}

                      & $s^\alpha_0 d^\alpha_1 \hat{\fk}$    &  $d^\alpha_1 \hat{\fk}$  &  $d^\alpha_1 \hat{\fk}$         &  $s^\alpha_0 d^\alpha_1 d^\alpha_1 \hat{\fk}$     &     &$s^\alpha_0 d^\alpha_i d^\alpha_1 \hat{\fk}$  & & $s^\alpha_0 d^\alpha_{\n^\alpha} d^\alpha_1 \hat{\fk}$               \\  \cline{2-5} \cline{7-7} \cline{9-9}
     
    \multicolumn{1}{r}{}   &   \multicolumn{1}{c}{$\vdots$}     & \multicolumn{1}{r}{} &  \multicolumn{1}{r}{}     & \multicolumn{1}{r}{}       &  \multicolumn{1}{l}{$\ddots$}    & \multicolumn{1}{r}{} & \multicolumn{1}{r}{}  & \multicolumn{1}{c}{$\vdots$}        \\ \cline{2-5} \cline{7-7} \cline{9-9}

   $\Lambda$  &   $=:H$   &  $x$ &   $x'$ &    $s^\alpha_0 d^\alpha_1 d^\alpha_i \hat{\fk}$ &     &  $s^\alpha_0 d^\alpha_i d^\alpha_i \hat{\fk}$ & &    $s^\alpha_0 d^\alpha_{\n^\alpha} d^\alpha_i \hat{\fk}$      \\   \cline{2-5} \cline{7-7} \cline{9-9}
      
    \multicolumn{1}{r}{}   &   \multicolumn{1}{c}{$\vdots$}     & \multicolumn{1}{r}{} &  \multicolumn{1}{r}{}     & \multicolumn{1}{r}{}       &  \multicolumn{1}{l}{}    & \multicolumn{1}{r}{} & \multicolumn{1}{r}{$\ddots$}  & \multicolumn{1}{c}{$\vdots$}         \\ \cline{2-5} \cline{7-7} \cline{9-9}

      &   $s^\alpha_0 d^\alpha_{\n^\alpha+1} \hat{\fk}$   &  $d^\alpha_{\n^\alpha+1} \hat{\fk}$ &  $d^\alpha_{\n^\alpha+1} \hat{\fk}$  &     $s^\alpha_0 d^\alpha_1 d^\alpha_{\n^\alpha+1} \hat{\fk}$     &  $\cdots$    & $s^\alpha_0 d^\alpha_i d^\alpha_{\n^\alpha+1} \hat{\fk}$ &  $\cdots$  &   $s^\alpha_0 d^\alpha_{\n^\alpha} d^\alpha_{\n^\alpha+1} \hat{\fk}$      \\   \cline{2-5} \cline{7-7} \cline{9-9}      
    \end{tabular}\end{center}  
    \end{table}
  
By filling $L'$ with a filler $\widehat{L'}$ we get a filler $H:=d^\alpha_i(\widehat{L'})$ of the sphere which is ``empty face'' of $L'$. The Glenn table shown in Table~\ref{daggerconditionsproof6} makes it easy to see that $H$ is a filler of $\squish^\alpha(x)$ with $d^\alpha_1 H=x'$. This shows $(1)$ holds. We have now shown $(1) \Leftrightarrow (2)$ and $(1) \Leftrightarrow (3),$ completing the proof.
\end{proof}

\begin{proposition} \label{simequivprop} Let $X$ be a inner-Kan $k$-simplicial set. Then the homotopy rel. boundary relation $\sim$ on the $\n$-cells of $x$ which is given by the equivalent conditions of Proposition~\ref{daggericonditions} is an equivalence relation.
\end{proposition}
\begin{proof} Characterization $(2)$ of $\sim$ in Proposition~\ref{daggericonditions} clearly gives a reflexive and symmetric relation. 

We now show $\sim$ is transitive using characterization $(1)$. Suppose $x \sim x'$ and $x' \sim x''$. Then we have a filler $H$ of $\squish^\alpha(x)$ and a filler $H'$ of $\squish^\alpha(x')$ such that $d_1^\alpha(H)=x'$ and $d_1^\alpha(H')=x''.$ Let $L$ be a $\Lambda^{\n+\e^\alpha+\e^\alpha}_{\{d^\alpha_1\}}$-horn with $d^\beta_i L := s^\alpha_0 s^\alpha_0 d^\beta_i x$ for $\beta \neq \alpha$ and $$d^\alpha L := [H,-,H',s^\alpha_0 s^\alpha_0 d^\alpha_1 x, s^\alpha_0 s^\alpha_0 d^\alpha_2 x, \ldots,  s^\alpha_0 s^\alpha_0 d^\alpha_{\n^\alpha} x].$$ We need to check the $k$-simplicial identities for $L$. At the same time we will verify the identites that check that the ``empty face'' of $L$, which (slightly abusing notation) we write as $d^\alpha_1 L $, matches $\squish^\alpha(x)$ except on its $d^\alpha_1$ face. First we give the $d^\alpha d^\alpha$ portion of the Glenn table for $L$, to help check these identities for $L$.

\begin{table}[H]\caption{The $d^\alpha d^\alpha$ part of the Glenn table for $L$}\label{transitivetable1} \begin{center}
    \begin{tabular}{ r | l || !{\vrule width 2pt} l | l | l | l  | l | l  !{\vrule width 2pt} }
    \cline{2-6} \cline{8-8}

                      & $H$  & $x$   &   $x'$    & $s^\alpha_0 d^\alpha_1 x$             & $s^\alpha_0 d^\alpha_2 x$         & $\cdots$ & $s^\alpha_0 d^\alpha_{\n^\alpha} x$     \\ \cline{2-6}\cline{8-8}
     
      $\Lambda \!$  &   &  $x$   & $x''$       & $s^\alpha_0 d^\alpha_1 x$            &   $s^\alpha_0 d^\alpha_2 x$           &        &  $s^\alpha_0 d^\alpha_{\n^\alpha} x$    \\ \cline{2-6} \cline{8-8}

                      & $H'$    &  $x'$  &  $x''$         &$s^\alpha_0 d^\alpha_1 x$        &    $s^\alpha_0 d^\alpha_2 x$        & & $s^\alpha_0 d^\alpha_{\n^\alpha} x$            \\  \cline{2-6} \cline{8-8}
     
                     &  $ s^\alpha_0 s^\alpha_0 d^\alpha_1 x$        &  $d^\alpha_0 s^\alpha_0 s^\alpha_0 d^\alpha_1 x$ &  $d^\alpha_1 s^\alpha_0 s^\alpha_0 d^\alpha_1 x$         &$d^\alpha_2 s^\alpha_0 s^\alpha_0 d^\alpha_1 x$       &    $d^\alpha_3 s^\alpha_0 s^\alpha_0 d^\alpha_1 x$      & & $d^\alpha_{\n^\alpha} s^\alpha_0 s^\alpha_0 d^\alpha_1 x$              \\  \cline{2-6} \cline{8-8}

                     &  $ s^\alpha_0 s^\alpha_0 d^\alpha_2 x$        &  $d^\alpha_0 s^\alpha_0 s^\alpha_0 d^\alpha_2 x$ &  $d^\alpha_1 s^\alpha_0 s^\alpha_0 d^\alpha_2 x$         &$d^\alpha_2 s^\alpha_0 s^\alpha_0 d^\alpha_2 x$       &    $d^\alpha_3 s^\alpha_0 s^\alpha_0 d^\alpha_2 x$       & & $d^\alpha_{\n^\alpha} s^\alpha_0 s^\alpha_0 d^\alpha_2 x$              \\  \cline{2-6} \cline{8-8}
          
    \multicolumn{1}{r}{}   &   \multicolumn{1}{c}{$\vdots$}     & \multicolumn{1}{r}{} &  \multicolumn{1}{r}{}     & \multicolumn{1}{r}{}       & \multicolumn{1}{r}{} & \multicolumn{1}{r}{$\ddots$}  & \multicolumn{1}{c}{$\vdots$}         \\ \cline{2-6} \cline{8-8}

                     &  $ s^\alpha_0 s^\alpha_0 d^\alpha_{\n^\alpha} x$        &  $d^\alpha_0 s^\alpha_0 s^\alpha_0 d^\alpha_{\n^\alpha} x \!$ &  $d^\alpha_1 s^\alpha_0 s^\alpha_0 d^\alpha_{\n^\alpha} x$         &$d^\alpha_2 s^\alpha_0 s^\alpha_0 d^\alpha_{\n^\alpha} x$       &    $d^\alpha_3 s^\alpha_0 s^\alpha_0 d^\alpha_{\n^\alpha} x$     & $\cdots$ & $d^\alpha_{\n^\alpha} s^\alpha_0 s^\alpha_0 d^\alpha_{\n^\alpha} x$              \\  \cline{2-6} \cline{8-8}
      
    \end{tabular}\end{center}  
    \end{table}
    
The simplicial identities for $d^\alpha_i d^\alpha_j L $ where $j\leq 2$ and $i\leq 1$ are immediately seen to hold from the six entries at the top left of Table~\ref{transitivetable1}. Next the identities $d^\alpha_i d^\alpha_0 L = d^\alpha_0 d^\alpha_{i+1} L$  and $d^\alpha_i d^\alpha_1 L = d^\alpha_1 d^\alpha_{i+1} L$ with $i>1$ assert $$s^\alpha_0 d^\alpha_i x = d^\alpha_0 s^\alpha_0 s^\alpha_0 d^\alpha_i x  \ \ \  \ \ s^\alpha_0 d^\alpha_i x = d^\alpha_1 s^\alpha_0 s^\alpha_0 d^\alpha_i x$$ which follow from $ d^\alpha_0 s^\alpha_0 = d^\alpha_1 s^\alpha_0=\id.$ The identities $d^\alpha_i d^\alpha_2 L = d^\alpha_2 d^\alpha_{i+1} L$ with $i>1$ assert $$s^\alpha_0 d^\alpha_i x = d^\alpha_2 s^\alpha_0 s^\alpha_0 d^\alpha_i x  $$ which can be seen from the simplicial identities: $$ d^\alpha_2 s^\alpha_0 s^\alpha_0 d^\alpha_i =s^\alpha_0 d^\alpha_1 s^\alpha_0 d^\alpha_i=s^\alpha_0 d^\alpha_i.$$ We have left to consider the simplicial identites  $d^\alpha_i d^\alpha_j L = d^\alpha_{j-1} d^\alpha_i$ with $j > i >2$.  These assert $$ d^\alpha_i s^\alpha_0 s^\alpha_0 d^\alpha_{j-2} x = d^\alpha_{j-1} s^\alpha_0 s^\alpha_0 d^\alpha_{i-2} x $$ We apply simplicial identites to each side to move the $d^\alpha$ terms left, using the fact that $j-2>0$ and $i-2>0$ to see this statement is equivalent to $$ d^\alpha_i d^\alpha_{j}s^\alpha_0 s^\alpha_0  x = d^\alpha_{j-1} d^\alpha_{i} s^\alpha_0 s^\alpha_0 x$$ which clearly holds.

The identities for $d^\beta_i d^{\beta'}_j L$ with $\beta, \beta' \neq \alpha$  follow  from the definition $d^{\beta'}_i L = s^\alpha_0 s^\alpha_0 d^{\beta'}_i x$ and the fact that $d^\beta_i$ and $d^{\beta'}_j$ commute with $s^\alpha_0$. We have left to show the identities of the form $d^\beta_i d^\alpha_j L = d^\alpha_j d^\beta_i L$ for $\beta \neq \alpha$, which assert $$d^\beta_i d^\alpha_j L=d^\alpha_j s^\alpha_0 s^\alpha_0 d^\beta_i x.$$ 
First consider $j=0$. We have $d^\beta_i d^\alpha_0 L = d^\beta_i H$ and since $H$ matches $\squish^\alpha(x)$ and thus $s^{\alpha}_0 x$ on $d^\beta$ faces we have   $d^\beta_i s^\alpha_0 x=s^\alpha_0  d^\beta_i x.$ 
So $$d^\beta_i d^\alpha_0 L=s^\alpha_0  d^\beta_i x=d^\alpha_0 s^\alpha_0 s^\alpha_0 d^\beta_i x= d^\alpha_0 d^\beta_i L$$ using the fact that $d^\alpha_0 s^\alpha_0=\id.$ The cases $j=1$ and $j=2$ follow by a similar argument, using the fact that $d x =d x' = d x''.$ For $j>2$ we have $$d^\beta_i d^\alpha_j L = d^\beta_i s^\alpha_0 s^\alpha_0 d^\alpha_{j-2} x = s^\alpha_0 s^\alpha_0 d^\alpha_{j-2}d^\beta_i  x =  d^\alpha_js^\alpha_0 s^\alpha_0 d^\beta_i  x =d^\alpha_j d^\beta_i L.$$ This completes the verification of the $k$-simplicial identities for $L$. 

Because $X$ is inner-Kan, there is a filler $\hat{L}$ of $L$. The ``empty face'' $d^\alpha_1 L$ matches $\squish^\alpha(x)$ except for its $d^\alpha_1$ face, which is $x''$. So $H'':=d^\alpha_1 \hat{L}$ is a filler of  $\squish^\alpha(x)$ with $d^\alpha_1 H''= x''$. Thus $x \sim x''$, and we conclude $\sim$ is transitive and an equivalence relation.  
\end{proof}

\subsection{$\alpha$-homotopy compared with $\sim$ }

We now justify our choice terminology in calling the relation $\sim$ ``homotopy rel. boundary.''

\begin{definition}
Let $X$ be a $k$-simplicial set. Let $\mathbf{n} = [\mathbf{n}^1, \mathbf{n}^2, \ldots, \mathbf{n}^k]$ be a $k$-fold index, with $\dim(\n) > 0 $. Let $x,x':\Delta^k[\n]\ra X$ be two $\mathbf{n}$-cells in $X$, with $dx = dx'$ and let $1 \leq \alpha \leq k$ be such that  $\n^\alpha>0.$  Then $x \sim_\alpha x'$ if there is a map $h:\Delta^k[\mathbf{n}]\times I_\alpha \ra X$  such that:
\begin{itemize}\item $h|_{d\Delta^k[\n]\times I_\alpha }=dx \circ \pi =dx'\circ \pi$, where $\pi: d\Delta^k[\n]\times I_\alpha \ra  d\Delta^k[\n]$ is the projection map.
\item  $h|_0 =x$ and $h|_1=x'$, where $h|_0$ and $h|_1$ are the restrictions of $h$ along the maps $i_0, i_1: \Delta^k[\mathbf{n}]\ra \Delta^k[\mathbf{n}]\times I_\alpha $ induced by the two maps $\Delta^k[\mathbf{0}]\rightrightarrows I_\alpha$.
\end{itemize}
\end{definition}

Let $I_\alpha:=\Delta^k[\e_\alpha]=\Delta^k[0,\ldots,0, 1,0, \ldots 0]. $ 
\begin{definition} \label{homotopydef}Let $x, x'$ be $\n$-cells in a $k$-simplicial set $X$. A \emph{$\alpha$-homotopy} from $x$ to $x'$ is a map 
\end{definition}
$h:\Delta^k[\n]\times I_\alpha \ra X$  such that:
\begin{itemize} 
\item  $h|_0 =x$ and $h|_1=x'$, where $h|_0$ and $h|_1$ are the restrictions of $h$ along the maps $i_0, i_1: \Delta^k[\n]\ra \Delta^k[\n]\times I_\alpha $ induced by the two maps $\Delta^k[\mathbf{0}]\rightrightarrows I_\alpha.$\end{itemize}
If also $h|_{d\Delta^k[\n]\times I_\alpha }=dx \circ \pi =dx'\circ \pi$, where $\pi: d\Delta^k[\n]\times I_\alpha \ra  d\Delta^k[\n]$ is the projection map, then we say that $h$ is \emph{constant on the boundary}. 

First we consider the simplicial case, setting $k=1$. Recall the presentation of the simplicial prism $\Delta[n] \times \Delta[1]$ as a colimit of $(n+1)$-simplicies (details are given for instance in \cite{JT99}). The prism $h=\Delta[n]\times \Delta[1]$ is constructed from $n+1$ copies of $\Delta[n+1]$, denoted $h_0, \ldots, h_n$, with the $(i+1)$-face of $h_i$ glued to the $(i+1)$-face of $h_{i+1}$. Then the inclusions $i_0,i_1 :\Delta[n] \rightrightarrows \Delta[n]\times \Delta[1]$, correspond to the inclusion of the $0$th face in $h_0$ and the $n$th face in $h_n$. The projection $\Delta[n]\times \Delta[1]\ra \Delta[n]$ is given on $h_i$ by the codegeneracy map $\hat{s}_i.$

More formally, $h=\Delta[n]\times\Delta[1]$ is isomorphic to the coequalizer of $$\coprod_{1\leq i \leq n} \Delta[n] \rightrightarrows \coprod_{0\leq j \leq n} \Delta[n+1]$$ where the two maps respectively send the $i$th copy of $\Delta[n]$ to either the $i$th face of the $(i-1)$th copy of $\Delta[n+1]$ or the $i$th face of $i$th face of the $i$th copy of $\Delta[n+1]$.

Now we consider the multi-simplicial case. Using Proposition~\ref{boxysimplicies} and Proposition~\ref{boxytimes} we consider a $k$-simplicial prism 
\begin{align*}\Delta^k[\n]\times\Delta^k[\e^\alpha] & \cong (\Delta[\n^1]\square \cdots \square \Delta[\n^k]) \times (\Delta[0]\square \cdots \Delta[0]\square \Delta[1] \square\Delta[0]\cdots \Delta[0]) \\ 
 &\cong   \Delta[\n^1]\square \cdots \square (\Delta[\n^\alpha] \times \Delta[1])\square \cdots \square \Delta[\n^k]
\end{align*}
Because $\square$ preserves colimits (Proposition~\ref{boxytimes}), the presentation of $\Delta[\n^\alpha] \times \Delta[1]$ given above gives rise to a presentation of the above multi-simplicial prism $\Delta^k[\n]\times\Delta^k[\e^\alpha]$. This prism is formed by a colimit of $\n^\alpha+1$ copies $g_0, \ldots, g_{\n^\alpha}$ of  $\Delta^k[\n+\e^\alpha],$ with the $i$th face of $g_{i-1}$ glued to the $i$th face of $g_i$.  
Then a homotopy $h:\Delta^k[\n]\times\Delta^k[\e^\alpha]$ is equivalent to a collection $h_0, \ldots, h_{\n^\alpha}$ of $(\n+\e^\alpha)$-cells of $X$ such that $d^\alpha_i h_{i-1} = d^\alpha_i h_i$ for $1\leq i \leq \n^\alpha$.

\begin{proposition} \label{alphahomotopyproposition} Let $X$ be a inner-Kan $k$-simplicial set with $\n$ cells $x$ and $x'$, with $\n^\alpha>0$. Then $x \sim_\alpha x'$ if and only if $x \sim x'$. 
\end{proposition}
\begin{proof}
First we assume $x \sim x'.$ Using $(1)$ of Proposition~\ref{daggericonditions} let $H$ be a filler of the squished horn $\squish^\alpha (x)$ such that $d^\alpha_1H = x'$. Now consider the trivial homotopy $h$ of $x'$ in $X$, given by $\Delta^k[\n]\times \Delta^k[\e^\alpha]\stackrel{\pi}{\ra}\Delta^k[\m]\stackrel{x'}{\ra}X.$
The cell $h_i$ of this homotopy is given by $s^\alpha_i x'.$  Clearly this homotopy is constant on its boundary. Because $d H$ matches $d(s^\alpha_0 x')$ except on face $d^\alpha_0$ we can replace $h_0$ in the homotopy $h$ by $\hat{\fk}'$, yielding a homotopy $h'$ (constant on the boundary) from $x$ to $x'$, showing $x \sim_\alpha x'$.

On the other hand, suppose $x \sim_\alpha x'$ and $h$ is a homotopy constant on its boundary from $x$ to $x'$. Consider the bottom cell $h_0$ of $h$. It follows from the fact that $h$ is constant on the boundary that $d(h_0)$ matches $d(s^\alpha_0 x)$ except possibly on the face $d^\alpha_1$. This makes $s^\alpha_0 x$ and $h_0$ fillers of the same inner horn, since $\n^\alpha>0$ so we can conclude $x \sim d^\alpha_1 h_0$ by the fact that $X$ satisfies $(1)$ of Proposition~\ref{daggericonditions}. By a similar argument, we can show $d_i h_{i}  \sim_\alpha d_{i+1} h_i$ for $0 < i \leq \n^\alpha$. So we have $$x=d_0 h_0 \sim d_1 h_0 = d_1 h_1 \sim \ldots \sim d_{\n^\alpha+1} h_{\n^\alpha} = x' .$$ By the transitive property of $\sim$ shown in Proposition~\ref{simequivprop}, we conclude $x \sim x'.$
\end{proof}
\begin{remark} In the case that $\n^\alpha=0$, we must take a different definition of  $\sim_\alpha$ in order for the statement in Proposition~\ref{alphahomotopyproposition} to hold. We must insist that the homotopy is ``invertible'' in an appropriate sense. If the appropriate definition is made, all of the $\sim_\alpha$ relations will be the same as long as $\dim(\n)>0.$ For $0$ cells, the appropriate definition of $\sim_\alpha$ is the relation of being connected by an invertible $1$-cell with multi-index $\e^\alpha$. However, these relations are in general different, and no preferred definition of $\sim$ is possible without further conditions on $X$. 
\end{remark}

\subsection{The homotopy category operation}

Recall from Definition~\ref{coskeletaldef} that a $k$-simplicial set is said to be $m$-subcoskeletal if the natural map $X \ra \Cosk^m(X)$ is a monomorphism. Proposition~\ref{coskeletallemma} shows that $X$ is $m$-subcoskeletal if and only if every $d\Delta^k[\n]$-sphere in $X$ has at most one filler when $\dim(\n)>m.$

\begin{proposition} \label{ostarconds} If $X$ is a inner-Kan $k$-simplicial set and is $m$-subcoskeletal then:
\begin{enumerate} 
\item $X$ is $(m+1)$-coskeletal.
\item  $X$ is $(m+1)$-reduced. 
\end{enumerate}
\end{proposition} 
\begin{proof} $(2)$ follows from Lemma~\ref{coskeletaltofiller}, and $(1)$ follows by Lemma~\ref{fillerstocoskeletal1}.
\end{proof}

\begin{definition} A $k$-fold simplicial set $X$ is called \emph{homotopically $m$-reduced} if it is inner-Kan and $m$-subcoskeletal, and for $\n$-cells $x$, $x'$ in $X$ with $\dim(\n)=m$, $x\sim x'$ implies $x=x'.$ \end{definition}

\begin{proposition} \label{daggerness} Let $X$ be inner-Kan $k$-simplicial set, and let $m\geq k$. Then $X$ is \emph{homotopically $m$-reduced} if and only if it is $m$-reduced.
\end{proposition}
\begin{proof}
First assume $X$ is homotopically $m$-reduced. Then in the case of a horn $H$ of type $\Lambda^{\n}(d^\alpha_i)$ in $X$ with $\dim(\n)=m+1$ if we have two fillers $\hat{H}$ and $\hat{H}'$ of $H$, we have from $(3)$ in Proposition~\ref{daggericonditions} that $d^\alpha_i \hat{H}\sim d^\alpha_i\hat{H}'$ from which we conclude using the hypothesis that $d^\alpha_i \hat{H}= d^\alpha_i\hat{H}'$ and thus $d\hat{H}= d\hat{H}'$. Then since $X$ is $m$-subcoskeletal we have $\hat{H}=\hat{H}'.$ For horns of dimension $m+2$ and greater in $X$ we have unique fillers by Proposition~\ref{ostarconds}. 

Now suppose $X$ is $m$-reduced. If we have $d\Delta^k[\n]$-sphere with $\dim(\n)=m>k$ then there is inner horn contained within it, since $\n^\alpha \geq 2$ for some $\alpha$ by the pigeonhole principle. Any filler of the $d\Delta^k[\n]$-sphere is also a filler of this inner horn, and is therefore unique, thus $X$ is $m$-subcoskeletal. The condition for $m$-dimensional cells that $x \sim x'$ implies $x= x'$ follows directly from description $(3)$ of $\sim$ from Proposition~\ref{daggericonditions}.
\end{proof} 
Let $k\cat{-inKan}$ denote the full subcategory of $\cat{Set}_{\Delta^k}$ on the inner-Kan $k$-simplicial sets, and $(m,k)\cat{-red}$ denotes the full subcategory of inner-Kan $k$-simplicial sets which are $m$-reduced. The inclusion $k\cat{-inKan} \ra (m,k)\cat{-red}$ has a left adjoint $h_m$ which we will now construct. First we give a truncated version of $h_m:$
\begin{definition} Define $h_m X|^m_0$, to be a truncated $k$-simplicial set in $\cat{Set}_{\Delta^k|^m_0}$ given by: 
\begin{itemize} 
\item $h_m X|^m_0$ is identical to $\tr^m X$ on cells of dimension $<m.$
\item Let $\m$ be a multi-index of dimension $m$. The the $\m$-cells of $h_m X|^m_0$ are the equivalence classes of $\m$-cells of $X$ under $\sim.$ The face maps are well-defined because $x \sim x'$ implies $dx = d x'$.  
\end{itemize}
Note that we have a natural map $\tr^m X \ra  h_m X|^m_0$ which is the identity below dimension $m$ and sends a cell of dimension $m$ to its $\sim$-equivalence class.
\end{definition}
We choose a section $s$ for the natural map $p : \tr^m X \ra h_m X|^m_0$, which can be given by any choice of representatives for equivalence classes of $m$-cells under $\sim$.
\begin{definition}
We will say two maps $Y \ra \tr^m X$ are \emph{$m$-homotopic} if their compositions with $p$ are equal. Composing with the section $s$ gives an equivalence between maps $Y\ra h^\alpha_m X|^m_0$ and equivalence classes of maps $Y \ra \tr^m X.$ For maps $f,f':Y \ra X$ we will say $f$ is \emph{$m$-homotopic} to $f'$ if $\tr^m f$ is $m$-homotopic to $\tr^m f'$ and write $f \sim_m f'.$ Note that clearly $f \sim_m f'$ if and only if $f$ and $f'$ agree on $m-1$ cells and lower and $f(x) \sim_m f'(x)$ for every $m$-cell $x$ in $Y$.
\end{definition}
\begin{definition} \label{homotopycatdef} The $k$-fold simplicial set $h_m X$ is defined by:
$$ (h_m)_{\n} := \Hom(\Delta^k[\n],X)/\sim_{m}.$$ 
\end{definition} 
\begin{proposition} \label{univpropprop} $ (h_m)_{\n}$ has the universal property that there is a natural isomorphism $$\Hom(Y, (h_m)_{\n}) \cong \Hom(Y,X)/\sim_{m}.$$ 
\end{proposition}
\begin{proof} Applying the definition, it's not hard to see that the map $Y \ra (h_m)_{\n}$ can be described as homotopy class of maps $f:\Sk^m Y \ra X$ such that for every cell $y$ in $Y$ there a representative of the class $f|_y$ as a map $\Sk^m y \ra X$ that in extends to a map $\Sk^{m+1} y \ra X.$ On the other hand a map in $\Hom(Y,X) / \sim_{m}$ is a homotopy class of map $\Sk^m Y$ such that some representative of the map $\Sk^m Y \ra X$ extends to $\Sk^{m+1}Y.$

This description gives an injective natural map  $ \Hom(Y,X)/\sim_{m}\ra \Hom(Y, (h_m)_{\n}).$ For this map to fail to be surjective would be for there to be a homotopy class of maps $\Sk^m Y\ra X$ that can be represented cell-by-cell, but not globally, by a map which is extendable to the $m+1$ skeleton. Lemma~\ref{liftlemma} below shows this is impossible.\end{proof}
\begin{lemma}\label{liftlemma}
If $\dim(\n) = m+1$ and $f,f' :\Sk^m \Delta^k[\n]= d\Delta^k[\n] \ra X$ are $m$-homotopic,  then $f$ is commutative (i.e. extends to $\Delta^k[\n]$) if and only if $f'$ is.
\end{lemma}
\begin{proof} The relation of $m$-homotopy for $f$ and $f'$ is equivalent to the statement that $f$ and $f'$ agree on $m-1$ cells and lower and $f(x) \sim_m f'(x)$ for each nondegenerate $m$-cell $x$ in $Y : = d\Delta^k[\n] $. Clearly it suffices to consider the case in which $f$ and $f'$ agree except on one face, which we take to be the face $d^\beta_i \Sk^m \Delta^k[\n] $ in the sphere. Since $d^\beta_i f(Y)\sim d^\beta_i f'(Y)$, and $f(Y)$ and $f'(Y)$ agree on their subhorn of type $\Lambda^{\n}(d^\beta_i)$, the statement in this follows immediately from condition $(2)$ for $\sim$ in Proposition~\ref{daggericonditions}.
\end{proof}

\begin{proposition} \label{hprops} Let $X$ be a inner-Kan $k$-simplicial set. Then:
\begin{enumerate}
\item $h_m X$ is an inner-Kan $k$-simplicial set.
\item $h_m X$ is $m$-reduced.
\item If $X$  is $m$-reduced, then the natural map $X\ra h_m X$ is an isomorphism.
\end{enumerate}
\end{proposition}
\begin{proof} First note that taking a representative of the map $\id: h_m X \ra h_m X$ viewed as an equivalence class of maps $\Hom(h_m X, X)$ gives a section $s: h_m X \ra X$ of the natural map $p: X \ra h_m X.$  Then $(1)$ is immediate since any horn in $h_m X$ can be lifted to a horn in $X$. 

For $(2)$, observe that $h_m X$ is $m$-subcoskeletal follows from \ref{univpropprop} since for a simplex $\Delta^k[\n]$ with $\dim(\n)>m$ the $\sim_m$ class of a map $f:\Delta^k[\n]\ra X$ depends only on its restriction to $\Sk^m (\Delta^k[\n]) \ra X,$ which in particular depends only on the restriction of $f$ to $d\Delta^k[\n].$ If $x'\sim x'$ are $\n$-cells in $h_m X$ with $\dim(\n) = m$, then using the definition of $\sim$, let $\fk$ be a inner horn in $h_m X$ having fillers $H, H'$ with $d^\alpha_i H=dx $ and $d^\alpha_i H' = dx'$. By the existence of the section $s$, we see that the fillers $s(H)$ and $s(H')$ of $s(\fk)$  witness the fact that $s(x)\sim s(x').$ The natural map $X\ra h_m X$ by its construction thus send $s(x)$ and $s(x')$ to the same cell, thus $x =x'.$ This shows that $h_m$ is $m$-reduced.

For $(3)$ note that if $X$ is $m$-reduced, then $f \sim_m f' $ for maps $f,f':\Delta^k[\n] \ra X$ if and only if  $\tr^m f= \tr^m f'$. Since $X$ is $m$-subcoskeletal, this relation is trivial, showing that $X \ra h_m X$ is an isomorphism. 
\end{proof}
\begin{corollary}\label{hisadjoint} $h_m$ is left adjoint to the inclusion $k\cat{-inKan} \ra (m,k)\cat{-red}$
\end{corollary}
\begin{proof}The natural map $X \ra h_m X$ that provides the unit of the adjunction. If $X$ is $m$-reduced then by $(2)$ in Proposition~\ref{hprops}, this map is an isomorphism. The inverse of this map provides the counit of the adjunction. The verifications of the zig-zag identities for the unit and counit are straightforward and are left to the reader.
\end{proof}

\begin{definition} For a inner-Kan simplicial set $X$, we define $\Bic(X)=\Bic(h_2 X).$ For an inner-Kan $k$-simplicial set $Y$ we define $\mbox{Vdc}(Y)=\mbox{Vdc}(h_2 Y)$. 
\end{definition}

\section{The nerve of a VDC with thin structure}
\subsection{The functor $\partial_*$ and the prismatic identities}
\begin{definition} \label{diagonaldef}
The \emph{diagonal} functor $\partial : \cat{Set}_{\Delta^2}\ra \cat{Set}_{\Delta}$ is given by $\partial(X)_{n}=X_{nn}$. This functor has an adjoint on both the left and the right, though we will only consider the right adjoint, $\partial_*$, which is given by $$\partial_*( X)_{mn} = \Hom(\Delta[m] \times \Delta[n], X) $$ with face maps given by $$d_i f = f \circ (\hat{d}^i \times \id) \ \ \ \ \ \ \ \ \ \ \ \ \delta_i f = f \circ (\id \times \hat{d}^i)$$  and degeneracies given by  $$s_i f = f \circ (\hat{s}^i \times \id) \ \ \ \ \ \ \ \ \ \ \ \ \varsigma_i f = f \circ (\id \times \hat{s}^i).$$
\end{definition}
Note that the zeroth row and column of $\partial_* (X)$ can be identified with $X$. For this section, we will also need to understand the \emph{prisms} of $X$, i.e. the elements of $\partial_*(X)_{1n}$ and $\partial_*(X)_{m1}.$   A map $f: \Delta[m]\times \Delta[1]\ra X$ in $\partial_*(X)_{m1}$ is a $2$-homotopy as defined in Definition~\ref{homotopydef}. Recall that the description of $\Delta[m]\times\Delta[1]$ as a colimit of $\Delta[m+1]$-simplicies allows us to describe a homotopy $f$ by giving cells $h_0 f, h_1 f,\ldots h_m f $ subject to the condition that $d_i h_{i-1}f = d_i h_i f$ for $1\leq i \leq m.$ In this way, we can view $h_i$ as a map $\partial_*(X)_{m1} \ra X_{m+1}.$ We write the cells making up a homotopy $f$ using the notation $$hf:=[h_0f, h_1f,\ldots, h_m f].$$

View $[n]$ as the totally ordered set $0<1< \ldots <n.$ Then $N([n])$ is canonically identified with $\Delta[n].$ Since $N$ preserves limits as a functor from the category of partially ordered sets to $\cat{Set}_{\Delta}$, we have $$\Delta[n]\times \Delta[1]\cong N([n])\times N([1])\cong N([n]\times [1]).$$
Since $N$ is full and faithful, a map $\Delta[m]\ra \Delta[n]\times \Delta[1]$ is the same as a map $[m] \ra [n]\times[1]$. The map $h_i$ is dual to the map $\hat{h}_i :[m]\ra [n]\times [1]$ given by the sequence $$\hat{h}_i =((0,0),(1,0),\ldots,(i,0), (i,1),\ldots (n,1)).$$ 

The face and degeneracy maps in the first direction $d_i$ and $s_i$ take a prism $f$ in $\partial_*(X)_{m1}$ to another prism, and make $\partial_*(X)_{\bullet 1}$ itself a simplicial set. These maps are dual to the obvious maps:
\begin{align*} \hat{d}_i:[n]\times [1] \ra [n+1]\times [1] \\ \hat{s}_i:[n]\times [1] \ra [n-1]\times [1]. \end{align*}

With this description, it's easy to check the \emph{prismatic identities}, which describe how the decomposition of $f$ into its nontrivial cells relates to the decomposition of the prisms $d_i f$ and $s_i f:$

\begin{align}h_i d_j &=   \begin{dcases*}
 			     d_j h_{i+1} & if $i \geq j$ \\
			     d_{j+1} h_i  & if $i <j$
   \end{dcases*} \\
h_i s_j &=   \begin{dcases*}
 			     s_j h_{i-1} & if $i > j$ \\
			     s_{j+1} h_i  & if $i \leq j$
   \end{dcases*} \end{align}

\subsection{The bisimplicial set $h_2 \partial_* N(\T)$}

Let $X$ be a $2$-reduced inner-Kan simplicial set. Then for  $\partial_* X $, $\sim$ is the equality relation for $(0,2)$ and $(2,0)$ cells, since $(\partial_* X)_{\bullet 2}$ and $(\partial_* X)_{2\bullet }$ are isomorphic to $X$. 

Let $\T$ be a $(2,1)$-category. We consider the relation  $\sim$ for $(1,1)$ cells of $\partial_* N(\T)$. Using the presentation of $\Delta[1]\times\Delta[1]$ as a colimit of two $\Delta[2]$ cells joined along their $d_1$ face, together with the definition of $N(\T)_2,$ an element  of $(\partial_* N(\T))_{11}$ is given by a diagram in $\T$ of the form shown in Figure~\ref{xelementpicture}:
\begin{figure}[H]
\begin{center}
\begin{tikzpicture}[scale=1.8,auto]

\begin{scope}

\node (10) at (1,1) {$b'$};
\node (00) at (0,1) {$a$};
\node (11) at (1,0) {$c$};
\node (01) at (0,0) {$b$};
\node[rotate=-135] at (.33,.33){$\Rightarrow$};
\node at (.44,.22) {$\eta$};
\node[rotate=45] at (.68,.68){$\Rightarrow$};
\node at (.79,.57) {$\eta'$};
\path[->] (00) edge node[midway]{$f'$}(10);
\path[->] (00) edge node[midway,swap]{$f$}(01);
\path[->] (01) edge node[midway,swap]{$h$}(11);
\path[->] (10) edge node[midway]{$h'$}(11);
\path[->] (00) edge node[midway, shift={(-.32,.15)}]{$g$}(11);
\end{scope}

\end{tikzpicture}
\caption{\label{xelementpicture} An element $[\eta', \eta]$ of $(\partial_* N(\T))_{11}$}    \end{center}
\end{figure}
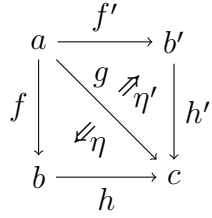 We denote the element of $(\partial_* N(\T))_{11}$ shown above in Figure~\ref{xelementpicture} by $[\eta', \eta]$. For such an element, we will call the composition $ \mathrm{comp}([\eta', \eta]):=\eta' \bullet \eta^{-1}$ the \emph{interior composition}. 

\begin{proposition} \label{sim1prop} $x \sim x'$ if and only if $dx = d x' $ and $\mathrm{comp}(x)=\mathrm{comp}(x')$. 
\end{proposition}
\begin{proof} We will work with condition $(1)$ for $\sim$ in Proposition~\ref{daggericonditions}, with $\alpha=1$

 Take $x = [\eta', \eta]$, and suppose $x \sim x'$. We seek to explicitly describe the horn  $\squish^1 x$. First we will describe $s_0 x $ using the prismatic identities. $$h s_0 x = [s_1 h_0 x , \ s_0 h_0 x,\ s_0 h_1 x].$$   We write the Glenn table for each cell $h_i s_0 x$ below. 
\begin{figure}[H]
	\begin{subfigure}[b]{0.49\textwidth} \centering\caption*{$h_0 s_0 x = s_1 h_0 x$}
		\begin{tabular}{ r | l || l | l | l | }
     \cline{2-5}
	    	      &  $s_0h'=\rho_{h'}$   &         $h'$      &  $h'$             &    $\id_{b'}$                \\  \cline{2-5} 
	 	$(1)$    &   $\eta'$     &  $ h'$     &        $g$        &  $f'$                       \\ \cline{2-5}
  		 $h_0 d_1 s_0 x$ &         $\eta'$    &   $h'$  & $g$               &     $f'$     \\   \cline{2-5}
    		  $h_0 d_2 s_0 x$    &     $s_1f' = \lambda_{f'}$     &    $\id_{b'}$         &                  $f'$       &  $f'$                       \\   \cline{2-5}
      \cline{2-5}
     \end{tabular} 
      \end{subfigure}
		 \begin{subfigure}[b]{0.49\textwidth}  \centering \caption*{$h_1 s_0 x = s_0 h_0 x$}
			
					\begin{tabular}{ r | l || l | l | l | }
			     \cline{2-5}
				    $h_0 d_0 s_0x$      &  $\eta'$   &         $h'$      &  $g$             &    $f'$                \\  \cline{2-5} 
						 $(1)$              &   $\eta'$             &  $ h'$     &        $g$        &  $f'$                       \\ \cline{2-5}
						 $(2)$               &   $s_0g = \rho_g$    &         $g$  & $g$               &     $\id_a$     \\   \cline{2-5}
 						  $h_1 d_2 s_0x$    &     $s_0f' = \rho_{f'}$     &    $f'$         &                  $f'$       &  $\id_a$                       \\   \cline{2-5}
    			 \cline{2-5}
    			 \end{tabular}\end{subfigure}
    \end{figure}
\begin{table}[H]\caption*{$ h_2 s_0 x = s_0 h_1 x $}\begin{center}
\begin{tabular}{ r | l || l | l | l | }
     \cline{2-5}
  $h_1 d_0 s_0 x$        &  $\eta$   &         $h$      &  $g$             &    $f$                \\  \cline{2-5} 
  $h_1 d_1 s_0 x$   &   $\eta$             &  $ h$     &        $g$        &  $f$                       \\ \cline{2-5}
  $(2)$ &         $s_0 (g) = \rho_g$    &   $g$  & $g$               &     $\id_a$     \\   \cline{2-5}
       &     $s_0f= \rho_f$     &    $f$         &                  $f$       &  $\id_a$                       \\   \cline{2-5}
    \cline{2-5}
    \end{tabular} \end{center}
    \end{table}    

These three Glenn tables show the image of every non-singular $1$ and $2$ cell in $\Sk^2 s_0 x: \Sk^2( \Delta[2]\times \Delta[1])\ra X,$ based on our description of $\Sk^2( \Delta[2]\times \Delta[1])$ as a colimit. To $\squish^1(x)$, by definition the $\Lambda[2^{\, 1}, 1]$-horn inside $s_0 x$, we remove the interior faces together with the faces making up $d_1 s_0 x$, and the $1$ cell in the interior of this face, which by simplicial and prismatic identities, is 
 $$d_1 h_0 d_1 s_0 x = d_1 h_1 d_1 s_0 x= d_1 d_1 h_0 s_0 x = d_1 d_1 h_1 s_0 x=d_1 d_2 h_1 s_0 x = d_1 d_2 h_2 s_0 x.$$ The following Glenn tables show the parts of the $s_0x$ which are also in $\squish^1 (x).$

\begin{figure}[H]
	\begin{subfigure}[b]{0.49\textwidth} \centering
\caption*{$ h_0 s_0 x = s_1 h_0 x $}
\begin{tabular}{ r | l || l | l | l | }
     \cline{2-5}
                          &  $s_0h' = \rho_{h'}$   &         $h'$                &  $h'$             &    $\id_{b'}$                \\  \cline{2-5} 
  $(1)$                    &         &         $ h'$                  &          &  $f'$                       \\ \cline{2-5}
        $h_0 d_1 s_0 x$           &                  &   $h'$  &          &     $f'$     \\   \cline{2-5}
      $h_0 d_2 s_0 x$    &     $s_1f' =\lambda_{f'}$     &    $\id_{b'}$         &                  $f'$       &  $f'$                       \\   \cline{2-5}
    \cline{2-5}
    \end{tabular} 
    \end{subfigure}
	\begin{subfigure}[b]{0.49\textwidth} \centering
	\caption*{$ h_1 s_0 x = s_0 h_0 x $}
\begin{tabular}{ r | l || l | l | l | }
     \cline{2-5}
    $h_0 d_0 s_0x$      &  $\eta'$   &         $h'$      &  $g$             &    $f'$                \\  \cline{2-5} 
 $(1)$    &              &  $ h'$     &               &  $f'$                       \\ \cline{2-5}
 $(2)$   &               &   $g$  &                      &     $\id_a$     \\   \cline{2-5}
   $h_1 d_2 s_0x$        &     $s_0 f'=\rho_{f'}$     &    $f'$         &                  $f'$       &  $\id_a$                       \\   \cline{2-5}
    \cline{2-5}
    \end{tabular} 
    \end{subfigure}
   	\end{figure}
\begin{table}[H]\caption*{$ h_2 s_0 x = s_0 h_1 x $}\begin{center}
\begin{tabular}{ r | l || l | l | l | }
     \cline{2-5}
  $h_1 d_0 s_0 x$        &  $\eta$   &         $h$      &  $g$             &    $f$                \\  \cline{2-5} 
  $h_1 d_1 s_0 x$   &                &  $ h$     &        &  $f$                       \\ \cline{2-5}
  $(2)$ &           &   $g$  &              &     $\id_a$     \\   \cline{2-5}
       &     $s_0 f = \rho_f$     &    $f$         &                  $f$       &  $\id_a$                       \\   \cline{2-5}
    \cline{2-5}
    \end{tabular} \end{center}
    \end{table}

Suppose we have $x'$ with $x\sim x'$. Then by Proposition~\ref{hprops} the map $\squish^1(x) \ra \partial_* X$ has an extension $H$ to $\Delta[2] \ra \Delta[1]$ with $d_1 H = x'$. Then from the above tables, we see that the interior $d_2 h_2 H = d_2 h_1 H$ is a morphism $g' \Rightarrow g \circ \id_a$ in $\T$. By the invertibility of the unitors in $\T$, we can write this $2$-morphism as $\rho_{g} \bullet \beta$ where $\beta: g' \Rightarrow g$. We will see in a series of steps using the fact that $N(\T)$ has unique fillers for inner $3$-horns that the choice of $\beta$ uniquely determines $H$. 

Then $d (h_2 H) = [\eta, Y, \rho_g \bullet \beta, \Id_f]$ where $Y$ is necessarily the unique $2$-morphism in $\T$  making the $3$-cell condition true: $$\alpha_{ h,f,\id_a}\bullet(h \rhd \rho_f )\bullet Y =(\eta \lhd \id_a )\bullet \rho_{g} \bullet \beta.$$ Applying the compatibility of $\alpha$ and $\rho$ (\textbf{B13}) on the left and the naturality of $\rho$ (\textbf{B7}) on the right this is equivalent to: $$\rho_{h \circ f} \bullet Y =\rho_{h \circ f} \bullet \eta \bullet \beta.$$ Thus $Y= \eta \bullet \beta= d_1 h_2 H $ A similar argument shows $\eta' \bullet \beta = d_1 h_1 H= d_1 h_0 H.$ Then $$d h_0 H = [ s_0h', \eta' \bullet \beta, - , s_1 f']$$ thus $h_0 H = s_1 (\eta' \bullet \beta)$ by uniquess of fillers for inner $3$-horns. We have now described every cell in $\Sk^2 H$, and we give the ``filled'' Glenn tables for $h_i H$ below:
\begin{figure}[H]	\begin{subfigure}[b]{0.49\textwidth} \centering
\caption*{$ h_0 H$}
\begin{tabular}{ r | l || l | l | l | }
     \cline{2-5}
          &  $s_0h'=\rho_{h'}$   &         $h'$      &  $h'$             &    $\id_{b'}$                \\  \cline{2-5} 
 $(1)$    &   $\eta'\bullet \beta$     &  $ h'$     &        $g'$        &  $f'$                       \\ \cline{2-5}
   $h_0 d_1 H$ &    $\eta'\bullet \beta$    &   $h'$  & $g'$               &     $f'$     \\   \cline{2-5}
      $h_0 d_2 H$    &     $s_1f' = \lambda_{f'}$     &    $\id_{b'}$         &                  $f'$       &  $f'$                       \\   \cline{2-5}
    \cline{2-5}
    \end{tabular} 
    \end{subfigure}
\begin{subfigure}[b]{0.49\textwidth} \centering \caption*{$ h_1 H$}
\begin{tabular}{ r | l || l | l | l | }
     \cline{2-5}
    $h_0 d_0 H$      &  $\eta'$   &         $h'$      &  $g$             &    $f'$                \\  \cline{2-5} 
 $(1)$              &   $\eta'\bullet \beta$                &  $ h'$     &        $g'$        &  $f'$                       \\ \cline{2-5}
 $(2)$               &   $\rho_g \bullet \beta$    &         $g$  & $g'$               &     $\id_a$     \\   \cline{2-5}
   $h_1 d_2 H$    &     $s_0f' = \rho_{f'}$     &    $f'$         &                  $f'$       &  $\id_a$                       \\   \cline{2-5}
    \cline{2-5}
    \end{tabular}
    \end{subfigure}\end{figure}
\begin{table}[H]\caption*{$ h_2 H$}\begin{center}
\begin{tabular}{ r | l || l | l | l | }
     \cline{2-5}
  $h_1 d_0 H$        &  $\eta$   &         $h$      &  $g$             &    $f$                \\  \cline{2-5} 
  $h_1 d_1 H$   &   $\eta \bullet \beta$             &  $ h$     &        $g'$        &  $f$                       \\ \cline{2-5}
  $(2)$ &         $ \rho_g \bullet \beta$    &   $g$  & $g'$               &     $\id_a$     \\   \cline{2-5}
       &     $s_0f= \rho_f$     &    $f$         &                  $f$       &  $\id_a$                       \\   \cline{2-5}
    \cline{2-5}
    \end{tabular} \end{center}
    \end{table}   
We deduce that $h_0 x' = h_0 d_1H =\eta' \bullet \beta$ and $h_0 x' = h_0 d_1 H$, so $x'=[\eta'\bullet \beta, \eta \bullet \beta]$ Thus we have $$\mathrm{comp}(x')=(\eta' \bullet \beta) \bullet (\eta\bullet \beta)^{-1} = \eta' \bullet \eta^{-1}=\mathrm{comp}(x).$$  The fact that  $dx = d x'$ is also apparent.

Conversely suppose $x'=[\theta',\ \theta]$ is a square in $\partial_* (N(\T))_{11}$ with $dx=dx'$ and $$\mathrm{comp}(x)=\eta' \bullet \eta^{-1} = \theta' \bullet \theta^{-1} =\mathrm{comp}(x').$$   Since we showed that there is a unique filler $H$ of $\squish^1(x)$ for every $\beta$, we can take $\beta = (\eta')^{-1} \bullet \theta'.$ Then if $H$ is corresponding filler of $\squish^1(x)$, we have 
\begin{align*}d_1 H&=[\eta' \bullet (\eta')^{-1} \bullet \theta',\  \eta \bullet  (\eta')^{-1} \bullet \theta'   ] \\ 
									 &=[\theta', \ \theta] = x'. \end{align*}
Thus $x \sim x'$.
\end{proof}
Observe that if $(T,\chi)$ is an algebraic $2$-reduced inner-Kan simplicial set, then $h_2 \partial_* T$ is $2$-reduced inner-Kan by Proposition~\ref{hprops}. Furthermore, since the only inner horns of dimension $2$ in $h_2 \partial_* T$ are of the form $\Lambda[2^{\, 1}, 0]$ and $\Lambda[0,2^{\, 1}]$, and the zeroth row and column of $h_2 \partial_* T$ are identified with $T$, $\chi$ induces a natural algebraic structure on $h_2 \partial_* T$, which we will also denote by $\chi.$

\begin{theorem} We have two functors from small $(2,1)$-categories to algebraic $2$-reduced inner-Kan bisimplicial sets, given  by $\T  \ra h_2 \partial_* N(\T)$ and $\T \ra N(\ES(\llcorner\T \lrcorner  )).$ There is a strict natural isomorphism $V$ between these functors. 
\end{theorem}
\begin{proof} Since both  $h_2 \partial_* N(\T)$ and $N(\ES(\llcorner\T \lrcorner ))$ are $3$-coskeletal, to give the isomorphism $h_2 \partial_* N(\T) \ra N(\ES(\llcorner\T \lrcorner ))$ it suffices to give an isomorphism $$V|^3_0 :h_2 \partial_* N(\T)|^3_0 \ra N(\ES(\llcorner\T \lrcorner ))|^3_0.$$ The zeroth row and column of each bisimplicial set can be identified with $N(\T),$ so we give the map for $(1,1),$ $(2,1)$, and $(1,2)$ cells. Note that the fact that $V$ is strict is immediate from this part of the definition.

An element $(h_2 \partial_* N(\T))_{11}$ is given by a $\sim$ equivalence class of squares in $\partial_* N(\T)$, which by Proposition~\ref{sim1prop} is given by $1$-morphisms $a \stackrel{f}{\ra} b \stackrel{g}{\ra}c $ and $a \stackrel{f'}{\ra} b' \stackrel{g'}{\ra}c $ together with a morphism $\Theta:g\circ f \Rightarrow g' \circ f'$ which is the interior composition of the cells in the equivalence class. The same data defines a square of $\ES(\llcorner\T \lrcorner )$, which is an element of $N(\ES(\llcorner\T \lrcorner ))_{11}.$
This correspondence defines $V|^3_0$ for $(1,1)$ cells.
 
Since $h_2 \partial_* N(\T)$ and $N(\ES(\llcorner\T \lrcorner ))$ are $3$-subcoskeletal, in order to check that $V|^3_0$ is well-defined and an isomorphism for $(2,1)$ cells, we must check that a $d\Delta^2[2,1]$-sphere $s$ in $h_2 \partial_* N(\T)$ is commutative if and only if $V(s)$ is commutative in $N(\ES(\llcorner\T \lrcorner)).$ Since $V|^3_0$ is an isomorphism on $2$-cells and lower, we can take $V|^3_0(s)$ to be the following arbitrary $d\Delta^2[2,1]$-sphere in $N(\ES(\llcorner\T \lrcorner)):$
$$[ \Sigma,\ \Pi,\ \Theta \ \ | \ \ \eta' ,\ \eta ]$$ with Glenn table:

\begin{table}[H]\begin{center}
    \subfloat
    {\begin{tabular}{ r | l || l | l !{\vrule width 2pt} l | l  |}
    \cline{2-6}

        &  $\Sigma:h \circ g \Rightarrow h' \circ g''$  &  $h$  &   $g''$   &     $h' $    &    $g$       \\ \cline{2-6} 
     
        &  $\Pi:h \circ i \Rightarrow i' \circ f'$  &  $h$  &   $f'$   &     $i' $    &    $i$     \\ \cline{2-6} 
 
        &  $\Theta:g''\circ f \Rightarrow g' \circ f'$  &  $g''$  &   $f'$   &     $g' $    &    $f$      \\ \cline{2-6}  \end{tabular}}
\subfloat
{
\begin{tabular}{ r | l || l | l | l !{\vrule width 2pt}  }
    \cline{2-5}

    &  $\eta': i \Rightarrow h'\circ g' $     &     $h' $    &    $i'$   &  $g' $        \\ \cline{2-5} 
     
 &  $\eta: i \Rightarrow g \circ f$    &     $g $    &    $i$   &  $f $      \\ \cline{2-5} 
 
  \end{tabular}
}
 \end{center}   \end{table}

This sphere is commutative in $N(\ES(\llcorner\T \lrcorner ))$ if and only if it meets the commutativity condition: 
$$(\Sigma \boxminus \Theta) \lact \eta \ract \eta' =\Pi.$$ 
Applying the definition of $\boxminus$, $\lact$, and $\ract$ in $\ES(\llcorner\T \lrcorner )$ from Section~\ref{edgesymsec} we see that this condition is equivalent to the following identity in $\T$:
\begin{equation}((\eta')^{-1}\lhd f')\bullet \alpha_{h',g',f'} \bullet (h' \rhd \Theta)\bullet \alpha^{-1}_{h',g'',f}\bullet (\Sigma\lhd f) \bullet \alpha_{h,g,f} \bullet (h \rhd \eta)=\Pi. \label{NEScomcond} \end{equation}
To check if $s$ is commutative in $h_2 \partial_* N(\T)$, by Lemma~\ref{liftlemma} we can check if any lift $\hat{s}$ to $\partial_* N(\T)$ found by taking a representative of each $2$ cell making up $s$. We take $$\hat{s}= [[\Id_{h\circ g},\Sigma],\ [\Id_{h\circ i}, \Pi],\ [\Id_{g'' \circ f}, \Theta] \ \ | \ \ \eta',\ \eta]$$

Using the prismatic identities, we can calculate that if $S$ is a filler of $\hat{s}$, then 
\begin{align*}  
d h_0 S &= [\eta',\ X,\ \Pi,\ \Theta]\\
d h_1 S &= [\Sigma,\ X ,\ Y ,\ \Id_{g''\circ f}] \\
d h_2 S &= [\Id_{h\circ g},\ \Id_{h\circ i},\ Y,\ \eta]
\end{align*}
We calculate $Y$ using the commutativity condition for $3$-cells in $N(\T)$, applied to $h_2 S$:
\begin{align*}
\alpha_{h,g,f} \bullet (h \rhd \eta) \bullet \Id_{h\circ i } &= ( \Id_{h\circ g} \lhd f) \bullet Y \\
\alpha_{h,g,f}\bullet (h\rhd \eta) &= Y
\end{align*}
Next we calculate $X$ using the commutativity condition for $h_1 S:$
 \begin{align*}
\alpha_{h',g'',f} \bullet (h' \rhd \Id_{g''\circ f}) \bullet X &= ( \Sigma \lhd f) \bullet \alpha_{h,g,f} \bullet (h \rhd \eta) \\
X &= \alpha^{-1}_{h',g'',f} \bullet ( \Sigma \lhd f) \bullet \alpha_{h,g,f} \bullet (h \rhd \eta)
\end{align*}
Finally we see that the filler $S$ exists and thus $\hat{s}$ and therefore $s$ is commutative if and only if $h_0 S=[\eta',\ X,\ \Pi,\ \Theta]$ is commutative in $N(\T)$, which holds if and only if 
 \begin{align*}
\alpha_{h',g',f'} \bullet (h' \rhd \Theta) \bullet \alpha^{-1}_{h',g'',f} \bullet (\Sigma \lhd f) \bullet \alpha_{h,g,f} \bullet (h \rhd \eta) &= ( \eta' \lhd f') \bullet \Pi \\
( (\eta')^{-1} \lhd f')\bullet \alpha_{h',g',f'} \bullet (h' \rhd \Theta) \bullet \alpha^{-1}_{h',g'',f} \bullet (\Sigma \lhd f) \bullet \alpha_{h,g,f} \bullet (h \rhd \eta) &=  \Pi
\end{align*}
This matches the condition of Equation~\ref{NEScomcond}, showing $s$ is commutative if and only if $V|^3_0(s)$ is commutative, which ensures $V|^3_0$ is well defined and bijective for $(2,1)$-cells.

 By symmetry, the same holds for $(1,2)$ cells, showing $$V|^3_0:h_2 \partial_* N(\T)|^3_0 \ra N(\ES(\llcorner\T \lrcorner ))|^3_0$$ is an isomorphism. Thus $$V:=\cosk^3 V|^3_0:h_2 \partial_* N(\T) \ra N(\ES(\llcorner\T \lrcorner ))$$ is an isomorphism. To show the naturality of $V$, consider a functor $F:\T \ra \mathcal{S}$ of $(2,1)$-categories. The only nontrivial condition we must check is for $(1,1)$-cells, that for $x=[\eta',\ \eta]$  we have $V (h_2 \partial_* N(F)(x))= N(\ES(\llcorner F\lrcorner))(V(x))$. This works out to the statement that $$[N(F)(\eta'),\ N(F)(\eta)] = [\rho_g \bullet \eta' ,\ \rho_g \bullet \eta]$$ and $[\eta',\ \eta]$ have the same interior composition, i.e. $$(\rho_g \bullet \eta)^{-1} \bullet (\rho_g \bullet \eta') = \eta^{-1} \bullet \eta.$$
 \end{proof}
The isomorphism $V: h_2\partial_*(N(\T)) \ra N(\ES(\llcorner\T \lrcorner ))$ together with the natural isomorphisms $u:N\Bic \ra \Id$ and $U:\Bic N \ra \Id$ from Chapter~\ref{bicchapter} gives us an natural isomorphism:

\begin{multline}W:\Vdc(h_2 \partial_* (T,\chi)) \stackrel{\cong}{\ra} \Vdc(h_2\partial_* N(\Bic(T))) \\ \stackrel{\cong}{\ra}  \Vdc(N(\ES(\llcorner\T \lrcorner ))) \stackrel{\cong}{\ra}\ES(\llcorner\Bic(T,\chi)\lrcorner). \end{multline}
$W$ can be seen to be strict using the fact that $V$ is strict by applying Theorem~\ref{vdcsummary}. 
\subsection{Thin structures for bisimplicial sets}

\begin{proposition} If $X$ is inner-Kan then $\partial_* X$ is inner-Kan. \label{partialkanprop}
\end{proposition}
\begin{proof} 
In \cite{Lur09}, an \emph{inner fibration} in $\cat{Set}_\Delta$ is a map which has the right lifting property with respect to inner horn inclusions, and a map is \emph{inner anodyne} if it has the left lifting property with respect to inner fibrations. In particular, if $X$ is a inner-Kan simplicial set, the map $X \ra \Delta[0]$ is an inner fibration, and if $f:A\ra B$ is anodyne, any map $A \ra X$ extends along $f$ to a map $B\ra X$. Equivalently, the map $f^\star: [B,X]\ra [A,X]$ has a section. Corollary~$2.3.2.4$ of \cite{Lur09}, which is due to Joyal (\cite{Joy08}), asserts that if $i:A\ra A'$ is inner-anondyne and $j:B\ra B'$ is a cofibration (i.e. an injective map), then the induced map $$ \left(A \times B'\coprod_{A\times B} A'\times B \right)\ra A'\times B'$$ is inner-anodyne. Applying this fact where $i:\Lambda^n_i \ra \Delta[n]$ is an inner-horn inclusion and $j:d\Delta[m]\ra \Delta[m]$ is the natural map, we get that the inclusion $$k:\left(\Lambda^n_i \times \Delta[m]\coprod_{\Lambda^n_i\times d\Delta[m] } \Delta[n]\times d\Delta[m] \right)\ra \Delta[n]\times\Delta[m] $$ is inner anodyne. It's not hard to see that $k$ is isomorphic to $\partial{\hat{k}}$ where $\hat{k}$ is the horn inclusion $\Lambda[n^{\ i},m]\ra \Delta[m,n]$. Thus if $X$ is an inner-Kan simplicial set, the map $$(\partial \hat{k})^\star : [\partial\Lambda[n^{\ i},m],X] \ra [\partial \Delta[n,m], X]$$ has a section. Applying the adjunction of $\partial$ and $\partial_*$ we conclude that the map $$\hat{k}^\star:[\Lambda[n^{\ i},m],\partial_* X] \ra [\Delta[n,m], \partial_*X]$$ has a section. By a similar argument, the natural map $$[\Lambda[n,m^{\ i}],\partial_* X] \ra [\Delta[n,m], \partial_*X]$$ with $0<i<m$ has a section. We conclude that $\partial_* X$ has fillers for inner horns, as was to be shown.
 \end{proof}
\begin{definition} \label{algbithin} A \emph{algebraic $2$-reduced inner-Kan bisimplicial set with thin structure} or ADCT consists of:
\begin{enumerate}
\item An algebraic $2$-reduced inner-Kan bisimplicial set $(X,\chi)$
\item An algebraic $2$-reduced inner-Kan simplicial set $(T,\tau)$
\item A strict morphism of algebraic $2$-reduced inner-Kan bisimplicial sets $t:h_2 \partial_* (T,\tau) \ra X$ which is an isomorphism when restricted to zeroth row and column. 
\end{enumerate} 
A morphism  $F:(X,T,t) \ra (Y,S,s)$ is given by the obvious commutative square. Such a morphism is \emph{strict} if the maps $X \ra Y$ and $T \ra S$ are strict.
\end{definition}
Note that by the universal property of $h_2$, since $X$ is $2$-reduced inner-Kan, a map $t:h_2\partial_* T \ra X$ is equivalent to a map $t': \partial_* T \ra X.$ 

We will show in this section that ADCT's correspond to VDC's with thin structure. First let $(X,T,t)$ be a ADCT . We have $$\Vdc(t):\Vdc(h_2 \partial_* (T,\tau)) \ra \Vdc(X),$$ which is strict by Theorem~\ref{vdcsummary} and the fact that $t$ is strict. Precomposing with the strict isomorphism $W^{-1}$ constructed above we get a map $$\ES(\llcorner\Bic(T,\tau)\lrcorner)\ra \Vdc(X).$$ Clearly this functor is strict and induces an isomorphism on the horizontal and vertical $(2,1)$-categories so this gives $(\Vdc(X),\Bic(T))$ the structure of a VDC with thin structure. It is easy to see (using the naturality of $W$) that this construction sends morphisms of ADCT's to functors which preserve thin structure. Thus this construction is a functor from the category of algebraic $2$-reduced inner-Kan bisimplicial set with thin structure to the category of small VDC's with thin structure. 

\begin{definition}  The construction given above is a functor from the category of algebraic $2$-reduced inner-Kan bisimplicial set with thin structure to the category of small VDC's with thin structure. We denote this functor by $\Vdcth$.
\end{definition}

In the opposite direction, if $(\D,\T,t)$ is a VDC with thin structure, then we have a strict functor $T: \ES(\llcorner\T\lrcorner)\ra \D$. Then applying $N$ we have a strict morphism $$N(T):N(\ES(\llcorner\T \lrcorner ))\ra N(\D).$$ Then precomposing with $V$ as constructed above yields a strict map $$ h_2\partial_*(N(\T))\ra N(\D)$$ which is a thin structure $N(\D)$ making it into a ADCT . Using the naturality of $V$, it is immediate that this construction sends strictly identity-preserving functors which preserve thin structure to morphisms of ADCT's. Thus this construction given above is a functor $N_t$ from  the category of small VDC's with thin structure and their strictly identity-preserving functors to the category of algebraic $2$-reduced inner-Kan bisimplicial set with thin structure.

\begin{theorem}\label{NVdcthm} As described in this section, the functors $N_t$ and $\Vdcth$ give an equivalence of categories between from the category of small VDC's with thin structure and strictly identity preserving functors to the category of ADCT's. Furthermore $N_t$ and $\Vdcth$ preserve strictness, and the isomorphisms $u$ and $U$ which give the equivalence are strict.
\end{theorem}

\begin{proof}  The fact that $N_t$ and $\Vdcth$ preserve strictness follow from the same statements for $N$ and $\Vdcth$ as given in Theorem~\ref{vdcsummary}. The construction of the strict isomorphisms $u$ and $U$ from the isomorphisms of the same names in Chapter~\ref{vdcchapter} is straightforward and left to the reader.
\end{proof}
The following corollary follows from Theorem~\ref{NVdcthm} and Conjecture~\ref{foldesconj}, so is itself (mildly) conjectural:
\begin{corollary}\label{foldingcor} $N_t\circ \ES$ and $\fold \circ \Vdcth$  give an equivalence of categories between from the category of small fancy bicategories and strictly identity preserving functors to the category of ADCT's. Furthermore these functors preserve strictness, and the isomorphisms  which give the equivalence are strict.
\end{corollary}

\section{The vertically trivial nerve for pseudo-double-categories and bicategories} 
\begin{definition}A pseudo-double category consists of a category $\D_0$ called the \emph{category of objects and vertical morphisms} and a category $\D_1$ called the \emph{category of horizontal morphisms and squares}, together with
\begin{itemize}
 \item \emph{left frame} and \emph{right frame} functors $L,R:\D_1\ra\D_0.$ 
 \item an \emph{identity} functor $U:\D_0\ra \D_1 $
 \item and a \emph{horizontal composition functor} $\odot:\D_1 \underset{\D_0}{\times} \D_1 \ra \D_1$
 \end{itemize}
  Also, we have additional structural data in the form of natural isomorphisms:\TabPositions{3.4cm}
\begin{itemize}
\item An \emph{associator}: \tab $\mathfrak{a}: \ \ M \odot (N \odot P) \ra (M \odot N )  \odot P  $
\item A \emph{right unitor}: \tab $\mathfrak{r}: \ \ M \ra M\odot U_A  $
\item A \emph{left unitor}: \tab $\mathfrak{l}: \ \ M \ra U_A\odot M $
\end{itemize}
We insist that the components of $\mathfrak{a}$, $\mathfrak{r}$ and $\mathfrak{l}$ \emph{horizontal-globular}, meaning that their left and right frames are identities. This data satisfies axioms asserting the compatibility of the unitors and the associator and the pentagon identity for the associator. 
\end{definition}
We draw squares in a pseudo-double category just as we draw them in a VDC:
\begin{center}
\begin{tikzpicture}[scale=1.8,auto]
\begin{scope}
\node (10) at (1,1) {$b$};
\node (00) at (0,1) {$a$};
\node (11) at (1,0) {$b'$};
\node (01) at (0,0) {$a'$};
\node[rotate=-45] at (.5,.5){$\Rightarrow$};
\node[scale=.8] at (.6,.65){$\Theta$};
\path[->] (00) edge node[midway]{$s(\Theta)$}(10);
\path[->] (00) edge node[midway,swap]{$L(\Theta)$}(01);
\path[->] (01) edge node[midway,swap]{$t(\Theta)$}(11);
\path[->] (10) edge node[midway]{$R(\Theta)$}(11);
\end{scope}
\end{tikzpicture}
\end{center}
\begin{definition} A pseudo-double category is called \emph{vertically trivial} if it has only identity vertical morphisms.\end{definition}
\begin{definition} Let $\D$ be a pseudo-double category. The \emph{horizontal bicategory} of $\D$ is the bicategory of objects, horizontal $1$-morphisms, and horizontal-globular $2$-morphisms. The associator and unitors of $\D$ can be used to give the associator and unitors of $\D$, since the components are assumed to be globular. In the other direction, we can construct a vertically trivial pseudo-double category $\D$ from a bicategory $\B$, taking the objects, horizontal morphisms , and squares of  $\D$ to be the $0$-morphisms, $1$-morphisms, and $2$-morphisms of $\B$, respectively. These constructions give an equivalence between vertically trivial pseudo-double categories and bicategories.
\end{definition}
\begin{remark} Many popular bicategories are horizontal bicategories of a pseudo-category in a natural way. For instance, the category $\mathbf{\mathcal{M}od}$ of rings, bimodules, and bimodule morphisms can be enlarged to a double category that includes ring homomorphisms as the vertical $1$-morphisms.
\end{remark}

\begin{definition} \label{pseudothin} A  \emph{fancy pseudo-double category} is a pseudo-double  category $\widetilde{D}$ together a vertically trivial pseudo-double category $\overline{\D}$ such that $\overline{\D}_1$  is a groupoid and a strict functor  $t_\D:\overline{\D} \rightarrow \widetilde{\D}$ which is an isomorphisms for objects and horizontal $1$-morphisms. 
\end{definition}
\begin{definition} A VDC will be called \emph{vertically $1$-trivial} if its  vertical double category $\mathcal{V}$ has only identity $2$-morphisms.
\end{definition}
If $\D$ is vertically $1$-trivial, the actions $\lact$ and $\ract$ of vertical $2$-morphisms on squares are trivial, and the vertical composition of $1$-morphisms $\circ_v$ and of squares $\boxminus$ are strictly associative and strictly unital. We construct a fancy pseudo-double category $\trunc(\D)$ from a vertically $1$-trivial VDC $\D$:

\begin{itemize} 
\item  $\widetilde{\trunc(\D)}_0$ is the category of objects of $\D$ and vertical morphisms
\item $\widetilde{\trunc(\D)}_1$ is the category of horizontal morphisms and squares, with composition given by $\boxminus$. The identity of this category for a horizontal $1$-morphism $f$ is given by $\ID_f$
\item The functor $L:\widetilde{\trunc(\D)}_1\ra\widetilde{\trunc(\D)}_0$ and $R:\widetilde{\trunc(\D)}_1\ra\widetilde{\trunc(\D)}_0$ take squares to their source and target (left and right) vertical morphisms, and take horizontal $1$-morphisms to their source and target objects, respectively.
\item The functor $U: \widetilde{\trunc(\D)}_0\ra\widetilde{\trunc(\D)}_1$ is given on objects by the horizontal identity $\id$ and on (vertical) morphisms $p$ by $\ID_p$ 
\item The horizontal composition functor $\odot$ is given for horizontal morphisms in $\widetilde{\trunc(\D)}_1$ by the composition $\circ_h$ of $H$, and for squares between horizontal morphisms by $\boxvert$
\item The associator $\mathfrak{a}$ is given for horizontal morphisms $f, g, h$ in $\widetilde{\trunc(\D)}_1$ by letting $\mathfrak{a}_{f,g,h}: h \circ (g \circ f) \ra (h\circ g)\circ f$ be defined to be $\ID_{(h\circ g)\circ f} \uact \alpha_{h,g,f}.$ 
\item The right unitor $\mathfrak{r}$ is given  for horizontal $f:a\ra b$ by letting $\mathfrak{r}_{f}: f \ra f\circ \id_a$ be defined to be $\ID_{f\circ \id_a} \uact \rho_{f}.$ 
\item  The left unitor $\mathfrak{l}$ is given  for horizontal $f:a\ra b$ by letting $\mathfrak{l}_{f}: f \ra  \id_b\circ f$ be defined to be $\ID_{ \id_b\circ f} \uact \lambda_{f}.$ 
\item $\overline{\trunc(\D)}$ is $\mathcal{H}$ viewed as a vertically trivial pseudo-double category. The map $t_\D:\overline{\trunc(\D)} \ra  \widetilde{\trunc(\D)} $  is given for a $2$-morphism  in $\mathcal{H}$  (i.e. squares in $\BBT$) by $t_\D(\beta)=\ID_g\uact \beta.$
\end{itemize}

\begin{proposition}\label{swappyprop}For horizontal $1$-morphisms $f,g$ and horizontal $\beta:f \Rightarrow g$ in a VDC, $\ID_g \uact \beta =\ID_f \dact \beta^{-1}$. \end{proposition}
\begin{proof} By middle associativity (Lemma~\ref{midass}) and the strict unitality of $\boxminus$, we have:

\begin{align*} \ID_g \uact \beta &=   \ID_f  \boxminus (\ID_g \uact \beta) \\ 
& = \ID_f \uact \beta \\ 
&= (\ID_f \dact \Id_f)\boxminus (\ID_g \uact \beta) \\
&= (\ID_f \dact (\Id_f \circ \beta^{-1}))\boxminus \ID_g \\
&= \ID_f \dact \beta^{-1}. 
\end{align*}
\end{proof}
From this, we see that each of our definitions could have equivalently been given using $\dact$ instead of $\uact$, so the apparent asymmetry of the construction is illusory.

We can verify the pseudo-double category axioms for $\widetilde{\trunc(\D)}$ from corresponding VDC axioms for $\D$:
\begin{itemize}
\item The category axioms for $\widetilde{\trunc(\D)}_0$ and $\widetilde{\trunc(\D)}_1$ follow from the strict associativity and unitality of $\circ_v$ and $\boxminus$, since the associators and unitors of $V$ are necessarily identities. 
\item The functoriality $L$ and $R$ follow from the values of the specified left and right vertical morphisms for squares composed by $\boxvert$, and for $\ID_f$.
\item The functoriality of $U$ follows from the compatibility of $\ID$ with $\circ_v$ and $\id_v$, (\textbf{VDC6}) and (\textbf{VDC5}).
\item The functoriality of $\odot$  follows from square interchange an the compatibility of $\ID$ with $\circ_h$, (\textbf{VDC11}) and (\textbf{VDC5}).
\item The naturality of $\mathfrak{a}$ follows from (\textbf{VDC9}), the compatibility of $\alpha_h$ and $\boxvert$, and Proposition~\ref{swappyprop}.
\item The naturality of $\mathfrak{r}$ and $\mathfrak{l}$ follows from (\textbf{VDC10}), the compatibility of $\rho_h$ and $\lambda_h$ and $\boxvert$, and Proposition~\ref{swappyprop}.
\item The globularity of the components of $\mathfrak{a}, \mathfrak{l},$ and $\mathfrak{r}$ holds directly by definition.
\item The compatibility axioms of $\mathfrak{a}$ and $\mathfrak{r}$ follows from the compatibility of $\alpha_h$ and $\rho_h$ and $\lambda_h$.
\item The pentagon identity for $\mathfrak{a}$ follows from the pentagon identity for $\alpha_h.$
\item That $t_\D$ is a strict functor of pseudo-double categories . The only non-trivial thing to check is  functoriality for squares. For horizontal composition, this follows from $\textbf{VDC}$, the interchange of $\boxvert$ and $\uact.$ The functoriality of $t_\D$ with respect to vertical composition of squares can be shown using \ref{swappyprop}. Let $\alpha:f\Rightarrow g$ and $\beta:g \Rightarrow h.$ we have:
\begin{align*} t_\D(\beta)\bullet t_\D(\alpha) &=   (\ID_h \uact \beta)  \boxminus (\ID _g\uact \alpha) \\ 
& = (\ID_g \dact \beta^{-1})  \boxminus (\ID_g \uact \alpha) \\ 
&= \ID_g \boxminus (\ID_g \uact \alpha)\dact \beta^{-1} \\
&= (\ID_g\dact \beta^{-1}) \uact  \alpha\\
&= (\ID_g\uact \beta) \uact  \alpha \\
&= \ID_g\uact (\beta \uact  \alpha)\\
&= t_\D (\beta\bullet \alpha)
\end{align*}

\end{itemize} 
Note we have only listed  most important axioms used, the other VDC axioms are used throughout.

In the other direction, let $\D$ be a pseudo-double category with thin structure. We can construct a VDC $\mathrm{Triv}(\D) $ as follows:

\begin{itemize} 
\item The vertical $(2,1)$-category $V$ of $\mathrm{Triv}(\D)$ is $\widetilde{\D}_0$, considered as a bicategory with only identity $2$-morphisms
\item The horizontal $(2,1)$-category $H$ of $\mathrm{Triv}(\D)$ is $\overline{\D}$, viewed as a bicategory.
\item The squares of $\mathrm{Triv}(\D)$ are  the morphisms of $\widetilde{\D}_1$ with left and right (source and target) vertical $1$-morphisms given by $L$ and $R$, and top and bottom (source and target) horizontal $1$-morphisms given by the source and target maps of $\widetilde{\D}_1.$
\item The compostion $\boxminus$ is given by the composition of $\widetilde{\D}_1$, while the compostion $\boxvert$ is given by $\odot$
\item The actions $\lact$ and $\ract$ are trivial, since we have only identity vertical $2$-morphisms. The action $\Theta\uact \beta$ and $\Theta \dact \beta$ are given by the composition of $\widetilde{\D}_1$, with $\Theta \uact \beta = \Theta \boxminus t(\beta)$ and  $\Theta \dact \beta = t(\beta^{-1}) \boxminus  \Theta$.
\item The pseudo identity $\ID$ is given for horizontal $1$-morphisms by the identity of $\widetilde{\D}_1$, and for vertical 2-morphisms by $U$. 
\end{itemize}

We leave it to the reader to check that this construction meets the axioms for a Verity double category.

The two constructions we have given are easily seen to be mutually inverse, up to strict isomorphism. We have not defined the notion of a functor of pseudo-double categories, or of fancy pseudo-double categories, but the appropriate definition is of course the one that makes  $\mathrm{Triv}$ and $\trunc$ into equivalences of categories between the category of $1$-trivial VDC's with thin structure and the category of pseudo-double categories.
\begin{definition} A VDC will be called \emph{vertically $0$-trivial} or simply \emph{vertically trivial} if it has only identity vertical $1$-morphism and $2$-morphisms.\end{definition}

Clearly $\mathrm{Triv}$ and $\trunc$ take a vertically trivial VDC to a fancy pseudo-double category with only identities for vertical morphisms. Such pseudo-double categories are equivalent to fancy bicategories , so we abuse notation slightly and take $\mathrm{Triv}$ and $\trunc$ to also indicate the corresponding equivalences between the categories of vertically trivial VDC's and the category of  fancy bicategories.
\begin{definition} An algebraic $2$-reduced inner-Kan bisimplicial set is \emph{vertically $1$-trivial} $X_v = X_{\bullet 0}$ is one-dimensional, equivalently if the face and degeneracy maps of $X_v$ are identities for dimensions greater than $1$. It will be called \emph{vertically $0$-trivial} or simply \emph{vertically trivial} if $X_v$ is $0$-dimensional, equivalently if all the face and degeneracy maps of $X_v$ are identities.
\end{definition}
\begin{theorem}  $\trunc \circ \Vdc$ and $N \circ \mathrm{Triv}$ give an equivalence of categories between the category of vertically trivial algebraic $2$-reduced inner-Kan bisimplicial sets and the category of small fancy bicategories and strictly identity-preserving functors. These equivalences preserve strictness, and the natural isomorphisms exhibiting the equivalence are strict, so that these functors are also equivalences between the strict versions of these categories.
\end{theorem}
\begin{proof}\label{verticallytrivialtheorem} The functors $N$ and $\Vdc$ as constructed in Chapter~\ref{vdcchapter} take a vertically trivial algebraic $2$-reduced inner-Kan bisimplicial set to a small vertically trivial VDC, and vice versa, giving an equivalence of categories by Theorem~\ref{vdcsummary}. The result follows since $\mathrm{Triv}$ and $\trunc$ are equivalences.
\end{proof}
\begin{remark} A similar statement to Theorem~\ref{verticallytrivialtheorem} holds between the category of vertically $1$-trivial algebraic $2$-reduced inner-Kan bisimplicial sets and the category of small pseudo-double categories with thin structure.
\end{remark}

\chapter{The $\Gamma$-nerve for symmetric monoidal groupoids\label{gammachapter}}
\section{Braided and symmetric monoidal groupoids}
It is well-known that a monoidal category, usually defined as a category equipped with a ``tensor'' operation $\otimes$ between objects and morphisms which satisfies  certain axioms,
is equivalent to a bicategory with a single object.

We find it convenient to \emph{define} a small monoidal category to be a  bicategory with a single object, and a \emph{small monoidal groupoid} as a $(2,1)$-category with a single object, denoted $\Box.$
To preserve continuity with our previous notation and avoid tediously introducing a set of new notation that is redundant with what we already have, we will retain the notations we made for the operations of a $(2,1)$-category. For instance, the monoidal product, usually written $A \otimes B$ will be written $A \circ B$ in our notation. However, we shift our usage of the words ``object'' and ``morphism'' to match the usual usage in a small monoidal category, calling the $1$-morphisms of our $(2,1)$-category \emph{objects} and the $2$-morphisms \emph{morphisms}.
\begin{definition} A \emph{small braided monoidal groupoid} is a small monoidal groupoid $\CalC$, together with for every pair $A,B$ of objects a morphism $\gamma_{A,B}: A \circ B \ra B \circ A$, called the \emph{braiding} satisfying the following axioms:
\begin{itemize}
\item \emph{compatibility of $\gamma$ with $\rho$ and $\lambda$ }
\begin{description}
\item[BMG1.]  For all $A$ in $\CalC,$ \tab  $\gamma_{A, \id_{\Box}} \bullet \rho_{A} =\lambda_{A}$
\end{description}
\item \emph{naturality of $\gamma$}
\begin{description}
\item[BMG2.]  For all $A$, $B$, $C$ and $f:A\ra B,$   \tab $\gamma_{B,C} \bullet (f \lhd C)  = (C \rhd f) \bullet \gamma_{A,C}$
\item[BMG3.] For all $A$, $B$, $C$ and $f:B\ra C,$   \tab $\gamma_{A,C} \bullet (A \rhd f)  = (f \lhd A) \bullet \gamma_{A,B}$
\end{description}
\item\emph{hexagon identities}
\begin{description}
\item[BMG4.] For all $A,$ $B,$ and $C$  \tab  $\alpha_{A,C,B} \bullet \gamma_{C \circ B, A} \bullet \alpha_{C,B,A}$

\tab $\quad = (\gamma_{C,A}\lhd B)\bullet\alpha_{C,A,B}\bullet (C \rhd \gamma_{B,A})$
\item[BMG5.] For all $A,$ $B,$ and $C$  \tab  $ \alpha^{-1}_{B,A,C} \bullet \gamma_{C ,B \circ A} \bullet \alpha^{-1}_{C,B,A} $

\tab $\quad= (B\rhd \gamma_{C,A}) \bullet\alpha^{-1}_{B,C,A}\bullet (\gamma_{C,B}\rhd A)$

\end{description}

\end{itemize}

\end{definition}
\begin{definition} A \emph{small symmetric monoidal groupoid} is a small braided monoidal groupoid $\CalC$ such that the braiding $\gamma$ satisfies the following additional axiom:
 \begin{description}
\item[SM.]  For all $A,B$ in $\CalC,$ \tab  $\gamma_{B,A}\bullet\gamma_{A,B}=\Id_{A\circ B}$
\end{description}
\end{definition}
\begin{definition}
A \emph{braided} functor  between braided or small symmetric monoidal groupoids is a functor $F:\CalC \ra \CalD$ between the underlying small monoidal groupoids such that the following axiom holds
 \begin{description}
\item[BMGFun.]  For all $A,B$ in $\CalC,$ \tab  $\gamma_{F(A),F(B)}\bullet\phi_{A,B}=\phi_{B,A}\bullet F(\gamma_{A,B})$
\end{description}
\end{definition}
\begin{remark}\label{redunremark}
The axioms \textbf{SM} and \textbf{BMG2} together can be easily seen to imply \textbf{BMG3}. Similarly, \textbf{SM} and \textbf{BMG4} imply \textbf{BMG5}. Thus if we check $\textbf{SM}$ for a small symmetric monoidal category, we can skip \textbf{BMG3} and \textbf{BMG5}. 
\end{remark}

\section{$\Gamma$-sets \label{gammasection}}
\begin{definition}[Segal] An object in the \emph{category of finite pointed cardinals}, denoted by $\Gamma^{\op}$ is given by a (possibly empty) set $\{1,2,\ldots, n\}$, together with a special ``marked'' point $\star$. A morphism in this category is any (not necessarily order preserving) map $f$ with $f(\star)=\star$. We denote an object $\{1,\ldots,n, \star \}$ of $\Gamma$ by $\overline{n}$. A map $f$ is denoted by $[f(1)f(2)\ldots f(n)]$, e.g. $[31\star32\star].$

The \emph{Segal category} $\Gamma$ is the opposite category of this category.\end{definition}
We will find it convenient to give a direct description of $\Gamma$, which is equivalent to the above definition:
\begin{definition}  \label{gammadef} $\Gamma$ is the category consisting of objects $\underline{n} = \{1,\ldots, n\}$, which correspond to the object $\overline{n}$ in the category of finite pointed sets. A map $\underline{m} \ra \underline{n}$ is a map $g$ from $\underline{m}$ to the power set of $\underline{n}$, such that the $g(i)$ are disjoint. For a map $f$ in the category of finite pointed cardinals, the corresponding morphism in $\Gamma$ is denoted $\hat{f}$ and given by $\hat{f}(i)= f^{-1}(i)$. Note that $f^{-1}(\star)= \underline{n} \setminus \bigcup_{i\in \underline{m}} \hat{f}(i).$ Composition of morphisms $\underline{m} \stackrel{g}{\ra} \underline{n} \stackrel{h}{\ra} \underline{k}$ is given by the formula:
$$\hat{h} \circ \hat{g} (i) =\bigcup_{j\in \hat{g}(i)} \hat{h}(j).$$ 
 \end{definition}
 We make $\Gamma$ into a dimensional category by letting $\dim(\underline{n})=n$.

We adopt a notation for specifying a map $\hat{g}(i) :\underline{m}\ra \underline{n}$ which is best described by example. Let $\hat{g}:\underline{3} \ra \underline{4}$ be given by $\hat{g}(1)=\{1,3\}$, $\hat{g}(2)=\{\}$, and $\hat{g}(3)=\{2\}$. Then we write $$g= [13,-,2].$$  In this notation, order of the elements between commas is irrelevant, but we will adopt the convention of always writing these in order. If context is clear, we abbreviate by writing the same map as $(13)\emptyset 2.$
\begin{definition}The following are notations for special maps of  $\Gamma$ whose domain is a fixed object $\underline{n}$, which is implicit. Collectively, we call these maps the \emph{generators} of $\Gamma$.
\begin{itemize} 
\item \emph{coskip} maps $\hat{k}_i:\underline{n}\ra \underline{n+1}$ given by $\hat{k}_i= [1,2,\ldots,\widehat{i+1},\ldots, n+1 ]$, dual to the \emph{skip} maps $$k_i = [1,\ldots, i, \star, i+1,\ldots, n]$$
\item \emph{comerge} maps  $\hat{m}_i:\underline{n}\ra \underline{n+1}$ given by $\hat{m}_i= [1,2,\ldots,i(i+1),\ldots, n+1]$, dual to the \emph{merge} maps $$m_i = [1,\ldots, i, i, \ldots, n]$$
\item \emph{coinsert} maps $\hat{s}_i:\underline{n} \ra \underline{n-1}$ given by $\hat{s}_i = [1,2, \ldots,i,-,i+1,\ldots, n-1]$, dual to the \emph{insert} maps $$s_i=[1,2,\ldots, \widehat{i+1},\ldots,n]$$
\item \emph{coswap} maps $\hat{w}_i:\underline{n} \ra \underline{n}$ given by $\hat{w}_i = [1,2, \ldots, i-1, i+1, i, i+2,\ldots, n]$, dual to the \emph{swap} maps $$w_i = [1,\ldots, i-1, i+1,i, i+2, \ldots, n]$$
\end{itemize} 
\end{definition}
\begin{proposition} \label{gammarelations}
The generators of $\Gamma$ satisfy the following relations
\begin{itemize} 
\item \emph{coskip and coinsert relations}
\begin{enumerate}
\item $\hat{k}_j\hat{k}_i= \hat{k}_i \hat{k}_{j-1}$, \ $i<j$
\item $\hat{s}_j\hat{s}_i =  \hat{s}_i \hat{s}_{j+1}$, \ $ i \leq j $
\item $\hat{s}_j \hat{k}_i =  \begin{cases}
 			     \id & i=j  \\
			     \hat{k}_i \hat{s}_{j-1}  & i <j \\
			     \hat{k}_{i-1} \hat{s}_j & i>j			      
   \end{cases} $\end{enumerate}
\item \emph{comerge relations}
\begin{enumerate}[resume]
\item $\hat{m}_j \hat{m}_i =				\hat{m}_i \hat{m}_{j-1}$, \ $ i < j $
\item $\hat{m}_j \hat{k}_i =  \begin{cases}
 			     \hat{k}_i\hat{k}_i = \hat{k}_{i+1}\hat{k}_i & i=j-1 \\
			     \hat{k}_i \hat{m}_{j-1}  & i <j-1 \\
			     \hat{k}_{i+1} \hat{m}_{j} & i>j-1			      
   \end{cases} $
\item  $\hat{s}_j \hat{m}_i =  \begin{cases}
 			     \id & i=j,j + 1 \\
			     \hat{m}_i \hat{s}_{j-1}  & i <j \\
			     \hat{m}_{i-1} \hat{s}_j & i>j+1			      
   \end{cases} $
\end{enumerate}

\item \emph{coswap relations}
\begin{enumerate}[resume]
\item  $\hat{w}_j \hat{k}_i =  \begin{cases}
 			     \hat{k}_{i+1} & i=j-1  \\
			     \hat{k}_{i-1} & i =j \\
			      \hat{k}_i \hat{w}_{j-1} & i<j-1		\\
			     \hat{k}_i \hat{w}_j & i>j			      
   \end{cases} $
\item $\hat{s}_j \hat{w}_i =  \begin{cases}
 			     \hat{s}_{j-1} & i=j  \\
			     \hat{s}_{j+1} &  i=j+1 \\
			     \hat{w}_i \hat{s}_j & i < j\\
			     \hat{w}_{i-1} \hat{s}_j &i	>j+1	      
   \end{cases} $
\item $\hat{m}_j \hat{w}_i =  \begin{cases}
 			     \hat{w}_{i+1}\hat{w}_{i} \hat{m}_{j+1} & i=j  \\
			     \hat{w}_{i}\hat{w}_{i+1} \hat{m}_{j-1}     & i =j-1 \\
			     \hat{w}_i \hat{m}_j &  i<j-1 \\                     
           \hat{w}_{i+1}\hat{m}_j & i>j
   \end{cases} $  
\item $\hat{w}_i \hat{m}_i= \hat{m}_i$
\end{enumerate}
\item \emph{symmetric group relations}
\begin{enumerate}[resume]
\item ${\hat{w}_i}^2 =\id$ 
\item $\hat{w}_j \hat{w}_i = \hat{w}_i \hat{w}_j$ \ \ \ $i \neq j-1, j+1$
\item $\hat{w}_i \hat{w}_{i+1} \hat{w}_i = \hat{w}_{i+1} \hat{w}_i \hat{w}_{i+1}$ 
   
\end{enumerate}
\end{itemize}
\end{proposition}
To verify these identities for oneself is easier than to read a proof for them, consequently they are left as an exercise.

\begin{definition}\label{cofulldef}
A map $\hat{f}:\underline{m}\ra\underline{n}$ in $\Gamma$ will be called \emph{coincreasing} if every element of $\hat{f}(j)$ is greater than every element of $\hat{f}(i)$ whenever $j > i.$ If also each $\hat{f}(i)$ has at most one element, we say $\hat{f}(j)$ is \emph{strictly coincreasing}.  Equivalently, $\hat{f}$ is (strictly) coincreasing if $f$ is (strictly) increasing when restricted to $f^{-1}(\overline{n} \setminus \star)$, in which case we say $f$ is \emph{(strictly) increasing}. The \emph{image} of $\hat{f}$ denoted by $\mathrm{Img}_{\hat{f}}:=\bigcup_{i \in \underline{m}} \hat{f}(i).$ We will say $\hat{f}$ is \emph{cofull} if $\mathrm{Img}_{\hat{f}}= \underline{n}.$ Equivalently $\hat{f}$ is cofull if $f$ is \emph{full}, which we take to mean that $f(i) \neq \star$ for $i \neq \star.$ 
\end{definition}
Note that (co)full and (co)increasing morphisms form subcategories in $\Gamma^{\op}$ and $\Gamma$. A morphism $\hat{f}$ in $\Gamma$ is an epimorphism if each $\hat{f}(i)$ has at most one element and $\hat{f}$ is cofull, or equivalently if its dual $f$ is monic. $\hat{f}$ is monic if each $\hat{f}(i)$ is nonempty, equivalently if its dual $f$ is epic. Epicness and monicness for morphism in $\Gamma^{\op}$ of course correspond to surjectivity and injectivity as maps of sets. 

\begin{proposition} A map $f:\overline{m} \ra \overline{m'''}$ in $\Gamma^{\op}$ has a unique factorization $f=SMWK$, where:
\begin{enumerate}
\item $K:\overline{m} \ra \overline{m'}$ is epic and is strictly increasing.
\item $W:\overline{m'} \ra \overline{m'}$ is an isomorphism.
\item $M:\overline{m'} \ra \overline{m''}$ is increasing, full, and epic.
\item $S:\overline{m''} \ra \overline{m'''}$ is increasing, full, and monic.
\item $W$ preserves order within the inverse image of points in $m''$ under $M$ , in the sense that if $M(i)=M(j)$, and $i<j$, then $W^{-1}(i) < W^{-1}(j)$.
\end{enumerate}
We call a map $k$-like if it meets the condition we have imposed on $K$ in $(1)$ above, similarly for $w$, $m,$ and $s$ and conditions $(2)$, $(3)$, $(4)$, respectively. 
\end{proposition}
\begin{proof}
First, we claim $f$ has a unique factorization $f=f'K$, where $K$ is $k$-like and $f'$ is full. Since $f'$ is full, we must have $K^{-1}(\star)= f^{-1}$. There is a unique $k$-like morphism $K$ with domain $\overline{m}$ that has this property, with codomain $\overline{m'}$ where $m'$ is the cardinality of $\overline{m} \setminus f^{-1}(\star)$, and $K$ sends $f^{-1}(\star)$ to $\star$ and the $i$th element of  $\overline{m} \setminus f^{-1}{\star}$, to $i\in \overline{m'}$. $K$ clearly has a left inverse $L$, which is full, and we can let $f'= fL.$ Uniqueness for $f'$ comes from the fact that $S$ is epic, i.e. right-cancellative.  

For the factorization of $f'$, we note that the the subcategory of full morphisms in $\Gamma^{\op}$ is isomorphic to the category $\cat{Fin}$ of finite cardinal numbers and set maps, by an equivalence given by removing the marked point. Applying this equivalence it suffices to factor the image $g$ of $f'$ under this equivalence uniquely as $g=SMW$ where $W$ is an isomorphism, $M$ is increasing and epic, and $S$ is increasing and monic, and $W$ preserves order within the inverse image of points under $M$ . 

First we claim we have a unique factorization $g= g'W$, where $g'$ is increasing, $W$, is an isomorphism, and $W$ preserves order within the inverse image of $g'$. Take $W$ to be the permutation that sorts the points $i$ of the domain $m'=\{1,\ldots,n\}$ into the order in which they are arranged first by the value of $g(i)$ then by the value of $i$. It is easy to check that this is the unique map such that $g W^{-1}$ is increasing, that $W$ preserves order within the inverse image of $g'$. Finally we use the usual epi-monic factorization in $\Delta$ to factor the increasing map $g'$ into a epic increasing map $M$ and a monic increasing map $I$. The condition that $W$ preserves order within the inverse image of points under $g'=SM$ is equivalent to the condition that $W$  preserves order within the inverse image of points under $M$, since $S$ is monic and so the non-empty inverse images of points under $M$ and $SM$ are the same. 
\end{proof}

\begin{corollary}[Four-way factorization] \label{fourwayfactorization} A map $f:\underline{m} \ra \underline{m'''}$ in $\Gamma$ has a unique factorization $f=\hat{K}\hat{W}\hat{M}\hat{S}$, where:
\begin{enumerate}
\item $\hat{K}:\underline{m''} \ra \underline{m'''}$ is coincreasing and each $\hat{K}(i)$ has exactly one element.
\item $\hat{W}:\underline{m''} \ra \underline{m''}$ has the form $[W_1, W_2, \ldots, W_{m''}].$ Equivalently, $W$ is an isomorphism, also equivalently, it is epi and monic. Necessarily in this case $i \ra W_i$ is a permutation of $\underline{m''}$, which we call $\sigma_{\hat{W}}.$
\item $\hat{M}:\underline{m'} \ra \underline{m''}$ is coincreasing, cofull, and monic.
\item $\hat{S}:\underline{m} \ra \underline{m'}$ is coincreasing, cofull, and epic.
\item $\hat{W}$ \emph{preserves order within} $\hat{M}$ meaning that the restriction $\sigma_{\hat{W}}|_{\hat{M}(i)}$ is order-preserving as a map to $\underline{m''}$.
\end{enumerate}
We call a map $\hat{k}$-like if it meets the condition we have imposed on $\hat{K}$ in $(1)$ above, similarly for $\hat{w}$, $\hat{m},$ and $\hat{s}$ and conditions $(2)$, $(3)$, $(4)$, respectively. 
\end{corollary}
\begin{lemma}\label{breakuplemma}~

\begin{itemize}
\item Every $\hat{k}$-like map $f$ has a unique factorization of the form $f=\hat{k}_{i_j}\hat{k}_{i_{j-1}}\ldots \hat{k}_{i_0}$ with $i_0 < i_1 \leq \ldots < i_j$. 
\item Every $\hat{m}$-like map $g$ has a unique factorization of the form  $g=\hat{m}_{i_j}\hat{m}_{i_{j-1}}\ldots \hat{m}_{i_0}$ with $i_0 < i_1< \ldots < i_j$
\item Every $\hat{s}$-like map $h$ has a unique factorization of the form  $h=\hat{s}_{i_j}\hat{s}_{i_{j-1}}\ldots \hat{s}_{i_0}$ with $i_0 > i_1 >\ldots > i_j$
\item We can choose a preferred way of writing each $\hat{w}$-like map $f'$ as a composition of elements $\hat{w}_i.$
\end{itemize}
\end{lemma}
\begin{proof}
For the $\hat{m}$-like map $f:\underline{n} \ra \underline{m}$ let $i_0, i_1,\ldots, i_j$ be the elements of $\underline{m}$ such that $i_k$ and $i_k+1$ are in the same $f(p)$. It is easy to see $f=\hat{m}_{i_j}\hat{m}_{i_{j-1}}\ldots\hat{m}_{i_0},$ this explicit formula for the factorization of a $\hat{m}$-like map at the same time guarantees such a factorization is unique. The $\hat{k}$ and $\hat{s}$ case can likewise be easily handled by explicitly giving the factorization. 

A $\hat{w}$-like map $f'$ is equivalent to the permutation $\sigma_{f'}$ and it is well-known that the symmetric group can be generated by transpositions of adjacent elements, i.e., elements $\sigma_{\hat{w}_i}$. We officially use the Bubble Sort algorithm to give a preferred way of writing a such a permutation as a product of transpositions. This equivalently gives us a preferred way of of writing $f'$ as a product of $\hat{w}_i$ maps.
\end{proof}
\begin{corollary} \label{uniquefactorizationcorollary} Any map $f$ in $\Gamma$ can be uniquely factored as  $$f = \hat{k}_{i^0_{j_0}}\ldots \hat{k}_{i^0_0}\hat{w}_{i^1_{j_1}}\ldots \hat{w}_{i^1_0}\hat{m}_{i^2_{j_2}}\ldots \hat{m}_{i^2_0}\hat{s}_{i^3_{j_3}}\hat{s}_{{i}^3_{{j_3}-1}}\ldots \hat{s}_{i^3_0}$$ such that
 \begin{enumerate}
 \item $i^0_0<  \ldots < i^0_{j_0}$ and $i^2_0 <  \ldots < i^2_{j_1}$ and $i^3_0 >  \ldots > i^3_{j_1}$
 \item $\sigma_{\hat{w}_{i^1_{j_1}}}\ldots \sigma_{\hat{w}_{i^1_0}}$ is a preferred factorization of a permutation
 \item For each $i^2_k$ the permutation $\sigma_{\hat{W}}$ has the property that $$\sigma_{\hat{W}}(i^2_k) < \sigma_{\hat{W}}(i^2_k+1)$$ for all $0\leq k \leq j_2,$ where $\hat{W}:\hat{w}_{i^1_{j_1}}\ldots\hat{w}_{i^1_0}$ 
\end{enumerate}
\end{corollary}
\begin{proof} First we apply Corollary~\ref{fourwayfactorization} to find the unique four-way factorization $f =\hat{K}\hat{W}\hat{M}\hat{S}$, then we apply Lemma~\ref{breakuplemma} to factor these four maps. We must only show that condition $(3)$ above holds if and only if $\hat{W}$ preserves order within $\hat{M}.$ This can be easily seen given the explicit description of the factorization of $\hat{M}$ given in the proof of Lemma~\ref{breakuplemma}. 
\end{proof}
\begin{proposition} A $\Gamma$-set, i.e. a functor $X:\Gamma^{\op} \ra \cat{Set}$, is equivalent to a set $X_n$ for each $n\geq 0$, with operators \begin{align*} k_j&: X_n \ra X_{n-1} \ \ 0\leq j\leq n-1 \\  w_j&: X_n \ra X_n \ \ 1\leq j\leq n-1 \\ m_j&: X_n \ra X_{n-1} \ \ 1\leq j \leq n-1 \\ s_j&:X_n \ra X_{n+1} \ \ 0 \leq j \leq n \end{align*} satisfying the ``opposite versions'' of the relations in Proposition~\ref{gammarelations}, e.g. $k_i k_j= k_{j-1}k_i $  instead of $\hat{k}_j\hat{k}_i= \hat{k}_i \hat{k}_{j-1}$ for $i<j.$ A morphism $X\ra Y$ between $\Gamma$-sets is equivalent to a map $X_i\ra Y_i$ for all $i$, commuting with the operators $k_j, w_j, m_j,$ and $s_j.$
\end{proposition}
\begin{proof}
Proposition~\ref{gammarelations} ensures that if we let $X_n := X(\underline{n})$ and $X(\hat{k}_i)= k_i$, and similarly for $s_i, m_i,$ and $w_i$, we get the data specified in the proposition. To give the inverse operation, constructing a functor $X :\Gamma^{\op} \ra \cat{Set}$ from this data, we set $X(\underline{n}):= X_n.$ To construct the map $X(f)$  we first use Corollary~\ref{uniquefactorizationcorollary} to factor $$f=\hat{k}_{i^0_{j_0}}\ldots \hat{k}_{i^0_0}\hat{w}_{i^1_{j_1}}\ldots \hat{w}_{i^1_0}\hat{m}_{i^2_{j_2}}\ldots \hat{m}_{i^2_0}\hat{s}_{i^3_{j_3}}\ldots \hat{s}_{i^3_0}$$ and then set $$X(f) :=s_{i^3_0}\ldots s_{{i}^3_{j_3}}m_{i^2_0}\ldots m_{{i}^2_{j_2}}w_{i^1_0}\ldots w_{{i}^1_{j_1}}k_{i^0_0}\ldots k_{{i}^0_{j_0}}.$$ 

If $g$ is a map such that $g\circ f$ is defined, we must show that $X(f)\circ X(g)= X(g\circ f).$ Use Corollary~\ref{uniquefactorizationcorollary} to make the factorization $$g=\hat{k}_{i^4_{j_4}}\ldots \hat{k}_{i^4_0}\hat{w}_{i^5_{j_5}}\ldots \hat{w}_{i^5_0}\hat{m}_{i^6_{j_6}}\ldots \hat{m}_{i^6_0}\hat{s}_{i^7_{j_7}}\ldots \hat{s}_{i^7_0}.$$ It is enough to show that the relations of Proposition~\ref{gammarelations} suffice to move from  $$\hat{k}_{i^4_{j_4}}\ldots \hat{k}_{i^4_0}\hat{w}_{i^5_{j_5}}\ldots \hat{w}_{i^5_0}\hat{m}_{i^6_{j_6}}\ldots \hat{m}_{i^6_0}\hat{s}_{i^7_{j_7}}\ldots \hat{s}_{i^7_0}\hat{k}_{i^0_{j_0}}\ldots \hat{k}_{i^0_0}\hat{w}_{i^1_{j_1}}\ldots \hat{w}_{i^1_0}\hat{m}_{i^2_{j_2}}\ldots \hat{m}_{i^2_0}\hat{s}_{i^3_{j_3}}\ldots \hat{s}_{i^3_0}$$ to the factorization of $g \circ f$ described in Corollary~\ref{uniquefactorizationcorollary}. 

It is not hard collect the $\hat{k}$,$\hat{w}$,  $\hat{m}$, and $\hat{s}$ terms together and to put them in the proper order. We must only check that the relations satisfy to move from $$\hat{w}_{i_j} \ldots \hat{w}_{i_0} \hat{m}_{k_l}\ldots \hat{m}_{k_0}$$ meeting conditions $(1)$ and $(2)$ of Corollary~\ref{uniquefactorizationcorollary}, to a factorization also meeting $(3)$ of Corollary~\ref{uniquefactorizationcorollary}, that $\sigma_{\hat{W}}:=\sigma_{\hat{w}_{i_j} \ldots \hat{w}_{i_0}}$ has the property that $$\sigma_{\hat{W}}(k_p) < \sigma_{\hat{W}}(k_p+1)$$ for all $0\leq p \leq l.$

To show this, take the smallest $p$ for which $$\sigma_{\hat{W}}(k_p) > \sigma_{\hat{W}}(k_p+1).$$ We can factor $\hat{W}$ as  $\hat{W}=\hat{W}'\hat{W}''$ where $W'$ preserves order within $\hat{M}:= \hat{m}_{k_l}\ldots \hat{m}_{k_0},$ and $W''$ preserving $\hat{M}$ in the sense that $\hat{W}''\hat{M}=\hat{M},$ which ensures that $\sigma_{W''}$ permutes elements only with the sets $\hat{M}(p).$ Then we can clearly factor $$W''=\hat{w}_{q_r}\ldots \hat{w}_{q_0}$$ where each $\hat{w}_{q_i}$ has the property that $q_i$ and $q_i+1$ are in the same $\hat{M}(p)$ set. By the explicit description of the factorization of a $\hat{m}$-like map given in the proof of Lemma~\ref{breakuplemma}, it follows that for each $q_i$, $\hat{m}_{q_i}$ also appears in our factorization of $\hat{M}$, i.e. $q_i=k_j$ for some $j$. In particular, let $q_0=k_r$. 

Use symmetric group laws to move from  $\hat{w}_{i_j} \ldots \hat{w}_{i_0} \hat{m}_{k_l}\ldots \hat{m}_{k_0}$ to a preferred factorization of $\hat{W'}$ followed by $$\hat{w}_{q_r}\ldots \hat{w}_{q_0}\hat{m}_{k_l}\ldots \hat{m}_{k_0}.$$ Then use the relation  $\hat{m}_j \hat{m}_i =\hat{m}_i \hat{m}_{j-1}$ for  $ i < j $ to move $\hat{m}_{k_r}$ to the front of the $\hat{m}$-terms, yielding $$\hat{w}_{q_r}\ldots \hat{w}_{q_0}\hat{m}_{k_r}\hat{m}_{k_l -1}\ldots\hat{m}_{k_{r+1}-1}\hat{m}_{k_{r-1}}\ldots \hat{m}_{k_0}.$$ Then since $q_0=k_r$ the relation $\hat{w}_i \hat{m}_i=\hat{m}_i$ eliminates the $\hat{w}_{q_0}$ term. Then we can apply the $ \hat{m}_j \hat{m}_i$ relation to move $\hat{m}_{k_r}$ back into place, giving us $$\hat{w}_{q_r}\ldots \hat{w}_{q_1}\hat{m}_{k_r}\hat{m}_{k_l -1}\ldots\hat{m}_{k_{r+1}-1}\hat{m}_{k_{r-1}}\ldots \hat{m}_{k_0}.$$ In this manner we can eliminate each $\hat{w}_{q_i}$ term, leaving us with a preferred factorization of $\hat{W'}$ followed by $\hat{m}_{k_l}\ldots \hat{m}_{k_0}.$ The fact that $\hat{W'}$ preserves order within $\hat{M'}$ guarantees that this factorization meets condition $(3)$ of Corollary~\ref{uniquefactorizationcorollary}. This completes the proof that $X(f)\circ X(g)= X(g\circ f).$

The final statement, that under the equivalence we have constructed a morphism $X\ra Y$ between $\Gamma$-sets is equivalent to a map $X_i\ra Y_i$ for all $i$, commuting with the operators $k_j, w_j, m_j,$ and $s_j$ requires only that every map in $\Gamma$ factors in some way into maps $\hat{k}_j, \hat{w}_j, \hat{m}_j,$ and $\hat{s}_j$, which is shown in  Corollary~\ref{uniquefactorizationcorollary}.
\end{proof}

\section{Inner-Kan $\Gamma$-sets \label{inKangamma}}
 First we adopt some notation specific to $\Gamma$-sets. If $X$ is a $\Gamma$-set, we denote the operator $X(i)$ given by the map $i$ of $\Gamma$ by $\gamma_{i}$. For instance, if it is clear from context we are specifying an operator with a given source, we can use notations like $\gamma_{(13)\emptyset 2}$ to specify operators.
 
 As before, we sometimes want to name a generic cell in $X$ together cells of that can be obtained from that cell with $\gamma_i$ operators. If we give an $n$-cell a name $y_{12\ldots n}$, we will write $y_i$ for $\gamma_i(y_{12\ldots n})$, for instance $y_{(13)\emptyset 2}.$

Recall that a \emph{coface} map in $\Gamma$ is a monic map which increases dimension by $1$. The following proposition justifies our use of Glenn tables in this section.
\begin{proposition}
$\Gamma$ is an excellent dimensional category, as defined in Definition~\ref{cofacesystemdef}.
\end{proposition}
\begin{proof}
It's straightforward to check that $\Gamma$ is a good dimensional category.

$\Gamma$ in fact satisfies a slightly stronger condition than being excellent. If $h$ and $h'$ are distinct coface maps in $\Gamma$, there is a pullback square 
\begin{center}
\begin{tikzcd}
  c'      \arrow{r}{h} & c            \\
c'''   \arrow{u}{g} \arrow{r}{g'}     & c'' \arrow{u}{h'}  
\end{tikzcd}
\end{center}
where $g$ and $g'$ are coface maps. The equivalent dual statement concerning pushouts in $\Gamma^{\op}$, the category of pointed finite cardinals is easy to check and is left to the reader.
%
\end{proof}

We adopt an ordered special coface system for $\Gamma$, for dimensions $\leq 4$ by choosing an ordered list of representatives for each equivalence class of coface maps under the equivalence relation of differing by precomposition with an isomorphism. For $\underline{4}$, we use the following list of cofaces: $$(234,\ (12)34, \ 1(23)4, \ 12(34),\ 123, \ 413, \ (24)13, \ 2(14)3,\ 24(13), 241).$$ For $\underline{3}$ we use $$(23,(12)3, 1(23),12, 13, (13)2),$$  and for $\underline{2}$, we use  $(2,(12),1).$  $\underline{1}$ has a unique coface. The strange convention for $\underline{4}$ has been chosen to make the universal Glenn table given in Table~\ref{gammaglenn4} follow a somewhat more coherent pattern.
 
\begin{definition} \label{innerkangammadef} Let $\underline{m}$ be an element of $\Gamma$ and let $P$ be a partition of the elements $\{1,2,\ldots, m\}$ making up $\underline{m}$ into two nonempty parts, $P_0$ and $P_1$. A coface map $\underline{m-1} \ra \underline{m}$ is said to \emph{respect $P$} if it does not merge of elements from $P_0$ and $P_1$ together, i.e. if it is not of the form $i_0i_1\ldots (i_j i_{j+1})\ldots i_{m-1}$ where $i_j$ and $i_{j+1}$ are in different parts of $P$.

The \emph{universal inner horn} of dimension $m$ associated with $P$, $\Lambda^m_P$ is formed by removing from $\Gamma[m]$ the cofaces which do not respect $P$. That is, $\Lambda^m_P$ is the sub-presheaf of the Yoneda presheaf $\Gamma[m]$ consisting of maps which factor through a coface that respects $P$. Note that a skip map respects any partition.
\end{definition}
Many universal inner horns of the same dimension are isomorphic by maps induced by automorphisms of $\underline{m}.$ For dimension $2$, the only horns are $\Lambda^2_{1|2}$, and $\Lambda^2_{2|1}$ with $1|2$ (for instance) abbreviating the partition $\{\{1\},\{2\}\}.$ These horns are identical by definition. In dimension $3$, all of $\Lambda^3_{1|23}$, $\Lambda^3_{2|13}$, and $\Lambda^3_{3|12}$ are isomorphic. Note that the horn does not depend on the order of the sets in our partition, but we adopt the convention of writing the smaller set first, or if equal, the set containing $1$. Last, in dimension $4$, we have two universal inner horns, up to isomorphism: $\Lambda^3_{1|234}$ and $\Lambda_{12|34}.$

\begin{definition}We choose representatives of the isomorphism classes of universal inner horns, called \emph{(special) universal inner horns}, which are those universal inner horns of the form $\Lambda^m_{12\ldots k|(k+1)\ldots(m)}.$
\end{definition}
\begin{definition} 
If $X$ is a $\Gamma$-set, a \emph{(special) inner horn in $X$} is a map from a (special) universal inner horn to $X$.

A $\Gamma$-set $X$ will be called \emph{inner-Kan} if $X(\underline{0})$ has a unique element, which we denote $\Box$, and all special inner horns (equivalently all inner horns) in $X$ have a filler, i.e. an extension along the natural inclusion $\Lambda^m_P \ra \Gamma[m].$
\end{definition} 
\begin{definition}An \emph{algebraic} inner-Kan $\Gamma$-set is an inner-Kan $\Gamma$-set $X$ equipped with a choice $\chi$ of a filler for each inner horn in $X$. A morphism of algebraic inner-Kan $\Gamma$-sets is just a morphism of the underlying $\Gamma$-sets, and a \emph{strict} morphism is one that preserves the algebraic structure the obvious way, sending preferred fillers to preferred fillers.
\end{definition}
\begin{proposition}\label{mcdprop}
Recall that the \emph{minimal complementary dimension} of a horn $H$ is the smallest dimension among the cells ``missing'' from $H$. The minimum complementary dimension of the horn $\Lambda^{m+n}_{12\ldots n|(n+1)\ldots (n+m)}$, where $n \leq m$, is $m.$
\end{proposition}
\begin{proof} A cell represented by a map $f$ with target $\underline{m+n}$ factors through some face of $\Lambda^{m+n}_{12\ldots n|(n+1)\ldots (n+m)}$ if it fails to be cofull (and thus factors through a face of $\Lambda^{m+n}_{12\ldots n|(n+1)\ldots (n+m)}$ given by a coskip map) or if two elements in the same part of the partition $$12\ldots n|(n+1)\ldots (n+m)$$ are in some $f(p)$. Given this description, it is easy to see that one of the smallest missing cells has the form $$(1(n+1))(2(n+2))\ldots (n(2n))(2n)(2n+1)\ldots(n+m)$$ which has dimension $m.$
\end{proof}
\begin{definition} \label{nreddef} An inner-Kan $\Gamma$-set $X$ will be called $n$-reduced if every horn in $X$ of minimum complementary dimension $n$ or greater has a unique filler. Note that this is consistent with our definition of the term $n$-reduced for simplicial sets and $k$-fold simplicial sets.
\end{definition}

The universal sphere of dimension $3$ is given by $$[23, (12)3, 1(23), 12, 13, (13)2],$$ with:
\begin{align*} 
\Lambda^2_{1|23} &= [23, -, 1(23), 12, 13, -] \\
\Lambda^2_{2|13} &= [23, -, -, 12, 13, (13)2] \\
\Lambda^2_{3|12} &= [23, (12)3,-, 12, 13, -] 
\end{align*}
So we get an inner horn by removing any two of the three inner faces $(12)3,$ $1(23),$ and $(13)2$.

The universal generalized Glenn table for $d\Gamma[3]$ is given in Table~\ref{gammaglenn3}

\begin{table}[H]  \begin{center}\caption{\label{gammaglenn3} The universal Glenn table for $d\Gamma[3]$}

    \begin{tabular}{ r | l || l | l | l | }
     \cline{2-5}

     &$23$   & $3$    &    $(23)$     &$2$            \\ \cline{2-5}
     
    &  $(12)3$ &$3$   &  $(123)$      & $(12) $        \\ \cline{2-5}  

    & $1(23)$     &   $(23)$  &     $(123)$   & $1$         \\   \cline{2-5}  
    
     & $12$     &   $2$  &     $(12)$      &    $1$    \\   \cline{2-5}   
     
     & $13$     &   $3$  &     $(13)$      &    $1$    \\   \cline{2-5}  
   
     &$(13)2$   & $2$  &    $(123)$     &$(13)$         \\ \cline{2-5}
    
    \end{tabular}\end{center}
    \end{table}

To give the universal generalized Glenn table for $d\Gamma[4],$
\begin{table}[H]  \begin{center}\caption{\label{gammaglenn4} The universal Glenn table for $d\Gamma[4]$}

    \begin{tabular}{ r | l || l | l | l | l z l | l | }
     \cline{2-8}

     &$234$   & $34$  &    $(23)4$     &$2(34)$           &$23$   & $24$ & $(24)3$  \\ \cline{2-8}
     
    &  $(12)34$ &$34$   &  $(123)4$   & $(12)(34) $   &     $(12)3$    &$(12)4$  & $(124)3$   \\ \cline{2-8}  

    & $1(23)4$     &   $(23)4$  &     $(123)4$      &    $1(234)$  & $1(23)$   &$14$ & $(14)(23)$    \\   \cline{2-8}  
    
     & $12(34)$     &   $2(34)$  &     $(12)(34)$      &    $1(234)$  & $12$   &$1(34)$ & $(134)2$    \\   \cline{2-8}   
     
     & $123$     &   $23$  &     $(12)3$      &    $1(23)$  & $12$   &$13$ & $(13)2$    \\   \cmidrule[1pt]{2-8}  
   
     &$413$   & $13$  &    $(14)3$     &$4(13)$           &$41$   & $43$ & $(34)1$  \\ \cline{2-8}
     
    &  $(24)13$ &$13$   &  $(124)3$   & $(24)(13) $   &     $(24)1$    &$(24)3$  & $(234)1$   \\ \cline{2-8}  

    & $2(14)3$     &   $(14)3$  &     $(124)3$      &    $2(134)$  & $2(14)$   &$23$ & $(23)(14)$    \\   \cline{3-8}  
    
     & $24(13)$     &   $4(13)$  &     $(24)(13)$      &    $2(134)$  & $24$   &$2(13)$ & $(123)4$    \\   \cline{2-8}   
     
     & $241$     &   $41$  &     $(24)1$      &    $2(14)$  & $24$   &$21$ & $(12)4$    \\   \cline{2-8}   
 
    \end{tabular}\end{center}
    \end{table}
 
\section{Symmetric quasimonoids vs inner-Kan $\Gamma$-sets \label{orbsec}}
\begin{definition}\label{phidefinition}
We give a definition of Segal's functor $\phi:\Delta \ra \Gamma$. We let $\phi([n])=\underline{n}$, which we identify with the set $$\mbox{Pair}(n):=\{(0,1), (1,2),\ldots, (n-1,n)\}$$  of pairs of consecutive integers in $[n]$ (we identify $k\in\underline{n}$ with $(k-1,k)$ in $\mbox{Pair}(n)$). Then for a map $f:[n]\ra[m]$ we say that a pair $(i, i+1)$ in $[n]$ \emph{covers} a pair $(j,j+1)$ in $[m]$ under $f$ if $f(i)\leq j \leq j+1\leq f(i+1)$. Then we define $\phi(f)$ to take a pair $(i,i+1)$ in $\mbox{Pair}(n)$ to the set of pairs it covers under $f$, which is a subset of $\mbox{Pair}(n)$. Using the identification of $\mbox{Pair}(n)$ with $\underline{n}$, this makes $\phi$ a functor from $\Delta$ to $\Gamma$.
$\phi$ induces a pullback functor $\phi^*:\cat{Set}_\Gamma \ra \cat{Set}_\Delta.$ 
\end{definition}
It is easy to see $\phi(\hat{s}_i)=\hat{s}_i,$ $\phi(\hat{d}_0)=\hat{k}_0$, $\phi(\hat{d}_{\max})=\hat{k}_{\max-1}$, and for all other values of $j$, $\phi(\hat{d}_{j})=\hat{m}_{j}.$ This tells us the relation between the operators of a $\Gamma$-set $X$ and the operators of a $\phi^\star(X)$, in particular if $x \in \phi^\star(X)_n$ then $$dx= [k_0 x, m_1 x, \ldots, m_{n-1} x, k_{n-1} x].$$

\begin{conjecture} \label{innerKanconjecture} If $X$ is an inner-Kan $\Gamma$-set, then $\phi^{\star}(X)$ is an inner-Kan simplicial set. \end{conjecture}
A combinatoral proof of this conjecture is likely feasible, though it seems non-trivial. We now give a definition that makes sense in full generality only if Conjecture~\ref{innerKanconjecture} is true.

\begin{definition}If $X$ is a $\Gamma$-set will call the $\Fl_nX:=h_n \phi^{\star}(X)$ the \emph{$n$-flattening} of $X$. (Recall that $h_n$ denotes the ``homotopy $n$-category'' operation on a quasicategory as defined in Definition~\ref{homotopycatdef}).     \end{definition}
For our purposes, we will only use $\Fl_2 X.$ Recall from the construction of $h_i$ in Section~\ref{homotopysection} that $h_2$ depends only on $\tr_3(X),$ and the only the inner-Kan property for $n\leq 4$ is used in the construction. So in order for $\Fl_2X$ to be a well-defined $2$-reduced inner-Kan simplicial set, we will only need to prove the following:
\begin{proposition}If $X$ is inner-Kan, then $\phi^{\star}(X)$ has fillers for inner horns of dimension $4$ or lower.\end{proposition}
\begin{proof} For an inner $2$-horn in $\phi^{\star}(X)$ of the form $[a, - , b]$, the same horn $[a, - , b]$ in $X$ is of type $\Lambda^2_{1|2}$. If we choose a filler $\chi(a,b)$ for this horn, the same cell considered as cell of $\phi^{\star}(X)$ is a filler for the original horn.

Now consider a $\Lambda^3_1$-horn $\frakh$ in $\phi^{\star}(X).$
\begin{table}[H]  \begin{center}\caption{\label{lambda31horn} The  $\Lambda^3_1$-horn $\frakh$ in $\phi^{\star}(X).$ }

    \begin{tabular}{ r | l || l | l | l | }
     \cline{2-5}

                 &$\frakh_{23}$     & $\frakh_3$    &    $\frakh_{(23)}$     &$\frakh_2$            \\ \cline{2-5}
     
    $\Lambda$    &             &$\frakh_3$   &  $\frakh_{(123)}$      & $\frakh_{(12)} $        \\ \cline{2-5}  

                 & $\frakh_{1(23)}$  &   $\frakh_{(23)}$  &     $\frakh_{(123)}$   & $\frakh_1$         \\   \cline{2-5}  
    
                 & $\frakh_{12}$     &   $\frakh_2$  &     $\frakh_{(12)}$      &    $\frakh_1$    \\   \cline{2-5}   
     
    \end{tabular}\end{center}
    \end{table}
    This horn is filled by the cell that fills the following $\Lambda^3_{1|23}$-horn in $X$:
\begin{table}[H]  \begin{center}\caption{\label{31hornfiller} A $\Lambda^3_{1|23}$-horn in $X$ providing the filler of the $\Lambda^3_1$-horn $\frakh$.}

    \begin{tabular}{ r | l || l | l | l | }
     \cline{2-5}

                 &$\frakh_{13}$   & $\frakh_3$    &    $\frakh_{(23)}$     &$\frakh_2$            \\ \cline{2-5}
     
    $\Lambda$    &     &$\frakh_3$   &  $\frakh_{(123)}$      & $\frakh_{(12)} $        \\ \cline{2-5}  

                 & $\frakh_{1(23)}$     &   $\frakh_{(23)}$  &     $\frakh_{(123)}$   & $\frakh_1$         \\   \cline{2-5}  
    
                 & $\frakh_{23}$     &   $\frakh_2$  &     $\frakh_{(12)}$      &    $\frakh_1$    \\   \cline{2-5}   
     
                 & $\chi(\frakh_3,\frakh_2)$     &   $\frakh_3$  &     $\gamma_{(12)}\chi(\frakh_3,\frakh_2)$      &    $\frakh_1$    \\   \cline{2-5}  
    
     $\Lambda$   &   & $\frakh_2$  &    $\frakh_{(123)}$     &$\frakh_{(13)}$         \\ \cline{2-5}
    
    \end{tabular}\end{center}
    \end{table}
The proof for $\Lambda^3_2$-horns in $\phi^{\star}(X)$ is similar, except that it uses a $\Lambda^3_{3|12}$-horn in $X$.

Next we consider a $\Lambda^4_1$-horn $[\frakH_{234},-,\frakH_{1(23)4}, \frakH_{12(34)},\frakH_{123}]$ in $\phi^{\star}(X)$. The following $\Lambda_{1|234}$ in $X$ gives a filling cell for this horn:
\begin{table}[H]  \begin{center}\caption{\label{lambda41proof} The $\Lambda_{1|234}$-horn in $X$ providing the filler of the $\Lambda^4_1$-horn $H$ in $\phi^{\star}(X).$}

    \begin{tabular}{ r | l || l | l | l | l z l | l | }
     \cline{2-8}

     &$\frakH_{234}$   & $\frakH_{34}$  &    $\frakH_{(23)4}$     &$\frakH_{2(34)}$           &$\frakH_{23}$   & $\frakH_{24}$ & $\frakH_{(24)3}$  \\ \cline{2-8}
     
   $\Lambda$ &      &$\frakH_{34}$   &  $\frakH_{(123)4}$   & $\frakH_{(12)(34)} $   &     $\frakH_{(12)3}$    &$\chi(\frakH_{4},\frakH_{(12)})$  &  $\gamma_{(12)3} \mathfrak{r}$   \\ \cline{2-8}  

    & $\frakH_{1(23)4}$     &   $\frakH_{(23)4}$  &     $\frakH_{(123)4}$      &    $\frakH_{1(234)}$  & $\frakH_{1(23)}$   &$\frakH_{14}$ & $\frakH_{(14)(23)}$    \\   \cline{2-8}  
    
     & $\frakH_{12(34)}$     &   $\frakH_{2(34)}$  &     $\frakH_{(12)(34)}$      &    $\frakH_{1(234)}$  & $\frakH_{12}$   &$\frakH_{1(34)}$ & $\frakH_{(134)2}$    \\   \cline{2-8}   
     
     & $\frakH_{123}$     &   $\frakH_{23}$  &     $\frakH_{(12)3}$      &    $\frakH_{1(23)}$  & $\frakH_{12}$   &$\frakH_{13}$ & $\frakH_{(13)2}$    \\   \cmidrule[1pt]{2-8}  
   
     &$=:\mathfrak{p}$   & $\frakH_{13}$  &    $\gamma_{(12)3}\mathfrak{p}$     &$\gamma_{1(23)}\mathfrak{p}$           &$\frakH_{41}$   & $\frakH_{43}$ & $\frakH_{(34)1}$  \\ \cline{2-8}
     
    &  $=:\mathfrak{r}$ &$\frakH_{13}$   &  $\gamma_{(12)3} \mathfrak{r}$          &  $\gamma_{1(23)} \mathfrak{r}$      &       $\gamma_{(12)3} \mathfrak{q}$     &$\frakH_{(24)3}$  & $\frakH_{(234)1}$   \\ \cline{2-8}  

 $\Lambda$   &    & $\gamma_{(12)3}\mathfrak{p}$   &   $\gamma_{(12)3} \mathfrak{r}$          &    $\frakH_{2(134)}$  &$\gamma_{1(23)}\mathfrak{q} $  &$\frakH_{23}$ & $\frakH_{(23)(14)}$    \\   \cline{2-8}  
    
 $\Lambda$    &   &  $\gamma_{1(23)}\mathfrak{p}$   &   $\gamma_{1(23)} \mathfrak{r}$           &    $\frakH_{2(134)}$  & $\frakH_{14}$   &$\frakH_{2(13)}$ & $\frakH_{(123)4}$    \\   \cline{2-8}   
     
     & $=:\mathfrak{q}$     &   $\frakH_{41}$  &     $\gamma_{(12)4} \mathfrak{q}$      &    $\gamma_{1(23)}\mathfrak{q} $  & $\frakH_{24}$   &$\frakH_{21}$ & $\chi(\frakH_{4},\frakH_{(12)})$    \\   \cline{2-8}   
 
    \end{tabular}\end{center}
    \end{table}
In Table~\ref{lambda41proof}, $\mathfrak{p}$ is defined to be the filler of the $\Lambda^3_{2|13}$-horn $$ \frakH:=[ \frakH_{13}  ,\     -    ,\ - ,\ \frakH_{41} ,\ \frakH_{43} ,\ \frakH_{(43)1}]$$
Then we similarly define $\mathfrak{q}$ and then $\mathfrak{r}$ to be fillers of $\Lambda^3_{2|13}$-horns, as given in the table. The cell which fills the horn given in Table~\ref{lambda41proof} also fills $\frakH$ (when viewed as a cell of $\phi^{\star}(X)$), finishing the $\Lambda^4_1$ case. The $\Lambda^4_3$ case is symmetrical to the $\Lambda^4_1$ case, and follows by a similar argument, using a $\Lambda^4_{4|123}$-horn.

Finally, consider a $\Lambda^4_2$-horn  $H:=[\frakH_{234},\frakH_{(12)34},-,  \frakH_{12(34)},\frakH_{123}]$ in $\phi^{\star}(X)$. We use a $\Lambda^4_{12|34}$-horn to construct a filling cell for this horn.
\begin{table}[H]  \begin{center}\caption{\label{lambda42proof} The $\Lambda_{12|34}$-horn in $X$ providing the filler of the $\Lambda^4_2$-horn $H$ in $\phi^{\star}(X).$}
    \begin{tabular}{ r | l || l | l | l | l z l | l | }
     \cline{2-8}

     &$\frakH_{234}$   & $\frakH_{34}$  &    $\frakH_{(23)4}$     &$\frakH_{2(34)}$           &$\frakH_{23}$   & $\frakH_{24}$ & $\frakH_{(24)3}$  \\ \cline{2-8}
     
    &  $\frakH_{(12)34}$    &$\frakH_{34}$   &  $\frakH_{(123)4}$   & $\frakH_{(12)(34)} $   &     $\frakH_{(12)3}$    &$\frakH_{(12)4}$  &  $\frakH_{(124)3}$   \\ \cline{2-8}  

 $\Lambda$   &     &   $\frakH_{(23)4}$  &     $\frakH_{(123)4}$      &    $\frakH_{1(234)}$  & $\frakH_{1(23)}$   &$\chi(\frakH_4,\frakH_1)$ & -    \\   \cline{2-8}  
    
     & $\frakH_{12(34)}$     &   $\frakH_{2(34)}$  &     $\frakH_{(12)(34)}$      &    $\frakH_{1(234)}$  & $\frakH_{12}$   &$\frakH_{1(34)}$ & $\frakH_{(134)2}$    \\   \cline{2-8}   
     
     & $\frakH_{123}$     &   $\frakH_{23}$  &     $\frakH_{(12)3}$      &    $\frakH_{1(23)}$  & $\frakH_{12}$   &$\frakH_{13}$ & $\frakH_{(13)2}$    \\   \cmidrule[1pt]{2-8}  
   
     &$=:\mathfrak{p}$   & $\frakH_{13}$  &    $\gamma_{(12)3}\mathfrak{p}$     &$\gamma_{1(23)}\mathfrak{p}$           &$\frakH_{41}$   & $\frakH_{43}$ & $\frakH_{(34)1}$  \\ \cline{2-8}
     
 $\Lambda$   &    &$\frakH_{13}$   &   $\frakH_{(124)3}$         &  -   &       $\gamma_{(12)3} \mathfrak{q}$     &$\frakH_{(24)3}$  & $\frakH_{(234)1}$   \\ \cline{2-8}  

 $\Lambda$   &    & $\gamma_{(12)3}\mathfrak{p}$   &    $\frakH_{(124)3}$        &    $\frakH_{2(134)}$  &$\gamma_{1(34)}\mathfrak{q} $  &$\frakH_{23}$ & -  \\   \cline{2-8}  
    
 $\Lambda$    &   &  $\gamma_{1(23)}\mathfrak{p}$   &   -        &    $\frakH_{2(134)}$  & $\frakH_{24}$   &$\frakH_{2(13)}$ & $\frakH_{(123)4}$    \\   \cline{2-8}   
     
     & $=:\mathfrak{q}$     &   $\frakH_{41}$  &     $\gamma_{(12)3} \mathfrak{q}$      &    $\gamma_{1(23)}\mathfrak{q} $  & $\frakH_{24}$   &$\frakH_{21}$ & $\frakH_{(12)4}$    \\   \cline{2-8}   
     
     \end{tabular}\end{center}
    \end{table}    

Similarly to above, we define $\mathfrak{p}$ and then $\mathfrak{q}$ to be fillers of $\Lambda^3_{2|13}$-horns, as given in the table. Note since the universal $\Lambda^4_{12|34}$-horn is missing iterated the faces of the form $(24)(13)$ and $(23)(14)$, as these cells of $\Gamma[4]$ do not factor through any face which respects the partition $12|34$, Table~\ref{lambda42proof} completely defines a $\Lambda^4_{12|34}$-horn in $X$. As before, the cell which fills this horn also fills $\frakH$ (when viewed as a cell of $\phi^{\star}(X)$), finishing the $\Lambda^4_2$ case.
\end{proof}

\begin{lemma}\label{partialfiller} Let $X$ be an inner-Kan $\Gamma$-set and let $x_{123}$ be a $3$-cell of $X$ with $$dx_{123} = [x_{23},\ x_{(12)3}, \ x_{1(23)},\ x_{12}, \ x_{13},\ x_{(13)2}].$$ Then if $y$ is a $2$-cell with $dy = [x_3, \gamma_{(12)}y, x_1]$, so that $$\mathfrak{h} := [x_{23},\ - , \ x_{1(23)},\ x_{12}, \ y,\ -]$$ is a $\Lambda_{1|23}$-horn  in $X$, then there is a filler $x'_{123}$ of $\mathfrak{h}$ such that $\gamma_{(12)} x'_{123} = \gamma_{(12)} x_{123}.$ That is, $x'_{123}$ is a filler of $$ [x_{23},\ x_{(12)3} , \ x_{1(23)},\ x_{12}, \ y,\ -].$$
\end{lemma}
\begin{proof} Let $\widehat{\mathfrak{h}}$ be a filler of $\mathfrak{h}.$

\begin{table}[H]  \begin{center}\caption{The $\Lambda^4_{1|234}$-horn defining $x'$}

    \begin{tabular}{ r | l || l | l | l | l z l | l | }
     \cline{2-8}

    & $x_{123}$     &   $x_{23}$  &     $x_{(12)3}$      &    $x_{1(23)}$  & $x_{12}$   &$x_{13}$ & $x_{(13)2}$    \\   \cline{2-8}   
  
  $\Lambda$  & $=:x'_{123}$        &   $x_{23}$  &     $x_{(12)3}$      &    $x_{1(23)}$  & $x_{12}$   &$y$ &    $\gamma_{(12)}\mathfrak{q}$     \\   \cline{2-8}   

    & $s_0 x_{(12)3}$     &  $x_{(12)3}$ &     $x_{(12)3}$     &    $s_0 x_{(123)}$  & $s_0 x_{(12)}$   &$s_0 x_3 $ & $x_{3(12)}$    \\   \cline{2-8}  
    
     & $s_0 x_{1(23)}$     &   $x_{1(23)}$  &    $x_{1(23)}$     &    $s_0 x_{(123)}$  & $s_0 x_1$   &$s_0 x_{(23)}$ & $x_{(23)1}$    \\   \cline{2-8}  
     
     & $s_0 x_{12}$     &   $x_{12}$  &    $x_{12}$     &    $s_0 x_{(12)}$           & $s_0 x_1$     &$s_0 x_2$ & $x_{21}$    \\   \cmidrule[1pt]{2-8}  
   
     &$s_1 x_{32}$   & $s_0 x_2$  &    $x_{32}$     &$x_{32}$           &$s_1 x_3 $   & $x_{32}$ & $s_1 x_{(23)}$   \\ \cline{2-8}
     
    &  $\mathfrak{q}$           &$s_0 x_2$   &  $ \gamma_{(12)}\mathfrak{q}$   & $\gamma_{(23)}\mathfrak{q}$   &      $\gamma_{(12)}\mathfrak{p}$    &$x_{(13)2}$  & $s_1 x_{(123)}$   \\ \cline{2-8}  

  $\Lambda$  &         &    $x_{32}$  &     $\gamma_{(12)}\mathfrak{q}$      &   $x_{1(23)}$  & $\gamma_{(23)}\mathfrak{p}$   &$x_{12}$& $x_{(12)3}$   \\   \cline{2-8}  
    
   $\Lambda$  &       &$x_{32}$     &     $\gamma_{(23)}\mathfrak{q}$      &    $x_{1(23)}$  & $x_{13}$   &$x_{12}$ & $x_{(12)3}$    \\   \cline{2-8}   
     
     & $\mathfrak{p}$     &   $s_1 x_3 $  &    $\gamma_{(12)}\mathfrak{p}$     &$\gamma_{(23)}\mathfrak{p}$     & $x_{13}$   &$s_1 x_1$ & $y$    \\   \cline{2-8}   
 
    \end{tabular}\end{center}
    \end{table}

The $3$-cell $\mathfrak{p}$ is defined to be any filler of the $\Lambda^3_{2|13}$-horn $$[ s_1 x_3  ,   \  -,   \  -,\   x_{13} ,\  s_1 x_1,\   y],$$ which we verify is a bona fide horn by writing its Glenn table in Table~\ref{Ahornmat}.
     \begin{table}[H]\begin{center}\caption{\label{Ahornmat}The generalized Glenn table defining $   \mathfrak{p} $}
      \begin{tabular}{ r | l || l | l | l | }
     \cline{2-5}

     &$s_1 x_2$   & $s_0\Box$   &    $x_2$     &   $x_2$          \\ \cline{2-5}
     
   $\Lambda$ &  $-$ &    $s_0 \Box$   &   $\gamma_{(12)}y$      &  $x_{(02)}$        \\ \cline{2-5}  

 $\Lambda$   & $-$     &   $x_2$  &       $\gamma_{(12)}y$       & $x_0$         \\   \cline{2-5}  
    
     & $x_{02}$     &   $x_2$  &     $x_{(02)}$      &    $x_0$    \\   \cline{2-5}   
     
     & $s_1 x_0$     &   $s_0\Box$  &     $x_0$      &    $x_0$    \\   \cline{2-5}  
   
     &$s_1 x_{(12)}$   & $x_2$  &    $\gamma_{(12)}y$     &$x_0$         \\ \cline{2-5}
    
    \end{tabular}\end{center}
    \end{table}
Similarly,    $\mathfrak{q}$  is defined to be any filler of the $\Lambda^3_{2|13}$-horn $$[ s_0 x_3  ,   \  -,   \  -,\   \gamma_{(23)}   \mathfrak{p} ,\  x_{(24)3},\  s_1 x_{(234)}],$$ though the verification that this makes a horn is left to the reader.

Filling the above $\Lambda^4_{1|234}$-horn gives $x'_{123}$ as desired, as the $\gamma_{(12)34}$ face of the filler.
\end{proof}
 
The simplicial set $\Fl_2 X$ has some extra structure that comes from the $\Gamma$-set structure of $X$. First, the element $\underline{2}$ of $\Gamma$ has a unique nontrivial automorphism $21$, yielding a automorphism $\gamma_{21}$ of the set of $2$-cells of $X$. 
\begin{proposition}
$\gamma_{21}$ induces an automorphism  $\sigma$ of the set of $2$-cells $\Fl_2 X$, with $\sigma^2=\id$ and $$ d (\sigma x) = [ d_3 x,\ d_2 x,\ d_1 x].$$ If $a$ is a $1$-cell, $\sigma s_0 a = s_1 a$.
\end{proposition}
\begin{proof}
The $2$-cells of $\Fl_2 X$ are  equivalence classes of $2$ cells in $\phi^{\star} X$ under the homotopy rel. boundary relation $\sim$ from Section~\ref{homotopysection}, so we must show that $\gamma_{21}$ respects $\sim$ in the sense that if $x\sim x'$ are $2$-cells in $X(\underline{2})=(\phi^\star X)_2$, then $\gamma_{21} x\sim \gamma_{21}x'$.

Taking $2$ cells $x \sim x'$, we use condition $(3)$ of Proposition~\ref{daggericonditions} to see that there is an inner horn $\mathfrak{h}$ in $\phi^\star X$ (without loss of generality a $\Lambda^3_1$-horn) with fillers $\hat{\mathfrak{h}}$ and $\hat{\mathfrak{h}}'$ such that \begin{align*}d_1 \hat{\mathfrak{h}}&=\gamma_{(12)3}\hat{\mathfrak{h}}=x\\ d_1 \hat{\mathfrak{h}}'&=\gamma_{(12)3}\hat{\mathfrak{h}}=x'. \end{align*}Then  $\gamma_{321} \hat{\mathfrak{h}}$ and $\gamma_{321} \hat{\mathfrak{h}}'$ are fillers of the same  $\Lambda^3_2$-horn, with \begin{align*}d_2 \gamma_{321}\hat{\mathfrak{h}}&=\gamma_{1(23)}\gamma_{321}\hat{\mathfrak{h}}=\gamma_{21}x \\ d_2 \gamma_{321}\hat{\mathfrak{h}}&=\gamma_{1(23)}\gamma_{321}\hat{\mathfrak{h}}'=\gamma_{21}x'.\end{align*} We conclude from condition $(3)$ of Proposition~\ref{daggericonditions} that $\gamma_{21} x\sim \gamma_{}x'.$
\end{proof}

\begin{lemma}[Reversal Law]\label{reversal} Let $\fraks=[a,\  b, \ c, \ d ]$ be a commutative sphere in $\Fl_2 X.$ Then $$\sigma \fraks : = [\sigma d,\ \sigma c, \ \sigma b, \ \sigma a]$$ is commutative. 
\end{lemma}
\begin{proof} Let $\widetilde{a}_{12}$, $\widetilde{b}_{12}$, $\widetilde{c}_{12}$, and $\widetilde{d}_{12}$ be representatives  in $X(\underline{2})=(\phi^\star X)_2$ for the $\sim$ equivalence classes $a$, $b$, $c$, and $d$ respectively. The fact that  $\fraks=[a,\  b, \ c, \ d ]$ is commutative in $\Fl_2 X$ implies that there exist $\widetilde{e}_{12}$ and $\widetilde{f}_{12}$ such that $$\mathfrak{p}:=[\widetilde{a}_{12},\ \widetilde{b}_{12},\ \widetilde{c}_{12},\ \widetilde{d}_{12},\ \widetilde{e}_{12},\ \widetilde{f}_{12}]$$ is commutative in $X$. 
Then the commutativity of the sphere $$\gamma_{321}\mathfrak{p}= \mathfrak{p}=[\widetilde{d}_{21},\ \widetilde{c}_{21},\ \widetilde{b}_{21},\ \widetilde{a}_{21},\ \widetilde{e}_{21},\ \widetilde{f}_{12}]$$ ensures that $[\sigma d,\ \sigma c, \ \sigma b, \ \sigma a]$ is commutative in $\Fl_2 X$.
\end{proof}
\begin{lemma}[Transposition Law] \label{magicsword} Suppose $\fraks_0=[a, \ b, \ c, \ d]$ and $\fraks_1=[e,\  b, \ \sigma f, \ \sigma d]$ are commutative in $\Fl_2(X)$. Then $$\fraks_2:=[\sigma a,\ f,\ c, \ e]$$ is a commutative sphere.
\end{lemma}
\begin{proof}Let $\widetilde{a}_{12}$, $\widetilde{b}_{12}$, $\widetilde{c}_{12}$, $\widetilde{d}_{12}$, $\widetilde{e}_{12}$, and $\widetilde{f}_{12}$ be representatives  in $X(\underline{2})=(\phi^\star X)_2$ for the $\sim$ equivalence classes $a$, $b$, $c$, $d$, $e$, and $f$ respectively.

We apply Lemma~\ref{partialfiller} to $[\widetilde{a}_{12},\ \widetilde{b}_{12},\ \widetilde{c}_{12},\ \widetilde{d}_{12},\ \widetilde{e}_{12},\ -]$ and $[\widetilde{e}_{12},\ \widetilde{b}_{12},\ \widetilde{f}_{21},\ \widetilde{d}_{21},\ \widetilde{a}_{12},\ -],$ getting fillers $\mathfrak{p}$ and $\mathfrak{q}$, with 
\begin{align*} d\mathfrak{p} &= [\widetilde{a}_{12},\ \widetilde{b}_{12},\ \widetilde{c}_{12},\ \widetilde{d}_{12},\ \widetilde{e}_{12},\ , \widetilde{f}'_{12}] \\ 
							 d\mathfrak{q} &= [\widetilde{e}_{12},\ \widetilde{b}_{12},\ \widetilde{f}_{21},\ \widetilde{d}_{21},\ \widetilde{a}_{12},\ \widetilde{c}'_{21}].\end{align*}
Then \begin{align*}
d(\gamma_{132} \mathfrak{p})  &= [\widetilde{a}_{21},\ \widetilde{f}'_{12},\   \widetilde{c}_{12},\  \widetilde{e}_{12},\ \widetilde{d}_{12},\ \widetilde{b}_{12}]\\
d(\gamma_{213} \mathfrak{q}) &= [\widetilde{a}_{12},\ \widetilde{b}_{12},\ \widetilde{c}'_{12},\ \widetilde{d}_{12},\    \widetilde{e}_{12},\ ,  \widetilde{f}_{12}] 
\end{align*}
Now we define a $3$-cell $\mathfrak{r}$ in $X$ using the following horn:
\begin{table}[H]  \begin{center}\caption{\label{swordprooftable} The $\Lambda^4_{123|4}$-horn defining $\mathfrak{r}$}

    \begin{tabular}{ r | l || l | l | l | l z l | l | }
     \cline{2-8}

    & $\mathfrak{p}$   &            $ \widetilde{a}_{12}$  &     $ \widetilde{b}_{12}$      &    $\widetilde{c}_{12}$   &  $\widetilde{d}_{12}$   &$\widetilde{e}_{12}$ & $\widetilde{f}'_{12}$    \\   \cline{2-8}   
  
   & $\gamma_{021} \mathfrak{q}$                 &   $\widetilde{a}_{12}$  &     $\widetilde{b}_{12}$    &    $\widetilde{c}'_{12}$    & $\widetilde{d}_{12}$    &$\widetilde{e}_{12}$ & $\widetilde{f}_{12}$    \\   \cline{2-8}   

    &   $s_0 \widetilde{b}_{12}$                &  $\widetilde{b}_{12}$ &     $\widetilde{b}_{12}$     &    $s_0 \widetilde{b}_{(12)}=s_0 \widetilde{c}_{(12)}$  & $s_0 \widetilde{b}_1$   &$s_0 \widetilde{b}_2$ & $\widetilde{b}_{21}$    \\   \cline{2-8}  
    
  $\Lambda$ &                     & $\widetilde{c}_{12}$       &     $\widetilde{c}'_{12}$    &    $s_0 \widetilde{c}_{(12)}$  & $s_0 \widetilde{d}_1=s_0 \widetilde{c}_1$   &$s_0 \widetilde{c}_2$ & $\widetilde{c}_{21}$    \\   \cline{2-8}  
     
   & $s_0 \widetilde{d}_{12}$         &   $\widetilde{d}_{12}$  &    $\widetilde{d}_{12}$     &    $s_0  \widetilde{d}_{(12)}$           & $s_0 \widetilde{d}_1=s_0 \widetilde{c}_1$     &$s_0 \widetilde{d}_2$ & $\widetilde{d}_{21}$    \\   \cmidrule[1pt]{2-8}  
   
     &$s_1 \widetilde{a}_{21}$   & $s_0 \widetilde{d}_2=s_0 \widetilde{a}_1$  &    $\widetilde{a}_{21}$     &$\widetilde{a}_{21}$           &$s_1 \widetilde{b}_2 =s_1 \widetilde{a}_2$   & $\widetilde{a}_{21}$ & $s_1 \widetilde{c}_2=\widetilde{a}_{(12)}$   \\ \cline{2-8}
     
     $\Lambda$           &           &$s_0 \widetilde{d}_2=s_0\widetilde{f}_2$     &  $\widetilde{f}_{12}$   &  $\widetilde{f}'_{12}$      &      $s_1 \widetilde{f}_1=s_1 \widetilde{e}_{(12)}$    &$\widetilde{f}_{12}$  & $s_1 \widetilde{b}_{(12)}=s_1 \widetilde{f}_{(12)}$   \\ \cline{2-8}  

    $\Lambda$    &  $=:\mathfrak{r}$       &    $\widetilde{a}_{21}$  &     $\widetilde{f}_{12}$      &   $\widetilde{c}_{12}$  & $\widetilde{e}_{12}$   &$\widetilde{d}_{12}$& $\widetilde{b}_{12}$   \\   \cline{2-8}  
    
           &  $\gamma_{102}\mathfrak{p}$            & $\widetilde{a}_{21}$     &     $\widetilde{f}'_{12}$      &    $\widetilde{c}_{12}$  & $\widetilde{e}_{12}$   &$\widetilde{d}_{12}$ & $\widetilde{b}_{12}$    \\   \cline{2-8}   
     
     & $s_2 \widetilde{e}_{12}$    &     $s_1 \widetilde{b}_2=s_1 \widetilde{e}_2$  &    $s_1 \widetilde{e}_{(12)}$     &$\widetilde{e}_{12}$     & $\widetilde{e}_{12}$   &$s_1 \widetilde{e}_1$ & $\widetilde{e}_{12}$    \\   \cline{2-8}   
 
    \end{tabular}\end{center}
    \end{table}
    
    Taking equivalence classes under $\sim$, we see that $\mathfrak{r}$ ensures the $[ \sigma a,\  f, \ c,\ e],$ is commutative in $\Fl_2(X)$ by Proposition~\ref{univpropprop}, proving the lemma.
\end{proof}
In order to apply Lemma~\ref{magicsword} conveniently, we introduce a special tabular notation for its application. A \emph{transposition law table} lists two $3$-cells used in the hypothesis of Lemma~\ref{magicsword}, followed by the sphere (marked by the symbol $\utimes$) which we can conclude is commutative:

\begin{table}[H]  \begin{center}\caption{A transposition law table}

    \begin{tabular}{ r | l || l | l | l | l | }
     \cline{2-6}

     &$\fraks_1$   & $e$  &    $b$     &$\sigma f$           &$\sigma d$      \\ \cline{2-6}
     
    &  $\fraks_0$ &$a$   &  $b$   & $c$   &     $d$        \\ \cline{2-6}  

   $\utimes$     &     &  $ \sigma a$ &   $ f$       &    $c$  & $e$        \\   \cline{2-6}  
     
    \end{tabular}\end{center}
    \end{table}

Unlike the simplicial Glenn table, there is no simple trick for remembering the pattern for the transposition law table that allows verification of its correctness. Nevertheless, the pattern of what-matches-what and what-matches-what-up-to-$\sigma$ is small and symmetrical enough to be easily verified.

\begin{definition} \label{symquasidef} A \emph{$2$-reduced symmetric quasimonoid} is a $2$-reduced inner-Kan simplicial set $X$ with a unique $0$-cell, together with a \emph{symmetrizing involution} $\sigma$ of its $2$-cells  with $$d(\sigma x_{12}) =[x_1, x_{(12)}, x_2]$$ and $\sigma s_0= s_1$ which satisfies the Reversal Law and the Transposition Law, Lemmas~\ref{reversal} and \ref{magicsword}. An \emph{algebraic $2$-reduced symmetric quasimonoid} is a $2$-reduced symmetric quasimonoid, together with an algebraic structure on the underlying simplicial set $X$.

A morphism of $2$-reduced symmetric quasimonoids is a map of simplicial sets preserving the involution $\sigma$. A morphism of algebraic $2$-reduced symmetric quasimonoids is called \emph{strict} if it is strict as a morphism of the underlying simplicial sets.
\end{definition}
So we have shown the functor $\Fl_2$ takes any inner-Kan $\Gamma$-set to a $2$-reduced symmetric quasimonoid $(\Fl_2 X,\sigma).$ It is easy to see that this construction is functorial, so we have a constructed a functor $\Fl_2$ from inner-Kan $\Gamma$-sets to $2$-reduced symmetric quasimonoids. Also, it is easy to see that if $X$ is algebraic, then the fillers for $\Lambda^2_{1|2}$-horns in $X$ provide fillers for $\Lambda^2_1$-horns in $\Fl_2 X$, giving the latter an algebraic structure.
\subsection{The functor $\Orb$}
Having constructed the functor $\Fl_2$, we now construct an inverse $\mbox{Orb}$ to this construction, taking a $2$-reduced symmetric quasimonoids to $\Gamma$-sets.  We first define the truncation $\mbox{Orb}|^3_0$ as a $\Gamma|^3_0$-set. For $i<3$, let $\mbox{Orb}|^3_0(Y)(\underline{i}):=Y(i)$. We define the symmetrizing involution of $\Orb(Y)$ by $$\gamma_{21}:= \sigma:\mbox{Orb}(Y)(\underline{1})\ra\mbox{Orb}(Y)(\underline{1}).$$ The $\Gamma$ relations for $i<3$ are easily seen to hold by the simplicial relations for $Y$ together with the hypotheses $d(\sigma x_{12}) =[x_1, x_{(12)}, x_2]$ and $\sigma s_0 =s_1.$

Having defined the truncation $\tr^2(\mbox{Orb}|^3_0(Y))$, we define the $3$ cells of $\mbox{Orb}|^3_0(Y)$ as a subset of the $3$-cells of $\cosk^2\tr^2(\mbox{Orb}|^3_0(Y)),$ i.e. as a subset of $d\Gamma[\underline{3}]$-spheres in $\tr^2(\mbox{Orb}|^3_0(Y)).$ 
\begin{definition} We define a sphere $\fraks:=[a_{12},\ b_{12},\ c_{12},\ d_{12},\ e_{12},\ f_{12}]$ of $\mbox{Orb}|^3_0(Y)$ to be commutative (i.e. in $\mbox{Orb}|^3_0(Y)(3))$ if and only if  both of\begin{align*}\Delta_{\Orb}(\fraks,0)&:=[a_{12},\ b_{12},\ c_{12},\ d_{12}] \\ \Delta_{\Orb}(\fraks,1)&:=[e_{12},\ b_{12},\ \sigma f_{12},\ \sigma d_{12}]\end{align*} are  commutative. Note that the Transposition Law then ensures $$\Delta_{\Orb}(\fraks,2):=[\sigma a_{12},\ f_{12},\ c_{12},\ e_{12}]$$ is also commutative.
\end{definition}

To check this is really a sphere (equivalently, that our definition follows the $\Gamma$-relations for merge and insert operators), we make the generalized Glenn table for $$[a_{12},\ b_{12},\ c_{12},\ d_{12},\ e_{12},\ f_{12}]:$$ 
\begin{table}[H]  \begin{center}\caption{\label{gammaglenny} The Glenn table for $[a_{12},\ b_{12},\ c_{12},\ d_{12},\ e_{12},\ f_{12}]$}

    \begin{tabular}{ r | l || l | l | l | }
     \cline{2-5}

     &$a_{12}$   & $a_2$    &    $a_{(12)}$     &$a_{1}$            \\ \cline{2-5}
     
    &  $b_{12}$ &$b_2$   &  $b_{(12)}$      & $b_1 $        \\ \cline{2-5}  

    & $c_{12}$     &   $c_2$  &     $c_{(12)}$   & $c_1$         \\   \cline{2-5}  
    
     & $d_{12}$     &   $d_2$  &     $d_{(12)}$      &    $d_1$    \\   \cline{2-5}   
     
     & $e_{12}$     &   $e_2$  &     $e_{(12)}$      &    $e_1$    \\   \cline{2-5}  
   
     &$f_{12}$   & $f_2$  &    $f_{(12)}$     &$f_1$         \\ \cline{2-5}
    
    \end{tabular}\end{center}
    \end{table}
For this to be a sphere, we must have certain relations between the faces as can be deduced by comparing this table to the universal Glenn table for $d\Gamma[3]$, Table~\ref{gammaglenn3}. For instance, we must have $a_2=b_2=e_2$. These relations all follow from the fact that $[a_{12},b_{12},c_{12},d_{12}]$ and $[e_{12},b_{12},f_{21},d_{21}]$ are both spheres in $Y$, as can be seen by making the Glenn tables for these spheres, which we leave to the reader.

In order to have defined insert maps $s_j$ and the flip maps $w_j$ for $\mbox{Orb}(Y)$ we note that
 \begin{align*}
 d(s_0 a_{12})&= [a_{12},\ a_{12},\ s_0 a_{(12)},\ s_0 a_1, \ s_0 a_2,\  a_{21}]    \\
 d(s_1 a_{12})&= [s_0 a_2,\ a_{12} ,\ a_{12},\ s_1 a_1,\  a_{12},\ s_1 a_{(12)}] \\
 d(s_2 a_{12})&= [s_1 a_{2},\ s_1 a_{(12)},\ a_{12}, \ a_{12}, \ s_1 a_1, a_{12}] 
 \end{align*}
These identities define the $s_i$ maps for $2$-cells in so that the $\Gamma$-identities involving $s_i$ maps hold, provided these spheres are all commutative according for our definition. For $s_0 a_{12}$ this is equivalent to showing both of the following spheres are commutative in $Y:$ \begin{align*}
&[a_{12},\ a_{12},\ s_0 a_{(12)},\ s_0 a_1] \\ 
&[s_0 a_{2},\ a_{12},\ a_{12},\ s_1 a_1].
\end{align*} These spheres are filled by the cells $s_0 a_{12}$ and $s_1 a_{12}$ respectively. For $s_1 a_{12}$ and $s_2 a_{12}$, we must similarly show the following spheres are commutative in $Y$:
\begin{align*}
&[s_0 a_2,\ a_{12} ,\ a_{12},\ s_1 a_1] \\
&[a_{12},\ a_{12}, \ s_0 a_{(12)}, \ s_0a_1] \\
&[s_1 a_{2},\ s_1 a_{(12)},\ a_{12}, \ a_{12}] \\
&[s_1 a_{1},\ s_1 a_{(12)},\ a_{12}, \ a_{12}] \\
\end{align*}
which are filled by $s_1 a_{12}$, $s_0 a_{12}$, $s_2 a_{12}$, and $s_2 a_{12}$ respectively.

If $x$ is a $3$-cell with $d x= [a, \ b,\  c,\  d,\  e,\  f]$, then 
\begin{align*}
d w_0 x&= d \gamma_{213} x = [e,\ b,\ \sigma f,\ \sigma d, \ a,\ \sigma c]\\
d w_1 x&= d \gamma_{132} x = [\sigma a,\ f,\ c,\ e, \ d,\ b]
\end{align*}
These identities define the $w_i$ maps for the $3$-cell $x$ so that the $\Gamma$-identities involving $w_i$ hold, provided these spheres are both commutative according for our definition. This commutativity is equivalent to the commutativity of the following four spheres in $Y$:
\begin{align*}
&[e, \ b,\ \sigma f,\ \sigma d]\\
&[a,\ b,\ c,\ d ]\\
&[\sigma a,\ f,\ c,\ e]\\
&[d,\ f,\ \sigma b,\ \sigma e]
\end{align*}
The first two of these spheres are seen to be commutative directly from the hypothesis that $[a, \ b,\ c,\  d,\  e,\  f]$ is commutative in $\mbox{Orb}|^3_0(Y)$. The commutativity of the third sphere then follows from applying the Transposition Law to the first two spheres, and commutativity of the last sphere follows from applying the Reversal Law to the first sphere.

We have now defined $\mbox{Orb}|^3_0(Y)$ as a $\Gamma|^3_0$-set. It is easy to see this construction is functorial. We define $\mbox{Orb}=\cosk^3 \mbox{Orb}|^3_0.$
\begin{proposition} Let $Y$ be a $2$-reduced symmetric quasimonoid. Then $\mbox{Orb}(Y)$ is $2$-reduced inner-Kan $\Gamma$-set. 
\end{proposition}
\begin{proof} Lemma~\ref{coskeletaltofiller} ensures from the fact that $\mbox{Orb}(Y)$ is $3$-coskeletal that every horn in $\mbox{Orb}(Y)$ of minimal complementary dimension $4$ or higher has a unique filler. Using Proposition~\ref{mcdprop} to calculate complementary dimension of horns, we have left to check the following (isomorphism classes of) horns: $\Lambda^2_{1|2},$ $\Lambda^3|_{1|23},$ $\Lambda^4_{1|234},$ $\Lambda^4_{12|34}$, $\Lambda^5_{12|345}$ and $\Lambda^6_{123|456}$.

It is clear for $2$-horns in $\mbox{Orb}(Y)$, since $2$-horns and fillers in $Y$ and $\mbox{Orb}(Y)$ are the same thing. By the same observation, if $Y$ has an algebraic structure, we get an algebraic structure on $\mbox{Orb}(Y)$. For $3$-horns, a $\Lambda^3_{1|23}$-horn in $\mbox{Orb}(Y)$ has the form of two horns $[a,\ -,\ c,\ d]$ and $[e,\  b,\ -,\ \sigma d]$, in $Y$ and a filler is equivalent to a filler of both these horns. So since $Y$ has unique fillers for $3$-horns, so does $\mbox{Orb}(Y).$ 

Next, consider the following $\Lambda^4_{1|234}$, horn in $\mbox{Orb}(Y)$. 
\begin{table}[H]  \begin{center}\caption{\label{gammaglennhorn4} A $\Lambda^4_{1|234}$-horn $\mathfrak{H}_{1234}$ in $\Orb(Y)$.}

    \begin{tabular}{ r | l || l | l | l | l z l | l | }
     \cline{2-8}
     &$\frakH_{234}$   & $\frakH_{34}$  &    $\frakH_{(23)4}$     &$\frakH_{2(34)}$           &$\frakH_{23}$   & $\frakH_{24}$ & $\frakH_{(24)3}$  \\ \cline{2-8}
     
  $\Lambda$  & $\frakH_{(12)34}$   &$\frakH_{34}$   &  $\frakH_{(123)4}$   & $\frakH_{(12)(34)} $   &     $\frakH_{(12)3}$    &$\frakH_{(12)4}$  & $\frakH_{(124)3}$   \\ \cline{2-8}  

    & $\frakH_{1(23)4}$     &   $\frakH_{(23)4}$  &     $\frakH_{(123)4}$      &    $\frakH_{1(234)}$  & $\frakH_{1(23)}$   &$\frakH_{14}$ & $\frakH_{(14)(23)}$    \\   \cline{2-8}  
    
     & $\frakH_{12(34)}$     &   $\frakH_{2(34)}$  &     $\frakH_{(12)(34)}$      &    $\frakH_{1(234)}$  & $\frakH_{12}$   &$\frakH_{1(34)}$ & $\frakH_{(134)2}$    \\   \cline{2-8}   
     
     & $\frakH_{123}$     &   $\frakH_{23}$  &     $\frakH_{(12)3}$      &    $\frakH_{1(23)}$  & $\frakH_{12}$   &$\frakH_{13}$ & $\frakH_{(13)2}$    \\   \cmidrule[2pt]{2-8}  
   
     &$\frakH_{413}$   & $\frakH_{13}$  &    $\frakH_{(14)3}$     &$\frakH_{4(13)}$           &$\frakH_{41}$   & $\frakH_{43}$ & $\frakH_{(34)1}$  \\ \cline{2-8}
     
    &  $\frakH_{(24)13}$ &$\frakH_{13}$   &  $\frakH_{(124)3}$   & $\frakH_{(24)(13)} $   &     $\frakH_{(24)1}$    &$\frakH_{(24)3}$  & $\frakH_{(234)1}$   \\ \cline{2-8}  

 $\Lambda$   &   $\frakH_{2(14)3}$  &   $\frakH_{(14)3}$  &     $\frakH_{(124)3}$      &    $\frakH_{2(134)}$  & $\frakH_{2(14)}$   &$\frakH_{23}$ & $\frakH_{(23)(14)}$    \\   \cline{2-8}  
    
 $\Lambda$    &  $\frakH_{24(13)}$   &   $\frakH_{4(13)}$  &     $\frakH_{(24)(13)}$      &    $\frakH_{2(134)}$  & $\frakH_{24}$   &$\frakH_{2(13)}$ & $\frakH_{(123)4}$    \\   \cline{2-8}   
     
     & $\frakH_{241}$     &   $\frakH_{41}$  &     $\frakH_{(24)1}$      &    $\frakH_{2(14)}$  & $\frakH_{24}$   &$\frakH_{21}$ & $\frakH_{(12)4}$    \\   \cline{2-8}   
    \end{tabular}\end{center}
    \end{table}
A filler of such a horn must be unique if it exists by Lemma~\ref{coskeletaltofiller}, so we need only show that the filler does exist. We must check that the spheres $\frakH_{(12)34}$, $\frakH_{2(14)3}$, and $\frakH_{24(13)}$ are commutative in $\Orb(Y)$, given the commutativity of the other seven spheres listed in Table~\ref{gammaglennhorn4}. This works out to checking the commutativity six spheres in $Y$, as follows:

\begin{table}[H]
 \caption{Commutativity of $\Delta_{\Orb}(\frakH_{(12)34},0)$ in $Y$ ~\label{4hornproof1}}\begin{center}
  \scalebox{1}{  \begin{tabular}{ r | l || l | l | l | l | }
     \cline{2-6}

             &    $\Delta_{\Orb}(\frakH_{234},0)$     & $\frakH_{34}$  &    $\frakH_{(23)4}$     &$\frakH_{2(34)}$     &$\frakH_{23}$     \\ \cline{2-6}
     
 $\Lambda$   & $\Delta_{\Orb}(\frakH_{(12)34},0)$     & $\frakH_{34}$   &  $\frakH_{(123)4}$   &$\frakH_{(12)(34)}$    &     $\frakH_{(12)3}$        \\   \cline{2-6}  

            &  $\Delta_{\Orb}(\frakH_{1(23)4},0)$     & $\frakH_{(23)4}$  &     $\frakH_{(123)4}$      &    $\frakH_{1(234)}$  & $\frakH_{1(23)}$       \\   \cline{2-6}  
            &   $\Delta_{\Orb}(\frakH_{12(34)},0)$    & $\frakH_{2(34)}$  &     $\frakH_{(12)(34)}$      &    $\frakH_{1(234)}$  & $\frakH_{12}$    \\ \cline{2-6}

            &  $\Delta_{\Orb}(\frakH_{123},0)$        & $\frakH_{23}$  &     $\frakH_{(12)3}$      &    $\frakH_{1(23)}$  & $\frakH_{12}$   \\   \cline{2-6} 
    \end{tabular}}\end{center}
    \end{table}
    
\begin{table}[H]
 \caption{Commutativity of  $ \Delta_{\Orb}(\frakH_{24(13)},1)$ in $Y$ ~\label{4hornproof2}}\begin{center}
  \scalebox{1}{  \begin{tabular}{ r | l || l | l | l | l | }
     \cline{2-6}

           &    $\sigma \Delta_{\Orb}(\frakH_{123},2)$     &$\frakH_{31}$  &    $\frakH_{(23)1}$    & $\frakH_{2(13)}$     & $\frakH_{23}$    \\ \cline{2-6}
                 
        & $\Delta_{\Orb}(\frakH_{(24)13},2)$   &$\frakH_{31}$    &  $\frakH_{(234)1}$   &   $\frakH_{(24)(13)}$           & $\frakH_{(24)3}$         \\   \cline{2-6}  

       &   $\sigma \Delta_{\Orb}(\frakH_{1(23)4},0)$   &$\frakH_{(23)1}$    &    $\frakH_{(234)1}$   &  $\frakH_{4(123)}$    &  $\frakH_{4(23)}$    \\ \cline{2-6}
   
   $\Lambda$  &   $ \Delta_{\Orb}(\frakH_{24(13)},1)$  &   $\frakH_{2(13)}$   &     $\frakH_{(24)(13)}$       &  $\frakH_{4(123)}$    &    $\frakH_{42}$   \\   \cline{2-6} 
           &  $\sigma \Delta_{\Orb}(\frakH_{234},1)$      &  $\frakH_{23}$ &     $\frakH_{(24)3}$       &     $\frakH_{4(23)}$  &  $\frakH_{42}$   \\   \cline{2-6} 
    \end{tabular}}\end{center}
    \end{table}  
\begin{table}[H]
 \caption{Commutativity of $\Delta_{\Orb}(\frakH_{(12)34},1)$ in $Y$ ~\label{4hornproof3}}\begin{center}
  \scalebox{1}{  \begin{tabular}{ r | l || l | l | l | l | }
     \cline{2-6}

           &    $\sigma\Delta_{\Orb}(\frakH_{241},1)$     &$\frakH_{24}$  &    $\frakH_{(12)4}$    &$\frakH_{1(24)}$    &$\frakH_{12}$      \\ \cline{2-6}
             
           &  $\sigma \Delta_{\Orb}(\frakH_{24(13)},1)$ (See Table~\ref{4hornproof2})  &   $\frakH_{24}$  &     $\frakH_{(123)4}$      &    $\frakH_{(13)(24)}$   & $\frakH_{(13)2}$       \\   \cline{2-6}  
                 
 $\Lambda$ & $\Delta_{\Orb}(\frakH_{(12)34},1)$   &$\frakH_{(12)4}$   &  $\frakH_{(123)4}$   &  $\frakH_{3(124)}$    &     $\frakH_{3(12)}$          \\   \cline{2-6}  

          &   $\sigma \Delta_{\Orb}(\frakH_{(24)13},0)$   &$\frakH_{1(24)}$   &     $\frakH_{(13)(24)}$   &   $\frakH_{3(124)}$   & $\frakH_{31}$    \\ \cline{2-6}
  
           &  $\sigma\Delta_{\Orb}(\frakH_{123},1)$      &  $\frakH_{12}$ &     $\frakH_{(13)2}$       &     $\frakH_{3(12)}$  &  $\frakH_{31}$   \\   \cline{2-6} 
    \end{tabular}}\end{center}
    \end{table}    
    
\begin{table}[H]
 \caption{Commutativity of  $ \Delta_{\Orb}(\frakH_{2(14)3},1)$ in $Y$ ~\label{4hornproof4}}\begin{center}
  \scalebox{1}{  \begin{tabular}{ r | l || l | l | l | l | }
     \cline{2-6}

           &    $\Delta_{\Orb}(\frakH_{123},0)$     &$\frakH_{23}$  &    $\frakH_{(12)3}$    & $\frakH_{1(23)}$     & $\frakH_{12}$    \\ \cline{2-6}
                 
    $\Lambda$    &  $ \Delta_{\Orb}(\frakH_{2(14)3},1)$    &$\frakH_{23}$    &  $\frakH_{(124)3}$   &   $\frakH_{(14)(23)}$           & $\frakH_{(14)2}$         \\   \cline{2-6}  

       &   $\sigma \Delta_{\Orb}(\frakH_{(12)34},1)$  (See Table~\ref{4hornproof3}) & $\frakH_{(12)3}$  &   $\frakH_{(124)3}$ &  $\frakH_{4(123)}$    &  $\frakH_{4(12)}$    \\ \cline{2-6}
   
      &   $\sigma \Delta_{\Orb}(\frakH_{1(23)4},1)$  &   $\frakH_{1(23)}$    &     $\frakH_{(14)(23)}$      & $\frakH_{4(123)}$ &    $\frakH_{41}$   \\   \cline{2-6} 
   
      &  $\sigma \Delta_{\Orb}(\frakH_{241},2)$      & $\frakH_{12}$     &    $\frakH_{(14)2}$      &     $\frakH_{4(12)}$&  $\frakH_{41}$ \\   \cline{2-6} 
    \end{tabular}}\end{center}
    \end{table}  

\begin{savenotes}
\begin{table}[H]
 \caption{Commutativity of $\Delta_{\Orb}(\frakH_{2(14)3},0)$ in $Y$ ~\label{4hornproof5}}\begin{center}
  \scalebox{1}{  \begin{tabular}{ r | l || l | l | l | l | }
     \cline{2-6}

           &    $\Delta_{\Orb}(\frakH_{413},1)$     &$\frakH_{43}$  &   $\frakH_{(14)3}$     &$\frakH_{1(34)}$    &$\frakH_{14}$      \\ \cline{2-6}

 & $\Delta_{\Orb}(\frakH_{(12)34},2)$ \footnote{Follows by the Transposition law from the commutativity of $\Delta_{\Orb}(\frakH_{(12)34},0)$ and $\Delta_{\Orb}(\frakH_{(12)34},1)$, which were shown above in Tables~\ref{4hornproof1}~and~\ref{4hornproof3}}    &$\frakH_{43}$   &   $\frakH_{(124)3}$      & $\frakH_{(12)(34)}$     &     $\frakH_{(12)4}$          \\  \cline{2-6}  
 
  $\Lambda$  &  $\Delta_{\Orb}(\frakH_{2(14)3},0)$ &   $\frakH_{(14)3}$  &     $\frakH_{(124)3}$      &    $\frakH_{2(134)}$   & $\frakH_{2(14)}$       \\   \cline{2-6} 
  
          &   $\Delta_{\Orb}(\frakH_{12(34)},1)$   & $\frakH_{1(34)}$   &     $\frakH_{(12)(34)}$  &    $\frakH_{2(134)}$  & $\frakH_{21}$    \\ \cline{2-6}
  
           &  $\Delta_{\Orb}(\frakH_{241},2)$      &  $\frakH_{14}$ &     $\frakH_{(12)4}$     &     $\frakH_{2(14)}$   &  $\frakH_{21}$   \\   \cline{2-6} 
    \end{tabular}}\end{center}
    \end{table}   
\end{savenotes}

\begin{table}[H]
 \caption{Commutativity of $\Delta_{\Orb}(\frakH_{24(13)},0)$ in $Y$ ~\label{4hornproof6}}\begin{center}
  \scalebox{1}{  \begin{tabular}{ r | l || l | l | l | l | }
     \cline{2-6}
           &    $\Delta_{\Orb}(\frakH_{413},0)$     & $\frakH_{13}$  &    $\frakH_{(14)3}$     &$\frakH_{4(13)}$           &$\frakH_{41}$      \\ \cline{2-6}
             
          & $\Delta_{\Orb}(\frakH_{(24)13},0)$      &$\frakH_{13}$   &  $\frakH_{(124)3}$   & $\frakH_{(24)(13)} $   &     $\frakH_{(24)1}$       \\  \cline{2-6}  
 
          &   $\Delta_{\Orb}(\frakH_{2(14)3},0)$ (See Table~\ref{4hornproof5})  & $\frakH_{(14)3}$   &  $\frakH_{(124)3}$   &    $\frakH_{2(134)}$  &  $\frakH_{2(14)}$     \\ \cline{2-6}  
          
      $\Lambda$      &  $\Delta_{\Orb}(\frakH_{24(13)},0)$ &   $\frakH_{4(13)}$  &     $\frakH_{(24)(13)}$      &    $\frakH_{2(134)}$   & $\frakH_{24}$       \\   \cline{2-6} 
  
           &  $\Delta_{\Orb}(\frakH_{241},0)$      &  $\frakH_{41}$ &     $\frakH_{(24)1}$     &     $\frakH_{2(14)}$   &  $\frakH_{24}$   \\   \cline{2-6} 
    \end{tabular}}\end{center}
    \end{table}
    
Next, we consider a $\Lambda^4_{12|34}$-horn in $\Orb(Y):$  
\begin{table}[H]  \begin{center}\caption{A $\Lambda^4_{12|34}$-horn $\frakH_{1234}$ in $\Orb(Y)$.}

    \begin{tabular}{ r | l || l | l | l | l z l | l | }
     \cline{2-8}
     &$\frakH_{234}$   & $\frakH_{34}$  &    $\frakH_{(23)4}$     &$\frakH_{2(34)}$           &$\frakH_{23}$   & $\frakH_{24}$ & $\frakH_{(24)3}$  \\ \cline{2-8}
     
  & $\frakH_{(12)34}$   &$\frakH_{34}$   &  $\frakH_{(123)4}$   & $\frakH_{(12)(34)} $   &     $\frakH_{(12)3}$    &$\frakH_{(12)4}$  & $\frakH_{(124)3}$   \\ \cline{2-8}  

   $\Lambda$   & $\frakH_{1(23)4}$     &   $\frakH_{(23)4}$  &     $\frakH_{(123)4}$      &    $\frakH_{1(234)}$  & $\frakH_{1(23)}$   &$\frakH_{14}$ & -    \\   \cline{2-8}  
    
     & $\frakH_{12(34)}$     &   $\frakH_{2(34)}$  &     $\frakH_{(12)(34)}$      &    $\frakH_{1(234)}$  & $\frakH_{12}$   &$\frakH_{1(34)}$ & $\frakH_{(134)2}$    \\   \cline{2-8}   
     
     & $\frakH_{123}$     &   $\frakH_{23}$  &     $\frakH_{(12)3}$      &    $\frakH_{1(23)}$  & $\frakH_{12}$   &$\frakH_{13}$ & $\frakH_{(13)2}$    \\   \cmidrule[1pt]{2-8}  
   
     &$\frakH_{413}$   & $\frakH_{13}$  &    $\frakH_{(14)3}$     &$\frakH_{4(13)}$           &$\frakH_{41}$   & $\frakH_{43}$ & $\frakH_{(34)1}$  \\ \cline{2-8}
     
 $\Lambda$   &  $\frakH_{(24)13}$ &$\frakH_{13}$   &  $\frakH_{(124)3}$   & - &     $\frakH_{(24)1}$    &$\frakH_{(24)3}$  & $\frakH_{(234)1}$   \\ \cline{2-8}  

 $\Lambda$   &   $\frakH_{2(14)3}$  &   $\frakH_{(14)3}$  &     $\frakH_{(124)3}$      &    $\frakH_{2(134)}$  & $\frakH_{2(14)}$   &$\frakH_{23}$ & -    \\   \cline{2-8}  
    
 $\Lambda$    &  $\frakH_{24(13)}$   &   $\frakH_{4(13)}$  &    -     &    $\frakH_{2(134)}$  & $\frakH_{24}$   &$\frakH_{2(13)}$ & $\frakH_{(123)4}$    \\   \cline{2-8}   
     
     & $\frakH_{241}$     &   $\frakH_{41}$  &     $\frakH_{(24)1}$      &    $\frakH_{2(14)}$  & $\frakH_{24}$   &$\frakH_{21}$ & $\frakH_{(12)4}$    \\   \cline{2-8}   
    \end{tabular}\end{center}
    \end{table}
To fill this horn, we must find a $p$ and a $q$ such that the following eight spheres are commutative in $Y$:
\begin{align*}
&\Delta_{\Orb}(\frakH_{1(23)4},0) =[  \frakH_{(23)4},\ \frakH_{(123)4},\ \frakH_{1(234)},\ \frakH_{1(23)}  ]   \\
&\Delta_{\Orb}(\frakH_{1(23)4},1) =[ \frakH_{14},\ \frakH_{(123)4},\ \sigma p,\ \frakH_{(23)1}] \\
&\Delta_{\Orb}(\frakH_{(24)13},0)=[  \frakH_{13},\ \frakH_{(124)3},\ q ,\ \frakH_{(24)1}    ] \\
&\Delta_{\Orb}(\frakH_{(24)13},1)=[  \frakH_{(24)3},\ \frakH_{(124)3},\ \frakH_{1(234)} ,\ \frakH_{1(24)}    ] \\
&\Delta_{\Orb}(\frakH_{2(14)3},0)=[  \frakH_{(14)3},\ \frakH_{(124)3},\ \frakH_{2(134)},\ \frakH_{2(14)}] \\
&\Delta_{\Orb}(\frakH_{2(14)3},1)=[  \frakH_{23},\ \frakH_{(124)3},\     p  ,\ \frakH_{(14)2}] \\
&\Delta_{\Orb}(\frakH_{24(13)},0)=[ \frakH_{4(13)},\  q   ,\ \frakH_{2(134)}, \frakH_{24}  ] \\
&\Delta_{\Orb}(\frakH_{24(13)},1)=[ \frakH_{2(13)},\  q   ,\ \frakH_{4(123)}, \frakH_{42}  ]
\end{align*}
First we define the cells $\sigma p$ and $q$ by filling the horns $$[ \frakH_{14},\ \frakH_{(123)4},\ -,\ \frakH_{(23)1}] $$ and $$[  \frakH_{13},\ \frakH_{(124)3},\ -,\ \frakH_{(24)1}]$$ respectively, ensuring that $\Delta_{\Orb}(\frakH_{1(23)4},1)$ and $\Delta_{\Orb}(\frakH_{(24)13},0)$ are commutative. We use Glenn tables to check the commutativity of two more spheres on this list:
\begin{table}[H]
 \caption{Commutativity of $\Delta_{\Orb}(\frakH_{1(23)4},0)$ in $Y$ ~\label{42hornproof1}}\begin{center}
  \scalebox{1}{  \begin{tabular}{ r | l || l | l | l | l | }
     \cline{2-6}

             &    $\Delta_{\Orb}(\frakH_{234},0)$     & $\frakH_{34}$  &    $\frakH_{(23)4}$     &$\frakH_{2(34)}$     &$\frakH_{23}$     \\ \cline{2-6}
     
              & $\Delta_{\Orb}(\frakH_{(12)34},0)$     & $\frakH_{34}$   &  $\frakH_{(123)4}$   &$\frakH_{(12)(34)}$    &     $\frakH_{(12)3}$        \\   \cline{2-6}  

   $\Lambda$    &  $\Delta_{\Orb}(\frakH_{1(23)4},0)$     & $\frakH_{(23)4}$  &     $\frakH_{(123)4}$      &    $\frakH_{1(234)}$  & $\frakH_{1(23)}$       \\   \cline{2-6}  
            &   $\Delta_{\Orb}(\frakH_{12(34)},0)$    & $\frakH_{2(34)}$  &     $\frakH_{(12)(34)}$      &    $\frakH_{1(234)}$  & $\frakH_{12}$    \\ \cline{2-6}

            &  $\Delta_{\Orb}(\frakH_{123},0)$        & $\frakH_{23}$  &     $\frakH_{(12)3}$      &    $\frakH_{1(23)}$  & $\frakH_{12}$   \\   \cline{2-6} 
    \end{tabular}}\end{center}
    \end{table}

\begin{table}[H]
 \caption{Commutativity of  $\Delta_{\Orb}(\frakH_{(24)13},1)$ in $Y$ ~\label{42hornproof2}}\begin{center}
  \scalebox{1}{  \begin{tabular}{ r | l || l | l | l | l | }
     \cline{2-6}

             &    $\Delta_{\Orb}(\frakH_{234},2)$     & $\frakH_{43}$  &    $\frakH_{(24)3}$     &$\frakH_{2(34)}$     &$\frakH_{24}$     \\ \cline{2-6}
     
              & $\Delta_{\Orb}(\frakH_{(12)34},2)$     & $\frakH_{43}$   &  $\frakH_{(124)3}$   &$\frakH_{(12)(34)}$    &     $\frakH_{(12)4}$        \\   \cline{2-6}  

   $\Lambda$    &  $\Delta_{\Orb}(\frakH_{(24)13},1)$     & $\frakH_{(24)3}$  &    $\frakH_{(124)3}$      &    $\frakH_{1(234)}$  & $\frakH_{1(24)}$       \\   \cline{2-6}  
            &   $\Delta_{\Orb}(\frakH_{12(34)},0)$    & $\frakH_{2(34)}$  &     $\frakH_{(12)(34)}$      &    $\frakH_{1(234)}$  & $\frakH_{12}$    \\ \cline{2-6}

            &  $\Delta_{\Orb}(\frakH_{124},0)$        & $\frakH_{24}$  &       $\frakH_{(12)4}$    &    $\frakH_{1(24)}$  & $\frakH_{12}$   \\   \cline{2-6} 
    \end{tabular}}\end{center}
    \end{table}
    
    With these spheres known to be commutative, we can use our calculations above for the $\Lambda^4_{1|234}$ case to finish the proof. For the last two faces $\frakH_{2(14)3}$ and $\frakH_{24(13)}$ the four spheres of $Y$ associated with them, the argument for $\Lambda^4_{1|234}$ shows these could in fact be shown to be commutative even if we did not know that $\frakH_{(12)34}$ was commutative.
    
    For a horn $H:\Lambda^5_{12|345}\ra \mbox{Orb}(Y)$, since $\mbox{Orb}(Y)$ is $3$-coskeletal, it is enough to extend $\tr^3 H:\tr^3\Lambda^5_{12|345}\ra \tr^3 \mbox{Orb}(Y)$ to $\tr^3 \Gamma[5].$ This extension is unique if it exists because $\mbox{Orb}(Y)$ is $2$-subcoskeletal. First we extend $H$ to the cell $\mathfrak{p}:=(13)(24)5$ of $\Gamma[5]$. The $4$-sphere $d\mathfrak{p}$ is contained within (i.e. factors through) $\Lambda^5_{12|345}$,  so to extend $H$ to the cell $\mathfrak{p}$, we must show $H(d\mathfrak{p})$  in $\mbox{Orb}(Y)$.
    
There is a map $$\frakH:\Lambda^4_{1|234}\stackrel{13(24)5}{\lra} \Lambda^5_{12|345}\stackrel{H}{\lra}\mbox{Orb}(Y)$$ where the first map is induced by $13(24)5:\underline{4}\ra\underline{5}$. Furthermore, the sphere $H(d\mathfrak{p})$ factors through $\frakH$, so filling $\frakH$ in $\Orb(Y)$ (as a $\Lambda^4_{1|234}$-horn, this was proved possible above) shows $H(d\mathfrak{p})$ is commutative in $\Orb(Y)$. The same argument works for each $3$-cell in $\Lambda^5_{12|345},$ finishing the extension of $\tr^3H$ to $\tr^3\Gamma[5]$ as desired.

The $\Lambda^6_{123|456}$ case follows by a similar argument to the $\Lambda^5_{12|345}$ case, for instance to extend a  $\Lambda^5_{123|456}$-horn to the cell $\mathfrak{p}:=(14)(25)(36)$ we use a map  $\Lambda^4_{12|345}\ra \Lambda^5_{123|456}$ induced by the map $1245(36):\underline{4}\ra\underline{5}.$
\end{proof}

\begin{theorem}\label{OrbFlequiv} The functors $\Orb$ and $\Fl_2$ are inverse equivalences of categories between the category of $2$-reduced (algebraic) inner-Kan symmetric quasimonoids and the category of $2$-reduced (algebraic) inner-Kan $\Gamma$-sets. Likewise, the functors give an  inverse equivalences of categories between the category of algebraic $2$-reduced inner-Kan symmetric quasimonoids with strict morphisms and the category of algebraic $2$-reduced inner-Kan $\Gamma$-sets with strict morphisms.
\end{theorem}
\begin{proof} We claim that, if $X$ is a $2$-reduced inner-Kan $\Gamma$-set then the relation $\sim$ of homotopy rel. boundary for $2$ cells in $\phi^{\star}(X)$, which is used in Section~\ref{homotopysection} to construct $\Fl_2(X)=h_2 \phi^{\star}(X),$ is the identity relation.

Suppose $a_{12} \sim a'_{12}$ in $\phi^{\star}(X).$ Then by part 1 of Lemma~\ref{daggericonditions}, we conclude that the sphere $[a_{12}, a'_{12}, s_0 a_{(12)},\ s_0{a_1}]$ is commutative in $\phi^{\star}(X),$ meaning there is $3$-cell $\frakp$ in $X$ with $$d\frakp=[a_{12}, a'_{12}, s_0 a_{(12)},\ s_0{a_1},\ e,\ f]$$ for some $e$ and $f$. Then we apply Lemma~\ref{partialfiller} to see that there is a cell $\frakp'$ with $$d\frakp' =  [a_{12}, a'_{12}, s_0 a_{(12)},\ s_0{a_1},\ s_0{a_2},\ f'].$$ Then $d\frakp'$ and $s_0 a$ are both fillers of the $\Lambda^3_{1|23}$-horn $$[a_{12},\ - ,\ s_0 a_{(12)},\ s_0 a_1,\ s_0a_2,\ - ]$$ so from the uniqueness of fillers for $3$-horns in $X$, since $$d(s_0a)= [a_{12},\ a_{12},\ s_0 a_{(12)},\ s_0 a_1, s_0 a_2, a_{21}]$$ we conclude that $a_{12}=a'_{12}$.

Thus, if $X$ is $2$-reduced, by the construction of $h_2$ the $2$-cells of $\Fl_2$ are equivalence classes of $2$-cells in $X$ under the trivial relation $\sim$, so these sets of $2$ cells can be naturally identified with each other, and we treat them as being the same. 

To show the natural equivalence $\Fl_2 \Orb (Y) \cong Y$, we must only show the equivalence up to $3$-cells, since both sides are $3$-coskeletal. Both sides are $2$-reduced, thus $2$-subcoskeletal by Lemma~\ref{fillerstosubcoskeletal}. By the above remarks and the construction of $\Orb$ above, we have identified the sets of $0,1,$ and $2$-cells of the two spaces. Further, the symmetrizing involution $\sigma$ of $Y$ is $\gamma_{21}$ for $\Orb(Y)$ which induces the symmetrizing involution $\sigma$ for $\Fl_2 \Orb(Y),$ showing these two involutions agree. We have left only to show that a sphere of  $\Fl_2 \Orb (Y)$ is commutative if and only if it is commutative as a sphere of $Y$.

A sphere $[a,\ b,\ c,\ d]$ of $\Fl_2 \Orb (Y)$ is commutative if and only if there is an $e$ and an $f$ such that $[a,\ b, \ c, \ d]$ and $[e,\ b,\ \sigma f,\ \sigma d]$ are commutative in $Y$. So we must show that such an $e$ and such an $f$ exist for any commutative sphere $[a,\ b, \ c, \ d].$  To construct such an $e$ and $f$, choose an algebraic structure $\chi$ for $Y$, then let $e= \chi(d_0a, d_2 d)$ and then let $\sigma f$ be defined by the filler of $[e,\ b, - , \sigma d].$

To show the natural equivalence $\Orb \Fl_2 X \cong X,$ we must similarly show a $d\Gamma[\underline{3}]$-sphere in $\Orb \Fl_2 X$ is commutative if and only if the same sphere in $X$ is commutative. $$[a,\ b,\ c,\ d,\ e,\ f]$$ is a commutative sphere in $\Orb \Fl_2 X$ if and only if $[a,\ b,\ c,\ d]$ and $[e,\  b,\  \sigma f,\ \sigma d]$ are commutative in $\Fl_2 X$. If $[a,\ b,\ c,\ d,\ e,\ f]$ is commutative in $X$, and $\frakp$ is the filling cell, then the $\frakp$ fills $[a,\ b,\ c,\ d]$ when viewed as a cell of $\Fl_2 X$, and $\gamma_{01} \frakp$ fills $[e,\  b,\  \sigma f,\ \sigma d].$ In the other direction, if $[a,\ b,\ c,\ d]$  and $[e,\  b,\  \sigma f,\ \sigma d]$ are commutative in $\Fl_2 X$, then by Lemma~\ref{partialfiller}, there commutative spheres in $X$ of the forms \begin{align*} \frakp&:=[a,\ b,\ c,\ d,\ e,\ f']\\ \frakq&:=[e,\  b,\  \sigma f,\ \sigma d,\  a,\ \sigma c'].\end{align*} Then Table~\ref{swordprooftable} above gives a horn which proves that the sphere $$\frakr:=[\sigma a,\ f,\ c, e,\ d,\ b]$$ is commutative in $X$, and then we have
$$d(\gamma_{23} \frakr)=[a,\ b,\ c,\ d,\ e, \ f]$$ showing this sphere is commutative in $X$.

The statement for the algebraic case is immediate since both $\Fl_2$ and $\Orb$ clearly preserve strictness for morphisms.
\end{proof}

\section{The small symmetric monoidal category $\Sym(X)$. \label{Symsection} }
Let $X$ be a $2$-reduced symmetric quasimonoid. Then viewing $X$ as a simplicial set with a unique $0$-cell, we have from Chapter~\ref{bicchapter} a $(2,1)$-category $\Bic(X)$, also having a unique object. 
\subsection{The braiding of $\Sym(X)$}
We now construct a small symmetric monoidal category $\Sym(X)$ with underlying monoidal category $\Bic(X)$.

Only one additional piece of structure is necessary to define $\Sym(X)$:
\begin{definition} The braiding of $\Sym(X)$, given by morphisms $\gamma_{A,B} : A \circ B \ra B \circ A,$ is defined by $\gamma_{A,B}:= \underline{\sigma \chi(A,B)}.$
\end{definition}
\subsection{Verification of the symmetric monoidal axioms for $\Sym(X)$}
We now verify the symmetric monoidal axioms for $\Sym(X)$. Recall from Remark~\ref{redunremark} that we need only check \textbf{BMG1}, \textbf{BMG2}, \textbf{BMG4}, and \textbf{SM}. First we establish a preliminary lemma:
\begin{lemma}\label{hatsigmalemma} Suppose $f:A \ra B$ is a $2$-morphism of $\Sym(X)$, considered as a $2$-cell of $X.$ Then $\sigma f = \widehat{f}.$ 
\end{lemma}
\begin{proof} $$d(\Delta_{\wedge}(f))= [s_1 B,\ \widehat{f},\ f,\ s_0 B].$$ Therefore, by the Matching Lemma, it suffices to show that $$[s_1 B,\ \sigma f, f, s_0 B]$$ is commutative. The following transposition law table verifies this commutativity:
\begin{table}[H]  \begin{center}

    \begin{tabular}{ r | l || l | l | l | l | }
     \cline{2-6}

     &$s_1f$   & $s_0 B$  &    $f$     &$f$           &$\sigma \Id_\id=\Id_\id$      \\ \cline{2-6}
     
    &  $s_1 f$ &$s_0B$  & $f$   & $f$   &     $\Id_\id$        \\ \cline{2-6}  

   $\utimes$     &     &  $ \sigma s_0 B = s_1B$ &   $\sigma f$       &    $f$  & $s_0 B$        \\   \cline{2-6}  
     
    \end{tabular}\end{center}
    \end{table}
\end{proof}
We will consistently use Lemma~\ref{hatsigmalemma} in this section to avoid all use of the $\widehat{f}$ notation, always writing it as $\sigma f$ instead.

\begin{proposition}[\textbf{BMG1} compatibility of the braiding and the unitors] 
\end{proposition}
\begin{proof} Since $\gamma_{A,B}= \underline{\sigma \chi(A,B)}$ and $\lambda_A =  \underline{s_1 A}$, by Lemma~\ref{bulletylemma}\footnote{asserting that $\underline{x\bullet \eta } = \underline{x} \bullet \eta$}, it suffices to show that   $$(\sigma \chi_{A, I}) \bullet \rho_{A} =s_1 A.$$  We have $$d ( \Delta_{\bullet} (\sigma \chi(A,\id), \rho_A) ) =[\sigma \chi(A, \id) ,\  \sigma \chi(A, \id) \circ \rho_{A},\ \rho_A, \Id_{A} ]$$ so by the Matching Lemma, it suffices to show the following sphere is commutative:  $$d ( \Delta_{\bullet} (\sigma \chi(A,\id), \rho_A) ) =[\sigma \chi(A, \id) ,\  \sigma \chi(A, \id) \circ \rho_{A},\ \rho_A, \Id_{A}].$$ We conclude the proof with the following transposition law table verifying this commutivity:

\begin{table}[H]  \begin{center}
   \begin{tabular}{ r | l || l | l | l | l | }
     \cline{2-6}

     &$s_0 \Id_A $   & $\Id_A$  &    $\Id_A$     &$s_0 A$           &$\sigma\Id_\id= \Id_\id$      \\ \cline{2-6}
     
    &  $S_1$ &$ \chi(A,\id)$   &  $\Id_A$   & $\rho_A$   &     $\Id_\id$        \\ \cline{2-6}  

   $\utimes$     &     &  $ \sigma \chi(A,\id)$ &   $\sigma (s_0 A) = s_1 A$       &    $\rho_A$  & $\Id_A$        \\   \cline{2-6}  
     
    \end{tabular}\end{center}
    \end{table}
\end{proof}
\begin{proposition}[\textbf{BMG2} naturality of the braiding in the first variable] \label{symnatfirst} Let $A$, $B$, and $C$ be objects of $\Sym(X)$, and let $f:A\ra B$ be a morphism. Then $$\gamma_{B,C} \bullet (f \lhd C)  = (C \rhd f) \bullet \gamma_{A,C}.$$
\end{proposition}
\begin{proof}  Using Lemma~\ref{bulletylemma}, it suffices to show $$\sigma (\chi(B,C)) \bullet (f \lhd C)  = (C \trhd f) \bullet \gamma_{A,C}.$$ We have $$ d(\Delta_{\bullet}(C \trhd f,\gamma_{A,C})= [C \trhd f,\  (C \trhd f)\bullet \gamma_{A,C},\ \gamma_{A,C},\ \Id_B].$$ By the Matching Lemma, it is enough to show the sphere $$[C \trhd f,\  (\sigma \chi(B,C))\bullet (f \lhd C),\ \gamma_{A,C},\ \Id_B]$$ is commutative. The following Glenn table proof verifies this commutativity:

\begin{table}[H] \caption{~}\label{isopt1}\begin{center}

    \begin{tabular}{ r | l || l | l | l | l | }
     \cline{2-6}

      & $\Delta_{\trhd}(C,f)$    &   $\Id_{C}$  &  $\chi(C,A)$       & $C \trhd f $           &$\sigma f$        \\ \cline{2-6}

   $\odot$       &      (Table~\ref{isopt2})                 &  $\Id_C$                        &          $\sigma \chi(A,C)$    &                 $(\sigma \chi(B,C)) \bullet (f \lhd C)$                  &      $\sigma f$   \\   \cline{2-6}  

           &   $\Delta_{-}(\sigma \chi(A,C))$       &    $\chi(C,A)$     & $\sigma \chi(A,C)$                  &   $\gamma_{A,C}$                                   & $\Id_A$         \\   \cline{2-6}  
        
    $\Lambda$       &         & $C\trhd f$                       &    $(\sigma \chi(B,C)) \bullet (f \lhd C)$      &     $\gamma_{A,C}$   &$\Id_B$      \\ \cline{2-6}

     &  $s_0(\sigma f)$                                 &   $\sigma f$                        &    $\sigma f$          &$\Id_A$     &$\Id_B$    \\   \cline{2-6} 
    \end{tabular}\end{center}
    \end{table}
 \begin{table}[H] \caption{~}\label{isopt2}\begin{center}
\scalebox{.97}{
    \begin{tabular}{ r | l || l | l | l | l | }
     \cline{2-6}

      & $s_1(\sigma\chi(B,C))$    &   $s_0 C$   &  $\sigma\chi(B,C)$       &$\sigma\chi(B,C)$          &$s_1 B$        \\ \cline{2-6}

   $\Lambda$      &        &  $s_0 C=\Id_C$         &          $\sigma \chi(A,C)$    &         $(\sigma \chi(B,C)) \bullet (f \lhd C)$                  &      $\sigma f$   \\   \cline{2-6}  

           &    $\sigma (\Delta_{\lhd}(f,C))  $           &    $\sigma\chi(B,C)$         & $\sigma\chi(A,C)$      & $f\lhd C$    &           $f$        \\  \cline{2-6}

           & $\Delta_{\bullet}( \sigma \chi(B,C),f\lhd C )$     &  $\sigma \chi(B,C)$       &         $(\sigma \chi(B,C)) \bullet (f \lhd C)$    &   $f\lhd C$    &        $s_0 B$ \\   \cline{2-6}  
  
     &  $\Delta_{\wedge}(f)$                                 &   $s_1 B$                        &      $\sigma f$        &  $f$  &  $s_0 B$  \\   \cline{2-6} 
    \end{tabular}}\end{center}
    \end{table}\end{proof}

\begin{proposition}[\textbf{BMG4} first hexagon identity] \label{firsthexiden} Let $A,B,C$ be objects of $\Sym(X).$ Then $$ \alpha_{A,C,B} \bullet \gamma_{C \circ B, A} \bullet \alpha_{C,B,A} = (\gamma_{C,A}\lhd B)\bullet\alpha_{C,A,B}\bullet (C \rhd \gamma_{B,A}).$$ 
\end{proposition}
\begin{proof}
By Lemma~\ref{bulletylemma}, it suffices to show  $$ \widetilde{\alpha}_{A,C,B} \bullet \gamma_{C \circ B, A} \bullet \alpha_{C,B,A} = (\gamma_{C,A}\tlhd B)\bullet\alpha_{C,A,B}\bullet (C \rhd \gamma_{B,A}).$$ Applying the Matching Lemma to $$d(\Delta_{\bullet}(\gamma_{C,A}\tlhd B,\alpha_{C,A,B}\bullet (C \rhd \gamma_{B,A}))),$$ we see that it suffices to show the following sphere is commutative: $$[\gamma_{C,A}\tlhd B ,\ \widetilde{\alpha}_{A,C,B} \bullet \gamma_{C \circ B, A} \bullet \alpha_{C,B,A} ,\ \alpha_{C,A,B}\bullet (C \rhd \gamma_{B,A}) ,\ \Id_B].$$
The following table proof, involving both Glenn and transposition law tables, verifies this commutativity. Note that where necessary we omit subscripts to save space.

\begin{table}[H] \caption{~ }\label{axHex1pt1}\begin{center}

 \scalebox{.92}{   \begin{tabular}{ r | l ||  l | l | l |l | }
     \cline{2-6}

                    & $\Delta_{\tlhd}(\gamma_{C,A} , B)$        &$\gamma_{C,A}$   &  $\gamma_{C,A}\tlhd B$      &   $\chi(C\circ A,B)$           &$s_1 B$     \\ \cline{2-6}
     
            $\odot$  &  (Table \ref{axHex1pt2}) & $\gamma_{C,A}$    & $\widetilde{\alpha}_{A,C,B} \bullet \gamma_{C \circ B, A} \bullet \alpha_{C,B,A} $   &     $\widetilde{\alpha}_{C,A,B}\bullet (C \rhd \gamma_{B,A})$     &        $s_1 B$   \\   \cline{2-6}  

      $\Lambda$        &        &$\gamma_{C,A}\tlhd B$    &$\widetilde{\alpha}_{A,C,B} \bullet \gamma_{C \circ B, A} \bullet \alpha_{C,B,A} $      &   $\alpha_{C,A,B}\bullet (C \rhd \gamma_{B,A})$    &$ \Id_B$          \\   \cline{2-6}  
     
    & $\Delta_{-}(\widetilde{\alpha}_{C,A,B}\bullet (C \rhd \gamma_{B,A}))$    & $\chi(C\circ A,B)$    &   $\widetilde{\alpha}_{C,A,B}\bullet (C \rhd \gamma_{B,A})$   &     $\alpha_{C,A,B}\bullet (C \rhd \gamma_{B,A})$    &    $ \Id_B$\\ \cline{2-6}

     &  $s_2(\Id_B)$   & $s_1 B$ &   $s_1 B$  &  $\Id_B$& $\Id_B$   \\   \cline{2-6} 
    \end{tabular}}\end{center}
    \end{table}
\begin{table}[H] \caption{~ }\label{axHex1pt2}\begin{center}
 \scalebox{.95}{   \begin{tabular}{ r | l ||  l | l | l |l | }
     \cline{2-6}

     $\odot_1$ &   (Table \ref{axHex1pt31})   & $\sigma \chi(A,C)$   &  $\chi(C,A)$      &   $\gamma_{C,A}$       &$\Id_A$     \\ \cline{2-6}
     
     $\odot_3$         &  (Table \ref{axHex1pt33}) & $\sigma \chi(A,C)$     & $\chi(C,B\circ A) $   &     $\widetilde{\alpha}_{A,C,B} \bullet \gamma_{C \circ B, A} \bullet \alpha_{C,B,A} $  &       $\sigma\chi(B,A)$   \\   \cline{2-6}  

      $\odot_2$     &   (Table \ref{axHex1pt32})      &$\chi(C,A)$   &$\chi(C,B\circ A) $   &   $\widetilde{\alpha}_{C,A,B}\bullet (C \rhd \gamma_{B,A})$    &    $\sigma\chi(B,A)$      \\   \cline{2-6}  
     
   $\Lambda$  &           & $\gamma_{C,A}$       &    $\widetilde{\alpha}_{A,C,B} \bullet \gamma_{C \circ B, A} \bullet \alpha_{C,B,A} $         &   $\widetilde{\alpha}_{C,A,B}\bullet (C \rhd \gamma_{B,A})$                    &      $s_1 B$   \\   \cline{2-6} 
     
                &   $s_1(\sigma\chi(B,A))$     & $\Id_A$    & $\sigma\chi(B,A)$    &     $\sigma\chi(B,A)$     &        $s_1 B$ \\ \cline{2-6}
    \end{tabular}}\end{center}
    \end{table}
\begin{table}[H] \caption{Commutativity of $\odot_1$ in Table~\ref{axHex1pt2} }\label{axHex1pt31}\begin{center}
 \scalebox{.95}{       \begin{tabular}{ r | l || l | l | l | l | }
     \cline{2-6}

     &$s_1 (\sigma \chi(C,A))$     & $\Id_A$  &    $\sigma \chi(C,A)$      &$\sigma \chi(C,A)$            &$s_1 C=\sigma \Id_C$      \\ \cline{2-6}
     
    &  $\Delta_{-}(\sigma\chi(C,A))$     &$ \chi(A,C)$   &  $\sigma \chi(C,A)$                     &  $\gamma_{C,A}$        &     $\Id_C$        \\ \cline{2-6}  

   $\utimes$     &     & $\sigma \chi(A,C)$   &  $\chi(C,A)$      &   $\gamma_{C,A}$       &$\Id_A$   \\   \cline{2-6}  
     
    \end{tabular}}\end{center}
    \end{table}
\begin{table}[H] \caption{Commutativity $\odot_2$ in Table~\ref{axHex1pt2}}\label{axHex1pt32}\begin{center}
 \scalebox{.93}{   \begin{tabular}{ r | l ||  l | l | l |l | }
     \cline{2-6}

         &   $\Delta_{\widetilde{\alpha}}(C,A,B)$   & $\chi(C,A)$   &  $\chi(C,A \circ B)$      &  $ \widetilde{\alpha}_{C,A,B}$      &$\chi(A,B)$  \\ \cline{2-6}
   
     $\Lambda$   &       &$\chi(C,A)$   &$\chi(C,B\circ A) $   &   $\widetilde{\alpha}_{C,A,B}\bullet (C \rhd \gamma_{B,A})$    &    $\sigma\chi(B,A)$      \\   \cline{2-6}     
    
        & $\Delta_{\rhd}(C,\gamma_{B,A})$ & $\chi(C,A\circ B)$      & $\chi(C,B\circ A) $   &     $C \rhd \gamma_{B,A} $  &       $\gamma_{B,A}$   \\   \cline{2-6}

       & $\Delta_{\bullet}( \widetilde{\alpha}_{C,A,B} ,\  C \rhd \gamma_{B,A}) $     &$\widetilde{\alpha}_{C,A,B} $     &     $\widetilde{\alpha}_{C,A,B}\bullet (C \rhd \gamma_{B,A})$       &   $C \rhd \gamma_{B,A}$                    &      $s_0 B$   \\   \cline{2-6} 
     
                &   $\Delta_{-}(\sigma\chi(B,A))$     & $\chi(A,B)$   & $\sigma\chi(B,A)$ &     $\gamma_{B,A}$     &        $s_0 B$ \\ \cline{2-6}
    \end{tabular}}\end{center}
    \end{table}
\begin{table}[H] \caption{Commutativity of $\odot_3$ in Table~\ref{axHex1pt2}, step 1 }\label{axHex1pt33}\begin{center}
 \scalebox{.95}{       \begin{tabular}{ r | l || l | l | l | l | }
     \cline{2-6} 
     &$\sigma\Delta_{\widetilde{\alpha}}(C,B,A)$  &$\sigma\chi(B,A)$  &  $\sigma \widetilde{\alpha}_{C,B,A}$  &   $\sigma \chi(C,B\circ A) $   &$\sigma \chi(C,B)$                  \\ \cline{2-6}
     
 $\odot$     & (Table \ref{axHex1pt34})    & $\chi(A,C)$   &  $\sigma \widetilde{\alpha}_{C,B,A}$  &   $\widetilde{\alpha}_{A,C,B} \bullet \gamma_{C \circ B, A} \bullet \alpha_{C,B,A} $      &$\chi(C,B)$   \\   \cline{2-6}  
  
 $\utimes$      &  &  $\sigma \chi(A,C)$     & $\chi(C,B\circ A) $   &     $\widetilde{\alpha}_{A,C,B} \bullet \gamma_{C \circ B, A} \bullet \alpha_{C,B,A} $  &       $\sigma\chi(B,A)$     \\ \cline{2-6}  
    \end{tabular}}\end{center}
    \end{table}
\begin{table}[H] \caption{Commutativity of $\odot_3$ in Table~\ref{axHex1pt2}, step 2 }\label{axHex1pt34}\begin{center}
 \scalebox{.93}{   \begin{tabular}{ r | l ||  l | l | l |l | }
     \cline{2-6}

         &   $\Delta_{\widetilde{\alpha}}(A,C,B)$   &  $\chi(A,C)$   &  $\chi(A,C \circ B)$      &  $ \widetilde{\alpha}_{A,C,B}$      &$\chi(C,B)$  \\ \cline{2-6}
   
     $\Lambda$   &       & $\chi(A,C)$   &  $\sigma \widetilde{\alpha}_{C,B,A}$  &   $\widetilde{\alpha}_{A,C,B} \bullet \gamma_{C \circ B, A} \bullet \alpha_{C,B,A} $      &$\chi(C,B)$   \\   \cline{2-6}    
    
      $\odot$  &  (Table \ref{axHex1pt35}) & $\chi(A,C\circ B)$      &  $\sigma \widetilde{\alpha}_{C,B,A}$  &    $\gamma_{C \circ B, A} \bullet \alpha_{C,B,A} $    &       $s_0(C\circ B)$   \\   \cline{2-6}

       & $\Delta_{\bullet}( \widetilde{\alpha} ,\  \gamma \bullet \alpha ) $     &$\widetilde{\alpha}_{A,C,B} $     &       $\widetilde{\alpha}_{A,C,B} \bullet \gamma_{C \circ B, A} \bullet \alpha_{C,B,A} $    &   $\gamma_{C \circ B, A} \bullet \alpha_{C,B,A} $                    &      $s_0B$   \\   \cline{2-6} 
     
                &   $s_0(\chi(C,B))$     & $\chi(C,B)$   & $\chi(C,B)$ &     $s_0(C\circ B)$     &        $s_0B$   \\ \cline{2-6}
    \end{tabular}}\end{center}
    \end{table}
\begin{table}[H] \caption{Commutativity of $\odot_3$ in Table~\ref{axHex1pt2}, step 3 }\label{axHex1pt35}\begin{center}
 \scalebox{.93}{   \begin{tabular}{ r | l ||  l | l | l |l | }
     \cline{2-6}

         &   $\Delta_{-}(\sigma\chi(C\circ B, A))$   &  $\chi(A,C\circ B)$   &  $\sigma\chi(C\circ B, A)$      &  $\gamma_{C \circ B, A}$      &$s_0(C\circ B)$  \\ \cline{2-6}
   
     $\Lambda$   &        & $\chi(A,C\circ B)$      &  $\sigma \widetilde{\alpha}_{C,B,A}$  &    $\gamma_{C \circ B, A} \bullet \alpha_{C,B,A} $    &       $s_0(C\circ B)$   \\   \cline{2-6}  
    
      &  $\sigma\Delta_{\widehat{-}}(\widetilde{\alpha}_{C,B,A}) $ (See below)  & $\sigma\chi(C\circ B, A).$    &  $\sigma \widetilde{\alpha}_{C,B,A}$  &   $ \alpha_{C,B,A}  $   &       $s_0(C\circ B)$   \\   \cline{2-6}

       & $\Delta_{\bullet}(\gamma_{C \circ B, A} ,\  \alpha_{C,B,A}) $     &$\gamma_{C \circ B, A}$     &         $\gamma_{C \circ B, A} \bullet \alpha_{C,B,A} $  &   $ \alpha_{C,B,A}  $                    &      $s_0s_0\Box$   \\   \cline{2-6} 
     
                &   $s_0s_0(C\circ B)$    & $s_0(C\circ B)$   & $s_0(C\circ B)$&     $s_0(C\circ B)$     &        $s_0s_0\Box$   \\ \cline{2-6}
    \end{tabular}}\end{center}
    \end{table}

For face $2$ in Table~\ref{axHex1pt35}, recall Lemma~\ref{overlinelemma}, which shows that $$\overline {\widetilde{\alpha}_{C,B,A}}= \sigma(\underline{\widetilde{\alpha}_{C,B,A}  })=\sigma \alpha_{C,B,A}.$$    
\end{proof}
\begin{proposition}
[\textbf{SM} symmetric axiom] Let $A,B$ be objects of $\Sym(X).$ Then $$\gamma_{A,B}\bullet \gamma_{B,A}=\Id_{B\circ A}.$$
\end{proposition}
\begin{proof} By Lemma~\ref{bulletylemma}, it suffices to show $\sigma\chi(A,B)\bullet \gamma_{B,A}=\chi(B,A).$ We have $$d\Delta(\chi(A,B), \gamma_{B,A})=[\sigma \chi(A,B),\ \sigma \chi(A,B) \bullet \gamma_{B,A} ,\ \gamma_{B,A}, \Id_A]$$ so by the Matching Lemma it suffices to show $$[\sigma \chi(A,B),\ \chi(B,A) ,\ \gamma_{B,A}, \Id_A]$$ is commutative. The following Transposition Law table shows this commutativity.
\begin{table}[H] \begin{center}
 \scalebox{1}{       \begin{tabular}{ r | l || l | l | l | l | }
     \cline{2-6}

     &$s_1 (\sigma \chi(B,A))$     & $\Id_A$  &    $\sigma \chi(B,A)$      &$\sigma \chi(B,A)$            &$s_1 B=\sigma \Id_B$      \\ \cline{2-6}
     
    &  $\Delta_{-}(\sigma\chi(B,A))$     &$ \chi(A,B)$   &  $\sigma \chi(B,A)$                     &  $\gamma_{B,A}$        &     $\Id_B$        \\ \cline{2-6}  

   $\utimes$     &     & $\sigma \chi(A,B)$   &  $\chi(B,A)$      &   $\gamma_{B,A}$       &$\Id_A$   \\   \cline{2-6}  
     
    \end{tabular}}\end{center}
    \end{table}
\end{proof}

We have now shown $\Sym(X)$ is a small symmetric monoidal category. We now show that $\Sym$ is functorial. If $F:X\ra Y$ is a morphism of $2$-reduced symmetric quasimonoids, we must check that $\Bic(F)$ is a braided functor. We must check
\begin{proposition}
[Compatibility of the braiding with the distributors] Let $A, B$ be $0$-cells of $X$. Then $$\phi_{B,A}\bullet F(\gamma_{A,B})=\gamma_{F(A),F(B)}\bullet \phi_{A,B}.$$
\end{proposition}
\begin{proof}By Lemma~\ref{bulletylemma} it suffices to show $$F(\chi(B,A))\bullet F(\gamma_{A,B})=\sigma \chi(F(A),F(B))\bullet \phi_{A,B}.$$ Since $F$ is a morphism of symmetric quasimonoids, $F$ respects $\sigma$, so we have $F(\sigma\chi(A,B))=\sigma F(\chi(A,B)).$ Thus we will show $ F(\sigma\chi(A,B))=F(\chi(B,A))\bullet F(\gamma_{A,B})$ and $\sigma F(\chi(A,B))=\sigma \chi(F(A),F(B))\bullet \phi_{A,B}.$ 

For  $\sigma F(\chi(A,B))=F(\chi(B,A))\bullet F(\gamma_{A,B})$ we can apply the Matching Lemma to $$d(\Delta_{\bullet}(F(\chi(B,A)),\ F(\gamma_{A,B})))=[   F(\chi(B,A)),\ F(\chi(B,A))\bullet F(\gamma_{A,B}))           ,\   F(\gamma_{A,B}) ,\ Id_{F(A)}     ]$$ to see that it suffices to show that the following sphere is commutative $$[   F(\chi(B,A)),\ F( \sigma \chi(A,B))           ,\   F(\gamma_{A,B}) ,\ Id_{F(A)}=F(\Id_A)      ].$$ This sphere is filled by $F(\Delta_{-}(\sigma \chi(A,B))).$

To show  $\sigma F(\chi(A,B))=\sigma \chi(F(A),F(B))\bullet \phi_{A,B}$  we can apply the Matching Lemma to $d(\Delta_{\bullet}(\sigma \chi(F(A),F(B)),\ \phi_{A,B}))$ to see that it suffices to show that the following sphere is commutative $$[ \sigma \chi(F(A),F(B)),\ \sigma F(\chi(A,B))          ,\   \phi_{A,B},\ Id_{F(A)}=F(\Id_A)      ].$$ This sphere is filled by $\sigma(\Delta_{\stackrel{\wedge}{-}}(\sigma F(\chi(A,B)) )).$
\end{proof}

\section{Symmetric structure on the nerve of a symmetric monoidal category. }
Let $\CalC$ be a small symmetric monoidal groupoid. By considering the underlying monoidal category of $\CalC$ as a $(2,1)$-category with one object and taking the Duskin nerve, we get a $2$-reduced quasimonoid $N(\CalC).$ To give $N(\CalC)$ the structure of a $2$-reduced symmetric quasimonoid, we must define a involution $\sigma$ on $2$-cells of $N(\CalC)$. 

\begin{definition}
Recall a $2$-cell in $N(\CalC)$ is given by $(C,B,A \vbar  f)$ where $f:B\Rightarrow C\circ A.$ We define the symmetric structure of  $N(\CalC)$ by $$\sigma(C,B,A \vbar f):=(A,B,C \vbar \gamma_{C,A}\bullet f).$$
\end{definition}
The fact that $\sigma$ is an involution follows from the symmetric law $\gamma_{A,C}\bullet \gamma_{C,A} =\Id_{C\circ A}. $ Likewise the fact that $\sigma s_0 A=s_1 A$ follows immediately from the compatibility of the braiding with the unitors, $\gamma_{A, I} \bullet \rho_{A} =\lambda_{A}.$ To check that $\sigma$ makes $\CalC$ a symmetric $2$-reduced quasimonoid, we have left to check that the Reversal Law and Transposition Law hold.

The following Lemma gives identity for symmetric monoidal categories which is easily seen to follow formally from the axioms:
\begin{lemma}
\begin{equation}\label{hexlike1} \gamma_{A\circ B,C} \bullet \alpha_{A,B,C}\bullet \gamma_{B\circ C,A}=      (C \rhd \gamma_{B,A})    \bullet \alpha^{-1}_{C,B,A}\bullet (\gamma_{B,C} \lhd  A)
\end{equation}   
\end{lemma}

\begin{proposition}[Reversal Law for $N(\CalC)$]
If we have a $3$-cell $x_{0123}$ in $N(\CalC)$, meeting the $3$-cell condition:
\begin{equation}\label{3cellconditionagain}\alpha_{x_{23},x_{12},x_{01}}\bullet(x_{23}\rhd x_{012})\bullet x_{023}=(x_{123}\lhd x_{01})\bullet x_{013}\end{equation}
then the sphere $\sigma x_{0123}= [ \sigma x_{012} ,\ \sigma x_{013} ,\  \sigma x_{023},\  \sigma x_{123} ]$  is commutative, meaning that
\begin{equation*}\alpha_{x_{01},x_{12},x_{23}}\bullet(x_{01}\rhd   (\gamma_{x_{23},x_{12}} \bullet x_{123})     ) \bullet \gamma_{x_{13}, x_{01}}  \bullet x_{013}=((\gamma_{x_{12}, x_{01}}\bullet x_{012})\lhd x_{23})\bullet \gamma_{x_{23}, x_{02}} \bullet x_{023}.\end{equation*}
\end{proposition}
\begin{proof}
This statement is easily seen to be equivalent to its converse using the symmetric law, so we will work backwards, starting with the $3$-cell condition, Equation~\ref{3cellcondition}. (Alternatively, the steps below are reversible so the proof may be read backwards).

 First apply the naturality of $\gamma$ to both sides, yielding 
 \begin{equation*}\alpha_{x_{01},x_{12},x_{23}}\bullet \gamma_{x_{12}\circ x_{23}, x_{01}} \bullet(  (\gamma_{x_{23},x_{12}} \bullet x_{123})  \lhd  x_{01}  )   \bullet x_{013}=\gamma_{x_{23}, x_{01}\circ x_{12}} \bullet (x_{23}\rhd (\gamma_{x_{12}, x_{01}}\bullet x_{012}) ) \bullet x_{023}.\end{equation*}
 Then we apply interchange and move $\gamma_{x_{23}, x_{01}\circ x_{12}}$ to the left to get:
  \begin{equation*}\begin{split}
  &\gamma_{x_{01}\circ x_{12},x_{23}} \bullet \alpha_{x_{01},x_{12},x_{23}}\bullet \gamma_{x_{12}\circ x_{23}, x_{01}} \bullet  (\gamma_{x_{23},x_{12}} \lhd  x_{01}) \bullet (x_{123}  \lhd  x_{01})    \bullet x_{013} \\ &\quad=(x_{23}\rhd \gamma_{x_{12}, x_{01}})\bullet (x_{23}\rhd x_{012})  \bullet x_{023}\end{split}\end{equation*}
  The first three terms match the left hand side of Equation~\ref{hexlike1}. Making the substitution, we get
    \begin{equation*} \begin{split}
    &(x_{23} \rhd \gamma_{x_{12},x_{01}})    \bullet \alpha^{-1}_{x_{23},x_{12},x_{01} }\bullet (\gamma_{x_{12},x_{23} } \lhd  x_{01}) \bullet  (\gamma_{x_{23},x_{12}} \lhd  x_{01}) \bullet (x_{123}  \lhd  x_{01})    \bullet x_{013} \\ &\quad =(x_{23}\rhd \gamma_{x_{12}, x_{01}})\bullet (x_{23}\rhd x_{012})  \bullet x_{023} \end{split}
    \end{equation*}
  Finally use the symmetric law to make cancellations and move $\alpha$ to the other side, yielding 
      \begin{equation*}
          (x_{123}  \lhd  x_{01})    \bullet x_{013}= \alpha_{x_{23},x_{12},x_{01} }\bullet (x_{23}\rhd x_{012})  \bullet x_{023}.\end{equation*}
\end{proof}
\begin{proposition}[Transposition law for $N(\CalC)$]
Suppose we have the following two commutative spheres of $N(\CalC):$ 
$$[y_{012}  ,\ x_{023},\  \sigma z_{012}      ,\ \sigma x_{012}       ]$$
$$[x_{123},\ x_{023},\ x_{013},\ x_{012}]$$
The commutativity of these spheres means the following equations of morphisms in $\CalC$ hold:
\begin{align}
\label{firstsword3cell}\alpha_{x_{23},x_{01},x_{12}}\bullet(x_{23}\rhd (\gamma_{x_{12},x_{01}} \bullet x_{012}))\bullet x_{023}&=(y_{012}\lhd x_{12})\bullet \gamma_{x_{12},y_{02}} \bullet z_{012} \\
\label{secondsword3cell}\alpha_{x_{23},x_{12},x_{01}}\bullet(x_{23}\rhd x_{012})\bullet x_{023}&=(x_{123}\lhd x_{01})\bullet x_{013} 
\end{align}
    Then the sphere $$[\sigma x_{123},\ z_{012},\ x_{013},\ y_{012}]$$ is commutative, meaning that the following equation holds:
     \begin{equation}\label{thirdsword3cell}
       \alpha_{x_{12},x_{23},x_{01}}\bullet (y_{012}\lhd x_{12}) \bullet z_{012}=((\gamma_{x_{23},x_{12}}\bullet x_{123}) \lhd x_{01})\bullet x_{013}.
      \end{equation}

\end{proposition}
\begin{proof} 
Starting with Equation~\ref{firstsword3cell}, we apply interchange and the naturality of the braiding then move a $\gamma$-term to the left side, yielding
\begin{equation*}
\gamma_{x_{23}\circ x_{01},x_{12}} \bullet  \alpha_{x_{23},x_{01},x_{12}} \bullet(x_{23}\rhd \gamma_{x_{12},x_{01}}) \bullet (x_{23}\rhd x_{012})\bullet x_{023}=(y_{012}\lhd x_{12}) \bullet z_{012} 
\end{equation*}
We can make a substitution to the first three terms on the left hand side, using the first hexagon identity, yielding:
 \begin{equation*}
 \alpha^{-1}_{x_{12},x_{23},x_{01}} \bullet (\gamma_{x_{23},x_{12}}\lhd x_{01})\bullet\alpha_{x_{23},x_{12},x_{01}} \bullet (x_{23}\rhd x_{012})\bullet x_{023}=(y_{012}\lhd x_{12}) \bullet z_{012} 
 \end{equation*}
 The last $3$ terms on the left hand side match the left hand side of Equation~\ref{secondsword3cell}. We make the substitution, yielding 
   \begin{equation*}
   \alpha^{-1}_{x_{12},x_{23},x_{01}} \bullet (\gamma_{x_{23},x_{12}}\lhd x_{01})\bullet(x_{123}\lhd x_{01})\bullet x_{013}=(y_{012}\lhd x_{12}) \bullet z_{012}.
   \end{equation*}
   A few simplifications yield Equation~\ref{thirdsword3cell}, finishing the proof.
\end{proof}
To check that $N$ is functorial, we must show that if $F:\CalC \ra \CalD$ is a strictly identity-preserving symmetric monoidal functor of small symmetric monoidal groupoids, then $N(F)$ preserves $\sigma$. Recall that $$N(F)(C,B,A\vbar f)= (F(C),F(B),F(A) \vbar\phi_{C,A}\bullet F(f)).$$ So the functoriality of $N$ is equivalent to showing for all $2$-cells $(C,B,A \vbar f)$ of $N(\CalC)$ that
$$\phi_{A,C}\bullet F(\gamma_{C,A}\bullet f)=\gamma_{F(C),F(A)}\bullet \phi_{C,A}\bullet F(f).$$ This is follows immediately from \textbf{BMGFun}.
\begin{theorem}\label{NSymequiv}
$N$ and $\Sym$ are inverse equivalences of categories between the category of small symmetric monoidal groupoids and strictly identity-preserving braided functors and the category of  $2$-reduced algebraic symmetric quasimonoids. Furthermore, $N$ and $\Sym$ preserve strictness, and the natural isomorphisms $u':N \Sym\cong \Id $ and $U': \Sym N \cong \Id$ exhibiting the equivalence are strict, thus $N$ and $\Sym$ are also inverse equivalences of categories between the category of $2$-reduced algebraic quasimonoids and strict morphisms and the category of small symmetric monoidal groupoids categories and strict braided functors.
\end{theorem}
\begin{proof}  Let $u':N \Sym\cong \Id $ and $U': \Sym N \cong \Id$ be defined on the underlying simplicial sets and small $(2,1)$-categories from the isomorphisms $u:N \Bic\cong \Id $ and $U: \Bic N \cong \Id.$ We must only show the isomorphisms are actually symmetric in the sense that $u'$ preserves $\sigma$ and $U': \Bic N \cong \Id$ satisfies \textbf{BMGFun} in order to promote these isomorphisms to isomorphism. Recall that $u(C,B,A \vbar f)= \underline{f}.$ So we must check that $\underline{\sigma f}=\gamma_{C,A}\bullet \underline{f}.$ By Lemma~\ref{bulletylemma}, it is enough to show $\sigma f=\sigma{\chi(C,A)}\bullet \underline{f}.$ Applying the matching lemma to $d\Delta_\bullet(\sigma \chi(C,A), \underline{f})$, we see that it suffice to show the sphere $$[ \sigma\chi(C,A) ,\ \sigma f ,\   \underline{f}, \Id_C   ]$$ is commutative. This sphere is filled by $\sigma \Delta_{\stackrel{\wedge}{-}}(f).$

For $U$ is defined for a morphism $f:A\ra B$ by  $U(f)=(B,A,\id_{\Box}\vbar \rho_B \bullet f).$ The symmetric axiom \textbf{BMGFun} for $U$ works out to
$$\underline{(B, A \circ B, A \vbar \gamma_{A,B}\bullet \Id_{A\circ B})}=(B\circ A,A\circ B,\id_{\Box}\vbar \rho_{B\circ A}\bullet\gamma_{A,B}),$$
To show this identity, by the Matching Lemma applied to $\Delta_{-}((B, A \circ B, A \vbar \gamma_{A,B})    )$, it suffices to show the following sphere is commutative in $N(\CalC)$:
$$[(B,B\circ A,A \vbar \Id_{B\circ A}),\, (B, A \circ B, A \vbar \gamma_{A,B}),\, (B\circ A,A\circ B ,\id_{\Box}\vbar \rho_{B\circ A}\bullet\gamma_{A,B}) ,\, (A,A,\id_{\Box} \vbar \rho_A)]$$
This commutativity of this sphere is equivalent to the following $3$-cell condition:
$$\alpha_{ B, A   ,\id_{\Box}} \bullet (B\rhd \rho_A)\bullet \gamma_{A,B}= (\Id_{B\circ A}\lhd \id_{\Box})\bullet \rho_{B\circ A}\bullet\gamma_{A,B}$$
This identity can be proven by applying the compatibility of $\alpha$ and $\rho$ to the left hand side.

Algebraic structure and strictness depend only on the underlying simplicial structure, so the second statement is immediate.
\end{proof}

\begin{theorem} \label{gammasummary} $\Orb \circ N$ and $\Sym\circ \Fl_2$ are inverse equivalences of categories between the category of small symmetric monoidal groupoids and strictly identity-preserving braided functors and the category of  $2$-reduced algebraic inner-Kan $\Gamma$-sets. Also, they are inverse equivalences of inverse equivalences of categories between the category of small symmetric monoidal groupoids and strict braided functors and the category of  $2$-reduced algebraic inner-Kan $\Gamma$-sets and strict morphisms.
\end{theorem}
\begin{proof}This follows from Theorem~\ref{OrbFlequiv} and Theorem~\ref{NSymequiv}. \end{proof}

We will denote the map $\Orb \circ N$ by $N_\Gamma$, which we view as a nerve for symmetric monoidal groupoids which is valued in $\Gamma$-sets. 
\begin{remark}
The functor which forgets algebraic structures on $\Gamma$-sets is an equivalence of categories, thus $N_\Gamma$ is an equivalence of categories from symmetric monoidal groupoids to (non-algebraic) $2$-reduced inner-Kan $\Gamma$-sets.
\end{remark}

\section{The $\Gamma$ nerve}

The following explicit description of $N_\Gamma(\C)$ can be inferred from our definitions:
\begin{itemize}
\item $N_\Gamma(\C)(\underline{0})$ has a unique object, denoted $\Box.$
\item A $\underline{1}$-cell of $N_\Gamma(\C)$ is an object of $\C$. 
\item A $\underline{2}$-cell of $N_\Gamma(\C)$ is a quadruple $(x_1, x_{(12)},x_2 \vbar x_{12})$ where $x_1,$  $x_{(12)},$ and $x_2$ are objects of $\C$ and $x_{12}: x_{(12)}\ra x_2 \circ x_1$ is a morphism.
\item A $\underline{3}$-cell fo $N_\Gamma(\C)$ is a commuting diagram:
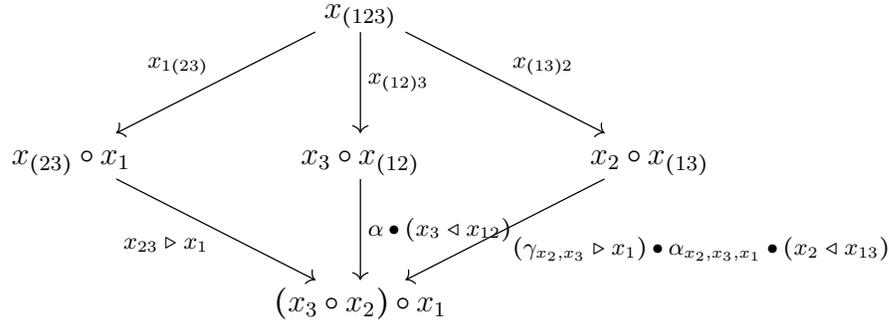
\begin{center}
\begin{tikzcd}[column sep=large, row sep=large]
 ~ & x_{(123)}   \arrow[swap]{dl}{x_{1(23)}} \arrow{d}{x_{(12)3}} \arrow{dr}{x_{(13)2}}     &    \\
x_{(23)}\circ x_{1}  \arrow[swap]{dr}{x_{23} \rhd x_1}  & x_{3} \circ x_{(12)}    \arrow{d}{\alpha\bullet (x_3\lhd x_{12})}    & x_{2} \circ x_{(13)} \arrow{dl}{( \gamma_{x_2,x_3}\rhd x_1)\bullet \alpha_{x_2,x_3,x_1}\bullet (x_2 \lhd x_{13})}  \\
 & (x_3\circ x_2)\circ x_1 &   
\end{tikzcd}
\end{center}
\item A $N_\Gamma(\C)$ is $3$-coskeletal.\end{itemize}

\begin{remark} \label{mandellrem}
A direct construction of $N_\Gamma$ can be obtained from a construction given by Mandell in \cite{Man10}. Mandell shows how a small symmetric monoidal category $\C$ gives rise to a functor $\mathcal{K}(\C):\Gamma^\op\ra \cat{Cat}$. This construction is used by Mandell to construct a $K$-theory functor associating a symmetric monoidal category to a special $\Gamma$-space. Taking objects levelwise in $\mathcal{K}(\C)$ yields a $\Gamma$-set which is easily seen to be isomorphic to $N_\Gamma(\C)$ since it is $3$-coskeletal and matches $N_\Gamma(\C)$ up to $3$-cells. 
\end{remark}
\begin{remark} \label{extendrem}
Both our explicit description of $N_\Gamma(\C)$ above and the construction of $N_\Gamma(\C)$ in Remark~\ref{mandellrem} allow $\C$ to be an arbitrary symmetric monoidal category.
\end{remark}

\subsection{Kan $\Gamma$-sets}
\begin{definition} \label{kandef}
Recall that $\Lambda^{\underline{2}}_{\{\hat{k}_0\}}$ and $\Lambda^{\underline{2}}_{\{\hat{k}_1\}}$ denote the universal horns obtained from removing the first and last faces from $\Gamma[2]$, i.e. the coface maps $\hat{k}_0$ and $\hat{k}_1$. An inner-Kan $\Gamma$-set $X$ will be called $\emph{Kan}$ if every horn of type $\Lambda^{\underline{2}}_{\{\hat{k}_0\}}$ and $\Lambda^{\underline{2}}_{\{\hat{k}_1\}}$ has a filler in $X$.
\end{definition}
It may seem strange that we have defined the notion of Kan $\Gamma$-sets to only require an extra filling condition for ``outer horns'' of dimension $2$, whereas a Kan simplicial set is required to have fillers for outer horns in every dimension. However, the following theorem, which is an easy corollary of Theorem~1.3 in \cite{Joy02a}, shows that these extra filling conditions are redundant in the case of simplicial sets:
\begin{theorem}[Joyal]
If a simplicial set is inner-Kan and has fillers for the outer horns $\Lambda^2_0$ and $\Lambda^2_2$, then it is Kan.
\end{theorem}
\begin{definition}
A symmetric monoidal groupoid $\C$ is called \emph{grouplike} or a \emph{ Picard groupoid} if for all objects $c\in \C$ there is a $c^{-1}$ such that $c \circ c^{-1}$ (and therefore $c^{-1}\circ c$) are isomorphic to the identity object $\id_\Box$.
\end{definition}
\begin{proposition}\label{Kanprop}
The symmetric monoidal groupoid $\C$ is grouplike if and only if $N_{\Gamma}(\C)=\Orb\circ \Gamma (\C)$ is Kan.
\end{proposition}
\begin{proof}
A filler of the horn $[-,\id_{\Box}, c]$ in $N_{\Gamma}(\C)$ provides a left inverse to $c$, which serves as well as a right inverse since $\C$ is symmetric, showing $\C$ is grouplike if $N_{\Gamma}(\C)$ is Kan. 

In the other direction, a $\Lambda^2_{\{\hat{k}_0\}}$-horn  $N_{\Gamma}(\C)$ of the form $[-,c,b]$ may be filled by the cell $(c\circ b^{-1}, c, b \vbar c \lhd \iota_b$ where $\iota_b: b^{-1}\circ b \ra \id_{\Box}$ is an isomorphism exhibiting $b^{-1}$ as the left inverse of $b$. A similar argument works for $\Lambda^{\underline{2}}_{\{\hat{k}_1\}}$-horns, showing that $N_{\Gamma}(\C)$ is Kan if $\C$ is grouplike.
\end{proof}

\subsection{Summary and extensions}
We have a commuting diagram:

\begin{center}
\begin{tikzcd}
\mbox{Sym. Mon. Groupoids} \arrow[xshift=-0.7ex,swap]{d}{N_\Gamma=\Orb \circ N}     \arrow[leftrightarrow]{r}{=} & \mbox{Sym. Mon. Groupoids} \arrow[xshift=-0.7ex,swap]{d}{N} \arrow{r}{u}    &  \mbox{$(2,1)$-categories} \arrow[xshift=-0.7ex,swap]{d}{N}   \\
\mbox{$2$-red. inner-Kan $\Gamma$-sets} \arrow[xshift=0.7ex,swap]{u}{\Sym \circ h_2\circ \phi^*}  \arrow[yshift=0.7ex]{r}{h_2\circ \phi^*}      & \mbox{$2$-red. sym. quasimonoids} \arrow[xshift=0.7ex,swap]{u}{\Sym} \arrow{r}{u'}  \arrow[yshift=-0.7ex]{l}{\Orb} & \mbox{$2$-red. inner-Kan $\Delta$-sets} \arrow[xshift=0.7ex,swap]{u}{\Bic} 
\end{tikzcd}
\end{center}
where $u$ is the functor that takes the underlying monoidal category of a symmetric monoidal category, viewed as a $(2,1)$-category with one object, and $u'$ forgets the symmetrizing involution $\sigma$ of a symmetric quasimonoid.

In the Kan case, we have:

\begin{center}
\begin{tikzcd}
\mbox{Gr. Sym. Mon. Groupoids} \arrow[xshift=-0.7ex,swap]{d}{N_\Gamma=\Orb \circ N}     \arrow[leftrightarrow]{r}{=} & \mbox{Gr. Sym. Mon. Groupoids} \arrow[xshift=-0.7ex,swap]{d}{N} \arrow{r}{u}    &  \mbox{$(2,0)$-categories} \arrow[xshift=-0.7ex,swap]{d}{N}   \\
\mbox{$2$-red. Kan $\Gamma$-sets} \arrow[xshift=0.7ex,swap]{u}{\Sym \circ h_2\circ \phi^*}  \arrow[yshift=0.7ex]{r}{h_2\circ \phi^*}      & \mbox{$2$-red. Kan sym. quasimonoids} \arrow[xshift=0.7ex,swap]{u}{\Sym} \arrow{r}{u'} \arrow[yshift=-0.7ex]{l}{\Orb} & \mbox{$2$-red. Kan $\Delta$-sets} \arrow[xshift=0.7ex,swap]{u}{\Bic} 
\end{tikzcd}
\end{center}

Finally, our discussion of $N_\Gamma$ would not be complete without some further discussion of the result of applying $N_\Gamma$ to an arbitrary symmetric monoidal category.

\begin{definition}
There is an endofunctor $\mbox{pr}_2$ of $\Delta$-sets that makes a simplicial set $2$-subcoskeletal and $3$-coskeletal \footnote{See Definition~\ref{coskeletaldef}}. It can be defined by letting $\mbox{pr}_2(X)(n)$ be equivalence classes of maps $f:\Delta[n]\ra X$ under the equivalence relation $f \sim f'$ if $\tr_2(f)= \tr_2(f')$ for $n\leq 3$, and defining higher cells by letting $\mbox{pr}_2(X)$ be $3$-coskeletal.
\end{definition}
\begin{proposition}
If $X$ is a $2$-reduced inner-Kan $\Gamma$-set, then $h_2\phi^*(X)=\mbox{pr}_2\phi^*(X)$
\end{proposition}
\begin{proof}
We showed in the proof of Theorem~\ref{OrbFlequiv} that if $X$ is a $2$-reduced inner-Kan $\Gamma$-set then the relation $\sim$ of homotopy rel. boundary for $2$ cells in $\phi^{\star}(X)$, which is used in Section~\ref{homotopysection} to construct $h_2,$ is the identity relation. The claim then follows directly from the construction of $h_2$.
\end{proof}

\begin{proposition}The following diagram is commutative up to isomorphism:
\begin{center}
\begin{tikzcd}
\mbox{Sym. Mon. Categories} \arrow{d}{N_\Gamma}     \arrow{r}{u} & \mbox{Bicategories} \arrow{d}{N}   \\
\mbox{$\Gamma$-sets}\arrow{r}{\mbox{pr}_2 \circ \phi^*}   &  \mbox{$\Delta$-sets}
\end{tikzcd}
\end{center}
\end{proposition}
\begin{proof}
The isomorphism between the two compositions is trivial to see except for $3$-cells, where the crucial issue is to check whether a commuting $3$-sphere in the Duskin nerve of $u(\C)$ can be lifted to a commuting $3$-sphere in $N_\Gamma(\C).$ In other words, we must show that a commuting diagram
\begin{center}
\begin{tikzcd}[column sep=large, row sep=large]
 ~ & x_{(123)}   \arrow[swap]{dl}{x_{1(23)}} \arrow{d}{x_{(12)3}}     &    \\
x_{(23)}\circ x_{1}  \arrow[swap]{dr}{x_{23} \rhd x_1}  & x_{3} \circ x_{(12)}    \arrow{d}{\alpha\bullet (x_3\lhd x_{12})}    \\
 & (x_3\circ x_2)\circ x_1 &   
\end{tikzcd}
\end{center}
giving a $3$-cell in the Duskin nerve can be completed to a $3$-cell in $N_\Gamma(\C)$:
\begin{center}
\begin{tikzcd}[column sep=large, row sep=large]
 ~ & x_{(123)}   \arrow[swap]{dl}{x_{1(23)}} \arrow{d}{x_{(12)3}} \arrow{dr}{x_{(13)2}}     &    \\
x_{(23)}\circ x_{1}  \arrow[swap]{dr}{x_{23} \rhd x_1}  & x_{3} \circ x_{(12)}    \arrow{d}{\alpha\bullet (x_3\lhd x_{12})}    & x_{2} \circ x_{(13)} \arrow{dl}{( \gamma_{x_2,x_3}\rhd x_1)\bullet \alpha_{x_2,x_3,x_1}\bullet (x_2 \lhd x_{13})}  \\
 & (x_3\circ x_2)\circ x_1 &   
\end{tikzcd}
\end{center}
This is always possible for instance letting $x_{(13)}=x_1\circ x_3$ and the map $x_{(13)}$ be the identity. Then the map at the lower right of the above diagram is invertible and the map $$x_{(13)2}:x_{(13)2}\ra x_{(13)}\circ x_2$$ may be filled in appropriately so that the diagram commutes.
\end{proof}

\begin{definition}
A \emph{stratified} simplicial set is simplicial set $X$ together with a subset $tX$ of the cells of $X$, called the \emph{thin} cells, such that $tX$ contains all degenerate cells and no $0$-cells of $X$.

Following an idea of Street \cite{Str87}, Verity defines in \cite{Ver} a horn-extension condition for stratified simplicial sets that, and calls stratified simplicial sets satisfying the condition \emph{weak complicial sets}. These weak complicial sets provide a model for weak $\infty$-categories.
\end{definition}

In \cite{Gur09}, Gurski uses the notion of weak complicial sets to contextualize Duskin's characterization in \cite{Dus02} of the nerve of an arbitrary bicategory. In particular, if $\C$ is a bicategory and $N(\C)$ is made into a stratified simplicial set such that all $3$-cells and higher in $N(\C)$ are marked as thin, along with those $1$ and $2$-cells which are associated respectively to invertible $1$ and $2$-morphisms of $\C,$ then $N(\C)$ a weak complicial set. This leads us to conjecture:

\begin{conjecture}
There is an appropriate notion of a weak complicial $\Gamma$-set, which models symmetric monoidal $\infty$-categories, and has the property that $N_\Gamma(\C)$ has the structure of a weak complicial $\Gamma$-set for all symmetric monoidal categories $\C$.
\end{conjecture}

\chapter{The globular nerve of a fancy bicategory \label{thetachapter}}
\section{The category $\Theta_2$}
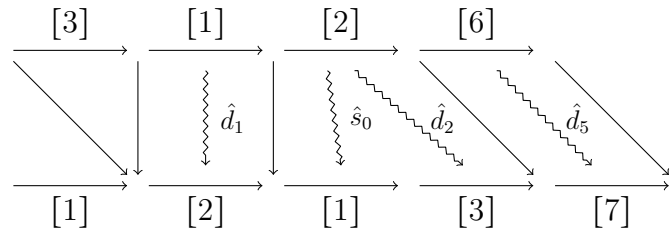
\begin{figure}
\begin{center}
\begin{tikzpicture}[scale=1.8,auto]

\begin{scope}

\node (10) at (1,1) {};
\node (00) at (0,1) {};
\node (11) at (1,0) {};
\node (01) at (0,0) {};
\node[rotate=-45] at (.5,.5){};
\node[scale=.8] at (.6,.65){};

\node (20) at (2,1) {};
\node (21) at (2,0) {};
\node (30) at (3,1) {};
\node (31) at (3,0) {};
\node (40) at (4,1) {};
\node (41) at (4,0) {};
\node (51) at (5,0) {};
\node[rotate=-45] at (1.5,.5){};
\node[scale=.8] at (1.6,.65){};

\path[->] (00) edge node[midway]{$[3]$}(10);
\path[->] (01) edge node[midway,swap]{$[1]$}(11);
\path[->] (10) edge (11);
\path[->] (10) edge node[midway]{$[1]$}(20);
\path[->] (11) edge node[midway,swap]{$[2]$}(21);
\path[->] (20) edge (21);
\path[->] (00) edge (11);
\path[->] (40) edge (51);
\path[->] (30) edge (41);
\path[->] (20) edge node[midway]{$[2]$}(30);
\path[->] (21) edge node[midway,swap]{$[1]$}(31);
\path[->] (30) edge node[midway]{$[6]$}(40);
\path[->] (31) edge node[midway,swap]{$[3]$}(41);
\path[->] (41) edge node[midway,swap]{$[7]$}(51);
\draw [->,
line join=round,
decorate, decoration={
    zigzag,
    segment length=4,
    amplitude=.9,post=lineto,
    post length=2pt
}]  (1.5,.85)  -- (1.5,.15);
\draw [->,
line join=round,
decorate, decoration={
    zigzag,
    segment length=4,
    amplitude=.9,post=lineto,
    post length=2pt
}]  (2.4,.85)  -- (2.5,.15);
\draw [->,
line join=round,
decorate, decoration={
    zigzag,
    segment length=4,
    amplitude=.9,post=lineto,
    post length=2pt
}]  (2.6,.85)  -- (3.4,.15);
\draw [->,
line join=round,
decorate, decoration={
    zigzag,
    segment length=4,
    amplitude=.9,post=lineto,
    post length=2pt
}]  (3.65,.85)  -- (4.35,.15);
\node[scale=.85] at (1.7,.5){$\hat{d}_1$};
\node[scale=.85] at (2.65,.5){$\hat{s}_0$};
\node[scale=.85] at (3.25,.5){$\hat{d}_2$};
\node[scale=.85] at (4.25,.5){$\hat{d}_5$};
\end{scope}

\end{tikzpicture}\caption{A morphism in $\Theta_2=\Delta\wreath\Delta$ \label{thetamorphismexample}}\end{center}\end{figure}

The category $\Theta_n$ was introduce by Joyal in \cite{Joy97} in an attempt to give a definition for $n$-categories. Rezk uses $\Theta_n$ in a different way in \cite{Rez10} to give a Cartesian model category for $(\infty,n)$ categories. We will work with a combinatorial definition of $\Theta_n$ which is due to Berger \cite{Ber07}.
\begin{definition} \label{wreathdef}
Let $\CalC$ be a small category. $\Gamma\wreath \CalC$ is a category obtained from $\CalC$, with
\begin{itemize}
\item An object of $\Gamma\wreath \CalC$ is a tuple $(\underline{n}, c_1,\ldots, c_{n})$ where $\underline{n}$ is an element of the Segal $\Gamma$ category as defined in Chapter~\ref{gammachapter}, and $c_1,\ldots, c_{n}$ are objects of $c$ 
\item A morphism $f:(\underline{n}, c_1,\ldots, c_{n})\ra (\underline{m}, c'_1,\ldots, c'_{m})$ consists of a map $\Ty(f):\underline{n}\ra\underline{m}$ in $\Gamma,$ called the \emph{type} of $f$, together with a morphism of $f_{ji}$ in $\CalC$ from  $c_i$ to $c'_j$ whenever $j\in f(i).$ These morphisms of $\CalC$  called the \emph{components} of $f$.
\item $g\circ f$ is defined by $\Ty(g\circ f)=\Ty(g)\circ\Ty(f).$ If $k$ is in $(g\circ f)(i)$ then there is a unique $j\in f(i)$ with $k\in g(j)$. We let $(g\circ f)_{ki}=g_{kj}\circ f_{ji}.$
\end{itemize}
\end{definition}
It is not (as the notation $\wreath$ perhaps suggests) possible to form a wreath product of two categories in a reasonable way. A slight generalization is given below, however:
\begin{definition} \label{simpwreathdef}
Suppose we have a category $\mathcal{B}$ and a functor $\varphi:\mathcal{B}\ra \Gamma.$ Then we define a category $\mathcal{B}\wreath \CalC$ as follows
\begin{itemize}
\item An object of $\mathcal{B}\wreath \CalC$ is a tuple $(b, c_1,\ldots, c_{n})$ where $b$ is an object of $\mathcal{B}$ with $\varphi(b)=\underline{n}$,  and $c_1,\ldots, c_{n}$ are objects of $c$ 
\item A morphism $f:(b, c_1,\ldots, c_{n})\ra (b', c'_1,\ldots, c'_{m})$ consists of a map $\Ty(f):b\ra b'$ in $\mathcal{B},$ called the \emph{type} of $f$, together with a morphism of $f_{ji}$ in $\CalC$ from  $c_i$ to $c'_j$ whenever $j\in \varphi(f)(i).$ As before, these morphisms of $\CalC$  called the \emph{components} of $f$.
\item  $\Ty(g\circ f)=\Ty(g)\circ\Ty(f).$ If $k$ is in $\varphi(g\circ f)(i)$ then there is a $j\in \varphi(f)(i)$ with $k\in \varphi(g)(j)$. Let $(g\circ f)_{ki}=g_{kj}\circ f_{ji}.$
\end{itemize}
\end{definition}
In particular, we have a canonical map $\phi:\Delta\ra\Gamma$ described in Definition~\ref{phidefinition}, allowing us to define the \emph{simplicial wreath product} $\Delta\wreath \CalC.$ We define $\Theta_1:=\Delta$ and $\Theta_n:=\Delta\wreath \Theta_{n-1}$. We will be chiefly concerned with $\Theta_2= \Delta\wreath \Delta$. An object of $\Theta_2$ is an element $[n]$ of $\Delta$ together objects $[c_1],[c_2],\ldots, [c_{n}]$ of $\Delta$. We denote such an object by $\langle c_1,c_2,\ldots,c_{n}\rangle.$ In Figure~\ref{thetamorphismexample} we give an example of a morphism of $\Theta_2$ from $\langle 3, 1, 2, 6\rangle $ to $\langle 1, 2, 1, 3, 7\rangle$.

To further clarify the nature of the category $\Theta_2$, there is an alternate way of defining $\Theta_2$ due to Berger in \cite{Ber02} and independently to Makkai and Zawadowski in \cite{MZ01} and which we consider for instance the object $\langle 1, 0, 2 \rangle$ to be the free strict bicategory on the diagram below:
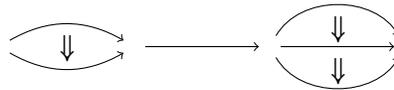
\begin{center}
\begin{tikzpicture}[scale=1.8,auto]

\begin{scope}

\node (10) at (1,1) {};
\node (00) at (0,1) {};

\node (20) at (2,1) {};
\node (30) at (3,1){};

\path[->] (10) edge node[midway]{}(20);

\draw[->] (20) to [out=-60,in=-120] node(h)[midway]{} (30) ;
\draw[->] (20) to [out=60,in=120] node(i')[midway]{} (30) ;
\path[->] (20) edge node(i)[midway]{}(30);
\node[rotate=-90] (aretb) at ($(i)+(0,-.26)$){$\Rightarrow$};
\node[rotate=-90] (aretc) at ($(i')+(0,-.26)$){$\Rightarrow$};

\draw[->] (00) to [out=-30,in=-150] node(g)[midway]{} (10) ;
\draw[->] (00) to [out=30,in=150] node(f)[midway]{} (10) ;

\node[rotate=-90] (areta) at ($(f)+(0,-.26)$){$\Rightarrow$};
\end{scope}
\end{tikzpicture}\end{center}
Then $\Theta_2$ can be viewed as the category of such free strict bicategories, with morphisms given by strict bicategory functors. Dual to this, we can think of $\Theta_2^{\op}$ as the category of ``pasting diagrams'', with morphisms given by ways of composing part of the diagram to get a new diagram. We will need only the combinatorial definition of $\Theta_2$ given above, however.

Recall that we can denote a morphism $g$ in $\Delta$ by $[g(0),\ldots,g(n)]$ or $g(0)g(1)\ldots g(n)$. For instance, the coface map $\hat{d}_1: [3]\ra[4]$ is denoted $0234$. We will call this the \emph{standard} notation for $g$. We define an alternate ``stars and bars'' notation for morphisms in $\Delta$ whereby $g=2335: [4]\ra [6]$ is denoted by $$\star\star|\star||\star \star|\star.$$ The general rule is that $g(i)=j$ if there are $j$ stars to the left of the $i$th bar (counting from $0$).

In this notation, $j\in \phi(f)(i)$ if and only if the $j$th star (counting from $0$) is between the $i$th and $(i+1)$st bar in this representation of $f$. To denote a morphism in $\Delta \wreath \Delta$ of type $f$ we start with the stars-and-bars notation for $f$ and replace the $j$th star by a notation for the component $f_{ji}$, as long as this star lies between two bars. The stars before the first bar and after the last bar are thus left in place.

This notation tells us the type of $f$ and its components at the same time. Commas between consecutive components and a $\cdot$ symbol between consecutive bars are used to make the notation easier to read. We give our notation for the morphism given in Figure~\ref{thetamorphismexample} as an example:
$$    \star|\cdot|01|001,013|0123467|  $$

\begin{definition} \label{innerdef}
We will call a map $f: [m]\ra [n]$ of $\Delta$ \emph{inner} if $0$ and $n$ are in the image of $f$. Note the faces that can be removed from a simplex $\Delta[n]$ to make an inner horn correspond to the inner coface maps of $\Delta$. We call a map of $\Theta_2$ \emph{inner} if its type is inner and each of its components are inner. By $\Lambda_{\{\mathfrak{f}\}}^{c}$ we denote the \emph{universal inner horn} formed by removing an inner coface map $\mathfrak{f}$ from the representable presheaf $\Theta_2[c].$ A \emph{inner horn} in a $\Theta_2$-set $X$ is a map from a universal inner horn to $X$. $X$ is called \emph{inner-Kan} if every inner horn $H$ in $X$ has a \emph{filler}, i.e. an extension of $H$ along the canonical inclusion of the universal inner horn into a representable presheaf.
\end{definition}
\begin{definition}
An \emph{algebraic} inner-Kan $\Theta_2$-set is an inner-Kan $\Theta_2$ set $X$ together with set $\chi$ of preferred filler for every inner horn in $X$. A map of algebraic inner-Kan $\Theta_2$-sets is called \emph{strict} if it preserves these preferred fillers in the obvious sense.
\end{definition}
\begin{definition}
An inner-Kan $\Theta_2$-set $X$ is called \emph{$2$-reduced} if every inner horn of minimal complementary dimension \footnote{See Definition~\ref{mincompdimdef}} $2$ or greater in $X$ has a unique filler.
\end{definition}

\section{Coface maps in $\Theta_2$ \label{tedioussec}}
\begin{definition} \label{thetacofacedef} The \emph{dimension} of an object $\langle c_1,\ldots, c_{n}\rangle$ in $\Theta_2$ is $n+\sum_{1\leq i \leq n} c_i $. A \emph{coface} map in $\Theta_2$ is a monic map which increases dimension by $1$. 
\end{definition}

\begin{definition}
A \emph{Reedy category} is a small dimensional category $\CalC$ together with two \emph{wide} (meaning that they include every object) subcategories $\CalC^{+}$ and $\CalC^{-}$ such that:
\begin{itemize}
\item Every non-identity morphism in $\CalC^+$ raises dimension and every non-identity morphism $\CalC^-$ lowers dimension
\item Every morphism in $f$ in $\CalC$ factors uniquely as $f=f^{+} f^-$ where $f^+$ is in $\CalC^+$ and $f^-$ is in $\CalC^-$.
\end{itemize}
\end{definition}
$\Theta_n$ was shown to be a Reedy category Berger in \cite{Ber02}, with $\Theta_2^-$ being the subcategory of epimorphisms and $\Theta_2^+$ being the subcategory of monomorphisms. See \cite{BR11} for a different proof using the combinatorial definition of $\Theta_2$ used above.
\begin{definition}
If we have a morphism $$f: \langle c_1, \ldots, c_{n}\rangle \ra \langle c'_1, \ldots, c'_{m}\rangle,$$ the \emph{$i$th part} $f_{\bullet i}$ of $f$ is the map $$|f_{j i}, f_{(j+1)i}, \ldots, f_{ki}|:\langle c_i\rangle \ra \langle c'_j, \ldots, c'_{k}\rangle,$$ where $f_{ji},\ldots f_{ki}$ are of all the components of the form $f_{-\, i}.$ 
\end{definition}
We will sometimes denote a morphism letting the names of its parts stand in for lists of components, denoting a map $f$ by $$\star \ldots \star |f_{\bullet 1}|f_{\bullet 2}|\ldots |f_{\bullet n}|\star \ldots \star.$$

\begin{proposition}\label{thetafactorprop}
For every morphism $f:a\ra c$ in $\Theta_2$ with $\dim a < \dim c$ which is not a coface map, there distinct coface map $d'$ and $d''$ and and morphisms $f'$ and $f''$ such that $f=d'\circ f'=d''\circ f''.$ 
\end{proposition}
\begin{proof} Taking the epi-monic factorization $f=f^+ f^-$, note that $f^+$ increases dimension by the Reedy axioms. Furthermore, if $f^-$ is non-identity, then it decreases dimension, thus $f^+$ increases dimension by $2$ or more. On the other hand, if $f^-$ is the identity then $f=f^+$, and the fact that $f$ is not a coface map ensures that $f^+$ increases dimension by $2$ or more. Thus it suffices to prove that any monic which increases dimension $2$ or more factors into two non-identity monic maps in two ways, $f=hg$ and $f=h'g'$ with $h \neq h'$. If we know this, we can easily see by induction that any monic map can in fact be entirely factored into a series of coface maps in at least two ways such that the last coface map in the factorization is distinct between the two factorizations.

Let $f$ be a non-identity monomorphism  $$f: \langle c_1, \ldots, c_{n}\rangle \ra \langle c'_1, \ldots, c'_{m}\rangle$$ which is not a coface map. Since $f$ increases dimension, and does not increase it by $1$, it increases dimension by at least $2.$

First suppose the type of $f$ is not inner. Without loss of generality suppose $1$ is not in the image of the type of $f$, so that  $f$ has no $f_{1j}$ component for any $j$. Then if $$\underbrace{\star \ldots \star}_\text{$k$ times}  |f_{\bullet 1}|f_{\bullet 2}|\ldots |f_{\bullet n}|\underbrace{\star \ldots \star}_\text{$l$ times}.$$
we know $k\geq 1$. We can factor $f=hg$ and $f=h'g'$:
\begin{align*}
g&= \star|\id|\id|\ldots|\id|: \langle c_1, \ldots, c_{n}\rangle \ra \langle 0,c_1, \ldots, c_{n}\rangle \\
h&=|\underbrace{0,0,\ldots ,0}_\text{$k$ times}|f_{\bullet 1}|f_{\bullet 2}|\ldots |f_{\bullet n}|\underbrace{\star \ldots \star}_\text{$l$ times}:\langle 0,c_1, \ldots, c_{n}\rangle\ra \langle c'_1, \ldots, c'_{m}\rangle.\\
g'&=\underbrace{\star \ldots \star}_\text{$k-1$ times}  |f_{\bullet 1}|f_{\bullet 2}|\ldots |f_{\bullet n}|\underbrace{\star \ldots \star}_\text{$l$ times}: \langle c_1, \ldots, c_{n}\rangle \ra \langle c'_1, \ldots,c'_{m}\rangle \\
h'&= \star  |\id|\id|\ldots |\id|:\langle c'_1, \ldots, c'_{m}\rangle\ra \langle c'_1, \ldots, c'_{m}\rangle.
\end{align*}
Since $\diff(f)>1$, one of the following must clearly hold:
\begin{itemize}
\item $k>1$
\item $l>0$
\item Some $f_{\bullet i}$ is not an identity map
\item $c'_1> 0$
\end{itemize}
If any of these hold, $h$ is not the identity map, and if any of the first three hold, $g'$ is not the identity map. In the case $c'_1>0$, however, $g'$ may be the identity map. In this case we revise our definition, defining 
\begin{align*}
g'&= \star|\id|\id|\ldots|\id|: \langle c_1, \ldots, c_{n}\rangle \ra \langle 0,c_1, \ldots, c_{n}\rangle \\
h'&=|1,\underbrace{0,\ldots ,0}_\text{$k-1$ times}|f_{\bullet 1}|f_{\bullet 2}|\ldots |f_{\bullet n}|\underbrace{\star \ldots \star}_\text{$l$ times}:\langle 0,c_1, \ldots, c_{n}\rangle\ra \langle c'_1, \ldots, c'_{m}\rangle.\\
\end{align*}
Now $g'$ and $h'$ are clearly non-identity and $h'\neq h$. This completes the case where the type of $f$ is not inner.

If the type of $f$ is inner, then $$\diff(f) \ = \sum_{1\leq i\leq n} \mbox{diff} (f_{\bullet i}) $$ where $$\diff(g):=\dim(\mbox{target}(g)) -\dim(\mbox{source}(g)).$$ Since $\diff(f)\geq 2$, and each $\diff(f_{\bullet i})\geq 0$ as they are monic, either  $\diff(f_{\bullet i}) >0$ and  $\diff(f_{\bullet j}) >0$ for two distinct parts of $f$, or  $\diff(f_{\bullet i}) >2$ for some $i$.

First we consider the first case, in which  $\diff(f_{\bullet i}) >0$ and  $\diff(f_{\bullet j}) >0.$  Let $f_{k i}, \ldots , f_{li}$ be the components of form $f_{{-} \, i}$ of $f$. Then we can factor $f=hg$ where  $h$ is obtained from $f$ by replacing $f_{\bullet i}$ by identity maps, meaning that
$$h:  \langle c_1, \ldots ,c_{i-1}, c'_{k}, c'_{k+1}, \ldots, c'_{l},c_{i+1}, \ldots, c_{n}\rangle \ra \langle c'_1, \ldots, c'_{m}\rangle$$ with the part of $h$ associated with $c'_{j} \in (c'_k, \ldots, c'_l)$ being the identity map $\id: \langle c'_{j}\rangle \ra \langle c'_{j}\rangle$ and the part of $h$ associated with $c_{j}$ being given by $f_{\bullet j}.$ The map $$g:\langle c_1, \ldots, c_{n}\rangle \ra  \langle c_1,\ldots ,c_{i-1}, c'_{k}, c'_{k+1}, \ldots, c'_{l},c_{i+1}, \ldots, c_{n}\rangle$$ is defined by $g_{\bullet j}= |\id|$ for $j\neq i$ and $g_{\bullet i}=f_{\bullet i}.$ It's easy to see $f=hg$ with $h$ and $g$ monic. If we do the same construction in the opposite way, switching the role of $i$ and $j$, we get a distinct factorization $f=h'g'.$

Finally consider the case where $\diff(f_{\bullet i}) >2$ for some $i$. It's easy to see how any factorization of $f_{\bullet i}$ can be used to give a factorization of $f$, so we need only show that $f_{\bullet i}$ itself has two distinct non-trivial factorizations. This reduces to proving the proposition in the case where $f$ has one part, so that $f$ has the form $\langle c_1\rangle \ra \langle c'_1, \ldots, c'_{m}\rangle.$

We consider three cases, in which $f$ has one, two, or more than two components. If $f$ has one component, that component is a monic map in $\Delta$ which increases dimension by $2$ or more. This case relies on the fairly obvious fact that such a map has at least two non-trivial factorizations; details are left to the reader.

If $f$ has two components, then $c'_1+c'_2 > m$. Consider the sequence $$f_{11}(1)+f_{21}(1), f_{11}(2)+f_{21}(1), \ldots ,f_{11}(m)+ f_{21}(m).$$ First suppose $f_{11}$ and $f_{21}$ are surjective, then  $f_{11}(1)+f_{21}(1)=0$ and $f_{11}(m)+ f_{21}(m)=c'_1+c'_2.$ Since there are $m+1$ terms in the sequence, there must be some $i$ such that $$f_{11}(i)+f_{21}(i)+1 < f_{11}(i+1)+f_{21}(i+1)$$ and therefore since $f_{11}$ and $f_{21}$ are surjective, we must have $f_{11}(i+1)=f_{11}(i)+1$ and $f_{21}(i+1)=f_{21}(i)+1.$ Then let $h,h': \langle c_1+1 \rangle \ra \langle c'_1, \ldots, c'_{m}\rangle$ be defined to be 
\begin{align*}
h&:= | f_{11}(1)f_{11}(2)\ldots f_{11}(i) f_{11}(i) f_{11}(i+1)\ldots f_{11}(m),\ \\ &\quad \quad f_{21}(1)f_{21}(2)\ldots f_{21}(i) f_{21}(i+1) f_{21}(i+1)\ldots f_{21}(m)|\\
h'&:= | f_{11}(1)f_{11}(2)\ldots f_{11}(i) f_{11}(i+1) f_{11}(i+1)\ldots f_{11}(m),\ \\   &\quad \quad f_{21}(1)f_{21}(2)\ldots f_{21}(i) f_{21}(i) f_{21}(i+1)\ldots f_{21}(m)|.
\end{align*}
The maps $h$ and $h'$ are distinct and non-identity, and letting $$g=g'=|\hat{d}_{i+1}|,$$ we have $f=hg=h'g'.$ 

In the case where $f_{11}$ or $f_{21}$ are surjective, without loss of generality assume $f_{11}$ is not surjective, with $i$ not in the image of $f_{11}$. Then we can factor $f_{11}= \hat{d}_i f'_{11}$. Likewise there is some $j$ for which $f_{11}(j)<i$ and $f_{11}(j+1)>i$, and we can factor $f''_{11} \hat{d}_{j}$ with $$f''_{11}(k)= \begin{cases} 
      f_{11}(k) & k<j\\
      i & k=j \\
      f_{11}(k-1) & k>j.
   \end{cases}$$  Then we have $$f=|f_{11}, f_{21}| = |\hat{d}_i| \id| \circ | f'_{11}, f_{21}| = |f''_{11}, f_{21} \hat{s}_j|\circ |\hat{d}_j|    $$
which are two distinct factorizations of $f$ by monics. 

If $f$ has at least three components $f_{11},\ldots, f_{m1}.$ Then let $h_1g_1$ be the epi-monic factorization of $|f_{11},\ldots, f_{(m-1)1}|$ and let $h_2g_2$ be the epi-monic factorization of $f_{m1}$. Then $$f=|h_1|h_2|\circ |g_1,g_2|$$ where $|g_1,g_2|$ is the map with one part and components given by the components of $g_1$ followed by $g_2$ as the last component. Note that $|h_1|h_2|$ is clearly monic (since it has two parts, each of which is monic) and both $|h_1|h_2|$ and $|g_1,g_2|$ are clearly non-identity. So letting $h= |h_1|h_2|$ and $g= |g_1,g_2|$ gives the factorization $f=hg$ as desired. The same construction can be done in the opposite way, ``breaking off'' $f_{11}$ instead of $f_{m1}$, leading to a distinct factorization $f=h'g'$.
\end{proof}
\begin{corollary}
$\Theta_2$ is a good dimensional category.
\end{corollary}
\begin{proof} Every morphism $f:a\ra c$ in $\Theta_2$ with $\dim a < \dim c$ either is a coface map or factors through a coface map by Proposition~\ref{thetafactorprop}.
\end{proof}
\begin{corollary}
For $\Theta_2$-sets, every universal horn obtained by removing a single face from a sphere is nice, i.e. each has a minimal complementary dimension one less than its dimension.
\end{corollary}
\begin{proof}This follows from the fact that each cell of $\Theta_2[c]$ of dimension $\dim(c)-2$ or smaller factors through two different coface maps by Proposition~\ref{thetafactorprop}.
\end{proof}
\begin{definition}Let $c=\langle c_1,c_2, \ldots,c_{n}\rangle$ and let $1\leq l \leq l' \leq m' \leq m\leq n$. We define the \emph{$(l',m')$-restriction}  $p_{l'm'}$ to be the cartesian projection $$p: [c_{l}]\times [c_{l+1}]\times \ldots \times [c_m] \ra [c_{l'}]\times [c_{l'+1}]\times \ldots \times [c_{m'}].$$ 

We call a subset of the form $S \subseteq  [c_{l}]\times [c_{l+1}] \times \ldots \times [c_m]$ a \emph{profile of kind $(l,m)$} of $c.$ We also allow an \emph{empty profile} $\emptyset$, which is the empty subset of the empty product. Taxing $S$ of kind $(l,m)$ and $S'$ of kind $(l',m')$ to be arbitrary profiles of $c$, we say $S' \sqsubseteq S$ if $ l \leq l' \leq m' \leq m$ and $S' \subseteq p_{l'm'}(S).$ For the same $S,S'$, let $(l'',l''+1, \ldots, m'')$ be the intersection of the intervals $(l,l+1,\ldots, m)$ and $(l',l'+1,\ldots,m')$, if non-empty. We define $S \sqcap S' :=\emptyset$ if this intersection is empty, otherwise $S \sqcap S':= p_{l''m''}S\cap p_{l''m''}S'.$
\end{definition}
\begin{fact}\label{proffact}
If $S,S'$ are profiles of $c$, then $S \sqsubseteq S'\sqcap S''$ if and only if $S\sqsubseteq S'$ and $S \sqsubseteq S''$.
\end{fact}

\begin{definition}
We call a map $f:c\ra c'$ in $\Theta_2$ \emph{primary} if $c$ is of the form $\langle c_1\rangle$, so that $f$ is of the form $\langle f_{l1}, \ldots, f_{m1} \rangle$. The \emph{image} of the primary map $f$ is the $(l,m)$ profile  $$\image(f):= \{ (f_{l1}(i),\ldots, f_{m1}(i))\ |\ i \in [c_l]  \}\subseteq [c'_l]\times [c'_{l+1}]\times \ldots \times [c'_m].$$ For an arbitrary map $$f:\langle c_1,\ldots, c_{n}\rangle \ra \langle c'_1,\ldots, c'_m \rangle $$ in $\Theta_2$, we define 
$$\image(f):=\image(f_{\bullet 1})\times \image(f_{\bullet 2})\times \ldots \times \image(f_{\bullet n})\subseteq [c'_l]\times [c'_{l+1}]\times \ldots \times [c'_m]$$ where $l = \Ty(f)(0)$ and $m=\Ty(f)(k')-1$. 
\end{definition}

\begin{proposition}\label{imfactorprop}
Let $f,g$ be maps in $\Theta_2$.  Then $g$ factors through $f$ as $g=f\circ h$ if and only if:
\begin{enumerate}
\item $\Ty(g)$ factors through $\Ty(f)$
\item $\image(g)\sqsubseteq\image(f).$ 
\end{enumerate}
\end{proposition}
\begin{proof}

It's not hard to see that the following criteria for factorization follows from the definition of $\Theta_2$: the map $g$ factors through $f$ as $g=f\circ h$ if and only if
\begin{enumerate}\setcounter{enumi}{2}

\item $\Ty(g)$ factors through $\Ty(f)$
\item For every non-empty part $f_{\bullet i}=|f_{li}, \ldots, f_{mi}|$ of $f$, we have that  $|g_{lj},g_{(l+1)j},\ldots,g_{mj}|$ factors through $f_{\bullet i}$ whenever there is a $j$ such that the components $g_{lj},g_{(l+1)j},\ldots,g_{mj}$ exist.
\end{enumerate}
Note that if $(3)$ holds there and $g_{lj}$ exists then all of $g_{lj},\ldots, g_{mj}$ exist as components of $g$.  In this case we say that $g$ \emph{hits} $f_{\bullet i}$, otherwise we say $g$ \emph{misses} $f_{\bullet i}.$
We must show that $(4)$ above is equivalent to $(2)$ in case $\Ty(g)$ factors through $\Ty(f).$ 

From the definition $\image(f):=\image(f_{\bullet 1})\times \image(f_{\bullet 2})\times \ldots \times \image(f_{\bullet m'})$ it follows that $\image (g)\sqsubseteq \image (f)$ if and only if for every part $f_{\bullet i}=|f_{li}, \ldots, f_{mi}|$ we have either $g$ misses $f_{\bullet i}$ or
$$p_{lm}(\image(g))=\image(|g_{lj},g_{(l+1)j},\ldots,g_{mj}|)\subseteq p_{lm}(\image(f))=\image(f_{\bullet i}).$$
So we must show that  $|g_{lj},g_{(l+1)j},\ldots,g_{mj}|$ factors through $f_{\bullet i}$ if and only if $$\image(|g_{lj},g_{(l+1)j},\ldots,g_{mj}|)\subseteq \image(f_{\bullet i}).$$ 
Let $\langle c_j \rangle$ and $\langle c_i\rangle$ be the domains of the primary maps $g':=|g_{lj},g_{(l+1)j},\ldots,g_{mj}|$ and $f_{\bullet i}$ respectively. If $\image(g')\subseteq \image(f_{\bullet i})$ then we can define a map $h:\langle c_j \rangle \ra \langle c_i\rangle$ to be the map with a unique component $h_{11}$ such that $h_{11}(k)=k'$ where $k'$ is the smallest value for which $$(g_{lj}(k), g_{(l+1)j}(k),\ldots,g_{mj}(k))=(f_{li}(k'),f_{(l+1)i}(k'),\ldots f_{mi}(k'))$$
It's easy to see $g'=f_{\bullet i} h.$ The other direction, stating if such an $h$ exists, then $\image(g')\subseteq \image(f_{\bullet i}),$ is immediate.
\end{proof}


We can characterize the coface maps in $\Theta_2$ with target $\langle c_1, \ldots, c_{n} \rangle$:
\begin{itemize}
\item If $c_1=0$ there is a coface map $$\hat{d}_{\star|}:=\star|\id|\id|\ldots|\id|$$ whose type is $\hat{d}_0= 12\ldots k$ and whose components are identities. Likewise if $c_{n}=0$ there is a coface map $$\hat{d}_{|\star}:=|\id|\id|\ldots|\id|\star$$ whose type is $\hat{d}_n$ and whose components are identities.
\item There is a coface map $$\hat{d}^i_{m}:=|\id| \id| \ldots |\id| \hat{d}_m|\id| \ldots |\id|        $$ whose type is the identity and whose components are identities, its $i$th component which is $\hat{d}_m.$ 
\item  Let $f:[c_i+c_{i+1}]\ra [c_i]\times[c_{i+1}]$ be a strictly increasing map of partially ordered sets, where $ [c_i]\times[c_{i+1}]$ is given the product partially ordered set structure where $(a,b)\leq (a',b')$ if $a\leq a'$ and $b\leq b'$. Then there is a face map $$h:=  {\hat{d}}^i_{|f|}:\langle c_1, c_2, \ldots, c_{i-1}, c_i+c_{i+1}, c_{i+2},\ldots, c_{n} \rangle \ra \langle c_1, c_2, \ldots, c_{n}\rangle$$ whose type is $\hat{d}_{i+1}$ and whose components are identities except $h_{ii}$ and $h_{(i+1)i}$ where $$(h_{ii}(m),h_{(i+1)i}(m))=(f_1(m), f_2(m)):=f(m).$$ We call $f$ the \emph{core} of $ \hat{d}^i_{|f|}$.
\end{itemize}
\begin{definition}
Consider a functor $F:\mathcal{C} \ra \mathcal{D}$ $d \in \mathcal{D}$ and an object $d\in \mathcal{D}$. The \emph{comma category} $(d/ F)$ is the category whose objects are arrows $f:F(c) \ra d$ and whose morphisms from $f:F(c)\ra d$ to $g:F(c')\ra d$ are morphisms $c \ra c'$ in $\mathcal{C}$ which make the obvious triangle commutative. $F$ is called \emph{cofinal} if for every $d \in \mathcal{D}$ the category $(d/F)$ is non-empty and connected.
\end{definition}

Let $c=\langle c_1, \ldots, c_{n}\rangle \in \Theta_2$ and let $\mathfrak{f}:c'\ra c$ be a face map. Recall that the category of objects $O(c,\mathfrak{f})$ of the universal horn $ \Lambda_{\{\mathfrak{f}\}}^{c} \subseteq d\Theta_2[c]$ is the category of morphisms to $c' \ra c$ in $\Theta_2$ which factor through a face map other than $\mathfrak{f}$, and commuting triangles between them. The Glenn category $G(c,\mathfrak{f})$ of $\Lambda_{\{\mathfrak{f}\}}^{c}$ is the subcategory of $O(c,\mathfrak{f})$ containing those $f:c'\ra c$ which are either a coface maps or a compositions of two coface maps, with morphisms between these given by precompositions by coface maps.

The following proposition, along with Proposition~\ref{cofinalprop2}, justifies the use of Glenn tables for $\Theta_2$-sets:
\begin{proposition} \label{cofinalprop} A horn of the form $\Lambda_{\{\mathfrak{f}\}}^{c}$ is tabular i.e the inclusion $I:G(c,\frakf) \ra O(c,\frakf)$ is cofinal.
\end{proposition}
\begin{proof}
Let $f:d \ra c$ be an object of $O(c,\mathfrak{f})$. The category $(f/I)$ has the following description:
\begin{itemize}
\item An object of $(f/I)$ is a factorization $f=h\circ g$ where $h$ is either a coface map other than $\mathfrak{f}$ or a composition of two coface maps other than $\mathfrak{f}$ (by Proposition~\ref{thetafactorprop}, the latter includes all $h$ which increase dimension by $2$). 
\item A morphism $h\circ g \ra h'\circ g'$ is a factorization $h'= h \circ h''$ where $h''$ is a coface map.
\end{itemize}
Note that since the coface maps of $\Theta_2$ are monic, the factorization $f = h \circ g$ is uniquely determined by $h$. Thus we will think of the objects of $(f/I)$ as being coface maps or composites of two coface maps (call these \emph{subcoface maps}) through which $f$ factors, without mention of $g$.  There is a morphism from $h$ to every composite $h \circ h''$ through which $f$ also factors. For brevity, we call a map through which $f$ factors \emph{$f$-valid}.

To show that $(f/I)$ is non-empty is trivial, by definition the fact that $f\in O(c,\mathfrak{f})$ means that $f$ factors through a coface map other than $\mathfrak{f.}$

We must show that $(f/I)$ is connected. Let $\smileeq$ be the equivalence relation on objects of $(f/I)$ generated by the morphisms. We must show $\smileeq$ is the complete relation. Clearly every $f$-valid subcoface map is $\smileeq$ to some $f$-valid coface, so it suffices to show any two $f$-valid coface maps are related by $\smileeq$ to each other. We must consider each possibly for two $f$-valid coface maps $h$ and $h'$, and in each case show $h \smileeq h'$. First we handle the ``easy cases'':
\begin{enumerate}
\item $h$ and $h'$ are both of the form $\hat{d}_{\star|}$, $\hat{d}_{|\star},$ or $\hat{d}^i_{m}.$
\item Either $h=\hat{d}_{\star|}$ and $h'= \hat{d}^i_{|f|}$ with $i>1$ or  $h=\hat{d}_{|\star}$ and  $h'= \hat{d}^i_{|f|}$ with $i<n-1.$ 
\item $h=\hat{d}^i_m$  and $h'=\hat{d}^j_{|g|}$ with $i\neq j, j+1$
\item $h=\hat{d}^i_{|f|}$  and $h'=\hat{d}^j_{|g|}$ with $i\neq j-1, j, j+1.$
\end{enumerate}
In each of these cases, there is a pullback square of the form 
\begin{center}
\begin{tikzcd}
  c'      \arrow{r}{h} & c            \\
c'''   \arrow{u}{g} \arrow{r}{g'}     & c'' \arrow{u}{h'}  
\end{tikzcd}
\end{center}
where $g$ and $g'$ are coface maps. Note that such a square of monic maps is a pullback square if and only if it is the case that a map $f$ factors through $h$ and $h'$ if and only if it factors through the pullback map $hg=h'g'$. This in turn implies that if $h$ and $h'$ are $f$-valid, then so is $hg=h'g'$, and $h \smileeq hg=h'g'\smileeq h'$ showing $h\smileeq h'$. By Proposition~\ref{imfactorprop} and Fact~\ref{proffact}, a map $e$ is the pullback map of $h,h'$ if and only if $\Ty(e)$ is the pullback map of $\Ty(h)$ and $\Ty(h')$ and $\image(e)=\image (h)\sqcap \image(h')$. Table~\ref{pullbacktable} gives the pullback squares in some possible cases of $(1)$-$(4)$, and each of these can be easily checked using these criteria.
Each non-listed case is symmetrical to a listed case and follows similarly.

\begin{table}[H] \caption{\label{pullbacktable}~}\begin{center}
    \begin{tabular}{c c c c c }
$h$       &   $h'$   &conditions &  g & g'       \\  \hhline{-----}  
$\hat{d}_{\star|}$ &$\hat{d}_{|\star}$ & &$\hat{d}_{|\star}$ & $\hat{d}_{\star|}$  \\
$\hat{d}_{\star|}$ &   $\hat{d}^i_{m}$ &  & $\hat{d}^{i-1}_m$& $ \hat{d}_{\star|}$  \\
 $\hat{d}^i_{m}$ & $\hat{d}^j_{m}$ & $i\neq j$ & $\hat{d}^j_{m}$     & $\hat{d}^i_{m}$   \\
 $\hat{d}^i_{m}$ & $\hat{d}^i_{m'}$ &  $m<m'$     &$\hat{d}^i_{m'-1}$ &  $\hat{d}^i_{m}$  \\
$\hat{d}_{\star|}$  & $\hat{d}^i_{|f|}$ &  $i>1$   & $\hat{d}^{i-1}_{|f|}$ &  $\hat{d}^i_{|f|}\hat{d}_{\star|}$\\
$\hat{d}^i_m$  & $\hat{d}^j_{|g|}$ &  $i<j$     & $\hat{d}^j_{|g|}$&$\hat{d}^i_m$   \\
$\hat{d}^i_m$  & $\hat{d}^j_{|g|}$ &  $i>j+1$    &$\hat{d}^j_{|g|}$ & $\hat{d}^{i-1}_m$    \\
$\hat{d}^i_{|f|}$ & $\hat{d}^j_{|g|}$ &  $i<j-1$  &  $\hat{d}^{j-1}_{|g|}$ & $\hat{d}^i_{|f|}$  
    \end{tabular}\end{center}
    \end{table}

We are left with the following cases
\begin{enumerate}[resume]
\item Either $h=\hat{d}_{\star|}$ and $h'= \hat{d}^1_{|f|}$  or  $h=\hat{d}_{|\star}$ and  $h'= \hat{d}^{n-1}_{|f|}.$ 
\item $h=\hat{d}^i_m$  and $h'=\hat{d}^j_{|g|}$ with $i= j$ or $i= j+1$
\item  $h=\hat{d}^i_{|g|}$  and $h'=\hat{d}^j_{|g'|}$ with $i=j$
\item  $h=\hat{d}^i_{|g|}$  and $h'=\hat{d}^j_{|g'|}$ with $i=j-1$ and $i=j+1.$
\end{enumerate}
For case $(5)$, first suppose $h=\hat{d}_{\star|}$ and $h'= \hat{d}^1_{|f|}$ with $c_2=0$. In this case $d_{\star|}d_{\star|}$ can be seen to be the pullback map of $h$ and $h'$, letting us conclude $h\smileeq h'$ by the same argument we used for cases $(1)$-$(4)$. Otherwise if $c_2>0$ it's easy to see that $\hat{d}^2_0$ and $\hat{d}^2_1$ are $f$-valid since the $\Ty(f)$ must factor through $\hat{d}_0\hat{d}_0$, and any such map can easily be seen to factor through $\hat{d}^2_0$ and $\hat{d}^2_1$. At least one of these two maps is not $\frakf$, without loss of generality $\hat{d}^2_0$. Then 
 $h=\hat{d}_{\star|}\smileeq \hat{d}^2_0$ by case $(1)$ and $h'= \hat{d}^1_{|f|}\smileeq \hat{d}^1_0$ by case $(6)$ which is shown below. We conclude that $h \smileeq h'$. The case where $h=\hat{d}_{|\star}$ and  $h'= \hat{d}^{n-1}_{|f|}$ is similar. 

The other three cases require a more complicated combinatorial proof, especially for $(7)$. We first work towards a lemma that will enable us to analyse these cases.

Consider coface maps $\hat{d}^i_{|g|}$  and $\hat{d}^i_{|g'|}.$   Because $$\image(\hat{d}^i_{|g|})=[c_1]\times [c_2]\times \ldots \times [c_{i-1}] \times \image{g} \times [c_{i+2}] \times \ldots \times [c_{n}]$$ where $\image{g}=p_{i,i+1} \image(\hat{d}^i_{|g|}) \subseteq [c_i]\times [c_{i+1}],$ if $\Ty(f)$ factors through $\Ty(\hat{d}^i_{|g|})=\hat{d}_{i+1}$ we have $\image{f}\sqsubseteq \image(\hat{d}^i_{|g|})$ and equivalently $\hat{d}^i_{|g|}$ is $f$-valid if and only if either $p_{i,i+1}$ is undefined or $p_{i,i+1} \image(f) \subseteq \image(g)$. For the purposes of this proof, we treat $p_{l,m} \image(f)=\emptyset$ whenever it is undefined. With this convention, call the fact above the \emph{$f$-validity Criterion}:
\begin{itemize}
\item Suppose $\Ty(f)$ factors through $\Ty(\hat{d}^i_{|g|})=\hat{d}_{i+1}.$ Then $\hat{d}^i_{|g|}$ is $f$-valid if and only if  $p_{i,i+1} \image(f) \subseteq \image(g).$
\end{itemize}

The set $\image(g)\sqcap \image(g')$ is necessarily totally ordered in the partially ordered set $[c_i]\times [c_{i+1}]$  and contains $(0,0)$ and $(c_i, c_{i+1})$. If $k$ is the cardinality of this set, there is a map monic map $$e^i_{|g|,|g'|}  :\langle c_1,\ldots,c_{i-1}, k-1 ,c_{i+2},\ldots, c_{n}\rangle\ra \langle c_1, \ldots, c_{n}\rangle $$ such that  $\image(   e^i_{|g|,|g'|}    )= \image(\hat{d}^i_{|g|})\sqcap \image(\hat{d}^i_{|g'|}).$ If $k= c_i+ c_{i+1}$, then $e^i_{|g|,|g'|}$ increases dimension by $2$. By Proposition~\ref{imfactorprop} and Fact~\ref{proffact} if $\hat{d}^i_{|g|}$  and $\hat{d}^i_{|g'|}$ are $f$-valid, then $\image(f)\sqsubseteq  \image(   e^i_{|g|,|g'|}    )$     and therefore $e^i_{|g|,|g'|}$ is $f$-valid and is an object in $(f/I)$. Since $\image( e^i_{|g|,|g'|} )\sqsubseteq \image(\hat{d}^i_{|g|})$ and      $\image( e^i_{|g|,|g'|} )\sqsubseteq \image(\hat{d}^i_{|g'|})$, we have that $e^i_{|g|,|g'|}$ factors through $\hat{d}^i_{|g|}$  and $\hat{d}^i_{|g'|}.$ Thus $\hat{d}^i_{|g|}\smileeq   e^i_{|g|,|g'|}\smileeq \hat{d}^i_{|g'|}.$ We will call what we have proven the Neighbor Rule:
\begin{itemize}
\item Suppose $\hat{d}^i_{|g|}\neq \frakf$  and $\hat{d}^i_{|g'|}\frakf$ are $f$-valid, thus objects of $(f/I)$. If $\image(g)\sqcap \image(g')$ has cardinality $c_i+c_{i+1}$, i.e. one less than the maximum cardinality for a totally ordered subset of $[c_i+c_{i+1}]$, then $\hat{d}^i_{|g|}\smileeq \hat{d}^i_{|g'|}.$ 
\end{itemize}
Now consider the case $(6)$. Since $h'=\hat{d}^j_{|g|}$ is $f$-valid, we have  $p_{i,i+1} \image(f) \subseteq \image(g)$,
and since $h=\hat{d}^i_m$ is $f$-valid, $$p_{i,i+1} \image(f) \subseteq \{0,1,\ldots, \hat{m}, \ldots, c_{i}\}\times [c_{i+1}].$$ Thus $$p_{i,i+1} \image(f) \subseteq \image(g)\setminus (\{m\}\times [c_{i+1}]).$$ We proceed by induction on the cardinality $|\image(g)\sqcap (\{m\}\times [c_{i+1}])|.$   First consider the case $|\image(g)\sqcap (\{m\}\times [c_{i+1}])|=1$, so that $$| \image(g)\setminus (\{m\}\times [c_{i+1}])|=c_{i}+c_{i+1}.$$ In this case there  is a map monic map $$e^i_{|g|,m}  :\langle c_1,\ldots,c_{i-1}, c_{i}+c_{i+1}-1 ,c_{i+2},\ldots, c_{n}\rangle\ra \langle c_1, \ldots, c_{n}\rangle $$ such that  $$\image(   e^i_{|g|,m}    )= [c_1]\times \ldots \times [c_{i-1}] \times (\image(g)\setminus (\{m\}\times [c_{i+1}])) \times [c_{i+2}]\times \ldots [c_{n}].$$
The map $e^i_{|g|,m}$ increases dimension by $2$ and is $f$-valid, and since $\image(e^i_{|g|,m})\sqsubseteq \image(\hat{d}^j_{|g|})$ and $\image(e^i_{|g|,m})\sqsubseteq \hat{d}^i_m$ it follows that $e^i_{|g|,m}$ factors through $\hat{d}^i_m$ and $\hat{d}^j_{|g|}$ therefore $\hat{d}^i_{|g|}\smileeq   e^i_{|g|,m}\smileeq \hat{d}^i_m.$

Now suppose $|\image(g)\sqcap (\{m\}\times [c_{i+1}])|=q>1$, and assume as our induction hypothesis that $ \hat{d}^j_{|g|}\smileeq \hat{d}^i_m$ whenever $|\image(g)\sqcap (\{m\}\times [c_{i+1}])|< q$. The map $g:[c_{i}+c_{i+1}]\ra [c_{i}]\times [c_{i+1}]$ has the form $$[(0,0),\ldots, (m-1,k),(m,k),(m,k+1),\ldots, (m,l),(m+1,l),\ldots, (c_i,c_{i+1}) ],$$ with $p_{i,i+1}\image(f)$ contained within $\image(g)\sqcap (\{0,1,\ldots, \hat{m}, \ldots, c_{i}\}\times [c_{i+1}]).$ 
\begin{align*}g'&=[(0,0),\ldots, (m-1,k),(m-1,k+1),(m,k+1),\ldots, (m,l),(m+1,l),\ldots, (c_i,c_{i+1}) ]\\
g''&= [(0,0),\ldots, (m-1,k),(m,k),\ldots, (m,l-1) ,(m+1,l-1),(m+1,l),\ldots, (c_i,c_{i+1}) ]
\end{align*}
Clearly $p_{i,i+1}\image(f)\sqsubseteq \image(g')$ and $p_{i,i+1}\image(f)\sqsubseteq \image(g')$ so $\hat{d}^i_{|g'|}$ and $\hat{d}^i_{|g'|}$ are $f$-valid, and $g'\neq g''$ so at least one of $\hat{d}^i_{|g'|}\neq \frakf$ or $\hat{d}^i_{|g''|}\neq\frakf$ holds. Assuming without loss of generality that $\hat{d}^i_{|g'|}\neq \frakf$, then $\hat{d}^i_{|g'|}$ is in $(f/I)$ with $\hat{d}^i_{|g'|}\smileeq \hat{d}^i_{|g|}$ by the Neighbor Rule, and $d^i_{|g'|}\smileeq d^i_m$ by the induction hypothesis. Therefore $\hat{d}^i_{|g|}\smileeq \hat{d}^i_m$, completing case $(6)$.

In case $(7),$  we take $h=\hat{d}^i_{|g|}$  and $h'=\hat{d}_{|g'|}$ in $(f/I).$ We first introduce a notation for specifying an arbitrary strictly increasing map $e:[c_i+c_{i+1}]\ra [c_i]\times [c_{i+1}]$, which is associated a face map $\hat{d}^i_{|e|}.$ The map $e$ is given by a maximal ordered subset of elements of $[c_i ]\times [c_{i+1}],$ which can be pictured  as a path through the grid $[c_i]\times [c_{i+1}].$ This path can be specified by a word consisting of letters $\mathtt{r}$ and $\mathtt{u}$, with the number of $\mathtt{r}$'s equal to $c_i$ and the number of $\mathtt{u}$'s equal to $c_{i+1}$, telling us the moves ``right'' and ``up'' that trace out the path. Figure~\ref{interiorpath} gives an example of this notation.
\begin{figure} \caption{ \label{interiorpath} The map $e:[6]\ra[4]\times [2]$ given by $[(0,0), (1,0), (1,1),(2,1),(2,2),(3,2),(4,2)].$ This map is denoted $\mathtt{rururr}.$ }
\begin{center}
\begin{tikzpicture}[scale=1.4,auto]
\begin{scope}

\foreach \x in {0,1,2,3,4}
\foreach \y in {0,1,2}
{

\fill (\x,\y) circle (0.03cm);
}

\draw (0,0)--(1,0)--(1,1)--(2,1)--(2,2)--(3,2)--(4,2);

\end{scope}
\end{tikzpicture}
\end{center}
\end{figure}
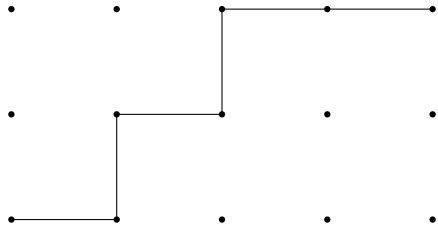 
If $\hat{d}^i_{|e|}$is $f$-valid then $p_{i,i+1}\image(f)\subseteq \image(e).$ We find it helpful to specify this subset within the notation for $e$. We do this by adding symbols $\upspoon$ to the word specifying $e$, with each $\upspoon$ signifying an element $(r,u)$ in $p_{i,i+1}\image(f)$ with $r$ being the number of $\mathtt{r}$'s to the left of this $\upspoon$ and $u$ being the number of $\mathtt{u}$'s to the left of the $\upspoon.$ We call this word the \emph{spoon notation for $e$, $\spoon(e)$}. We call the pieces of  $\spoon(e)$ between successive $\upspoon$'s  (or before the first $\upspoon$ or after the last $\upspoon$) the \emph{segments} of $\spoon(e)$. 
\begin{figure} \caption{ \label{interiorpath2}  $e$ from Figure~\ref{interiorpath} with subset $p_{i,i+1}\image(f)=\{(1,0),(2,1),(4,2)\}\subseteq \image(e),$ which are circled above. The notation $\spoon(e)= \mathtt{r}\alddiv \mathtt{ur}\alddiv\mathtt{urr}\alddiv$ specifies $e$ together with this subset. }
\begin{center}
\begin{tikzpicture}[scale=1.4,auto]
\begin{scope}

\foreach \x in {0,1,2,3,4}
\foreach \y in {0,1,2}
{

\fill (\x,\y) circle (0.03cm);
}

\draw (0,0)--(1,0)--(1,1)--(2,1)--(2,2)--(3,2)--(4,2);
\draw (1,0) circle (.1cm);
\draw (2,1) circle (.1cm);
\draw (4,2) circle (.1cm);
\end{scope}
\end{tikzpicture}
\end{center}
\end{figure} 

If $\hat{d}^i_{|e|}$ is $f$-valid then a word $w$ of letters $\mathtt{u}$, $\mathtt{r}$, and $\upspoon$ specifies a map $e'$ with $\hat{d}^i_{|e'|}$ being $f$-valid and $\spoon(e')=w$ and if and only if $w$ and $\spoon(e)$ have the same number of $\upspoon$'s and the corresponding segments of $\spoon(e)$ and $w$ each have the same number of $\mathtt{u}$'s and $\mathtt{r}$'s. In this case we say  $w$ is $e$-\emph{aligned}. So $\spoon$  is a bijective correspondence between $e$-aligned words and maps $e'$ such that $\hat{d}^i_{|e'|}$ is $f$-valid.  

Suppose $\hat{d}^i_{|e|}$, $\hat{d}^i_{|e'|}$ are $f$-valid. It's not hard to see $|\image(e)\sqcap \image(e')|=c_i +c_{i+1}$ if $\spoon(e')$ is obtained from  $\spoon(e)$ by swapping an $\mathtt{r}$ with an adjacent $\mathtt{u}.$ Note that this swap must be done within a segment to preserve $f$-validity. We call such a move a \emph{basic swap}. If $w$ is obtained from $\spoon(e)$ by a basic swap, then $w$ is $e$-aligned so that $\spoon^{-1}(w)$ exists with $\hat{d}^i_{|\spoon^{-1}(w)|}$ being $f$-valid and  $\hat{d}^i_{|e|}\smileeq \hat{d}^i_{|\spoon^{-1}(w)|}$ by the Neighbor Rule, so long as $\hat{d}^i_{|\spoon^{-1}(w)|}\neq \frakf$. 

It's clear that any two $e$-aligned words are related by a series of basic swaps. For $f$-valid coface maps $\hat{d}^i_{|g|}$ and $\hat{d}^i_{|g'|}$ if $\frakf$ is not of the form $\hat{d}^i_{|\mathfrak{g}|}$, then we can relate $\spoon(g)$ and $\spoon(g')$ by a series of basic swaps, with each intermediate word $w$ having the property that $\hat{d}^i_{|\spoon^{-1}(w)|} \in (f/I).$ By the above discussion, we can conclude $\hat{d}^i_{|g|}\smileeq \hat{d}^i_{|g'|}$.

Next we consider the case  $\frakf=\hat{d}^i_{|\mathfrak{g}|}.$ We consider the following cases for $\hat{d}^i_{|g|}$ and $\hat{d}^i_{|g'|}$, bearing in mind that $\spoon(g')$ is $g$-aligned and vice versa:
\begin{enumerate}[label=\roman*.]
\item $\spoon(g)$ and $\spoon(g')$ have two or more segments which contain at least one $\mathtt{r}$ and at least one $\mathtt{u}$
\item $\spoon(g)$ and $\spoon(g')$ have a segment which has either two or more instances of $\mathtt{r}$ or two or more instances of $\mathtt{u}$
\item $\spoon(g)$ and $\spoon(g')$ have only one segment with two letters, which are a $\mathtt{r}$ and a $\mathtt{u}$ 
\item Every segment of $\spoon(g)$ and $\spoon(g')$ has exactly one letter. 
\end{enumerate}
In case $(\mathrm{i})$ it's not hard to see that one can move from $\spoon(g)$ to $\spoon(g')$ by a series of basic swaps, while avoiding $\spoon(\mathfrak{g})$, showing that  $\hat{d}^i_{|g|}\smileeq \hat{d}^i_{|g'|}.$ Details are left to the reader.

In case $(\mathrm{ii})$, without loss of generality assume $\spoon(g)$ has a segment which has either two or more instances of $\mathtt{r}$. Then there is some $m$ for which $p_{i,i+1}\image(f) \cap (\{m\}\times [c_{i+1}])=\emptyset.$ Then $\hat{d}^i_m$ is $f$-valid and $\hat{d}^i_m\neq \frakf$ by the  assumption that $\frakf=\hat{d}^i_{|\mathfrak{g}|}.$ So $\hat{d}^i_m\in (f/I)$ with case $(5)$ proved above showing  $\hat{d}^i_{|g|}\smileeq \hat{d}^i_m \smileeq  \hat{d}^i_{|g'|}.$ 

In case $(\mathrm{iii})$, $\spoon(g)$ and $\spoon(g')$ are either identical or directly related by a basic swap, showing $\hat{d}^i_{|g|}\smileeq \hat{d}^i_{|g'|}$ by the Neighbor Rule. 

Finally in case $(\mathrm{iv})$ the maps $g$ and $g'$ are necessarily identical. We have shown $\hat{d}^i_{|g|}\smileeq \hat{d}^i_{|g'|}$ in every case, completing case $(7).$

For case $(8)$ we have $f$-valid maps $h=\hat{d}^{i-1}_{|g|}$  and $h'=\hat{d}^{i}_{|g'|}$ in $(f/I).$ Since $\Ty(f)$ factors through $\Ty(h)=\hat{d}_i$ and $\Ty(h')=\hat{d}_{i+1}$, it factors through the map $\hat{d}_i\hat{d}_i$ which skips both $i$ and $i+1$. Thus $f$ either has no components of the form $f_{(i-1)-},$ $f_{i-}$ or $f_{(i+1)-}$ in which case $p_{i-1,i+1} \image(f)=\emptyset$ by convention or there is an $l$ such that  $f_{(i-1)l},$ $f_{il}$, and $f_{(i+1)l}$ exists. In either case $p_{i-1,i+1} \image(f)$ is a totally ordered subset of $[c_{i-1}]\times [c_i]\times [c_{i+1}]$. 

We now construct sequences
 \begin{align*}
k \ra \alpha_k&:[c_{i-1}+c_i+c_{i+1}]\ra [c_{i-1}] \\ 
k \ra \beta_k&:[c_{i-1}+c_i+c_{i+1}]\ra [c_{i}] \\ 
 k\ra   \gamma_k&:[c_{i-1}+c_i+c_{i+1}]\ra [c_{i+1}] 
\end{align*}
such that for all $0\leq k \leq c_{i-1}+c_i+c_{i+1}:$
\begin{enumerate}[label=\roman*.]
\item $\alpha_k+\beta_k+\gamma_k=k$
\item $(\alpha_k,\beta_k)\in \image(g)$
\item $(\beta_k,\gamma_k)\in \image(g')$
\item $(\alpha_k,\beta_k,\gamma_k) \leq (s_{i-1},s_i,s_{i+1})$ in the partial order of $ [c_{i-1}]\times [c_i]\times [c_{i+1}]$ for all $ (s_{i-1},s_i,s_{i+1})\in p_{i-1,i+1} \image(f)$ such that $k\leq s_{i-1}+s_i+s_{i+1}$. In particular if $$k=\alpha_k+\beta_k+\gamma_k= s_{i-1}+s_i+s_{i+1}$$ this implies that  $(\alpha_k,\beta_k,\gamma_k) = (s_{i-1},s_i,s_{i+1})$
\end{enumerate}
We proceed inductively, defining $(\alpha_0,\beta_0,\gamma_0) = (0,0,0)$, which clearly meets the conditions above. To define $(\alpha_m,\beta_m,\gamma_m)$ for $m>0$, assume that we have defined  $(\alpha_k,\beta_k,\gamma_k)$ meeting conditions (i)-(iv) above for each $k<m$.

First consider the case where $(\alpha_{m-1}+1,\beta_{m-1})\in \image(g)$. Then either the triple $$(\alpha_{m-1}+1,\beta_{m-1}, \gamma_{m-1})$$ meets each of our conditions, in which case we define  $$(\alpha_m,\beta_m,\gamma_m)=(\alpha_{m-1}+1,\beta_{m-1}, \gamma_{m-1})$$ or there is an $ (s_{i-1},s_i,s_{i+1})\in p_{i-1,i,i+1}\image(f)$ such that  $m\leq s_{i-1}+s_i+s_{i+1}$ and $$(\alpha_{m-1}+1,\beta_{m-1}, \gamma_{m-1})\nleq (s_{i-1},s_i,s_{i+1}).$$ By the induction hypothesis, $$(\alpha_{m-1},\beta_{m-1}, \gamma_{m-1}) \leq  (s_{i-1},s_i,s_{i+1})$$ thus we conclude $s_{i-1}=\alpha_{m-1}$. By the $f$-validity Criterion we know that $p_{i-1,i}\image(f)\subseteq \image(g)$ thus $ (s_{i-1},s_i)=(\alpha_{m-1}, s_i)\in \image(g)$. If $s_i>\beta_{m-1}$ then $(\alpha_{m-1}+1,\beta_{m-1})$ and $(\alpha_{m-1}, s_i)$ are elements of $\image(g)$ which are incomparable in the partial order on $[c_{i-1}]\times [c_i]$ which is impossible since $\image(g)$  is a totally ordered subset. We conclude $s_i = \beta_{m-1}$. The fact that $$m=\alpha_{m-1}+\beta_{m-1}+\gamma_{m-1}+1\leq s_{i-1}+s_i+s_{i+1}=\alpha_m+\beta_m+s_{i+1}$$ lets us conclude that $\gamma_{m-1}<s_{i+1}$. Since $$(\beta_{m-1},\gamma_{m-1})<(s_i, s_{i+1})=(\beta_{m-1},s_{i+1})$$ are both in  $p_{i,i+1}\image(f)\subseteq \image(g')$, and $\image(g')$ is a maximal totally ordered subset of $[c_i]\times [c_{i+1}],$ we have that $(\beta_{m-1},\gamma_{m-1}+1)\in \image(g')$. We therefore define  $$(\alpha_m,\beta_m,\gamma_m)=(\alpha_{m-1},\beta_{m-1}, \gamma_{m-1}+1),$$ which is easily seen to meet conditions (i)-(iv).

We considered the case where $(\alpha_{m-1}+1,\beta_{m-1})\in \image(g)$ above, and symmetrically if  $(\beta_{m-1},\gamma_{m-1}+1)\in \image(g')$ we use a similar argument to define $(\alpha_m,\beta_m,\gamma_m)$. Now suppose neither of these cases holds. We must have \begin{align*}(\alpha_{m-1},\beta_{m-1}+1)&\in \image(g)\\ (\beta_{m-1}+1,\gamma_{m-1})&\in \image(g').\end{align*} In this case we define
 $$(\alpha_m,\beta_m,\gamma_m)=(\alpha_{m-1},\beta_{m-1}+1, \gamma_{m-1}),$$ which clearly meets conditions (i)-(iii). To check condition (iv) suppose this condition fails so that we have $(s_{i-1},s_i,s_{i+1})\in p_{i-1,i+1}\image(f)$  such that   $m\leq s_{i-1}+s_i+s_{i+1}$ and $$(\alpha_{m-1},\beta_{m-1}+1, \gamma_{m-1})\nleq (s_{i-1},s_i,s_{i+1}).$$ Since by the induction hypothesis  $$(\alpha_{m-1},\beta_{m-1}, \gamma_{m-1})\leq (s_{i-1},s_i,s_{i+1})$$ we must have $s_i=\beta_{m-1},$ in which case either $s_{i-1}> \alpha_{m-1}$ or $s_{i+1}> \gamma_{m-1}$. Without loss of generality suppose  $s_{i-1}> \alpha_{m-1}$. Then by the $f$-validity Criterion  $p_{i-1,i}\image(f)\subseteq \image(g)$ thus $ (s_{i-1},s_i)=(s_{i-1},\beta_{m-1})\in \image(g)$. However, this pair is incomparable with $(\alpha_{m-1},\beta_{m-1}+1)\in \image(g)$, contradicting the fact that $\image(g)$ is totally ordered. This contradiction shows that (iv) must hold.
 
  Having constructed $\alpha,$  $\beta$ and $\gamma$ inductively, let $e$ be the map of with $\Ty(e)=\hat{d}_i\hat{d}_i$ that has the form $$e=|\id| \cdots |\alpha,\beta,\gamma|\id | \cdots |\id|.$$ Condition (i) ensures $e$ is monic while conditions (ii) and (iii) ensure that $\image(e)\sqsubseteq \image(\hat{d}^{i-1}_{|g|})$ and $\image(e)\sqsubseteq \image(\hat{d}^{i}_{|g'|})$ respectively. Condition (iv) ensures that $\image(f)\sqsubseteq \image(e)$. Thus since $e$ increases dimension by $2$, we have that $e \in (f/I)$, and by Proposition~\ref{thetafactorprop} we have that $e$ factors through $\hat{d}^{i-1}_{|g|}$ and $\hat{d}^{i}_{|g'|}.$ Therefore $\hat{d}^{i-1}_{|g|} \smileeq e \smileeq \hat{d}^{i}_{|g'|}.$ This completes the last case, we have shown that for each $h,h' \in (f/I)$ we have $h\smileeq h'$, showing that $(f/I)$ is connected. 
\end{proof}
Moving from the horn case to the sphere case, let $c=\langle c_1, \ldots, c_{n}\rangle \in \Theta_2$ and let $\mathfrak{f}:c'\ra c$ be a coface map. The category of objects $O(c)$ of the universal sphere $d\Theta_2[c]$ is the category of morphisms $c' \ra c$ in $\Theta_2$ which factor through an object of dimension less than $\dim(c)$, and commuting triangles between them. The Glenn category $G(c)$ of $d\Theta_2[c]$ is the subcategory of $O(c)$ containing those $f:c'\ra c$ which are either a coface maps or a compositions of two coface maps, with morphisms between these given by precompositions by coface maps.
\begin{proposition} \label{cofinalprop2} A horn of the form $\Lambda_{\emptyset}^c=d\underline{c}$ is tabular, i.e the inclusion $I:G(c) \ra O(c)$ is cofinal.
\end{proposition} 
\begin{proof} The proof is the same as for Proposition~\ref{cofinalprop} except that we don't have to worry about $\frakf$, which simplifies the proof considerably.
\end{proof}
\section{Universal Glenn tables for $\Theta_2$-sets}
Propositions~\ref{cofinalprop} and \ref{cofinalprop2} allow us to use the Glenn tables to define horns and spheres in a $\Theta_2$-set. First we give the universal Glenn tables for the objects $c$ for which we want to make horns $\Lambda_{\{\mathfrak{f}\}}^{c}$ and $d\Theta_2[c],$ which is for $\dim(c)=3,4.$  Note that we will not endeavour to give an ordering of the cofaces of all objects in $\Theta_2$, but our preferred orderings for objects of dimension $3$ and $4$ are implicitly defined by the order of the universal Glenn tables for $3$ and $4$ dimensional objects in $\Theta_2$ which we list below. Note that the orderings defined below for the cofaces of $3$ dimensional objects are observed in the tables for $4$ dimensional objects. These tables are presented with greyed cells inserted in some of the rows. This is done to make all but one of these tables match the simplicial Glenn table pattern, making the pattern of these universal tables easier to remember.

To make these tables easier to read, we color $\langle 1 \rangle$ cells blue and $\langle 0, 0 \rangle$ cells yellow.  The symbol $\wedgedot$ is used to indicate row corresponds to an inner face, which can be removed to make an inner horn.

\begin{table}[H] \caption{Universal Glenn table for $\langle 2\rangle$ }\begin{center}
    \begin{tabular}{ r | l || l | l | }
   \hhline{~---}

                          &\onecell $|12|$       &   $|2|$   &     $|1|$           \\  \hhline{~---}

         $\wedgedot$        &\onecell $|02|$        & $|2|$  & $|0| $                      \\  \hhline{~---}
        
                              & \onecell   $|01|$    &  $|1|$   &   $|0|$             \\ \hhline{~---}
           
    \hhline{~---}
    \end{tabular}\end{center}
    \end{table}
    
    \begin{table}[H] \caption{Universal Glenn table for $\langle 1, 0\rangle$ }\begin{center}
        \begin{tabular}{ r | l || l | l | l|}
        \hhline{~----}

              &\zzcell $|1 | 0|$       &   $\star|0|$   &     $|1,0|$ &$|1|\star$          \\   \hhline{~----}

                  &\zzcell $|0|0|$        & $\star|0|$  & $|0,0| $      &$|0|\star$                \\   \hhline{~----}
            
   $\wedgedot$          &\onecell    $|01,00|$    &  $|1,0|$   &   $|0,0|$          &           \blank \\ \hhline{~----}
                                  
                & \onecell   $|01|\star$    &  $|1|\star$   &   $|0|\star$ &   \blank      \\ \hhline{~----}
               
        \end{tabular}\end{center}
        \end{table}
        
            \begin{table}[H] \caption{Universal Glenn table for $\langle 0, 1\rangle$ }\begin{center}
                \begin{tabular}{ r | l || l | l | l|}
             \hhline{~----}

                                      &\onecell $\star|01|$       &   \blank  &     $\star|1|$ &$\star|0|$          \\   \hhline{~----}

                          $\wedgedot$                &\onecell $|00,01|$        & \blank & $|0,1| $      &$|0,0|$                \\  \hhline{~----}
                    
                                          & \zzcell   $|0|1|$    &  $\star|1|$   &   $|0,1| $           &     $|0|\star$     \\ \hhline{~----}
                                          
                                       & \zzcell  $|0|0|$    & $\star|0|$   &   $|0,0|$ &   $|0|\star$       \\ \hhline{~----}
                       
                \end{tabular}\end{center}
                \end{table}
                
          \begin{table}[H] \caption{Universal Glenn table for $\langle 3 \rangle$ }\begin{center}
           \begin{tabular}{ r | l || l | l | l|}
                                  \hhline{~----}
                    &$|123|$       & \onecell   $|23|$  &   \onecell   $|13|$ &\onecell $|12|$          \\  \hhline{~----}
           
             $\wedgedot$        & $|023|$        & \onecell $|23|$    &\onecell  $|03|$      &\onecell $|02|$                \\  \hhline{~----}
            
            $\wedgedot$         &    $|013|$    & \onecell  $|13|$ &  \onecell  $|03|$            &  \onecell    $|01|$     \\ \hhline{~----}
                  
                   &    $|012|$    & \onecell $|12|$      &  \onecell $|02|$  &  \onecell  $|01|$       \\ \hhline{~----}
                                \end{tabular}\end{center}
                                \end{table}
          \begin{table}[H] \caption{Universal Glenn table for $\langle 2,0 \rangle$ }\begin{center}
           \begin{tabular}{ r | l || l | l | l|l|}
                                  \hhline{~-----}
                    &$|12|0|$       & \zzcell   $|2|0|$  &  \zzcell    $|1|0|$ & \onecell $|12,00|$  &  \onecell  $|12|\star$     \\  \hhline{~-----}
           
            $\wedgedot$         & $|02|0|$        &\zzcell  $|2|0|$    &\zzcell  $|0|0|$      &\onecell $|02,00|$   & \onecell  $|02|\star$           \\  \hhline{~-----}
            
                     &    $|01|0|$    &\zzcell  $|1|0|$      & \zzcell   $|0|0|$            &   \onecell   $|01,00|$& \onecell  $|01|\star$   \\ \hhline{~-----}
                  
        $\wedgedot$           &    $|012,000|$    & \onecell $|12,00|$      & \onecell  $|02,00|$  & \onecell $|01,00|$ &  \blank  \\ \hhline{~-----}
                   
                   &    $|012|\star$    &\onecell  $|12|\star$      & \onecell  $|02|\star$  &  \onecell  $|01|\star$    &  \blank  \\ \hhline{~-----}
                                \end{tabular}\end{center}
                                \end{table}
          \begin{table}[H] \caption{Universal Glenn table for $\langle 0,2 \rangle$ }\begin{center}
           \begin{tabular}{ r | l || l | l | l|l|}
                                \hhline{~-----}
                    &$\star|012|$       &   \blank  & \onecell     $\star|12|$ &\onecell $\star|02|$  &  \onecell  $\star|01|$     \\   \hhline{~-----}
           
           $\wedgedot$          & $|000,012|$        & \blank    &\onecell  $|00,12|$      &\onecell$|00,02|$   &  \onecell  $|00,01|$           \\  \hhline{~-----}
            
                     &    $|0|12|$    &\onecell   $\star|12|$     &\onecell   $|00,12|$            &   \zzcell   $|0|2|$& \zzcell  $|0|1|$   \\ \hhline{~-----}
                  
           $\wedgedot$        &    $|0|02|$    &\onecell $\star|02|$     & \onecell  $|00,02|$ &\zzcell  $|0|2|$ & \zzcell  $|0|0|$ \\   \hhline{~-----}
                   
                   &    $|0|01|$    &\onecell  $\star|01|$      &  \onecell  $|00,01|$    &  \zzcell  $|0|1|$    & \zzcell  $|0|0|$ \\  \hhline{~-----}
                                \end{tabular}\end{center}
                                \end{table}
          \begin{table}[H] \caption{Universal Glenn table for $\langle 1,1 \rangle$ }\begin{center}
           \begin{tabular}{ r | l || l | l | l|l|l|}
                                  \hhline{~------}
                    &$|1|01|$       &  \onecell  $ \star|01|$  &\onecell   $|11,01|$&\blank & \zzcell $|1|1|$  &   \zzcell $|1|0|$     \\   \hhline{~------}
           
                     & $|0|01|$        &\onecell  $\star|01|$    &\blank &\onecell  $|00,01|$      & \zzcell $|0|1|$   &  \zzcell $|0|0|$          \\    \hhline{~------}
            
         $\wedgedot$            &  $|011,001|$    & \onecell  $|11,01|$ &  \blank   & \onecell   $|01,01|$& \blank &\onecell  $|01,00|$  \\  \hhline{~------}
                  
       $\wedgedot$            &   $|001,011|$    &\blank    & \onecell  \ \circletext{$|01,11|$} \ &\onecell  $|01,01|$ & \onecell  \ \circletext{$|00,01|$} \   &\blank\\  \hhline{~------}
                   
                   &    $|01|1|$    & \zzcell $|1|1|$      &  \zzcell  $|0|1|$    &  \blank   & \onecell  $|01,11|$ &\onecell  $|01|\star$\\ \hhline{~------}
                  
                   &    $|01|0|$    &\zzcell  $|1|0|$      & \zzcell   $|0|0|$        &\onecell   $|01,00|$&   \blank &\onecell $|01|\star$\\  \hhline{~------}
                                \end{tabular}\end{center}
                                \end{table}
   Note that this table for $\langle 1, 1\rangle$ does not fit the usual simplicial Glenn Table pattern. But the pattern of this universal table can be remembered by noting that it would fit the usual pattern if the circled cells were swapped.                             
                                
          \begin{table}[H] \caption{Universal Glenn table for $\langle 1,0,0 \rangle$ }\begin{center}
           \begin{tabular}{ r | l || l | l | l|l|}
                                 \hhline{~-----}                            
                    &$|1|0|0|$       &  \zzcell   $\star|0|0|$  & \zzcell  $|1,0|0|$&\zzcell  $|1|0,0|$ &\zzcell $|1|0|\star$      \\   \hhline{~-----}
           
                     &$|0|0|0|$       &\zzcell  $\star|0|0|$   &\zzcell $|0,0|0|$ &\zzcell  $|0|0,0|$      &\zzcell $|0|0|\star$           \\    \hhline{~-----}
            
          $\wedgedot$           &    $|01,00|0|$    & \zzcell    $|1,0|0|$ &\zzcell  $|0,0|0|$  &   \onecell $|01,00,00|$&  \onecell $|01,00|\star$  \\ \hhline{~-----}
                  
          $\wedgedot$         &    $|01|0,0|$    &\zzcell $|1|0,0|$  &\zzcell  $|0|0,0|$   &   \onecell $|01,00,00|$ &  \onecell $|01|\star\star$  \\  \hhline{~-----}
                   
                   &    $|01|0|\star$    &\zzcell $|1|0|\star$    &\zzcell $ |0|0|\star$     & \onecell$|01,00|\star$  &  \onecell $|01|\star\star$\\  \hhline{~-----}

                                \end{tabular}\end{center}
                                \end{table}
          \begin{table}[H] \caption{Universal Glenn table for $\langle 0,1,0 \rangle$ }\begin{center}
           \begin{tabular}{ r | l || l | l | l|l|l|}
                                \hhline{~------}
                    &$\star|01|0|$       &    \blank  &\zzcell   $\star|1|0|$&\zzcell  $\star|0|0|$ & \onecell $\star|01,00|$ &  \onecell$\star|01|\star$   \\  \hhline{~------}
           
            $\wedgedot$         &$|00,01|0|$       & \blank  &\zzcell $|0,1|0|$ & \zzcell $|0,0|0|$    & \onecell $|00,01,00|$   &  \onecell $|00,01|\star$     \\  \hhline{~------}
            
                     &    $|0|1|0|$    &  \zzcell   $\star |1|0|$ &\zzcell  $|0,1|0|$  &   \blank & \zzcell  $|0|1,0|$ &\zzcell  $|0|1|\star$ \\ \hhline{~------}
                  
                   &    $|0|0|0|$    &\zzcell $\star|0|0|$ &\zzcell  $|0,0|0|$   &  \blank & \zzcell  $|0|0,0|$&\zzcell  $|0|0|\star$ \\ \hhline{~------}
                   
           $\wedgedot$        &    $|0|01,00|$    &  \onecell $\star|01,00|$    & \onecell $|00,01,00|$    &\zzcell  $|0|1,0|$  & \zzcell  $|0|0,0|$& \blank \\ \hhline{~------}
                   &    $|0|01|\star$    & \onecell $\star|01|\star$    & \onecell $|00,01|\star$    &\zzcell  $|0|1|\star$  & \zzcell  $|0|0|\star$& \blank \\ \hhline{~------}
                                \end{tabular}\end{center}
                                \end{table}
          \begin{table}[H] \caption{Universal Glenn table for $\langle 0,0,1 \rangle$ }\begin{center}
           \begin{tabular}{ r | l || l | l | l|l|}
                              \hhline{~-----}
                    &$\star|0|01|$       &  \onecell  $\star\star|01|$  &  \onecell $\star|00,01|$& \zzcell $\star|0|1|$ & \zzcell $\star|0|0|$      \\   \hhline{~-----}
           
                    $\wedgedot$ &$|0,0|01|$       & \onecell $\star\star|01|$   & \onecell $|00,00,01|$ &\zzcell  $|0,0|1|$      &\zzcell $|0,0|0|$           \\    \hhline{~-----}
            
                    $\wedgedot$ &    $|0|00,01|$    &  \onecell   $\star|00,01|$ & \onecell $|00,00,01|$  &   \zzcell $|0|0,1|$& \zzcell  $|0|0,0|$  \\  \hhline{~-----}
                  
                   &    $|0|0|1|$    &\zzcell $\star|0|1|$  &\zzcell $|0,0|1|$    &  \zzcell  $|0|0,1|$ & \zzcell  $|0|0|\star$  \\  \hhline{~-----}
                   
                   &    $|0|0|0|$    &\zzcell $\star|0|0|$    &\zzcell $|0,0|0|$       &\zzcell $|0|0,0|$  & \zzcell  $|0|0|\star$\\  \hhline{~-----}

                                \end{tabular}\end{center}
                                \end{table}

\section{The structure of $\TBic(X)$}
There is a full and faithful functor $\psi:\Delta \ra \Theta_2$ given by sending $[n]$ to $$\langle \underbrace{0,0,\ldots, 0}_\text{$n$ times} \rangle.$$ This induces a functor $\psi^\star:\cat{Set}_{\Theta_2}\ra \cat{Set}_{\Delta}.$ It is easy to see that if $X$ is ($2$-reduced) inner-Kan, then so is $\psi^\star(X).$

Let $(X,\chi)$ be an algebraic $2$-reduced inner-Kan $\Theta_2$-set. Note that the only universal inner horn of minimal complementary dimension $1$ or lower is of type $\langle 0,0\rangle$. 

We define a fancy bicategory \footnote{See Definition~\ref{fancydef}}  $ \TBic(X):$
\begin{definition}
The objects and $1$-morphisms of $\widetilde{\TBic(X)}$ are the $\langle\ \rangle$-cells and $\langle 0 \rangle$ cells of $X$, which are the same as the objects and $1$-morphisms of $\overline{\TBic(X)}:=\Bic(\psi^{\star}(X)).$ A morphism from $f$ to $g$ in $\widetilde{\TBic(X)}$ is a $\langle 1 \rangle$ cell $\eta$ in $X$ such that $\theta_{|0|} \eta=f$ and $\theta_{|1|}\eta = g.$ 
\end{definition}
\begin{definition}[Identity for $1$-morphisms in $\widetilde{\TBic(X)}$] For a $1$-morphism $f$ of $\widetilde{\TBic(X)}$ the identity for $f$, which we denote $\ID_f$, is defined by $\ID_f:=\theta_{|00|}f$.
\end{definition}
\begin{definition}[$\bullet$ in $\widetilde{\TBic(X)}$]Let $\eta:f \Rightarrow g$ and $\theta:g\Rightarrow h$ be two $2$-morphisms in $\widetilde{\TBic(X)}.$ The following $\langle 2\rangle $-horn defines $\theta\bullet \eta$ in $\widetilde{\TBic(X)}$. 
\begin{table}[H] \caption{The $\langle 2\rangle $-horn $\hat{\Theta}_{\bullet}(\theta,\eta)$ defining $\theta\bullet \eta$ }\begin{center}
    \begin{tabular}{ r | l || l | l | }
\hhline{~---}

                          &\onecell $\theta$       &   $h$   &     $g$           \\  \hhline{~---}

                      $\Lambda$        &\onecell $=:\theta\bullet \eta$        & $h$  & $f $                      \\  \hhline{~---}
        
                              &  \onecell  $\eta$    &  $g$   &   $f$             \\ \hhline{~---}
           
  \hhline{~---}
    \end{tabular}\end{center}
    \end{table}

\end{definition}

\begin{definition}[Right whiskering $\rhd$ in $\widetilde{\TBic(X)}$]Let $\eta:f \Rightarrow g:a\ra b$ be a $2$-morphism and $h:b\ra c$ be a $1$-morphism in $\widetilde{\TBic(X)}$. The following $\langle 1,0\rangle $-horn defines $h\rhd \eta.$ 

    \begin{table}[H] \caption{The $\langle 1,0\rangle $-horn $\hat{\Theta}_{\rhd}(h,\eta)$ defining $h\rhd \eta$ }\begin{center}
        \begin{tabular}{ r | l || l | l | l|}
\hhline{~----}
                              &\zzcell $\chi(h,g)$       &   $h$   &     $h\circ g$ &$g$          \\  \hhline{~----}

                                  &\zzcell $\chi(h,f)$        & $h$  & $h\circ f$      &$f$                \\   \hhline{~----}
            
                        $\Lambda$            &\onecell  $=:h\rhd\eta$    &    $h\circ g$    &  $h\circ f$ &           \blank \\ \hhline{~----}
                                  
                               &  \onecell  $\eta$    &  $g$   &   $f$ &   \blank      \\ \hhline{~----}
        \end{tabular}\end{center}
        \end{table}
\end{definition}
\begin{definition}[Left whiskering $\lhd$ in $\widetilde{\TBic(X)}$]Let $\eta:g \Rightarrow h:b\ra c$ be a $2$-morphism  and $f:a\ra b$ be a $1$-morphism in $\widetilde{\TBic(X)}$. The following $\langle 0,1\rangle $-horn defines $\eta \lhd f.$

            \begin{table}[H] \caption{The $\langle 0,1\rangle $-horn $\hat{\Theta}_{\lhd}(\eta,f)$ defining $\eta\lhd f$ }\begin{center}
                \begin{tabular}{ r | l || l | l | l|}
                \hhline{~----}

                                      &$\eta$ \onecell      &   \blank  &     $h$ &$g$          \\  \hhline{~----}

                          $\Lambda$                & $=:\eta\lhd f$ \onecell       &\blank  & $h\circ f$     &$g\circ f$                \\   \hhline{~----}
                    
                                          &    $\chi(h,f)$ \zzcell  &  $h$   &   $h\circ f$          &     $f$     \\ \hhline{~----}
                                          
                                       &    $\chi(g,f)$ \zzcell   & $g$   &  $g\circ f$   &   $f$       \\ \hhline{~----}
                       
                \end{tabular}\end{center}
                \end{table}
                \end{definition}
\begin{definition}[Thin structure map $t$ ]Let $\beta:f \Rightarrow g:a\ra b$ be a $2$-morphism of $\overline{\TBic(X)}$. The following $\langle 0,1\rangle $-horn defines $t(\beta).$ 
            \begin{table}[H] \caption{The $\langle 0,1\rangle $-horn $ \hat{\Theta}_{t}(\beta)$ defining $t(\beta)$}\begin{center}
                \begin{tabular}{ r | l || l | l | l|}
                 \hhline{~----}

                                      &$\ID_g$  \onecell     &   \blank  &     $g$ &$g$          \\   \hhline{~----}

                          $\Lambda$                & $t(\beta)$  \onecell      &\blank  & $g$     &$f$                \\   \hhline{~----}
                    
                                          &    $s_0 g$  \zzcell  &  $g$   &   $g$          &     $\id_a$     \\ \hhline{~----}
                                          
                                       &    $\beta$ \zzcell   & $g$   &  $f$   &   $\id_a$       \\ \hhline{~----}
                       
                \end{tabular}\end{center}
                \end{table}
                \end{definition}
\begin{definition}
The unitors and associators of $\widetilde{\TBic(X)}$ are defined by $t(\rho_f)$, $t(\lambda_f)$, and $t(\alpha_{h,g,f})$ respectively.
\end{definition}
\section{Verification that $\TBic(X)$ is a fancy bicategory}
 Let $X$ be a $2$-reduces inner-Kan $\Theta_2$-set. As in Chapters~\ref{bicchapter} and \ref{vdcchapter}, the following analogue of Lemma~\ref{lemma1} is the main tool we use for checking the fancy bicategory axioms for $\TBic(X):$
\begin{lemma}[Matching Lemma] Suppose we have two commutative $d\Theta[c]$-spheres $S$, $S'$ in $X$, such that every corresponding face, except for possibly a single inner face, of these spheres match. Then $S=S'$ and in particular they indeed match on this putatively non-matching inner face.
\end{lemma}
\begin{proof} By the uniqueness of fillers for inner horns in $X$.
\end{proof}

\subsection{Axioms for vertical composition in $\widetilde{\TBic(X)}.$}
\begin{proposition}[\textbf{B1}, identity for $1$-morphisms] For all $f,g$ and $\eta:f\Rightarrow g$ in $\widetilde{\TBic(X)},$ we have $\eta \bullet \ID_f =\ID_g \bullet \eta = \eta$.
\end{proposition}
\begin{proof}\begin{align*}
d(\Theta_{\bullet}(\eta,\ID_f))&=[\eta,\ \eta\bullet\ID_f,\ \ID_f]\\
d(\theta_{001}\eta)&=[\eta,\ \eta,\ \ID_f] \\
d(\Theta_{\bullet}(\ID_g,\eta))&=[\ID_g,\ \ID_g\bullet\eta,\ \eta] \\
d(\theta_{011}\eta)&=[\ID_g,\ \eta,\ \eta]
\end{align*}
So by the Matching Lemma $\eta \bullet \ID_f =\eta$ and $\ID_g \bullet \eta = \eta.$
\end{proof}
\begin{proposition}[\textbf{B2} associativity of vertical composition] For all $f \stackrel{\eta}{\Rightarrow} g \stackrel{\theta}{\Rightarrow} h \stackrel{\iota}{\Rightarrow}i,$ in $\widetilde{\TBic(X)}$, we have $\iota \bullet (\theta \bullet \eta) = (\iota \bullet \theta) \bullet \eta.$
\end{proposition}
\begin{proof}By the Matching Lemma applied to $\Theta_{\bullet}(\iota, \theta \bullet \eta)$ it suffices to show the following $\langle 2 \rangle$ sphere is commutative: $$[\iota, (\iota\bullet \theta)\bullet \eta, \theta \bullet \eta].$$ The following $\langle 3\rangle$-horn in $\widetilde{\TBic(X)}$ verifies this commutativity:
          \begin{table}[H] \caption{ $\langle 3 \rangle$-horn}\begin{center}
           \begin{tabular}{ r | l || l | l | l|}
                                  \hhline{~----}
                    &$\Theta_{\bullet}(\iota,\theta)$       & \onecell   $\iota$  &   \onecell   $\iota\bullet\theta$ &\onecell $\theta$          \\  \hhline{~----}
           
             $\Lambda$        &  & \onecell $\iota$    &\onecell  $(\iota\bullet \theta \bullet) \eta  $    &\onecell $\theta \bullet \eta$                \\  \hhline{~----}
            
                 &  $\Theta_{\bullet}(\iota\bullet\theta,\eta) $    & \onecell  $\iota\bullet \theta$ &  \onecell   $(\iota\bullet \theta \bullet) \eta$          &  \onecell    $\eta$     \\ \hhline{~----}
                  
                   &    $\Theta_{\bullet}(\iota\bullet \theta, \eta)$    & \onecell $\theta$      &  \onecell $\theta\bullet\eta$  &  \onecell  $\eta$       \\ \hhline{~----}
                                \end{tabular}\end{center}
                                \end{table}
\end{proof}
\begin{definition}[Thin structure map inverse $t^{-1}$ ]Let $\beta:f \Rightarrow g:a\ra b$ be a $2$-morphism of $\overline{\TBic(X)}$. The following $\langle 0,1\rangle $-horn defines $t^{-1}(\beta).$ 
            \begin{table}[H] \caption{The $\langle 0,1\rangle $-horn $ \hat{\Theta}_{t^{-1}}(\beta)$ defining $t^{-1}(\beta)$}\begin{center}
                \begin{tabular}{ r | l || l | l | l|}
    \hhline{~----}
                          &$\ID_g$ \onecell      &   \blank  &     $g$ &$g$          \\  \hhline{~----}

                          $\Lambda$                &\onecell $=:t^{-1}(\beta)$        &\blank  & $f$     &$g$                \\   \hhline{~----}
                    
                                          &  \zzcell $\beta$    &  $g$   &   $f$          &     $\id_a$     \\ \hhline{~----}
                                          
                                       &    \zzcell  $s_0 g$   & $g$   &  $g$   &   $\id_a$       \\ \hhline{~----}
                       
                \end{tabular}\end{center}
                \end{table}
                \end{definition}
\subsection{Strict functoriality of $t$ and the $\langle 0,1\rangle$ and $\langle 1,0\rangle$ cell criteria.}
\begin{proposition}[$t$ preserves identities for $1$-morphisms] Let $f:a\ra b$ be a $1$-morphism of $\TBic(X).$ Then $t(\Id_g)=\ID_g$
\end{proposition}
\begin{proof}We have \begin{align*}d(\theta_{|01|\cdot|}(\ID_g))&=[\ID_g,\ \ID_g,\ s_0 g,\ s_0 g=\Id_g] \\ d(\Theta_{t}(\Id_g))&=[\ID_g,\ t(\Id_g),\ s_0 g,\ \Id_g]. \end{align*} The Matching Lemma then shows $t(\Id_g)=\ID_g.$
\end{proof}
\begin{proposition}\label{inversesphere} Let $\beta:f \Rightarrow g:a\ra b$ be a $2$-morphism of $\overline{\TBic(X)}$. Then $t(\beta)$ is invertible and $t(\beta)^{-1}=t^{-1}(\beta).$
\end{proposition}
\begin{proof} By the Matching Lemma applied to $\Delta_{\bullet}(t(\beta),t^{-1}(\beta))$ it suffices to show the following $d\Theta_2\langle 0,0,0 \rangle$-sphere\footnote{$\Theta_2\langle 0,0,0 \rangle$ is an abbreviation for the representable presheaf $\Theta_2[\langle 0,0,0 \rangle]$} is commutative $$[t(\beta),\ \Id_g ,\ t^{-1}(\beta)].$$ The following horn verifies this commutativity:
          \begin{table}[H] \caption{ $\langle 0,2 \rangle$-horn}\begin{center}
           \begin{tabular}{ r | l || l | l | l|l|}
                                   \hhline{~-----}
                             
                    &$\theta_{|000|} g$       &   \blank  &   \onecell   $\ID_g$ & \onecell  $\ID_g$ &   \onecell   $\ID_g$  \\   \hhline{~-----}
           
           $\Lambda$       &     & \blank    & \onecell $t(\beta)$      & \onecell $\ID_g$   & \onecell  $ t^{-1}(\beta)$           \\    \hhline{~-----}
            
                     &   $\Theta_{t}(\beta)$   &  \onecell $\ID_g$     &  \onecell $t(\beta)$            & \zzcell     $\Id_g$& \zzcell  $\beta$   \\  \hhline{~-----}
                  
                          & $ \theta_{|00|\cdot|}(g)$   & \onecell $\ID_g$     & \onecell  $\ID_g$ &\zzcell  $\Id_g$ & \zzcell  $\Id_g$ \\  \hhline{~-----}
                   
                   &    $\Theta_{t^{-1}}(\beta)$     & \onecell $\ID_g$      &  \onecell $  t^{-1}(\beta)$      &\zzcell   $\beta$     &\zzcell   $\Id_g$ \\  \hhline{~-----}
                                \end{tabular}\end{center}
                                \end{table}
\end{proof}
Because of Proposition~\ref{inversesphere}, we will make no further use of the notation $t^{-1}(\beta)$, but we will use $\Theta_{t^{-1}}(\beta)$ frequently.
\begin{proposition}Let $\beta:f \Rightarrow g:a\ra b$ be a $2$-morphism of $\overline{\TBic(X)}$. Then $t(\beta^{-1})=t(\beta)^{-1}.$
\end{proposition}
\begin{proof} By the Matching Lemma applied to $\Theta_{t}(\beta^{-1})$, it suffices to show the $d\Theta_2\langle 0,1\rangle$-sphere $$[\ID_f,\ t(\beta)^{-1},\ \Id_f,\ \beta^{-1}]$$ is commutative. The following $d\Theta_2\langle 0, 0, 1\rangle $-sphere demonstrates this commutativity: 
          \begin{table}[H] \caption{ $\langle 0,0,1 \rangle$-horn }\begin{center}
           \begin{tabular}{ r | l || l | l | l|l|}
                                   \hhline{~-----}

                    &$\theta_{|00|0|}(\beta)$       &   \onecell  $\ID_g$ &  \onecell $\ID_f$&\zzcell  $\beta$  &\zzcell $\beta$      \\   \hhline{~-----}
           
                    & $\Theta_{t^{-1}}(\beta)$    & \onecell    $\ID_g$ & \onecell $ t(\beta)^{-1}$  &\zzcell    $\beta$& \zzcell  $\Id_g$  \\ \hhline{~-----}
                    
            $\Lambda$   &             & \onecell $\ID_f$   & \onecell$ t(\beta)^{-1}$ &\zzcell  $\Id_f$      &\zzcell $\beta^{-1}$           \\   \hhline{~-----}
                   &    $s_0\beta$    &\zzcell $\beta$  &\zzcell $\beta$    &  \zzcell $\Id_f$ & \zzcell$\Id_{\id_a}$ \\  \hhline{~-----}
                   
                   &    $\Delta_{\bullet}(\beta,\beta^{-1})$    &\zzcell$\beta$    &\zzcell$\Id_g$       &\zzcell$\beta^{-1}$  & \zzcell $\Id_{\id_a}$\\ \hhline{~-----}

                                \end{tabular}\end{center}
                                \end{table}
\end{proof}
\begin{proposition}[$t$ preserves $\bullet$]Let $\beta:f \Rightarrow g:a\ra b$ and $\gamma:g \Rightarrow h:a\ra b$  be a $2$-morphisms of $\overline{\TBic(X)}$. Then $t(\gamma \bullet \beta )=t(\gamma)\bullet t(\beta).$
\end{proposition}
\begin{proof} By the Matching Lemma applied to $\Theta_{\bullet}(t(\gamma),t(\beta))$, it suffices to show the $d\Theta_2\langle 2 \rangle$-sphere $$[t(\gamma), t(\gamma\bullet \beta), t(\beta)]$$ is commutative. The following Glenn table proof demonstrates this commutativity: 
          \begin{table}[H] \caption{$\langle 0,2 \rangle$-horn}\begin{center}
           \begin{tabular}{ r | l || l | l | l|l|}
                                   \hhline{~-----}

                    &$\theta_{|000|} g$       &   \blank  &  \onecell $\ID_g$ & \onecell $\ID_g$ &  \onecell   $\ID_g$    \\   \hhline{~-----}
           
               $\Lambda$      &                     & \blank    & \onecell $t(\gamma)$      & \onecell $t(\gamma\bullet \beta)$   & \onecell $t(\beta)$           \\   \hhline{~-----}
            
                     &    $\Theta_{t^{-1}}(\gamma^{-1})$    & \onecell $\ID_g$  & \onecell $t(\gamma)=t(\gamma^{-1})^{-1}$           & \zzcell    $\gamma^{-1}$&\zzcell  $\Id_g$   \\  \hhline{~-----}
                  
           $\odot$        &   (Table~\ref{thomproof2})    &  \onecell  $\ID_g$    &  \onecell  $t(\gamma\bullet \beta)$  &\zzcell $\gamma^{-1}$ & \zzcell $\beta$ \\ \hhline{~-----}
                   
                   &    $\Theta_{t}(\beta)$   &  \onecell  $\ID_g$      &   \onecell  $t(\beta)$     &   \zzcell $\Id_g$    & \zzcell $\beta$ \\  \hhline{~-----}
                                \end{tabular}\end{center}
                                \end{table}
          \begin{table}[H] \caption{$\langle 0,0,1 \rangle$-horn \label{thomproof2} }\begin{center}
           \begin{tabular}{ r | l || l | l | l|l|}
                         \hhline{~-----}

                    &    $\theta_{|00|0|}(\gamma)$    & \onecell    $\ID_h$ & \onecell  $\ID_g$ &\zzcell $\gamma$   &\zzcell $\gamma$        \\  \hhline{~-----}
           
                     &$\Theta_{t}(\gamma\bullet\beta)$  & \onecell $\ID_h$   &\onecell  $t(\gamma\bullet \beta)$  &\zzcell  $\Id_h$      &\zzcell $\gamma\bullet \beta$          \\   \hhline{~-----}
            
                    $\Lambda$ &                & \onecell $\ID_g$    & \onecell  $t(\gamma\bullet \beta)$  & \zzcell  $\gamma^{-1}$ & \zzcell  $\beta$  \\\hhline{~-----}
                  
                   &    $\Delta_{\bullet}(\gamma,\gamma^{-1})$     &\zzcell $\gamma$    & \zzcell $\Id_h$    &\zzcell  $\gamma^{-1}$  & \zzcell  $\Id_{\id_a}$  \\ \hhline{~-----}
                   
                   &    $\Delta_{\bullet}(\gamma,\beta)$    &\zzcell $\gamma$    &\zzcell $\gamma\bullet \beta$       &\zzcell  $\beta$   &\zzcell   $\Id_{\id_a}$\\ \hhline{~-----}

                                \end{tabular}\end{center}
                                \end{table}
\end{proof}
\begin{lemma}[Commutativity criterion for $d\Theta_2\langle 0, 1 \rangle$-spheres]
\label{02celllemma} Let $X$ be a $2$-reduced inner-Kan $\Theta_2$-set. Let $$S=[\theta,\ \eta,\ \gamma ,\ \beta]$$ be a $d\Theta_2\langle 0, 1 \rangle$-sphere in $X$, with the following Glenn table:
            \begin{table}[H] \caption{The $d\Theta_2\langle 0, 1\rangle$-sphere $S$}\begin{center}
                \begin{tabular}{ r | l || l | l | l|}
                 \hhline{~----}

                                &$\theta$  \onecell     &   \blank  &     $h'$ &$h$          \\  \hhline{~----}

                       & $\eta$\onecell        &\blank  & $g'$     &$g$    \\   \hhline{~----}
                    
                                    &  \zzcell  $\gamma$    &  $h'$ & $g'$ &     $f$     \\ \hhline{~----}
                                          
                                    &  \zzcell  $\beta$    & $h$   &  $g$   &   $f$       \\  \hhline{~----}
                       
                \end{tabular}\end{center}
                \end{table}
 Then $S$ is commutative if and only if $\eta=t(\underline{\gamma})^{-1}\bullet (\theta \lhd f) \bullet t(\underline{\beta}).$ 
\end{lemma}
\begin{proof} By the uniqueness of fillers for inner $\langle 0, 1 \rangle$-horns in $X$, it is enough to show  $$S=[\theta,\ t(\underline{\gamma})^{-1}\bullet (\theta \lhd f) \bullet t(\underline{\beta}),\ \gamma ,\ \beta]$$ is commutative in $X$. The following Glenn table proof demonstrates this commutativity.
          \begin{table}[H] \caption{$\langle 0,2 \rangle$-horn \label{3lemmaproof1} }\begin{center}
          \scalebox{.95}{ \begin{tabular}{ r | l || l | l | l|l|}
                                \hhline{~-----}
                    &$\theta_{|001|}(\theta)$       &   \blank  & \onecell    $\theta$   &\onecell $\theta$   &  \onecell  $\ID_h$     \\   \hhline{~-----}
           
                     & $\Theta_{\bullet}(t(\underline{\gamma})^{-1}\bullet (\theta \lhd f) ,  t(\underline{\beta}))$        & \blank    &\onecell  $t(\underline{\gamma})^{-1}\bullet (\theta \lhd f)$      &\onecell$t(\underline{\gamma})^{-1}\bullet (\theta \lhd f) \bullet t(\underline{\beta})$   &  \onecell  $ t(\underline{\beta})$           \\  \hhline{~-----}
            
                 $\odot_2$     &  (Table~\ref{3lemmaproof3})     &\onecell   $\theta$       &\onecell   $t(\underline{\gamma})^{-1}\bullet (\theta \lhd f)$          &   \zzcell   $\gamma$& \zzcell  $\chi(h,f)$   \\ \hhline{~-----}
                  
           $\Lambda$        &       &\onecell $\theta$     & \onecell $t(\underline{\gamma})^{-1}\bullet (\theta \lhd f) \bullet t(\underline{\beta})$ &\zzcell  $\gamma$ & \zzcell  $\beta$ \\   \hhline{~-----}
                   
                 $\odot_1$  &  (Table~\ref{3lemmaproof2})        &\onecell  $\ID_h$      &  \onecell  $t(\underline{\beta})$    &  \zzcell  $\chi(h,f)$    & \zzcell  $\beta$ \\  \hhline{~-----}
                                \end{tabular}}\end{center}\end{table}
          \begin{table}[H] \caption{ $\langle 0,0,1 \rangle$-horn showing $\odot_1$ from Table~\ref{3lemmaproof1} \label{3lemmaproof2} }\begin{center}
           \begin{tabular}{ r | l || l | l | l|l|}
                              \hhline{~-----}
                    &$\theta_{|00|0|}(\chi(h,f))$       &  \onecell  $\ID_{h}$ &  \onecell $\ID_{h\circ f}$& \zzcell $\chi(h,f)$  & \zzcell $\chi(h,f)$       \\   \hhline{~-----}
           
                    $\Lambda$ &       &\onecell  $\ID_h$      &  \onecell  $t(\underline{\beta})$    &  \zzcell  $\chi(h,f)$    & \zzcell  $\beta$     \\    \hhline{~-----}
            
                    &    $\Theta_{t}(\underline{\beta})$    &  \onecell  $\ID_{h\circ f}$ & \onecell $t(\underline{\beta})$   &   \zzcell $\Id_{h\circ f}$& \zzcell  $\underline{\beta}$  \\  \hhline{~-----}
                  
                   &    $s_0 \chi(h,f)$    &\zzcell  $\chi(h,f)$  &\zzcell $\chi(h,f)$   &  \zzcell $\Id_{h\circ f}$ & \zzcell  $\Id_{f}$ \\  \hhline{~-----}
                   
                   &    $\Delta_{-}(\beta)$    &\zzcell $\chi(h,f)$    &\zzcell $\beta$        &\zzcell $\underline{\beta}$  & \zzcell  $\Id_{f}$ \\  \hhline{~-----}

                                \end{tabular}\end{center}
                                \end{table}
          \begin{table}[H] \caption{$\langle 0,2 \rangle$-horn showing $\odot_2$ from Table~\ref{3lemmaproof1} \label{3lemmaproof3} }\begin{center}
          \scalebox{1}{ \begin{tabular}{ r | l || l | l | l|l|}
                                \hhline{~-----}
                    &$\theta_{|011|}(\theta)$       &   \blank  & \onecell    $\ID_h'$   &\onecell $\theta$   &  \onecell  $ \theta $     \\   \hhline{~-----}
           
                     & $\Theta_{\bullet}(t(\underline{\gamma})^{-1} ,  \theta \lhd f)$        & \blank    &\onecell  $t(\underline{\gamma})^{-1} $      &\onecell$t(\underline{\gamma})^{-1}\bullet (\theta \lhd f)$    &  \onecell  $ \theta \lhd f$           \\  \hhline{~-----}
                     
                $\odot$      & (Table~\ref{3lemmaproof4})   &\onecell   $\ID_{h'}$       &\onecell   $ t(\underline{\gamma})^{-1}$         &   \zzcell   $\gamma$& \zzcell  $\chi(h',f')$  \\   \hhline{~-----}
                
               $\Lambda$    &      &\onecell   $\theta$       &\onecell   $t(\underline{\gamma})^{-1}\bullet (\theta \lhd f)$          &   \zzcell   $\gamma$& \zzcell  $\chi(h,f)$   \\ \hhline{~-----}
                  
                 & $\Theta_{\lhd}(\theta,f)$      &\onecell  $ \theta $      &  \onecell $ \theta \lhd f$     &  \zzcell  $\chi(h',f')$    & \zzcell  $\chi(h,f)$ \\  \hhline{~-----}
                                \end{tabular}}\end{center}\end{table}
               
          \begin{table}[H] \caption{ $\langle 0,0,1 \rangle$-horn showing $\odot$ from Table~\ref{3lemmaproof3} \label{3lemmaproof4} }\begin{center}
           \begin{tabular}{ r | l || l | l | l|l|}
                              \hhline{~-----}
                    &$\theta_{|00|0|}(\chi(h',f'))$       &  \onecell  $\ID_{h'}$ &  \onecell $\ID_{h'\circ f'}$& \zzcell $\chi(h',f')$  & \zzcell $\chi(h',f')$       \\   \hhline{~-----}
           
                    $\Lambda$ &      &\onecell   $\ID_{h'}$       &\onecell   $ t(\underline{\gamma})^{-1}$         &   \zzcell   $\gamma$& \zzcell  $\chi(h',f')$  \\   \hhline{~-----}
            
                    &    $\Theta_{t^{-1}}(\underline{\gamma})$    &  \onecell  $\ID_{h'\circ f'}$ & \onecell $ t(\underline{\gamma})^{-1}$   &   \zzcell $\underline{\gamma}$& \zzcell 
                      $\Id_{h'\circ f'}$ \\  \hhline{~-----}
                  
                   &      $\Delta_{-}(\gamma)$  &\zzcell $\chi(h',f')$  &\zzcell $\gamma$   &  \zzcell $\underline{\gamma}$  & \zzcell  $\Id_{f'}$ \\  \hhline{~-----}
                   
                   &    $s_0 \chi(h',f')$    &\zzcell $\chi(h',f')$   &\zzcell$\chi(h',f')$       &\zzcell  $\Id_{h'\circ f'}$  & \zzcell  $\Id_{f'}$ \\  \hhline{~-----}

                                \end{tabular}\end{center}
                                \end{table}
\end{proof}
The following verification of the naturality of $\rho$ is a bit out of place, but we need it to continue our discussion of commutativity criteria.
\begin{proposition}[\textbf{B7} Naturality of $\rho$] \label{natofrho} Let $\eta:f\Rightarrow g:a\ra b$ in $\widetilde{\TBic(X)}$. Then $$t(\rho_g)^{-1}\bullet (\eta\lhd \id_a )\bullet t(\rho_f)= \eta.$$
\end{proposition}
\begin{proof} Apply Proposition~\ref{02celllemma} to the sphere $$d(\theta_{|01|\cdot|}(\eta))= [\eta,\ \eta,\ \Id_f,\ \Id_g]$$ yielding $$t(\underline{\Id_g})^{-1}\bullet (\eta\lhd \id_a )\bullet t(\underline{\Id_f})= \eta.$$ We are done once we recall that $\underline{\Id_g}=\rho_g$  and $\underline{\Id_f}=\rho_f$ by definition.\end{proof}
\begin{corollary}
\label{02cellcor} Let $X$ be a $2$-reduced inner-Kan $\Theta_2$-set. Let $$S=[\theta,\ \eta,\ \gamma ,\ \beta]$$ be a $d\Theta_2\langle 0, 1 \rangle$-sphere in $X$ as in Lemma~\ref{02celllemma} except with the additional condition that $\gamma$ and $\beta$  are $2$-morphisms of $\overline{\TBic(X)}.$ Then $S$ is commutative if and only if $\theta=t(\gamma)^{-1}\bullet \eta  \bullet t(\beta).$ 
\end{corollary}
\begin{proof}By Lemma~\ref{02celllemma} $S$ is commutative if and only if \begin{equation}\label{02cellcoreq} \eta=t(\underline{\gamma})^{-1}\bullet (\theta \lhd \id) \bullet t(\underline{\beta}).\end{equation} By Lemma~\ref{bulletylemma}, $$\underline{\gamma}=\underline{\Id_{h'}\bullet \gamma}=\underline{\Id_{h'}}\bullet \gamma=\rho_{h'}\bullet \gamma.$$ So Equation~\ref{02cellcoreq} is equivalent to $$\theta=t(\rho_{h'} \bullet \gamma )^{-1}\bullet (\theta \lhd \id) \bullet t(\rho_{h}\bullet \beta).$$ The naturality of $\rho$ as shown in Proposition~\ref{natofrho} shows this is equivalent to $\eta=t(\gamma)^{-1}\bullet \theta  \bullet t(\beta)$ as was to be shown.
\end{proof}
\begin{lemma}
\label{10cellcor} Let $X$ be a $2$-reduced inner-Kan $\Theta_2$-set. Let $$S=[\theta,\ \eta,\ \gamma ,\ \beta]$$ be a $d\Theta_2\langle 0, 1 \rangle$-sphere in $X$ as in Corollary~\ref{02cellcor}, with $\gamma$ and $\beta$ being $2$-morphisms in $\overline{\TBic(X)}.$ Then the $d\Theta_2\langle 2\rangle$-sphere $$\widehat{S}:=[ \widehat{\gamma} ,\ \widehat{\beta},\ \eta,\ \theta]$$ is commutative if and only if $S$ is commutative, which is equivalent by Corollary~\ref{02cellcor} to the condition  $$\eta=t(\gamma)^{-1}\bullet \theta  \bullet t(\beta).$$
\end{lemma}
\begin{proof}       Consider the following Glenn Table:
   \begin{table}[H] \caption{ $\langle 0,1,0 \rangle$-horns \label{SShattable} }\begin{center}
           \begin{tabular}{ r | l || l | l | l|l|l|}
                                \hhline{~------}
                    &$\theta_{|\cdot|01|}(\theta)$       &    \blank  &\zzcell   $\widehat{\Id_{h'}}$&\zzcell $\widehat{\Id_{h}}$& \onecell $\theta$  &  \onecell$\theta$   \\  \hhline{~------}
           
            $\Lambda_0$         & $\widehat{S}$ & \blank  &\zzcell  $\widehat{\gamma}$ & \zzcell $\widehat{\beta}$    & \onecell $\eta$  &  \onecell $\theta$    \\  \hhline{~------}
            
                     &    $\Delta_{\wedge}(\gamma)$    &  \zzcell   $\widehat{\Id_{h'}}$ &\zzcell  $\widehat{\gamma}$  &   \blank & \zzcell   $\gamma$ &\zzcell  $\Id_{h'}$ \\ \hhline{~------}
                  
                   &    $\Delta_{\wedge}(\beta)$    &\zzcell $\widehat{\Id_h}$ &\zzcell  $\widehat{\beta}$   &  \blank & \zzcell  $\beta$ &\zzcell  $\Id_h$ \\ \hhline{~------}
                   
           $\Lambda_1$        &  $S$    &  \onecell $\theta$    & \onecell $\eta$    &\zzcell  $\gamma$  & \zzcell  $\beta$& \blank \\ \hhline{~------}
                   &    $\theta_{|01|\cdot|}(\theta)$    & \onecell $\theta$   & \onecell $\theta$    &\zzcell  $\Id_{h'}$  & \zzcell  $\Id_h$& \blank \\ \hhline{~------}
                                \end{tabular}\end{center}
                                \end{table}
We can make two $\langle 0,1,0\rangle$-horns in $X$ from Table~\ref{SShattable}. Removing $\Lambda_0$, we get a horn showing that $\hat{S}$ is commutative if $S$ is commutative, and removing $\Lambda_1$ yields a horn which shows that $S$ is commutative if $\hat{S}$ is commutative. 
\end{proof}
\begin{lemma}[Commutativity criterion for $d\Theta_2\langle 1, 0 \rangle$-spheres]\label{10celllemma}
Let $X$ be a $2$-reduced inner-Kan $\Theta_2$-set. Let $$S=[\gamma,\ \beta,\ \theta ,\ \eta]$$ be a $d\Theta_2\langle 1, 0 \rangle$-sphere in $X$, with the following Glenn table:
            \begin{table}[H] \caption{The $d\Theta_2\langle 1, 0\rangle$-sphere $S$}\begin{center}
                \begin{tabular}{ r | l || l | l | l|}
                 \hhline{~----}
              &  \zzcell  $\gamma$    & $h$   &  $g'$   &   $f'$    \\ \hhline{~----}
                       
                 &  \zzcell  $\beta$    &  $h$ & $g$ &     $f$           \\  \hhline{~----}
              & $\theta$\onecell         & $g'$     &$g$     &\blank    \\   \hhline{~----}                      
                 
                    &$\eta$  \onecell     &      $f'$ &$f$   &  \blank         \\  \hhline{~----}
                \end{tabular}\end{center}
                \end{table}
 Then $S$ is commutative if and only if $\theta=t(\underline{\gamma})^{-1}\bullet (h \rhd \eta) \bullet t(\underline{\beta}).$ \end{lemma}
 \begin{proof} By the uniqueness of fillers for inner $\langle 1, 0 \rangle$-horns in $X$, it is enough to show  $$S=[\gamma ,\ \beta,\  t(\underline{\gamma})^{-1}\bullet (h\rhd \eta) \bullet t(\underline{\beta}),\ \eta]$$ is commutative in $X$. The following Glenn table proof demonstrates this commutativity.
     \begin{table}[H] \caption{\label{10lemmaproof1} $\langle 2,0 \rangle$-horn}\begin{center}
         \scalebox{.94}{  \begin{tabular}{ r | l || l | l | l|l|}
                                    \hhline{~-----}
                  $\odot_1$    &  (Table~\ref{10lemmaproof2})     & \zzcell    $\gamma$  &  \zzcell    $\chi(h,f')$ & \onecell $t(\underline{\gamma})^{-1}$  &  \onecell  $\ID_{f'}$     \\  \hhline{~-----}
             
              $\Lambda$         &       &\zzcell   $\gamma$    &\zzcell  $\beta$       &\onecell $ t(\underline{\gamma})^{-1}\bullet (h\rhd \eta) \bullet t(\underline{\beta})$   & \onecell  $\eta$           \\  \hhline{~-----}
              
             $\odot_2$          & (Table~\ref{10lemmaproof3})      &\zzcell  $\chi(h,f')$      & \zzcell   $\beta$            &   \onecell   $ (h\rhd \eta) \bullet t(\underline{\beta})$& \onecell  $\eta$   \\ \hhline{~-----}
                       
  &    $\Theta_{\bullet}(t(\underline{\gamma})^{-1}, (h\rhd \eta) \bullet t(\underline{\beta}) )$    & \onecell $t(\underline{\gamma})^{-1}$      & \onecell  $ t(\underline{\gamma})^{-1}\bullet (h\rhd \eta) \bullet t(\underline{\beta})$  & \onecell $ (h\rhd \eta) \bullet t(\underline{\beta})$ &  \blank  \\ \hhline{~-----}
                     
                     &    $\theta_{|011|}\eta$    &\onecell  $\ID_{f'}$      & \onecell   $\eta$   &  \onecell  $\eta$    &  \blank  \\ \hhline{~-----}
                                  \end{tabular}}\end{center}
                                  \end{table}
  \begin{table}[H] \caption{ $\langle 1,0,0 \rangle$-horn showing $\odot_1$ from Table~\ref{10lemmaproof1} \label{10lemmaproof2} }\begin{center}
           \begin{tabular}{ r | l || l | l | l|l|}
                                 \hhline{~-----}                            
                    &$\Delta_{\stackrel{\wedge}{-}}(\gamma)$       &  \zzcell   $\widehat{\Id_h}$  & \zzcell  $\overline{\gamma}$&\zzcell  $\gamma$ &\zzcell $\chi(h,f')$      \\   \hhline{~-----}
           
                     &$s_2\chi(h,f')$       &\zzcell  $\widehat{\Id_h}$    &\zzcell $\widehat{\Id_{h\circ f'}}$ &\zzcell  $\chi(h,f')$    &\zzcell  $\chi(h,f')$          \\    \hhline{~-----}
            
           &    Lemma~\ref{10cellcor}   & \zzcell    $\overline{\gamma}$ &\zzcell  $\widehat{\Id_{h\circ f'}}$   &   \onecell $t(\underline{\gamma})^{-1}$ &  \onecell $\ID_{h\circ f'}$  \\ \hhline{~-----}
                  
          $\Lambda$         &     & \zzcell    $\gamma$  &  \zzcell    $\chi(h,f')$ & \onecell $t(\underline{\gamma})^{-1}$  &  \onecell  $\ID_{f'}$     \\  \hhline{~-----}
                   
                   &    $\theta_{|00|0|}\chi(h,f')$    &\zzcell $\chi(h,f')$    &\zzcell $\chi(h,f')$    & \onecell $\ID_{h\circ f'}$ &  \onecell $\ID_{f'}$\\  \hhline{~-----}

                                \end{tabular}\end{center}
                                \end{table}
 \begin{table}[H] \caption{$\langle 2,0 \rangle$-horn showing $\odot_2$ from Table~\ref{10lemmaproof1} \label{10lemmaproof3}}\begin{center}
           \begin{tabular}{ r | l || l | l | l|l|}
                                  \hhline{~-----}
                    &$\Theta_{\rhd}(h,\eta)$       & \zzcell  $\chi(h,f')$  &  \zzcell    $\chi(h,f)$ & \onecell $h\rhd \eta$    &  \onecell  $\eta$     \\  \hhline{~-----}
           
            $\Lambda$         &      &\zzcell  $\chi(h,f')$      & \zzcell   $\beta$            &   \onecell   $ (h\rhd \eta) \bullet t(\underline{\beta})$& \onecell  $\eta$      \\  \hhline{~-----}
            
              $\odot$       & (Table~\ref{10lemmaproof4})      &\zzcell  $\chi(h,f)$       & \zzcell    $\beta$            &   \onecell   $t(\underline{\beta})$ & \onecell  $\ID_f$   \\ \hhline{~-----}
    
              &    $\Theta_{\bullet}(h\rhd \eta,t(\underline{\beta}) )$    & \onecell $h\rhd \eta$      & \onecell   $ (h\rhd \eta) \bullet t(\underline{\beta})$ & \onecell $t(\underline{\beta})$ &  \blank  \\ \hhline{~-----}
                   
                   &    $\theta_{|001|}\eta$    &\onecell  $\eta$      & \onecell  $\eta$  & $\ID_f$ \onecell     &  \blank  \\ \hhline{~-----}
                                \end{tabular}\end{center}
                                \end{table}

  \begin{table}[H] \caption{ $\langle 1,0,0 \rangle$-horn showing $\odot$ from Table~\ref{10lemmaproof3} \label{10lemmaproof4} }\begin{center}
           \begin{tabular}{ r | l || l | l | l|l|}
                                 \hhline{~-----}                            
                    &$s_2\chi(h,f)$      &  \zzcell   $\widehat{\Id_h}$  & \zzcell  $\widehat{\Id_{h\circ f}}$&\zzcell  $\chi(h,f)$ &\zzcell $\chi(h,f)$    \\   \hhline{~-----}
           
                     &   $\Delta_{\stackrel{\wedge}{-}}(\beta)$     &\zzcell  $\widehat{\Id_h}$    &\zzcell $\overline{\beta}$ &\zzcell   $\beta$    &\zzcell  $\chi(h,f)$          \\    \hhline{~-----}
            
           &    Lemma~\ref{10cellcor}   & \zzcell    $\widehat{\Id_{h\circ f}}$ &\zzcell  $\overline{\beta}$  &   \onecell  $t(\underline{\beta})$  &  \onecell $\ID_{h\circ f'}$  \\ \hhline{~-----}
                  
          $\Lambda$         &     &\zzcell  $\chi(h,f)$       & \zzcell    $\beta$            &   \onecell   $t(\underline{\beta})$ & \onecell  $\ID_f$    \\  \hhline{~-----}
                   
                   &    $\theta_{|00|0|}\chi(h,f)$    &\zzcell $\chi(h,f)$   &\zzcell $\chi(h,f)$    & \onecell $\ID_{h\circ f'}$ &  \onecell $\ID_{f'}$\\  \hhline{~-----}

                                \end{tabular}\end{center}
                                \end{table}
 \end{proof}
 \begin{proposition}[$t$ preserves $\lhd$] Let $\eta:g\Rightarrow h:b\ra c$ and $f:a\ra b$ in $\overline{\TBic(X)}$. Then $t(\eta)\lhd f = t(\eta \lhd f).$ \label{tpreservesrhd}\end{proposition}
 \begin{proof} By Lemma~\ref{02celllemma} it suffices to show the $\langle 0,1\rangle$ sphere $$[t(\eta),\ \ID_{g\circ f},\ \eta \tlhd f,\ \chi_{g,f}]$$ is commutative. The following $\langle 0 , 0 , 1 \rangle$-horn verifies this commutativity.

\begin{table}[H] \caption{ $\langle 0,0,1 \rangle$-horn }\begin{center}
            \begin{tabular}{ r | l || l | l | l|l|}
                               \hhline{~-----}
                     &Corollary~\ref{02cellcor}       &  \onecell   $t(\eta)$  &  \onecell $\ID_{g}$& \zzcell $\eta$  & \zzcell $\Id_g$       \\   \hhline{~-----}
            
                     $\Lambda$ &      &\onecell   $t(\eta)$       &\onecell   $\ID_{g\circ f}$         &   \zzcell   $\eta \tlhd f$& \zzcell  $\chi(g,f)$  \\   \hhline{~-----}
             
                     &    $\theta_{|0|00|}\chi(g,f)$    &  \onecell  $\ID_{g}$ & \onecell  $\ID_{g\circ f}$    &   \zzcell $\chi(g,f)$ & \zzcell 
                      $\chi(g,f)$ \\  \hhline{~-----}
                   
                    &      $\Delta_{\tlhd}(\eta,f)$  &\zzcell $\eta$   &\zzcell $\eta \tlhd f$  &  \zzcell $\chi(g,f)$  & \zzcell  $\widehat{\Id_f}$ \\  \hhline{~-----}
                    
                    &    $s_1 \chi(g,f)$    &\zzcell $\Id_g$   &\zzcell$\chi(g,f)$      &\zzcell  $\chi(g,f)$ & \zzcell  $\widehat{\Id_f}$ \\  \hhline{~-----}
 
                                 \end{tabular}\end{center}
                                 \end{table}\end{proof}
 \begin{proposition}[$t$ preserves $\rhd$] Let $\eta:f\Rightarrow g:a\ra b$ and $h:b\ra c$ in $\overline{\TBic(X)}$. Then $h\rhd t(\eta) = t(h\rhd\eta).$\end{proposition} \begin{proof} This follows by an argument symmetrical to the proof of Proposition~\ref{tpreservesrhd}.
  \end{proof}
  
  We have now shown that $t$ is a strict functor, since we have shown that it preserves identity for $1$-morphisms, $\bullet$, $\rhd$, and $\lhd$, and it preserves the unitors and the associator by definition. We are left to verify the rest of the bicategory axioms for $\widetilde{\TBic(X)}.$ Of these, \textbf{B4},\textbf{B12}-\textbf{B15}, and \textbf{B17} all follow immediately from the definition of the unitors and associators for $\widetilde{\TBic(X)}$ together with the fact that $t$ is a strict functor and the corresponding axioms hold in $\overline{\TBic(X)}.$ The naturality of $\lambda$, \textbf{B8} follows by using Lemma~\ref{10celllemma}  in a similar argument to proof of the naturality of $\rho$ (\textbf{B7}) in Proposition~\ref{natofrho}. Also recall that \textbf{B3} is the axiom guarantees the existence of inverses, which does not apply in the bicategory case. This leaves the interchange axioms \textbf{B5},\textbf{B6}, and \textbf{B16} and the naturality of the associator axioms \textbf{B9}--\textbf{B11}.
\subsection{Interchange axioms for $\widetilde{\TBic(X)}$}
\begin{proposition}[\textbf{B5} Interchange of $\rhd$ and $\bullet$] \label{rhdinterchange} Let $f \stackrel{\eta}\Rightarrow g \stackrel{\theta}\Rightarrow h: a \to b$ and $i:b \to c$ in $\widetilde{\TBic(X)}$,  then $(i \rhd \theta) \bullet (i \rhd \eta)=i \rhd (\theta \bullet \eta).$
\end{proposition}
\begin{proof}
Applying the Matching Lemma to $\Theta_{\bullet}(i\rhd \theta, i\rhd \eta)$, we see that it suffices to show that the following $d\Theta_2\langle 2\rangle$-sphere is commutative: $$[  i \rhd \theta  ,\ i \rhd (\theta \bullet \eta),\ i \rhd \eta ].$$ The following $\langle 2,0\rangle$-horn shows this commutativity:

\begin{table}[H] \caption{ $\langle 2,0 \rangle$-horn }\begin{center}
\begin{tabular}{ r | l || l | l | l|l|}
                                  \hhline{~-----}
&$\Theta_{\rhd}(i,\theta)$     & \zzcell   $\chi(i, h)$  &  \zzcell  $\chi(i, g)$ & \onecell  $i \rhd \theta$     &  \onecell  $\theta$     \\  \hhline{~-----}
           
 & $\Theta_{\rhd}(i,\theta\bullet\eta)$        &\zzcell  $\chi(i, h)$    &\zzcell   $\chi(i, f)$     &\onecell $i \rhd (\theta \bullet \eta)$  & \onecell  $\theta \bullet \eta$           \\  \hhline{~-----}
            
 &  $\Theta_{\rhd}(i,\eta)$      &\zzcell   $\chi(i, g)$     & \zzcell  $\chi(i, f)$          &   \onecell   $i \rhd \eta$ & \onecell  $\eta$   \\ \hhline{~-----}
                  
$\Lambda$           &       & \onecell $i \rhd \theta$      & \onecell  $i \rhd (\theta \bullet \eta)$  & \onecell $i \rhd \eta$ &  \blank  \\ \hhline{~-----}
                   
&    $\Theta_{\bullet}(\theta, \eta)$    &\onecell  $\theta$      & \onecell  $\theta \bullet \eta$  &  \onecell  $\eta$    &  \blank  \\ \hhline{~-----}
                                \end{tabular}\end{center}
                                \end{table}
\end{proof}
\begin{proposition}[\textbf{B6} Interchange of $\lhd$ and $\bullet$] Let $g \stackrel{\eta}\Rightarrow h \stackrel{\theta}\Rightarrow i: b \to c$ and $f:a \to b$ in $\widetilde{\TBic(X)}$,  then $(\theta\lhd f) \bullet (\eta\lhd f)=(\theta \bullet \eta)\lhd f.$
\end{proposition}
\begin{proof}Follows by a similar argument to the proof of \textbf{B5} in Proposition~\ref{rhdinterchange} above, using a $\langle 0,2\rangle$-horn.
\end{proof}
\begin{proposition}
[\textbf{B16} Full Interchange]
Let $f\stackrel{\eta}{\Rightarrow} g:a\ra b$ and $h\stackrel{\theta}{\Rightarrow} i:b\ra c$ in $\widetilde{\TBic(X)}.$ Then $(i\rhd \eta) \bullet (\theta \lhd f)=(\theta \lhd g)\bullet (h \rhd \eta).$
\end{proposition}
\begin{proof}
Applying the Matching Lemma to $\Theta_{i\rhd \eta, \theta \lhd f}$, we see that it suffices to show that the following $d\Theta_2\langle 2\rangle$-sphere is commutative: $$[  i \rhd \eta  ,\ (\theta \lhd g)\bullet (h \rhd \eta),\ \theta \lhd f ].$$ The following $\langle 1,1 \rangle$-horn shows this commutativity:
   \begin{table}[H] \caption{Universal Glenn table for $\langle 1,1 \rangle$ }\begin{center}
           \begin{tabular}{ r | l || l | l | l|l|l|}
                                  \hhline{~------}
                    &$\Theta_{\lhd}(\theta, g)$       &  \onecell  $ \theta$  &\onecell  $\theta \lhd g$ &\blank & \zzcell $\chi(i,g)$  &   \zzcell  $\chi(h,g)$     \\   \hhline{~------}
           
                     & $\Theta_{\lhd}(\theta, f)$        &\onecell  $\theta$    &\blank &\onecell   $\theta \lhd f$      & \zzcell $\chi(i,f)$  &  \zzcell $\chi(i,f)$        \\    \hhline{~------}
            
            &  $\Theta_{\bullet}(\theta\lhd g,\ h \rhd \eta)$    & \onecell  $\theta \lhd g$ &  \blank   & \onecell   $(\theta \lhd g)\bullet (h \rhd \eta)$& \blank &\onecell  $h \rhd \eta$  \\  \hhline{~------}
                  
       $\Lambda$     &      &\blank    & \onecell \ \circletext{$i \rhd \eta$} \ &\onecell  $(\theta \lhd g)\bullet (h \rhd \eta)$ & \onecell \ \circletext{$\theta \lhd f $} \  &\blank\\  \hhline{~------}
                   
                   &    $\Theta_{\rhd}(i,\eta)$    & \zzcell $\chi(i,g)$      &  \zzcell  $\chi(i,f)$    &  \blank   & \onecell  $i \rhd \eta$ &\onecell  $\eta$\\ \hhline{~------}
                  
                   &    $\Theta_{\rhd}(h,\eta)$    &\zzcell  $\chi(h,g)$      & \zzcell   $\chi(h,f)$        &\onecell   $h \rhd \eta$&   \blank & \onecell $\eta$ \\  \hhline{~------}
                                \end{tabular}\end{center}
                                \end{table}
\end{proof}
\subsection{ Naturality of the associator in $\TBic$}
\begin{proposition}[\textbf{B9} Naturality of the associator in the first argument] \label{natassfirst}
Let $ f\stackrel{\eta}{\Rightarrow} g:a \to b$ and $b \stackrel{h}\to c \stackrel{i}\to d$ in $\widetilde{\TBic(X)}$, then  $$ i \rhd (h \rhd \eta)\bullet=  t(\alpha_{i,h,g})^{-1} \bullet            ((i \circ h) \rhd \eta) \bullet t( \alpha_{i,h,f}).$$ \end{proposition}
\begin{proof}By Lemma~\ref{10celllemma} this is equivalent to the commutativity of the $\langle 1,0\rangle$ sphere $$S=[\widetilde{\alpha}_{i,h,g} ,\  \widetilde{\alpha}_{i,h,f},\  i \rhd (h \rhd \eta),\ \eta].$$                        
          \begin{table}[H] \caption{$\langle 1,0,0 \rangle$-horn }\begin{center}
           \begin{tabular}{ r | l || l | l | l|l|}
                                 \hhline{~-----}                            
                    &$\Delta_{\widetilde{\alpha}}(i,h,g)$  &  \zzcell   $\chi(i,h)$  & \zzcell  $\chi(i, h\circ g)$&\zzcell  $\widetilde{\alpha}_{i,h,g}$ &\zzcell $\chi(h,g)$      \\   \hhline{~-----}
           
                     &$\Delta_{\widetilde{\alpha}}(i,h,f)$       &\zzcell $\chi(i,h)$    &\zzcell $\chi(i, h\circ f)$ &\zzcell  $\widetilde{\alpha}_{i,h,f}$      &\zzcell $\chi(h,f)$           \\    \hhline{~-----}
            
               &   $\Theta_{\rhd}(i, h\rhd \eta)$ & \zzcell    $\chi(i, h\circ g)$&\zzcell  $\chi(i, h\circ f)$  &   \onecell  $ i \rhd (h \rhd \eta)$&  \onecell $h \rhd \eta$  \\ \hhline{~-----}
                  
   $\Lambda$   &       &\zzcell $\widetilde{\alpha}_{i,h,g}$  &\zzcell  $\widetilde{\alpha}_{i,h,f}$   &\onecell  $ i \rhd (h \rhd \eta)$&  \onecell $\eta$  \\  \hhline{~-----}
                   
                   &     $\Theta_{\rhd}(h, \eta)$    &\zzcell $\chi(h,g)$    &\zzcell $\chi(h,f)$       & \onecell$h \rhd \eta$  &  \onecell  $\eta$ \\  \hhline{~-----}

                                \end{tabular}\end{center}
                                \end{table}

\end{proof}
\begin{proposition}[\textbf{B10} Naturality of the associator in the second argument] Let  $a\stackrel{f}{\ra}b$ and $g\stackrel{\eta}{\Rightarrow}h:b\ra c$ and $c\stackrel{i}{\ra}d$ in $\widetilde{\TBic(X)}$, then $$i \rhd (\eta \lhd f)=t(\alpha_{i,h,f})^{-1}\bullet ((i \rhd \eta) \lhd f) \bullet t(\alpha_{i,g,f}).$$ \end{proposition}
\begin{proof}By Lemma~\ref{02celllemma} this is equivalent to the commutativity of the $\langle 0,1\rangle$ sphere $$S=[i \rhd \eta ,\  i \rhd (\eta \lhd f),\ \widetilde{\alpha}_{i,h,f},\ \widetilde{\alpha}_{i,g,f}].$$    

   \begin{table}[H] \caption{ $\langle 0,1,0 \rangle$-horn}\begin{center}
           \begin{tabular}{ r | l || l | l | l|l|l|}
                                \hhline{~------}
                    & $\Theta_{\rhd}(i, \eta)$   &    \blank  &\zzcell    $\chi(i,h)$ &\zzcell $\chi(i,g)$& \onecell  $i \rhd \eta$   &  \onecell$\eta$   \\  \hhline{~------}
           
                    &  $\Theta_{\rhd}(i, \eta\lhd f)$  & \blank  &\zzcell  $\chi(i,h\circ f)$ & \zzcell $\chi(i,g\circ f)$   & \onecell $i \rhd (\eta \lhd f)$  &  \onecell $\eta \lhd f$    \\  \hhline{~------}
            
                     &    $\Delta_{\widetilde{\alpha}}(i,h,f)$    &  \zzcell   $\chi(i,h)$ &\zzcell  $\chi(i,h\circ f)$  &   \blank & \zzcell  $\widetilde{\alpha}_{i,h,f}$  &\zzcell  $\chi(h,f)$   \\ \hhline{~------}
                  
                   &    $\Delta_{\widetilde{\alpha}}(i,h,g)$   &\zzcell  $\chi(i,g)$ &\zzcell  $\chi(i,g\circ f)$    &  \blank & \zzcell $\widetilde{\alpha}_{i,g,f}$ &\zzcell $\chi(g,f)$  \\ \hhline{~------}
                   
           $\Lambda$        &    &  \onecell $i \rhd \eta$    & \onecell $i \rhd (\eta \lhd f)$    &\zzcell  $\widetilde{\alpha}_{i,h,f}$  & \zzcell  $\widetilde{\alpha}_{i,g,f}$& \blank \\ \hhline{~------}
                   &  $\Theta_{\lhd}(i,\eta )$   & \onecell $\eta$   & \onecell $\eta \lhd f$    &\zzcell $\chi(h,f)$  & \zzcell  $\chi(g,f)$ & \blank \\ \hhline{~------}
                                \end{tabular}\end{center}
                                \end{table}\end{proof}
\begin{proposition}[\textbf{B11} Naturality of the associator in the last argument] Let $a\stackrel{f}{\ra}b\stackrel{g}{\ra}c$ and $h\stackrel{\eta}{\Rightarrow}i:c\ra d$ in $\widetilde{\TBic(X)}$, then $$ \eta \lhd (g \circ f)=t(\alpha_{i,g,f})^{-1} \bullet((\eta \lhd g) \lhd f) \bullet t(\alpha_{h,g,f}).$$ \end{proposition} \begin{proof} This follows by an argument symmetrical to that used for the naturality of the associator in the first argument in the proof of Proposition~\ref{natassfirst}. This argument uses a $\langle 0,0,1\rangle$-horn.
\end{proof}
\subsection{The construction $\TBic$ is functorial }
Let $F:X\ra Y$ be a morphism of $2$-reduced inner-Kan $\Theta_2$-sets. To show the construction $\TBic$ is functorial, we must define a functor of fancy bicategories $\TBic(X)\ra \TBic(Y).$ We have functor $$\overline{\TBic(F)}:=\Bic(\psi^{\star}F):\Bic(\psi^{\star}X)\ra\Bic(\psi^{\star}Y),$$ also defining a map $\widetilde{\TBic(X)}\ra \widetilde{\TBic(Y)}$ for objects and $1$-morphisms. For $2$-morphisms $\widetilde{\TBic(F)}:\widetilde{\TBic(X)}\ra \widetilde{\TBic(Y)}$ is defined by $F$ restricted to $\langle 1 \rangle$ cells. The distributor is $t(\phi_{g,f})$ and the unitor is $t(\upsilon_a)$ where $\phi_{g,f}$ and $\upsilon_a$ are the distributor and unitor of $\overline{\TBic(F)}.$  By this definition $\TBic(F)$ clearly respects $t$ and makes a functor of fancy bicategories if we can show $\widetilde{\TBic(F)}$ is a functor. 

The compatibility of the distributor of $\TBic(F)$ and the associator and the compatibility of the unitor of $\TBic(F)$ and the unitors of $\widetilde{\TBic(X)}$ and $\widetilde{\TBic(Y)}$  (\textbf{TFun5},\textbf{TFun6}, and \textbf{TFun7}) follow from the same axioms for $\Bic(\psi^{\star}F).$ The fact that $\TBic(F)$ preserves identity (\textbf{TFun1}) follows from the fact that $F$ preserves $\theta_{|00|}$.

For \textbf{TFun2} we have
\begin{align*}d(F(\Theta_{\bullet}(\theta,\eta)))&=[F(\theta),\ F(\theta\bullet \eta),\ F(\eta)] \\ d(\Theta_{\bullet}(F(\theta),F(\eta)))&= [F(\theta),\ F(\theta)\bullet F(\eta),\ F(\eta)] \end{align*}
this shows $F(\theta\bullet \eta )=F(\theta)\bullet F(\eta)$ for $2$-morphisms of $\widetilde{\TBic(X)}$, which proves this axiom.

For \textbf{TFun3} we have $$d(F(\Theta_{\rhd}(h,\eta)))=[F(\chi(h,g)),\ F(\chi(h,f)),\ F(h\rhd \eta),\ F(\eta)].$$ From Lemma~\ref{10celllemma} we conclude 
\begin{align*}
F(h\rhd \eta)&=t(\underline{F(\chi(h,g))})^{-1}\bullet (F(h)\rhd F(\eta))\bullet t(\underline{F(\chi(h,f))})\\
&=t(\phi_{h,g})^{-1}\bullet (F(h)\rhd F(\eta))\bullet t(\phi_{h,f})
\end{align*}
This proves \textbf{TFun3}. The $\lhd$ case \textbf{TFun4} is similar. The fact the construction given above preserves composition (i.e. a functor) follows trivially from the fact that $\Bic$ is a functor.

\section{The $\Theta_2$ nerve of a fancy bicategory}
The idea of defining a $\Theta_n$-set nerve of a strict $n$-category is due to Berger \cite{Ber02}. We will define a generalization of this concept in the case $n=2$, defining a functor $N_{\theta}$ taking a small fancy bicategory $\B$ to a $\Theta_2$-set.
 
 We first define the $2$-truncation $N_\theta(\B)|^2_0$ as a $\Theta_2|^2_0$-set. We let the $\langle \rangle$, $\langle 0\rangle$ and $\langle 0,0\rangle$ cells be the $0$, $1$ and $2$ cells of the Duskin nerve $N(\olb)$, respectively, with appropriate face maps induced by the full and faithful embedding $\psi:\Delta \ra \Theta_2$. The $\langle 1 \rangle$ cells of $N_\theta(\B)|^2_0$ be the $2$-morphisms of $\wtb$ with $\theta_{|0|} \eta= \mbox{source}(\eta)$ and $\theta_{|1|}\eta=\mbox{target}(\eta)$. Recalling that we identify the $1$-morphisms of $\olb$ and $\wtb$, we define $\theta_{|00|}f=\Id_f.$ It is straightforward to check that these maps generate $\Theta_2|^2_0$ and that the above definition defines a $\Theta_2|^2_0$-set. We define $\Theta_2|^3_0$ as a subsheaf of $\tr^3\cosk^2(N_\theta(\B)|^2_0)$. To do this, we must define which $3$-spheres will be commutative in $N_\theta(\B)|^3_0.$
 
 A $d\Theta_2\langle 2 \rangle$-sphere $[\pi,\theta,\eta]$ will be commutative if $\theta=\pi \bullet \eta$. A $d\Theta_2\langle 1,0\rangle$-sphere with Glenn table given below in Table~\ref{10sphereS} is commutative if $\theta=t_\B(\gamma)^{-1}\bullet (h \rhd \eta)\bullet t_\B(\beta).$
             \begin{table}[H] \caption{The $d\Theta_2\langle 1, 0\rangle$-sphere $S$ \label{10sphereS}}\begin{center}
                 \begin{tabular}{ r | l || l | l | l|}
                  \hhline{~----}
               &  \zzcell  $(h,g',f' \vbar \gamma)$    & $h$   &  $g'$   &   $f'$    \\ \hhline{~----}
                        
                  &  \zzcell  $(h,g,f \vbar \beta)$    &  $h$ & $g$ &     $f$           \\  \hhline{~----}
               & $\theta$\onecell         & $g'$     &$g$     &\blank    \\   \hhline{~----}                      
                  
                     &$\eta$  \onecell     &      $f'$ &$f$   &  \blank         \\  \hhline{~----}
                 \end{tabular}\end{center}
                 \end{table}
A $d\Theta_2\langle 0,1\rangle$-sphere with Glenn table given below in Table~\ref{01sphereS} is commutative is commutative if $\eta=t_\B(\gamma)^{-1}\bullet (\theta\lhd f) \bullet t_\B(\beta).$ 
            \begin{table}[H] \caption{The $d\Theta_2\langle 0, 1\rangle$-sphere $S$ \label{01sphereS}}\begin{center}
                \begin{tabular}{ r | l || l | l | l|}
                 \hhline{~----}

                                &$\theta$  \onecell     &   \blank  &     $h'$ &$h$          \\  \hhline{~----}

                       & $\eta$\onecell        &\blank  & $g'$     &$g$    \\   \hhline{~----}
                    
                                    &  \zzcell  $(h',g',f \vbar \gamma)$    &  $h'$ & $g'$ &     $f$     \\ \hhline{~----}
                                          
                                    &  \zzcell  $(h,g,f \vbar \beta)$    & $h$   &  $g$   &   $f$       \\  \hhline{~----}
                       
                \end{tabular}\end{center}
                \end{table} 
 We must check that each sphere that can be obtained from a $2$ dimensional or smaller cell by map in $\Theta_2$ is commutative.  It is not hard to check that any such operator factors through $\theta_{|001|}$, $\theta_{|011|}$, $\theta_{|00|0|}$, $\theta_{|0|00|}$, $\theta_{|\cdot|01|}$, $\theta_{|01|\cdot|}$ or maps in the image of a simplicial degeneracy map under $\psi$. The fact that the last of these produce commutative spheres follows from the fact that the  simplicial degeneracies give commutative spheres in the Duskin nerve which was shown in Chapter~\ref{bicchapter}. 
 
 For $\eta:f\Rightarrow g$ the sphere given by $\theta_{|001|}\eta$ is $[\eta,\eta,\Id_f]$ which can immediately be seen to satisfy the commutativity condition $\eta=\eta\bullet \Id_f$ for $\langle 2\rangle$ cells. The other cases are also easily verified, with the $\theta_{|\cdot|01|}$ and $\theta_{|01|\cdot|}$ case using the naturality of the unitors of $\wtb$. 
 
 Finally we define $N_{\theta}(\B)=\Cosk^3(N_{\theta}(\B)|^3_0).$
 
 \begin{proposition} $N_{\theta}(\B)$ is $2$-reduced inner-Kan.
  \end{proposition}
  \begin{proof}
  Each inner horn of dimension $3$ is immediately seen to have a unique filler from our definition. The only $2$-dimensional inner horn is for $\langle 0,0\rangle$, which is part of image of $\psi$. These horns have preferred fillers of the form $(g, g\circ f, f \vbar \Id_{g\circ f})$ as in the Duskin nerve. Lemma~\ref{coskeletaltofiller} ensures that inner horns of dimension $5$ and higher have unique fillers.
  
The fact that an inner $\langle 0,0,0,0\rangle$-horn in $N_{\theta}(\B)$ has a unique filler follows from the fact that the Duskin nerve  $N(\overline{\B})$ is $2$-reduced inner-Kan. 

Next consider $\langle 1,0,0\rangle$-horns. The following equations give the commutativity conditions for the five faces of a $\langle 1, 0 , 0 \rangle$ cell $x_{|01|0|0|}$. Note that for the calculations in this proof, we conflate a $\langle 0, 0\rangle$ cell with its interior $2$-morphism in $\overline{\B}$.
\begin{align}
\alpha_{ x_{\star\star|0|}, x_{\star|0|\star},x_{|1|\star \star}} \bullet (x_{\star\star|0|}\rhd x_{|1|0|\star})\bullet x_{|1,0|0|} &=(x_{\star|0|0|}\lhd x_{|1|\star\star})\bullet x_{|1|0,0|}  \label{100nerve1} \\
\alpha_{ x_{\star\star|0|}, x_{\star|0|\star},x_{|0|\star \star}} \bullet (x_{\star\star|0|}\rhd x_{|0|0|\star})\bullet x_{|0,0|0|} &=(x_{\star|0|0|}\lhd x_{|0|\star\star})\bullet x_{|0|0,0|} \label{100nerve2} \\
(\Lambda) \quad \quad x_{|01,00,00|}&= t_\B(x_{|1,0|0|})^{-1}\bullet (x_{\star\star|0|}\rhd x_{|01,00|\star}   )\bullet t_\B(x_{|0,0|0|}) \label{100nerve3} \\
(\Lambda) \quad \quad x_{|01,00,00|}&= t_\B(x_{|1|0,0|})^{-1}\bullet (x_{\star|0,0|}\rhd x_{|01|\star\star}   )\bullet t_\B(x_{|0|0,0|})  \label{100nerve4} \\
x_{|01,00|\star}&= t_\B(x_{|1|0|\star})^{-1}\bullet (x_{\star|0|\star}\rhd x_{|01|\star\star}   )\bullet t_\B(x_{|0|0|\star})  \label{100nerve5}
\end{align}
We must check that the third or fourth equation holds given the other four equations. Equivalently, we must show Equation~\ref{100nerve3} and Equation~\ref{100nerve4} are equivalent, given the other three equations. We show this by showing the right hand side of the two equations are equal, given the other three equations. We start with the right had side of \ref{100nerve3}:

\begin{align*}
 &t_\B(x_{|1,0|0|})^{-1}\bullet (x_{\star\star|0|}\rhd x_{|01,00|\star}   )\bullet t_\B(x_{|0,0|0|})  \end{align*}
 We substitute for $x_{|01,00| \star}$ using Equation~\ref{100nerve5} then use interchange for $\rhd$:
 \begin{align*}
=& t_\B(x_{|1,0|0|})^{-1}\bullet (x_{\star\star|0|}\rhd ( t_\B(x_{|1|0|\star})^{-1}\bullet (x_{\star|0|\star}\rhd x_{|01|\star\star}   )\bullet t_\B(x_{|0|0|\star}))  )\bullet t_\B(x_{|0,0|0|})  \\
=&t_\B(x_{|1,0|0|})^{-1}\bullet t_\B(x_{\star\star|0|}\rhd  x_{|1|0|\star})^{-1} \bullet (x_{\star\star|0|}\rhd (x_{\star|0|\star}\rhd x_{|01|\star\star}))   \bullet t_\B(x_{\star\star|0|}\rhd    x_{|0|0|\star}) \bullet t_\B(x_{|0,0|0|})  \\
\end{align*}
If we apply $t_\B$ to both sides of Equations~\ref{100nerve1} and \ref{100nerve2} we can make a substitution for the first two and last two terms, then apply the naturality of the associator to make a simplification:
\begin{align*}
=&t_\B(x_{|1|0,0|})^{-1}\bullet t_\B(x_{\star|0|0|}\lhd x_{|1|\star\star})^{-1} \bullet t_\B(\alpha_{ x_{\star\star|0|}  , x_{\star|0|\star}, x_{|1|\star \star}}) \bullet (x_{\star\star|0|}\rhd (x_{\star|0|\star}\rhd x_{|01|\star\star}))   \\ &\quad \bullet  t_\B(\alpha_{  x_{\star\star|0|}, x_{\star|0|\star},  x_{|0|\star \star} })^{-1} \bullet t_\B(x_{\star|0|0|}\lhd x_{|0|\star\star})\bullet t_\B(x_{|0|0,0|})\\
=&t_\B(x_{|1|0,0|})^{-1}\bullet t_\B(x_{\star|0|0|}\lhd x_{|1|\star\star})^{-1} \bullet ((x_{\star\star|0|}\circ  x_{\star|0|\star})\rhd x_{|01|\star\star})  \bullet t_\B(x_{\star|0|0|}\lhd x_{|0|\star\star})\bullet t_\B(x_{|0|0,0|})\\
\end{align*}
Finally we use the full interchange law and make a cancellation:
\begin{align*}
=&t_\B(x_{|1|0,0|})^{-1}\bullet t_\B(x_{\star|0|0|}\lhd x_{|1|\star\star})^{-1} \bullet t_\B(x_{\star|0|0|}\lhd x_{|1|\star\star})  \bullet (x_{\star|0,0|}\rhd x_{|01|\star\star}   )    \bullet t_\B(x_{|0|0,0|})\\
=&t_\B(x_{|1|0,0|})^{-1}\bullet (x_{\star|0,0|}\rhd x_{|01|\star\star}   )    \bullet t_\B(x_{|0|0,0|}).
\end{align*}

For $\langle 0,1,0\rangle$ we have the following commutativity conditions for the faces of a $\langle 0,1,0\rangle$ cell $x_{|0|01|0|}.$

\begin{align}
x_{\star|01,00|}&= t_\B(x_{\star|1|0|})^{-1}\bullet ( x_{\star\star|0|}\rhd x_{\star|01|\star})\bullet t_\B(x_{\star|0|0|})     \label{010nerve1} \\
 (\Lambda) \quad \quad  x_{|00,01,00|}&= t_\B(x_{|0,1|0|})^{-1}\bullet ( x_{\star\star|0|}\rhd x_{|00,01|\star})\bullet t_\B(x_{|0,0|0|}) \label{010nerve2}  \\
\alpha_{  x_{\star\star|0|}  , x_{\star|1|\star}, x_{|0|\star \star}} \bullet (x_{\star\star|0|}\rhd x_{|0|1|\star})\bullet x_{|0,1|0|} &=(x_{\star|1|0|}\lhd x_{|0|\star\star})\bullet x_{|0|1,0|}  \label{010nerve3} \\
\alpha_{      x_{\star\star|0|}, x_{\star|0|\star},x_{|0|\star \star}} \bullet (x_{\star\star|0|}\rhd x_{|0|0|\star})\bullet x_{|0,0|0|} &=(x_{\star|0|0|}\lhd x_{|0|\star\star})\bullet x_{|0|0,0|} \label{010nerve4}   \\
(\Lambda) \quad \quad  x_{|00,01,00|}&= t_\B(x_{|0|1,0|})^{-1}\bullet ( x_{\star|01,00|}\lhd x_{|0|\star\star} )\bullet t_\B(x_{|0|0,0|})    \label{010nerve5}     \\
x_{|00,01|\star}&= t_\B(x_{|0|1|\star})^{-1}\bullet (  x_{\star|01|\star}\lhd x_{|0|\star\star} )\bullet t_\B(x_{|0|0|\star}) \label{010nerve6}
\end{align}
We must show Equation~\ref{010nerve2} and Equation~\ref{010nerve5} are equivalent, given the other four equations. We show this by showing the right hand side of the two equations are equal, given the other four equations. We start with the right had side of \ref{010nerve2} then substitute for $ x_{|00,01|\star}$ using Equation~\ref{010nerve6}, then use interchange for $\rhd$:
\begin{align*} 
&t_\B(x_{|0,1|0|})^{-1}\bullet ( x_{\star\star|0|}\rhd x_{|00,01|\star})\bullet t_\B(x_{|0,0|0|})\\
=&t_\B(x_{|0,1|0|})^{-1}\bullet ( x_{\star\star|0|}\rhd  ( t_\B(x_{|0|1|\star})^{-1}\bullet (  x_{\star|01|\star}\lhd x_{|0|\star\star} )\bullet t_\B(x_{|0|0|\star})  )  )\bullet t_\B(x_{|0,0|0|})\\
=&t_\B(x_{|0,1|0|})^{-1}\bullet  t_\B(  x_{\star\star|0|}\rhd    x_{|0|1|\star})^{-1}\bullet ( x_{\star\star|0|}\rhd ( x_{\star|01|\star}\lhd x_{|0|\star\star} ) )\bullet t_\B(x_{\star\star|0|}\rhd x_{|0|0|\star})   \bullet t_\B(x_{|0,0|0|})\\
 \end{align*}
If we apply $t_\B$ to both sides of Equations~\ref{010nerve3} and \ref{010nerve4} we can make a substitution for the first two and last two terms, then apply the naturality of the associator and then interchange for $\lhd$ simplifications:
\begin{align*}
=&t_\B(x_{|0|1,0|})^{-1} \bullet t_\B(x_{\star|1|0|}\lhd x_{|0|\star\star})^{-1} \bullet t_\B(\alpha_{     x_{\star\star|0|}     , x_{\star|1|\star},x_{|0|\star \star} }) \bullet ( x_{\star\star|0|}\rhd ( x_{\star|01|\star}\lhd x_{|0|\star\star} ) ) \\ & \quad \quad \bullet t_\B(\alpha_{  x_{\star\star|0|}  , x_{\star|0|\star}, x_{|0|\star \star}})^{-1} \bullet t_\B(x_{\star|0|0|}\lhd x_{|0|\star\star})  \bullet t_\B(x_{|0|0,0|}) \\
=&t_\B(x_{|0|1,0|})^{-1} \bullet t_\B(x_{\star|1|0|}\lhd x_{|0|\star\star})^{-1} \bullet  ( ( x_{\star\star|0|}\rhd  x_{\star|01|\star} )\lhd x_{|0|\star\star} )  \bullet t_\B(x_{\star|0|0|}\lhd x_{|0|\star\star})  \bullet t_\B(x_{|0|0,0|}) \\
=&t_\B(x_{|0|1,0|})^{-1} \bullet   (( t_\B(x_{\star|1|0|})^{-1} \bullet ( x_{\star\star|0|}\rhd  x_{\star|01|\star})  \bullet t_\B(x_{\star|0|0|}))\lhd x_{|0|\star\star} )  \bullet t_\B(x_{|0|0,0|}) 
 \end{align*}
Finally we make a substitution using Equation~\ref{010nerve1}, yielding the right hand side of \ref{010nerve2} as desired:
$$=t_\B(x_{|0|1,0|})^{-1} \bullet   (x_{\star|01,00|}\lhd x_{|0|\star\star} )  \bullet t_\B(x_{|0|0,0|}).$$

The $\langle 0,0,1\rangle$ case is similar to the $\langle 1,0,0\rangle$ case above. 

For the $\langle 2,0\rangle$ case we have the following commutativity conditions for the faces of a $\langle 2,0 \rangle$ cell $x_{|012|0|}.$
\begin{align}
                   x_{|12,00|}&= t_\B(x_{|2|0|})^{-1}\bullet (x_{\star \star|0|}\rhd x_{|12|\star} )\bullet t_\B(x_{|1|0|})   \label{20nerve1}            \\
(\Lambda)\quad \quad x_{|02,00|}&= t_\B(x_{|2|0|})^{-1}\bullet (x_{\star \star|0|}\rhd x_{|02|\star})\bullet t_\B(x_{|0|0|})  \label{20nerve2}         \\
                  x_{|01,00|}&= t_\B(x_{|1|0|})^{-1}\bullet (x_{\star \star|0|}\rhd x_{|01|\star})\bullet t_\B(x_{|0|0|})     \label{20nerve3}           \\
(\Lambda)\quad \quad  x_{|02,00|}&=x_{|12,00|}\bullet x_{|01,00|}                                                       \label{20nerve4}          \\
                    x_{|02|\star}&=x_{|12|\star}\bullet x_{|01|\star}                                                   \label{20nerve5}
\end{align}
As usual we must show Equations~\ref{20nerve2} and \ref{20nerve4} are equivalent given the other three. We do this by showing the right hand sides of these equations are equal. We start with the right hand side of Equation~\ref{20nerve4} and make substitutions using Equation~\ref{20nerve1} and Equation~\ref{20nerve3}.
\begin{align*}
&x_{|12,00|}\bullet x_{|01,00|}   \\
=&t_\B(x_{|2|0|})^{-1}\bullet (x_{\star \star|0|}\rhd x_{|12|\star} )\bullet t_\B(x_{|1|0|})  \bullet t_\B(x_{|1|0|})^{-1}\bullet (x_{\star \star|0|}\rhd x_{|01|\star})\bullet t_\B(x_{|0|0|})    \\
=&t_\B(x_{|2|0|})^{-1}\bullet (x_{\star \star|0|}\rhd x_{|12|\star} )\bullet (x_{\star \star|0|}\rhd x_{|01|\star})\bullet t_\B(x_{|0|0|})  
\end{align*}
Then we use interchange for $\rhd$ and make a substitution using Equation~\ref{20nerve5} to finish the proof:
\begin{align*}
=&t_\B(x_{|2|0|})^{-1}\bullet (x_{\star \star|0|}\rhd (x_{|12|\star} \bullet x_{|01|\star}))\bullet t_\B(x_{|0|0|})\\
=&t_\B(x_{|2|0|})^{-1}\bullet (x_{\star \star|0|}\rhd x_{|02|\star})\bullet t_\B(x_{|0|0|})\\
\end{align*}
The $\langle 0,2\rangle$ case is similar.

For the $\langle 1,1\rangle$ case we  have the following commutativity conditions for the faces of a $\langle 1,1 \rangle$ cell $x_{|01|01|}.$
\begin{align}
x_{|11,01|}&= t_\B(x_{|1|1|})^{-1}\bullet (x_{\star|01|}\lhd x_{|1|\star\star} ) \bullet t_\B(x_{|1|0|}) \label{11nerve1}    \\
x_{|00,01|}&= t_\B(x_{|0|1|})^{-1}\bullet (x_{\star|01|}\lhd x_{|0|\star\star} ) \bullet t_\B(x_{|0|0|})  \label{11nerve2}  \\
(\Lambda)\quad \quad x_{|01,01|}&= x_{|11,01|}\bullet x_{|01,00|}     \label{11nerve3}   \\
(\Lambda)\quad \quad x_{|01,01|}&= x_{|01,11|}\bullet x_{|00,01|}     \label{11nerve4}         \\
x_{|01,11|}&=t_\B(x_{|1|1|})^{-1}\bullet (x_{\star \star |0|}\rhd x_{|01|\star} ) \bullet t_\B(x_{|0|1|}) \label{11nerve5}    \\
x_{|01,00|}&=t_\B(x_{|1|0|})^{-1} \bullet (x_{\star \star |0|}\rhd x_{|01|\star}) \bullet t_\B(x_{|0|0|})  \label{11nerve6} 
\end{align}
We must show Equations~\ref{11nerve3} and \ref{11nerve4} are equivalent given the other four equations. We show the right hand sides of these equations are equal. We start with the right hand side of Equation~\ref{11nerve3}, make substitutions using Equations~\ref{11nerve2} and \ref{11nerve5} then apply full interchange:
\begin{align*}
 &x_{|11,01|}\bullet x_{|01,00|}  \\
 =& t_\B(x_{|1|1|})^{-1}\bullet (x_{\star|01|}\lhd x_{|1|\star\star} ) \bullet t_\B(x_{|1|0|})\bullet t_\B(x_{|1|0|})^{-1} \bullet (x_{\star \star |0|}\rhd x_{|01|\star}) \bullet t_\B(x_{|0|0|}) \\
  =& t_\B(x_{|1|1|})^{-1}\bullet (x_{\star|01|}\lhd x_{|1|\star\star} ) \bullet (x_{\star \star |0|}\rhd x_{|01|\star}) \bullet t_\B(x_{|0|0|}) \\
  =& t_\B(x_{|1|1|})^{-1}\bullet (x_{\star \star |0|}\rhd x_{|01|\star} ) \bullet (x_{\star|01|}\lhd x_{|0|\star\star} ) \bullet t_\B(x_{|0|0|}) \\
  =& t_\B(x_{|1|1|})^{-1}\bullet (x_{\star \star |0|}\rhd x_{|01|\star} )\bullet t_\B(x_{|1|0|}) \bullet t_\B(x_{|1|0|})^{-1} \bullet (x_{\star|01|}\lhd x_{|0|\star\star} ) \bullet t_\B(x_{|0|0|})\\
\end{align*}
Finally substitutions using Equations~\ref{11nerve1} and \ref{11nerve6} yield the right hand side of Equation~\ref{11nerve4} as desired: $$=x_{|01,11|}\bullet x_{|00,01|}.$$

The last inner horns we must check are $\langle 3 \rangle$-horns.  We have the following commutativity conditions for the faces of a $\langle 4 \rangle$ cell $x_{|0123|}.$
\begin{align}
x_{|13|}= x_{|23|}\bullet x_{|12|}\\
(\Lambda) \quad \quad x_{|03|}= x_{|23|}\bullet x_{|02|} \label{3nerve2}  \\
(\Lambda) \quad \quad x_{|03|}= x_{|13|}\bullet x_{|01|}  \label{3nerve3} \\
x_{|02|}= x_{|12|}\bullet x_{|01|}
\end{align}
We must show Equations~\ref{3nerve2} and \ref{3nerve3} are equivalent given the other two equations. This is easily seen using the associativity of $\bullet$.
\end{proof}
\subsection{$N_\theta$ is functorial}
Let $\B$ and $\C$ be small fancy bicategories, and let  $F$ be a functor between these:
\begin{center}
\begin{tikzcd}
\olb \arrow{r}{\olf} \arrow{d}{t_\B} & \olc \arrow{d}{t_\C} \\
\wtb \arrow{r}{\wtf}                  & \wtc
\end{tikzcd}
\end{center}
We first define a functor $N_\theta(F)|^3_0:N_\theta(\B)|^3_0\ra N_\theta(\C)|^3_0.$

The Duskin nerve $N(\overline{F})$ defines a map $\psi^{\star}N_{\theta}(\B)\ra \psi^{\star} N_\theta(\C)$, defining $N_\theta(F)$ on $\langle \rangle$, $\langle 0\rangle$, $\langle 0,0\rangle$, and $\langle 0,0,0\rangle$ cells. On $\langle 1\rangle$ cells $N_\theta(F)$ is given by the map of $2$-morphisms given by $\wtf$. We must check that $N_\theta$ is defined for $\langle 2 \rangle$, $\langle 1,0\rangle$ and $\langle 0, 1 \rangle$ cells. We must check that $N_\theta(F)$ as defined above sends commutative spheres of these types in $N_\theta(\B)|^3_0$ to commutative spheres in $N_\theta(\C)|^3_0.$

A commutative $d\Theta_2\langle 2 \rangle$-sphere $S$ in $N_\theta(\B)|^3_0$ and its image $N_\theta(F)(S)$ have the following form: \begin{align*}S&:=[\theta, \theta\bullet \eta, \eta]\\ N_\theta(F)(S)&=[\wtf(\theta),\ \wtf(\theta\bullet \eta),\ \wtf(\eta)].\end{align*}
The commutativity of $N_\theta(F)(S)$ follows from the functoriality of $\wtf$ with respect to $\bullet$, \textbf{BFun2}.

A commutative $d\Theta_2\langle 1,0\rangle$-sphere $S$ in $N_\theta (\B)|^3_0$ and its image $N_\theta(F)|^3_0(S)$ by definition have the following form:
\begin{align*}
S&:=[(h,g',f'\vbar \gamma),\ (h,g,f\vbar \beta),\ t_\B(\gamma)^{-1}\bullet (h \rhd \eta) \bullet t_\B(\beta) ,\ \eta]\\ 
N_\theta(F)|^3_0(S)&= [(\olf(h),\olf(g'),\olf(f')\vbar \phi_{h,f'}\bullet \olf(\gamma)), \\ & \quad  (\olf(h),\olf(g),\olf(f)\vbar \phi_{h,f}\bullet  \olf(\beta)),\     \wtf(t_\B(\gamma)^{-1}\bullet (h \rhd \eta) \bullet t_\B(\beta)) ,\ \wtf(\eta)]
\end{align*}
$N_\theta(F)|^3_0(S)$ is commutative if and only if $$ t_\C(   \phi_{h,f'}\bullet \olf(\gamma)     )^{-1} \bullet (\olf(h) \rhd \wtf(\eta))\bullet t_\C(  \phi_{h,f}\bullet  \olf(\beta))  ) = \wtf(t_\B(\gamma)^{-1}\bullet (h \rhd \eta) \bullet t_\B(\beta)) $$
This is easily checked using the fact that $\olf$ and $\wtf$ commute with $t_\B$ and $t_\C$, and the naturality of the unitor of $\wtf$, which is $t_\C(\phi_{g,f}).$ 

The $\langle 0,1\rangle$ case is similar to the $\langle 0,1\rangle$ case. This defines $N_\theta(F)|^3_0$, and we define $$N_\theta(F):=\cosk^3(  N_\theta(F)|^3_0  ).$$ It is trivial to check that $N_\theta$ is functorial (respects composition of functors) given the fact that the Duskin nerve is functorial.
\section{$N_\theta$ and $\TBic$ are inverse equivalences of categories}
\subsection{The isomorphism $u: X \cong N_{\theta}( \TBic(X))$}
Since both $X$ and $N_{\theta}( \TBic(X))$ are $3$-coskeletal, it is enough to give an isomorphism $u':X|^3_0 \ra N_{\theta}( \TBic(X))|^3_0.$
This isomorphism is defined for $\langle \rangle$, $\langle 0 \rangle$, $\langle 0,0 \rangle$, and $\langle 0,0,0\rangle$ cells by the isomorphism $u:\psi^\star X \cong N(\Bic(\psi^\star X))=N_\theta(\overline{\TBic(X)}).$ defined in Section~\ref{bicisosec}. For $\langle 0,0\rangle$ cells, this map has the form $u'(x_{|0|0|})=(x_{\star|0|}, x_{|0,0|},x_{|0|\star} \vbar \underline{x_{|0|0|}}).$

The $\langle 1\rangle$-cells of $X$ and $N_{\theta}( \TBic(X))$ are identical, so we define $u'$ to be the identity for these cells. For a $\langle 2\rangle$ cell $x_{|012|}$ in $X$ we must have $$d(u'(x_{|012|}))=[x_{|12|},\ x_{|02|},\ x_{|01|}].$$ To show that $u'$ can be defined and is bijective on $\langle 2 \rangle$ cells, we must show this sphere commutes if and only if $x_{|012|}$  is a cell in $X$, or equivalently $$dx_{|012|}= [x_{|12|},\ x_{|02|},\ x_{|01|}]$$ is commutative in $X$. By definition of $N_\theta$, $d(u'(x_{|012|}))$ is commutative if and only if  $x_{|02|}=x_{12} \bullet x_{|02|}$, in $\widetilde{\TBic(X)}$. By the definition of $\bullet$ and the Matching Lemma for $X$, the sphere $d(x_{|012|})$ is also commutative if and only if $x_{|02|}=x_{12} \bullet x_{|02|}.$ This shows that $u'$ is well-defined and bijective for $\langle 2 \rangle $ cells.

Similarly, for a $\langle 1,0\rangle$ cell $x_{|01|0|}$ in $X$, we must have 
$$d(u'(x_{|01|0|})) = [ ( x_{\star|0|}, x_{|1,0|}, x_{|1|\star} \vbar \underline{x_{|1|0|}} ),\       ( x_{\star|0|}, x_{|0,0|}, x_{|0|\star} \vbar \underline{x_{|0|0|}} )              ,\ x_{|01,00|},\ x_{|01|\star} ]    $$
This sphere commutes if and only if $$   x_{|01,00|}    =\underline{x_{|1|0|}}\bullet (x_{\star|0|} \rhd x_{|01|\star}) \bullet  \underline{x_{|0|0|}} $$ in $\widetilde{\TBic(X)}$  which by Lemma~\ref{10celllemma} holds if and only if the sphere $$d(x_{|01|0|}  ) = [x_{|1|0|}, \ x_{|0|0|},\ x_{|01,00| }, x_{|01|\star}       ]$$
is commutative in $X$. This allows $u$ to be defined bijectively for $\langle 1,0\rangle$ cells.

The $\langle 0,1\rangle$ case is similar, making $u'$ an isomorphism. Finally we define $u= \cosk^3 u'$. 

The only non-trivial case to check for the naturality of $u$ is for $\langle 0, 0\rangle$ cells, which follows immediately from the fact that $u: X \cong N(\Bic(X))$ as defined in Section~\ref{bicisosec} is natural.

\subsection{The isomorphism $U$}
We now construct a natural isomorphism: $$U:\B  \cong   \TBic(  N_\theta(\B))    ,  $$ which is given by a square \begin{center}
\begin{tikzcd}
\olb \arrow{r}{\overline{U}} \arrow{d}{t_\B} & \overline{\TBic(  N_\theta(\B))} \arrow{d}{t'}) \\
\wtb \arrow{r}{\widetilde{U}}                  &  \widetilde{\TBic(  N_\theta(\B))}
\end{tikzcd}\end{center}

The map $\overline{U}$ is given by the strict isomorphism $U: \olb \cong \Bic(N(\mathcal{\olb}))=\overline{\TBic(  N_\theta(\B))}.$  Note that $\olb$ and $ \widetilde{\TBic(  N_\theta(\B))}$ have identical objects, $1$-morphisms, and $2$-morphisms, so we define $\widetilde{U}$ to be given by these identity maps with trivial unitor. We must check that this makes the square above commute, and that $ \widetilde{\TBic(  N_\theta(\B))}$ and $\wtb$ have the same vertical composition $\bullet$, whiskerings $\rhd$ and $\lhd,$ and identity for $1$-morphisms. 

The equivalence of the identity for $1$-morphisms in the two categories is trivial and left to the reader. To see that the square commutes, recall that by definition for $\beta:f \Rightarrow g:a\ra b$ in  $\overline{U}$ we have $\overline{U}(\beta)=\rho'_g \bullet \beta$ where $\rho'_g$ is the unitor of $\overline{U}$. So $t'(\overline{U}(\beta))=t'(\rho_{g'}\bullet \beta)$ makes the following sphere $\Theta_{t}(\rho_{g'}\bullet \beta)$ in $N_\theta  (\B)$ commutative:
$$[\Id_g, \ t'(\overline{U}(\beta)),\ s_0g =\rho'_g,\ \rho_{g'}\bullet \beta].$$ By the $\langle 0,1\rangle$ cell condition for $N_\theta(\B)$, this means
\begin{align*}t'(\overline{U}(\beta))&= t_\B(\rho'_g)^{-1} \bullet (\Id_g \lhd \id_a) \bullet t_\B(\rho'_g\bullet \beta)\\
t'(\overline{U}(\beta))&= t_\B(\beta)=\widetilde{U}(t_\B(\beta))
\end{align*}
showing that the square commutes.

The equivalence of the $\bullet$, $\lhd$ and $\rhd$ of $ \widetilde{\TBic(  N_\theta(\B))}$ and $\wtb$ are all immediately seen from the commutativity conditions for $\langle 2\rangle$, $\langle 0, 1\rangle$, and $\langle 1, 0\rangle$ cells in  $N_\theta(\B)$ respectively; details are left to the reader. The fact that $U$ respects composition and is natural follows trivial from the fact that $\overline{U}$ respects composition and is natural, which was shown in Section~\ref{summarysection}. $U$ is an isomorphism since both $\overline{U}$ and $\widetilde{U}$ are isomorphisms. 

\begin{theorem} \label{joyalthm}
The constructions $N_\theta$, and  $\TBic$ are inverse equivalences of categories between the category of $2$-reduced inner Kan algebraic $\Theta_2$-sets and the category of small fancy bicategories and strictly identity preserving functors. These constructions preserve strictness, and the natural isomorphisms exhibiting the equivalence are strict, these constructions also give an inverse equivalences of categories between the category of $2$-reduced inner-Kan algebraic $\Theta_2$-sets and strict morphisms and the category of small  fancy bicategories and strict functors.
\end{theorem}
\begin{remark}
The functor which forgets algebraic structures on $\Theta_2$-sets is an equivalence of categories, thus $N_{\Theta_2}$ is an equivalence of categories from fancy bicategories to (non-algebraic) $2$-reduced inner-Kan $\Theta_2$-sets.
\end{remark}
\chapter{Epilogue \label{epilogue}}
\section{Generalizations}
We have given several equivalences between categorical structures and presheaf categories with certain horn-filling conditions. Generalizing these presheaf categories allows us to suggest some definitions for certain types of higher categorical structures. 
\subsection{Verity $n$-fold categories}
The \emph{dimension} of an object $\mathbf{n}:=[n_1, \ldots, n_k]$ in $\Delta^m$ is given by $\dim(\mathbf{n})=\sum_{1\leq i \leq m} n_i.$

A \emph{(inner) coface} map in $\Delta^m$ is given by a map which is an (inner) coface map in $\Delta$ in component and the identity map in each other component. A \emph{universal inner horn} is obtained by removing an inner coface map from a representable presheaf (i.e. simplex). Note that any universal inner horn obtained from $\mathbf{n}$ in $\Delta^m$ is clearly nice, i.e. has minimal complementary dimension\footnote{See Definition~\ref{mincompdimdef}} $\dim(\mathbf{n})-1.$

As usual an inner horn in a $\Delta^m$-set $X$ is a map from a universal inner horn to $X$, and $X$ is called \emph{inner-Kan} if every inner horn $H$ in $X$ has a filler, i.e. an extension of $H$ along the natural inclusion of the universal horn in its corresponding representable presheaf. $X$ is called \emph{$k$-reduced} if every inner horn of minimal complementary dimension $k$ or greater has a unique filler.

\begin{definition}
A \emph{Verity $(\infty, n)$-fold category} is a inner-Kan $\Delta^n$-set. A \emph{small Verity $n$-fold category} is a $n$-reduced inner-Kan $\Delta^n$-set.
\end{definition}
\begin{definition}
We call a $\Delta^n$-set $X$ \emph{horizontally trivial} if each $k$-simplicial set $$X_{m_{1} \ldots m_{n-k-1} 0\underbrace{\bullet \ldots \bullet}_\text{$k$-times}  }$$ is trivial, i.e. has only identity maps. A \emph{Barwick-like small fancy $(\infty,n)$-category} is a horizontally trivial small Verity $(\infty, n)$-fold category, and a \emph{Barwick-like small fancy $n$-category} is a horizontally trivial small Verity $n$-fold category.
\end{definition}
The reference to Barwick is to his notion of ``$n$-fold Segal spaces'', defined in \cite{Bar05}. Refer to \cite{BSP11} for an exposition of this definition, or to \cite{Lur09b}. Our definition makes it so that there is a close analogy:
\begin{itemize}
\item Quasicategories are to (not necessarily complete) Segal spaces as Barwick-like small fancy  $(\infty,n)$-categories are to (not necessarily complete) $n$-fold Segal spaces. 
\end{itemize}
In Chapter~\ref{bicatchapter} we showed a Barwick-like fancy $2$-category is equivalent to a small fancy bicategory, so our definition gives a generalization of this concept. 
\begin{remark}
In Chapter~\ref{bicatchapter} we compared small fancy bicategories to \emph{vertically} trivial bisimplicial sets, whereas we have adopted the opposite convention above, using horizontally trivial bisimplicial sets. The reason for this notational inconsistency stems from our conformance in that chapter to the usual usage of ``vertical'' and ``horizontal'' 1-morphisms for pseudo-double categories. We beg our gentle readers' forbearance.
\end{remark}
\subsection{Verity $n$-fold categories with thin structure}
Recall Definition~\ref{algbithin}:
\begin{definition} An \emph{$2$-reduced inner-Kan bisimplicial set with thin structure} consists of:
\begin{enumerate}
\item An $2$-reduced inner-Kan  bisimplicial set $X$
\item An $2$-reduced inner-Kan simplicial set $T$
\item A morphism of $2$-reduced inner-Kan bisimplicial sets $t:h_2 \partial_* T \ra X$ which is an isomorphism when restricted to zeroth row and column. 
\end{enumerate} 
A morphism  $F:(X,T,t) \ra (Y,S,s)$ is given by the obvious commutative square. 
\end{definition}
The non-reduced version of this definition, called a \emph{inner-Kan bisimplicial set with thin structure}, is found by relaxing the condition that $T$ and $X$ are $2$-reduced and replacing $t$ with a map $\partial_* T \ra X$.

\begin{definition}  Let $f$ be a morphism $\underline{m}\ra \underline{n}$ of the Segal category $\Gamma$ defined in Definition~\ref{gammadef}. We define a \emph{compilation along $f$} functor $D_f: \cat{Set}_{\Delta^m} \ra \cat{Set}_{\Delta^n}$ as follows: $$D_f(X)_{i_1\ldots i_n}:= \Hom\left( \left(\prod_{k\in f(1)} \Delta[i_k]\right)\square \left(\prod_{k\in f(2)} \Delta[i_k]\right)\square \ldots \square  \left(\prod_{k\in f(m)} \Delta[i_k]\right), X \right).$$ This $D_f$ is functorial in an  obvious way.
\end{definition}
\begin{example}~

\begin{itemize}
\item $D_{10}: \cat{Set}_{\Delta^2} \ra \cat{Set}_{\Delta^2}$ is the functor that switches the two indicies of a bisimplicial set.
\item $D_{0\emptyset}$ and $D_{\emptyset 0}$ are functors $\cat{Set}_{\Delta^2} \ra \cat{Set}_{\Delta}$ which give the zeroth column and zeroth row of a bisimplicial set.
\item $D_{0}$ and $D_{1}$ are the functors $\mathbf{const}^2_1$ and $\mathbf{const}^2_2$ from $\cat{Set}_{\Delta}$ to $\cat{Set}_{\Delta^2}$ as defined in Definition~\ref{constdef}, which make a bisimplicial set which is constant in one direction.
\item $D_{(01)}$ is the diagonal functor $\partial_\star: \cat{Set}_{\Delta}\ra \cat{Set}_{\Delta^2}$ defined in Definition~\ref{diagonaldef}.
\end{itemize}
\end{example}
If $f$ and $g$ are composable morphisms of $\Gamma$, there is a natural isomorphism $D_gD_f \cong D_{g\circ f}$, which can be used to make a $D$ a bicategory functor from $\Gamma$, viewed as a bicategory with trivial $2$-morphisms, to $\cat{Cat}$, the category of small categories. We call $D$ the \emph{simplicial compilation functor}

Let $F$ be a bicategory functor from a small category $\CalC,$ viewed as a category with trivial $2$ morphisms, to $\cat{Cat}.$ Recall the \emph{Grothendieck construction} $\CalC \int F$ is a category whose objects are pairs $(c, x)$ where $c \in \CalC$ and $x \in G(c)$. A morphism $(c,x)\ra (c',x')$ is a pair $(f,\varphi)$ where $f:c\ra c'$ in $\CalC$ and $\phi:F(f)(x)\ra x'$. For $$(c,x)\stackrel{(f,\varphi)}{\ra}(c',x')\stackrel{(f',\varphi')}{\ra}(c'',x'')$$ the composition is defined by $$(f'\circ f,\ \ \varphi'\circ F(f')(\varphi)\circ (\phi_{f',f})_x)$$
where $\phi_{f',f}$ is the unitor of $F.$ The following diagram clarifies the map on the right above:
\begin{center}
\begin{tikzcd}
F(g\circ f)(x) \arrow{r}{(\phi_{f',f})_x } & F(f')\circ F(f)(x) \arrow{r}{F(f')(\varphi)} & F(f')(x') \arrow{r}{\varphi'} & x''
\end{tikzcd}
\end{center}
There is an obvious forgetful functor $p: \CalC \int F \ra \CalC.$ The category of sections $\mbox{Sec}(\CalC \int F)$ is the category of functors $S:\CalC \ra \CalC \int F$ such that $p\circ S$ is the identity, and natural transformations between them. Explicitly, an element $S$ of $\mbox{Sec}(\CalC \int F)$ consists of
\begin{itemize}
\item For each $c \in \CalC$ an object $S(c)$ in $F(c)$
\item For each morphism $f:c \ra c'$ in $\CalC$ a map $S(f) : F(f)(S(c))\ra S(c')$ 
\end{itemize}
such that $S(\id_a)$ is the unitor of $\upsilon_a$ of $F$ and for $f:c \ra c'$ and $f': c' \ra c''$, the following ``composition condition" diagram commutes
\begin{center}
\begin{tikzcd}
F(f'\circ f)(S(c)) \arrow{r}{(\phi_{f',f})_{S(c)}} \arrow{rrd}[below]{S(g\circ f)} & F(f')\circ F(f)(S(c)) \ \arrow{r}{F(f')(S(f))} &\ F(f')(S(c')) \arrow{d}{S(f')} \\
&            & S(c'')
\end{tikzcd}
\end{center}
Recall the wreath product operation $\Delta \wreath$ from Definition~\ref{wreathdef}. Let $\star$ be the category with one object and an identity morphism. Then $\Delta \wreath \star$ is canonically equivalent to the image of the Segal functor $\phi:\Delta \ra \Gamma$. The image of this functor is the subcategory consisting of morphisms $f:\underline{m} \ra \underline{n}$ which are coincreasing (meaning that if $i > j$ then each element of $f(i)$ is greater than each element of $f(j)$), and coconsecutive (meaning that the image $\bigcup_{i\in \underline{m}} f(i)$ of $f$ is a set of consecutive integers.)

Restricting the map $\phi:\Delta \ra \Gamma$ to the full subcategory $\Delta|^n_0$ on $[0], [1],\ldots, [n]$, we get an operation $\Delta|^n_0 \wreath $ with $\Delta|^n_0 \wreath \star$ being the full subcategory of the image of $\phi$ on $\underline{0},\underline{1}, \ldots \underline{n}$. So the simplicial compilation functor $D$ can be restricted to a bicategory functor $D|^n_0: \Delta|^n_0\wreath \star \ra \cat{Cat}$.

Let $\mathcal{Q}$ denote the full subcategory of $\mbox{Sec}(\Delta|^2_0\wreath \star \int D|^2_0)$ on sections $S$  for which each $S(\underline{1})$ and $S(\underline{2})$ are inner-Kan. Then an object $S\in \mathcal{Q}$ is yields a inner-Kan bisimplicial set with thin structure, with $X:=S(\underline{2})$ and $T:= S(\underline{1}).$ The map $t:\partial_\star T \ra X$ is given by $$S((01)):D_{(01)}(T)=\partial_\star T \ra X$$ and the condition that this map is an isomorphism when restricted to zeroth row and column is ensured by applying the composition condition for objects of $\mbox{Sec}(\Delta|^2_0\wreath \star \int D|^2_0)$ to the identities $$\id_{\underline{1}}=0\emptyset \circ (01)=\emptyset 0 \circ (01)$$ in $\Delta|^2_0\wreath \star.$ In fact it is not hard to check that this gives an equivalence between $\mathcal{Q}$ and the category of inner-Kan bisimplicial set with thin structure.

Similarly let $\mathcal{R}$ denote the full subcategory of $\mbox{Sec}(\Delta|^2_0\wreath \star \int D|^2_0)$ on sections $S$  for which each $S(\underline{1})$ and $S(\underline{2})$ are $2$-reduced inner-Kan. An object of $\mathcal{R}$ yields a $2$-reduced inner-Kan bisimplicial set with thin structure in the same manner, except that the map $t: h_2 \partial_\star T \ra X$ is found by taking the adjoint to the map $$S((01)):D_{(01)}(T)=\partial_\star T \ra X,$$ using the fact that $h_2$ is left adjoint to the inclusion of the category of $2$-reduced inner-Kan bisimplicial sets in the category of inner-Kan bisimplicial sets. Again, this gives an  equivalence between $\mathcal{R}$ and the category of $2$-reduced inner-Kan bisimplicial set with thin structure.

We now are in a position to suggest a generalization to the notion of an inner-Kan bisimplicial set with thin structure:
\begin{definition}
The category of \emph{($k$-reduced) inner-Kan $n$-simplicial sets with thin structure} is the full subcategory of  $\mbox{Sec}(\Delta|^n_0\wreath \star \int D|^n_0)$ on sections $S$  for which each $S(\underline{i})$ is ($k$-reduced) inner-Kan for $i>0$. 
\end{definition}

If the mildly conjectural Corollary~\ref{foldingcor} holds, then $n$-reduced inner-Kan $n$-simplicial sets with thin structure are another generalization of small fancy bicategories, along with a Barwick-like fancy $n$-categories. It is unclear if these are equivalent for $n>2$. 

Similarly, we suggest that  inner-Kan $n$-simplicial sets with thin structure are a possible model for small fancy $(\infty,n)$-categories.

\subsection{Joyal-like small fancy $n$-categories }
\begin{definition}
For $[n]\in \Delta$ we define $\dim([n])=n$. For $\Delta\wreath \CalC$ where $\CalC$ is a dimensional category we define $\dim([c_1,\ldots, c_m])= \dim(c_1)+\dim(c_2)+\ldots +\dim(c_m)+m.$ This inductively defines $\dim$ for $\Theta_n =\Delta\wreath \Delta \wreath \ldots \wreath \Delta$. A \emph{coface} map in $\Theta_n$ is a monic map which increases dimension by $1$.
\end{definition}

Recall that a coface map $f:[n]\ra[m]$ in    $\Delta=\Theta_1$ is called \emph{inner} if $0$ and $m$ are in the image of $f$. Inductively, we say a coface map $f$ in $\Theta_n =  \Delta \wreath \Theta_{n-1}$ is \emph{inner} if the type of $f$ is inner and all components of $f$ are inner. This definition is consistent with Definition~\ref{innerdef} in the case $n=2$. 

A \emph{universal inner horn} is $\Theta_n$-set obtained by removing an inner coface from a sphere\footnote{See Definition~\ref{mincompdimdef}}. An inner-horn in a $\Theta_n$-set $X$ is a map from a universal inner horn to $X$. A $\Theta_n$-set $X$ is called inner-Kan if every inner horn in $X$ has a filler. An inner-Kan $\Theta_n$-set $X$ is called \emph{$k$-reduced} if every inner horn in $X$ of minimal complementary dimension $k$ or greater has a unique filler.

\begin{definition}
We call an $n$-reduced inner-Kan $\Theta_n$-set a \emph{Joyal-like fancy $n$-categories}. An inner-Kan $\Theta_n$-set will be called a \emph{Joyal small fancy $(\infty,n)$-category}. 
\end{definition} 
By Theorem~\ref{joyalthm}, the above definition gives a generalization of small fancy bicategories.

\subsection{Barwick-like  small symmetric monoidal fancy $n$-categories}
A \emph{simplicial coface} map in $\Gamma\times\Delta^n$ is given by a map which is an coface map in a $\Delta$ component and the identity map in each other component.  A \emph{$\Gamma$-coface} is a coface map in a $\Gamma$ component and the identity map in each other component. A \emph{universal inner horn} is a $\Gamma \times \Delta^n$ obtained by either removing a simplicial inner coface map from a representable presheaf (i.e. simplex), or by removing a set of $\Gamma$-coface maps which correspond to a set of coface maps which are removed from a representable $\Gamma$-set make an inner universal horn, as described in Definition~\ref{innerkangammadef}.

As usual an inner horn in a $\Gamma\times\Delta^m$-set $X$ is a map from a universal inner horn to $X$, and $X$ is called \emph{inner-Kan} if every inner horn $H$ in $X$ has a filler. $X$ is called \emph{$k$-reduced} if every inner horn of minimal complementary dimension $k$ or greater has a unique filler.

\begin{definition}
We call a $\Gamma\times \Delta^n$-set $X$ \emph{horizontally trivial} if each $k$-simplicial set $$X_{m_{0} \ldots m_{n-k-1}0 \underbrace{\bullet \ldots \bullet}_\text{$k$-times} }$$ is trivial, i.e. has only identity maps. A \emph{Barwick-like small symmetric monoidal fancy $(\infty,n)$-category} is a horizontally trivial inner-Kan $\Gamma\times \Delta^n$-set. A \emph{Barwick-like small symmetric monoidal fancy $n$-category} is a $n+1$-reduced vertically trivial inner-Kan $\Gamma\times \Delta^n$-set. 
\end{definition}

%
\section{Model category speculations}
Dimitri Ara \cite{Ara12} and Harry Gindi \cite{Gin12} independently construct a model structure on the category $\cat{Set}_{\Theta_n}$ which is Quillen equivalent to Rezk's \cite{Rez10} model structure for complete Rezk spaces, which is a model structure on the category of $\Theta_n$-spaces and is equivalent to other model categories for $(\infty,n)$-categories. In this construction, Ara and Gindi modify an idea of Joyal and Cisinski. However, an inner-Kan $\Theta_2$-set is not necessarily fibrant in this model structure. For instance, if $\B$ is the bicategory consisting of two parallel $1$-morphisms and an inverse pair  of $2$-morphisms between them, Remark 5.27 of \cite{Ara12} shows that the $\Theta_2$-set $N_\theta(\ulcorner \B \urcorner)$ is not fibrant in this model structure. In light of Theorem~\ref{gammasummary}, this is not surprising, because we expect inner-Kan $\Theta_2$-sets to model fancy $(\infty,2)$-categories and not $(\infty,2)$-categories.

This leads us to conjecture:
\begin{conjecture} A $\Theta_2$-set $X$ is fibrant in the Ara-Gindi model structure if and only if it is a $2$-reduced inner-Kan and $\TBic(X)$ is a complete fancy bicategory. 
\end{conjecture}
\begin{conjecture}\label{araconj}
More speculatively, we suggest that there is a model structure $M_{\Theta_n}$ on $\cat{Set}_{\Theta_n}$ such that the inner-Kan $\Theta_n$-sets are the fibrant objects and such that the two adjunctions given in \cite{Ara12} between $\cat{Set}_{\Theta_n}$ and $\Theta_n$-spaces are Quillen equivalences between $M_{\Theta_n}$ and Rezk's model structure $\mathrm{Se}_C$ whose fibrant objects are his (not necessarily complete) ``Segal objects''.
\end{conjecture}
One of the adjunctions mentioned in Conjecture~\ref{araconj}  is easy to define: 
\begin{definition} Let $\C$ be a small category with a terminal object $0$. There is a pair of functors $$i_\C: \C \leftrightarrow \C \times \Delta : p_\C$$ where $p_\C$ is the obvious projection and $i_\C(x)= (x,0)$. The pullbacks $$p_\C^*:\cat{Set}_{\C}\leftrightarrow \cat{Set}_{\C \times \Delta}:i_\C^*$$ are an adjoint pair, with $i_\C^*$ being the right adjoint. 
\end{definition}
Note that the category $\cat{Set}_{\C \times \Delta}$ can equivalently be viewed as the category of $\C$-spaces, i.e. presheaves of simplicial sets on $\C$.  In the case $\Theta_n=\C$ the adjoint pair $(i^*_{\Theta_n},p^*_{\Theta_n} )$ is one of two Quillen equivalences given by Ara and Gindi.


In the case of $\Gamma$-sets, Theorem~\ref{gammasummary} provides support the idea that inner-Kan $\Gamma$-sets are a model for  small symmetric monoidal $(\infty,0)$-categories. We propose the following analogy:
\begin{itemize}
\item Quasicategories are to Segal categories as inner-Kan $\Gamma$-sets are to special $\Gamma$-spaces.
\end{itemize}
A \emph{Segal category} is a Segal space whose "space of objects" $X_0$ is discrete. See section 5 of \cite{JT07} for details on the relation between quasicategories and Segal categories.

A \emph{special $\Gamma$-space}, defined in \cite{Seg74}, is a presheaf of topological spaces on $\Gamma$ meeting certain conditions. Note that these are simply called $\Gamma$-spaces by Segal. Special $\Gamma$-spaces provide a model for grouplike $E_\infty$ spaces where the multiplication is ``geometrically'' defined, that is, only defined up to coherent homotopy. These $\Gamma$-spaces are used as models for the notion of small symmetric monoidal $(\infty,0)$-categories, for instance see remark 2.4.2.2 in \cite{Lur11}.
This leads us to conjecture:
\begin{conjecture} There is a model structure $M_{\Gamma}$ on $\cat{Set}_\Gamma$ whose fibrant objects are the inner-Kan $\Gamma$-sets and such that  $(i^*_{\Gamma},p^*_{\Gamma} )$ is a Quillen equivalence between $M_{\Gamma}$ and the model structure on $\Gamma$-spaces whose fibrant objects are the special $\Gamma$-spaces, which was first defined in \cite{BF77}.
\end{conjecture}
\section{The singular $\Gamma$-set of a pointed topological space}
If inner-Kan $\Gamma$-sets model $E_\infty$-spaces, then Proposition~\ref{Kanprop} supports the idea that Kan $\Gamma$-sets\footnote{see Definition~\ref{kandef}} serve as a model for grouplike $E_\infty$-spaces. In the next subsection, we hypothesize that Kan $\Gamma$-spaces have a close relation to connected infinite loop spaces, which are another model for grouplike $E_\infty$-spaces.
\subsection{Kan $\Gamma$-sets vs. infinite loop spaces}
The overall idea of this subsection, as well as the constructions of our realization and singular $\Gamma$-set functors are analogous to corresponding ideas from Segal's famous paper \cite{Seg74}, but with $\Gamma$-sets replacing Segal's $\Gamma$-spaces.

Let $\cat{Top}_\star$ denote the category of pointed topological spaces and let $(K,e)\in \cat{Top}_\star.$ 
\begin{definition}
The functor $R_{(K,e)}:\Gamma \ra \cat{Top_\star}$, also denoted $R_K$, is defined as follows
\begin{itemize}
\item $R_{K}(\underline{n}) = K^n$
\item For $f: \underline{n}\ra \underline{m}$ we define $$R_K(f)(k_1,\ldots, k_n)_j=\begin{cases} k_i & j\in f(i) \\ e & j \notin \bigcup_{i\in \underline{n}} f(i).
   \end{cases} $$
\end{itemize}
\end{definition}
There is a unique-up-to-isomorphism colimit-preserving \emph{$K$-realization} functor $\cat{Set}_\Gamma \ra \cat{Top}_\star$ which is isomorphic to $R_K$ when precomposed with the Yoneda embedding, and we abuse notation slightly by also denoting this functor $R_K$. This functor has a right adjoint \emph{$K$-singular} functor $\Sing_K$ which is given by the formula $$\Sing_K(X)(\underline{n})=\cat{Top}_\star(R_K(\Gamma[n]),X).$$

We will focus on the functors $R_{(S^1,0)}$ and $\Sing_{S^1}$ where $(S^1,0)$ is the standard pointed circle. The following proposition serves to contextualize this choice:
\begin{proposition}
There is a natural isomorphism of functors between $R_{(S^1,0)}:{\cat{Set}_\Gamma}_\star \ra \cat{Top}_\star$ and the functor $X\ra (|\phi^\star(X)|, s_0\Box)$, where $(S^1,0)$ is the standard pointed circle and ${\cat{Set}_\Gamma}_\star$ denotes the subcategory of $\cat{Set}_\Gamma$ consisting of $\Gamma$-sets $Y$ such that $Y(\underline{0})$ consists of a single point denoted $\Box$.
\end{proposition}
\begin{proof} This may be proved by constructing the isomorphism explicitly for representable presheaves, which is left to the reader as an exercise.
\end{proof}
\begin{definition}
A \emph{delooping} of pointed topological space $(Y,e)$ is a pointed space $(Y_1,e_1)$ together with a weak homotopy equivalence $(Y,e) \simeq \Omega (Y_1,e_1).$ We say $Y$ is an \emph{infinite loop space} if it has a delooping and its delooping has a delooping and so forth ad infinitum. 
\end{definition}
Let $(Y,e)$ be a pointed, connected topological space. The $\Gamma$-set $\Sing_{S^1}(Y)$ in general need not be either inner-Kan or Kan. For instance, applying $R_{S^1}$ to universal horn inclusion $\Lambda^2_{1|2}\ra  \Gamma[2]$ we get the canonical inclusion $S^1 \vee S^1 \ra S^1 \times S^1$. Applying the adjunction between $R_{S^1}$ and $\Sing_{S^1}$ we see that $\Sing_{S^1}(Y)$ has fillers for $\Lambda^2_{1|2}$-horns if and only if any pointed map $S^1 \vee S^1 \ra Y$ has an extension to $S^1\times S^1$, which holds if and only if $\pi_1(Y)$ is abelian, which is necessary for $Y$ to have a delooping. Based on this fact, we make the following fairly speculative conjecture:
\begin{conjecture}
Let $(Y,e)$ be a pointed, connected topological space. Then $Y$ is an infinite loop space if and only if $\Sing_{S^1}(Y)$ is Kan. Likewise, if $X$ is a $\Gamma$-set, then $R_{S^1}(X)$ is an infinite loop space if and only if $X$ is Kan.
\end{conjecture}

\section{Fancy bicategory theory \label{fancysec}}
While our motivation for introducing fancy bicategories was to compare with certain presheaf categories generalizing the notion of quasicategories, fancy bicategories have some interesting utility in their own right. We described in Definition~\ref{fancificationdef} how both the category of bicategories and weak functors and the category of strict bicategories and strict functors are subcategories of the category of fancy bicategories in a natural way:
\begin{definition} \label{fancificationdef2}
For a bicategory $\B$ the \emph{complete fancification $\llcorner\B\lrcorner$} takes $\widetilde{\llcorner\B\lrcorner}:=\B$ and $\overline{\llcorner\B\lrcorner}$ to be the subcategory of $\B$ consisting of all objects and $1$-morphisms and all invertible $2$-morphisms. We also denote the same fancy bicategory as $F^{\downarrow}(\B)$. 

For a strict bicategory $\B$, the \emph{sparse fancification $\ulcorner\B\urcorner$ } takes $\widetilde{\ulcorner\B\urcorner}:=\B$  and $\overline{\ulcorner\B\urcorner}$ to be the objects and $1$-morphisms of $\B$ together with identity $2$-morphisms. In this case there is an obvious canonical functor $\iota_\B: \ulcorner\B\urcorner \ra \llcorner\B\lrcorner.$ We also denote the same fancy bicategory as $F^{\uparrow}(\B)$. 
\end{definition}
We will see in this section that this unification of  weak bicategory theory and strict bicategory theory can be extended further, defining notions that unify strict and weak transformations and modifications of functors, along with a concept of $2$-limit for fancy bicategories that unifies two standard notions $2$-limit.

\subsection{Bicategory transformations and modifications}
Let $\B$, $\C$ be bicategories, and let $F,G:\B\ra \C$ be bicategory functors.  A \emph{pseudonatural transformation}, $u:F\Rightarrow G$ consists of a $1$-morphism \emph{components} $u_b:F(b)\ra G(b)$ for each object $b\in \B$ and an invertible \emph{naturalizer} $2$-morphism $u_f$ for each $1$-morphisms $f:b\ra b'$ of $\B$, fitting in the following square:
\begin{center}
\begin{tikzpicture}[scale=1.8,auto]
\begin{scope}
\node (10) at (1,1) {$F(b')$};
\node (00) at (0,1) {$F(b)$};
\node (11) at (1,0) {$G(b')$};
\node (01) at (0,0) {$G(b)$};
\node[rotate=-135] at (.5,.5){$\Rightarrow$};
\node at (.35,.65){$u_f$};
\path[->] (00) edge node[midway,scale=.8]{$F(f)$}(10);
\path[->] (00) edge node[midway,swap]{$u_b$}(01);
\path[->] (01) edge node[midway,swap,scale=.8]{$G(f)$}(11);
\path[->] (10) edge node[midway]{$u_{b'}$}(11);
\end{scope}
\end{tikzpicture}
\end{center}
A pseudonatural transformation must satisfy relations expressing the compatibility of the naturalizers with the unitors and distributors of $F$ and $G$, as well as a naturality condition with respect to $2$-morphisms of $\B$.

For transformations $u,v:F\Rightarrow G$, a \emph{modification} $\frakm:u \Rightarrow v$ consists of a \emph{component} $2$-morphism $\frakm_b :u_b \Rightarrow v_b$ which commutes with the naturalizers of $u$ and $v$ in the obvious sense. These notions of pseudonatural transformation and modifications make $\Fun(\B,\C)$ into a strict bicategory. If $F$ and $G$ are strict functors of bicategories, a \emph{strict transformation} $u:F\Rightarrow G$ is called \emph{strict} if its naturalizer contains only identities. If $\B$ and $\C$ are strict bicategories, we write $\SFun(\B,\C)$ for the strict bicategory of strict functors, strict transformations, and modifications.

\subsection{Fancy transformations and modifications}

Let $\B$, $\C$ be fancy bicategories. We have defined a functor $F:\B\ra \C$  in Definition~\ref{fancydef} as consisting of bicategory functors $\olf:\olb \ra \olc$ and $\wtf:\wtb \ra \wtc$, such that the following diagram commutes strictly:
\begin{center}
\begin{tikzcd}
\olb \arrow{r}{\olf} \arrow{d}{t_\B} & \olc \arrow{d}{t_\C} \\
\wtb \arrow{r}{\wtf}                  & \wtc
\end{tikzcd}
\end{center}

We can define transformations and modifications between these fancy bicategory functors giving the fancy functors $\B\ra\C$  the structure of a fancy bicategory   $\{\B,\C \}$. 

Let $F,G: \B \ra \C$ be fancy bicategory functors. A \emph{fancy natural transformation} $u:F\Rightarrow G$ is a transformation $\overline{u}:\overline{F} \Rightarrow \overline{G}.$ If we apply $t_\C$ to the components and naturalizer of $\overline{u}$, we get a fancy natural transformation $\widetilde{u}:\widetilde{F} \Rightarrow \widetilde{G}.$ These fancy natural transformations are the $1$-morphisms of  $\{\B,\C \}$. A $2$-morphism $u\Rightarrow v$ of $\overline{\{\B,\C \}}$ is called a \emph{thin modification} and is given by a modification $\overline{u} \Rightarrow \overline{v},$  whereas a $2$-morphism $u\Rightarrow v$ of $\widetilde{ \{\B,\C\}}$ is called a \emph{thick modification} and is given by a modification $\widetilde{u} \Rightarrow \widetilde{v}.$ Applying $t_\C$ to components gives a map from thin modifications to thick modifications, giving the thin structure map $t:\overline{ \{\B,\C\}}\ra \widetilde{ \{\B,\C\}}.$ 
\begin{definition}
Fancy bicategory functors have an obvious composition, which can be extended to give a fancy bicategory functor $$\circ :\{\mathcal{A},\B\} \times \{\B,\C\} \ra \{\mathcal{A},\C\}.$$
In particular, $F\in \{\A,\B\}$ and $G\in\{\B\ra \C\}$  induce functors denoted
\begin{align*}
\presuper{\circ}{F} :\{\B,\C\} \ra \{\A,\C\} \\
G^\circ : \{\A,\B\} \ra \{\A,\C\}.
\end{align*}
\end{definition}

\subsection{Fancy categories}
\begin{definition}
A (small) fancy category $B$ consists of a (small) groupoid $\overline{B}$ and a (small) category $\widetilde{B}$ together with a map $t_B: \overline{B} \ra \widetilde{B}$ which is an isomorphism on the set of objects. A \emph{functor} $F$ of fancy categories is a square:
\begin{center}
\begin{tikzcd}
\overline{B} \arrow{r}{\olf} \arrow{d}{t_B} & \overline{C} \arrow{d}{t_C} \\
\widetilde{B} \arrow{r}{\wtf}                  & \widetilde{C}
\end{tikzcd}
\end{center}
A \emph{thin} transformation $u:F \Rightarrow G$ is a natural transformation $\overline{u}: \overline{F}\Rightarrow \overline{G}$ and a \emph{thick} transformation is a natural transformation $\widetilde{u}: \widetilde{F}\Rightarrow \widetilde{G}.$   We define a fancy bicategory $\TCat$ with objects small fancy categories and $1$-morphisms functors between them. $2$-morphisms of $\overline{\TCat}$ are thin natural transformations and $2$-morphisms of $\widetilde{\TCat}$ are thick natural transformations. Applying $t_C$ to components gives a map from thin transformations $F\Rightarrow G$ to thick transformations $F\Rightarrow G$, giving the thin structure map $t$ of $\TCat$.
\end{definition}
If $\B$ is a fancy bicategories, and $a,b$ are objects in $\B$, then its easy to see that the morphisms $\B(a,b)$ can be given the structure of a fancy category in a natural way.  In fact, the construction $(a,b)\ra \B(a,b)$ can be extended to give a functor $\fancyhom_\B: \B^{\op}\times \B \ra \TCat$ where $\B^{\op}$ is the $1$-cell dual of $\B$, reversing the $1$-morphisms of $\B$ but not either set of $2$-morphisms.

\begin{definition}
Given a small category $B$, the \emph{complete fancification $\llcorner B\lrcorner$} is constructed by taking $\widetilde{\llcorner B\lrcorner}:=B$ and $\overline{\llcorner B\lrcorner}$ to be the subcategory of $B$ consisting of all objects and all invertible morphisms. This operation extends to a strictly full and faithful functor of fancy bicategories,  $$\cfancya:\llcorner\cat{Cat}\lrcorner \ra \TCat .$$ A fancy category is called \emph{complete} if it is isomorphic to a fancy bicategory in the image of the complete fancification functor.

Given a small category $B$, the \emph{sparse fancification $\ulcorner B\urcorner$ } is formed by taking $\widetilde{\ulcorner B\urcorner}:=B$   and $\overline{\ulcorner B\urcorner}$ to be the objects  together with identity morphisms. This extends to a strictly full and faithful functor of fancy bicategories, denoted $$\sfancya:\ulcorner\cat{Cat}\urcorner \ra \TCat.$$ A fancy category $B$ will be called \emph{sparse} if it is isomorphic to a fancy bicategory in the image of this functor, equivalently, if $\overline{B}$ has only identity $2$-morphisms.
\end{definition}
\subsection{Basic fancified bicategory theory}
\begin{definition}~
\begin{itemize}
\item A $1$-morphisms $f:a\ra b$ in a fancy bicategory $\B$ is called an \emph{equivalence} if it is an equivalence in $\overline{\B}.$ 
\item Let $F: \B \ra \C$ be a functor of fancy bicategories. $F$ induces a strict natural transformation $$\fancyhom_\B \Rightarrow \fancyhom_\C \circ  (F^{\op}\times F),$$
in particular there is a functor of fancy categories $\B(a,b)\ra \C(F(a),F(b))$. We say $F$ is \emph{strictly fully faithful} if this transformation is an isomorphism.
\item The category of fancy bicategories has a terminal object, denoted $\star$, where $\overline{\star}$ and $\widetilde{\star}$ both consist of a single object, identity $1$-morphism and identity $2$-morphism.
\item If $b$  in a bicategory or fancy bicategory $\B$, there is a functor $\lambda_b^\B:\star \ra \B$ which sends the unique object to $b$, where $\star$ denotes either the terminal bicategory or terminal fancy bicategory.
\item For fancy bicategories $\B$ and $\C$ there is an obvious cartesian product $\B \times \C$.
\item  If $\A$, $\B$, $\C$ are fancy bicategories, we have a strict isomorphism of fancy bicategories $$\fancyadj_{\A,\B,\C}:\{\A\times \B, \C \}\ra \{\A, \{\B,\C\}\}.$$
\item  If $\A$, $\B$, $\C$ are bicategories, we have a strict isomorphism of bicategories $$\Adj_{\A,\B,\C}:(\A\times \B, \C )\ra (\A, (\B,\C)).$$
\item  If $\A$, $\B$, $\C$ are strict bicategories, we have a strict isomorphism of bicategories $$\Adjst_{\A,\B,\C}:[\A\times \B, \C ]\ra [\A, [\B,\C]].$$
\end{itemize}
\end{definition}

The complete and sparse fancification constructions $\B \ra F^{\downarrow}(\B) = \llcorner \B \lrcorner$ and $\B \ra F^{\uparrow}(\B) = \ulcorner \B \urcorner$ can be extended to a functors from a ``fancy tricategory'' of bicategories to the fancy tricategory of fancy bicategories. We will not explore this point, but in Definition~\ref{fancificationdef} we noted that these two  constructions can be promoted to apply to bicategory functors or to strict functors of strict bicategories, respectively. We can apply these constructions to transformations and modifications as well, as we describe below.

Let  $\B$ and $\C$  be bicategories. Let $(\B,\C)$ denote the strict bicategory of (weak) functors, transformations, and modifications from $\B$ to $\C$. If $\B$ and $\C$ are strict bicategories, let $[\B,\C]$ denote the strict bicategory of strict functors, strict transformations, and modifications from $\B$ to $\C$, with $i_{\B,\C}:[\B,\C]\ra(\B,\C)$ being the obvious canonical strict functor. There is a canonical strict isomorphism of fancy bicategories :
\begin{itemize}
\item $F^{\downarrow }_{\B,\C}:\llcorner (\B,\C)\lrcorner \stackrel{F \ra \llcorner F \lrcorner }{\lra} \{\llcorner \B \lrcorner, \llcorner \C \lrcorner\}$
\end{itemize} 
If $\B$ and $\C$ are strict bicategories, there is a canonical strict isomorphism of fancy bicategories:
\begin{itemize}
\item $F^{\uparrow }_{\B,\C}:\ulcorner [\B,\C] \urcorner \stackrel{F \ra \ulcorner F \urcorner }{\lra} \{\ulcorner \B \urcorner, \ulcorner \C \urcorner\}$
\end{itemize}

\begin{proposition}~ \label{fancyfacts}
\begin{itemize}
\item If $\B$ is a bicategory then a $1$-morphism $f$ of $\B$  an equivalence in $\llcorner \B \lrcorner$ if and only if it is an equivalence in $\B$. If $\B$ is strict, then $f$ is an equivalence in $\ulcorner \B \urcorner$ if and only if it is an \emph{isomorphism} in $\B$.
\item $\cfancya$ and $\sfancya$ are strictly fully faithful
\item If $F$ is a strictly fully faithful functor of fancy bicategories, so is $F^\circ$
\item A strict isomorphism of fancy bicategories is strictly fully faithful
\item The composition of two strictly fully faithful functors is strictly fully faithful
\item If $\B$ is a strict bicategory and $\C$ is a fancy bicategory then  $\presuper{\circ}{(\iota_\B)}:[\llcorner \B \lrcorner,\C]\ra[\ulcorner \B \urcorner,\C]$ is a strict isomorphism
\item  Let $\B$ be a bicategory. There is a strict isomorphism: 
\begin{align*}\fancyhom_{\llcorner \B \lrcorner} &\cong f^{\downarrow} \circ \llcorner \Hom_\B \lrcorner     \end{align*}
\item If $\B$ is strict, there is a strict isomorphism:
\begin{align*}
\fancyhom_{\ulcorner \B \urcorner} &\cong f^{\uparrow} \circ \ulcorner \Hom_\B \urcorner    
           \end{align*}
\item Let $F:\A\times \B \ra \C$ be a bicategory functor. There is a strict isomorphism:
$$\fancyadj_{\llcorner \A\lrcorner,\llcorner \B \lrcorner,\llcorner \C \lrcorner}(\llcorner F \lrcorner ) \cong F^{\downarrow}_{\B,\C} \circ \llcorner \Adj_{\A,\B,\C}(F)\lrcorner   $$
\item If $\A$, $\B$, $\C$, and $F$ are strict, there is a strict isomorphism:
$$\fancyadj_{\ulcorner \A\urcorner,\ulcorner \B \urcorner,\ulcorner \C \urcorner}(\ulcorner F \urcorner ) \cong F^{\uparrow}_{\B,\C} \circ \ulcorner \Adjst_{\A,\B,\C}(F)\urcorner $$
\item Let $F:\A \times \B \ra \C$ and $G: \C \ra \D$ be a functor of fancy bicategories. Then
$$\fancyadj_{\A,\B,\D}(G\circ F)=G^\circ \circ \fancyadj_{\A,\B,\C }(F)$$
\item Let $\B$ be a bicategory and $b$ an object in $\B$. Then $\lambda^{\llcorner \B \lrcorner}_b= \llcorner \lambda^\B_b\lrcorner.$ If $\B$ is strict then  $\lambda^{\ulcorner \B \urcorner}_b= \ulcorner \lambda^\B_b\urcorner.$ 
\end{itemize}
\end{proposition}
\begin{proof} These are individually easy to check and are left to the reader.
\end{proof}
\subsection{$2$-limits}
The categorical concept of limit has several generalizations to bicategory theory. Setting aside the various kinds of lax and colax limits, there are three natural kinds of $2$-limit, corresponding to three levels of strictness. The last of these to be defined was Kelly's pseudolimits, defined in \cite{Kel89} which also gives details on the other definitions and gives various earlier references. See also \cite{Lac10} for a more recent exposition and \cite{nla13} for a concise overview. In this section, we show that the fully strict and fully weak notions are both special cases of a general ``fancy 2-limit''.

For this discussion, we break from our usual practice of using calligraphic letters for $2$-dimensional categorical structures and use Roman letters for bicategories to distinguish them from fancy bicategories. The chief facts and isomorphisms used in our analysis come from Proposition~\ref{fancyfacts}, but we will also use some other more obvious facts and isomorphisms without mention, for instance the equality $\llcorner G \lrcorner \circ \llcorner F \lrcorner = \llcorner G \circ F\lrcorner$  for composable bicategory functors $F$ and $G$. 

Let $K$ and $D$ be bicategories, and $J: D \ra \cat{Cat}$ and $F:D \ra K$ be functors. A \emph{$J$-weighted $2$-limit of $F$} is an object $L\in K$ and a natural equivalence (i.e. an equivalence in the bicategory $(K^{\op},\cat{Cat})$)  between the two functors $K^{\op} \ra \cat{Cat}$ which send an $X \in K$ to the left and right side of the equation below:
$$K(X,L) \cong \Hom_{(D,\cat{Cat})}(J, K(X, F - ) )$$

If $K$ and $D$ are strict bicategories, and $J$ and $F$ are strict functors, a \emph{$J$-weighted strict $2$-limit of $F$} is an object $L\in K$ and a \emph{strict} natural \emph{isomorphism} in $[K^{\op},\cat{Cat}]$ between the functors sending $X\in K$  to the left and right side of:
$$K(X,L) \cong \Hom_{[D,\cat{Cat}]}(J, K(X, F - ) )$$
In the same case where $K$, $D$, $J$, and $F$ are strict, a \emph{$J$-weighted pseudo $2$-limit of $F$}
is an object $L\in K$ and a strict natural isomorphism in $[K^{\op},\cat{Cat}]$ between the functors sending $X\in K$  to the left and right side of:
$$K(X,L) \cong \Hom_{(D,\cat{Cat})}(J, K(X, F - ) )$$

\begin{definition}
Let $\K$ and $\D$ be fancy bicategories, and $J: \D \ra \cat{Cat}$ and $F:\D \ra \K$ be functors. A \emph{$J$-weighted fancy $2$-limit of $F$} is an object $L\in \K$ and a natural equivalence (i.e. an equivalence in the fancy bicategory $\{\K^{\op},\TCat\}$)  between the two functors $\K^{\op} \ra \TCat$ which send an $X \in \K$ to the left and right side of the equation below:
$$\K(X,L) \cong \Hom_{(\D,\TCat)}(J, \K(X, F - ) )$$
Explicitly, this is an equivalence $$\fancyhom_\K \circ (\id \times \lambda_L^\K)\cong \fancyhom_{\{\D,\TCat\}}\circ \left( \lambda_J^{\{\D,\TCat\}^\op} \times \fancyadj(\fancyhom_\K \circ (\id \times F)    )   \right)$$
\end{definition}

\begin{theorem} \label{gen2lim}
Let $K$ and $D$ be bicategories, and $J: D \ra \cat{Cat}$ and $F:D \ra K$ be functors. A \emph{$J$-weighted $2$-limit of $F$} is equivalent to a $f^{\downarrow} \circ \llcorner J \lrcorner$-weighted fancy $2$-limit of $\llcorner F \lrcorner: \llcorner D \lrcorner \ra \llcorner K \lrcorner$. 
\end{theorem}
\begin{proof}
 A  $f^{\downarrow} \circ \llcorner J \lrcorner$-weighted fancy $2$-limit of $\llcorner F \lrcorner$ is an object $L$ in $\llcorner K \lrcorner$ (or equivalently in $K$) and an equivalence:
  $$\fancyhom_{\llcorner K \lrcorner }  \circ (\id \times \lambda_L^{\llcorner K \lrcorner})\cong \fancyhom_{\{\llcorner D \lrcorner \ , \ \TCat\}}\circ \left( \lambda_{f^{\downarrow} \circ \llcorner J \lrcorner}^{\{\llcorner D\lrcorner\ ,\ \TCat\}^\op} \times \fancyadj(\fancyhom_{\llcorner K\lrcorner } \circ (\id \times \llcorner F\lrcorner )    )   \right).$$

Starting with the left hand side, we have strict isomorphisms:
\begin{align*}
\fancyhom_{\llcorner K \lrcorner }  \circ (\id \times \lambda_L^{\llcorner K \lrcorner})
\cong&f^{\downarrow}\circ \llcorner\Hom_{ K  }\lrcorner  \circ (\id \times\llcorner \lambda_L^{ K} \lrcorner)\\
\cong&f^{\downarrow}\circ \llcorner\Hom_{ K  }  \circ (\id \times \lambda_L^{ K} )\lrcorner.
\end{align*}
For the right hand side, we have strict isomorphisms
\begin{align*}
 &\fancyhom_{\{\llcorner D \lrcorner \ , \ \TCat\}}\circ \left( \lambda_{f^{\downarrow} \circ \llcorner J \lrcorner}^{\{\llcorner D\lrcorner\ ,\ \TCat\}^\op} \times \fancyadj(\fancyhom_{\llcorner K\lrcorner } \circ (\id \times \llcorner F\lrcorner )    )   \right) \\
\cong& \fancyhom_{\{\llcorner D \lrcorner \ , \ \TCat\}}\circ \left( (({f^{\downarrow}}^\circ)^\op \circ\lambda_{ \llcorner J \lrcorner}^{\{\llcorner D\lrcorner\ ,\ \llcorner \cat{Cat} \lrcorner \}^\op} )\times \fancyadj(f^{\downarrow}\circ\llcorner\Hom_{ K  }\lrcorner  \circ (\id \times \llcorner F\lrcorner )    )   \right) \\
\cong& \fancyhom_{\{\llcorner D \lrcorner \ , \ \TCat\}}\circ \left( (({f^{\downarrow}}^\circ)^\op \circ (F^{\downarrow}_{D,\cat{Cat}})^\op \circ \lambda_{ J }^{\llcorner (D \ ,\  \cat{Cat})\lrcorner ^\op } )\times ({f^{\downarrow}}^\circ \circ \fancyadj(\llcorner\Hom_{ K  }\circ (\id \times  F )  \lrcorner  )  ) \right)\\
\cong& \fancyhom_{\{\llcorner D \lrcorner \ , \ \TCat\}}\circ  \\  &\quad \quad \quad ( ({f^{\downarrow}}^\circ \circ F^{\downarrow}_{D,\cat{Cat}})^\op \circ \llcorner\lambda_{ J }^{ (D \ ,\  \cat{Cat}) ^\op }\lrcorner )\times ({f^{\downarrow}}^\circ \circ F^{\downarrow}_{D,\cat{Cat}} \circ \llcorner\Adj(\Hom_{ K  }\circ (\id \times  F )   )\lrcorner   ) \\
\cong& \fancyhom_{\{\llcorner D \lrcorner \ , \ \TCat\}}\circ  \\  &\quad \quad \quad \left( ({f^{\downarrow}}^\circ \circ F^{\downarrow}_{D,\cat{Cat}})^\op \times ({f^{\downarrow}}^\circ \circ F^{\downarrow}_{D,\cat{Cat}}) \right) \circ \llcorner\lambda_{ J }^{ (D  , \cat{Cat}) ^\op }\times \Adj(\Hom_{ K  }\circ (\id \times  F )   )\lrcorner   \\
\end{align*}
Since ${f^{\downarrow}}^\circ \circ F^{\downarrow}_{D,\cat{Cat}}$ is strictly fully faithful, we continue with a strict isomorphism:
\begin{align*}
\cong& \fancyhom_{\llcorner( D  ,  \cat{Cat})\lrcorner}  \circ \llcorner\lambda_{ J }^{ (D  ,  \cat{Cat}) ^\op }\times \Adj(\Hom_{ K  }\circ (\id \times  F )   )\lrcorner   \\
\cong&f^{\downarrow}\circ \llcorner\Hom_{( D  ,  \cat{Cat})} \lrcorner \circ \llcorner\lambda_{ J }^{ (D  ,  \cat{Cat}) ^\op }\times \Adj(\Hom_{ K  }\circ (\id \times  F )   )\lrcorner \\
=&f^{\downarrow}\circ \llcorner\Hom_{( D  ,  \cat{Cat})}  \circ \lambda_{ J }^{ (D  ,  \cat{Cat}) ^\op }\times \Adj(\Hom_{ K  }\circ (\id \times  F )   )\lrcorner. 
\end{align*}
So giving our   fancy $2$-limit  is the same as giving an equivalence in $\{\llcorner K\lrcorner ,\TCat\}$
$$f^{\downarrow}\circ \llcorner\Hom_{ K  }  \circ (\id \times \lambda_L^{ K} )\lrcorner \cong f^{\downarrow}\circ \llcorner\Hom_{( D  ,  \cat{Cat})}  \circ \lambda_{ J }^{ (D , \cat{Cat}) ^\op }\times \Adj(\Hom_{ K  }\circ (\id \times  F )   )\lrcorner. $$
Since $f^{\downarrow}$ is fully faithful, this the same as giving an equivalence in $\{\llcorner\K \lrcorner, \llcorner\cat{Cat}\lrcorner\}$
$$ \llcorner\Hom_{ K  }  \circ (\id \times \lambda_L^{ K} )\lrcorner \cong \llcorner\Hom_{( D  ,  \cat{Cat})}  \circ \lambda_{ J }^{ (D ,  \cat{Cat}) ^\op }\times \Adj(\Hom_{ K  }\circ (\id \times  F )   )\lrcorner$$
Applying the isomorphism $F^{\downarrow}_{D, \cat{Cat}}$ this is the same as giving an equivalence in the fancy bicategory $\llcorner(D, \cat{Cat})\lrcorner$ between the functors:
$$\Hom_{ K  }  \circ (\id \times \lambda_L^{ K} )\cong \Hom_{( D  ,  \cat{Cat})}  \circ \lambda_{ J }^{ (D ,  \cat{Cat}) ^\op }\times \Adj(\Hom_{ K  }\circ (\id \times  F )   )$$
This is the same as giving an equivalence between these same functors in the \emph{bicategory} $(D,\cat{Cat})$ which is exactly the data needed to exhibit $L$ as a $J$-weighted $2$-limit of $F$.\end{proof}
\begin{theorem}
Let $K$ and $D$ be strict bicategories, and $J: D \ra \cat{Cat}$ and$F:D \ra K$ be strict  functors. A \emph{$J$-weighted strict $2$-limit of $F$} is equivalent to a $f^{\uparrow} \circ \ulcorner J \urcorner$-weighted fancy $2$-limit of $\ulcorner F \urcorner: \ulcorner D \urcorner \ra \ulcorner K \urcorner$. 
\end{theorem}
\begin{proof}
 A  $f^{\uparrow} \circ \ulcorner J \urcorner$-weighted fancy $2$-limit of $\ulcorner F \urcorner$ is an object $L$ in $\ulcorner K \urcorner$ (or equivalently in $K$) and an equivalence:
  $$\fancyhom_{\ulcorner K \urcorner }  \circ (\id \times \lambda_L^{\ulcorner K \urcorner})\cong \fancyhom_{\{\ulcorner D \urcorner \ , \ \TCat\}}\circ \left( \lambda_{f^{\uparrow} \circ \ulcorner J \urcorner}^{\{\ulcorner D\urcorner\ ,\ \TCat\}^\op} \times \fancyadj(\fancyhom_{\ulcorner K\urcorner } \circ (\id \times \ulcorner F\urcorner )    )   \right).$$
A similar argument to the one used in the proof of Theorem~\ref{gen2lim} shows this is the same an equivalence in $\{\ulcorner K\urcorner ,\ulcorner \cat{Cat} \urcorner \}$ between:
$$\ulcorner\Hom_{ K  }   (\id \times \lambda_L^{ K} )\urcorner \cong  \ulcorner\Hom_{[ D  ,  \cat{Cat}]}  \circ \lambda_{ J }^{ [D , \cat{Cat}]^\op }\times \Adj(\Hom_{ K  }\circ (\id \times  F )   )\urcorner. $$
Applying the isomorphism $F^{\downarrow}_{D, \cat{Cat}}$ this is the same as giving an equivalence in the fancy bicategory $\ulcorner[D, \cat{Cat}]\urcorner$ between the functors:
$$\Hom_{ K  }  \circ (\id \times \lambda_L^{ K} )\cong \Hom_{[ D  ,  \cat{Cat}]}  \circ \lambda_{ J }^{ [D ,  \cat{Cat}] ^\op }\times \Adjst(\Hom_{ K  }\circ (\id \times  F )   )$$
This is the same as a (strict) \emph{isomorphism} between the same functors in the bicategory $[D, \cat{Cat}]$, which is the data needed to exhibit $L$ as a strict $J$-weighted $2$-limit of $F$.
\end{proof}
\printbibliography
\end{document}